\documentclass[a4paper,10pt,reqno]{amsart}
\usepackage{amsmath,amsfonts,amsthm,amssymb,xcolor}
\usepackage[T1]{fontenc}
\usepackage{csquotes}
\usepackage{graphicx}
\usepackage{pstricks}
\usepackage{lmodern}
\usepackage{calc}
\usepackage{mathabx}
\usepackage{mathrsfs}
\usepackage{ulem}
\usepackage{hyperref}
\usepackage{stmaryrd}

\usepackage{cancel}

\usepackage{scalerel}

\newcommand{\pe}{\mathbin{\scaleobj{0.7}{\tikz \draw (0,0) node[shape=circle,draw,inner sep=0pt,minimum size=8.5pt] {\footnotesize $=$};}}}

\newcommand{\pl}{\mathbin{\scaleobj{0.7}{\tikz \draw (0,0) node[shape=circle,draw,inner sep=0pt,minimum size=8.5pt] {\footnotesize $<$};}}}
\newcommand{\pg}{\mathbin{\scaleobj{0.7}{\tikz \draw (0,0) node[shape=circle,draw,inner sep=0pt,minimum size=8.5pt] {\footnotesize $>$};}}}
\newcommand{\ple}{\mathbin{\scaleobj{0.7}{\tikz \draw (0,0) node[shape=circle,draw,inner sep=0pt,minimum size=8.5pt] {\footnotesize $\leqslant$};}}}
\newcommand{\pge}{\mathbin{\scaleobj{0.7}{\tikz \draw (0,0) node[shape=circle,draw,inner sep=0pt,minimum size=8.5pt] {\footnotesize $\geqslant$};}}}

\usepackage{accents}

\usepackage{enumerate}

\numberwithin{equation}{section}
\usepackage{tikz}

 \usepackage{caption} 
\newcommand{\vtila}{\widehat{v}^{(\la)}}
\newcommand{\wtila}{\widehat{w}^{(\la)}}

\newcommand{\vla}{v^{(\la)}}
\newcommand{\wla}{w^{(\la)}}

\newcommand{\id}{\mbox{Id}}

\newcommand{\1}{\mathbf{1}}

 \newcommand{\ov}{\overline}
\newcommand{\ct}{\mathcal T}

\newcommand{\eps}{\varepsilon}

\newcommand{\dis}{\displaystyle}

\newcommand{\B}{\mathcal{B}}

\newcommand{\cq}{\mathcal{Q}}

\newcommand{\C}{\mathbb C}

\newcommand{\R}{\mathbb R}
\newcommand{\N}{\mathbb N}

\newcommand{\T}{\mathbb T}

\newcommand{\cb}{\mathcal B}
\newcommand{\cac}{\mathcal C}
\newcommand{\cd}{\mathcal D}

\newcommand{\cf}{\mathcal F}

\newcommand{\ci}{\mathcal I}
\newcommand{\cj}{\mathcal J}

\newcommand{\cl}{\mathcal L}
\newcommand{\cm}{\mathcal M}

\newcommand{\cs}{\mathcal S}
\newcommand{\cw}{\mathcal W}

\newcommand{\calt}{\mathcal T}

\newcommand{\frakc}{\mathfrak{c}}

\newcommand{\scret}{\mathscr{T}}

\newcommand{\al}{\alpha}

\newcommand{\ep}{\varepsilon}

\newcommand{\ga}{\gamma}
\newcommand{\gga}{\Gamma}
\newcommand{\ka}{\kappa}
\newcommand{\la}{\lambda}

\newcommand{\om}{\omega}

\newcommand{\si}{\sigma}

\newcommand{\vp}{\varphi}

\newcommand{\tn}{|\!|\!|}

\newtheorem{theorem}{Theorem}[section]
\newtheorem*{acknowledgements}{Acknowledgements}

\newtheorem{corollary}[theorem]{Corollary}

\newtheorem{definition}[theorem]{Definition}

\newtheorem{lemma}[theorem]{Lemma}
\newtheorem{notation}[theorem]{Notation}

\newtheorem{proposition}[theorem]{Proposition}

\theoremstyle{remark}
\newtheorem{remark}[theorem]{Remark}

\pgfdeclareshape{crosscircle}
{
  \inheritsavedanchors[from=circle] 
  \inheritanchorborder[from=circle]
  \inheritanchor[from=circle]{north}
  \inheritanchor[from=circle]{north west}
  \inheritanchor[from=circle]{north east}
  \inheritanchor[from=circle]{center}
  \inheritanchor[from=circle]{west}
  \inheritanchor[from=circle]{east}
  \inheritanchor[from=circle]{mid}
  \inheritanchor[from=circle]{mid west}
  \inheritanchor[from=circle]{mid east}
  \inheritanchor[from=circle]{base}
  \inheritanchor[from=circle]{base west}
  \inheritanchor[from=circle]{base east}
  \inheritanchor[from=circle]{south}
  \inheritanchor[from=circle]{south west}
  \inheritanchor[from=circle]{south east}
  \inheritbackgroundpath[from=circle]
  \foregroundpath{
    \centerpoint%
    \pgf@xc=\pgf@x%
    \pgf@yc=\pgf@y%
    \pgfutil@tempdima=\radius%
    \pgfmathsetlength{\pgf@xb}{\pgfkeysvalueof{/pgf/outer xsep}}%
    \pgfmathsetlength{\pgf@yb}{\pgfkeysvalueof{/pgf/outer ysep}}%
    \ifdim\pgf@xb<\pgf@yb%
      \advance\pgfutil@tempdima by-\pgf@yb%
    \else%
      \advance\pgfutil@tempdima by-\pgf@xb%
    \fi%
    \pgfpathmoveto{\pgfpointadd{\pgfqpoint{\pgf@xc}{\pgf@yc}}{\pgfqpoint{-0.707107\pgfutil@tempdima}{-0.707107\pgfutil@tempdima}}}
    \pgfpathlineto{\pgfpointadd{\pgfqpoint{\pgf@xc}{\pgf@yc}}{\pgfqpoint{0.707107\pgfutil@tempdima}{0.707107\pgfutil@tempdima}}}
    \pgfpathmoveto{\pgfpointadd{\pgfqpoint{\pgf@xc}{\pgf@yc}}{\pgfqpoint{-0.707107\pgfutil@tempdima}{0.707107\pgfutil@tempdima}}}
    \pgfpathlineto{\pgfpointadd{\pgfqpoint{\pgf@xc}{\pgf@yc}}{\pgfqpoint{0.707107\pgfutil@tempdima}{-0.707107\pgfutil@tempdima}}}
  }
}
\makeatother

\definecolor{gr}{rgb}   {0.,   0.69,   0.23 }
\definecolor{bl}{rgb}   {0.,   0.5,   1. }
\definecolor{mg}{rgb}   {0.85,  0.,    0.85}
\definecolor{yl}{rgb}   {0.8,  0.7,   0.}
\definecolor{or}{rgb}  {0.7,0.2,0.2}
\definecolor{marron}{rgb}{0.64,0.16,0.16}

\colorlet{symbols}{black!90!black}
\colorlet{symbolsb}{black!90!black}
\colorlet{testcolor}{green!60!black}

\usetikzlibrary{shapes.misc}
\usetikzlibrary{shapes.symbols}
\usetikzlibrary{decorations}
\usetikzlibrary{decorations.markings}

\def\drawx{\draw[-,solid] (-3pt,-3pt) -- (3pt,3pt);\draw[-,solid] (-3pt,3pt) -- (3pt,-3pt);}
\tikzset{
	root/.style={circle,fill=testcolor,inner sep=0pt, minimum size=2mm},
	dot/.style={circle,fill=black,inner sep=0pt, minimum size=1mm},
	var/.style={circle,fill=black!10,draw=black,inner sep=0pt, minimum size=2mm},
	dotred/.style={circle,fill=black!50,inner sep=0pt, minimum size=2mm},
	generic/.style={semithick,shorten >=1pt,shorten <=1pt},
	dist/.style={ultra thick,draw=testcolor,shorten >=1pt,shorten <=1pt},
	testfcn/.style={ultra thick,testcolor,shorten >=1pt,shorten <=1pt,<-},
	testfcnx/.style={ultra thick,testcolor,shorten >=1pt,shorten <=1pt,<-,
		postaction={decorate,decoration={markings,mark=at position 0.6 with {\drawx}}}},
	kprime/.style={semithick,shorten >=1pt,shorten <=1pt,densely dashed,->},
	kprimex/.style={semithick,shorten >=1pt,shorten <=1pt,densely dashed,->,
		postaction={decorate,decoration={markings,mark=at position 0.4 with {\drawx}}}},
	kernel/.style={semithick,shorten >=1pt,shorten <=1pt,->},
	multx/.style={shorten >=1pt,shorten <=1pt,
		postaction={decorate,decoration={markings,mark=at position 0.5 with {\drawx}}}},
	kernelx/.style={semithick,shorten >=1pt,shorten <=1pt,->,
		postaction={decorate,decoration={markings,mark=at position 0.4 with {\drawx}}}},
	kernel1/.style={->,semithick,shorten >=1pt,shorten <=1pt,postaction={decorate,decoration={markings,mark=at position 0.45 with {\draw[-] (0,-0.1) -- (0,0.1);}}}},
	kernel2/.style={->,semithick,shorten >=1pt,shorten <=1pt,postaction={decorate,decoration={markings,mark=at position 0.45 with {\draw[-] (0.05,-0.1) -- (0.05,0.1);\draw[-] (-0.05,-0.1) -- (-0.05,0.1);}}}},
	kernelBig/.style={semithick,shorten >=1pt,shorten <=1pt,decorate, decoration={zigzag,amplitude=1.5pt,segment length = 3pt,pre length=2pt,post length=2pt}},
	rho/.style={dotted,semithick,shorten >=1pt,shorten <=1pt},
	renorm/.style={shape=circle,fill=white,inner sep=1pt},
	res/.style={circle,draw=symbols,inner sep=0pt,minimum size=1.2mm},
	labl/.style={shape=rectangle,fill=white,inner sep=1pt},
	xi/.style={circle,fill=symbols!10,draw=symbols,inner sep=0pt,minimum size=1.2mm},
	xiblack/.style={circle,fill=symbolsb,draw=symbolsb,inner sep=0pt,minimum size=1.2mm},
	xix/.style={crosscircle,fill=symbols!10,draw=symbols,inner sep=0pt,minimum size=1.2mm},
	xib/.style={circle,fill=symbols!10,draw=symbols,inner sep=0pt,minimum size=1.6mm},
	xibx/.style={crosscircle,fill=symbols!10,draw=symbols,inner sep=0pt,minimum size=1.6mm},
	not/.style={circle,fill=symbols,draw=symbols,inner sep=0pt,minimum size=0.5mm},
	notblack/.style={circle,fill=symbolsb,draw=symbolsb,inner sep=0pt,minimum size=0.5mm},
	>=stealth,
	}
\makeatletter
\def\DeclareSymbol#1#2#3{\expandafter\gdef\csname MH@symb@#1\endcsname{\tikz[baseline=#2,scale=0.15,draw=symbols]{#3}}\expandafter\gdef\csname MH@symb@#1s\endcsname{\scalebox{0.7}{\tikz[baseline=#2,scale=0.15,draw=symbols]{#3}}}}
\def\<#1>{\csname MH@symb@#1\endcsname}
\makeatother

\DeclareSymbol{circle}{0.5}{\draw (0,0.7) node[xi] {};}
\DeclareSymbol{line}{0.5}{\draw (0,0.2) node[not] {} -- (0,1.4) node[not] {};}

\DeclareSymbol{Psi}{0.5}{\draw (0,0) node[not] {} -- (0,1.5) node[xi] {};}
\DeclareSymbol{Psiblack}{0.5}{\draw (0,0) node[notblack] {} -- (0,1.5) node[xiblack] {};}
\DeclareSymbol{Psi2}{0.5}{\draw (-0.6,1.5) node[xi] {} -- (0,0) node[not] {} -- (0.6,1.5) node[xi] {};}
\DeclareSymbol{Psi3}{0.5}{\draw (-1,1.5) node[xi] {} -- (0,0) node[not] {}; \draw (0,1.5) node[xi] {} -- (0,0) node[not] {};\draw (1,1.5) node[xi] {} -- (0,0) node[not] {}}
\DeclareSymbol{IPsi2}{0}{\draw (0,1) -- (0.8,2.2) node[xi] {};\draw (0,-0.25) node[not] {} -- (0,1) node[not] {} -- (-0.8,2.2) node[xi] {};}
\DeclareSymbol{IPsi3}{0}{\draw (0,1) -- (1,2.2) node[xi] {};\draw (0,1) -- (0,2.2) node[xi] {};\draw (0,-0.25) node[not] {} -- (0,1) node[not] {} -- (-1,2.2) node[xi] {};}
\DeclareSymbol{PsiIPsi2}{0}{\draw (0,1) -- (0.8,2.2) node[xi] {};\draw (0,-0.25) node[not] {} -- (0,1) node[not] {} -- (-0.8,2.2) node[xi] {};\draw (0,-0.25) node[not] {} -- (0.8,1) node[xi] {};}

\DeclareSymbol{Psi2IPsi3}{0}{\draw (0,1) -- (1,2.2) node[xi] {};\draw (0,1) -- (0,2.2) node[xi] {};\draw (0,-0.25) node[not] {} -- (0,1) node[not] {} -- (-1,2.2) node[xi] {};\draw (0,-0.25) node[not] {} -- (1,0.75) node[xi] {};\draw (0,-0.25) node[res] {} -- (-1,0.75) node[xi] {};}

\DeclareSymbol{PsiIPsi3}{0}{\draw (0,1) -- (1,2.2) node[xi] {};\draw (0,1) -- (0,2.2) node[xi] {};\draw (0,-0.25) node[not] {} -- (0,1) node[not] {} -- (-1,2.2) node[xi] {};\draw (0,-0.25) node[not] {} -- (-1,0.75) node[xi] {};\draw (0,-0.25) node[res] {};}

\DeclareSymbol{PsiIPsi3nr}{0}{\draw (0,1) -- (1,2.2) node[xi] {};\draw (0,1) -- (0,2.2) node[xi] {};\draw (0,-0.25) node[not] {} -- (0,1) node[not] {} -- (-1,2.2) node[xi] {};\draw (0,-0.25) node[not] {} -- (-1,0.75) node[xi] {};\draw (0,-0.25) node[not] {};}

\DeclareSymbol{Psi2IPsi2}{0}{\draw (0,1) -- (0.8,2) node[xi] {};\draw (0,-0.25) node[not] {} -- (0,1) node[not] {} -- (-0.8,2) node[xi] {};\draw (0,-0.25) node[not] {} -- (-0.9,0.75) node[xi] {};\draw (0,-0.25) node[not] {} -- (0.9,0.75) node[xi] {};\draw (0,-0.25) node[res] {};}

\DeclareSymbol{Psi2IPsi2nr}{0}{\draw (0,1) -- (0.8,2) node[xi] {};\draw (0,-0.25) node[not] {} -- (0,1) node[not] {} -- (-0.8,2) node[xi] {};\draw (0,-0.25) node[not] {} -- (-0.9,0.75) node[xi] {};\draw (0,-0.25) node[not] {} -- (0.9,0.75) node[xi] {};\draw (0,-0.25) node[not] {};}

\DeclareSymbol{Psi2IPsi3nr}{0}{\draw (0,1) -- (1,2.2) node[xi] {};\draw (0,1) -- (0,2.2) node[xi] {};\draw (0,-0.25) node[not] {} -- (0,1) node[not] {} -- (-1,2.2) node[xi] {};\draw (0,-0.25) node[not] {} -- (1,0.75) node[xi] {};\draw (0,-0.25) node[not] {} -- (-1,0.75) node[xi] {};}

\date{\today}

\title[On the parabolic $\Phi_3^4$ model for the harmonic oscillator II]{On the parabolic $\Phi_3^4$ model for the harmonic oscillator II: global existence and invariant measures}

\begin{document}



\begin{abstract}
We establish an a priori bound for the dynamical parabolic $\Phi_3^4$ model with harmonic potential. This bound yields the global well-posedness of the equation and, via the Krylov-Bogoliubov method, the existence of an invariant measure, shown to be non-Gaussian. The argument builds on the strategy developed by Mourrat and Weber for the periodic $\Phi_3^4$ model, with substantial modifications to handle the non-compact geometry of $\R^3$ and the spectral framework imposed by the harmonic oscillator. We further prove that this measure is unique in the small-coupling regime. 
\end{abstract}

\maketitle

\begin{center}
{\large
Aur\'elien Deya\footnote{Universit\'e de Lorraine, CNRS, IECL, F-54000 Nancy, France. \\Email: {\tt aurelien.deya@univ-lorraine.fr}}, Reika Fukuizumi\footnote{Department of Mathematics,
School of fundamental science and engineering, Waseda University, 3-4-1, Okubo, Shinjuku-ku, Tokyo, Japan. \\Email: {\tt fukuizumi@waseda.jp}}, and Laurent Thomann\footnote{Universit\'e de Lorraine, CNRS, IECL, F-54000 Nancy, France. \\Email: {\tt laurent.thomann@univ-lorraine.fr}}}
\end{center}

 \tableofcontents
 

\section{Introduction}\label{sec:introd}
\subsection{Motivation and main result}

In this article, we continue our investigation initiated in \cite{DFT} regarding the parabolic~$\Phi_3^4$ model with a harmonic potential, which, after renormalization, can formally be written as
\begin{equation} \label{eq:intro}
\left\{
\begin{aligned}
&  \partial_t X + HX= -\la X^3 +\infty \cdot X+ \xi, \quad t>0, \quad x \in \R^3, \\
& X_0 = u,
\end{aligned}
\right.
\end{equation}
with $\la>0$, and where 
$$H:=-\Delta_{\R^3} +|x|^2$$
stands for the harmonic oscillator on $\R^3$, and $\xi$ denotes a space-time white noise defined 
on a complete filtered probability space $(\Omega,\cf,\mathbb{P})$, with (formal) covariance
$$\mathbb{E}\big[ \xi(t,x) \xi(s,y) \big] =\delta_{t-s} \delta_{x-y}.$$
 Equation~\eqref{eq:intro} can thus be seen as the gradient flow of the energy functional associated with the~$\Phi^4_3$ measure with harmonic confinement, or equivalently, in the language of quantum field theory, as the {\it stochastic quantization} equation of this measure.

\medskip

This equation belongs to the broad class of {\it singular stochastic partial differential equations}, a subject that has undergone spectacular development over the past decade. The term ``singular'' reflects the fact that the noise drives the solution into a space of distributions in which the cubic nonlinearity $-X^3$ has no classical meaning, so that a renormalization procedure is needed to give the equation a rigorous interpretation. The foundational contributions of Hairer~\cite{Hai14}, who introduced the theory of {\it regularity structures}, and of Gubinelli, Imkeller and Perkowski~\cite{GIP}, who developed the {\it paracontrolled calculus}, have provided systematic frameworks for treating such equations in the parabolic setting. 

\smallskip

Prototypical examples include the KPZ equation in one space dimension, the dynamical $\Phi^4_2$ and $\Phi^4_3$ models arising in Euclidean quantum field theory, and various singular parabolic or dispersive equations driven by space-time white noise. A common feature of these equations is the need for renormalization: formally divergent polynomial expressions in the solution must be replaced by finite Wick-ordered counterparts, whose precise form depends on the ambient function spaces and on the regularization scheme. In the present work, the operator $H = -\Delta + |x|^2$ replaces the flat Laplacian $-\Delta$, and the spectral theory of the harmonic oscillator, in particular the Hermite function basis and the associated Besov spaces~\cite{FI}, plays the role played by Fourier analysis in the periodic setting. 

\medskip

The study of equation~\eqref{eq:intro} is motivated by the physics of {\it Bose-Einstein condensates} (BEC). At positive temperature, thermal fluctuations in a BEC confined by a harmonic trap can be modeled by the {\it stochastic Gross-Pitaevskii equation} (SGPE)~\cite{GAF, GD}, a complex Ginzburg-Landau equation driven by space-time white noise:
\begin{equation}\label{eq:SGPE}
  (\beta - i\alpha)\,(\partial_t \psi + H\psi) = -\lambda\,|\psi|^2 \psi +\xi_\C,
  \quad t>0, \quad x\in\R^d,
\end{equation}
where $\alpha,\beta > 0$ control the relative strength of the dispersive and dissipative parts, and $\xi_\C$ is a space-time white noise.  Equation~\eqref{eq:intro} somehow corresponds to a real version of \eqref{eq:SGPE}, and our study can thus be considered as a first step towards a full understanding of the SGPE dynamics.

\medskip

The mathematical analysis of the SGPE was initiated by de~Bouard, Debussche and the second author ~\cite{dBDF18} in one dimension. The more challenging problem of space-time white noise forcing in dimension two was treated by the same authors~\cite{dBDF23}, who established global well-posedness via an inhomogeneous Wick renormalization and proved the existence of an invariant measure. The three-dimensional case, which is the focus of the present work and of our companion paper~\cite{DFT}, is substantially more singular: even restricting ourselves to the simpler case $\alpha=0$, the solution lives in a space of distributions of negative Sobolev regularity, and the renormalization requires the subtraction of two divergent constants rather than one.

In~\cite{DFT} we constructed the stochastic diagrams associated to~\eqref{eq:intro} and obtained local well-posedness results. The purpose of the present study is to show that the problem is globally well-posed with rougher initial conditions, and to establish a priori bounds on the solutions that are uniform in the initial datum. This in turn will allow us to apply the Krylov-Bogoliubov method to prove the existence of an invariant measure for the dynamics of~\eqref{eq:intro}.

As noted above, we restrict our discussion to $\alpha =0$ and thus focus on the formulation of the problem in \eqref{eq:intro}. Nevertheless, in light of previous results on the torus (see \cite{HIN}), we expect similar a priori estimates to hold when $\alpha \neq 0$ and $\beta$ is large enough.

\

Beyond existence of invariant measures, we address the question of uniqueness: for a sufficiently small coupling constant $\la>0$, we prove that the rescaled equation~\eqref{eq-la} below admits a unique invariant measure. We also establish that this measure is non-Gaussian. From the physical standpoint, the invariant measure would describe the statistical equilibrium of the system governed by \eqref{eq:SGPE}. Its uniqueness then would ensure that the long-time behaviour of the condensate is determined by the noise alone, regardless of the initial configuration (see \cite{dBDF18}).

\medskip

To place our results in perspective, let us recall the state of the art regarding the ergodic properties of the {\it standard} $\Phi^4_3$ equation, in which the harmonic oscillator $H$ is replaced by the Laplacian $-\Delta$:
\begin{equation}\label{eq:Phi43std}
  \partial_t X - \Delta X = -\la X^3 + \infty \cdot X + \xi, \quad t>0, \quad x \in \T^3 \text{ or } \R^3.
\end{equation}
On the torus $\T^3$, the {\it coming down from infinity} property -- a uniform-in-time bound on solutions independent of the initial datum -- was first established by Mourrat-Weber~\cite{MW} via paracontrolled calculus, and subsequently re-derived using regularity structures by Moinat-Weber~\cite{MoW}. This a priori bound is the cornerstone of the long-time analysis. Combined with the strong Feller property proved by Hairer-Mattingly~\cite{HM18} and with the support theorem established by Hairer-Sch{\"o}nbauer~\cite{HS22}, it yields the existence and uniqueness of an invariant measure as well as the exponential ergodicity of the semigroup (see~\cite[Theorem 2.2]{DHYZ25}); an analogous result for the $\Phi^4_2$ equation on the torus had been obtained earlier by Tsatsoulis-Weber~\cite{TW18}. In \cite{BDFT1, BDFT2}, Bailleul,  Dang, Ferdinand and  T\^o recently generalized these results to a $\Phi^4_3$ model 
on a closed manifold.

\smallskip

Extending all these previous considerations to infinite volume is substantially harder: compactness arguments are no longer available, the noise needs to be set up on weighted spaces, and the known strategies towards the strong Feller property fail on $\R^3$, to mention but a few obstacles. Nevertheless, global well-posedness and the construction of the $\Phi^4_3$ measure on $\R^3$ were achieved by Gubinelli-Hofmanov\'a~\cite{GH}. Very recently, Duch, Hairer, Yi and Zhao~\cite{DHYZ25} proved the global well-posedness of the $\Phi^4_3$ dynamic on $\R^3$ in a suitable weighted Besov space, and showed that at {\it high temperature}, {\it i.e.}, when $\la >0$ in \eqref{eq:Phi43std} is small enough, any two solutions driven by the same noise converge to each other exponentially fast. 
This allows them to characterize the infinite-volume \(\Phi^4_3\) measure as the unique invariant measure of the dynamics and, through stochastic quantization, to verify the Osterwalder-Schrader axioms, including invariance under Euclidean isometries and exponential decay of correlations. Their argument, however, does not cover the entire high-temperature regime up to the phase transition. In the low-temperature regime, one expects uniqueness to fail, with multiple coexisting invariant measures, in line with the phase coexistence phenomenon established both for the Gibbs measure and in the discrete setting (see \cite{FSS,CGW} and the references therein).

\

Our approach to the coming-down-from-infinity property and to the construction of an invariant measure for \eqref{eq:intro} broadly follows the strategy developed by Mourrat and Weber~\cite{MW} for the parabolic $\Phi_3^4$ model on the torus. The central point of the argument is to establish Lyapunov functionals on the fluctuation part of the solution, using $L^p$-norms for very large $p\geq 1$. In passing, one may draw a parallel with the strategy adopted by Burq-Tzvetkov~\cite{BurqTzv} to obtain global solutions for the wave equation with random initial data. See also~\cite{BKTV} and references therein for the Schr\"odinger case.

Adapting the strategy  of \cite{MW} to the framework dictated by the harmonic oscillator requires several preliminary steps and substantial modifications:

\smallskip

\noindent
$(i)$ We first refine the construction and the topology of the stochastic diagrams beyond what was done in our companion paper~\cite{DFT}. In particular, the regularity of the diagrams must now be controlled in stronger norms, suited to a global-in-time analysis. This is achieved through new estimates combined with an interpolation argument. 

\smallskip

\noindent
$(ii)$ Within this enhanced framework, the passage from the periodic setting $\T^3$ to the Euclidean space $\R^3$ with harmonic confinement still requires several non-trivial adaptations. In particular, the lack of the embedding $L^q(\T^3)\subset L^p(\T^3)$ (for $p\leq q$) in the whole-space framework must be compensated by suitable inclusions in harmonic Besov spaces (see Lemma~\ref{lem-bis}).

\smallskip

\noindent
$(iii)$ We further improve the bounds obtained in \cite{MW} by establishing estimates in more general Besov topologies (see Theorem~\ref{theo:impr-bou}), which are particularly useful in the remainder of the analysis.

\smallskip

\noindent
$(iv)$ Once estimates uniform in the initial condition, and polynomial bounds with respect to the stochastic diagrams, are available, the application of the Krylov-Bogoliubov criterion proceeds in a way entirely analogous to the torus case and yields the desired invariant measure. We then prove that this measure is non-Gaussian (see Theorem~\ref{theo:non-gauss}), using cumulant techniques inspired by Malliavin calculus (see~\cite{NP, nualart-book}).

\

Finally, the other main contribution of this work is the proof of uniqueness of the invariant measure when the coupling constant $\la >0$ in \eqref{eq:intro} is small enough. Following the ideas of~\cite{DHYZ25}, the argument here consists in showing that two solutions with different initial conditions converge exponentially fast to one another in $L^\infty$. This is achieved by linearizing the dynamics around their difference and performing a careful analysis of the corresponding linear equation.
\

\

The result of all these investigations can be summarized as follows. For $\alpha \in \R$, we denote by $\B^\alpha_{\infty, \infty}(\R^3)$ the Besov space based on the harmonic oscillator (see Section \ref{sect-Besov} for a precise definition).

\begin{theorem}\label{thm:main}

Let $(\xi^{(n)})$ be the smooth regularization of $\xi$ given by \eqref{regu-noise}. For every $\la>0$, there exists a sequence $(\frakc^{(\la),(n)})$ of \textit{deterministic} functions on $\R^3$ such that the following assertions hold:

\smallskip

\noindent
 $(i)$  Uniformly in $|x| \ll 2^{\frac{n}2}$, we have
\begin{equation}\label{expC}
\frakc^{(\la),(n)}(x)=   \frac{3\la \cdot 2^{n/2}}{8\sqrt{2}\,\pi^{3/2}}- \frac{3 \lambda |x|}{8\pi} +\frac{3 \lambda^2}{32\pi^2}   \ln\big(2^{-\frac{n}2}\langle x \rangle\big)+ \mathcal{O}(1).
\end{equation}

\noindent
$(ii)$ Let $\eps>0$ and $u \in \cb^{-\frac12-\eps}_{\infty,\infty}(\R^3)$.  The sequence $(X^{(\la),(n)})$ of solutions to the renormalized stochastic equation
\begin{equation} \label{appro-equ-introd}
\left\{
\begin{aligned}
&(\partial_t +H)  X^{(\la),(n)}  = -\lambda (X^{(\la),(n)})^3+ \frakc^{(\la),(n)}X^{(\la),(n)} +\xi^{(n)}, \quad t>0, \quad x \in \R^3, \\
& X^{(\la),(n)}_0=u,
\end{aligned}
\right.
\end{equation} 
converges almost surely to a limit solution $X^{(\lambda),u}$ in the space 
$$e^{-tH}u+\bigcap_{\eta>0}\cac\big( [0,T ] ;\cb^{-\frac12-\eta}_{\infty,\infty}(\R^3)\big),$$
for all $T>0$.

 Moreover, for every $\eps >0$, one has 
$$X^{(\lambda),u} - \<Psi> +\la \<IPsi3>   \in \mathcal{C}\big((0,T ]; \mathcal{B}_{\infty, \infty}^{\frac34-\eps}(\R^3)\big), $$
  where $\<Psi>$ and $\<IPsi3> $ are explicit and only depend on the noise $\xi$.
	
\

\noindent
$(iii)$ For every $\la >0$, the semigroup $(P^{(\la)}_t)_{t\geq 0}$ generated by $X^{(\lambda),u}$ admits at least one invariant measure $\rho$, supported in  $\B^{-\frac12-\eps}_{\infty, \infty}(\R^3)$: for all bounded Borel functions {$ \Phi: \B_{\infty,\infty}^{-\frac12-\eps}(\R^3) \to \R$}
$$ \int_{\cb^{-\frac12-\eps}_{\infty,\infty}(\R^3)}  \Phi(v) \rho^{(\la)}(dv)= \int_{\cb^{-\frac12-\eps}_{\infty,\infty}(\R^3)} \big(P^{(\la)}_t \Phi)(v)\rho^{(\la)}(dv).$$

Furthermore, for every $\la >0$, any invariant measure is non-Gaussian.

\

\noindent
$(iv)$ If $\la>0$ is small enough, then the invariant measure~$\rho^{(\lambda)}$ is unique. Moreover, for all $u\in \B^{-\frac12-\eps}_{\infty, \infty}(\R^3)$, every bounded Lipschitz function $\Phi:\B^{-\frac12-\eps}_{\infty, \infty}(\R^3)\to\R$ and all $T\geq 1$, it holds that
\begin{equation*} 
\big|\mathbb{E}\big[\Phi(X_T^{(\la),u})\big]-\int \Phi(v) \rho^{(\lambda)}(dv)  \big|\leq c \big\|\Phi\big\|_{\text{Lip}} e^{-c T} ,
\end{equation*}
for some $c>0$.
\end{theorem}

\begin{remark}
The sequence $\frakc^{(n)}$ in \eqref{appro-equ-introd} is explicitly given by the expression in \eqref{con}, and with this definition in mind, it becomes clear that the expansion \eqref{expC} has already been shown in our previous work \cite{DFT} (see \cite[Section 11]{DFT}). It is recalled here for clarity. 
\end{remark}

\begin{remark}
The low regularity	$\B^{-\frac12-\eps}_{\infty, \infty}(\R^3)$ of the dynamics is the natural one, as it matches that of the associated linear solution $\<Psi>$. We refer to Section \ref{Sect3} for more details.
\end{remark}

\begin{remark}
The spectral gap inequality for \eqref{eq:SGPE}, {\it i.e.}, the exponential decay in item $(iv)$ above, has been established in \cite{dBDF18} in the one-dimensional case.
The method there shows in addition that the exponential rate $c$ coincides with the dissipation coefficient $\beta$.
\end{remark}

\

The presence of the parameter $\lambda$ in front of the nonlinearity in \eqref{eq:intro} admits an enlightening reinterpretation in terms of a modulation of the {\it confining harmonic potential}.

\smallskip

Consider indeed, for every $\alpha \in \mathbb{D}:=\{2^{2\ell}, \, \ell\geq 0\}$, the modulated operator 
$$H_\alpha := -\Delta + \alpha^2 |x|^2$$
 on $\R^3$, as well as the corresponding equation
\begin{equation} \label{eq-alpha-renormalized0}
\left\{
\begin{aligned}
&(\partial_t + H_\alpha)\, Y^{(\alpha),(n)} = -\bigl(Y^{(\alpha),(n)}\bigr)^{\!3} 
+ \widetilde{\frakc}^{(\alpha),(n)} Y^{(\alpha),(n)} + \xi^{(\alpha),(n)}, \quad t>0, \quad x \in \R^3, \\
& Y^{(\alpha),(n)}(0) = v,
\end{aligned}
\right.
\end{equation} 
where $\xi^{(\alpha),(n)}$ is the regularization of $\xi$ given by \eqref{regu-noise-alpha}. In this setting, we have the following counterpart of Theorem~\ref{thm:main} (see Section \ref{subsubsec:besov-scaling} for the definition of the general Besov spaces $\cb^{\si}_{\infty,\infty}(H_{\alpha})$).
 
\begin{theorem}[Strong-confinement uniqueness]\label{thm:uniq-alpha}
Fix $\al\in \mathbb{D}$. Then there exists a sequence $(\widetilde{\frakc}^{(\alpha),(n)})$ of \textit{deterministic} functions on $\R^3$ such that the following assertions hold:

\smallskip
 
$(i)$  For every $\eps > 0$ and $v \in \cb^{-\frac12-\eps}_{\infty,\infty}(H_{\alpha})$, 
the sequence $(Y^{(\alpha),(n)})$ of solutions to \eqref{eq-alpha-renormalized0} 
converges almost surely to a limit process $Y^{(\alpha),v}$ in the space 
$$e^{-tH_{\alpha}}v+\bigcap_{\eta>0}\cac\big( [0,T] ;\cb^{-\frac12-\eta}_{\infty,\infty}(H_{\alpha})\big),$$
for every $T > 0$.

Moreover, for every $u\in \cb^{-\frac12-\eps}_{\infty,\infty}(H)$, it holds that
\begin{equation}\label{identilaw}
Y^{(\alpha),\cs_\al u} \stackrel{\text{(law)}}{=} \ct_\alpha\, X^{(\la),u} \qquad \text{with} \quad \la = \alpha^{-\frac12},
\end{equation}
and where we have set
$$(\cs_{\alpha}u)(x) := \alpha^{\frac14}\, u(\al^{\frac12}x), \qquad (\ct_\alpha Z)(t,x) \,:=\, \alpha^{\frac14}\, Z(\alpha t, \alpha^{\frac12} x).$$
 
\smallskip
 
$(ii)$   If $\alpha $ is large enough, then the semigroup $(Q^{(\alpha)}_t)_{t \geq 0}$ generated by 
$Y^{(\alpha)}$ admits a unique invariant measure $\nu^{(\alpha)}$, supported in 
$\cb^{-\frac12-\eps}_{\infty,\infty}(H_{\alpha})$. Moreover, $\nu^{(\alpha)}$ is non-Gaussian, 
and is given explicitly as the pushforward
\begin{equation}\label{pushforward-formula}
\nu^{(\alpha)}= (\cs_\al)_\sharp \rho^{(\la)}=  \rho^{(\la)} \circ (\cs_\al)^{-1}
\qquad \text{with} \quad \la = \alpha^{-\frac12},
\end{equation}
and where $\rho^{(\la)}$ is the unique invariant measure for \eqref{appro-equ-introd} 
provided by item~$(iv)$ of Theorem~\ref{thm:main}.
 
\smallskip
 
$(iii)$  For every $v \in \cb^{-\frac12-\eps}_{\infty,\infty}(H_{\alpha})$, 
every bounded Lipschitz function $\Phi : \cb^{-\frac12-\eps}_{\infty,\infty}(H_{\alpha}) \to \R$ 
and every $T \geq 1$,
\begin{equation}\label{expo-conv-alpha}
\bigl|\, \mathbb{E}[\Phi(Y^{(\alpha),v}_T)] -\int \Phi(v) \nu^{(\alpha)}(dv)    \, \bigr| 
\,\leq\, c_0\, \alpha^{1/4}\, \|\Phi\|_{\mathrm{Lip}}\, e^{-c \alpha T}.
\end{equation}
\end{theorem}

\

Items~$(ii)$-$(iii)$ of Theorem \ref{thm:uniq-alpha} thus establish that the {\it high-temperature} (or {\it small coupling}) regime $\lambda \to 0^+$ addressed in Theorem \ref{thm:main} (item~$(iv)$) corresponds, through \eqref{identilaw}, to the {\it strong confinement} regime $\alpha \to +\infty$ in \eqref{eq-alpha-renormalized0}.  \medskip

\color{black}

\begin{remark}

Our argument yields uniqueness only in the small-coupling regime. In view of the various $\Phi^4_3$ models studied so far (see the beginning of the introduction for appropriate references), the picture beyond this regime is known to depend sensitively on the underlying geometry.
\end{remark}

\

\subsection{Organization of the paper}
 The paper is organized as follows. In Section~\ref{Sect2}, we recall the definition of harmonic Besov spaces and set up the probability framework. In Section~\ref{Sect3}, we establish local well-posedness for equation~\eqref{eq:intro} with rough initial conditions, adapting the paracontrolled approach developed in~\cite{DFT}. Section~\ref{Sect4} concerns the derivation of a priori bounds and the globalization argument; we follow the general strategy of Mourrat and Weber~\cite{MW}, with the necessary adaptations to the whole-space setting. In Section~\ref{Sect5}, we prove the Markov property for the dynamics and establish the existence of an invariant measure via the Krylov-Bogoliubov method; we further show that this measure is non-Gaussian. Section~\ref{sec:uniq-inv} is devoted to the proof of uniqueness of the invariant measure in the small-coupling regime, as well as to the transition between~\eqref{appro-equ-introd} and \eqref{eq-alpha-renormalized0}. Finally, the appendices gather functional inequalities and the detailed analysis of higher-order stochastic diagrams.

\

\textit{Note that, in order to simplify the notation, we will only prove items $(i)$ to $(iii)$ of Theorem \ref{thm:main} for $\lambda=1$
 (see Sections \ref{Sect3}–\ref{sec:uniq-inv})}. The proof extends to any $\la >0$ with no additional difficulty. The role of $\lambda$ will be emphasized in Section \ref{sec:uniq-inv} through the high-temperature condition.

\

\begin{acknowledgements}
The authors thank Tristan Robert for enlightening discussions. A. Deya and L. Thomann  were partially supported by the ANR project "SMOOTH" ANR-22-CE40-0017. 
R. Fukuizumi was supported by JSPS KAKENHI Grant Number 23H01079. 
\end{acknowledgements}

\

\section{Setting and outline}\label{Sect2}

\subsection{Besov spaces associated with the harmonic oscillator}\label{sect-Besov}

We recall the definition of the harmonic Besov spaces that we introduced in \cite[Section~2]{DFT} (see also \cite{FI} and \cite[Section~13]{DFT} for further properties).  Let 
\[
\mathcal{A} := \Big\{ \xi \in \mathbb{R}_+ : \tfrac{3}{4} \le \xi \le \tfrac{8}{3} \Big\}.
\]
There exists a  function \(\chi \in \mathcal{C}_0^{\infty}(\mathcal{A})\), taking values in \([0,1]\), such that   for all \(\xi \geq \frac43\),
\begin{equation}\label{partition}
 \sum_{j=0}^{+\infty} \chi_j(\xi) = 1,
 \qquad \text{where} \quad
 \chi_j(\xi) := \chi(2^{-j}\xi), \quad j \ge 0.
\end{equation}

The associated Hermite multipliers \((\delta_j)_{j \ge 0}\) are defined by 
\begin{equation*} 
\delta_j u = \chi_j(\sqrt{H})u = \chi\!\left(\frac{\sqrt{H}}{2^j}\right)u, \quad j \ge 0.
\end{equation*} 
Setting
\(\theta(x) := \chi(\sqrt{|x|})\), we have \(\theta  \in \mathcal{C}_0^{\infty}(\mathbb{R})\) and
\[
\operatorname{Supp}\theta 
\subset \Big\{ \xi \in \mathbb{R} : (\tfrac{3}{4})^2 \le |\xi| \le (\tfrac{8}{3})^2 \Big\}.
\]

For \(j \ge 0\), we can equivalently write
\begin{equation*} 
 \delta_j = \theta\!\left(\frac{H}{2^{2j}}\right).
\end{equation*}

\medskip

The Besov spaces associated with the harmonic oscillator are defined, for \(1 \le p \le \infty\) and \(\sigma \in \mathbb{R}\), by
\begin{equation}\label{def-besov} 
\mathcal{B}^{\sigma}_{p,\infty}(\mathbb{R}^3)
= \Big\{ u \in \mathscr{S}'(\mathbb{R}^3) : 
\big(2^{j\sigma} \|\delta_j u\|_{L^p(\mathbb{R}^3)}\big)_{j \ge 0} \in \ell^\infty \Big\},
\end{equation}
and these spaces are endowed with the norm
\begin{equation}\label{def-theta}
\|u\|_{\mathcal{B}^{\sigma}_{p,\infty}(\mathbb{R}^3)} =
\displaystyle \sup_{j \ge 0} \| 2^{j\sigma} \delta_j u \|_{L^p(\mathbb{R}^3)}.
\end{equation}
For convenience, in the sequel, we denote
$$
\mathcal{B}^{\sigma}_{p} := \mathcal{B}^{\sigma}_{p,\infty}(\mathbb{R}^3).
$$
We refer to \cite[Section 13]{DFT} for some results on the spaces $\mathcal{B}^{\sigma}_{p}$. 
\medskip

For $1 \leq p \leq  \infty$, we denote the usual Lebesgue spaces by 
$$L^p=L^p(\R^3).$$

 We consider a Hilbertian basis $(\vp_k)_{k \geq 0}$ of $L^2(\R^3)$, composed of  eigenvectors of $H$ with eigenvalues~$(\la_k)_{k \geq 0}$: 
$$H\varphi_k=\lambda_k \varphi_k.$$ 
For more details on the spectral theory of the harmonic oscillator, we refer to \cite[Section 2.1]{DFT}. 

Finally, we denote the kernel of the exponential operator $e^{-tH}$ by
\begin{equation*}  
K_t(x,y):=\sum_{k\geq 0} e^{-t\lambda_k}\varphi_k(x)\varphi_k(y),
\end{equation*}
and note that it is explicitly given by the Mehler formula (see {\it e.g.} \cite[page 109]{Taylor}):
\begin{equation}  \label{mehler2}
K_t(x,y)  = (2\pi\sinh 2t)^{-\frac{3}{2}} \exp\left(- \frac{ \vert x-y\vert^2}{4\tanh t}-\frac{\tanh t}{4}\vert x+y\vert^2\right).
\end{equation} 
 
\subsection{Probability framework}

By classical convention, we represent the white noise $\xi$ as 
\begin{equation}\label{repres-noise}
\xi:=dW_t, \quad \text{where}  \;\;\;  W_t(x):=\sum_{k \geq 0} \beta^{(k)}_t  \vp_k(x)
\end{equation}
is a cylindrical Wiener process defined through a family $(\beta^{(k)})_{k\geq 0}$ of independent (real-valued) Brownian motions, on the probability space $(\Omega, \mathcal{F}, \mathbb{P})$ introduced above. We denote by  $(\mathcal{F}_t)_{t \geq 0}$ the natural filtration generated by these Brownian motions. 

\smallskip

Just as in \cite{DFT}, we consider the natural regularization $(\xi^{(n)})_{n\geq 1}$ of $\xi$ given for all $t\geq 0$ by
\begin{equation}\label{regu-noise}
 \xi^{(n)}_t:= \frac{dW^{(n)}_t}{dt}, \quad \text{where} \;\;\; W^{(n)}_t(x):=\sum_{k \geq 0} e^{-\eps_n \la_k} \beta^{(k)}_t  \vp_k(x),
\end{equation} 
where $\eps_n:=2^{-n}$. For every $n\geq 1$ and $s\geq 0$, we denote by $\<Psi>^{(n)}_{s,.}$ the so-called linear solution associated with our problem (starting at $s$), that is, the process satisfying
\begin{equation*} 
\left\{
\begin{aligned}
&(\partial_t +H)  \<Psi>^{(n)}_{s,.}  = \xi^{(n)}, \quad t>s, \quad x \in \R^3, \\
& \<Psi>^{(n)}_{s,s}=0.
\end{aligned}
\right.
\end{equation*} 
In other words, for all $0\leq s\leq t$ and $x\in \R^3$,
\begin{equation}\label{premluxo}
\<Psi>^{(n)}_{s,t}(x)=\int_{s}^t \big(e^{-(t-r)H}dW^{(n)}_r\big)(x)=\int_{s}^{t} \int_{\R^3} K_{t-r}(x,y)\, W^{(n)}(dr,dy). 
\end{equation}

\smallskip

In the sequel, we shall rely on the calculation rules in Wiener chaoses, and to this end, we introduce the multiple integrals $(I^W_k)_{k\geq 0}$ driven by $W$ (see \cite[Section~1.1.2]{nualart-book}). 
 With this notation, observe that \eqref{premluxo} can be naturally recast as (see \cite[Section 4.1]{DFT} for details) 
\begin{equation}\label{defi-f-n-ell}
\<Psi>^{(n)}_{0,t}(y)=I^W_1\big(F^{(n)}_{t,y}\big), \quad \text{with} \ F^{(n)}_{t,y}(s,w):=K_{t-s+\eps_n}(y,w) \1_{[0,t]}(s).
\end{equation}
Let us denote by $\cac^{(n)}$ the covariance function of $\<Psi>^{(n)}_{0,.}$, that is, for all $t_1,t_2\geq 0$ and $y_1,y_2\in \R^3$,
$$
\cac^{(n)}_{t_1,t_2}(y_1,y_2):= \mathbb{E}\Big[ \<Psi>^{(n)}_{0,t_1}(y_1) \<Psi>^{(n)}_{0,t_2}(y_2)\Big]=\langle F^{(n)}_{t_1,y_1} ,F^{(n)}_{t_2,y_2} \rangle_{L^2(\R \times \R^3)}.
$$
Given the expression of $F^{(n)}$ in \eqref{defi-f-n-ell}, it is readily checked that 
\begin{equation}\label{cova-c-n-l}
\cac^{(n)}_{t_1,t_2}(y_1,y_2)=\frac12\int^{t_2+t_1+2\varepsilon_n}_{|t_2-t_1|+2\varepsilon_n}d\si\,  K_{\si}(y_1,y_2).
\end{equation}

\

\

\begin{remark}  
With these notations, we can explain in detail  the case $\lambda=0$ in \eqref{appro-equ-introd}. In this case, the counterterms are trivial, $\frakc^{(0),(n)}=0$, and  we have
$$X^{(0),(n)}_t=e^{-t H}X_0+\<Psi>^{(n)}_{0,t}= e^{-t H}\Big(X_0-\int_{-\infty}^0  e^{rH}dW^{(n)}_r\Big)+\<Psi>^{(n)}_{-\infty,t}, $$
where $\<Psi>^{(n)}_{-\infty,.}$ is the unique stationary solution to the linear  equation 
$$(\partial_t +H)  \<Psi>^{(n)}_{-\infty,.}  = \xi^{(n)}, \quad t \in \R, \quad x \in \R^3.$$
Assume that $X_0 \in\B^{-\frac12-\eps}_{\infty, \infty}(\R^3)$. Then $X^{(0),(n)}$ converges almost surely to $X^{(0)}$ in the space 
$$ e^{-tH}X_0+\bigcap_{\eta>0}\cac\big( [0,T] ;\cb^{-\frac12-\eta}_{\infty,\infty}(\R^3)\big),$$
 where $X^{(0)}$ is given by 
\begin{equation}\label{deta0}
X^{(0)}_t= e^{-t H}\Big(X_0-\int_{-\infty}^0  e^{rH}dW_r\Big)+\<Psi>_{-\infty,t}. 
\end{equation}
Here $\<Psi>_{-\infty,.}$ is the  unique stationary solution to the linear  Langevin equation 
$$(\partial_t +H)  \<Psi>_{-\infty,.}  = \xi, \quad t \in \R, \quad x \in \R^3,$$
which yields the measure $\rho^{(0)}:= \mathscr{L}( \<Psi>_{-\infty,t})$ for all $t \in \R$. This law corresponds to the Gaussian free field (GFF), also characterized by $\rho^{(0)}=\mathscr{L}(\gamma)$, where 
$$\gamma : \om \mapsto \sum_{n \geq 0}\frac{g_n(\om) \varphi_n}{\sqrt{2\lambda_n}},$$
for a family $(g_n)_{n \ge 0}$ of independent, standard real-valued Gaussians. The GFF can be formally written as
$$d\rho^{(0)}=Z^{-1}\exp \Big( - \int_{\R^3} \big(|\nabla u|^2+|x|^2 |u|^2\big) dx \Big)du.$$
We refer to the introduction of \cite{PRT2} for more properties of this measure.

Now let $\Phi:\B^{-\frac12-\eps}_{\infty, \infty}(\R^3)\to\R$ be a bounded Lipschitz function. For all $t\in \R$, we have
$$ \mathbb{E}\big[  \Phi(\<Psi>_{-\infty,t})\big]   =  \int \Phi(v) \rho^{(0)}(dv).$$
Due to the exponential decay of the semigroup $e^{-tH}$, one deduces from \eqref{deta0} that for all $T\geq 1$,
\begin{equation*}
\Big|\mathbb{E}\big[\Phi(X_T^{(0)})\big]- \int \Phi(v) \rho^{(0)}(dv)  \Big|\leq c \big|\Phi\big|_{\text{Lip}} e^{-c T} ,
\end{equation*}
for some $c>0$. This is a direct consequence of the fact that when a linear operator has a strictly positive discrete spectrum, the transition semigroup of the associated linear SPDE converges exponentially fast to its unique Gaussian invariant measure (see {\it e.g.} \cite[Section 11.3]{DPZ} and \cite[Section~5.3]{hairer}). 
\end{remark} 

\

 \section{Local well-posedness for rough initial conditions}\label{Sect3}

We start with the local well-posedness analysis of equation~\eqref{eq:intro}. As in our previous work~\cite{DFT}, we first introduce the stochastic diagrams that encode the singular part of the dynamics, and then reformulate the problem as a fixed-point equation for a more regular remainder. Although this local analysis is similar in spirit to \cite[Section 3]{DFT}, three new elements distinguish it:

$(i)$ we refine the regularity results for the diagrams by using topologies that are better suited to the forthcoming globalization argument (compare Table \ref{table-reg} below with the one in~\cite[Section 3]{DFT});

$(ii)$ we show how to reduce the analysis to \textit{time-independent} renormalization constants, by controlling the difference with the original time-dependent quantities (see Lemma~\ref{lem:c});

$(iii)$ we rely on topologies with a singularity at zero (see the norm in \eqref{defino-x-ept}), which allow us to consider distributional-valued initial conditions.

\smallskip

Throughout this section, we take $\la=1$ (see the end of Section~\ref{sec:introd}).

\subsection{The $\Phi^4_3$ diagrams}\label{subsec:dia}

Recall that the white-noise approximation $(\xi^{(n)})_{n\geq 1}$ has been introduced in \eqref{regu-noise}, as well as the (approximated) linear solution 
\begin{equation}\label{def-luxost}
\<Psi>^{(n)}_{s,t}(x):=\int_{s}^t \big(e^{-(t-r)H}dW^{(n)}_r\big)(x). 
\end{equation}
Then we define the higher-order diagrams at the core of our analysis by the formulas: for all $0\leq s\leq t$,
\begin{equation*} 
\<Psi2>^{(n)}_{s,t}:=\big(\<Psi>^{(n)}_{s,t}\big)^2-\frakc^{\mathbf{1},(n)}_{s,t}, \quad  \<IPsi2>^{(n)}_{s,t}:=\int_s^t dr \, e^{-(t-r)H} \<Psi2>^{(n)}_{s,r}, 
\end{equation*}
\begin{equation*} 
\<Psi3>^{(n)}_{s,t}:=\big(\<Psi>^{(n)}_{s,t}\big)^3-3\frakc^{\mathbf{1},(n)}_{s,t} \<Psi>^{(n)}_{s,t}, \quad \<IPsi3>^{(n)}_{s,t}:=\int_s^t dr \, e^{-(t-r)H}\<Psi3>^{(n)}_{s,r},
\end{equation*}
\begin{equation}\label{ordres4} 
\<PsiIPsi3>^{(n)}_{s,t}:=\<Psi>^{(n)}_{s,t}\pe \<IPsi3>^{(n)}_{s,t}, \quad \quad \<Psi2IPsi2>^{(n)}_{s,t}:=\<Psi2>^{(n)}_{s,t} \pe \<IPsi2>^{(n)}_{s,t}-\frakc^{\mathbf{2},(n)}_{s,t},
\end{equation}
\begin{equation}\label{ord5}
\<Psi2IPsi3>^{(n)}_{s,t}:=\<Psi2>^{(n)}_{s,t} \pe \<IPsi3>^{(n)}_{s,t} -3\, \frakc^{\mathbf{2},(n)}_{s,t}\, \<Psi>^{(n)}_{s,t} ,
\end{equation}
where the deterministic sequences $(\frakc^{\mathbf{1},(n)}),(\frakc^{\mathbf{2},(n)})$ are respectively given by
\begin{equation}\label{constante-re}
\frakc^{\mathbf{1},(n)}_{s,t}(x):=\mathbb{E}\Big[ \big| \<Psi>^{(n)}_{s,t}(x)\big|^2\Big] \quad \text{and} \quad \frakc^{\mathbf{2},(n)}_{s,t}(x):=\mathbb{E}\Big[ \<Psi2>^{(n)}_{s,t}(x) \<IPsi2>^{(n)}_{s,t}(x)\Big].
\end{equation}

We denote this set of diagrams by
\begin{equation}\label{zn} 
Z_{s,.}(\xi^{(n)}):=\Big(\<Psi>^{(n)}_{s,.}, \<Psi2>^{(n)}_{s,.}, \<IPsi3>^{(n)}_{s,.},\<PsiIPsi3>^{(n)}_{s,.},\<Psi2IPsi2>^{(n)}_{s,.},\<Psi2IPsi3>^{(n)}_{s,.}\Big).
\end{equation}

\

Observe that, as an immediate consequence of the stationarity of the Brownian increments, one has, for all $s\geq 0$,
\begin{equation*}
W^{(n)}_{.+s}-W^{(n)}_s\stackrel{\text{(law)}}{=} W^{(n)}_{.} \quad \text{in} \ \ \cac(\R_+\times \R^3 ;\R).
\end{equation*}
 Therefore, since
\begin{equation*}
\<Psi>^{(n)}_{s,s+t}=\int_{s}^{s+t} e^{-(s+t-r)H}dW^{(n)}_r=\int_{0}^{t} e^{-(t-r)H}d\big(W^{(n)}_{s+r}-W^{(n)}_s\big)
=\int_{0}^{t} e^{-(t-r)H}dW^{(n)}_{r},
\end{equation*}
it follows that for all $s\geq 0$,
\begin{equation*}
\<Psi>^{(n)}_{s,s+.}\stackrel{\text{(law)}}{=} \<Psi>^{(n)}_{0,.} \quad \text{in} \ \ \cac(\R_+\times \R^3 ;\R).
\end{equation*}

Going back to the above definitions of the components of $Z(\xi^{(n)})$, we deduce the following useful stationarity property.

\begin{lemma}\label{lem:stati-z}
For every $n\geq 1$, the process $Z(\xi^{(n)})$ is stationary in time, in the sense that for every $s\geq 0$,
\begin{equation*}
Z_{s,s+.}(\xi^{(n)})\stackrel{\text{(law)}}{=} Z_{0,.}(\xi^{(n)}) \quad \text{in} \ \ \cac(\R_+\times \R^3 ;\R).
\end{equation*}
\end{lemma}

\color{black}

\

The following statement summarizes our main control results on these stochastic diagrams.

\begin{proposition}\label{prop:conv-arbre}
For all $0\leq s<T$ and $0<\eps<\frac14$, the sequence  $\big(Z_{s,.}(\xi^{(n)})\big)_{n\geq 1}$ defined by \eqref{zn} \color{black} converges almost surely in the space
\begin{align}
&\mathcal Z_{\eps,  [s,T]}
:=\cac([s,T]; \cb_{\infty}^{-\frac12-\eps}) \times \cac([s,T];\cb_{\infty}^{-1-\eps}) \times  \big(\cac([s,T];\cb_{\infty}^{\frac12-\eps})\cap \cac^{\frac14-\eps}([s,T]; \cb^\eps_{\infty})\big)\nonumber\\
&\hspace{2cm} \times 
\cac([s,T]; \cb^{-\eps} _{\infty})\times \cac([s,T];\cb^{-\eps}_{\infty}) \times \cac([s,T];\cb_{\infty}^{-\frac 56-\eps}) .\label{space-for-z}
\end{align}
We naturally denote its limit by
$$Z_{s,.}:=Z_{s,.}(\xi)=\Big(\<Psi>_{s,.}, \<Psi2>_{s,.}, \<IPsi3>_{s,.},\<PsiIPsi3>_{s,.},\<Psi2IPsi2>_{s,.},\<Psi2IPsi3>_{s,.}\Big).$$
and set  
	\begin{align}
&\big\|Z_{s,.}\big\|_{\mathcal Z_{\frac{\eps}{2},[s,T]}}:=\big\|\<Psi>_{s,.}\big\|_{\cac([s,T]; \cb_{\infty}^{-\frac12-\eps})} + \big\|\<Psi2>_{s,.}\big\|_{\cac([s,T]; \cb_{\infty}^{-1-\eps})}+\big\|\<IPsi3>_{s,.}\big\|_{\cac([s,T]; \cb_{\infty}^{\frac12-\eps})}+\big\|\<IPsi3>_{s,.}\big\|_{\cac^{\frac14-\eps}([s,T]; \cb_{\infty}^{\eps})} \nonumber\\ 
&\hspace{3cm}+\big\|\<PsiIPsi3>_{s,.}\big\|_{\cac([s,T]; \cb_{\infty}^{-\eps})}+\big\|\<Psi2IPsi2>_{s,.}\big\|_{\cac([s,T]; \cb_{\infty}^{-\eps})}+\big\|\<Psi2IPsi3>_{s,.}\big\|_{\cac([s,T]; \cb_{\infty}^{-\frac56-\eps})} .\label{K}
\end{align}

Moreover, it holds that for every $p\geq 1$
\begin{equation} \label{est:tree_uniform}
\sup_{s \ge 0} \mathbb{E}\Big[\big\|Z_{s,.}\big\|_{\mathcal Z_{\frac{\eps}{2},[s,s+1]}}^p\Big] < \infty.
\end{equation}
\end{proposition}

\smallskip

\begin{proof}

Using Lemma \ref{lem:stati-z}, the convergence and regularity results for the diagrams of order less than 3 (that is, $\<Psi>, \<Psi2>, \<IPsi3>$) can be immediately derived from \cite[Proposition 3.2]{DFT}.

\smallskip

By contrast, the convergence and regularity results for the fourth-order diagrams $\<PsiIPsi3>,\<Psi2IPsi2>$ and the fifth-order diagram $\<Psi2IPsi3>$ are expressed in different topologies than in \cite[Proposition 3.2]{DFT}, better suited to the subsequent globalization arguments.
This topological refinement stems from an interpolation procedure, developed in Appendices \ref{sec:diag-4th-order-1}, \ref{sec:diag-4th-order-2}, and \ref{sec:fifth-o-i}.

\end{proof}

\smallskip

\begin{table}
{\renewcommand{\arraystretch}{2}
\begin{tabular}{| c |c | c | c |c | c | c |c | c |  }
\hline 
\; $\tau $\; &  \ $\<Psi> $& $\<Psi2>$ &\;  $\<IPsi3>$\;  & \; $\<PsiIPsi3>$ \; & \; $\<Psi2IPsi2>$ \; & $\<Psi2IPsi3>$ \    \\
\hline 
 $\mathcal{E}_{\tau}$ &  $\cac_T \cb_{\infty}^{-\frac12-\eps}$  & $\cac_T\cb_{\infty}^{-1-\eps} \cap \cac^{\frac14-\eps}_T\cb_{\infty}^{-\frac32+\eps} $ & $\cac_T\cb_{\infty}^{\frac12-\eps}\cap \cac^{\frac14-\eps}_T \cb^\eps_{\infty}$    
 &   $ \cac_T\cb^{-\eps} _{\infty}$   &     $\cac_T\cb^{-{\eps}}_{\infty} $    &  $ \cac_T\cb_{\infty}^{ -\frac 56-\eps} $   \\[5pt]
 \hline 
\end{tabular}}

 \caption{The diagrams and their regularity : $\tau \in \mathcal{E}_{\tau}$.}\label{table-reg}
 \end{table}

\subsection{Renormalizing sequence}\label{subsec:ren-seq}

We explicitly define the renormalizing sequence $(\frakc^{(n)})$ in \eqref{appro-equ-introd} as
\begin{equation}\label{con}
\frakc^{(n)} := 3\frakc^{\mathbf{1},(n)} - 9\frakc^{\mathbf{2},(n)},
\end{equation}
with 
$$\frakc^{\mathbf{1},(n)}:=\frakc^{\mathbf{1},(n)}_{-\infty,0} \quad \text{and} \quad \frakc^{\mathbf{2},(n)}:=\frakc^{\mathbf{2},(n)}_{-\infty,0}.$$
These two formulas, formally interpreted through \eqref{constante-re}, must be rigorously read as
\begin{equation}\label{cstts-rigor-1}
\frakc^{\mathbf{1},(n)}(x):=\frac12 \int_{2\varepsilon_n}^{+\infty} d\sigma \, K_{\sigma}(x,x)
		\end{equation}
		and
		\begin{equation}\label{cstts-rigor-2}
\frakc^{\mathbf{2},(n)}(x):= \frac12 \int_0^{+\infty} ds \int dw \, K_{s}(x,w)
    \bigg( \int_{s+2\varepsilon_n}^{+\infty} d\sigma \, K_{\sigma}(x,w) \bigg)^2.
		\end{equation}

\

To compare these two (time-independent) quantities with the (time-dependent) sequences in~\eqref{constante-re}, we introduce the difference functions
$$
\overline{\frakc}^{\mathbf{1},(n)}_{s,t}:=\frakc^{\mathbf{1},(n)}-\frakc^{\mathbf{1},(n)}_{s,t}, \qquad \overline{\frakc}^{\mathbf{2},(n)}_{s,t}:=\frakc^{\mathbf{2},(n)}-\frakc^{\mathbf{2},(n)}_{s,t}
$$
and
$$\overline{\frakc}^{(n)}_{s,t} := 3\, \overline{\frakc}^{\mathbf{1},(n)}_{s,t} -9\, \overline{\frakc}^{\mathbf{2},(n)}_{s,t},$$
for all $0\leq s<t$. We also set (arbitrarily) $\overline{\frakc}^{\mathbf{1},(n)}_{t,t}=\overline{\frakc}^{\mathbf{2},(n)}_{t,t}=\overline{\frakc}^{(n)}_{t,t}=0$ for all $t\geq 0$.

\smallskip

For clarity, we introduce the notation: for all $s\geq 0$ and $s<T\leq s+1$,
\begin{equation}\label{llbrack}
\big\llbracket  f\big\rrbracket_{\eps,[s,T]}:= \max\Big(\sup_{s<t\leq T} |t-s|^{{\frac{9}{16}}} \big\| f_{s,t}\big\|_{L^\infty},\; \sup_{s<t\leq T} |t-s|^{\frac{7}{8}}  \big\| f_{s,t}\big\|_{\cb^{\frac12+\eps}_\infty}\Big).
\end{equation}

\

\begin{lemma}\label{lem:c}
For every $\eps>0$ small enough, it holds that
\begin{equation}\label{knack}
\sup_{n\geq 1}\sup_{s\geq 0}\Big( \big\llbracket  \overline{\frakc}^{\mathbf{1},(n)}_{s,.}\big\rrbracket_{\eps,[s,s+1]}+ \big\llbracket  \overline{\frakc}^{\mathbf{2},(n)}_{s,.}\big\rrbracket_{\eps,[s,s+1]}\Big) <\infty.
\end{equation}

\end{lemma}

\begin{proof}

\smallskip

\noindent
1. \textit{Bound for $\overline{\frakc}^{\mathbf{1},(n)}_{s,t}$.}
Recall the expression (see \cite[Section 4.2]{DFT}) 
\begin{align*}
&\frakc^{\mathbf{1},(n)}_{s,t}(x)=\mathbb{E}\Big[ \big|\<Psi>^{(n)}_{s,t}(x)\big|^2\Big]=\frac12\int_{2\eps_n}^{2(t-s+\eps_n)}d\si\,  K_\si(x,x),
\end{align*}
which, in combination with \eqref{cstts-rigor-1}, yields
\begin{align}
&\overline{\frakc}^{\mathbf{1},(n)}_{s,t}(x)=\frac12\int_{2(t-s+\eps_n)}^{+\infty}d\si\,  K_\si(x,x).\label{overline-frakc-1}
\end{align}
As a result, for all $0\leq s<t$ such that $|t-s|\leq 1$, one has
\begin{align*}
\big\|\overline{\frakc}^{\mathbf{1},(n)}_{s,t}\big\|_{\cb^0_\infty}\lesssim \big\|\overline{\frakc}^{\mathbf{1},(n)}_{s,t}\big\|_{L^\infty}&\lesssim \int_{t-s}^{+\infty}\frac{d\si}{\si^{\frac32}}\lesssim \frac{1}{|t-s|^{\frac12}}.
\end{align*}
On the other hand, using the expression in \eqref{mehler2}, we can compute explicitly
\begin{multline}\label{expli-H}
H\big(\overline{\frakc}^{\mathbf{1},(n)}_{s,t}\big)(x)=\\
=c \int_{2(t-s+\eps_n)}^{+\infty}\frac{d\si}{(\sinh 2\si)^{\frac{3}{2}}} \Big(6 \tanh(\si)-4\tanh^2(\si)|x|^2  +|x|^2\Big) \exp\left(- \tanh(\si) \vert x\vert^2\right),
\end{multline}
which gives, for all $0\leq s<t$ such that $|t-s|\leq 1$,
\begin{align*}
\big\|H\big(\overline{\frakc}^{\mathbf{1},(n)}_{s,t}\big)\big\|_{L^\infty}&\lesssim \int_{2(t-s)}^{+\infty}\frac{d\si}{(\sinh 2\si)^{\frac{3}{2}}}  \Big( \tanh(\si)+\frac{1}{\tanh(\si)}\Big)\\
&\lesssim \int_{2(t-s)}^{2}\frac{d\si}{(\sinh 2\si)^{\frac{3}{2}}}\frac{1}{\tanh(\si)}+\int_{2}^{+\infty}\frac{d\si}{(\sinh 2\si)^{\frac{3}{2}}} \\
&\lesssim \int_{2(t-s)}^{2}\frac{d\si}{\si^{\frac{5}{2}}}+1\lesssim \frac{1}{|t-s|^{\frac32}},
\end{align*}
and hence by \cite[Lemma 13.11]{DFT}
$$\big\|\overline{\frakc}^{\mathbf{1},(n)}_{s,t} \big\|_{\cb_\infty^{2-\eps}} \lesssim \big\| \overline{\frakc}^{\mathbf{1},(n)}_{s,t}\big\|_{\cw^{2,\infty}} \lesssim \frac{1}{|t-s|^{\frac32}}.$$
By interpolation, 
\begin{align*}
  \big\| \overline{\frakc}^{\mathbf{1},(n)}_{s,t} \big\|_{\cb_\infty^{\frac12+\eps}}&\lesssim \big\| \overline{\frakc}^{\mathbf{1},(n)}_{s,t} \big\|_{\cb_\infty^{0}}^{\frac{3-4\eps}{4-2\eps}} \big\| \overline{\frakc}^{\mathbf{1},(n)}_{s,t} \big\|_{\cb_\infty^{2-\eps}}^{\frac{1+2\eps}{4-2\eps}} \lesssim \frac{1}{|t-s|^{\frac{6+2\eps}{8-4\eps}}} ,
\end{align*}
for $\eps>0$ small enough, which is sufficient for our purpose. 

\

\smallskip

\noindent
2. \textit{Bound for $\overline{\frakc}^{\mathbf{2},(n)}_{s,t}$.}
Let $0\leq s<t$ such that $|t-s|\leq 1$. Recall the expression (see\cite[expression~(10.1)]{DFT}) 
\begin{align*}
\frakc^{\mathbf{2},(n)}_{s,t}(x)&=\mathbb{E}\Big[ \<Psi2>^{(n)}_{s,t}(x) \<IPsi2>^{(n)}_{s,t}(x)\Big]\\
&=\frac12 \int_s^{t} dr \int dw \, K_{t-r}(x,w)
    \bigg( \int_{t-r+2\varepsilon_n}^{t-r+2(r-s)+2\eps_n} d\sigma \, K_{\sigma}(x,w) \bigg)^2\\
& =\frac12 \int_0^{t-s} dr \int dw \, K_{r}(x,w)
    \bigg( \int_{r+2\varepsilon_n}^{r+2(t-s-r)+2\eps_n} d\sigma \, K_{\sigma}(x,w) \bigg)^2, 
\end{align*}
which entails, in combination with \eqref{cstts-rigor-2}, 
\begin{align}
&\overline{\frakc}^{\mathbf{2},(n)}_{s,t}(x)=\frac12 \int_0^{+\infty} dr \int dw \, K_{r}(x,w)
    \bigg( \int_{r+2\varepsilon_n}^{+\infty} d\sigma \, K_{\sigma}(x,w) \bigg)^2\nonumber\\
&\hspace{4cm}-\frac12 \int_0^{t-s} dr \int dw \, K_{r}(x,w)
    \bigg( \int_{r+2\varepsilon_n}^{r+2(t-s-r)+2\ep_n} d\sigma \, K_{\sigma}(x,w) \bigg)^2\nonumber\\
&=\frac12 \int_{t-s}^{+\infty} dr \int dw \, K_{r}(x,w)
    \bigg( \int_{r+2\varepsilon_n}^{+\infty} d\sigma \, K_{\sigma}(x,w) \bigg)^2\nonumber\\
&\hspace{1cm}+\frac12 \int_0^{t-s} dr \int dw \, K_{r}(x,w) \bigg[\bigg( \int_{r+2\varepsilon_n}^{+\infty} d\sigma \, K_{\sigma}(x,w) \bigg)^2-\bigg( \int_{r+2\varepsilon_n}^{r+2(t-s-r)+2\ep_n} d\sigma \, K_{\sigma}(x,w) \bigg)^2\bigg]\nonumber\\
&=\frac12 \int_{t-s}^{+\infty} dr \int dw \, K_{r}(x,w)
    \bigg( \int_{r+2\varepsilon_n}^{+\infty} d\sigma \, K_{\sigma}(x,w) \bigg)^2\nonumber\\
&\hspace{1cm}+\frac12 \int_0^{t-s} dr \int dw \, K_{r}(x,w) \bigg( \int_{r+2(t-s-r)+2\ep_n}^{+\infty} d\sigma \, K_{\sigma}(x,w) \bigg)\bigg( \int_{r+2\varepsilon_n}^{+\infty} d\eta \, K_{\eta}(x,w) \bigg)\nonumber\\
&\hspace{1cm}+\frac12 \int_0^{t-s} dr \int dw \, K_{r}(x,w) \bigg( \int_{r+2(t-s-r)+2\ep_n}^{+\infty} d\sigma \, K_{\sigma}(x,w)  \bigg)\bigg( \int_{r+2\varepsilon_n}^{r+2(t-s-r)+2\ep_n} d\eta \, K_{\eta}(x,w) \bigg)\nonumber\\
&=:A^{(n)}_{s,t}(x)+B^{(n)}_{s,t}(x)+C^{(n)}_{s,t}(x).\label{decompo-c-2-c-2}
\end{align}

\

As far as $A^{(n)}_{s,t}$ is concerned, since $\dis \int dw \, K_{r}(x,w)\lesssim 1$, one has
\begin{align*}
\big|A^{(n)}_{s,t}(x)\big|&\lesssim  \int_{t-s}^{+\infty} dr \int dw \, K_{r}(x,w)
    \bigg( \int_{r+2\varepsilon_n}^{+\infty} \frac{d\sigma}{\sinh(2\si)^{\frac32}}  \bigg)^2\\
		&\lesssim  \int_{t-s}^{1} dr 
    \bigg( \int_{r}^{+\infty} \frac{d\sigma}{\si^{\frac32}}  \bigg)^2  + 1  \lesssim  \frac{1}{|t-s|^\eps},
\end{align*} 
for any $\eps>0$. Then, observe that 
\begin{align*}
\big|C^{(n)}_{s,t}(x)\big| \lesssim  \big|B^{(n)}_{s,t}(x)\big|,
\end{align*}
and thus it is enough to bound $\big|B^{(n)}_{s,t}(x)\big|$. To this end, we can write 
\begin{align*}
\big|B^{(n)}_{s,t}(x)\big|&\lesssim   \int_0^{t-s} dr \int dw \, K_{r}(x,w) \bigg( \int_{r+2(t-s-r)}^{+\infty} \frac{d\sigma}{\sinh(2\si)^{\frac32}} \bigg)\bigg( \int_{r}^{+\infty} \frac{d\eta}{\sinh(2\eta)^{\frac32}}\bigg)\\
&\lesssim     \int_0^{t-s} \frac{dr}{(r+2(t-s-r))^{\frac12} r^{\frac12}}\lesssim     \int_0^{t-s} \frac{dr}{(t-s-r)^{\frac12} r^{\frac12}}   \lesssim   1,
\end{align*} 
using again the bound $\dis \int dw \, K_{r}(x,w)\lesssim 1$ for any $r>0$. 

Going back to \eqref{decompo-c-2-c-2}, we have thus shown that for  any $\eps>0$,
\begin{equation*}
\big\|\overline{\frakc}^{\mathbf{2},(n)}_{s,t}\big\|_{L^\infty} \lesssim \frac{1}{|t-s|^\eps}.
\end{equation*}

We now establish a uniform bound for $H\big(\overline{\frakc}^{\mathbf{2},(n)}_{s,t}\big)$. By the Minkowski estimate, one has 
 \begin{align}
\|A^{(n)}_{s,t}\|_{\mathcal{W}^{\frac12+2\eps,\infty}}&= \|H^{\frac14+\eps}A^{(n)}_{s,t}\|_{L^\infty}\nonumber\\
&\lesssim  \int_{t-s}^{+\infty} dr  \int_{r}^{+\infty} d\sigma  \int_{r}^{+\infty} d\eta \,\Big\| H^{\frac14+\eps}\Big(  \int dw \,    K_r(.,w) K_\sigma(.,w) K_\eta(.,w)\Big)\Big\|_{L^{\infty}}\label{Minko}
\end{align}
and by interpolation
 \begin{multline}\label{intter}
 \Big\| H^{\frac14+\eps}\Big( \int dw \,    K_r(.,w) K_\sigma(.,w) K_\eta(.,w)\Big)\Big\|_{L^{\infty}}\leq \\
 \leq \Big\|\int dw \,    K_r(.,w) K_\sigma(.,w) K_\eta(.,w)\Big\|^{\frac34-\eps}_{L^{\infty}} \Big\| H\Big( \int dw \,    K_r(.,w) K_\sigma(.,w) K_\eta(.,w)\Big)\Big\|^{\frac14+\eps}_{L^{\infty}}.
\end{multline}
For the first term in \eqref{intter}, we can write for any $r,\si,\eta>0$
\begin{equation}\label{firs}
 \int dw \,    K_r(x,w) K_\sigma(x,w) K_\eta(x,w)  \lesssim \frac{1}{(\sigma \eta)^{\frac32}} \int dw \,    K_r(x,w)   \lesssim \frac{1}{(\sigma \eta)^{\frac32}}.
\end{equation}
Let us turn to the control of the second term in \eqref{intter}. To begin with, observe that for any three functions $f,g,h$ on $\R^3$, we can expand $H\big(fgh\big)$ as
\begin{equation*}
H\big(fgh\big)= (H f) (g h)-2 \,\langle \nabla f, \nabla g\rangle\, h -2\, \langle \nabla f, \nabla h\rangle\, g -f \,\Big( (\Delta g)\,  h+2 \,\langle \nabla g,\nabla h\rangle+g\, (\Delta h) \Big) ,
\end{equation*}
and thus, we have here
\small
\begin{align*} 
&\Big| \int dw \, H\big(K_r(.,w) K_\sigma(.,w) K_\eta(.,w)\big)\Big|\lesssim  \int dw \, \big|(H_xK_r)(.,w) K_\sigma(.,w) K_\eta(.,w)\big| \\
& + \int dw \, \big|\langle \nabla_x K_r(.,w), \nabla_x K_\sigma(.,w)\rangle  K_\eta(.,w)+\langle \nabla_x K_r(.,w), \nabla_x K_\eta(.,w)\rangle  K_\sigma(.,w)\big|\\
&+ \int dw \, \Big[ \big|K_r(.,w) \Delta_xK_\sigma (.,w)  K_\eta(.,w)\big|+\big|K_r(.,w)  \langle \nabla_x K_\si(.,w),\nabla_x  K_\eta(.,w) \rangle \big|+\big|K_r(.,w)K_\sigma (.,w)   \Delta_x K_\eta(.,w)\big|\Big] .
\end{align*}
\normalsize
To bound each of these quantities, recall that, thanks to \cite[Lemma 2.3]{DFT}, we can rely on the following pointwise estimates: for all $x,w \in \R^3$ and $\sigma>0$,
 \begin{equation}\label{Bb0}
\big| \big( H_{x}K_{\sigma}\big)(x,w)\big|  \lesssim \sigma^{-\frac{5}{2}}  \exp\Big( - \frac{|x-w|^2}{8 \tanh(\sigma)}  \Big),
\end{equation}
 \begin{equation}\label{Bb1}
\big| \big( \Delta_{x}K_{\sigma}\big)(x,w)\big|+ \big| \big( |x|^2K_{\sigma}\big)(x,w)\big| \lesssim \sigma^{-\frac{5}{2}}  \exp\Big( - \frac{|x-w|^2}{8 \tanh(\sigma)}  \Big),
\end{equation}
and 
\begin{equation} \label{Bb2}
\big| \big( \nabla_{x}K_{\sigma}\big)(x,w)\big| \lesssim \sigma^{-2}  \exp\Big( - \frac{|x-w|^2}{8 \tanh(\sigma)}  \Big).
\end{equation}
Using  \eqref{Bb0}, \eqref{Bb1} and \eqref{Bb2}, together with the elementary bound
$$ \int dw \,  \exp\Big( - \frac{|x-w|^2}{8 \tanh(r)}  \Big)\lesssim (\tanh(r))^{\frac32}\lesssim r^{\frac32},$$
we deduce for $r,\si,\eta>0$
\begin{multline*}
\Big| \int dw \, H\big(K_r(.,w) K_\eta(.,w) K_\sigma(.,w)\big)\Big| \lesssim \\
\begin{aligned}
 &\lesssim \Big(\frac{1}{r^{\frac52}\sigma^{\frac32}\eta^{\frac32}}+\frac{1}{r^{2}\sigma^{2}\eta^{\frac32}}+\frac{1}{r^{2}\sigma^{\frac32}\eta^{2}}+ \frac{1}{r^{\frac32}\sigma^{\frac32}\eta^{\frac52}}+ \frac{1}{r^{\frac32}\sigma^{\frac52}\eta^{\frac32}}+  \frac{1}{r^{\frac32}\sigma^{2}\eta^{2}}\Big) \int dw \,  \exp\Big( - \frac{|x-w|^2}{8 \tanh(r)}  \Big)\\
&\lesssim   \frac{1}{r\sigma^{\frac32}\eta^{\frac32}}+\frac{1}{r^{\frac12}\sigma^{2}\eta^{\frac32}}+\frac{1}{r^{\frac12}\sigma^{\frac32}\eta^{2}}+ \frac{1}{\sigma^{\frac32}\eta^{\frac52}}+ \frac{1}{ \sigma^{\frac52}\eta^{\frac32}}+  \frac{1}{ \sigma^{2}\eta^{2}} .
\end{aligned}
\end{multline*}
This estimate together with \eqref{intter} and \eqref{firs} implies
 \begin{multline} \label{eqreff}
\Big\| H^{\frac14+\eps}\Big( \int dw \,    K_r(.,w) K_\sigma(.,w) K_\eta(.,w)\Big)\Big\|_{L^{\infty}} \lesssim \\
 \begin{aligned}
 &\lesssim \frac{1}{(\sigma \eta)^{\frac32(\frac34-\eps)}}    \Big(   \frac{1}{r\sigma^{\frac32}\eta^{\frac32}}+\frac{1}{r^{\frac12}\sigma^{2}\eta^{\frac32}}+\frac{1}{r^{\frac12}\sigma^{\frac32}\eta^{2}}+ \frac{1}{\sigma^{\frac32}\eta^{\frac52}}+ \frac{1}{ \sigma^{\frac52}\eta^{\frac32}}+  \frac{1}{ \sigma^{2}\eta^{2}}\Big)^{\frac14+\eps} \\
 &\lesssim \frac{1}{r^{\frac14+\eps}\sigma^{\frac32}\eta^{\frac32}}+\frac{1}{r^{\frac18+\frac{\eps}2}\sigma^{\frac{13}8+\frac{\eps}2}\eta^{\frac32}}+\frac{1}{r^{\frac18+\frac{\eps}2}\sigma^{\frac32}\eta^{\frac{13}8+\frac{\eps}2}}+ \frac{1}{\sigma^{\frac32}\eta^{\frac74+\eps}}+ \frac{1}{ \sigma^{\frac74+\eps}\eta^{\frac32}}+  \frac{1}{ \sigma^{\frac{13}8+\frac{\eps}2}\eta^{\frac{13}8+\frac{\eps}2}}.
 \end{aligned}
 \end{multline}
Coming back to \eqref{Minko} and by integration of the previous line
 \begin{equation*} 
\|A^{(n)}_{s,t}\|_{\mathcal{W}^{\frac12+2\eps,\infty}} \lesssim  \int_{t-s}^{+\infty} \frac{dr}{r^{\frac54+\eps}}\lesssim \frac{1}{|t-s|^{\frac14+\eps}},
\end{equation*}
which implies the result if $\eps>0$ is small enough.

\smallskip

We now bound the contribution of $B^{(n)}_{s,t}$. Similarly to  \eqref{Minko} we have 
 \begin{equation*}
\|B^{(n)}_{s,t}\|_{\mathcal{W}^{\frac12+2\eps,\infty}}\lesssim  \int_0^{t-s} dr  \int_{r+2(t-s-r)}^{+\infty} d\sigma  \int_{r}^{+\infty} d\eta \,\Big\| H^{\frac14+\eps}\Big(  \int dw \,    K_r(.,w) K_\sigma(.,w) K_\eta(.,w)\Big)\Big\|_{L^{\infty}} ,
\end{equation*}
and using \eqref{eqreff} we deduce 
 \begin{equation*} 
\|B^{(n)}_{s,t}\|_{\mathcal{W}^{\frac12+2\eps,\infty}} \lesssim  \int_0^{t-s}dr\int_{r+2(t-s-r)}^{+\infty} d\sigma \, \Big(\frac{1}{r^{\frac34+\eps}\sigma^{\frac32}}+\frac{1}{r^{\frac58+\frac{\eps}2}\sigma^{\frac{13}8+\frac{\eps}2}}+  \frac{1}{r^{\frac12} \sigma^{\frac74+\eps}} \Big).
\end{equation*}
Then, since $r+2(t-s-r) \geq t-s$ for $0 \leq r \leq t-s$, we obtain
 \begin{align*} 
\|B^{(n)}_{s,t}\|_{\mathcal{W}^{\frac12+2\eps,\infty}} \lesssim & \int_0^{t-s}dr\int_{t-s}^{+\infty} d\sigma \, \Big(\frac{1}{r^{\frac34+\eps}\sigma^{\frac32}}+\frac{1}{r^{\frac58+\frac{\eps}2}\sigma^{\frac{13}8+\frac{\eps}2}}+  \frac{1}{r^{\frac12} \sigma^{\frac74+\eps}} \Big)\\
\lesssim & \int_0^{t-s}dr\,   \Big(\frac{1}{r^{\frac34+\eps}(t-s)^\frac12}+\frac{1}{r^{\frac58+\frac{\eps}2}(t-s)^{\frac{5}8+\frac{\eps}2}}+  \frac{1}{r^{\frac12} (t-s)^{\frac34+\eps}} \Big) \\
\lesssim & \frac{1}{(t-s)^{\frac14+\eps}},
\end{align*}
which suffices for $\eps>0$ small. The study of the term $C^{(n)}_{s,t}$ is similar.
\end{proof}

\

Based on the expressions derived in the previous proof, the limits of the quantities $\overline{\frakc}^{\mathbf{1},(n)},\overline{\frakc}^{\mathbf{2},(n)}$ (as $n\to\infty$) are easily identified as follows. For $0\leq s<t$, we define the ($n$-independent) functions
\begin{equation}\label{def-c1-limit}
\overline{\frakc}^{\mathbf{1}}_{s,t}(x)
:=\frac12\int_{2(t-s)}^{+\infty}d\si\,  K_\si(x,x),
\end{equation}
and
\begin{align}\label{def-c2-limit}
\overline{\frakc}^{\mathbf{2}}_{s,t}(x)
&:=\frac12 \int_{t-s}^{+\infty} dr \int dw \, K_{r}(x,w)
    \bigg( \int_{r}^{+\infty} d\sigma \, K_{\sigma}(x,w) \bigg)^2\nonumber\\
&+\frac12 \int_0^{t-s} dr \int dw \, K_{r}(x,w) \bigg( \int_{r+2(t-s-r)}^{+\infty} d\sigma \, K_{\sigma}(x,w) \bigg)\bigg( \int_{r}^{+\infty} d\eta \, K_{\eta}(x,w) \bigg)\nonumber\\
&+\frac12 \int_0^{t-s} dr \int dw \, K_{r}(x,w) \bigg( \int_{r+2(t-s-r)}^{+\infty} d\sigma \, K_{\sigma}(x,w)  \bigg)\bigg( \int_{r}^{r+2(t-s-r)} d\eta \, K_{\eta}(x,w) \bigg)\nonumber\\
&=:A_{s,t}(x)+B_{s,t}(x)+C_{s,t}(x),
\end{align}
together with $\overline{\frakc}_{s,t}:=3\,\overline{\frakc}^{\mathbf{1}}_{s,t}-9\,\overline{\frakc}^{\mathbf{2}}_{s,t}$.  We also set (arbitrarily) $\overline{\frakc}^{\mathbf{1}}_{t,t}=\overline{\frakc}^{\mathbf{2}}_{t,t}=\overline{\frakc}_{t,t}=0$ for all $t\geq 0$.

\begin{proposition}\label{prop:c-convergence}
For every $\eps>0$ small enough, it holds that
$$ \sup_{s\geq 0}\Big(\big\llbracket  \overline{\frakc}^{\mathbf{1},(n)}_{s,.}-\overline{\frakc}^{\mathbf{1}}_{s,.}\big\rrbracket_{\eps,[s,s+1]}+\big\llbracket  \overline{\frakc}^{\mathbf{2},(n)}_{s,.}-\overline{\frakc}^{\mathbf{2}}_{s,.}\big\rrbracket_{\eps,[s,s+1]}\Big)\stackrel{n\to\infty}{\longrightarrow} 0.$$

\end{proposition}

\begin{proof}

Throughout, $0\leq s<t$ with $|t-s|\leq 1$ and $\eps_n=2^{-n}$, and we record
the following elementary estimate for further use: for any $\beta>1$, any $a>0$ and any $\theta\in(0,1]$,
\begin{equation}\label{slab-elem}
\int_{a}^{a+2\eps_n}\frac{d\si}{\si^{\beta}}
\;\leq\; \int_{0}^{2\eps_n}\frac{d\si}{(a+\si)^{\beta}}\;\leq\; \frac{1}{a^{\,\beta-(1-\theta)}}\int_{0}^{2\eps_n}\frac{d\si}{\si^{1-\theta}}
\;\lesssim\; \frac{\eps_n^{\theta}}{a^{\,\beta-(1-\theta)}}.
\end{equation}

\

\noindent
1. \textit{Convergence of $\overline{\frakc}^{\mathbf{1},(n)}-\overline{\frakc}^{\mathbf{1}}$.}
From \eqref{overline-frakc-1} and \eqref{def-c1-limit}, one has
\begin{equation*}
\overline{\frakc}^{\mathbf{1},(n)}_{s,t}(x)-\overline{\frakc}^{\mathbf{1}}_{s,t}(x)
=\frac12\int_{2(t-s+\eps_n)}^{+\infty}d\si\, K_\si(x,x)
-\frac12\int_{2(t-s)}^{+\infty}d\si\, K_\si(x,x)
=-\frac12\int_{2(t-s)}^{2(t-s+\eps_n)}d\si\, K_\si(x,x).
\end{equation*}
Using $\|K_\si(\cdot,\cdot)\|_{L^\infty}\lesssim \si^{-\frac32}$, we obtain for the $L^\infty$ norm
\begin{equation*}
\big\|\overline{\frakc}^{\mathbf{1},(n)}_{s,t}-\overline{\frakc}^{\mathbf{1}}_{s,t}\big\|_{L^\infty}
\lesssim \int_{2(t-s)}^{2(t-s)+2\eps_n}\frac{d\si}{\si^{\frac32}}
\lesssim \frac{\eps_n^{\frac{1}{16}}}{|t-s|^{\frac{9}{16}}}
\end{equation*}
where we have used \eqref{slab-elem} with $\beta=\frac32$, $a=2(t-s)$ and $\theta=\frac{1}{16}$ (so that $\beta-(1-\theta)=\frac{9}{16}$).
Likewise, using the same computation as in \eqref{expli-H}, we obtain for the $\cw^{2,\infty}$ (hence $\cb^{2-\eps}_\infty$) norm,
\begin{equation*}
\big\|H\big(\overline{\frakc}^{\mathbf{1},(n)}_{s,t}-\overline{\frakc}^{\mathbf{1}}_{s,t}\big)\big\|_{L^\infty}
\lesssim \int_{2(t-s)}^{2(t-s)+2\eps_n}\frac{d\si}{(\sinh 2\si)^{\frac32}}\Big(\tanh\si+\tfrac{1}{\tanh\si}\Big)
\lesssim \int_{2(t-s)}^{2(t-s)+2\eps_n}\frac{d\si}{\si^{\frac52}}
\lesssim \frac{\eps_n^{\frac{1}{16}}}{|t-s|^{\frac{25}{16}}}
\end{equation*}
where we have used \eqref{slab-elem} with $\beta=\frac52$, $a=2(t-s)$ and
$\theta=\frac{1}{16}$ (so that $\beta-(1-\theta)=\frac{25}{16}$).

\smallskip

Interpolating exactly as in the proof of Lemma~\ref{lem:c}, we obtain for $\eps>0$ small enough
\begin{equation*}
\big\|\overline{\frakc}^{\mathbf{1},(n)}_{s,t}-\overline{\frakc}^{\mathbf{1}}_{s,t}\big\|_{\cb^{\frac12+\eps}_\infty}
\lesssim \Big(\frac{\eps_n^{\frac{1}{16}}}{|t-s|^{\frac{9}{16}}}\Big)^{\frac{3-4\eps}{4-2\eps}}
\Big(\frac{\eps_n^{\frac{1}{16}}}{|t-s|^{\frac{25}{16}}}\Big)^{\frac{1+2\eps}{4-2\eps}}
\lesssim \frac{\eps_n^{\frac{1}{16}}}{|t-s|^{\frac{7}{8}}},
\end{equation*}
due to $\tfrac{9}{16}\cdot\tfrac{3-4\eps}{4-2\eps}+\tfrac{25}{16}\cdot\tfrac{1+2\eps}{4-2\eps}\leq \tfrac{7}{8}$ for any $\eps$ small enough.

\smallskip

We have thus established that 
$$\sup_{s\geq 0}\ \llbracket\overline{\frakc}^{\mathbf{1},(n)}_{s,.}-\overline{\frakc}^{\mathbf{1}}_{s,.}\rrbracket_{\eps,[s,s+1]}\lesssim \eps_n^{\frac{1}{16}}\stackrel{n\to\infty}{\longrightarrow} 0.$$

\smallskip

\noindent
2. \textit{Convergence of $\overline{\frakc}^{\mathbf{2},(n)}-\overline{\frakc}^{\mathbf{2}}$.}
Based on the decompositions \eqref{decompo-c-2-c-2} for $\overline{\frakc}^{\mathbf{2},(n)}_{s,t}$ and \eqref{def-c2-limit} for $\overline{\frakc}^{\mathbf{2}}_{s,t}$, it suffices to bound successively $A^{(n)}-A$, $B^{(n)}-B$ and $C^{(n)}-C$. We will only focus on bounding the first two differences, and let the reader check that $C^{(n)}-C$ can be treated with similar arguments.

\smallskip

To start with, let us write
\small
\begin{align*}
&\big(A^{(n)}_{s,t}-A_{s,t}\big)(x)=\frac12 \int_{t-s}^{+\infty} dr \int dw \, K_{r}(x,w)
    \bigg( \int_{r+2\varepsilon_n}^{+\infty} d\sigma \, K_{\sigma}(x,w) \bigg)^2\\
		&\hspace{1cm}-\frac12 \int_{t-s}^{+\infty} dr \int dw \, K_{r}(x,w)
    \bigg( \int_{r}^{+\infty} d\sigma \, K_{\sigma}(x,w) \bigg)^2\\
&=-\frac12 \int_{t-s}^{+\infty} dr\int_r^{r+2\varepsilon_n} d\sigma\int_{r}^{+\infty} d\eta \int dw \, K_{r}(x,w) K_{\sigma}(x,w) K_{\eta}(x,w) \\
		&\hspace{1cm}-\frac12 \int_{t-s}^{+\infty} dr\int_r^{r+2\varepsilon_n} d\sigma\int_{r+2\eps_n}^{+\infty} d\eta \int dw \, K_{r}(x,w) K_{\sigma}(x,w) K_{\eta}(x,w),
\end{align*}
\normalsize
and
\small
\begin{align*}
&\big(B^{(n)}_{s,t}-B_{s,t}\big)(x)=\frac12 \int_0^{t-s} dr \int dw \, K_{r}(x,w) \bigg( \int_{r+2(t-s-r)+2\ep_n}^{+\infty} d\sigma \, K_{\sigma}(x,w) \bigg)\bigg( \int_{r+2\varepsilon_n}^{+\infty} d\eta \, K_{\eta}(x,w) \bigg)\\
&\hspace{1cm}-\frac12 \int_0^{t-s} dr \int dw \, K_{r}(x,w) \bigg( \int_{r+2(t-s-r)}^{+\infty} d\sigma \, K_{\sigma}(x,w) \bigg)\bigg( \int_{r}^{+\infty} d\eta \, K_{\eta}(x,w) \bigg) \\
&=-\frac12 \int_0^{t-s} dr\int_{r+2(t-s-r)}^{r+2(t-s-r)+2\eps_n} d\sigma \int_{r+2\varepsilon_n}^{+\infty} d\eta\int dw \, K_{r}(x,w)  K_{\sigma}(x,w) K_{\eta}(x,w) \\
&\hspace{1cm}-\frac12 \int_0^{t-s} dr\int_{r+2(t-s-r)}^{+\infty} d\sigma  \int_{r}^{r+2\eps_n} d\eta\int dw \, K_{r}(x,w) K_{\sigma}(x,w)  K_{\eta}(x,w) =:\mathbf{B}^{1,(n)}_{s,t}(x)+\mathbf{B}^{2,(n)}_{s,t}(x).
\end{align*}
\normalsize

\

\noindent
\textit{(a) $L^\infty$ bound.} For $A^{(n)}-A$, using $\|K_\si(\cdot,\cdot)\|_{L^\infty}\lesssim \si^{-\frac32}$ and $\int dw\, K_r(x,w)\lesssim 1$, we obtain
\begin{equation*}
\big|A^{(n)}_{s,t}(x)-A_{s,t}(x)\big|
\lesssim \int_{t-s}^{+\infty}dr\,\Big(\int_r^{r+2\eps_n}\frac{d\si}{\si^{\frac32}}\Big)
\Big(\int_r^{+\infty}\frac{d\si}{\si^{\frac32}}\Big)
\lesssim \eps_n^{\frac18}\int_{t-s}^{+\infty}\frac{dr}{r^{\frac98}}
\lesssim \frac{\eps_n^{\frac18}}{|t-s|^{\frac18}},
\end{equation*}
where we have used \eqref{slab-elem} with $\beta=\frac32$, $a=r$ and $\theta=\frac{1}{8}$ (so that $\beta-(1-\theta)=\frac{5}{8}$).

\smallskip

For $B^{(n)}-B$, one has for similar reasons
\begin{align*}
&\big|B^{(n)}_{s,t}(x)-B_{s,t}(x)\big|\lesssim \\
&\lesssim \int_{0}^{t-s}dr\,\Big(\int_{r+2(t-s-r)}^{r+2(t-s-r)+2\eps_n}\frac{d\si}{\si^{\frac32}}\Big)\Big(\int_r^{+\infty}\frac{d\eta}{\eta^{\frac32}}\Big)+\int_{0}^{t-s}dr\,\Big(\int_{r+2(t-s-r)}^{+\infty}\frac{d\si}{\si^{\frac32}}\Big)\Big(\int_r^{r+2\eps_n}\frac{d\eta}{\eta^{\frac32}}\Big)\\
&\lesssim \eps_n^{\frac18}\bigg(\int_{0}^{t-s}\frac{dr}{(r+2(t-s-r))^{\frac58} r^{\frac12}}+\int_{0}^{t-s}\frac{dr}{ (r+2(t-s-r))^{\frac12} r^{\frac58}}\bigg)\lesssim \frac{\eps_n^{\frac18}}{|t-s|^{\frac18}},
\end{align*}
and we can thus conclude that
\begin{equation}\label{overline-c2-inf}
\big\| \overline{\frakc}_{s,t}^{\mathbf{2},(n)}-\overline{\frakc}_{s,t}^{\mathbf{2}}\big\|_{L^\infty} \lesssim \frac{\eps_n^{\frac18}}{|t-s|^{\frac18}}.
\end{equation}

\

\noindent
\textit{(b) $\cb^{\frac12+\eps}_\infty$ bound.} For $A^{(n)}-A$, we can rely on \eqref{eqreff} to assert that
\begin{align*}
&\|A^{(n)}_{s,t}-A_{s,t}\|_{\mathcal{W}^{\frac12+2\eps,\infty}}\lesssim \int_{t-s}^{+\infty} dr\int_r^{r+2\varepsilon_n} d\sigma\int_{r}^{+\infty} d\eta\\
& \bigg[\frac{1}{r^{\frac14+\eps}\sigma^{\frac32}\eta^{\frac32}}+\frac{1}{r^{\frac18+\frac{\eps}2}\sigma^{\frac{13}8+\frac{\eps}2}\eta^{\frac32}}+\frac{1}{r^{\frac18+\frac{\eps}2}\sigma^{\frac32}\eta^{\frac{13}8+\frac{\eps}2}}+ \frac{1}{\sigma^{\frac32}\eta^{\frac74+\eps}}+ \frac{1}{ \sigma^{\frac74+\eps}\eta^{\frac32}}+  \frac{1}{ \sigma^{\frac{13}8+\frac{\eps}2}\eta^{\frac{13}8+\frac{\eps}2}}\bigg]\\
&\lesssim \eps_n^{\frac18} \int_{t-s}^{+\infty} dr \bigg[\frac{1}{r^{\frac14+\eps}r^{\frac58}r^{\frac12}}+\frac{1}{r^{\frac18+\frac{\eps}2}r^{\frac34+\frac{\eps}2}r^{\frac12}}+\frac{1}{r^{\frac18+\frac{\eps}2}r^{\frac58}r^{\frac{5}8+\frac{\eps}2}}+ \frac{1}{r^{\frac58}r^{\frac34+\eps}}+ \frac{1}{ r^{\frac78+\eps}r^{\frac12}}+  \frac{1}{ r^{\frac34+\frac{\eps}2}r^{\frac{5}8+\frac{\eps}2}}\bigg]\\
&\lesssim \eps_n^{\frac18} \int_{t-s}^{+\infty} \frac{dr}{r^{\frac{11}{8}+\eps}}\lesssim \frac{\eps_n^{\frac18}}{|t-s|^{\frac{3}{8}+\eps}}.
\end{align*}
Along the same lines,
\small
\begin{align*}
&\|\mathbf{B}^{1,(n)}_{s,t}\|_{\mathcal{W}^{\frac12+2\eps,\infty}}\lesssim \int_0^{t-s} dr\int_{r+2(t-s-r)}^{r+2(t-s-r)+2\eps_n} d\sigma \int_{r+2\varepsilon_n}^{+\infty} d\eta\\
& \bigg[\frac{1}{r^{\frac14+\eps}\sigma^{\frac32}\eta^{\frac32}}+\frac{1}{r^{\frac18+\frac{\eps}2}\sigma^{\frac{13}8+\frac{\eps}2}\eta^{\frac32}}+\frac{1}{r^{\frac18+\frac{\eps}2}\sigma^{\frac32}\eta^{\frac{13}8+\frac{\eps}2}}+ \frac{1}{\sigma^{\frac32}\eta^{\frac74+\eps}}+ \frac{1}{ \sigma^{\frac74+\eps}\eta^{\frac32}}+  \frac{1}{ \sigma^{\frac{13}8+\frac{\eps}2}\eta^{\frac{13}8+\frac{\eps}2}}\bigg]\\
&\lesssim \eps_n^{\frac18} \int_{0}^{t-s} dr \bigg[\frac{1}{r^{\frac14+\eps}(r+2(t-s-r))^{\frac58}r^{\frac12}}+\frac{1}{r^{\frac18+\frac{\eps}2}(r+2(t-s-r))^{\frac34+\frac{\eps}2}r^{\frac12}}+\frac{1}{r^{\frac18+\frac{\eps}2}(r+2(t-s-r))^{\frac58}r^{\frac{5}8+\frac{\eps}2}}\\
&\hspace{2cm}+ \frac{1}{(r+2(t-s-r))^{\frac58}r^{\frac34+\eps}}+ \frac{1}{ (r+2(t-s-r))^{\frac78+\eps}r^{\frac12}}+  \frac{1}{ (r+2(t-s-r))^{\frac34+\frac{\eps}2}r^{\frac{5}8+\frac{\eps}2}}\bigg]\\
&\lesssim \eps_n^{\frac18} \int_{0}^{t-s} dr \bigg[\frac{1}{r^{\frac34+\eps}(r+2(t-s-r))^{\frac58}}+\frac{1}{r^{\frac58+\frac{\eps}2}(r+2(t-s-r))^{\frac34+\frac{\eps}2}}+ \frac{1}{r^{\frac12} (r+2(t-s-r))^{\frac78+\eps}}\bigg]\\
&\lesssim \frac{\eps_n^{\frac18}}{|t-s|^{\frac{3}{8}+\eps}},
\end{align*}
 \normalsize
and likewise  
\small
\begin{align*}
&\|\mathbf{B}^{2,(n)}_{s,t}\|_{\mathcal{W}^{\frac12+2\eps,\infty}}\lesssim \int_0^{t-s} dr\int_{r+2(t-s-r)}^{+\infty} d\sigma  \int_{r}^{r+2\eps_n} d\eta\\
& \bigg[\frac{1}{r^{\frac14+\eps}\sigma^{\frac32}\eta^{\frac32}}+\frac{1}{r^{\frac18+\frac{\eps}2}\sigma^{\frac{13}8+\frac{\eps}2}\eta^{\frac32}}+\frac{1}{r^{\frac18+\frac{\eps}2}\sigma^{\frac32}\eta^{\frac{13}8+\frac{\eps}2}}+ \frac{1}{\sigma^{\frac32}\eta^{\frac74+\eps}}+ \frac{1}{ \sigma^{\frac74+\eps}\eta^{\frac32}}+  \frac{1}{ \sigma^{\frac{13}8+\frac{\eps}2}\eta^{\frac{13}8+\frac{\eps}2}}\bigg]\\
&\lesssim \eps_n^{\frac18} \int_{0}^{t-s} dr \bigg[\frac{1}{r^{\frac14+\eps}(r+2(t-s-r))^{\frac12}r^{\frac58}}+\frac{1}{r^{\frac18+\frac{\eps}2}(r+2(t-s-r))^{\frac58+\frac{\eps}2}r^{\frac58}}+\frac{1}{r^{\frac18+\frac{\eps}2}(r+2(t-s-r))^{\frac12}r^{\frac34+\frac{\eps}2}}\\
&\hspace{2cm}+ \frac{1}{(r+2(t-s-r))^{\frac12}r^{\frac78+\eps}}+ \frac{1}{ (r+2(t-s-r))^{\frac34+\eps}r^{\frac58}}+  \frac{1}{ (r+2(t-s-r))^{\frac58+\frac{\eps}2}r^{\frac34+\frac{\eps}2}}\bigg]\\
&\lesssim \eps_n^{\frac18} \int_{0}^{t-s} dr \bigg[\frac{1}{r^{\frac78+\eps}(r+2(t-s-r))^{\frac12}}+\frac{1}{r^{\frac34+\frac{\eps}2}(r+2(t-s-r))^{\frac58+\frac{\eps}2}}+ \frac{1}{r^{\frac58} (r+2(t-s-r))^{\frac34+\eps}}\bigg]\\
&\lesssim \frac{\eps_n^{\frac18}}{|t-s|^{\frac{3}{8}+\eps}}.
\end{align*}
 \normalsize

\smallskip

Gathering the bounds on $\mathbf{B}^{1,(n)}$ and $\mathbf{B}^{2,(n)}$, together with the analogous estimate on $A^{(n)}-A$, we obtain
$$\big\|\overline{\frakc}^{\mathbf{2},(n)}_{s,t}-\overline{\frakc}^{\mathbf{2}}_{s,t}\big\|_{\cb^{\frac12+\eps}_\infty}\lesssim \frac{\eps_n^{\frac18}}{|t-s|^{\frac38+\eps}},$$
which, combined with \eqref{overline-c2-inf}, entails
$$\sup_{s\geq 0}\ \llbracket\overline{\frakc}^{\mathbf{2},(n)}_{s,.}-\overline{\frakc}^{\mathbf{2}}_{s,.}\rrbracket_{\eps,[s,s+1]}\lesssim \eps_n^{\frac{1}{8}}\stackrel{n\to\infty}{\longrightarrow} 0.$$
\end{proof}

\subsection{Reformulation of the dynamics}\label{Sect32}

\

\smallskip

 We fix two times $s,T_\star$ such that $0\leq s< T_\star\leq s+1$, together with $Z_{s,.}\in \mathcal Z_{\eps,  [s,T_\star]}$, $\overline{\frakc}_{s,.}:=3\overline{\frakc}^{\mathbf{1}}_{s,.}-9\overline{\frakc}^{\mathbf{2}}_{s,.}$ satisfying $\llbracket \overline{\frakc}^{\mathbf{i}}_{s,.}\rrbracket_{3\eps,[s,T_\star]}<\infty$ and $u\in \cb_{\infty}^{-\frac12-\eps}$.  The interval $[s,T_\star]$ is the range on which the diagrams $Z_{s,\cdot}$ are defined and controlled, whereas $T$ denotes the existence time of the solution itself; in particular $s\leq T\leq T_\star$. 

\medskip

In this setting, we denote by
$$ (v,w)=(v^u_{s,.},w^u_{s,.})$$
 the solution of the following auxiliary system: for all $(r,T)$ with $s\leq r\leq T\leq T_\star$,
 \begin{equation} \label{syst1}
\left\{
\begin{aligned}
& \partial_t v + Hv  = F(v+w;Z_{s,.})-\mu v , \quad s\leq r \leq T\\
& \partial_t w + Hw  = 3\, \<Psi2IPsi3>_{s,.}+G(v,w;Z_{s,.}{,\overline{\frakc}_{s,.}})+\mu v,\\
& (v,w)(s)=(0,u).
\end{aligned}
\right.
\end{equation}
where 
\begin{equation*} 
 F(v+w;Z_{s,.})_{r} =F(v_r+w_r;Z_{s,r}) = -3(v_r+w_r -\<IPsi3>_{s,r}) \pl  \<Psi2>_{s,r}
\end{equation*}
\begin{multline}\label{definition-g}
G(v,w;Z_{s,.} ,\overline{\frakc}_{s,.})_{r} =- (v_r+w_r)^3 -3\,  \mathrm{com}^{Z_{s,.}} (v,w)_{s,r}-3\, w_r \pe \<Psi2>_{s,r}  \\
-3(v_r+w_r-\<IPsi3>_{s,r}) \pg \<Psi2>_{s,r}+ P(v_r+w_r; Z_{s,r})+ {\overline{\frakc}_{s,r}}\big(\<Psi>_{s,r}-\<IPsi3>_{s,r}+v_r+w_r\big)
\end{multline}
with
\begin{align}
\mathrm{com}^{Z_{s,.}} (v,w)_{s,r}&: = \mathrm{com}^{Z_{s,.}}_1 (v,w)_{s,r} \pe  \<Psi2>_{s,r} +\mathrm{com}^{Z_{s,.}}_2 (v+w)_{s,r}, \label{comm1} \\
\mathrm{com}^{Z_{s,.}}_1 (v,w)_{s,t} &:=-3\bigg[ \int_s^t e^{-(t-r)H} [(v_r+w_r-\<IPsi3>_{s,r}) \pl  \<Psi2>_{s,r}] \, dr 
-(v_t+w_t-\<IPsi3>_{s,t}) \pl  \<IPsi2>_{s,t}\bigg], \nonumber \\ 
\mathrm{com}^{Z_{s,.}}_2 (v+w)_{s,t}&:= \big[ \pl, \pe\big] \Big(-3\big(v_t+w_t-\<IPsi3>_{s,t}\big), \<IPsi2>_{s,t}, \<Psi2>_{s,t}\Big), \nonumber
\end{align}
and 
$$ P(v+w; Z_{s,.})_{r}=\tau^{(0)}(Z_{s,r}) +\tau^{(1)}(Z_{s,r}) \cdot(v_r+w_r) +\tau^{(2)}(Z_{s,r}) \cdot (v_r+w_r)^2$$
where
\begin{align}
\tau^{(0)}(Z_{s,r}) &:= (\<IPsi3>_{s,r})^3 -3\Big[ \<Psi>_{s,r}  \pg (\<IPsi3>_{s,r})^2  + \<Psi>_{s,r}  \pl (\<IPsi3>_{s,r})^2 +  \<Psi>_{s,r} \pe [\<IPsi3>_{s,r} \pe \<IPsi3>_{s,r}] \nonumber\\
&\hspace{1cm} +2 \<IPsi3>_{s,r}\,\<PsiIPsi3>_{s,r}+2[\pl, \pe](\<IPsi3>_{s,r}, \<IPsi3>_{s,r}, \<Psi>_{s,r})\Big]-9\, \<IPsi3>_{s,r}\, \<Psi2IPsi2>_{s,r} \label{deft0}\\
\tau^{(1)}(Z_{s,r}) &:=  6 \Big[\<IPsi3>_{s,r} \pg \<Psi>_{s,r} + \<IPsi3>_{s,r} \pl \<Psi>_{s,r}+\<PsiIPsi3>_{s,r} \Big] - 3 \,(\<IPsi3>_{s,r})^2+9\, \<Psi2IPsi2>_{s,r},\label{deft1} \\
\tau^{(2)}(Z_{s,r}) &:= -3\, \<Psi>_{s,r} +3\<IPsi3>_{s,r} .\label{deft2}
\end{align}

\

The following identification result (at the level of the approximated equation) can then be proved with the same arguments as in \cite[Proposition 3.5]{DFT}. 
 
\begin{lemma}\label{lem:identifi-solu}
Fix $\eps>0$ small enough and $n\geq 1$. Let $Z^{(n)}$ be the set of diagrams introduced in Section \ref{subsec:dia}, $\overline{\frakc}^{(n)}$ be defined as in Section \ref{subsec:ren-seq}, and assume that for all $s\geq 0$ and $u\in \cb_\infty^{-\frac12-\eps}$, the system \eqref{syst1} with $Z_{s,.}:=Z^{(n)}_{s,.}$, $\overline{\frakc}_{s,.}:=\overline{\frakc}^{(n)}_{s,.}$ admits a unique solution $(v^{(n),u}_{s,.},w^{(n),u}_{s,.})$ on $[s,s+1]$. 

\smallskip

Then the process $X^{(n)}$ defined (recursively)  by $X^{(n)}_0= u \in \cb^{-\frac12-\eps}_{\infty}$ and for all $0\leq s\leq t\leq s+1$
 \begin{equation}\label{decompos-x-n}
X^{(n)}_t=\<Psi>^{(n)}_{s,t}-\<IPsi3>^{(n)}_{s,t}+v_{s,t}^{(n),X^{(n)}_s}+w_{s,t}^{(n),X^{(n)}_s}
\end{equation}
is the unique solution of the problem
\begin{equation} \label{eq-bis}
\left\{
\begin{aligned}
&(\partial_t +H)  X^{(n)}  = -(X^{(n)})^3+ \frakc^{(n)}X^{(n)} +\xi^{(n)}, \quad t>0, \quad x \in \R^3, \\
&  X^{(n)}_0= u\in \cb^{-\frac12-\eps}_{\infty}, 
\end{aligned}
\right.
\end{equation}  
where $\xi^{(n)}$ denotes the spatially regularized version of the noise $\xi$ introduced in \eqref{regu-noise}.

\end{lemma}

\color{black}

\smallskip

\begin{remark}
Observe that the parameter $\mu >0$ introduced in \eqref{syst1} does not modify in the end the sum $v+w$ in the decomposition \eqref{decompos-x}. This parameter will however prove useful in the derivation of a priori bounds for $v$ and $w$ in the next sections.
\end{remark}

\medskip

\subsection{Local solution starting at $s$}

We now study the well-posedness of the  system \eqref{syst1}. The mild form of the equation reads  
  \begin{equation}  \label{mild:v}
\left\{
\begin{aligned}
& v_t =   \int_s^t d\tau \, e^{-(t-\tau)(H+\mu)} \big[F(v_\tau +w_\tau; Z_{s,\tau} ) -\mu v_\tau\big]  , \\
& w_t = e^{-(t-s)H} u + \int_s^t d\tau\, e^{-(t-\tau)H} [3\ \<Psi2IPsi3>_{s,\tau}+G(v,w; Z_{s,.},\overline{\frakc}_{s,.})_\tau+\mu v_\tau].
\end{aligned}
\right.
\end{equation}
 Let $s< T\leq T_\star$. For every small $\eps>0$, we consider the space 
\begin{equation*} 
\mathcal{X}_{\eps, [s,T]} := \Big\{ (v,w) \in \mathcal{C}\big([s,T]; \mathcal{B}_{\infty}^{ {-\frac23 +\eps}}\big) \times  \mathcal{C}\big([s,T]; \mathcal{B}_{\infty}^{{ -\frac12-\eps}}\big) \ \text{such that} \ \big\|(v,w)\big\|_{\eps,[s,T]}<\infty  \Big\}
\end{equation*}
where the norm $\big\|(v,w)\big\|_{\eps,[s,T]}=\big\|(v,w)\big\|_{\mathcal{X}_{\eps,[s,T]}}$ is defined by 
\begin{align}
\big\|(v,w)\big\|_{\eps,[s,T]} &:=  \max \Big\{ \sup_{s< t \le T} (t-s)^{\frac{7}{12}+\frac{\eps}{2}} \big\| v_t\big\|_{\mathcal{B}_{\infty}^{\frac12 +2\eps}}, 
\sup_{s \le t \le T} \big\| v_t \big\|_{\mathcal{B}_{\infty}^{-\frac23 +\eps}}, 
\sup_{s< r<t \le T} (r-s)^{\frac{7}{12}+\frac{\eps}{2}} \frac{ \big\|v_t-v_r\big\|_{\cb_{\infty}^{\frac32\eps}}} { |t-r|^{\frac14 +\frac{\eps}{4}} }, \nonumber\\
&\hspace{0.5cm} \sup_{s< t \le T} (t-s)^{\frac{3}{4}+\frac{3}{2}\eps} \big\|w_t \big\|_{\mathcal{B}_{\infty}^{1+2\eps}}, 
\sup_{s \le t \le T} \big\|w_t \big\|_{\mathcal{B}_{\infty}^{-\frac12 -\eps}}, 
\sup_{s< r<t \le T} (r-s)^{\frac{3}{4}+\frac{3}{2}\eps} \frac{ \big \|w_t-w_r \big\|_{\cb_{\infty}^{\eps}}} {|t-r|^{\frac12 +\frac{\eps}{2}}} \Big\}.\label{defino-x-ept}
\end{align}

We denote by 
$$\gga_{Z{,\overline{\frakc}}, w_0}:=(\gga_{Z}^{\mathbf{v}},\gga^{\mathbf{w}}_{Z{,\overline{\frakc}},w_0})$$
the map derived from the system \eqref{mild:v}, that is
\begin{align*}
\gga^{\mathbf{v}}_{Z}[v,w]_{s,t}:&= \int_s^t d\tau \, e^{-(t-\tau)(H+\mu)} \big[F(v_\tau +w_\tau; Z_{s,\tau} ) -\mu v_\tau\big]  \\
\gga^{\mathbf{w}}_{Z,\overline{\frakc},u}[v,w]_{s,t} :&= e^{-(t-s)H}u + \int_s^t d\tau\, e^{-(t-\tau)H} \Big[3\ \<Psi2IPsi3>_{0,\tau}+G(v,w; Z_{s,.},\overline{\frakc}_{s,.})_\tau+\mu v_\tau\Big].
\end{align*}

We are now in a position to state our main result regarding local well-posedness of the system~\eqref{syst1}. As in our previous work \cite{DFT}, we adapt the argument developed by Mourrat and Weber~\cite{MW}.

\begin{proposition}\label{prop:fixed-point}
Fix $0\leq s< T_\star \leq s+1$, and let $T \in (s, T_\star]$. For every $\eps>0$ small enough, 
there exists $\nu>0$ such that the following assertions hold.
\smallskip

\noindent
$(i)$ For all $Z_{s,.}\in \mathcal Z_{\frac{\eps}{2},[s,T_\star]}$, $\overline{\frakc}_{s,.}$ such that $\big\llbracket  \overline{\frakc}_{s,.}\big\rrbracket_{3\eps,[s,T_\star]}<\infty$ and $u\in  \mathcal{B}_{\infty}^{-\frac12-\eps}$, the map $\Gamma_{Z,\overline{\frakc}, u}[.,.]_{s,.}$ is well defined from $\mathcal{X}_{\eps,[s,T]}$ to $\mathcal{X}_{\eps,[s,T]}$ 
and for every $(v,w)\in \mathcal{X}_{\eps,[s,T]}$, one has 
\begin{multline*}
\big\| \gga_{Z,\overline{\frakc},u}[v,w]_{s,.}\big\|_{\eps,[s,T]} \leq \\
\leq C \big( 1+\big\| Z_{s,.}\big\|_{\mathcal Z_{\frac{\eps}{2},[s,T_\star]}}^3 +\big\llbracket  \overline{\frakc}_{s,.}\big\rrbracket_{3\eps,[s,T_\star]}^3 \big)\big( 1+\big\|u\big\|_{  \mathcal{B}^{-\frac12-\eps}_{\infty}}\big)\big(1+(T-s)^\nu \big\| (v,w)\big\|_{\eps,[s,T]}^{3} \big),
\end{multline*}
for some universal constant $C>0$.
\smallskip

\noindent
$(ii)$ Given two sets of diagrams $Z_{s,.},Z'_{s,.}\in \mathcal Z_{\frac{\eps}{2},[s,T_\star]}$, two functions $ \overline{\frakc}_{s,.}, \overline{\frakc}'_{s,.}$ for which $\big\llbracket  \overline{\frakc}_{s,.}\big\rrbracket_{3\eps,[s,T_\star]}$, $\big\llbracket  \overline{\frakc}'_{s,.}\big\rrbracket_{3\eps,[s,T_\star]}<\infty$, and $u\in  \mathcal{B}_{\infty}^{-\frac12-\eps}$, one has for all two elements $(v,w),(v',w')\in \mathcal{X}_{\eps,[s,T]}$,
\begin{multline*}
\big\| \gga_{Z, \overline{\frakc},u}[v,w]_{s,.}-\gga_{Z', \overline{\frakc}', u}[v',w']_{s,.}\big\|_{\eps,T}\leq \\
 \leq C (T-s)^\nu \Big(1+\big\| Z\big\|_{\mathcal Z_{\frac{\eps}{2},[s,T_\star]}}^2+\big\| Z'\big\|_{\mathcal Z_{\frac{\eps}{2},[s,T_\star]}}^2+ \big\llbracket  \overline{\frakc}\big\rrbracket_{3\eps,[s,T_\star]}^2+ \big\llbracket  \overline{\frakc}'\big\rrbracket_{3\eps,[s,T_\star]}^2\Big) 
\Big( 1+  \big\| (v,w)\big\|^2_{ \eps,[s,T]} + \big\| (v',w')\big\|^2_{ \eps,[s,T]} \Big)\\
 \times \Big\{\big\| (v,w)-(v',w')\big\|_{ \eps,[s,T]} +\big\| Z_{s,.}-Z'_{s,.}\big\|_{\mathcal Z_{\frac{\eps}{2},[s,T_\star]}} +\big\llbracket \overline{\frakc}- \overline{\frakc}'\big\rrbracket_{3\eps,[s,T_\star]}\Big\},
\end{multline*}
for some universal constant $C>0$.
\end{proposition}

\

We give the proof below. Once it is established, the local existence of the solution to \eqref{mild:v} follows from classical PDE arguments.
\begin{corollary}\label{cor:fixed-point}
Fix $0\leq s< T_\star \leq s+1$ and $\eps>0$ small enough. Then there exist $\zeta\geq 1$ and $C\in (0,1)$ such that for all $Z_{s,.}\in \mathcal Z_{\frac{\eps}{2},[s,T_\star]}$, $\overline{\frakc}_{s,.}$ with $\big\llbracket  \overline{\frakc}_{s,.}\big\rrbracket_{3\eps,[s,T_\star]}<\infty$ and $u\in  \mathcal{B}_{\infty}^{-\frac12-\eps}$, the system~\eqref{mild:v} admits a unique solution $(v,w)$ in $\mathcal{X}_{\eps,[s,s+\delta]}$, where
\begin{equation}
\delta:=C\, \Big( 1+\big\| Z_{s,.}\big\|_{\mathcal Z_{\frac{\eps}{2},[s,T_\star]}} +\big\llbracket  \overline{\frakc}_{s,.}\big\rrbracket_{3\eps,[s,T_\star]} \big)^{-\zeta}\big( 1+\big\|u\big\|_{  \mathcal{B}^{-\frac12-\eps}_{\infty}}\Big)^{-\zeta} \in (0,1).
\end{equation}
\end{corollary}

\

In accordance with Lemma \ref{lem:identifi-solu}, we define for the sake of clarity:
\begin{definition}\label{def:rough-sol}
Fix $0\leq T_0<T$ and $\eps>0$ small enough. Consider a family of diagrams $(Z_{s,t})_{T_0\leq s<t\leq T}$ such that $Z_{s,.}  \in \mathcal Z_{\frac{\eps}{2},[s,T]}$, as well as a family of functions $(\overline{\frakc}_{s,t})_{T_0\leq s<t\leq T}$ such that $\llbracket \overline{\frakc}_{s,.}\rrbracket_{3\eps,[s,T]} <\infty$.

\smallskip

A path $X\in \cac([T_0,T];\cb_\infty^{-\frac12-\eps})$ is called a (rough) solution of \eqref{eq:intro} on $[T_0,T]$ associated with $(Z,\overline{\frakc})$ and starting at $u\in \cb_\infty^{-\frac12-\eps}$ if $X_{T_0}=u$ and if, for all $T_0\leq s <t\leq (s+1)\wedge T$, $X_t$ admits the decomposition
 \begin{equation}\label{decompos-x}
X_t=\<Psi>_{s,t}-\<IPsi3>_{s,t}+v_{s,t}^{X_s}+w_{s,t}^{X_s}
\end{equation}
where $(v_{s,.}^{X_s},w_{s,.}^{X_s}) \in \mathcal{X}_{\eps,[s,(s+1)\wedge T]}$ is a solution of the system \eqref{mild:v} associated with $(Z_{s,.},\overline{\frakc}_{s,.})$ and starting at $X_s$.

\end{definition}

\

Let us finally specialize the above results to our initial problem, using in particular the convergence properties established in Propositions \ref{prop:conv-arbre} and \ref{prop:c-convergence}. \medskip

 Let $(Z,Z^{(n)})$ and $(\overline{\frakc},\overline{\frakc}^{(n)})$ be defined as in Propositions \ref{prop:conv-arbre} and \ref{prop:c-convergence}. In addition, set
 \begin{equation}\label{def-normes}
\big\| \! \big\| Z\big\| \! \big\|_{\mathcal Z_{\eps,[\tau_0,\tau_0+1]}}:=\sup_{s\in [\tau_0,\tau_0+1]} \big\| Z_{s,.}\big\|_{\mathcal Z_{\frac{\eps}{2},[s,\tau_0+1]}} \quad \text{and} \quad \big\llbracket\!\! \big\llbracket  \overline{\frakc}\big\rrbracket \!\! \big\rrbracket_{\eps,[\tau_0,\tau_0+1]}:=\sup_{s\in [\tau_0,\tau_0+1]}\big\llbracket  \overline{\frakc}_{s,.}\big\rrbracket_{3\eps,[s,\tau_0+1]} .
 \end{equation}

\begin{corollary}\label{coro:main}
Fix $\eps>0$ small enough, $\tau_0\geq 0$ and $u\in \cb_{\infty}^{-\frac12-\eps}(\R^3)$. Then there exist $\zeta\geq 1$ and $C\in (0,1)$ such that the following assertions hold.

\smallskip

\noindent
$(i)$ Almost surely, the equation \eqref{eq:intro} associated with $(Z,\overline{\frakc})$ and $u\in \cb_\infty^{-\frac12-\eps}$ admits a unique (rough) solution $X$ on $[\tau_0,\tau_0+\delta_{\tau_0,u}]$, where
\begin{equation}
\delta_{\tau_0,u}:=C\, \Big( 1+\big\|u\big\|_{  \mathcal{B}^{-\frac12-\eps}_{\infty}}\Big)^{-\zeta}\Big( 1+\big\| \! \big\| Z\big\| \! \big\|_{\mathcal Z_{\eps,[\tau_0,\tau_0+1]}} +\big\llbracket\!\! \big\llbracket  \overline{\frakc}\big\rrbracket \!\! \big\rrbracket_{\eps,[\tau_0,\tau_0+1]} \Big)^{-\zeta}\in (0,1].
\end{equation}

\smallskip

\noindent
$(ii)$ The sequence $(X^{(n)})$ of approximated solutions with initial condition $u$ at $\tau_0$ converges to $X$ in the space $\cac([\tau_0,\tau_0+\delta_{\tau_0}'];\cb_\infty^{-\frac12-\eps})$, where
\begin{equation}
\delta_{\tau_0,u}':=C\,\Big( 1+\big\|u\big\|_{  \mathcal{B}^{-\frac12-\eps}_{\infty}}\Big)^{-\zeta} \inf_{n\geq 1}\Big(\Big( 1+\big\| \! \big\| Z^{(n)}\big\| \! \big\|_{\mathcal Z_{\eps,[\tau_0,\tau_0+1]}} +\big\llbracket\!\! \big\llbracket  \overline{\frakc}^{(n)}\big\rrbracket \!\! \big\rrbracket_{\eps,[\tau_0,\tau_0+1]}\Big)^{-\zeta}\Big).
\end{equation}
\end{corollary}
\color{black}

\

\subsection{Proof of Proposition \ref{prop:fixed-point}}\label{proo-1}

\

\smallskip

In order to simplify the presentation, we consider only the case $ s=0$ and $T_\star=1$, but the argument extends to all $0\leq s< T_\star \leq s+1$. Moreover, we will focus on the proof of item $(i)$, and we leave it to the reader to verify that the assertion in item $(ii)$ can be established using entirely similar arguments.

\smallskip

Let $\eps>0$ (small enough), $0<T \leq 1$, $Z_{0,.}\in \mathcal Z_{\frac{\eps}{2},[0,1]}$,  {$\overline{\frakc}_{0,.}$ such that $\big\llbracket  \overline{\frakc}_{0,.}\big\rrbracket_{3\eps,[0,1]}<\infty$,} $w_0 \in   \mathcal{B}_{\infty}^{{ -\frac12-\eps} }$ and $(v,w)\in \mathcal{X}_{\eps,[0,T]}$. We bound the two components $\gga^{\mathbf{v}}_Z$ and $\gga^{\mathbf{w}}_{Z,w_0}$ of $\gga_{Z,w_0}$ separately.

\smallskip

For clarity, we set in what follows
$$  \big\|Z\big\|_{\eps}:=\big\|Z_{0,.}\big\|_{\mathcal Z_{\frac{\eps}{2},[0,1]}}, \quad \big\llbracket \overline{\frakc}\big\rrbracket_{\eps}:=\big\llbracket  \overline{\frakc}_{0,.}\big\rrbracket_{3\eps,[0,1]} \quad \text{and} \quad \big\|(v,w)\big\|_{\eps,T}:=\big\|(v,w)\big\|_{\mathcal{X}_{\eps,[0,T]}}.$$

For $0\leq \theta, \gamma \leq 1$ and $\tau \in \R$ we define the space
$$ 
 \big\| v \big\|_{\mathcal{D}_T^{\theta, \tau, \gamma}} :=\max\Big\{  \sup_{0< t \le T} t^{\theta} \big\| v_t\big\|_{\mathcal{B}_{\infty}^{\tau}},  \quad \sup_{0\le t\le T} \big\|v_t\big\|_{\mathcal{B}_{\infty}^{\tau-2\theta}},
 \quad \sup_{0< s<t \le T} s^{\theta} \frac{ \big\|v_t-v_s\big\|_{\cb_{\infty}^{\tau-2\gamma}} } { |t-s|^{\gamma} }
 \Big \}.$$
Then let
$$V_{\eps, T} :=  \mathcal{D}_T^{\frac{7}{12}+\frac{\eps}{2}, \frac12 + 2\eps,  \frac{1}{4} + \frac{\eps}{4}}   \quad \text{and} \quad 
 W_{\eps,T} := \mathcal{D}_T^{\frac34+\frac32\eps, 1+2\eps, \frac12+\frac{\eps}{2}},$$ 
so that we can write
$$  \big\|(v,w)\big\|_{\eps, T}  =\max\Big(\big\|v\big\|_{V_{\eps, T} } ,\big\|w\big\|_{W_{\eps, T} }\Big). $$

\begin{lemma} \label{lem:basic}
The following estimates hold.
\begin{enumerate}[$(i)$]
\item For $\beta \ge \alpha$ and $0 \le \tau \le 1$, we have 
$$\big\|e^{-tH} f \big\|_{\mathcal{D}_T^{\frac{\beta-\alpha}{2}, \beta, \tau}} \lesssim \big\|f\big\|_{\mathcal{B}_{\infty}^{\alpha}}.$$   
\item Let $0 < \eta <1$, $\alpha<\beta $, $\alpha \le \gamma <\alpha-2\eta +2$, $\gamma \le \beta <\alpha+2$. 
For any $0<\delta < \frac{\beta-\alpha}{2}$, 
$$
\Big\| \displaystyle{\int_0^t d\tau \, e^{-(t-\tau)H} u_\tau } \Big\|_{\mathcal{D}_T^{\frac{\beta-\gamma}{2}, \beta, \delta}} 
\lesssim T^{\frac{\alpha-2\eta +2 -\gamma}{2}} \Big( \sup_{0<t \le T} t^{\eta} \big\|u_t\big\|_{\mathcal{B}_{\infty}^{\alpha}} \Big).$$ 
In particular, for all $-\frac32 <\al\leq -\frac23$, $0\leq \eta<\frac{\al}{2}+\frac43$ and $\eps>0$ small enough,
$$
\Big\| \displaystyle{\int_0^t d\tau \, e^{-(t-\tau)H} u_\tau } \Big\|_{V_{\eps,T}} 
\lesssim T^{\frac{\alpha-2\eta}{2}+\frac43-\frac{\eps}{2}} \Big( \sup_{0<t \le T} t^{\eta} \big\|u_t\big\|_{\mathcal{B}_{\infty}^{\alpha}} \Big)
$$ 
while for all $-1 <\al\leq -\frac12$, $0\leq \eta<\frac{\al}{2}+\frac54$ and $\eps>0$ small enough,
$$
\Big\| \displaystyle{\int_0^t d\tau \, e^{-(t-\tau)H} u_\tau } \Big\|_{W_{\eps,T}} 
\lesssim T^{\frac{\alpha-2\eta }{2}+\frac54+\frac{\eps}{2}} \Big( \sup_{0<t \le T} t^{\eta} \big\|u_t\big\|_{\mathcal{B}_{\infty}^{\alpha}} \Big).
$$ 
\item If $(v,w) \in \mathcal{X}_{\eps,T}$, then $v \in  V_{\eps, T} $ and 
$w \in W_{\eps,T}$.  Moreover,  
for $-\frac{2}{3}+\eps \le \gamma \le \frac12 +2\eps$, we have  
 $$\big\|v_t\big\|_{\cb_{\infty}^{\gamma}} \lesssim t^{-\frac12 (\gamma-\eps+\frac23)}  \big\| v \big\|_{V_{\eps,T}}, $$ 
 and $-\frac{1}{2}-\eps \le \gamma \le 1+2\eps$, we have  
 $$\big\|w_t\big\|_{\cb_{\infty}^{\gamma}} \lesssim t^{-\frac12 (\gamma+\eps+\frac12)}  \big\| w\big\|_{W_{\eps,T}}.$$ 
 \end{enumerate}
\end{lemma}

 In particular, for  $-\frac{1}{2}-\eps \le \gamma \le \frac12+2\eps$, we have  
 \begin{equation}\label{eqgamma}
  \big\|v_t\big\|_{\cb_{\infty}^{\gamma}} + \big\|w_t\big\|_{\cb_{\infty}^{\gamma}} \lesssim t^{-\frac12 (\gamma-\eps+\frac23)}    \big\|(v,w)\big\|_{\eps,T}.  
   \end{equation}
Using the bound $\big\|u\big\|_{L^{\infty}} \lesssim \big\|u\big\|_{\cb^{\frac{\eps}2}_{\infty}} $ and \eqref{eqgamma} with $\gamma =\frac{\eps}2$, we deduce
 \begin{equation}\label{eqinfi}
  \big\|v_t\big\|_{L^{\infty}} + \big\|w_t\big\|_{L^{\infty}} \lesssim t^{-\frac13+\frac{\eps}4}    \big\|(v,w)\big\|_{\eps,T}.  
   \end{equation}

\begin{lemma} \label{FG} Let $0<T\le 1$. 
Fix $\eps>0$ small enough. For every $(v,w) \in \mathcal{X}_{\eps,T}$  and every $0<t\le T$, there exist $C'({\eps}),C''(\eps)>0$ such that 
\begin{equation} \label{eq:F}
\big\| F(v_t+w_t; Z_{0,t})\big\|_{\cb_{\infty}^{-1-\eps}} \le C' \Big( \big\|Z\big\|_{\eps}+ t^{-\frac13} \big\|(v,w)\big\|_{{{\eps},T}}\Big) \big\|Z\big\|_{\eps},
\end{equation}
\begin{equation} \label{eq:G}
\big\| G(v,w; Z_{0,.}{,\overline{\frakc}_{0,.}})_t \big\|_{\cb_{\infty}^{-\frac12-2\eps}} \le C'' t^{-1+\frac{3\eps}4} \Big(1+\big\|(v,w)\big\|^3_{\eps,T}\Big)\Big(1+{\big\llbracket  \overline{\frakc}\big\rrbracket_{\eps}^3+}\big\|Z\big\|_{\eps}^3\Big).
\end{equation}
\end{lemma}

\begin{proof}
We first prove \eqref{eq:F}. By Proposition \ref{Prop-est-para} $(ii)$, for $0  \le s \le T$ 
\begin{eqnarray}  
\big\| F(v_s+w_s; Z_{0,s} ) \big\|_{\mathcal{B}_{\infty}^{-1-\eps}} 
&\lesssim &
\big\| (v_s+w_s-\<IPsi3>_{0,s}) \pl  \<Psi2>_{0,s} \big\|_{\mathcal{B}_{\infty}^{-1-\eps}}\nonumber \\
&  \lesssim &  \big\|v_s+w_s-\<IPsi3>_{0,s}\big\|_{L^{\infty}} \big\|  \<Psi2>_{0,s}\big\|_{\mathcal{B}_{\infty}^{-1-\eps}}\nonumber  \\
&\lesssim &\Big(\big\|v_s\big\|_{L^{\infty}}  +\big\|w_s\big\|_{L^{\infty}} 
+\big\|\<IPsi3>_{0,s}\big\|_{\cb_{\infty}^{\frac12-\eps}} \Big)\big\|  \<Psi2>_{0,s}\big\|_{\mathcal{B}_{\infty}^{-1-\eps}} \nonumber \\
&\lesssim& \Big(s^{-\frac13} \big\|(v,w)\big\|_{\eps,T} +\big\|Z\big\|_{\eps}\Big) \big\|Z\big\|_{\eps}\label{gav}
\end{eqnarray}
where in the last line we have used \eqref{eqinfi}. 

\medskip

 We now turn to \eqref{eq:G} and will bound each term in the decomposition \eqref{definition-g} separately.  

\smallskip

First, one has, by \eqref{eqinfi}, for $0\le s\le T,$
\begin{multline*}
\big\| \big(v_s+w_s\big) ^3 \big\| _{\cb_{\infty}^{-\frac12 -2\eps}} \lesssim \big\|\big(v_s+w_s\big)^3\big\|_{L^{\infty}(\R^3)} \lesssim \\ \lesssim \big\|v_s+w_s\big\|_{L^{\infty}(\R^3)}^3
\lesssim \Big(s^{-\frac13+\frac{\eps}4} \big\|(v,w)\big\|_{\eps,T}\Big)^3\lesssim  s^{-1+\frac{3\eps}4} \big\|(v,w)\big\|^3_{\eps,T}. 
\end{multline*}
Also by Proposition \ref{Prop-est-para} $(i)$, 
$$
\big\| w_s \pe \<Psi2>_{0,s}\big\| _{\cb_{\infty}^{-\frac12 -2\eps}} \lesssim \big\| w_s \pe \<Psi2>_{0,s}\big\| _{\cb_{\infty}^{\eps}} \lesssim\big\| w_s\big\|_{\cb_{\infty}^{1 +2\eps}} 
\big\| \<Psi2>_{0,s} \big\|_{\cb_{\infty}^{-1 -\eps}}\lesssim  s^{-\frac12 \big(\frac32 +3\eps \big)}\big\|w\big\|_{W_{\eps,T}} \big\|Z\big\|_{\eps},
$$
and in a similar way, by Proposition \ref{Prop-est-para} $(iii)$
\begin{align*}
\big\| \big(v_s+w_s-\<IPsi3>_{0,s}\big) \pg \<Psi2>_{0,s} \big\|_{\cb_{\infty}^{-\frac12 -2\eps}} &
\lesssim \big\| v_s+w_s-\<IPsi3>_{0,s}\big\|_{\cb_{\infty}^{\frac12-\eps}} \big\| \<Psi2>_{0,s} \big\|_{\cb_{\infty}^{-1 -\eps}}\\
&\lesssim s^{-\frac{7}{12}}  \big\|(v,w)\big\|_{\eps,T}\big\|Z\big\|_{\eps} +\big\|Z\big\|_{\eps}^2 .
\end{align*}

For the control of $\big[\pl, \pe\big](\<IPsi3>, \<IPsi3>, \<Psi>)$, which appears in $\tau^{(0)}$, we can appeal to \cite[Proposition~13.25]{DFT}, which gives 
$$\big\| \big[\pl, \pe\big](\<IPsi3>_{0,s}, \<IPsi3>_{0,s}, \<Psi>_{0,s})\big\|_{\cb_{\infty}^{-\frac12 -2\eps}}
\lesssim \big\| \<IPsi3>_{0,s}\big\|_{\cb_{\infty}^{\frac12 -\eps}}^2 \big\| \<Psi>_{0,s}\big\|_{\cb_{\infty}^{-\frac12 -\eps}}\lesssim \big\|Z\big\|_{\eps}^3.$$

Now, observe that as an easy consequence of the properties in Proposition \ref{Prop-est-para}, we obtain
$$ \big\|\tau^{(0)}\big\|_{\cb_{\infty}^{-\frac12 -2\eps}} \lesssim 1+\big\|Z\big\|_{\eps}^3, \quad \quad \big\|\tau^{(1)}\big\|_{\cb_{\infty}^{-\frac12 -\eps}} \lesssim  1+\big\|Z\big\|_{\eps}^2 \quad \quad \text{and}
\quad \quad  \big\|\tau^{(2)}\big\|_{\cb_{\infty}^{-\frac12 -\eps}} \lesssim \big\|Z\big\|_{\eps}.$$
We deduce in particular that
\begin{multline*}
\big\| \tau^{(1)}_s \cdot \big(v_s+w_s\big) \big\|_{\cb_{\infty}^{-\frac12 -2\eps}} \lesssim \big\| \tau^{(1)}_s \cdot \big(v_s+w_s\big) \big\|_{\cb_{\infty}^{-\frac12 -\eps}} 
\lesssim \big\|\tau^{(1)}_s\big\|_{\cb_{\infty}^{-\frac12 -\eps}} \big\|v_s+w_s\big\|_{\mathcal{B}_{\infty}^{\frac12+2\eps}} \\
\lesssim  \big( 1+\big\|Z\big\|_{\eps}^2\big)  s^{-\frac7{12} -\frac{\eps}2} \big\|(v,w)\big\|_{\eps,T}
\end{multline*}
and thanks to \cite[Corollary 13.21]{DFT} and \eqref{eqgamma}
\begin{align*}
\big\| \tau^{(2)}_s\cdot \big(v_s+w_s\big)^2 \big\|_{\cb_{\infty}^{-\frac12 -2\eps}} &\lesssim \big\|\tau^{(2)}_s\big\|_{\cb_{\infty}^{-\frac12-\eps}}  
\big\|\big(v_s+w_s\big)^2\big\|_{{\mathcal B}_{\infty}^{\frac12 +2\eps}}\\
&\lesssim \big\|\tau^{(2)}_s\big\|_{\cb_{\infty}^{-\frac12-\eps}}   \big\|v_s+w_s\big\|_{L^\infty} \big\|v_s+w_s\big\|_{{\mathcal B}_{\infty}^{\frac12 +2\eps}} \\
&\lesssim s^{-\frac{11}{12}-\frac{\eps}4} \big\|Z\big\|_{\eps} \big\|(v,w)\big\|_{\eps,T}^2 .
\end{align*}

\smallskip

As far as the term ${\overline{\frakc}_{0,r}}\big(\<Psi>_{0,r}-\<IPsi3>_{0,r}+v_r+w_r\big)$ is concerned, we can rely on Lemma \ref{lem:c} to assert that for every $\eps>0$ small enough,
\begin{align*}
\big\|{\overline{\frakc}_{0,r}}\big(\<Psi>_{0,r}-\<IPsi3>_{0,r}+v_r+w_r\big) \big\|_{\cb_{\infty}^{-\frac12 -2\eps}}
&\lesssim \big\|{\overline{\frakc}_{0,r}}\big\|_{\cb_{\infty}^{\frac12 +3\eps}}\big\|\<Psi>_{0,r}-\<IPsi3>_{0,r}+v_r+w_r \big\|_{\cb_{\infty}^{-\frac12 -2\eps}}\\
&\lesssim \frac{1}{r^{\frac78}} {\big\llbracket  \overline{\frakc}\big\rrbracket_{\eps}} \Big(\big\|Z\big\|_{\eps} +\big\|v_r \big\|_{\cb_{\infty}^{-\frac12 -\eps}}+\big\|w_r \big\|_{\cb_{\infty}^{-\frac12 -\eps}}\Big)\\
&\lesssim \frac{1}{r^{\frac{23}{24}}}  {\big\llbracket  \overline{\frakc}\big\rrbracket_{\eps}}\big(\big\|Z\big\|_{\eps} +\big\|(v,w)\big\|_{\eps,T}\big),
\end{align*}
due to \eqref{eqgamma}.

\smallskip

To achieve \eqref{eq:G}, it only remains to control the term $\big\| \mathrm{com}^{Z_{0,.}} (v+w)_{0,s}\big\|_{\cb_{\infty}^{-\frac12-2\eps}}$ defined in \eqref{comm1}.  To do so, we write
\begin{align*}
\big\| \mathrm{com}^{Z_{0,.}}(v,w)_{0,s}\big\|_{\cb_{\infty}^{-\frac12-2\eps}}& \lesssim  \big\|\mathrm{com}^{Z_{0,.}}_1(v,w)_{0,s} \pe \<Psi2>_{0,s} \big\| _{\cb_{\infty}^{-\frac12-2\eps}} +  \big\|\mathrm{com}^{Z_{0,.}}_2 (v+w)_{0,s}\big\| _{\cb_{\infty}^{-\frac12-2\eps}}\\
&\lesssim \big\|\mathrm{com}^{Z_{0,.}}_1(v,w)_{0,s} \big\| _{{\mathcal{B}}_{\infty} ^{1+2\eps}}  \big\| \<Psi2>_{0,s}\big\|_{\cb_{\infty}^{-1 -\eps}} + \big\|  \mathrm{com}_2^{Z_{0,.}} (v+w)_{0,s} \big\|_{\cb_{\infty}^\eps}.
\end{align*}
Then we combine the technical results of Corollary \ref{coro22} below, Proposition \ref{Prop-est-para} and Lemma~\ref{lem:basic}. 
This yields for every $\eps>0$ small enough,
\begin{align*}
\big\| \mathrm{com}^{Z_{0,.}}(v,w)_{0,s}\big\|_{\cb_{\infty}^{-\frac12-2\eps}}& \lesssim  \Big[  { s^{-\frac12}}   \big( 1+ \big\|Z\big\|_{\eps}^2 \big) \big( 1+ \big\|(v,w)\big\|_{ \eps,T}\big)\Big] 
\big\| \<Psi2>_{0,s}\big\|_{\cb_{\infty}^{-1 -\eps}}\\
&\hspace{1cm}+ \Big[\big\| v_s \big\|_{\cb_{\infty}^{\frac12 -\eps}} + \big\|w_s\big\|_{\cb_{\infty}^{\frac12-\eps}} + \big\|\<IPsi3>_{0,s}\big\|_{\cb_{\infty}^{\frac12-\eps}} \Big]\big\| \<IPsi2>_{0,s}\big\|_{\cb_{\infty}^{1-2\eps}} 
\big\| \<Psi2>_{0,s}\big\|_{\cb_{\infty}^{-1-\eps}}\\
&\lesssim \Big[    { s^{-\frac12}} \big( 1+ \big\|Z\big\|_{\eps}^2 \big) 
\big( 1+ \big\|(v,w)\big\|_{ \eps,T}\big)\Big]\big\|Z\big\|_{\eps}\\
&\hspace{2cm}+\Big[  s^{-\frac{7}{12} +\eps} \big\| (v,w) \big\|_{\eps,T} +\big\|Z\big\|_{\eps} \Big] \big\|Z\big\|_{\eps}^2.
\end{align*}
This completes the proof of \eqref{eq:G}.
\end{proof}

\begin{proof}[Proof of Proposition \ref{prop:fixed-point}]
{\it Bound on $\gga_{Z}^{\mathbf{v}}[v,w]$.} By Lemma \ref{lem:basic} $(ii)$ and Lemma \ref{FG} we obtain 
\begin{align*}
\big\|\gga_{Z}^{\mathbf{v}} [v,w] \big\| _{V_{\eps,T}} 
&=  \Big\| \int_0^t  ds\, e^{-(t-s)(H+\mu)}  F(v_s +w_s; Z_{0,s} )\Big\|_{V_{\eps,T}} \\ 
&  \lesssim  T^{\frac12 -\eps} \sup_{0 < t \le T} t^{\frac13} \big\| F(v_t +w_t; Z_{0,t} )\big\|_{\mathcal{B}_{\infty}^{-1-\eps}} \\
&\lesssim  T^{\frac12-\eps} \Big(\big\|(v,w)\big\|_{ \eps,T} +T^{\frac13} \big\|Z\big\|_{\eps}\Big) \big\|Z\big\|_{\eps}\\
& \lesssim  P_2\big(\big\|Z\big\|_{\eps}\big) T^{\frac12-\eps} \Big(1+\big\|(v,w)\big\|_{\eps,T}\Big) .
\end{align*}
On the other hand, similarly to \eqref{gav}, we have, for $Z_{0,.},Z'_{0,.} \in \mathcal Z_{\frac{\eps}{2},[0,1]}$, $(v,w), (v',w') \in \mathcal{X}_{\eps,T},$ 
\begin{align*}
\big\|F_Z(v,w) -  F_{Z'} (v',w') \big\|_{\cb_{\infty}^{-1-\eps}}
\le  P_1 (1+t^{-\frac13})
\big\{\big\|Z-Z'\big\|_{\eps} +\big\|(v,w) -(v',w')\big\|_{\eps,T} \big\}
\end{align*}
where $P_1$ is a positive, first-order polynomial expression of the variables 
$\big\|(v,w)\big\|_{\eps,T}$, $\big\|(v',w')\big\|_{\eps,T}$, $\big\|Z\big\|_{\frac{\eps}{2},1}$, and $\big\|Z'\big\|_{\frac{\eps}{2},1} $.  
Thus,  
\begin{equation*}
\big\|\gga^{\mathbf{v}}_{Z} [v,w]-\gga^{\mathbf{v}}_{Z'} [v,w] \big\|_{V_T^{\eps}} \lesssim 
  \bar{P}_1 T^{\frac12 -\eps} \Big(\big\|Z-Z'\big\|_{\eps} +\big\|(v,w) -(v',w')\big\|_{\eps,T}\Big)
\end{equation*}
with some $\tilde{P}_1$ and $\bar{P_1}$ positive, first-order polynomials of the variables 
$\big\|(v,w)\big\|_{\eps,T}$, $\big\|(v',w')\big\|_{\eps,T}$, $\big\|Z\big\|_{\eps}$, $\big\|Z'\big\|_{\eps} $.  

\smallskip

\textit{Bound on $\gga^{\mathbf{w}}[v,w]$.}  
Let us first focus on the estimates related to $\<Psi2IPsi3>$. For any $t \le T$, one has
\begin{eqnarray*}
&& \Big\|\int_0^t e^{-(t-s)H} \<Psi2IPsi3>_{0,s}\,  ds \Big\|_{\cb_{\infty}^{-\frac{1}{2}-\eps}} \lesssim \int_0^t |t-s|^{-\frac16}  \big\| \<Psi2IPsi3>_{0,s}\big\|_{\cb_{\infty}^{-\frac56-\eps}} ds 
\le T^{\frac56}  \big\| \<Psi2IPsi3>\big\|_{C_T \cb_{\infty}^{-\frac56-\eps}} \le T^{\frac56} \big\|Z\big\|_{\eps}, 
\end{eqnarray*}
while
\begin{eqnarray*}
&&  \Big\| \int_0^t e^{-(t-s)H} \<Psi2IPsi3>_{0,s}\,  ds \Big\|_{\cb_{\infty}^{1+2\eps}}
\lesssim  \int_0^t |t-s|^{-(\frac{11}{12} +\frac{3}{2} \eps)}  \big\| \<Psi2IPsi3>_{0,s}\big\|_{\cb_{\infty}^{-\frac56-\eps}} ds 
\lesssim t^{\frac1{12}-\frac32\eps} \big\|Z\big\|_{\eps} , 
\end{eqnarray*}
hence
\begin{align*}
t^{\frac{3}{4}+\frac{3}{2}\eps}  \Big\| \int_0^t e^{-(t-s)H} \<Psi2IPsi3>_{0,s}\,  ds \Big\|_{\cb_{\infty}^{1+2\eps}}\lesssim  T^{\frac34} \big\|Z\big\|_{\eps} . 
\end{align*}
Also, for $0\le s<t \le T,$
\begin{multline*}
\Big\| \int_0^t e^{-(t-r)H} \<Psi2IPsi3>_{0,r}\, dr -\int_0^s e^{-(s-r)H} \<Psi2IPsi3>_{0,r}\,  dr \Big\|_{\cb_{\infty}^{\eps}} \lesssim\\
\begin{aligned}
& \lesssim \Big\| \int_s^t e^{-(t-r)H} \<Psi2IPsi3>_{0,r}\,  dr \Big\|_{\cb_{\infty}^{\eps}} 
+\Big\| (\mathrm{Id}- e^{-(t-s)H}) \int_0^s e^{-(s-r)H} \<Psi2IPsi3>_{0,r}\,  dr \Big\|_{\cb_{\infty}^{\eps}} \\
&\lesssim \int_s^t |t-r|^{-\frac{5}{12}-\eps} \big\| \<Psi2IPsi3>_{0,r} \big\|_{\cb_{\infty}^{-\frac56-\eps}} dr 
+ |t-s|^{\frac12+\frac{\eps}{2}} \Big\| \int_0^s e^{-(s-r)H} \<Psi2IPsi3>_{0,r}\, dr \Big\|_{\cb_{\infty}^{1+2\eps}} \\
&\lesssim |t-s|^{\frac{7}{12} -\eps} \big\|Z\big\|_{\eps} +  |t-s|^{\frac12+\frac{\eps}{2}} s^{\frac{1}{12}-\frac{3}{2} \eps}  \big\|Z\big\|_{\eps}. 
\end{aligned}
\end{multline*}
Thus, 
\begin{eqnarray*}
 \sup_{0\le s<t \le T} s^{\frac34 +\frac32 \eps} |t-s|^{-\frac12 -\frac{\eps}{2}} 
\Big\| \int_0^t e^{-(t-r)H} \<Psi2IPsi3>_{0,r}\, dr -\int_0^s e^{-(s-r)H} \<Psi2IPsi3>_{0,r}\, dr \Big\|_{\cb_{\infty}^{\eps}}  \lesssim
T^{\frac34} \big\|Z\big\|_{\eps}.
\end{eqnarray*}

As for the estimates related to ${G}(v,w; Z_{0,.}{,\overline{\frakc}_{0,.}})_s$, we have, using Lemma \ref{lem:basic} $(ii)$ and Lemma~\ref{FG}, 
\begin{align*}
 \Big\| \int_0^t  ds\, e^{-(t-s)H}  {G}(v,w; Z_{0,.}{,\overline{\frakc}_{0,.}})_s   \Big\|_{W_{\eps,T}} & \lesssim T^{\frac{\eps}4} \Big( \sup_{0<t\le T} t^{1-\frac{3\eps}4} \big\| {G}(v,w; Z_{0,.}{,\overline{\frakc}_{0,.}})_t\big\|_{\cb_{\infty}^{-\frac12-2\eps}}\Big)  \\
&\lesssim T^{\frac{\eps}4}   \Big(1+\big\|(v,w)\big\|_{\eps,T}^3\Big)  \Big(1+{\big\llbracket  \overline{\frakc}\big\rrbracket_{\eps}^3+}\big\|Z\big\|_{\eps}^3\Big).
\end{align*}
 Finally, using Lemma \ref{lem:basic} $(ii)$, one has
$$\Big\| \int_0^t  ds\, e^{-(t-s)H} v_s  \Big\|_{W_{\eps,T}}\lesssim T^{\frac{11}{12}} \Big( \sup_{0<t\le T} \big\|v_t\big\|_{\cb_{\infty}^{-\frac23+\eps}}\Big)\lesssim T^{\frac{11}{12}}\big\|(v,w)\big\|_{\eps,T}.$$
Combining the above bounds, we obtain that 
\begin{multline*}
\big\|\gga_{Z,  {\overline{\frakc},}w_0}^{\mathbf{w}} [v,w]\big\|_{W_{\eps,T}}  \lesssim \big\| e^{-tH} w_0\big\|_{W_{\eps,T}} \\
\begin{aligned}
&\hspace{0.5cm}  +  \Big\|\int_0^t e^{(t-s)H} \<Psi2IPsi3>_{0,s} ds \Big\|_{W_{\eps,T}}+\Big\| \int_0^t  ds\, e^{-(t-s)H}  {G}(v,w; Z_{0,.}, {\overline{\frakc}_{0,.}})_s   \Big\|_{W_{\eps,T}} +\mu\, \Big\| \int_0^t  ds\, e^{-(t-s)H} v_s  \Big\|_{W_{\eps,T}}  \\
&\lesssim
\big\|w_0\big\|_{\cb_{\infty}^{-\frac12-\eps}} + T^{\frac34} \big\|Z\big\|_{\eps}+ T^{{\frac{\eps}{4}}} \big(1+\|(v,w)\|_{\eps,T}^3\big)  {\big(1+\big\llbracket  \overline{\frakc}\big\rrbracket_{\eps}^3+\big\|Z\big\|_{\eps}^3\big)}
 + T^{\frac{11}{12}}\|(v,w)\|_{\eps,T}
 \\
&\lesssim
\big\|Z\big\|_{\eps}+\|w_0\|_{\cb_{\infty}^{-\frac12-\eps}}  +T^{{\frac{\eps}{4}}}\big(1+\big\|(v,w)\big\|_{\eps,T}^3\big)  {\Big(1+\big\llbracket  \overline{\frakc}\big\rrbracket_{\eps}^3+\big\|Z\big\|_{\eps}^3\Big)}.
\end{aligned}
\end{multline*}

\smallskip

With similar arguments, we can check that 
\begin{multline*}
\big\|\gga^{\mathbf{w}}_{Z,{\overline{\frakc},} w_0} [v,w]-\gga^{\mathbf{w}}_{Z',  {\overline{\frakc}',} w'_0} [v,w] \big\|_{W_T^{\eps}}  \lesssim  \\
 \lesssim \tilde{P}_2 \Big(  \big\|w_0-w_0'\big\|_{\cb_{\infty}^{-\frac12-\eps}}  \Big) +T^{\nu} P_2 \Big( \big\|(v,w) -(v',w')\big\|_{\eps,T} {+\big\llbracket \overline{\frakc}- \overline{\frakc}'\big\rrbracket_{\eps}} +\big\|Z-Z'\big\|_{\eps} \Big),
\end{multline*}
for $\nu >0$. This completes the proof of Proposition \ref{prop:fixed-point}.
\end{proof}

\subsection{First commutator estimate}

  Here we state the counterpart to \cite[Proposition 2.2]{MW}. Recall the notation 
  $$\delta_{st}f=f_t-f_s.$$
  Then, we have

  \begin{proposition}\label{prop22}
 Let $\eps>0$ and $\beta \in(4\eps, 1+2\eps]$.   Then for all $p \geq 1$ large enough, $r_1,r_2 \geq p$, $0\leq s<T\leq s+1$, $(v,w) \in \mathcal{X}_{\eps,[s, T]}$ and $t \in [s,T)$, one has  
  \begin{multline*}
 \big\|\mathrm{com}^{Z_{s,.}}_1(v,w)_{s,t}  \big\|_{\B^{1+2\eps}_p} \lesssim \big\|Z_{s,.}\big\|_{\mathcal Z_{\frac{\eps}{2},[s,T]}}^2+\big\|Z_{s,.}\big\|_{\mathcal Z_{\frac{\eps}{2},[s,T]}} \int_s^t  \frac{d\tau}{(t-\tau)^{1+2\eps-\frac{\beta}2}} \Big(\big\|v_\tau\big\|_{\B^\beta_p}+\big\|w_\tau\big\|_{\B^\beta_p}\Big)+\\
 +\big\|Z_{s,.}\big\|_{\mathcal Z_{\frac{\eps}{2},[s,T]}}\int_s^t  \frac{d\tau}{(t-\tau)^{1+2\eps}} \Big(\big\|\delta_{\tau t}v\big\|_{L^{r_1}}+\big\|\delta_{\tau t}w\big\|_{L^{r_2}}\Big), 
   \end{multline*}
   where the implicit proportional constant only depends on $\eps$ and $p$,  and where the quantity $\big\|Z_{s,.}\big\|_{\mathcal Z_{\frac{\eps}{2},[s,T]}}$ has been introduced in \eqref{K}. 
   \end{proposition}
	
	We can deduce the following result which we have already used in the proof of Lemma \ref{FG}.

	\begin{corollary}\label{coro22}
  Let $\eps>0$ small enough and $0\leq s<T\leq s+1$. We have for  every $(v,w) \in \mathcal{X}_{\eps, [s,T]}$ and $t \in [s,T)$,
  \begin{equation*}
 \big\|\mathrm{com}^{Z_{s,.}}_1(v,w)_{s,t}  \big\|_{\B^{1+2\eps}_\infty} \lesssim \big\|Z_{s,.}\big\|_{\mathcal Z_{\frac{\eps}{2},[s,T]}}^2+(t-s)^{-\frac12}\big\|Z_{s,.}\big\|_{\mathcal Z_{\frac{\eps}{2},[s,T]}}\big\|(v,w)\big\|_{\eps,[s,T]}.
   \end{equation*}
   \end{corollary}

 \begin{proof}[Proof of Corollary \ref{coro22}] 
 We can apply the result of Proposition \ref{prop22} with $\beta=\frac12+2\eps$ and $p=r_1=r_2=\infty$ to obtain
   \begin{multline*}
 \big\|\mathrm{com}^{Z_{s,.}}_1(v,w)_{s,t}  \big\|_{\B^{1+2\eps}_\infty} \lesssim \big\|Z_{s,.}\big\|_{\mathcal Z_{\frac{\eps}{2},[s,T]}}^2+\big\|Z_{s,.}\big\|_{\mathcal Z_{\frac{\eps}{2},[s,T]}} \int_s^t  \frac{d\tau}{(t-\tau)^{\frac34+\eps}} \Big(\big\|v_\tau\big\|_{\B^{\frac12+2\eps}_\infty}+\big\|w_\tau\big\|_{\B^{\frac12+2\eps}_\infty}\Big)+\\
 +\big\|Z_{s,.}\big\|_{\mathcal Z_{\frac{\eps}{2},[s,T]}}\int_s^t  \frac{d\tau}{(t-\tau)^{1+2\eps}} \Big(\big\|\delta_{\tau t}v\big\|_{L^{\infty}}+\big\|\delta_{\tau t}w\big\|_{L^{\infty}}\Big).
   \end{multline*}
   Then we use \eqref{eqgamma} with $\gamma=\frac12+ 2\eps$, and by the definition of $\big\|(v,w)\big\|_{\eps,[s,T]}$, we deduce 
      \begin{multline*}
 \big\|\mathrm{com}^{Z_{s,.}}_1(v,w)_{s,t}  \big\|_{\B^{1+2\eps}_\infty} \lesssim \big\|Z_{s,.}\big\|_{\mathcal Z_{\frac{\eps}{2},[s,T]}}^2+\big\|Z_{s,.}\big\|_{\mathcal Z_{\frac{\eps}{2},[s,T]}}\big\|(v,w)\big\|_{\eps,[s,T]} \int_s^t  d\tau\, (t-\tau)^{-\frac34-\eps} \tau^{-\frac7{12}-\frac{\eps}2}+ \\
  +\big\|Z_{s,.}\big\|_{\mathcal Z_{\frac{\eps}{2},[s,T]}}\big\|(v,w)\big\|_{\eps,[s,T]}\int_s^t d\tau \, (t-\tau)^{-1-2\eps} \Big(\tau^{-\frac{7}{12}-\frac{\eps}2} (t-\tau)^{\frac14+\frac{\eps}4}  +  \tau^{-\frac{3}{4}-\frac{3\eps}2} (t-\tau)^{\frac12+\frac{\eps}2}  \Big) \\
   \lesssim \big\|Z_{s,.}\big\|_{\mathcal Z_{\frac{\eps}{2},[s,T]}}^2+(t-s)^{-\frac13-\frac{9\eps}4}\big\|Z_{s,.}\big\|_{\mathcal Z_{\frac{\eps}{2},[s,T]}}\big\|(v,w)\big\|_{\eps,[s,T]},
   \end{multline*}
   which yields the claim.
  \end{proof}

\medskip

 \begin{proof}[Proof of Proposition \ref{prop22}] The argument is very similar to the proof of \cite[Proposition  2.2]{MW}. In our case, we need to add a degree of freedom 
with the norms $\big\|\delta_{\tau t}v\big\|_{L^{r_1}}$ and $\big\|\delta_{\tau t}w\big\|_{L^{r_2}}$  on the right-hand side in the estimates, where $r_1, r_2\geq p$. As in \cite{MW} we introduce the commutator operator
 $$[e^{-t H }, \pl] : (f,g) \longmapsto e^{-t H}(f \pl g) -f \pl (e^{-t H}g),$$
 and we write the decomposition 
 \begin{multline}\label{eq-22}
 e^{-(t-\tau) H}  \big[(v_\tau+w_\tau- \<IPsi3>_{s,\tau}) \pl \<Psi2>_{s,\tau} \big]=\\
  =(v_\tau+w_\tau- \<IPsi3>_{s,\tau}) \pl \big[ e^{-(t-\tau) H} \<Psi2>_{s,\tau} \big]+ \big[e^{-(t-\tau)H}, \pl\big] \big(v_\tau+w_\tau  - \<IPsi3>_{s,\tau},   \<Psi2>_{s,\tau} \big).
  \end{multline}
The contribution of the second term of \eqref{eq-22} is estimated as in \cite{MW}. Actually, using \eqref{heat1} and the commutation result \cite[Lemma~13.26]{DFT}, we can prove that
  \begin{multline*}
  \Big\| \int_s^t d\tau \,  \big[e^{-(t-\tau)H}, \pl\big] \big(v_\tau+w_\tau  - \<IPsi3>_{s,\tau},   \<Psi2>_{s,\tau} \big)  \Big\|_{\B^{1+2\eps}_p} \lesssim \\
  \lesssim \big\|Z_{s,.}\big\|_{\mathcal Z_{\frac{\eps}{2},[s,T]}}^2 + \big\|Z_{s,.}\big\|_{\mathcal Z_{\frac{\eps}{2},[s,T]}}\int_s^t  \frac{d\tau}{(t-\tau)^{1+2\eps-\frac{\beta}2}} \Big(\big\|v_\tau\big\|_{\B^\beta_p}+\big\|w_\tau\big\|_{\B^\beta_p}\Big).
  \end{multline*}
  In order to treat the contribution of the first term of \eqref{eq-22}, we use the identity
   \begin{multline*}
 \Big[(v_t+w_t- \<IPsi3>_{s,t}) \pl \<IPsi2>_{s,t} \Big] -  \int_s^t d\tau \,(v_\tau+w_\tau- \<IPsi3>_{s,\tau}) \pl \big[ e^{-(t-\tau) H} \<Psi2>_{s,\tau} \big]  =\\
 =\int_s^t d\tau \, \big[ \delta_{\tau t} (v+w- \<IPsi3>_{s,.}  )   \big] \pl \big[ e^{-(t-\tau) H} \<Psi2>_{s,\tau} \big] .
    \end{multline*}
  Then 
  \begin{multline*}
  \Big\| \int_s^t d\tau \, \big[ \delta_{\tau t} (v+w- \<IPsi3>_{s,.}  )   \big] \pl \big[ e^{-(t-\tau) H} \<Psi2>_{s,\tau} \big]  \Big\|_{\B^{1+2\eps}_p} \lesssim \\
    \begin{aligned}
& \lesssim   \int_s^t  d\tau\,  \Big\| \big[ \delta_{\tau t} (v+w- \<IPsi3>_{s,.}  )   \big] \pl \big[  e^{-(t-\tau) H}    \<Psi2>_{s,\tau}  \big]   \Big\|_{\B^{1+2\eps}_p} \\
& \lesssim   \int_s^t d\tau\,  \Big\| \big[ \delta_{\tau t} v   \big] \pl \big[  e^{-(t-\tau) H}    \<Psi2>_{s,\tau}  \big]   \Big\|_{\B^{1+2\eps}_p} + \int_s^t  d\tau\, \Big\| \big[ \delta_{\tau t} w   \big] \pl \big[  e^{-(t-\tau) H}    \<Psi2>_{s,\tau}  \big]   \Big\|_{\B^{1+2\eps}_p}  \\
&\qquad \qquad + \int_s^t d\tau \,  \Big\| \big[ \delta_{\tau t}    \<IPsi3>_{s,.} \big] \pl \big[  e^{-(t-\tau) H}    \<Psi2>_{s,\tau}  \big]   \Big\|_{\B^{1+2\eps}_p}  .
   \end{aligned}
   \end{multline*}
  
Let $r_1\geq p$. By Proposition \ref{Prop-est-para} $(ii)$  and Lemma  \ref{lem:actisemi} we have, for $1/r_1+ 1/\ov{r}_1=1/p$
\begin{eqnarray}\label{eq-23}
   \int_s^t  d\tau\, \Big\| \big[ \delta_{\tau t}   v\big] \pl \big[  e^{-(t-\tau) H}    \<Psi2>_{s,\tau}  \big]   \Big\|_{\B^{1+2\eps}_p}   &\lesssim & \int_s^t d\tau\,  \big\|   \delta_{\tau t}   v  \big\|_{L^{r_1}} \Big\|   e^{-(t-\tau) H}    \<Psi2>_{s,\tau}    \Big\|_{\B^{1+2\eps}_{\ov{r}_1}}  \nonumber \\
   &\lesssim &  \int_s^t   \frac{d\tau}{(t-\tau)^{1+\frac{3\eps}2}}    \big\|   \delta_{\tau t}   v  \big\|_{L^{r_1}}    \big\|      \<Psi2>_{s,\tau}    \big\|_{\B^{-1-\eps}_{\ov{r}_1}}  .
 \end{eqnarray}
  Now we use Lemma \ref{lem-bis}: if $p=p(\eps)$ is large enough, there exists $3/\ov{r}_1< \eta<\eps/2$ such that  
  $$ \big\|      \<Psi2>_{s,\tau}    \big\|_{\B^{-1-\eps}_{\ov{r}_1}}     \lesssim \big\|      \<Psi2>_{s,\tau}    \big\|_{\B^{-1-\frac{\eps}2}_{\infty} } \leq \big\|Z_{s,.}\big\|_{\mathcal Z_{\frac{\eps}{2},[s,T]}},$$
  then using \eqref{eq-23} we deduce 
  $$    \int_s^t d\tau\,  \Big\| \big[ \delta_{\tau t}   v\big] \pl \big[  e^{-(t-\tau) H}    \<Psi2>_{s,\tau}  \big]   \Big\|_{\B^{1+2\eps}_p} \lesssim     \big\|Z_{s,.}\big\|_{\mathcal Z_{\frac{\eps}{2},[s,T]}} \int_s^t   \frac{d\tau}{(t-\tau)^{1+\frac{3\eps}2}}    \big\|   \delta_{\tau t}   v  \big\|_{L^{r_1}}     .     $$
  Similarly, for any $r_2 \geq p$
    $$    \int_s^t  d\tau\, \Big\| \big[ \delta_{\tau t}   w\big] \pl \big[  e^{-(t-\tau) H}    \<Psi2>_{s,\tau}  \big]   \Big\|_{\B^{1+2\eps}_p}  \lesssim     \big\|Z_{s,.}\big\|_{\mathcal Z_{\frac{\eps}{2},[s,T]}}\int_s^t    \frac{d\tau}{(t-\tau)^{1+\frac{3\eps}2}}    \big\|   \delta_{\tau t}   w  \big\|_{L^{r_2}}.     $$

  By Proposition \ref{Prop-est-para}  and Lemma  \ref{lem:actisemi} we have
\begin{eqnarray*}
   \int_s^t d\tau\,   \Big\| \big[ \delta_{\tau t}    \<IPsi3>_{s,.} \big] \pl \big[  e^{-(t-\tau) H}    \<Psi2>_{s,\tau}  \big]   \Big\|_{\B^{1+2\eps}_p}   &\lesssim & \int_s^t  d\tau\, \big\|   \delta_{\tau t}    \<IPsi3>_{s,.}  \big\|_{L^p} \Big\|   e^{-(t-\tau) H}    \<Psi2>_{s,\tau}     \Big\|_{\B^{1+2\eps}_\infty}  \\
   &\lesssim &\big\|Z_{s,.}\big\|_{\mathcal Z_{\frac{\eps}{2},[s,T]}}  \int_s^t     \frac{d\tau}{(t-\tau)^{1+\frac{3\eps}2}}\big\|   \delta_{\tau t }    \<IPsi3>_{s,.}   \big\|_{L^p}   \\
      &\lesssim &\big\|Z_{s,.}\big\|_{\mathcal Z_{\frac{\eps}{2},[s,T]}}^2,
  \end{eqnarray*}
  where in the last estimate we have used that
    $$   \big\|   \delta_{\tau t}    \<IPsi3>_{s,.}   \big\|_{L^p}  \lesssim  \big\| \delta_{\tau t}  \<IPsi3>_{s,.}  \big\|_{\B^{{\frac{\eps}2}}_{\infty}}  \lesssim \big\|Z_{s,.}\big\|_{\mathcal Z_{\frac{\eps}{2},[s,T]}}(t-\tau)^{\frac18}.$$
  We gather all the previous estimates and complete the proof.
   \end{proof}


\section{A priori bounds and globalization argument}\label{Sect4}

Having established the local well-posedness of the model, we now turn to the derivation of global a priori bounds. The key challenge is to obtain estimates on the solution $(v,w)$ that are  uniform with respect to initial data, which will  allow us to extend the local solution to a global one. As before, our strategy follows the approach of Mourrat and Weber~\cite{MW}, with the modifications required by the whole-space setting.

\

\subsection{Main result}

Our main \enquote{coming down from infinity} result can now be stated as follows.

\begin{theorem}  \label{thm:global_Lp} 
Fix $\eps,\nu>0$ small enough, $0\leq s <T\leq s+1$, $u\in \cb_\infty^{-\frac12-\eps}$, $Z_{s,.}\in \mathcal Z_{\frac{\eps}{2},[s,T]}$, $\overline{\frakc}_{s,.}:=3 \overline{\frakc}^{\mathbf{1}}_{s,.}-9 \overline{\frakc}^{\mathbf{2}}_{s,.}$ with $\big\llbracket  \overline{\frakc}^{\mathbf{1}}_{s,.}\big\rrbracket_{3\eps,[s,T]}, \big\llbracket  \overline{\frakc}^{\mathbf{2}}_{s,.}\big\rrbracket_{3\eps,[s,T]}<\infty$, and set
\begin{equation}\label{notation-mst}
M_{s,T}:=1+\big\|Z_{s,.}\big\|_{\mathcal Z_{\frac{\eps}{2},[s,T]}}+\big\llbracket  \overline{\frakc}^{\mathbf{1}}_{s,.}\big\rrbracket_{3\eps,[s,T]}+ \big\llbracket  \overline{\frakc}^{\mathbf{2}}_{s,.}\big\rrbracket_{3\eps,[s,T]}.
\end{equation}
Then {there exist $\ka\geq 1$ (depending only on $\eps,\nu$) and $\mu \geq 1$ (depending on $\eps$, $\nu$ and $M_{s,T}$)} such that if $(v^u_{s,.},w^u_{s,.})\in \mathcal{X}_{\eps,[s,T]} $ is a solution of the system \eqref{mild:v} starting at $u$ {(with parameter $\mu$)}, one has
\begin{equation}\label{boundun}
\max\Big(\sup_{t\in [s,T]} \big\|v^u_{s,t}\big\|_{\cb_{\infty}^{-\eps}},  \sup_{t\in [s,T]} (t-s)^{\frac38+\nu}\big\|v^u_{s,t}\big\|_{\cb_{\infty}^{\frac34-\eps}} \Big)\lesssim M^{{\ka}}_{s,T},
\end{equation}
as well as
\begin{equation}\label{bounddeux}
\max\Big(\sup_{t\in [s,T]}(t-s)^{\frac12+\nu}\big\|w^u_{s,t}\big\|_{\cb_{\infty}^{-\eps}}, \sup_{t\in [s,T]} (t-s)^{\frac32+\nu}\big\|w^u_{s,t}\big\|_{\cb_{\infty}^{\frac32-\eps}} \Big)\lesssim M^\ka_{s,T},
\end{equation}
where the proportional constants in \eqref{boundun} and \eqref{bounddeux} only depend on $\eps$ and $\nu$ (in particular, they do not depend on $u$).
\end{theorem}

\

The above estimates immediately provide us with the following fundamental control on the solution $X$ of \eqref{eq:intro} (recall that the notion of rough solution has been introduced in Definition \ref{def:rough-sol}).

\begin{corollary}\label{coro:bounx}
In the setting of Theorem \ref{thm:global_Lp}, there exists $\ka\geq 1$ such that any (rough) solution~$X^u$ of \eqref{eq:intro} on $[s,T]$ associated with $(Z,\overline{\frakc})$ and starting at $u\in \cb_\infty^{-\frac12-\eps}$ satisfies 
\begin{equation}\label{bounx}
\sup_{s < t \le T} (t-s)^{\frac12+\nu} \big\|X^u_t\big\|_{\B_{\infty}^{-\frac12-\eps}} \lesssim M^\ka_{s,T},
\end{equation}
where the proportional constant does not depend on $u$.
\end{corollary}

\

\begin{corollary}\label{bound_final} 
Fix $\eps>0$ small enough, $T>0$ and $u\in \cb_{\infty}^{-\frac12-\eps}(\R^3)$. Let $(Z,Z^{(n)})$ and $(\overline{\frakc},\overline{\frakc}^{(n)})$ be defined as in Propositions \ref{prop:conv-arbre} and \ref{prop:c-convergence}. 

\smallskip

Then the following assertions hold.

\smallskip

\noindent
$(i)$ Almost surely, the equation \eqref{eq:intro} associated with $(Z,\overline{\frakc})$ and $u\in \cb_\infty^{-\frac12-\eps}$ admits a unique (rough) solution~$X^u$ on $[0,T]$. Moreover, for every $0\leq s\leq T$, one has
\begin{equation}\label{dec-x}
X^u_. - \<Psi>_{s,.} + \<IPsi3>_{s,.} \in \mathcal{C}\big((s,(s+1)\wedge T]; \mathcal{B}_{\infty}^{\frac34-\eps}(\R^3)\big). 
\end{equation}

\smallskip

\noindent
$(ii)$ The sequence $(X^{(n),u})$ of approximated solutions on $[0,T]$ starting from $u$ converges to $X^u$ in the space $\cac([0,T];\cb_\infty^{-\frac12-\eps})$.

\smallskip

\noindent
$(iii)$ For every $p\geq 1$, it holds that
\begin{equation}\label{bounx-mom}
\sup_{s \ge 0}\, \mathbb{E} \Big[ \sup_{s < t \le s+1}  \sup_{u \in \B_{\infty}^{-\frac{1}{2}-\eps}} \Big( (t^{\frac12+\nu} \wedge 1) \big\|X^{u}_t\big\|_{\B_{\infty}^{-\frac12-\eps}} \Big)^p\Big] <\infty  
\end{equation}
as well as
\begin{equation}\label{bounx-mom-n}
\sup_{n\geq 1} \sup_{s \ge 0}\, \mathbb{E} \Big[ \sup_{s < t \le s+1}  \sup_{u \in \B_{\infty}^{-\frac{1}{2}-\eps}} \Big( (t^{\frac12+\nu} \wedge 1) \big\|X^{(n),u}_t\big\|_{\B_{\infty}^{-\frac12-\eps}} \Big)^p\Big] <\infty.  
\end{equation}
\end{corollary}

\

\begin{proof}[Proof of Corollary \ref{bound_final}] Recall the definition \eqref{def-normes}. 

\smallskip

\noindent
$(i)$ By iterating the result of Corollary \ref{coro:main}, we can construct an increasing sequence of times $(T_i)_{i\geq -1}\in [0,1]$ such that:

\smallskip 

\noindent
$\bullet$ $T_{-1}=T_0=0$.

\smallskip 

\noindent
$\bullet$ for all $i\geq 0$, there exists a unique solution $X^u$ of \eqref{eq:intro} on $[0,T_i]$ starting from $u$ at time $0$.

\smallskip

\noindent
$\bullet$ one has for all $i\geq 0$, $T_{i+1}:=\min \big(1,T_i+\delta_{i+1}\big)$, where
\begin{equation*}
\delta_{i+1}:=C\, \big( 1+\big\|X^u_{T_i}\big\|_{  \mathcal{B}^{-\frac12-\eps}_{\infty}}\big)^{-\zeta}\big( 1+\big\| \! \big\| Z\big\| \! \big\|_{\mathcal Z_{\eps,[0,2]}} +\big\llbracket\!\! \big\llbracket  \overline{\frakc}\big\rrbracket \!\! \big\rrbracket_{\eps,[0,2]} \big)^{-\zeta},
\end{equation*}
for universal constants $\zeta\geq 1$ and $C\in (0,1)$.

\smallskip

Thanks to \eqref{bounx}, we deduce that for all $i\geq 1$,
$$\delta_{i+1} \geq \frac{C' T_{i}^{\zeta'}}{\big(1+\big\| \! \big\| Z\big\| \! \big\|_{\mathcal Z_{\eps,[0,2]}} +\big\llbracket\!\! \big\llbracket  \overline{\frakc}\big\rrbracket \!\! \big\rrbracket_{\eps,[0,2]} \big)^{\zeta''}}\geq \frac{C' T_{1}^{\zeta'}}{\big(1+\big\| \! \big\| Z\big\| \! \big\|_{\mathcal Z_{\eps,[0,2]}} +\big\llbracket\!\! \big\llbracket  \overline{\frakc}\big\rrbracket \!\! \big\rrbracket_{\eps,[0,2]} \big)^{\zeta''}},$$
for constants $C',\zeta',\zeta''>0$. Note that this lower bound still depends on $u$ (through $T_1$), but it no longer depends on the successive values $X^u_{T_1}, X^u_{T_2}, X^u_{T_3},\ldots$.

\smallskip

As a result, there must exist $i_0\geq 0$ such that $T_{i_0}=1$, and $X^u$ is thus the unique solution of~\eqref{eq:intro} on the whole interval $[0,1]$.

\smallskip

We can repeat the argument on each interval $[\ell,\ell+1]$, $\ell \in \N$, and derive the existence of a unique solution on any interval $[0,T]$. Decomposition \eqref{dec-x} is then a straightforward consequence of Definition \ref{def:rough-sol}.

\smallskip

The assertion in item $(ii)$ is a direct consequence of Corollary \ref{coro:main} (item $(ii)$), while the bounds \eqref{bounx-mom} and~\eqref{bounx-mom-n} follow from the combination of Corollary \ref{coro:bounx}, Proposition \ref{prop:conv-arbre} and Lemma \ref{lem:c}.
\end{proof}

\color{black}

\

\

The first part of the proof of Theorem \ref{thm:global_Lp} follows the same general strategy as in \cite{MW}; indeed, most of the estimates therein carry over {\it mutatis mutandis}. However, a fundamental difference arises: at several points in \cite{MW}, the embedding $L^q(\T^3) \subset L^p(\T^3)$ (for $p \leq q$) is crucially used, but this embedding of course fails in the setting of $\R^3$. The obstruction is overcome by means of Lemma~\ref{lem-bis}, which provides a comparison between $\B^{\alpha}_{p,\infty}(\R^3)$ and $\B^{\beta}_{q,\infty}(\R^3)$ spaces in the regime $p \leq q$, $\alpha < \beta$.

\smallskip

Furthermore, in \cite{MW}, the parameters $\eps > 0$ (a small real number) and $p \geq 1$ (a large even integer) are chosen independently of one another. In the present setting, however, the two parameters are coupled: we first fix $\eps > 0$ sufficiently small ({\it e.g.} $\eps = 10^{-5}$), and then select $p = p(\eps)$ subject to the constraint $\dis p \geq \frac{100}{\eps}$.

\medskip

For clarity of exposition, we present the proof of Theorem \ref{thm:global_Lp} only in the case $s = 0$, as the argument for general $s \geq 0$ is entirely analogous. Accordingly, we fix $T\in (0,1]$ and write 
$$v := v^u_{0,\cdot} \; \text{ and } \;  w := w^u_{0,\cdot}$$
where $(v^u_{0,.},w^u_{0,.})\in \mathcal{X}_{\eps,[0,T]} $ is assumed to be a solution of the system \eqref{syst1} on $[0,T]$ starting at $u$. We also define
$$M:=M_{0,T}.$$

 
 \subsection{A priori estimate on $v$}

 \subsubsection{The bound on $v$}
  
 Here we obtain some  bounds on $v$ which slightly extend \cite[Theorem~3.1]{MW}, where only the case $r=p$ was needed.  Recall that the notation $\mu \geq 1$ stands for the parameter appearing in \eqref{syst1}.

\begin{theorem}\label{thm21}
Let \( \varepsilon > 0 \) be sufficiently small.  Let \( p, p', q, r \in [1,\infty] \) with \( p' \leq q \), \( p \leq q \), $p \leq r$, and let \( \theta \in [-3\varepsilon, 1-4\varepsilon) \) such that
\begin{align*}
\frac{\theta + 3\varepsilon}{2} + \frac{3}{2}\left( \frac{1}{p'} - \frac{1}{q} \right) &< 1, \\
\sigma := \frac{\theta+1+\varepsilon}{2} + \frac{3}{2}\left( \frac{1}{p} - \frac{1}{q} \right) &< 1,
\end{align*}
and define
\begin{equation}\label{underlinemu}
\underline{\mu} := \mu - 1 - \big( M \Gamma(1-\sigma) \big)^{\frac{1}{1-\sigma}},
\end{equation}
where $\Gamma$ denotes the Euler Gamma function. Then, for all \(0\leq  t < T \), we have:
\begin{equation}\label{est-v-thm21}
\big\| v_t \big\|_{\B^\theta_q} \lesssim   M \int_0^t d\tau \, e^{-\underline{\mu}(t-\tau)} (t-\tau)^{-\sigma} \Big( \big\| w_\tau \big\|_{L^r} + M \Big) . 
\end{equation}

Moreover, if \( \theta > \varepsilon \), and if $p \geq 1$ is large enough, then for
\begin{equation*}
\sigma' := \frac{1}{2} + \varepsilon
\quad \text{and} \quad
\underline{\mu}' := \mu - 1 - \left( M \Gamma(1-\sigma') \right)^{\frac{1}{1-\sigma'}},
\end{equation*}
we have, for all \( 0 \leq s \leq t < T \),
\begin{equation}\label{est-v-thm21-bis}
\big\| \delta_{st} v \big\|_{L^p} \lesssim (\mu+M) |t-s|^{\frac{\theta-\varepsilon}{2}} \big\| v_s \big\|_{\B^\theta_p} + M \int_s^t d\tau \, e^{-\underline{\mu}'(t-\tau)} (t-\tau)^{-\sigma'} \Big( \big\| w_\tau \big\|_{L^r} + M \Big) ,
\end{equation}
where \( \delta_{st}v = v_t - v_s \).

In all estimates, the implicit constants depend on \( \varepsilon \), \( p \), \( q \), and \( \theta \), but neither on \( \mu \geq 1 \) nor on~$M$.
\end{theorem}

The argument is very similar to the proof of \cite[Theorem 3.1]{MW}. In our case, we add a degree of freedom with the norm $\|w\|_{L^r}$ on the right-hand side in the estimates, where $r \geq p$.

\begin{proof}
Recall that 
\begin{equation*} 
F(v+w;Z_{0,.}) = -3(v+w -\<IPsi3>_{0,.}) \pl  \<Psi2>_{0,.}.
\end{equation*}
As in \cite{MW} we decompose $F(v+w;Z_{0,.})= -3(w -\<IPsi3>_{0,.}) \pl  \<Psi2>_{0,.} -3v\pl  \<Psi2>_{0,.}$. Let $r \geq p$, then we have the estimates 
\begin{equation}\label{bff}
\big\| (w -\<IPsi3>_{0,.}) \pl  \<Psi2>_{0,.} \big\|_{\B_p^{-1-2\eps}} \lesssim \big\| w -\<IPsi3>_{0,.} \big\|_{L^r} \big\|\<Psi2>_{0,.} \big\|_{\B_{\ov{r}}^{-1-2\eps}} ,
\end{equation}
where $1/r+1/\ov{r}=1/p$. Assume that $p=p(\eps) \geq 1$ is large enough, then by Lemma \ref{lem-bis}, if we choose $3/p<\eta<\eps$ we have 
\begin{equation}\label{bff2}
\big\|\<Psi2>_{0,.}\big\|_{\B_{\ov{r}}^{ -1-2\eps}} \lesssim  \big\|\<Psi2>_{0,.}\big\|_{\B_{\infty}^{ -1-2\eps+\eta}} \lesssim   \big\|\<Psi2>_{0,.}\big\|_{\B_{\infty}^{ -1-\eps}}  \lesssim   M.
 \end{equation}
Therefore, from \eqref{bff} and \eqref{bff2}
\begin{equation*}
\big\| (w -\<IPsi3>_{0,.}) \pl  \<Psi2>_{0,.} \big\|_{\B^{-1-\eps}_p}  \lesssim M \Big( \big\| w \big\|_{L^r}+M\Big).
\end{equation*}
By Proposition \ref{Prop-est-para} $(iii)$ we also have 
\begin{equation*}
\big\| v\pl  \<Psi2>_{0,.} \big\|_{\B^{-1-4\eps}_q} \lesssim   \big\|  v \big\|_{\B^{-3\eps}_{q}}   \big\|\<Psi2>_{0,.}\big\|_{\B^{{-1-\eps}}_{\infty}} \lesssim M \big\|  v \big\|_{\B^{-3\eps}_{q}}  .
\end{equation*}
The rest of the argument is similar to that of \cite{MW}.
\end{proof}

\subsubsection{On the parameters} The strategy is the following. First, we fix the parameter $\eps>0$ small enough, 
such that the local well-posedness result holds. Then $p= p(\eps) \geq 1$ has to be taken large enough. For convenience,~$p$ can be chosen as an even integer.  In the sequel, we fix
\begin{equation*}
\sigma <\frac78, \quad \quad \underline{\mu} := \mu-1-\big(8M \big)^8, \quad \text{where} \ \mu\geq \big(8M \big)^8.
\end{equation*}

\

In the procedure that follows, the three parameters $(\eps,p,\mu)$ will be chosen in the following order:

\smallskip

\noindent
$(1)$ we first choose $\eps>0$ small enough, depending only on universal constants,

\smallskip

\noindent
$(2)$ we then choose $p\geq 1$ large enough, possibly depending on $\eps$,

\smallskip

\noindent
$(3)$ we finally choose $\mu\geq 1$ large enough, depending on $\eps$, $p$ \textit{and also on $M$} (see condition~\eqref{condiC-bis}). 

\

For now, we only label the first condition
\begin{equation}\label{condiC}
\mu\geq 1+\big(8 M \big)^8.
\end{equation}
Note that under this condition, the positivity of $\underline{\mu}$ in \eqref{underlinemu} (and in \eqref{est-v-thm21}) is guaranteed for every $\sigma <\frac78$, which will be used several times in the proof.

\

 \subsection{A priori estimate on $\delta_{st}w$}\label{subsec:aprideltaw}
 
We now derive an estimate on the time increments $\delta_{st}w$ of the second component. This estimate, combined with the previous bound on $v$, will be iterated in the subsequent sections to close the a priori argument. 
  
\begin{theorem} \label{thm41}
{Let $\eps>0$ be small enough. Then for every $p \geq 1$ large enough, every $\mu$ satisfying~\eqref{condiC}} and for all $0\leq s \leq t  \leq T$, we have 
\begin{equation}\label{rhsthm41}
\big\|\delta_{s t} w\big\|_{L^p}     \lesssim {(\mu M)}^7 (t-s)^{\tfrac{1}{8}}
\Bigg[
1 + \big\| w_s \big\|_{\B^{1+4\eps}_p}
+ \Big( \int_0^t d\tau\, \big\|w_\tau\big\|^p_{\B^{1+4\varepsilon}_p}   \Big)^{\tfrac{1}{p}}+ \left( \int_0^t d\tau\,  \big\|w_\tau\big\|^{3p}_{L^{3p}} \right)^{\tfrac{1}{p}}
\Bigg] 
\end{equation}
where the implicit constant depends on $\eps$ and $p$, but neither on $M$ nor on $\mu$ satisfying \eqref{condiC}.
\end{theorem}

Recall that 
$$w_t = e^{-(t-s)H} w_s +  {\int_s^t d\tau\, e^{-(t-\tau)H} \Big[3\ \<Psi2IPsi3>_{0,\tau}+G(v,w; Z_{0,.}, {\overline{\frakc}_{0,.}})_\tau+\mu v_\tau\Big]},$$
where $G$ is defined in \eqref{definition-g}. We introduce 
\begin{equation}\label{defD1}
\widehat{\delta}_{st}w:=w_t-e^{-(t-s)H}w_s- 3\int_s^t d\tau\, e^{-(t-\tau)H} \  {\<Psi2IPsi3>_{0,\tau}},
\end{equation}
so that
\begin{equation}\label{defD2}
 \widehat{\delta}_{st}w= \int_{s}^t d\tau \, e^{-(t-\tau)H} \Big[ {G(v,w; Z_{0,.},\overline{\frakc}_{0,.})_\tau}+\mu v_\tau\Big].
\end{equation}
Observe that we then have 
\begin{equation}\label{defD3}
\delta_{st}w=\widehat{\delta}_{st}w+\big(e^{-(t-s)H}w_s-w_s\big)+ 3\int_s^t d\tau\, e^{-(t-\tau)H} \  \<Psi2IPsi3>_{0,\tau},
\end{equation}
and accordingly, for $p\geq 1$ large enough
\begin{align}
\big\| \delta_{st}w \big\|_{L^p} &\lesssim \big\|\widehat{\delta}_{st}w\big\|_{L^p}+\big\|\big(e^{-(t-s)H}w_s-w_s\big)\big\|_{L^p}+ \big\|\int_s^t d\tau\, e^{-(t-\tau)H} \  \<Psi2IPsi3>_{0,\tau}\big\|_{L^p}\nonumber\\
&\lesssim \big\|\widehat{\delta}_{st}w\big\|_{L^p}+(t-s)^{\frac18}\Big[\big\| w_s \big\|_{\B^{1+4\eps}_p}+M\Big],\label{del-to-delhat}
\end{align}
which shows that we only need to focus on the estimate of $\big\|\widehat{\delta}_{st}w\big\|_{L^p}$. 

\smallskip

To this end, decompose $ \widehat{\delta}_{st}w$ as
\begin{equation}\label{decomp-delha-r}
 \widehat{\delta}_{st}w= \int_{s}^t d\tau \, e^{-(t-\tau)H} \Big[\mathrm{com}^{{Z_{0,\cdot}}}_1(v,w)_{0,\tau} \pe \<Psi2>_{0,\tau} +\overline{\frakc}_{0,\tau}\big(\<Psi>_{0,\tau}-\<IPsi3>_{0,\tau}+v_\tau+w_\tau\big)+R_\mu(v,w; Z_{0,.},\overline{\frakc}_{0,.})_\tau\Big],
\end{equation}
where the expression of $R_\mu(v,w; Z_{0,.},\overline{\frakc}_{0,.})_\tau$ can be derived from \eqref{defD2}, namely
\begin{align*}
R_\mu(v,w; Z_{0,.},\overline{\frakc}_{0,.})_\tau&:=\mu v_\tau- (v_\tau+w_\tau)^3 -3\,  \mathrm{com}^{Z_{0,.}}_2 (v+w)_{0,\tau}-3\, w_\tau \pe \<Psi2>_{0,\tau}  \\
&\hspace{2cm}-3(v_\tau+w_\tau-\<IPsi3>_{0,\tau}) \pg \<Psi2>_{0,\tau}+ P(v_\tau+w_\tau; Z_{0,\tau}).
\end{align*}

\smallskip

It turns out that all the terms composing $R_\mu(v,w; Z_{0,.},\overline{\frakc}_{0,.})_\tau$ can be treated by following line by line the proofs of \cite[Lemmas 4.2, 4.6, 4.7 and 4.8]{MW}. The result stemming from this observation can be summarized as follows.

\begin{lemma}\label{lem:id-MW}
 Let $\eps>0$ be small enough. Then for every $p \geq 1$ large enough, every $\mu$ satisfying~\eqref{condiC} and for all $0\leq s \leq t  \leq T$, we have
  \begin{multline*}
\Big\| \int_s^t d\tau \,  e^{-{(t-\tau)}H} \Big[  R_\mu(v,w; Z_{0,.},\overline{\frakc}_{0,.})_\tau\Big]\Big\|_{L^p}\lesssim  \\
 \lesssim  (\mu M)^7 (t-s)^{\tfrac{1}{8}}
\Bigg[
1 
+ \Big( \int_0^t d\tau\, \big\|w_\tau\big\|^p_{\B^{1+4\varepsilon}_p}   \Big)^{\tfrac{1}{p}}+ \left( \int_0^t d\tau\,  \big\|w_\tau\big\|^{3p}_{L^{3p}} \right)^{\tfrac{1}{p}}
\Bigg] ,
    \end{multline*}
where the implicit constant depends on $\eps$ and $p$.
\end{lemma}

Let us now turn to the estimate of the two remaining terms in \eqref{decomp-delha-r}, starting with the term arising from $\mathrm{com}^{{Z_{0,\cdot}}}_1(v,w)_{0,\tau} \pe \<Psi2>_{0,\tau} $. The latter indeed requires a treatment different from that of the proof of \cite[Lemma 4.5]{MW}, since one can no longer bound $\|w \|_{L^{p}}$ by $\|w \|_{L^{3p}}$ in the present non-compact setting.

\

\begin{lemma}\label{lem43}
{Let $\eps>0$ be small enough. Then for every $p \geq 1$ large enough, every $r\geq p$, every~$\mu$ satisfying \eqref{condiC}} and for all $0 \leq t \leq T$,
we have
  \begin{equation*}
\big\|\mathrm{com}^{{Z_{0,.}}}_1(v,w)_{0,t}
\big\|_{\B^{1+2\eps}_p}
\lesssim {(\mu M)}^3
  + {(\mu M)}^2\int_0^t  \frac{d {\tau}}{(t- {\tau})^{\frac34+\eps}}
      \big\|w_{ {\tau}}\big\|_{\B^{\frac12+2\eps}_p}
  + {\mu M} \int_0^t   \frac{d  {\tau}}{(t- {\tau})^{1+2\eps}}
      \big\|\delta_{ {\tau} t}w\big\|_{L^r},
   \end{equation*}
   where the implicit constant depends on $\eps$ and $p$.
\end{lemma}

\begin{proof}
We use Proposition \ref{prop22} with $r_1=p$ and $r_2=r$, and $s=0$, and we obtain
  \begin{multline*}
 \big\|\mathrm{com}^{{Z_{0,.}}}_1(v,w)_{{0,}t}
 \big\|_{\B^{1+2\eps}_p}
 \lesssim M^2
 + M\int_0^t \frac{d\tau}{(t-\tau)^{\frac34+\eps}}
     \Big(\big\|v_\tau\big\|_{\B^{\frac12+2\eps}_p}+\big\|w_\tau\big\|_{\B^{\frac12+2\eps}_p}\Big)\\
 + M\int_0^t \frac{d\tau}{(t-\tau)^{1+2\eps}}
     \Big(\big\|\delta_{\tau t}v\big\|_{L^{p}}+\big\|\delta_{\tau t}w\big\|_{L^{r}}\Big).
   \end{multline*}
The terms involving $v$ can be absorbed using the estimates of  Theorem~\ref{thm21} applied with
$\theta = \frac{1}{2}+2\eps$ and $q = p$ {(which, owing to Condition \eqref{condiC}, ensures that $\underline{\mu},\underline{\mu}'\geq 0$ in \eqref{est-v-thm21}-\eqref{est-v-thm21-bis})}. The argument is similar to  the proof of \cite[Lemma~4.3]{MW}, and we refer to it for the details.
\end{proof}


\begin{lemma}\label{lem46}
{Let $\eps>0$ be small enough. Then for every $p \geq 1$ large enough, every $\mu$ satisfying~\eqref{condiC}} and for all $0\leq s \leq t  \leq T$, we have
  \begin{multline}\label{AAAAA}
\Big\| \int_s^t d\tau \,  e^{-{(t-\tau)}H} \Big[  \mathrm{com}^{{Z_{0,\cdot}}}_1(v,w)_{0,\tau} \pe \<Psi2>_{0,\tau}\Big]\Big\|_{L^p}  \lesssim {(\mu M)^3} (t-s)^{\frac18} \bigg[1+     \Big(\int_0^t d\tau \, \big\|w_{\tau}\big\|^{p}_{\B^{\frac12+2\eps}_p}\Big)^{\frac1p}     \bigg] +\\
 \qquad +{\mu M}(t-s)^{1-\frac1{6p}} \tn w\tn^{\frac12}_{p,t}\Big(\int_0^t d\tau\, \big\|w_{\tau}\big\|^{3p}_{L^{3p}}\Big)^{\frac1{6p}} 
    \end{multline}
where $ \tn w\tn_{p,t}$ is defined by 
\begin{equation}\label{tnwdef}
\tn w\tn_{p,t}  := \sup_{0\leq \tau' < \tau \leq t} \frac{ \big\|\widehat{\delta}_{\tau'\tau}w\big\|_{L^p} }{|\tau-\tau'|^{\frac18}},
\end{equation}
with $\widehat{\delta}$ given by \eqref{defD1}. Moreover, the implicit constant depends on $\eps$ and $p$.
\end{lemma}

\begin{proof} 
We use Proposition \ref{Prop-est-para} to write 
  \begin{equation*}
\Big\| \int_s^t d\tau \,  e^{-{(t-\tau)}H} \Big[  \mathrm{com}^{{Z_{0,\cdot}}}_1(v,w)_{0,\tau} \pe
     \<Psi2>_{0,\tau}\Big]\Big\|_{L^p} \lesssim M  \int_s^t  d\tau \, \big\|\mathrm{com}^{{Z_{0,\cdot}}}_1(v,w)_{{0,}\tau} \big\|_{\B^{1+2\eps}_p}.
    \end{equation*}
Next, we apply Lemma \ref{lem43} and $r=\frac{3p}{2}$. Then 
\begin{equation*}
\Big\|  \mathrm{com}^{{Z_{0,\cdot}}}_1(v,w)_{{0,}\tau}
\Big\|_{\mathcal{B}^{1+2\eps}_p}\lesssim
{(\mu M)}^3+ {(\mu M)}^2\!\int_0^\tau \frac{d\tau'}{(\tau-\tau')^{\frac{3}{4}+\eps}}\big\|w_{\tau'}\big\|_{\mathcal{B}^{\frac{1}{2}+2\eps}_p}+ {\mu M}\!\int_0^\tau \frac{d\tau'}{(\tau-\tau')^{1+2\eps}}  \big\|\delta_{\tau'\tau}w\big\|_{L^{\frac{3p}{2}}}.
\end{equation*}
Compared to \cite[Lemma 4.5]{MW}, only the treatment of the term 
\begin{equation}\label{termdelta}
 \int_0^\tau \frac{d\tau'}{(\tau-\tau')^{1+2\eps}} \big\|\delta_{\tau'\tau}w\big\|_{L^{\frac{3p}2}}  
 \end{equation}
  is different, and we stress below what changes.

Let $0 \leq \tau' \leq \tau $. By \eqref{heat1} and \eqref{heat2},  we have the bound
\begin{eqnarray*}
\Big|  \big\| \delta_{\tau'\tau}w  \big\|_{L^{\frac{3p}2}}-\big\| \widehat{\delta}_{\tau'\tau}w \big\|_{L^{\frac{3p}2}} \Big| &\leq &  \big\| (\delta_{\tau'\tau}-\widehat{\delta}_{\tau'\tau} ) w   \big\|_{L^{\frac{3p}2}}    \\
&\lesssim & \big\|  (1- e^{-(\tau-\tau')H}) w_{\tau'}\big\|_{L^{\frac{3p}2}}{+\Big\|\int_{\tau'}^{\tau} d\beta\, e^{-(\tau-\beta)H} \, \<Psi2IPsi3>_{0,\beta}\Big\|_{L^{\frac{3p}2}}}\\
 &\lesssim  &(\tau-\tau')^{\frac18}\big\|w_{\tau'}\big\|_{\B^{\frac12+\eps}_{\frac{3p}2}}{+\int_{\tau'}^{\tau} \frac{d\beta}{(\tau-\beta)^{\frac12+\eps}} \ \big\|\<Psi2IPsi3>_{0,\beta}\big\|_{\B^{-1}_{\frac{3p}2}}}\\
 &\lesssim  &(\tau-\tau')^{\frac18}\big\|w_{\tau'}\big\|_{\B^{\frac12+\eps}_{\frac{3p}2}}{+\int_{\tau'}^{\tau} \frac{d\beta}{(\tau-\beta)^{\frac12+\eps}} \ \big\|\<Psi2IPsi3>_{0,\beta}\big\|_{\B^{-\frac56-\eps}_{\infty}}}\\
  &\lesssim  &(\tau-\tau')^{\frac18}\big\|w_{\tau'}\big\|_{\B^{\frac12+2\eps}_{p}}{+(\tau-\tau')^{\frac18} M},
\end{eqnarray*}
where we used  Lemmas \ref{lem-bis} and \ref{lem:inclusion-besov}, assuming that $p=p(\eps)\geq1$ large enough. As a consequence, $ \big\|\delta_{\tau'\tau}w\big\|_{L^{\frac{3p}2}}$ may be replaced with $ \big\|\widehat{\delta}_{\tau'\tau}w\big\|_{L^{\frac{3p}2}}$ in \eqref{termdelta}.

Next, by interpolation and  \eqref{heat-Lp}
\begin{eqnarray*}
\big\| \widehat{\delta}_{\tau'\tau}w  \big\|_{L^{\frac{3p}2}}  &\leq &    \big\| \widehat{\delta}_{\tau'\tau}w  \big\|^{\frac12}_{L^{p}}  \big\| \widehat{\delta}_{\tau'\tau}w \big\|^{\frac12}_{L^{3p}}    \\
&\lesssim &    \big\| \widehat{\delta}_{\tau'\tau}w  \big\|^{\frac12}_{L^{p}}  \Big(      \big\| w_{\tau} \big\|^{\frac12}_{L^{3p}}+  \big\|w_{\tau'}\big\|^{\frac12}_{L^{3p}} \big) .
\end{eqnarray*}
These bounds allow us to proceed as in  \cite[Lemma 4.5]{MW}.
 \end{proof}

\

Returning to \eqref{decomp-delha-r}, let us finally consider the term arising from $\overline{\frakc}_{0,\tau}\big(\<Psi>_{0,\tau}-\<IPsi3>_{0,\tau}+v_\tau+w_\tau\big)$, which does not appear in \cite{MW}.
\begin{lemma}\label{lem:nouveau}
Let $\eps>0$ be small enough. Then for every $p \geq 1$ large enough, every $\mu$ satisfying~\eqref{condiC} and for all $0\leq s \leq t  \leq T$, we have
  \begin{align}
&\Big\| \int_s^t d\tau \,  e^{-{(t-\tau)}H} \Big[ \overline{\frakc}_{0,\tau}\big(\<Psi>_{0,\tau}-\<IPsi3>_{0,\tau}+v_\tau+w_\tau\big)\Big]\Big\|_{L^p}\lesssim M^3 (t-s)^{\tfrac{1}{8}}
\Bigg[
1 + \Big( \int_0^t d\tau\, \big\|w_\tau\big\|^p_{\B^{1+4\varepsilon}_p}   \Big)^{\tfrac{1}{p}}
\Bigg],\label{AAAAA-andouillette}
    \end{align}
where the implicit constant depends on $\eps$ and $p$.
\end{lemma}

\begin{proof}

\smallskip

\noindent
\textit{(i)} By applying the paraproduct rules (see Proposition \ref{Prop-est-para}), and for $p\geq 1$ large enough, we have on the one hand
\begin{align*}
\Big\| \int_s^t d\tau \,  e^{-{(t-\tau)}H} \big[ \overline{\frakc}_{0,\tau}\pl \<Psi>_{0,\tau}\big]\Big\|_{L^p}&\lesssim  \int_s^t \frac{d\tau}{(t-\tau)^{\frac14+\eps}} \,   \big\| \overline{\frakc}_{0,\tau}\pl \<Psi>_{0,\tau}\big\|_{\cb_p^{-\frac12-2\eps}}\\
&\lesssim  \int_s^t \frac{d\tau}{(t-\tau)^{\frac14+\eps}} \,   \big\| \overline{\frakc}_{0,\tau}\big\|_{L^\infty}\big\|  \<Psi>_{0,\tau}\big\|_{\cb_p^{-\frac12-2\eps}}\\
&\lesssim M \int_s^t \frac{d\tau}{(t-\tau)^{\frac14+\eps}\tau^{\frac{9}{16}}} \,   \big\|  \<Psi>_{0,\tau}\big\|_{\cb_\infty^{-\frac12-\eps}} \lesssim M^2 (t-s)^{\frac18},
\end{align*}
and on the other hand
\begin{align*}
\Big\| \int_s^t d\tau \,  e^{-{(t-\tau)}H} \big[ \overline{\frakc}_{0,\tau}\pge\<Psi>_{0,\tau}\big]\Big\|_{L^p}&\lesssim \int_s^t d\tau \,  \big\| \overline{\frakc}_{0,\tau}\pge\<Psi>_{0,\tau}\big\|_{\cb^\eps_p}\\
& \lesssim \int_s^t d\tau \,  \big\| \overline{\frakc}_{0,\tau}\big\|_{\cb^{\frac12+3\eps}_\infty}  \big\|\<Psi>_{0,\tau}\big\|_{\cb^{-\frac12-2\eps}_p} \\
& \lesssim M\int_s^t \frac{d\tau}{\tau^{\frac{7}{8}}} \,   \big\|\<Psi>_{0,\tau}\big\|_{\cb^{-\frac12-\eps}_\infty} \lesssim M^2 (t-s)^{\frac18},
\end{align*}
so that in the end
\begin{align*}
&\Big\| \int_s^t d\tau \,  e^{-{(t-\tau)}H} \big[ \overline{\frakc}_{0,\tau}\<Psi>_{0,\tau}\big]\Big\|_{L^p}\lesssim  M^2 (t-s)^{\frac{1}{8}}.
\end{align*}

\smallskip

\noindent
\textit{(ii)} One has, for all $p\geq 1$ large enough,
\begin{multline*}
\Big\| \int_s^t d\tau \,  e^{-{(t-\tau)}H} \Big[ \overline{\frakc}_{0,\tau}\big(-\<IPsi3>_{0,\tau}+v_\tau+w_\tau\big)\Big]\Big\|_{L^p}\lesssim \\
\begin{aligned}
&\lesssim  \int_s^t d\tau \, \Big\| \overline{\frakc}_{0,\tau}\big(-\<IPsi3>_{0,\tau}+v_\tau+w_\tau\big)\Big\|_{L^p}\lesssim  \int_s^t d\tau \, \big\| \overline{\frakc}_{0,\tau}\big\|_{L^\infty} \Big\|-\<IPsi3>_{0,\tau}+v_\tau+w_\tau\Big\|_{L^p}\\
&\lesssim  M\Big[\int_s^t \frac{d\tau}{\tau^{\frac{9}{16}}}  \Big\|\<IPsi3>_{0,\tau}\Big\|_{L^p}+\int_s^t \frac{d\tau}{\tau^{\frac{9}{16}}}  \big\|w_\tau\big\|_{\cb_p^{1+4\eps}}+\int_s^t \frac{d\tau}{\tau^{\frac{9}{16}}}  \big\|v_\tau\big\|_{L^p}\Big]\\
&\lesssim  M\Big[M(t-s)^{\frac18}+\bigg(\int_s^t \frac{d\tau}{\tau^{\frac{9}{16}\frac{p}{p-1}}}\bigg)^{1-\frac1p}  \bigg(\int_s^t d\tau\, \big\|w_\tau\big\|_{\cb_p^{1+4\eps}}^p\bigg)^{\frac1p}+\int_s^t \frac{d\tau}{\tau^{\frac{9}{16}}}  \big\|v_\tau\big\|_{L^p}\Big]\\
&\lesssim  M\Big[M(t-s)^{\frac18}+(t-s)^{\frac18} \bigg(\int_s^t d\tau\, \big\|w_\tau\big\|_{\cb_p^{1+4\eps}}^p\bigg)^{\frac1p}+\int_s^t \frac{d\tau}{\tau^{\frac{9}{16}}}  \big\|v_\tau\big\|_{L^p}\Big].
\end{aligned}
\end{multline*}
Finally, for the last term into brackets, we can use Theorem \ref{thm21} (with $\theta=\eps$ and $q=r=p$) to assert that
\begin{align*}
&\int_s^t \frac{d\tau}{\tau^{\frac{9}{16}}}  \big\|v_\tau\big\|_{L^p}\lesssim \int_s^t \frac{d\tau}{\tau^{\frac{9}{16}}}  \big\|v_\tau\big\|_{\cb^\eps_p}\lesssim M \int_s^t \frac{d\tau}{\tau^{\frac{9}{16}}}  \int_0^\tau \frac{d\tau'}{(\tau-\tau')^{\frac12+\eps}}\Big(M+ \big\|w_{\tau'}\big\|_{L^{p}}\Big),
\end{align*}
where we have used Condition \eqref{condiC} to ensure that $\underline{\mu} \geq 0$ in \eqref{est-v-thm21}. This entails, for all $p\geq 1$ large enough,
\begin{align*}
\int_s^t \frac{d\tau}{\tau^{\frac{9}{16}}}  \big\|v_\tau\big\|_{L^p}&\lesssim M^2 \int_s^t \frac{d\tau}{\tau^{\frac{9}{16}}}  \int_0^\tau \frac{d\tau'}{(\tau-\tau')^{\frac12+\eps}}+ M \int_s^t \frac{d\tau}{\tau^{\frac{9}{16}}}  \int_0^\tau \frac{d\tau'}{(\tau-\tau')^{\frac12+\eps}}\big\|w_{\tau'}\big\|_{\cb_p^{1+4\eps}}\\
&\lesssim M^2 (t-s)^{\frac18}+ M \int_s^t \frac{d\tau}{\tau^{\frac{9}{16}}}  \bigg(\int_0^\tau \frac{d\tau'}{(\tau-\tau')^{(\frac12+\eps)\frac{p}{p-1}}}\bigg)^{1-\frac1p}\bigg(\int_0^t d\tau'\, \big\|w_{\tau'}\big\|_{\cb_p^{1+4\eps}}^p\bigg)^{\frac1p}\\
&\lesssim M^2 (t-s)^{\frac18}+ M (t-s)^{\frac18}  \bigg(\int_0^t d\tau\, \big\|w_{\tau}\big\|_{\cb_p^{1+4\eps}}^p\bigg)^{\frac1p}.
\end{align*}

\

Combining the estimates in the above items \textit{(i)} and \textit{(ii)} immediately yields the bound \eqref{AAAAA-andouillette}.
\end{proof}

We now have all the ingredients to conclude the proof of Theorem \ref{thm41} as in \cite{MW}. Namely, combining the bounds obtained in Lemmas \ref{lem:id-MW}, \ref{lem46} and \ref{lem:nouveau}, we obtain that for every $p\geq 1$ large enough, every $\mu $ satisfying \eqref{condiC} and every $t\in [0,T]$,
\begin{align*}
&\tn w\tn_{p,t} \lesssim {(\mu M)}^7 \Bigg[1 + \Big( \int_0^t d\tau\, \big\|w_\tau\big\|^p_{\B^{1+4\varepsilon}_p}   \Big)^{\tfrac{1}{p}}+ \left( \int_0^t d\tau\,  \big\|w_\tau\big\|^{3p}_{L^{3p}} \right)^{\tfrac{1}{p}}\Bigg]+ {\mu M}\Big(\int_0^t d\tau\, \big\|w_{\tau}\big\|^{3p}_{L^{3p}}\Big)^{\frac1{6p}} \tn w\tn^{\frac12}_{p,t}
\end{align*}
where $\tn w\tn_{p,t} $ is the quantity introduced in \eqref{tnwdef}.
Using the elementary implication $x \leq a + \sqrt{b\,x} \implies x \leq 2( a + b)$ for all $x,a,b\geq 0$, we deduce that for every $p\geq 1$ large enough, every $\mu $ satisfying \eqref{condiC} and every $t\in [0,T]$,
\begin{align*}
&\tn w\tn_{p,t} \lesssim {(\mu M)}^7 \Bigg[1 + \Big( \int_0^t d\tau\, \big\|w_\tau\big\|^p_{\B^{1+4\varepsilon}_p}   \Big)^{\tfrac{1}{p}}+ \left( \int_0^t d\tau\,  \big\|w_\tau\big\|^{3p}_{L^{3p}} \right)^{\tfrac{1}{p}}\Bigg].
\end{align*}

The claimed estimate for $\|\delta_{st} w\|_{L^p}$ finally stems from the observation in \eqref{del-to-delhat}.

\


\subsection{A priori estimate on $w$ in Besov spaces}

We now combine the estimates on $v$ and on the time increments of $w$ to derive a bound on $w$ in Besov spaces. The result below is the counterpart to \cite[Theorem 5.1]{MW}, and the argument is the same. 

\begin{theorem}\label{thm51}
Let $\eps>0$ be small enough and $0 < \gamma <\frac43$. Then for every $p \geq 1$ large enough, every $\mu$ satisfying \eqref{condiC}, and for all $0<s<t<T$, we have 
 \begin{multline*}
\big\|w_t\big\|_{\B^\gamma_p} \lesssim \big\| e^{-(t-s) H}w_s \big\|_{\B^\gamma_p}+ {(\mu M)}^9 \bigg[1+   \bigg(\int_s^t dr\, \big\|w_r\big\|^{3p}_{L^{3p}}\bigg)^{\frac1p}+  \bigg(\int_s^t dr\, \big\|w_r\big\|^{p}_{\B^{1+4\eps}_p}\bigg)^{\frac1p}+ \big\|v_s\big\|^3_{\B^{-3\eps}_{2p}}\bigg]+ \\
+{(\mu M)}^3 (t-s)^{\frac12(1-\gamma)-3\eps-\frac1{3p}} \bigg[  1 +\big\|v_s\big\|^2_{\B^{-3\eps}_{2p}}  +   \bigg(\int_s^t  dr\, \big\|w_r\big\|^{3p}_{L^{3p}}\bigg)^{\frac2{3p}}     \bigg],
 \end{multline*}
 where the implicit constant depends only on $\eps$, $p$ and $\gamma$, but neither on $M$ nor on $\mu$ satisfying~\eqref{condiC}.
\end{theorem}

The next result is the analogue of \cite[Corollary 5.2]{MW}.
\begin{corollary}\label{coro52}
Let $\eps>0$ be small enough. Then for every $p \geq 1$ large enough and every $\mu$ satisfying \eqref{condiC}, there exists $\kappa >0$ depending only on $p$ such that for all $0<s<t<T$, we have 
\begin{multline*}
\bigg( \int_s^t dr \, \big\|w_r\big\|^p_{\B^{1+7\varepsilon}_p}  \bigg)^{\tfrac{1}{p}} 
 \lesssim \bigg( \int_s^t dr\, \big\|e^{- (r-s)H}w_s\big\|^p_{\B^{1+7\varepsilon}_p}  \bigg)^{\tfrac{1}{p}} +(\mu M)^\kappa \bigg[ 1+  \bigg( \int_s^t dr\,  \big\|w_r\big\|^{3p}_{L^{3p}}  \bigg)^{\tfrac{1}{p}} +  \big\|v_s\big\|^3_{\B^{-3\eps}_{2p}}\bigg],
\end{multline*}
 where the implicit constant depends only on $\eps$ and $p$,  but neither on $M$ nor on $\mu$.
\end{corollary}

\begin{proof}[Proof of Corollary \ref{coro52}]
The proof follows the same lines as \cite[Corollary 5.2]{MW}. There is, however, one change to make. Namely, the estimate (5.2) there has to be replaced by the following one: 
\begin{equation}\label{foll}
\big\|w \big\|_{\B^{1+4\eps}_p}  \leq     \big\|w \big\|^{\frac{3-3\nu}{3-\nu}}_{\B^{1+7\eps}_{p}}+    \big\|w \big\|^3_{L^{3p}},
\end{equation}
with $\nu= \frac{\eps}{2+13\eps}>0$, so that $\frac{3-3\nu}{3-\nu}<1$. Let us now prove \eqref{foll}. By the interpolation inequality~\eqref{interp}, we have 
\begin{eqnarray}
\big\|w \big\|_{\B^{1+4\eps}_p}  \leq   \big\|w \big\|_{\B^{1+6\eps}_p} &\leq &\big\|w \big\|^{1-\nu}_{\B^{1+\frac{13}2\eps}_{q}} \big\|w \big\|^\nu_{\B^{0}_{3p}} \nonumber\\
&\lesssim &\big\|w \big\|^{1-\nu}_{\B^{1+\frac{13}2\eps}_{q}} \big\|w \big\|^\nu_{L^{3p}}, \label{in74}
\end{eqnarray}
where $\nu= \frac{\eps}{2+13\eps}$ and $q= \frac{1-\nu}{1-\frac{\nu}3}p<p$. Then by \eqref{lplq}, if $p\geq 1$ is large enough, there exists $0< \eta < \frac{\eps}2$ such that 
\begin{equation*}
\big\|w \big\|_{\B^{1+\frac{13}2\eps}_{q}} \lesssim \big\|w \big\|_{\B^{1+\frac{13}2\eps+\eta}_{p}}  \lesssim \big\|w \big\|_{\B^{1+7\eps}_{p}}, 
\end{equation*}
and then \eqref{foll} follows by \eqref{in74} and the Young inequality.
\end{proof}

\
 

\subsection{A priori estimate on $w$ in Lebesgue spaces}

In this section, we obtain an a priori bound on $w$ in Lebesgue spaces. This is the step where we crucially exploit the defocusing sign of the cubic nonlinearity, which provides a coercive structure through integration by parts. The main result  is the following, and is the analogue of \cite[Theorem 6.1]{MW}. {Recall again that the notation $\mu \geq 1$ stands for the parameter appearing in \eqref{syst1}.}

\begin{theorem}\label{thm61}
 Let $\eps >0$ be small enough. Then for all $p\geq 1$ large enough, there exist $\mu_0, \kappa>0$ depending on $\eps$ and $p$ such that if 
 \begin{equation}\label{condiC-bis}
\mu \geq \max\big(1+\big(8 M\big)^8,\mu_0 M^{30p}\big),
 \end{equation}
we have for all $0\leq s \leq t \leq T$,
 \begin{equation*}
 \big\|w_t\big\|^{3p-2}_{L^{3p-2}}+\int_s^t dr \, \big\|w_r\big\|^{3p}_{L^{3p}} \lesssim \big\|w_s\big\|^{3p-2}_{L^{3p-2}}+(\mu M)^\kappa \bigg[1+\big\|v_s\big\|^{3p}_{\B^{-3\eps}_{2p}}+\int_s^t dr\, \big\|w_r \big\|^p_{\B^{1+6\eps}_p}\bigg],
 \end{equation*}
 where the implicit constant depends only on $\eps$ and $p$.
\end{theorem}

We proceed as follows: we test the second equation of the system  \eqref{syst1} against $w^{3p-3}$, where $p\geq 1$ is a large even integer. Recall the definition \eqref{definition-g} of $G$ and define $\widetilde{G}$ by 
$$G(v,w; Z_{0,.},\overline{\frakc}_{0,.})=-w^3+ \widetilde{G}(v,w; Z_{0,.},\overline{\frakc}_{0,.}).$$
Then, similarly to  \cite[Proposition 6.2]{MW} we obtain the relation 
\begin{multline}\label{ipp}
\frac{1}{3p-2} \Big(\big\|w_t \big\|^{3p-2}_{L^{3p-2}}-\big\|w_s \big\|^{3p-2}_{L^{3p-2}}\Big)+(3p-3) \int_s^t dr \, \big\| |\nabla w_r|^2 w_r^{3p-4}\big\|_{L^1}  \\
+(3p-3) \int_0^tdr\,  \big\| |x|^2 w^{3p-2}_r\big\|_{L^1}+\int_{s}^t dr \, \big\|w_r\big\|^{3p}_{L^{3p}} = \int_s^tdr\,  \langle  {3\, \<Psi2IPsi3>_{0,.}+\widetilde{G}(v,w; Z_{0,.},\overline{\frakc}_{0,.})}+\mu v , w^{3p-3} \rangle (r) .
\end{multline}

We decompose the right-hand side of \eqref{ipp} in the following way : 
\begin{equation*}
 {3\, \<Psi2IPsi3>_{0,.}+\widetilde{G}(v,w; Z_{0,.},\overline{\frakc}_{0,.})}+\mu v=\sum_{j=1}^6  {A_j(v,w; Z_{0,.},\overline{\frakc}_{0,.})}, 
\end{equation*}
with 
\begin{equation*}
 {A_1(v,w; Z_{0,.},\overline{\frakc}_{0,.})}:= -\big(3w^2 v+ 3wv^2+v^3\big)
\end{equation*}
\begin{equation*}
 {A_2(v,w; Z_{0,.},\overline{\frakc}_{0,.}):= -3 \mathrm{com}^{Z_0,.}_1(v,w)_{0,.} \pe \<Psi2>_{0,.}}
\end{equation*}
\begin{equation*}
 {A_3(v,w; Z_{0,.},\overline{\frakc}_{0,.}):= -3  w \pe \<Psi2>_{0,.}}   
\end{equation*}
\begin{equation*}
A_4(v,w; Z_{0,.},\overline{\frakc}_{0,.}):= \tau^{(2)}(Z_{0,.}) (v+w)^2
\end{equation*}
\begin{multline*}
 A_5(v,w; Z_{0,.},\overline{\frakc}_{0,.}):=-3\mathrm{com}^{Z_{0,.}}_2(v+w)_{0,.}-3(v+w-\<IPsi3>_{s,.} )\pg \<Psi2>_{s,.}\\
+\tau^{(0)}(Z_{0,.})+\tau^{(1)}(Z_{0,.})(v+w)+\mu v,
\end{multline*}
\begin{equation*}
{A_6(v,w; Z_{0,.},\overline{\frakc}_{0,.}):=3\, \<Psi2IPsi3>_{0,.}+\overline{\frakc}_{0,.}\big(\<Psi>_{0,.}-\<IPsi3>_{0,.}+v+w\big),}
\end{equation*}
 where the $\tau^{(j)}(Z_{0,.})$ are defined in \eqref{deft0}-\eqref{deft2}.

\

 As a consequence, from \eqref{ipp} we obtain the a priori bound 
\begin{equation*} 
 \big\|w_t \big\|^{3p-2}_{L^{3p-2}}+ \int_{s}^t dr\, \big\|w_r\big\|^{3p}_{L^{3p}} \lesssim \big\|w_s \big\|^{3p-2}_{L^{3p-2}}+  \sum_{j=1}^6 \int_s^t dr\,  \big| \langle {A_j(v,w; Z_{0,.},\overline{\frakc}_{0,.}} , w^{3p-3}\rangle(r) \big|  .
\end{equation*}

The rest of the proof of Theorem \ref{thm61} follows the same lines as \cite[Theorem 6.1]{MW}. To be more specific, the contributions of the terms $(A_j)_{j \notin \{4,6\}}$ can be treated exactly as in the latter reference -- by replacing the standard Besov spaces $B_p^{\alpha}(\R^3)$ therein with the spaces $\B_p^{\alpha}(\R^3)$ under consideration here -- and we only need to add a few details about the contributions of $A_4$ and $A_6$.

\

Let us focus on $A_4$ first. For this term, we need a few additional technicalities in comparison with \cite[Lemma 6.6]{MW}, since, in our situation, we simply cannot bound $\|w \|_{L^{3p-2}}$ by~$\|w \|_{L^{3p}}$ as on a torus.

 \begin{lemma}\label{lem:termea4}
 Let $\eps >0$ be small enough. Then for all $p\geq 1$ large enough, there exists an exponent $\kappa >0$ depending only on $\eps$ and $p$ such that  {for every $\mu$ satisfying \eqref{condiC}}, every $\delta \in (0,1]$ and all $0\leq s \leq t \leq T$, we have 
 \begin{multline}\label{618}
 \int_s^t dr \, \big|   \langle   A_4(v,w; Z_{0,.},\overline{\frakc}_{0,.}), w^{3p-3}\rangle(r) \big|      \lesssim   (\delta^{-1} M)^\kappa \Big[  1+\int_s^t  dr\,  \big\| w_r \big\|^p_{\B^{1+4\eps}_p} \Big] +\\
+ \delta \Big[   \big\| v_s\big\|^{3p}_{ \B^{-3\eps}_{2p}}    +   \int_s^t dr\,  \big\| w_r\big\|^{3p}_{L^{3p}} +  \sup_{s\leq r \leq t } \big\|w_r \big\|^{3p-2}_{L^{3p-2}}  \Big],
 \end{multline}
 where the implicit constant depends only on $\eps$ and $p$.
  \end{lemma}

  \begin{proof}
To alleviate notations, we simply write $\tau^{(2)}=\tau^{(2)}(Z_{0,.})$ in the proof. By definition \eqref{deft2}, we have $\tau^{(2)} \in \B^{-\frac12-\eps}_\infty $. 

\smallskip
 
\noindent
 {\it Step 1: Contribution of $\dis \int_s^t dr \,\Big|  \langle \tau^{(2)}w^2, w^{3p-3} \rangle(r) \Big|  $.}  As in \cite[Lemma 6.6]{MW} we can establish that
  \begin{equation*}
      \langle \tau^{(2)} w^2, w^{3p-3}\rangle \lesssim \delta \big\| w\big\|^{3p}_{L^{3p}} + \delta^{-\kappa}M^{2p} \big\|w \big\|^p_{\B^{1+4\eps}_p}.
   \end{equation*}
 
\noindent
  {\it Step 2:  Contribution of $\dis   \int_s^t dr \,\Big| \langle  \tau^{(2)}v^2, w^{3p-3} \rangle(r) \Big| $.} Let $0< \eta <\frac{\eps}4$.  We begin with the estimate 
 \begin{eqnarray*}
     \langle \tau^{(2)} v^2, w^{3p-3}\rangle =    \langle \langle x \rangle^\eta \tau^{(2)} , \langle x \rangle^{-\eta} v^2  w^{3p-3}\rangle &\lesssim & \big\|     \langle x \rangle^\eta  \tau^{(2)}  \big\|_{\B^{-\frac12-\frac{5\eps}4}_\infty}   \big\| \langle x \rangle^{-\eta}  v^2 w^{3p-3} \big\|_{\B^{\frac12+\frac32\eps}_1} \\
      &\lesssim &    M  \big\| \langle x \rangle^{-\eta}  v^2 w^{3p-3} \big\|_{\B^{\frac12+\frac32\eps}_1},
      \end{eqnarray*}
where we have used the estimate $\big\|     \langle x \rangle^\eta  {\tau^{(2)}}  \big\|_{\B^{-\frac12-\frac{5\eps}4}_\infty}   \lesssim \big\|     {\tau^{(2)}}  \big\|_{\B^{-\frac12-\eps}_\infty}  \lesssim M$, by \cite[Lemma 13.12]{DFT}.

 Then by Proposition \ref{Prop-est-para} $(v)$ 
 $$  \big\| \langle x \rangle^{-\eta} v^2 w^{3p-3} \big\|_{\B^{\frac12+\frac32\eps}_1} \lesssim \big\|v^2\big\|_{\B^{\frac12+\frac32\eps}_{3p-2}}   \big\| \langle x \rangle^{-\eta}  w^{3p-3} \big\|_{L^{\frac{3p-2}{3p-3}}} + \big\|v^2\big\|_{L^{\frac{3p}2}}  \big\| \langle x \rangle^{-\eta} w^{3p-3} \big\|_{\B_{\frac{3p}{3p-2}}^{\frac12+2\eps}}.  $$
 On the one hand, using Lemma \ref{lem-bis} we obtain   $\|v^2\|_{\B^{\frac12+\frac32\eps}_{3p-2}}    \lesssim  \|v^2\|_{\B^{\frac12+2\eps}_{\infty}}   $. Then we can prove that 
 \begin{equation}\label{114}
\big\| \langle x \rangle^{-\eta}  w^{3p-3} \big\|_{L^{\frac{3p-2}{3p-3}}}  \lesssim \big\|    w \big\|^{\frac{3p-6}2}_{L^{3p-2}}  \big\| w \big\|^{\frac{3p}2}_{L^{3p}},
 \end{equation}
 which is the analogous estimate to  \cite[(6.21)]{MW}. Namely, by the H\"older estimate, there exists $r \in (3p-2, 3p)$ such that 
 \begin{equation}\label{14}
\big\| w \big\|^{3p-3}_{L^{r}}  \lesssim \big\| w \big\|^{\frac{3p-6}2}_{L^{3p-2}}  \big\| w \big\|^{\frac{3p}2}_{L^{3p}}.
 \end{equation}
To be more precise, we use here the interpolation estimate
 $$\big\| w \big\|_{L^r} \leq \big\| w \big\|^{1-\theta}_{L^{p_0}} \big\| w \big\|^{\theta}_{L^{p_1}},$$
 with $p_0=3p-2$, $p_1=3p$, $\theta=\frac{3p}{2(3p-3)}$ and $\frac1r= \frac{1-\theta}{p_0}+ \frac{\theta}{p_1}$. 
 Set  $q=\frac{r}{3p-2} \in (1, \frac{3p}{3p-2})$, then 
  \begin{equation}\label{rhs}
\Big\| \langle x \rangle^{-\eta}  w^{3p-3} \Big\|^{\frac{3p-2}{3p-3}}_{L^{\frac{3p-2}{3p-3}}}  =\int   dx\,  \langle x \rangle^{-\eta \frac{3p-2}{3p-3}}   w^{3p-2}         \leq \Big(  \int  dx \,    \langle x \rangle^{-\eta q' \frac{3p-2}{3p-3}}      \Big)^{\frac1{q'}}  \Big(  \int dx  \,     w^{(3p-2)q}   \Big)^{\frac1{q}}.
 \end{equation}
 If $p\geq 1$ is large enough, then $q' \gg 1$ is large, so that the first integral in the right-hand side of \eqref{rhs} converges and we obtain 
   \begin{equation*}
\big\| \langle x \rangle^{-\eta}  w^{3p-3} \big\|_{L^{\frac{3p-2}{3p-3}}}           \lesssim   \big\|w \big\|^{3p-3}_{L^r} ,
 \end{equation*}
 which together with \eqref{14} implies \eqref{114}.

\smallskip
 
 On the other hand, by \cite[Lemma 13.12]{DFT}, 
 $$ \big\| \langle x \rangle^{-\eta} w^{3p-3} \big\|_{\B_{\frac{3p}{3p-2}}^{\frac12+2\eps}} \lesssim   \big\| w^{3p-3} \big\|_{\B_{\frac{3p}{3p-2}}^{\frac12+2\eps}}.  $$
 With these bounds, we can follow the argument of \cite[Lemma 6.6]{MW} to prove that $\dis \int_s^t  dr\, \Big| \langle  \tau^{(2)} v^2, w^{3p-3} \rangle(r) \Big|      $ is bounded by the right-hand side of \eqref{618}.
 \medskip
 
\noindent
{\it Step 3: Contribution of $\dis  \int_s^t dr\, \Big| \langle \tau^{(2)} vw, w^{3p-3} \rangle(r) \Big|    $.}  Let $0< \eta <\frac{\eps}4$, then as in step 2, we have the  estimate 
 \begin{eqnarray*}
     \langle \tau^{(2)} v w, w^{3p-3}\rangle =    \langle \langle x \rangle^\eta \tau^{(2)} , \langle x \rangle^{-\eta} v  w^{3p-2}\rangle &\lesssim & \big\|     \langle x \rangle^\eta  \tau^{(2)} \big\|_{\B^{-\frac12-\frac{5\eps}4}_\infty}   \big\| \langle x \rangle^{-\eta}  v w^{3p-2} \big\|_{\B^{\frac12+\frac32\eps}_1} \\
      &\lesssim &    M  \big\| \langle x \rangle^{-\eta}  v w^{3p-2} \big\|_{\B^{\frac12+\frac32\eps}_1}.
      \end{eqnarray*}
  Then 
 $$  \big\| \langle x \rangle^{-\eta} v w^{3p-2} \big\|_{\B^{\frac12+\frac32\eps}_1} \lesssim \big\|v\big\|_{\B^{\frac12+\frac32\eps}_{\infty}}   \big\| \langle x \rangle^{-\eta}  w^{3p-2} \big\|_{L^{1}} + \big\|v\big\|_{L^{{3p}}}  \big\| \langle x \rangle^{-\eta} w^{3p-2} \big\|_{\B_{\frac{3p}{3p-1}}^{\frac12+2\eps}}.  $$
 
 Similarly to the previous case, we can show that 
  \begin{equation*} 
\big\| \langle x \rangle^{-\eta}  w^{3p-2} \big\|_{L^{1}}  \lesssim \big\| w \big\|^{\frac{3p-4}2}_{L^{3p-2}}  \big\| w \big\|^{\frac{3p}2}_{L^{3p}}.
 \end{equation*}
Moreover, using \cite[Lemma 13.12]{DFT}, we obtain 
 $$ \big\| \langle x \rangle^{-\eta} w^{3p-2} \big\|_{\B_{\frac{3p}{3p-1}}^{\frac12+2\eps}} \lesssim  \big\| w^{3p-2} \big\|_{\B_{\frac{3p}{3p-1}}^{\frac12+2\eps}},$$
 and we can then complete the proof as in \cite[Lemma 6.6]{MW}.
  \end{proof}

	As for the contribution of $A_6$, we can rely on the following result.

  \begin{lemma}\label{lem-globali:a6}
 Let $\eps >0$ be small enough. Then for all $p\geq 1$ large enough, there exists an exponent $\kappa >0$ depending only on $\eps$ and $p$ such that for every {$\mu$ satisfying \eqref{condiC}}, every $\delta \in (0,1]$ and all $0\leq s \leq t \leq T$, we have 
\begin{multline}\label{estim-A6}
\int_s^t dr \, \big|\langle A_6(v,w; Z_{0,.},\overline{\frakc}_{0,.}), w^{3p-3}\rangle(r)\big|
\;\lesssim\;
(\delta^{-1}M)^{\kappa}\Big[1+\int_s^t dr\, \big\|w_r\big\|^p_{\B^{1+4\eps}_p}\Big]\\
+\;\delta\Big[\,1+\big\|v_s\big\|^{3p}_{\B^{-3\eps}_{2p}}
+\int_s^t dr\,\big\|w_r\big\|^{3p}_{L^{3p}}
+\sup_{s\leq r\leq t}\big\|w_r\big\|^{3p-2}_{L^{3p-2}}\Big],
\end{multline}
where the implicit constant depends only on $\eps$ and $p$.
  \end{lemma}

\begin{proof}

We treat the contributions of the diagram
$\<Psi2IPsi3>_{0,.}$ and the $\overline{\frakc}$ term
$\overline{\frakc}_{0,.}\big(\<Psi>_{0,.}-\<IPsi3>_{0,.}+v+w\big)$ separately.

\medskip

\noindent
\textit{Step 1: the diagram term $\<Psi2IPsi3>_{0,.}$.}
The term arising from $\<Psi2IPsi3>_{0,.}$ can be bounded as   
\begin{align*}
\big|\langle \<Psi2IPsi3>_{0,.},w^{3p-3}\rangle\big|& \leq M \big\|w^{3p-3}\big\|_{\cb_1^{1+\eps}}\leq c M \big\|w^{3p-4}\big\|_{L^{\frac{3p}{3p-4}}} \big\|w\big\|_{\cb_{\frac{3p}{4}}^{1+\eps}}\\
&\leq \Big( \delta^{\frac{3p-4}{3p}}\big\|w\big\|_{L^{3p}}^{3p-4}\Big) \Big(c M \delta^{-\frac{3p-4}{3p}} \big\|w\big\|_{\cb_{\frac{3p}{4}}^{1+\eps}}\Big) \\
&\leq \delta \big\|w\big\|_{L^{3p}}^{3p}+ \Big(c M \delta^{-\frac{3p-4}{3p}} \big\|w\big\|_{\cb_{\frac{3p}{4}}^{1+\eps}}\Big)^{\frac{3p}{4}}.
\end{align*}
Now, thanks to \eqref{lplq}, we have the 
embedding $\cb_p^{1+4\eps }\hookrightarrow\B^{1+\eps}_{\frac{3p}4}$, and we can  write 
$$\big\|w\big\|_{\cb_{\frac{3p}{4}}^{1+\eps}}^{\frac{3p}{4}} \leq c_p  \big( 1+\big\|w\big\|_{\cb_{\frac{3p}{4}}^{1+\eps}}^{p}\big) \leq c_p\big( 1+\big\|w\big\|_{\cb_{p}^{1+4\eps}}^{p}\big) $$
therefore we obtain
\begin{align}
\big|\langle  {\<Psi2IPsi3>_{0,.}},w^{3p-3}\rangle \big| \leq \delta \big\|w\big\|_{L^{3p}}^{3p}+c_p  M^{\frac{3p}{4}} \delta^{-\frac{3p-4}{4}} \big( 1+\big\|w\big\|_{\cb_{p}^{1+4\eps}}^{p}\big).\label{a6-1}
\end{align}
Integrating in time for $p \geq 1$ large enough, we obtain for the first term of $A_6$
\begin{equation}\label{A6-diag}
\int_s^t dr\,\big|\langle\<Psi2IPsi3>_{0,r},w_r^{3p-3}\rangle\big|
\leq \delta\int_s^t dr\,\big\|w_r\big\|_{L^{3p}}^{3p}
+ (\delta^{-1}M)^{\kappa}\Big[1+\int_s^t dr\,\big\|w_r\big\|_{\B^{1+4\eps}_p}^p\Big],
\end{equation}
which is an admissible contribution to the right-hand side of
\eqref{estim-A6}.

\

\noindent
\textit{Step 2: the $\overline{\frakc}$ term.}  We focus on the term with the lowest regularity, which is the contribution of $ \overline{\frakc}_{0,.}\<Psi>_{0,.} $. We have
\begin{equation*}
  \big|\langle \overline{\frakc}_{0,.}\<Psi>_{0,.},w^{3p-3}\rangle\big| \leq \big\|    \overline{\frakc}_{0,.}\<Psi>_{0,.}   \big\|_{\B^{-\frac12-\eps}_\infty}   \big\|w^{3p-3}\big\|_{\cb_1^{\frac12+\eps}}
   \lesssim  M \big\|   \overline{\frakc}_{0,.}   \big\|_{\B^{\frac12+3\eps}_\infty}  \big\|w^{3p-3}\big\|_{\cb_1^{\frac12+\eps}}.
    \end{equation*}

 Since $\big\llbracket  \overline{\frakc}_{0,.}\big\rrbracket_{3\eps,[0,T]}\leq M$, we have $\big\| \overline{\frakc}_{0,r}\big\|_{\B^{\frac12+3\eps}_\infty}\leq M\, r^{-\frac78}$ for all $0<r\leq T$, and we obtain, for all $0 <\delta\leq 1$,
\begin{align*}
\big|\langle \overline{\frakc}_{0,r}\<Psi>_{0,r},w^{3p-3}\rangle\big|& \leq c M^2\, r^{-\frac78}  \big\|w^{3p-3}\big\|_{\cb_1^{\frac12+\eps}}\\
& \leq     c M^2\, r^{-\frac78}   \big\|w^{3p-4}\big\|_{L^{\frac{3p-2}{3p-4}}} \big\|w\big\|_{\cb_{\frac{3p-2}{2}}^{\frac12+\eps}}\\
&\leq \Big( \delta^{\frac{3p-4}{3p-2}}\frac{r^{-\frac78}}{A} \big\|w\big\|_{L^{3p-2}}^{3p-4}\Big) \Big(c A M^2\delta^{-\frac{3p-4}{3p-2}} \big\|w\big\|_{\cb_{\frac{3p-2}{2}}^{\frac12+\eps}}\Big) ,
 \end{align*}
for some parameter $A>0$ to be specified later. Thanks to the Young inequality, we obtain 
\begin{align}
\big|\langle \overline{\frakc}_{0,r}\<Psi>_{0,r},w^{3p-3}\rangle\big|  & \leq  \delta \frac{r^{-\frac{7(3p-2)}{8(3p-4)}}}{A^{\frac{3p-2}{3p-4}}}\big\|w\big\|_{L^{3p-2}}^{3p-2}+ c \,A^{\frac{3p-2}{2}}M^{3p-2} \,\delta^{-\frac{3p-4}{2}} \big\|w\big\|_{\cb_{\frac{3p-2}{2}}^{\frac12+\eps}}^{\frac{3p-2}{2}}.\label{sugge-1}
\end{align}
Since the parameter $A$ will be fixed to a finite value in \eqref{defparaguay}, the factor $A^{\frac{3p-2}{2}}$ can henceforth be absorbed into the constant $c$, and we drop it in the sequel. By interpolation, we can check that
\begin{align*}
&\big\|w\big\|_{\cb_{\frac{3p-2}{2}}^{\frac12+\eps}}^{\frac{3p-2}{2}} \lesssim \big\|w\big\|_{L^{3p}}^{\frac{3p-6}{4}}\big\|w\big\|_{\cb_{p}^{(1+2\eps)\frac{3p-2}{3p+2}}}^{\frac{3p+2}{4}}\lesssim \big\|w\big\|_{L^{3p}}^{\frac{3p-6}{4}}\big\|w\big\|_{\cb_{p}^{1+2\eps}}^{\frac{3p+2}{4}}
\end{align*}
which yields by the Young inequality
\begin{align*}
c M^{3p-2}\delta^{-\frac{3p-4}{2}}\big\|w\big\|_{\cb_{\frac{3p-2}{2}}^{\frac12+\eps}}^{\frac{3p-2}{2}} &\leq \Big(\delta^{\frac{p-2}{4p}} \big\|w\big\|_{L^{3p}}^{\frac{3p-6}{4}}\Big) \Big(c M^{3p-2}\delta^{-\frac{6p^2-7p-2}{4p}} \big\|w\big\|_{\cb_{p}^{1+2\eps}}^{\frac{3p+2}{4}}\Big)\\
&\leq \delta \big\|w\big\|_{L^{3p}}^{3p}+c\big(\delta^{-1}M\big)^\ka\big\|w\big\|_{\cb_{p}^{1+2\eps}}^{p},
\end{align*}
for some large $\ka\geq 1$. Going back to \eqref{sugge-1}, we deduce
 \begin{align}
\big|\langle \overline{\frakc}_{0,r}\<Psi>_{0,r},w^{3p-3}\rangle\big|  & \leq  \delta \frac{r^{-\frac{7(3p-2)}{8(3p-4)}}}{A^{\frac{3p-2}{3p-4}}}\big\|w\big\|_{L^{3p-2}}^{3p-2}+\delta \big\|w\big\|_{L^{3p}}^{3p}+c\big(\delta^{-1}M\big)^\ka\big\|w\big\|_{\cb_{p}^{1+2\eps}}^{p}.
\end{align}

\

For  $p \geq 1$ large enough, the singularity in $r$ is integrable, and setting
\begin{equation}\label{defparaguay}
A:= \Big(\int_0^1 dr\, r^{-\frac{7(3p-2)}{8(3p-4)}}\Big)^{\frac{3p-4}{3p-2}},
\end{equation}
we obtain
\begin{equation*} 
\int_s^t dr\,\big|\langle \overline{\frakc}_{0,r}\<Psi>_{0,r},w^{3p-3}_r\rangle\big| 
\leq \delta\Big[\sup_{s\leq r\leq t}\big\|w_r\big\|^{3p-2}_{L^{3p-2}}+\int_s^t dr\,\big\|w_r\big\|^{3p}_{L^{3p}}\Big]
+ (\delta^{-1}M)^{\kappa}\Big[1+\int_s^t dr\,\big\|w_r\big\|_{\B^{1+4\eps}_p}^p\Big],
\end{equation*}
which is an admissible contribution to the right-hand side of \eqref{estim-A6}. 

\smallskip

The contribution of $\overline{\frakc}_{0,.}\<IPsi3>_{0,.}$ can be treated with the same arguments (the regularity of $\<IPsi3>$ being better than that of $\<Psi>$).

\smallskip

The contributions of $\overline{\frakc}_{0,.}v$ and $\overline{\frakc}_{0,.}w$ are easier to treat and we omit the details here.

\medskip

Collecting Steps 1 and 2 gives \eqref{estim-A6}.
\end{proof}


\subsection{Conclusion of the argument}

  We are now in a position to combine all the preceding estimates to close the a priori bounds. The following result is thus the counterpart of \cite[Theorem~7.1]{MW}.
	
 \begin{theorem}\label{thm411}
 Let $\eps, \nu >0$ be small enough. Then for all $p\geq 1$ large enough, there exists an exponent $\kappa \geq 1$ such that  for all $\mu$ satisfying
 \begin{equation}\label{condiC-ter}
\mu \geq \max\big(1+\big(8  M \big)^8,\mu_0 M^{30p}\big),
 \end{equation}
where $\mu_0=\mu_0(\eps,p)>0$ is the constant coming from Theorem \ref{thm61}, and for every $t\in [0,T]$, we have 
 \begin{equation} \label{thm411-1}
 \big\| v_t \big\|_{\B^{-3\eps}_{2p}}+t^{\frac12+\nu}\big\| w_t \big\|_{L^{3p-2}} \lesssim {(\mu M)}^\kappa,
 \end{equation}
 where the implicit constant depends only on $\eps$ and $p$.
  \end{theorem}

\subsubsection{Preliminary controls}

 \begin{lemma}\label{lemma411}
 Let $\eps >0$ be small enough. Then for all $p\geq 1$ large enough, there exists $\kappa \geq 1$ depending only on $\eps$ and $p$ such that {for all $\mu$ satisfying \eqref{condiC-ter}} and for all $s\leq t\in [0,T]$, we have 
  \begin{equation}\label{est-c}
 \big\|w_t\big\|^{3p-2}_{L^{3p-2}}+\int_s^t dr\,  \big\|w_r\big\|^{3p}_{L^{3p}}+\int_{s}^t dr\,   \big\|w_r\big\|^{p}_{\B^{1+7\eps}_p}    
 \lesssim {(\mu M)}^{\kappa} \Big( 1+  \big\|w_s\big\|^{3p-2}_{L^{3p-2}}+\big\|v_s\big\|^{3p}_{\B^{-3\eps}_{2p}} +\big\|w_s \big\|^{\frac{3p-2}3}_{\B^{1+7\eps}_p}\Big),
 \end{equation}
 where the implicit constant depends only on $\eps$ and $p$.
  \end{lemma}

\begin{proof}[Proof of Lemma \ref{lemma411}]
The proof follows the same lines as \cite[Lemma 7.2]{MW}, but there are two main changes to make: \medskip

1) The estimate  \cite[(7.5)]{MW} has to be replaced with the estimate \eqref{foll}.\medskip

2) In the second part of the proof, we choose here $\gamma=(1+ \frac{13}2\eps)(1-\frac1p)$. Then by interpolation 
  \begin{equation*}
\big\|w \big\|^p_{\B^{\gamma}_p} \lesssim   \big\|w \big\|^{p-1}_{\B^{1+\frac{13}2\eps}_{q}} \big\|w \big\|_{L^{3p}}, 
\end{equation*}
with $\dis q=\frac{3p-2}{3}$. Then, with \eqref{lplq} we obtain that for $p\geq 1$ large enough,
  \begin{equation*}
\big\|w \big\|^p_{\B^{\gamma}_p} \lesssim   \big\|w \big\|^{p-1}_{\B^{1+7\eps}_{p}} \big\|w \big\|_{L^{3p}} \lesssim  \big\|w \big\|^{\frac{3p-2}3}_{\B^{1+7\eps}_{p}}+  \big\|w \big\|^{3p-2}_{L^{3p}}.
\end{equation*}
This inequality allows us to conclude as in \cite[Lemma 7.2]{MW}.
\end{proof}

\begin{corollary} \label{coro:to-gron}
 Let $\eps >0$ be small enough. Then for all $p\geq 1$ large enough, there exist $\kappa \geq 1$ and $0< \sigma<p$ depending only on $\eps$ and $p$ such that for {all $\mu$ satisfying \eqref{condiC-ter}} and all $s\leq t\in [0,T]$, we have 
 \begin{equation*}
 \big\|w_t\big\|^{3p-2}_{L^{3p-2}}+\int_s^t dr\,  \big\|w_r\big\|^{3p}_{L^{3p}}+\int_{s}^t  dr\,   \big\|w_r\big\|^{p}_{\B^{1+7\eps}_p}    
 \lesssim {(\mu M)}^{\kappa} \Big( 1+  \big\|w_s\big\|^{3\sigma}_{L^{3p}}+\big\|v_s\big\|^{3p}_{\B^{-3\eps}_{2p}} +\big\|w_s \big\|^{\sigma}_{\B^{1+7\eps}_p}\Big),
 \end{equation*}
 where the implicit constant depends only on $\eps$ and $p$.
  \end{corollary}
  
\begin{proof}
By the H\"older estimate and \eqref{lplq}, for all $\gamma>1$
\begin{equation}\label{erq1}
\big\|w\big\|_{L^{3p-2}} \leq \big\|w\big\|^{\theta}_{L^{3p}} \big\|w\big\|^{1-\theta}_{L^{3}} \lesssim \big\|w\big\|^{\theta}_{L^{3p}} \big\|w\big\|^{1-\theta}_{\B^{\gamma}_{p}},
\end{equation}
with $\dis \theta=\frac{p(3p-5)}{(p-1)(3p-2)}<1$, which is close to 1 since $p \gg 1$. By the Young inequality, the bound~\eqref{erq1} implies for all $q\geq 1$
\begin{equation*}
\big\|w\big\|_{L^{3p-2}} \lesssim    \big\|w\big\|^{q\theta}_{L^{3p}}+ \big\|w\big\|^{(1-\theta)q'}_{\B^{\gamma}_{p}}
\end{equation*}
thus, for all $q\geq 1$ and $\ga >1$,
\begin{equation}\label{using425}
 \big\|w\big\|^{3p-2}_{L^{3p-2}} +\big\|w\big\|^{\frac{3p-2}3}_{\B^{\gamma}_{p}}\lesssim    \big\|w\big\|^{q\theta(3p-2)}_{L^{3p}}+ \big\|w\big\|^{(1-\theta)q'(3p-2)}_{\B^{\gamma}_{p}}+ \big\|w\big\|^{\frac{3p-2}3}_{\B^{\gamma}_{p}}.
\end{equation}
We choose $1<q<\infty$ such that $(1-\theta)q'=\frac13$, namely 
$$q'=\frac{(p-1)(3p-2)}6, \qquad q=\frac{(p-1)(3p-2)}{3p^2-5p-4}. $$
We check that $3p-2<q \theta (3p-2)< 3p$. Set 
\begin{equation}\label{def-sigma}
\sigma:=\frac{q \theta (3p-2)}3<p.
\end{equation}
   Then using \eqref{using425} 
\begin{equation}\label{bsigma}
 \big\|w\big\|^{3p-2}_{L^{3p-2}} +\big\|w\big\|^{\frac{3p-2}3}_{\B^{\gamma}_{p}}\lesssim    1+ \big\|w\big\|^{3\sigma}_{L^{3p}}+ \big\|w\big\|^{\sigma}_{\B^{\gamma}_{p}}.
\end{equation}
We now set $\gamma=1+7\eps$ and combine the previous estimate with \eqref{est-c} to obtain the desired bound.
\end{proof}

In order to obtain a time decay for $w$ we will use the following result. This variation of \cite[Lemma 7.3]{MW} allows us to avoid working with time sequences.

\begin{lemma}\label{lem:pointwise}
Let $\tau>0$, $\beta>1$, $c\geq 1$ and $A\ge 1$. Let $F, L:[0,\tau]\to[0,\infty)$ be continuous functions 
satisfying, for all $0\le s\le t \leq \tau$,
\begin{equation}\label{eq:H1}
\int_s^t dr \, F^\beta(r) \le c\,F(s), 
\end{equation}
and
\begin{equation}\label{eq:H2}
L(t)\le A F(s).
\end{equation}
Then, there exists $C=C(\beta)>0$ such that for every $t\in{(0,\tau]}$,
\[
L(t) \leq C A\, c^{\frac{1}{\beta-1}}\,t^{-\frac1{\beta-1}}.
\]
\end{lemma}

\begin{proof}

Let $\dis G(s):=\int_s^\tau dr\,F^\beta(r)$.
It satisfies $G(s)\le cF(s)<\infty$ for every $s\leq \tau$.
The function~$G$ is  continuously differentiable   and $G'=-F^\beta$.
From $G\le cF$ we obtain
$F^\beta\ge (G/c)^\beta$, hence $-G'\ge c^{-\beta}G^\beta$. Wherever $G>0$,
\begin{equation}\label{ineq:H2}
\frac{d}{ds}\,G^{1-\beta}=(\beta-1)\,G^{-\beta}(-G')\;\ge\;\frac{\beta-1}{c^\beta}.
\end{equation}
Integrating  \eqref{ineq:H2} over $(0,s)$ gives  {(since $G(0)^{1-\beta} \ge 0$)}
\[
G(s)^{1-\beta} \geq G(s)^{1-\beta}-G(0)^{1-\beta}\;\ge\;\frac{\beta-1}{c^\beta}\,s,
\]
that is, 
\begin{equation}\label{defgg}
G(s)\;\le\;\Big(\frac{c^\beta}{(\beta-1)\,s}\Big)^{\!\frac1{\beta-1}}.
\end{equation}

\smallskip
Fix $t\in(0,\tau)$.
By continuity, $F$ attains its minimum over $[t/2,t]$ at some point $r^\star \in [t/2,t]$. Since $F\geq 0$, \eqref{defgg} directly yields
 \[ F^\beta(r^\star)\;\le\;\frac{2}{t}\int_{\frac{t}2}^{t}dr \, F^\beta(r)
\;\le\;\frac{2}{t}\,G \Big(\frac t2\Big)
\;\le\;\frac{2}{t}\Big(\frac{2c^\beta}{(\beta-1)\,t}\Big)^{\!\frac1{\beta-1}}, \]
which in turn implies
\begin{equation}\label{Z-2}
F(r^\star)\;\le\;C(\beta)\,c^{\frac1{\beta-1}}\,t^{-\frac1{\beta-1}},
\qquad
C(\beta)=2^{\frac1\beta}\Big(\tfrac{2}{\beta-1}\Big)^{\!\frac{1}{\beta(\beta-1)}}.
\end{equation}

\smallskip
Since $r^\star\le t{\ \leq \tau}$, we apply~\eqref{eq:H2} with $s=r^\star$, and then by~\eqref{Z-2}, 
\begin{equation*} 
L(t)\;\le\;A F(r^\star) 
\;\le\;A C(\beta)\,c^{\frac1{\beta-1}}\,t^{-\frac1{\beta-1}},
\end{equation*}
which was the claim.
\end{proof}

\subsubsection{{Proof of Theorem \ref{thm411}}}

\

\smallskip

Consider $\sigma$ as defined in \eqref{def-sigma}. A direct computation shows that $\sigma-p \to -\frac{2}3$   as $p\to\infty$, which in turn implies  that 
\begin{equation}\label{dl}
\frac{\sigma}{(3p-2)(p-\sigma)}
\;\xrightarrow[p\to\infty]{}\;\tfrac{1}{2}.
\end{equation}
Let us set 
$$F(t):=   \big\|w_t\big\|^{3\sigma}_{L^{3p}}+\big\|w_t\big\|^{\sigma}_{\B^{1+7\eps}_{p}} , \qquad L(t):=\big\| w_t \big\|^{3p-2}_{L^{3p-2}},$$
and define 
\begin{equation}\label{def:stop-time-tau}
\tau := \inf \Big\{ t \geqslant 0 : F(t) \leq  1 \text{ or }   F(t)\leq \big\|v_t\big\|^{3p}_{\B^{-3\eps}_{2p}}  \Big\} \wedge T. 
\end{equation}
 
Then,  by Corollary \ref{coro:to-gron}, the continuous functions $F$ and $L$ satisfy the assumptions of Lemma~\ref{lem:pointwise} {on $[0,\tau]$}, with $\beta=p/\sigma>1$, $c={(\mu M)}^\kappa$ and some $A \lesssim {(\mu M)}^{\kappa}$. 

\smallskip

Therefore, {for every $t\in (0,\tau]$},
 \begin{equation*} 
\big\| w_t \big\|_{L^{3p-2}} \lesssim {(\mu M)}^{\kappa_1} t^{-\frac{\sigma}{(3p-2)(p-\sigma)}} \lesssim {(\mu M)}^{\kappa_1} t^{-\frac12-\nu},
 \end{equation*}
for some $\nu=\nu(p)>0$, where we have used \eqref{dl} to derive the second estimate. 

\smallskip

By inserting this bound into \eqref{est-v-thm21}, we also obtain {that for every $t\in (0,\tau]$}
 \begin{equation}\label{bouvinter} 
\big\| v_t \big\|_{\B^{-3\eps}_{2p}}  \lesssim {(\mu M)}^{\kappa_2}
 \end{equation}
which completes the proof of {\eqref{thm411-1} on $[0,\tau]$}. 

\

To extend the estimate from $[0,\tau]$ to $[0,T]$, observe that by definition of $\tau$, one has
$$F(\tau) \leq 1+\big\|v_t\big\|^{3p}_{\B^{-3\eps}_{2p}}.$$
Therefore, by Corollary \ref{coro:to-gron}, we deduce that for every $t\in [\tau,T]$,
\begin{align*}
\big\|w_t\big\|^{3p-2}_{L^{3p-2}}   &\lesssim {(\mu M)}^{\kappa} \Big( 1+  F(\tau)+\big\|v_\tau\big\|^{3p}_{\B^{-3\eps}_{2p}} \Big)\lesssim {(\mu M)}^{\kappa} \Big( 1+\big\|v_\tau\big\|^{3p}_{\B^{-3\eps}_{2p}} \Big),
\end{align*}
which, by \eqref{bouvinter}, entails that $\big\|w_t\big\|^{3p-2}_{L^{3p-2}}  \lesssim M^{\kappa_3}$ for every $t\in [\tau,T]$. Finally, for every $t\in [0,T]$, one has
 \begin{equation*} 
\big\| w_t \big\|_{L^{3p-2}}  \lesssim {(\mu M)}^{\kappa_4} t^{-\frac12-\nu},
 \end{equation*}
and we can insert this control into \eqref{est-v-thm21} to complete the proof of the statement.

\

\subsection{Improved estimates}

 The estimates in Theorem \ref{thm411} are sufficient to prove the coming-down-from-infinity property (that is, Corollary \ref{bound_final}) and to establish the globalization of the solution. However, to obtain the more general bounds in Theorem \ref{thm:global_Lp}, which will be required in the sequel, we need the refinements below.

\smallskip

\subsubsection{Refinement of Theorem \ref{thm51}}

We now state the following (a posteriori) improvement of Theorem \ref{thm51}.

\begin{theorem}\label{improv-thm51}
Fix $0<\ga <\frac32$ {and $\eps >0$ small enough (depending on $\ga$). Then for all $p\geq 1$ large enough, there exists $\ka\geq 1$ such that for all $\mu$ satisfying \eqref{condiC-ter}} and all $0\leq s<t\leq {T}$,
 \begin{align}
\big\|w_t\big\|_{\B^\gamma_p}& \lesssim \big\| e^{-(t-s) H}w_s \big\|_{\B^\gamma_p}+ {(\mu M)}^\ka \bigg[1+   \bigg(\int_s^t dr\,   \big\|w_r\big\|^{3p}_{L^{3p}}\bigg)^{\frac1p}+  \bigg(\int_s^tdr\,   \big\|w_r\big\|^{p}_{\B^{1+4\eps}_p}\bigg)^{\frac1p}\bigg].\label{boutheo51st}
 \end{align}

\end{theorem}

Just as Theorem \ref{thm51}, Theorem \ref{improv-thm51} can be proved along the same lines as \cite[Theorem 5.1]{MW}, except for the following improvement over \cite[Lemma 5.6]{MW} (recall that $\tau^{(2)}(Z_{0,r})$ was defined in~\eqref{deft2}).

\begin{lemma}\label{lem:saucisse}
Fix $0<\ga <\frac32$ {and $\eps >0$ be small enough (depending on $\ga$). Then for all $p\geq 1$ large enough, there exists $\ka\geq 1$ such that for all $\mu$ satisfying \eqref{condiC-ter}} and all $0\leq s<t\leq T$,
\begin{align*}
&\bigg\| \int_s^t dr\, e^{-(t-r)H} \big( \tau^{(2)}(Z_{0,r})\cdot (v_r+w_r)^2 \big) \bigg\|_{\cb^\ga_p}\lesssim {(\mu M)}^\ka \Bigg[ 1+\bigg(\int_s^t dr\, \big\|w_r\big\|_{L^{3p}}^{3p}\bigg)^{\frac1p}+\bigg(\int_s^t dr\, \big\|w_r\big\|_{\cb^{1+4\eps}_p}^{p}\bigg)^{\frac1p}\Bigg].
\end{align*}

\end{lemma}

\begin{proof}
According to Theorem~\ref{thm411}, one has first, for all $p\geq 1$ large enough,
\begin{equation}\label{bounve-1}
\sup_{r\in [0,{T}]} \Big( \big\| v_r \big\|_{\B^{-\eps}_{p}}+r^{\frac12+\eps}\big\| w_r \big\|_{L^{p}} \Big)\lesssim {(\mu M)}^\kappa
\end{equation}
and from here observe that 
\begin{align}
 &\big\|v_r\big\|_{\cb^{\frac12+\eps}_p}    \lesssim \int_0^r ds \, \Big\| e^{-(r-s)H}\Big(\big(v_s+w_s -\<IPsi3>_{0,s}\big) \pl  \<Psi2>_{0,s}\Big) \Big\|_{\cb^{\frac12+\eps}_p}\nonumber\\
&   \lesssim \int_0^r \frac{ds}{(r-s)^{\frac34+\frac{3\eps}{2}}} \, \Big\|\big(v_s+w_s -\<IPsi3>_{0,s}\big) \pl  \<Psi2>_{0,s} \Big\|_{\cb^{-1-2\eps}_p}\nonumber\\
&\lesssim M\int_0^r \frac{ds}{(r-s)^{\frac34+\frac{3\eps}{2}}} \, \big\|v_s+w_s -\<IPsi3>_{0,s}\big\|_{\cb^{-\eps}_p}\lesssim {(\mu M)}^\ka\int_0^r \frac{ds}{(r-s)^{\frac34+\frac{3\eps}{2}}} \frac{1}{s^{\frac12 {+\eps}}} \lesssim \frac{{(\mu M)}^\ka}{r^{\frac14+3\eps}}.\label{bounve-2}
\end{align}
Interpolating between \eqref{bounve-1} and \eqref{bounve-2} also yields
\begin{equation}\label{bounve-3}
\big\|v_r\big\|_{L^p}\lesssim \big\|v_r\big\|_{\cb^\eps_p} \lesssim \big\|v_r\big\|_{\cb^{-\eps}_p}^{\frac{1}{1+4\eps}}\big\|v_r\big\|_{\cb^{\frac12+\eps}_p}^{\frac{4\eps}{1+4\eps}}\lesssim \frac{{(\mu M)}^\ka}{r^{\eta}},
\end{equation}
where, for the sake of clarity, we have set $\eta:=\frac{\eps+12\eps^2}{1+4\eps}$.

\smallskip

Let us now write
\begin{align*}
&\bigg\| \int_s^t dr \, e^{-(t-r)H} \big( \tau^{(2)}(Z_{0,r})\cdot (v_r+w_r)^2 \big) \bigg\|_{\cb^\ga_p} \leq\\
&\leq \bigg\| \int_s^t dr\, e^{-(t-r)H} \big(  (v_r+w_r)^2\pl \tau^{(2)}(Z_{0,r})\big) \bigg\|_{\cb^\ga_p}+\bigg\| \int_s^t dr\, e^{-(t-r)H} \big( \tau^{(2)}(Z_{0,r})\ple (v_r+w_r)^2 \big) \bigg\|_{\cb^\ga_p}.
\end{align*}
On the one hand,  thanks to Proposition \ref{Prop-est-para} $(iii)$, we obtain that 
\begin{multline*}
\bigg\| \int_s^t dr\, e^{-(t-r)H} \big(  (v_r+w_r)^2\pl \tau^{(2)}(Z_{0,r})\big) \bigg\|_{\cb^\ga_p}\lesssim \\
\begin{aligned}
&\lesssim \int_s^t \frac{dr}{|t-r|^{\frac12(\ga+\frac12+2\eps)}} \big\| (v_r+w_r)^2\pl \tau^{(2)}(Z_{0,r}) \big\|_{\cb^{-\frac12-2\eps}_p}\\
&\lesssim M\int_s^t \frac{dr}{|t-r|^{\frac12(\ga+\frac12+2\eps)}} \big\| (v_r+w_r)^2\big\|_{\cb^{-\eps}_p}.
\end{aligned}
\end{multline*}
With \eqref{lplq} we can write for all $p \geq 1$ large enough, 
$$ \big\| (v_r+w_r)^2\big\|_{\cb^{-\eps}_p} \lesssim  \big\| (v_r+w_r)^2\big\|_{\cb^{0}_{\frac{3p}{2}}} \lesssim \big\| (v_r+w_r)^2\big\|_{L^{\frac{3p}{2}}}  \lesssim \big\| v_r\big\|^2_{L^{3p}}+ \big\| w_r\big\|^2_{L^{3p}},$$
then, thanks to \eqref{bounve-3},
\begin{multline*}
\bigg\| \int_s^t dr\, e^{-(t-r)H} \big(  (v_r+w_r)^2\pl \tau^{(2)}(Z_{0,r})\big) \bigg\|_{\cb^\ga_p}\lesssim \\
\begin{aligned}
&\lesssim M\int_s^t \frac{dr}{|t-r|^{\frac12(\ga+\frac12+2\eps)}} \big\| v_r\big\|_{L^{3p}}^2+M\int_s^t \frac{dr}{|t-r|^{\frac12(\ga+\frac12+2\eps)}} \big\| w_r\big\|_{L^{3p}}^2\\
&\lesssim {(\mu M)}^\ka\bigg[\int_s^t \frac{dr}{|t-r|^{\frac12(\ga+\frac12+2\eps)}r^{2\eta}} +\int_s^t \frac{dr}{|t-r|^{\frac12(\ga+\frac12+2\eps)}} \big\| w_r\big\|_{L^{3p}}^3\bigg]\\
&\lesssim {(\mu M)}^\ka\bigg[1+  \bigg(\int_s^t dr\, \big\|w_r\big\|_{L^{3p}}^{3p}\bigg)^{\frac1p}\bigg]
\end{aligned}
\end{multline*}
for $p\geq 1$ large enough, due to the fact that $\frac{\gamma}2+\frac14+\eps+2\eta<1$ for $\eps>0$ small enough. 

\

On the other hand,  thanks to Proposition \ref{Prop-est-para} (items $(i)$, $(ii)$ and $(v)$), and using \eqref{bounve-2} and~\eqref{bounve-3} as well, one has
\begin{multline*} 
\big\| \tau^{(2)}(Z_{0,r})\ple (v_r+w_r)^2 \big\|_{\cb^\eps_p} \lesssim M
\big\| (v_r+w_r)^2 \big\|_{\cb^{\frac12+2\eps}_p}\lesssim \\
\lesssim M \big\| v_r+w_r \big\|_{L^{3p}} \big\| v_r+w_r \big\|_{\cb^{\frac12+2\eps}_{\frac{3p}{2}}}
\lesssim {(\mu M)}^\ka  \Big( \frac{1}{r^\eta}+\big\| w_r \big\|_{L^{3p}}\Big) \Big(\frac{1}{r^{\frac14+4\eps}}+\big\| w_r \big\|_{\cb^{\frac12+2\eps}_{\frac{3p}{2}}}\Big).
\end{multline*}
Therefore by \eqref{heat1}
\begin{multline}\label{decominte}
\bigg\| \int_s^t dr\, e^{-(t-r)H} \big( \tau^{(2)}(Z_{0,r})\ple (v_r+w_r)^2 \big) \bigg\|_{\cb^\ga_p}\lesssim \\
\begin{aligned}
&\lesssim  \int_s^t \frac{dr}{|t-r|^{\frac{\gamma}2}}  \big\| \tau^{(2)}(Z_{0,r})\ple (v_r+w_r)^2 \big\|_{\cb^\eps_p} \\
&\lesssim {(\mu M)}^\ka \bigg[1+\int_s^t \frac{dr}{|t-r|^{\frac\ga 2}r^\eta}\big\| w_r \big\|_{\cb^{\frac12+2\eps}_{\frac{3p}{2}}}+\int_s^t \frac{dr}{|t-r|^{\frac{\ga}2 }r^{\frac14+4\eps}}\big\| w_r \big\|_{L^{3p}} \\
&\hspace{7cm}+\int_s^t \frac{dr}{|t-r|^{\frac{\gamma}2}}\big\| w_r \big\|_{L^{3p}}\big\| w_r \big\|_{\cb^{\frac12+2\eps}_{\frac{3p}{2}}}\bigg],
\end{aligned}
\end{multline}
due to the fact that $\frac{\gamma}2+\eta+\frac14+4\eps<1$ for $\eps>0$ small enough.

\smallskip

Then, using the fact that for $p$ large enough, one has (see Lemma \ref{lem:inclusion-besov})
$$\big\| w_r \big\|_{\cb^{\frac12+2\eps}_{\frac{3p}{2}}} \lesssim \big\| w_r \big\|_{\cb^{\frac12+3\eps}_{p}} \lesssim \big\| w_r \big\|_{\cb^{1+3\eps}_{p}} ,$$
the first two integrals in \eqref{decominte} can be bounded for $p\ge 1$ large enough as
\begin{multline*}
1+\int_s^t \frac{dr}{|t-r|^{\frac{\ga}2}r^\eta}\big\| w_r \big\|_{\cb^{\frac12+2\eps}_{\frac{3p}{2}}}+\int_s^t \frac{dr}{|t-r|^{\frac{\ga}2}r^{\frac14+4\eps}}\big\| w_r \big\|_{L^{3p}}\lesssim \\
\begin{aligned}
&\lesssim 1+\int_s^t \frac{dr}{|t-r|^{\frac{\ga}2}r^\eta}\big\| w_r \big\|_{\cb^{1+3\eps}_{p}}+\int_s^t \frac{dr}{|t-r|^{\frac{\ga}2}r^{\frac14+4\eps}}\big\| w_r \big\|_{L^{3p}}^3\\
&\lesssim 1+\bigg(\int_s^t dr\, \big\| w_r \big\|_{\cb^{1+3\eps}_{p}}^p\bigg)^{\frac1p}+\bigg(\int_s^t dr\, \big\| w_r \big\|_{L^{3p}}^{3p}\bigg)^{\frac1p}.
\end{aligned}
\end{multline*}

As for the last integral in \eqref{decominte}, we have for $p\geq 1$ large enough,
\begin{align}
&\int_s^t \frac{dr}{|t-r|^{\frac{\ga}2}}\big\| w_r \big\|_{L^{3p}}\big\| w_r \big\|_{\cb^{\frac12+2\eps}_{\frac{3p}{2}}}\lesssim \bigg(\int_s^t dr\, \big\| w_r \big\|_{L^{3p}}^{3p}\bigg)^{\frac{1}{3p}} \bigg(\int_s^t dr\, \big\| w_r \big\|_{\cb^{\frac12+2\eps}_{\frac{3p}{2}}}^p\bigg)^{\frac1p}.\label{lasest}
\end{align}
Since
$$\big\| w_r \big\|_{\cb^{\frac12+2\eps}_{\frac{3p}{2}}}\lesssim \big\| w_r \big\|_{\cb^{0}_{3p}}^{\frac12} \big\| w_r \big\|_{\cb^{1+4\eps}_{p}}^{\frac12} \lesssim \big\| w_r \big\|_{L^{3p}}^{\frac12} \big\| w_r \big\|_{\cb^{1+4\eps}_{p}}^{\frac12},$$
we deduce that
\begin{align*}
&\bigg(\int_s^t dr\, \big\| w_r \big\|_{\cb^{\frac12+2\eps}_{\frac{3p}{2}}}^p\bigg)^{\frac1p}\lesssim \bigg(\int_s^t dr\, \big\| w_r \big\|_{L^{3p}}^{\frac{p}{2}} \big\| w_r \big\|_{\cb^{1+4\eps}_{p}}^{\frac{p}{2}}\bigg)^{\frac1p}\\
&\hspace{1cm}\lesssim \bigg(\int_s^t dr\, \big\| w_r \big\|_{L^{3p}}^{3p}\bigg)^{\frac{1}{6p}}\bigg(\int_s^t dr\, \big\| w_r \big\|_{\cb^{1+4\eps}_{p}}^{\frac{3p}{5}}\bigg)^{\frac{5}{6p}}\lesssim \bigg(\int_s^t dr\, \big\| w_r \big\|_{L^{3p}}^{3p}\bigg)^{\frac{1}{6p}}\bigg(\int_s^t dr\, \big\| w_r \big\|_{\cb^{1+4\eps}_{p}}^{p}\bigg)^{\frac{1}{2p}}.
\end{align*}
Going back to \eqref{lasest}, we obtain
\begin{align*}
\int_s^t \frac{dr}{|t-r|^{\frac{\ga}2}}\big\| w_r \big\|_{L^{3p}}\big\| w_r \big\|_{\cb^{\frac12+2\eps}_{\frac{3p}{2}}}&\lesssim \bigg(\int_s^t dr\, \big\| w_r \big\|_{L^{3p}}^{3p}\bigg)^{\frac{1}{2p}} \bigg(\int_s^t dr\, \big\| w_r \big\|_{\cb^{1+4\eps}_{p}}^{p}\bigg)^{\frac{1}{2p}}\\
&\lesssim \bigg(\int_s^t dr\, \big\| w_r \big\|_{L^{3p}}^{3p}\bigg)^{\frac{1}{p}}+ \bigg(\int_s^t dr\, \big\| w_r \big\|_{\cb^{1+4\eps}_{p}}^{p}\bigg)^{\frac{1}{p}},
\end{align*}
which concludes the proof of the lemma.
\end{proof}

\

\subsubsection{Consequence: a more general control on $w$}

\begin{lemma}\label{lemmun}
Fix $\frac43<\ga <\frac32$ {and $\eps >0$ small enough (depending on $\ga$). Then for all $p\geq 1$ large enough, there exists $\ka\geq 1$ such that for all $\mu$ satisfying \eqref{condiC-ter} and  all $0\leq s<t\leq T$},
 \begin{align}
\int_s^t dr\, \big\|w_r\big\|_{\B^\gamma_p}^p& \lesssim  \big\| w_s \big\|_{L^{3p-2}}^{3p-2}+ \big\| w_s \big\|_{\B^\ga_p}^{\frac{3p-2}{3}}+ {(\mu M)}^\ka \bigg[1+ \int_s^t dr\,   \big\|w_r\big\|^{3p}_{L^{3p}}\bigg].\label{lemmun-estim}
 \end{align}
\end{lemma}

\smallskip

\begin{proof}[Proof of Lemma \ref{lemmun}]

\

\smallskip

\noindent
\textit{Step 1: a first simplification.} By integrating \eqref{boutheo51st} between $s$ and $t$, we immediately obtain that for some $c_0,\ka_0>0$,
 \begin{align}
\int_s^t dr\, \big\|w_r\big\|_{\B^\gamma_p}^p& \leq c_0 \int_s^t dr\, \big\| e^{-(t-r) H}w_s \big\|_{\B^\gamma_p}^p+c_0 {(\mu M)}^{\ka_0} \bigg[1+ \int_s^t dr\,   \big\|w_r\big\|^{3p}_{L^{3p}}\bigg]+c_0 {(\mu M)}^{\ka_0} \int_s^t dr \,   \big\|w_r\big\|^{p}_{\B^{1+4\eps}_p}.\label{resu-integr}
 \end{align}
Then observe that for $\eps>0$ small enough and $\eps p\geq 1$ large enough, one has, for $\nu_0:= \frac{1+5\eps}{\ga+\eps}<1$,
\begin{align}
&\big\|w_r\big\|_{\B^{1+4\eps}_p} \lesssim \big\|w_r\big\|_{\B^{-\eps}_p}^{1-\nu_0} \big\|w_r\big\|_{\B^{\ga}_p}^{\nu_0}\lesssim \big\|w_r\big\|_{\B^{0}_{3p}}^{1-\nu_0} \big\|w_r\big\|_{\B^{\ga}_p}^{\nu_0}\lesssim \big\|w_r\big\|_{L^{3p}}^{1-\nu_0} \big\|w_r\big\|_{\B^{\ga}_p}^{\nu_0}\lesssim \big\|w_r\big\|_{L^{3p}}^{3} +\big\|w_r\big\|_{\B^{\ga}_p}^{\frac{3\nu_0}{2+\nu_0}},\label{interponu0}
\end{align}
which yields, for some $c_1>0$,
 \begin{align*}
c_0 {(\mu M)}^{\ka_0} \int_s^t dr \,   \big\|w_r\big\|^{p}_{\B^{1+4\eps}_p}& \leq c_1 {(\mu M)}^{\ka_0} \int_s^t dr \,   \big\|w_r\big\|_{L^{3p}}^{3p}+c_1 {(\mu M)}^{\ka_0} \int_s^t dr \,   \big\|w_r\big\|_{\B^{\ga}_p}^{\frac{3\nu_0 }{2+\nu_0}p}\\
&\leq c_1 {(\mu M)}^{\ka_0} \int_s^t dr \,   \big\|w_r\big\|_{L^{3p}}^{3p}+c_1 {(\mu M)}^{\ka_0} \bigg(\int_s^t dr \,   \big\|w_r\big\|_{\B^{\ga}_p}^{p}\bigg)^{\frac{3\nu_0 }{2+\nu_0}}\\
&\leq c_1 {(\mu M)}^{\ka_0} \int_s^t dr \,   \big\|w_r\big\|_{L^{3p}}^{3p}+ \frac{2-2\nu_0}{2+\nu_0} \big(c_1 {(\mu M)}^{\ka_0}\big)^{\frac{2+\nu_0}{2-2\nu_0}}+\frac{3\nu_0 }{2+\nu_0} \int_s^t dr \,   \big\|w_r\big\|_{\B^{\ga}_p}^{p}.
 \end{align*}
Going back to \eqref{resu-integr}, we obtain that for some constants $c_2,\ka_1>0$,
 \begin{align*}
\int_s^t dr\, \big\|w_r\big\|_{\B^\gamma_p}^p& \leq c_2 \int_s^t dr\, \big\| e^{-(t-r) H}w_s \big\|_{\B^\gamma_p}^p+c_2 {(\mu M)}^{\ka_1} \bigg[1+ \int_s^t dr\,   \big\|w_r\big\|^{3p}_{L^{3p}}\bigg]+\frac{3\nu_0 }{2+\nu_0} \int_s^t dr \,   \big\|w_r\big\|_{\B^{\ga}_p}^{p},
 \end{align*}
and since $\frac{3\nu_0 }{2+\nu_0}<1$, we deduce that
 \begin{align}
\int_s^t dr\, \big\|w_r\big\|_{\B^\gamma_p}^p& \lesssim \int_s^t dr\, \big\| e^{-(t-r) H}w_s \big\|_{\B^\gamma_p}^p+ {(\mu M)}^{\ka_1} \bigg[1+ \int_s^t dr\,   \big\|w_r\big\|^{3p}_{L^{3p}}\bigg].\label{stepun}
 \end{align}

\

\noindent
\textit{Step 2: initial-value term.} Let us set $\beta:=\ga-\frac{\ga+\eps}{p}$. With this parameter in hand, one has
\begin{align*}
&\int_s^t dr\, \big\| e^{-(t-r) H}w_s \big\|_{\B^\gamma_p}^p\lesssim \bigg(\int_s^t \frac{dr}{|t-r|^{\frac{\ga-\beta}{2}p}}\bigg) \big\| w_s \big\|_{\B^\beta_p}^p \lesssim \big\| w_s \big\|_{\B^\beta_p}^p ,
\end{align*}
due to $\frac{\ga-\beta}{2}p=\frac{\ga+\eps}{2}<1$. Then, similarly to \eqref{interponu0}, one has for $\nu_1:=\frac{\beta+\eps}{\ga+\eps}=1-\frac1p<1$,
\begin{align*}
&\big\| w_s \big\|_{\B^\beta_p} \lesssim \big\| w_s \big\|_{\B^{-\eps}_p}^{1-\nu_1} \big\| w_s \big\|_{\B^\ga_p}^{\nu_1}\lesssim \big\| w_s \big\|_{L^{3p-2}}^{1-\nu_1} \big\| w_s \big\|_{\B^\ga_p}^{\nu_1}\lesssim \big\| w_s \big\|_{L^{3p-2}}^{\frac{3p-2}{p}}+ \big\| w_s \big\|_{\B^\ga_p}^{\frac{3p-2}{(2+\nu_1)p-2} \nu_1}.
\end{align*}
As a result,
\begin{align*}
&\int_s^t dr\, \big\| e^{-(t-r) H}w_s \big\|_{\B^\gamma_p}^p\lesssim \big\| w_s \big\|_{L^{3p-2}}^{3p-2}+ \big\| w_s \big\|_{\B^\ga_p}^{\frac{3p-2}{(2+\nu_1)p-2} \nu_1 p} ,
\end{align*}
and we can finally check that
$$\frac{3p-2}{(2+\nu_1)p-2} \nu_1 p=\frac{3p-2}{3p-3}(p-1)=\frac{3p-2}{3}.$$
Going back to \eqref{stepun}, we deduce the desired bound \eqref{lemmun-estim}.
\end{proof}

\

\

\begin{lemma}\label{lemmdeux}
Fix $\frac43<\ga <\frac32$. For every small enough $\eps>0$, there exists $\beta=\beta_{\ga,\eps}\in (0,1)$ such that {for all $p\geq 1$ large enough, all $\mu$ satisfying \eqref{condiC-ter}} and all $0\leq s<t\leq T$, one has
 \begin{align*}
 &\big\|w_t\big\|^{3p-2}_{L^{3p-2}}+\int_s^t dr\,  \big\|w_r\|^{3p}_{L^{3p}} \lesssim  {(\mu M)}^\ka \bigg[ 1+ \big\|w_s\big\|^{3p-2}_{L^{3p-2}}+\int_s^tdr\,  \big\|w_r \big\|^{\beta p}_{\B^{\ga}_p}\bigg],
 \end{align*}
where $\ka\geq 1$ only depends on $\ga$, $\eps$ and $p$.
\end{lemma}

\smallskip

\begin{proof}[Proof of Lemma \ref{lemmdeux}]

We know by Theorem \ref{thm61}  and Theorem~\ref{thm411} that there exist constants $c_0,\ka_0>0$ for which
 \begin{align}
 &\big\|w_t\big\|^{3p-2}_{L^{3p-2}}+\int_s^t dr\,  \big\|w_r\big\|^{3p}_{L^{3p}} \leq c_0 \big\|w_s\big\|^{3p-2}_{L^{3p-2}}+c_0 {(\mu M)}^{\ka_0}+c_0 {(\mu M)}^{\ka_0}\int_s^tdr\,  \big\|w_r \big\|^p_{\B^{1+6\eps}_p}.\label{lemmdeux-1}
 \end{align}
Then, just as in \eqref{interponu0}, one has for $\nu:= \frac{1+7\eps}{\ga+\eps}<1$,
\begin{align*}
&\big\|w_r\big\|_{\B^{1+6\eps}_p} \lesssim \big\|w_r\big\|_{L^{3p}}^{1-\nu} \big\|w_r\big\|_{\B^{\ga}_p}^{\nu},
\end{align*}
and so, for some constant $c_1>0$,
$$c_0 {(\mu M)}^{\ka_0}\big\|w_r\big\|_{\B^{1+6\eps}_p}^p \leq \big\|w_r\big\|_{L^{3p}}^{(1-\nu)p}\Big( c_0 {(\mu M)}^{\ka_0} \big\|w_r\big\|_{\B^{\ga}_p}^{\nu p}\Big) \leq \frac{1-\nu}{3} \big\|w_r\big\|_{L^{3p}}^{3p}+\frac{2+\nu}{3} \Big( c_0 {(\mu M)}^{\ka_0} \big\|w_r\big\|_{\B^{\ga}_p}^{\nu p}\Big)^{\frac{3}{2+\nu}}.$$
Thus, setting $\beta:=\frac{3\nu}{2+\nu}<1$, there exist constants $c_1,\ka_1>0$ such that
$$c_0 {(\mu M)}^{\ka_0}\big\|w_r\big\|_{\B^{1+6\eps}_p}^p \leq \frac{1}{3} \big\|w_r\big\|_{L^{3p}}^{3p}+c_1 {(\mu M)}^{\ka_1} \big\|w_r\big\|_{\B^{\ga}_p}^{\beta p}.$$
Going back to \eqref{lemmdeux-1}, we deduce that
 \begin{align*}
 &\big\|w_t\big\|^{3p-2}_{L^{3p-2}}+\int_s^t dr\,  \big\|w_r\big\|^{3p}_{L^{3p}} \leq 3c_0 \big\|w_s\big\|^{3p-2}_{L^{3p-2}}+3c_0 {(\mu M)}^{\ka_0}+3c_1 {(\mu M)}^{\ka_1} \int_s^tdr\,  \big\|w_r\big\|_{\B^{\ga}_p}^{\beta p}.
 \end{align*}
which immediately entails the desired result.
\end{proof}

We now have all the tools in hand to prove our main control result:

\begin{theorem}\label{theo:impr-bou}
{Fix $\eps >0$ small enough. Then for all $p\geq 1$ large enough, there exists $\ka\geq 1$ such that for all $\mu$ satisfying \eqref{condiC-ter} and all $t\in [0,T]$,}
\begin{equation}\label{theo:impr-bou-1}
\big\|w_t\big\|_{\cb_p^{\frac32-\eps}} \lesssim \frac{{(\mu M)}^\ka}{t^{\frac32+\nu}}. 
\end{equation}
\end{theorem}

\begin{proof}
Fix $\frac43<\ga<\frac32$ (eventually taking $\gamma=\tfrac32-\eps$). By Lemma \ref{lemmdeux}, one has
 \begin{align*}
 &\int_s^t dr\,  \big\|w_r\big\|^{3p}_{L^{3p}} \lesssim  {(\mu M)}^{\ka_0} \bigg[ 1+ \big\|w_s\big\|^{3p-2}_{L^{3p-2}}+\int_s^tdr\,  \big\|w_r \big\|^{\beta p}_{\B^{\ga}_p}\bigg],
 \end{align*}
which, combined with Lemma \ref{lemmun}, gives
 \begin{align*}
\int_s^t dr\, \big\|w_r\big\|_{\B^\gamma_p}^p& \lesssim  \big\| w_s \big\|_{L^{3p-2}}^{3p-2}+ \big\| w_s \big\|_{\B^\ga_p}^{\frac{3p-2}{3}}+ {(\mu M)}^{\ka_1} \bigg[1+ \int_s^t dr\,   \big\|w_r\big\|^{3p}_{L^{3p}} \bigg]\\
&\lesssim {(\mu M)}^{\ka_2} \Big[1+\big\| w_s \big\|_{L^{3p-2}}^{3p-2}+ \big\| w_s \big\|_{\B^\ga_p}^{\frac{3p-2}{3}}\Big]+ {(\mu M)}^{\ka_2} \int_s^tdr\,  \big\|w_r \big\|^{\beta p}_{\B^{\ga}_p}.
 \end{align*}
By summing these two estimates, we deduce the existence of constants $c_0,\ka_3>0$ such that
 \begin{align*}
 \int_s^t dr\,  \big\|w_r\big\|^{3p}_{L^{3p}}+\int_s^t dr\, \big\|w_r\big\|_{\B^\gamma_p}^p&\leq c_0 {(\mu M)}^{\ka_3} \Big[1+\big\| w_s \big\|_{L^{3p-2}}^{3p-2}+ \big\| w_s \big\|_{\B^\ga_p}^{\frac{3p-2}{3}}\Big]+ c_0 {(\mu M)}^{\ka_3} \int_s^tdr\,  \big\|w_r \big\|^{\beta p}_{\B^{\ga}_p}.
 \end{align*}
We now recall the following elementary form of Young's inequality: given $\beta\in (0,1)$, there exist $a\geq 0$ and $A>0$ such that for all $C\geq 0$ and $x\geq 0$, $C x^\beta \leq A  C^a+\frac12 x$. As a result, we can write here
$$ c_0 {(\mu M)}^{\ka_3} \int_s^tdr\,  \big\|w_r \big\|^{\beta p}_{\B^{\ga}_p} \leq c_1 {(\mu M)}^{\ka_4}+\frac12 \int_s^tdr\,  \big\|w_r \big\|^{ p}_{\B^{\ga}_p},$$
and we have thus shown that for some $\ka_5>0$,
 \begin{align}
 \int_s^t dr\,  \big\|w_r\big\|^{3p}_{L^{3p}}+\int_s^t dr\, \big\|w_r\big\|_{\B^\gamma_p}^p&\lesssim  {(\mu M)}^{\ka_5} \Big[1+\big\| w_s \big\|_{L^{3p-2}}^{3p-2}+ \big\| w_s \big\|_{\B^\ga_p}^{\frac{3p-2}{3}}\Big] \nonumber \\
 &\lesssim  {(\mu M)}^{\ka_5} \Big[1+ \big\|w_s\big\|^{3\sigma}_{L^{3p}}+ \big\|w_s\big\|^{\sigma}_{\B^{\gamma}_{p}}\Big], \label{shown}
 \end{align}
where in the last line we have used \eqref{bsigma}, and $\sigma$ is fixed as in \eqref{def-sigma}.

\medskip

In order to extend this control to every time $t\in [0,T]$, recall that by Theorem \ref{improv-thm51}, one has for all $0\leq s<t\leq T$
 \begin{align*}
\big\|w_t\big\|_{\B^\gamma_p}& \lesssim \big\| w_s \big\|_{\B^\gamma_p}+ {(\mu M)}^{\ka_6} \bigg[1+   \bigg(\int_s^t dr\,   \big\|w_r\big\|^{3p}_{L^{3p}}\bigg)^{\frac1p}+  \bigg(\int_s^tdr\,   \big\|w_r\big\|^{p}_{\B^{1+4\eps}_p}\bigg)^{\frac1p}\bigg], 
 \end{align*}
which, combined with \eqref{shown}, gives (we recall that $\sigma/p<1$)
 \begin{align*}
\big\|w_t\big\|_{\B^\gamma_p}& \lesssim {(\mu M)}^{\ka_7} \Big[1+\big\| w_s \big\|_{\B^\ga_p}+\big\| w_s \big\|_{L^{3p}}^{\frac{3\si}{p}}+ \big\| w_s \big\|_{\B^\ga_p}^{\frac{\si}{p}}\Big] \lesssim {(\mu M)}^{\ka_7} \Big[1+\big\| w_s \big\|_{L^{3p}}^3+ \big\| w_s \big\|_{\B^\ga_p}\Big],
 \end{align*}
 so that 
  \begin{align*}
\big\|w_t\big\|^\sigma_{\B^\gamma_p}&  \lesssim {(\mu M)}^{\ka_8} \Big[1+\big\| w_s \big\|_{L^{3p}}^{3\sigma}+ \big\| w_s \big\|^{\sigma}_{\B^\ga_p}\Big].
 \end{align*}

To conclude, we apply the same strategy as in the proof of Theorem \ref{thm411}. Namely we set
$$F(t):=   \big\|w_t\big\|^{3\sigma}_{L^{3p}}+\big\|w_t\big\|^{\sigma}_{\B^{\gamma}_{p}} , \quad\quad  L(t):=\big\| w_t \big\|^{\sigma}_{\B_p^{\gamma}},$$
and introduce the stopping time $\tau$ analogous to the one in \eqref{def:stop-time-tau}.

\smallskip

Then, by Lemma \ref{lem:pointwise}, one has on $[0,\tau]$,
 \begin{equation*} 
\big\| w_t \big\|_{\B_p^{\gamma}} \lesssim {(\mu M)}^{\kappa_9} t^{-\frac{1}{p-\sigma}}\lesssim {(\mu M)}^{\kappa_9} t^{-\frac32-\nu},
 \end{equation*}
where we used that  $\frac{1}{p-\si}\geq \frac32$ and $\frac{1}{p-\si} \to \frac32$ as $p\to\infty$. This concludes the proof. 
\end{proof}

\

\begin{proof}[Proof of Theorem \ref{thm:global_Lp} (for $s=0$)]

Fix $\eps>0$ small enough. By taking $p\geq 1$ large enough, fixing $\mu:= \max\big(1+\big(8M \big)^8,\mu_0 M^{30p}\big)$ (in accordance with condition \eqref{condiC-ter}) and using Lemma \ref{lem:inclusion-besov}, we obtain from Theorem \ref{thm411} that uniformly over $u \in \B_{\infty}^{-\frac{1}{2}-\eps}$,
\begin{equation}\label{thm411-1-bis}
\sup_{t\in [0,T]}\Big( \big\| v_t \big\|_{\B^{-\eps}_{\infty}}+ t^{\frac12+\nu}\big\| w_t \big\|_{\cb_{\infty}^{-\eps}} \Big)\lesssim M^\kappa.
\end{equation}
On the other hand, arguing as in \eqref{bounve-2}, we derive for every $t\in [0,T]$
\begin{align*}
\big\|v_t\big\|_{\cb^{\frac34-2\eps}_\infty} &\lesssim \int_0^t ds \, \Big\| e^{-(t-s)H}\Big(\big(v_s+w_s -\<IPsi3>_{s}\big) \pl  \<Psi2>_{s}\Big) \Big\|_{\cb^{\frac34-2\eps}_\infty}\\
&\lesssim \int_0^t \frac{ds}{(t-s)^{\frac78}} \, \Big\|\big(v_s+w_s -\<IPsi3>_{s}\big) \pl  \<Psi2>_{s} \Big\|_{\cb^{-1-2\eps}_\infty}\\
&\lesssim M\int_0^t \frac{ds}{(t-s)^{\frac78}} \, \big\|v_s+w_s -\<IPsi3>_{s}\big\|_{\cb^{-\eps}_\infty}\lesssim M^\ka\int_0^t \frac{ds}{(t-s)^{\frac78}} \frac{1}{s^{\frac12+\nu}}\lesssim \frac{M^\ka}{t^{\frac38+\nu}} ,
\end{align*}
where we have used \eqref{thm411-1-bis} to control $\big\|v_s+w_s\big\|_{\cb^{-\eps}_\infty}$. As a result,
\begin{align*}
\sup_{t\in [0,1]} t^{\frac38+\nu}\big\|v_t\big\|_{\cb^{\frac34-2\eps}_\infty} &\lesssim M^\ka,
\end{align*}
which, combined with \eqref{theo:impr-bou-1} and Lemma \ref{lem:inclusion-besov}, concludes the proof of Theorem \ref{thm:global_Lp}.
\end{proof}


\section{Markov property and existence of an invariant measure} \label{Sect5}

 In this section, we exploit the a priori bounds established in the preceding section to prove the existence of an invariant measure for the dynamics~\eqref{eq:intro}. We also show that any such invariant measure cannot be Gaussian.

\smallskip

With the representation \eqref{repres-noise} of the noise in mind, we recall that the notation $(\mathcal{F}_t)_{t \geq 0}$ stands for the natural filtration generated by the family $(\beta^{(k)})_{k\geq 0}$. Now, given the unique solution $X_t=X^u_t$ of \eqref{eq:intro} obtained in Corollary \ref{bound_final}, let us denote by $(P_t)_{t\ge 0}$ the family of operators defined for all $\Phi \in B_b(\B_{\infty}^{-\frac12-\eps};\R)$ by
\begin{equation}\label{transi}
P_t \Phi (u): = \mathbb{E} \big[ \Phi \big(X_t^u\big)\big], \quad u \in \B_{\infty}^{-\frac12-\eps}, \quad t\ge 0,
\end{equation}
where $B_b(\B_{\infty}^{-\frac12-\eps};\R)$ is the space of bounded Borel functions $\Phi: \B_{\infty}^{-\frac12-\eps} \to \R$.

\subsection{Markov property}

We first ensure that the Markov property is satisfied at the level of the approximate solution ${X^{(n),u}}$, that is, for the solution of \eqref{appro-equ-introd}. For the sake of clarity, we denote by $\mathcal{J} (u , t_0; \frakc^{(n)},z)_t$ the solution of \eqref{mild} at time $t$ with initial data $u$ at $t_0$ and (regular) driver $z$, {\it i.e.}, $f:=\mathcal{J} (u, t_0;\frakc^{(n)}, z)$ satisfies
\begin{equation} \label{mild}
f_t = e^{-(t-t_0)H}u -\int_{t_0}^t ds\, e^{-(t-s)H} (f_s)^3   +\int_{t_0}^tds\,  e^{-(t-s)H} \frakc^{(n)} f_s  + z_t. 
\end{equation}

\begin{lemma}
The approximation ${X^{(n),u}}$ satisfies, for any $\Phi \in \mathcal{C}_b(\B_{\infty}^{-\frac12-\eps})$,
\begin{equation}\label{lemmaMarkov}
 \mathbb{E} \Big[\Phi\big({X^{(n),u}_{t+h}}\big) | \mathcal{F}_t\Big] = \big(P^{(n)}_h\Phi\big) \big({X^{(n),u}_{t}}\big),
 \end{equation} 
where  
$$\big(P^{(n)}_h\Phi\big) (u):=\mathbb{E} \Big[ \Phi \big({X^{(n),u}_{h}}\big)\Big].$$
\end{lemma}

\begin{proof}
Recall that $\<Psi>^{(n)}_{t_1,t_2}$ is defined (in \eqref{def-luxost}) as
\begin{equation*} 
\<Psi>^{(n)}_{t_1,t_2}(x)=\int_{t_1}^{t_2} \big(e^{-(t_2-s)H}dW^{(n)}_s\big)(x)=\int_{t_1}^{t_2} \int_{\R^3} K_{t_2-s}(x,y)\, W^{(n)}(ds,dy). 
\end{equation*}
The approximating equation 
$$(\partial_t +H) {X^{(n),u}} =-\big({X^{(n),u}}\big)^3 +\frakc^{(n)} {X^{(n),u}}+\xi^{(n)}$$
with initial data $u\in \B_{\infty}^{-\frac12-\eps}$ can be written in the mild form 
\begin{equation*}  
{X^{(n),u}_t}= e^{-tH}u -\int_0^t ds\, e^{-(t-s)H} ({X^{(n),u}_s})^3   +\int_0^tds\,  e^{-(t-s)H} \frakc^{(n)} {X^{(n),u}_s}  + \<Psi>^{(n)}_{0,t}. 
\end{equation*}
For $h>0$, one has
\begin{multline*}
{X^{(n),u}_{t+h}}=e^{-h H} \Big\{e^{-tH}u -\int_0^t ds\, e^{-(t-s)H} ({X^{(n),u}_s})^3   +\int_0^t ds\, e^{-(t-s)H} \frakc^{(n)} {X^{(n),u}_s} \Big\} \\
 -\int_{t}^{t+h} ds\,e^{-(t+h-s)H} ({X^{(n),u}_s})^3  +\int_{t}^{t+h} ds\,e^{-(t+h-s)H} \frakc^{(n)} {X^{(n),u}_s}  + \<Psi>^{(n)}_{0,t+h}.
\end{multline*}
Note that $\<Psi>^{(n)}_{0,t+h}= e^{-h H} \<Psi>^{(n)}_{0,t} + \<Psi>^{(n)}_{t,t+h}$. Therefore we have 
\begin{align*}
{X^{(n),u}_{t+h}} &=e^{-h H} {X^{(n),u}_t} -\int_{t}^{t+h} ds\, e^{-(t+h-s)H} ({X^{(n),u}_s})^3  +\int_{t}^{t+h} ds\,e^{-(t+h-s)H} \frakc^{(n)} {X^{(n),u}_s} + \<Psi>^{(n)}_{t,t+h}\\
&=e^{-h H} {X^{(n),u}_t} -\int_{0}^{h} ds\, e^{-(h-s)H} ({X^{(n),u}_{t+s}})^3  +\int_{0}^{h} ds\,e^{-(h-s)H} \frakc^{(n)} {X^{(n),u}_{t+s}}  + \<Psi>^{(n)}_{t,t+h}
\end{align*}
or in other words, by uniqueness of the solution,
$${X^{(n),u}_{t+h}}=\cj\big({X^{(n),u}_t},0;\frakc^{(n)}, \<Psi>^{(n)}_{t,t+.}\big)_h.$$
Since ${X^{(n),u}_t}$ is $\mathcal{F}_t$ measurable and  $\<Psi>^{(n)}_{t,t+.}$ is independent of $\mathcal{F}_t$, we deduce that
$$
\mathbb{E} \big[\Phi({X^{(n),u}_{t+h}}) | \mathcal{F}_t\big] =\mathbb{E} \Big[\Phi\big(\cj\big({X^{(n),u}_t},0;\frakc^{(n)}, \<Psi>^{(n)}_{t,t+.}\big)_{h}\big) \big| \mathcal{F}_t\Big]=  \big(P^{(n)}_h\Phi\big) \big({X^{(n),u}_t}\big),
$$
with  
$$\big(P^{(n)}_h\Phi\big) (y)=\mathbb{E} \Big[\Phi \big(\mathcal{J}\big(y , 0; \frakc^{(n)},\<Psi>^{(n)}_{t,t+.}\big)_h\big)\Big]= \mathbb{E} \Big[\Phi \big(\mathcal{J}\big(y , 0;\frakc^{(n)}, \<Psi>^{(n)}_{0,.}\big)_{h}\big)\Big],$$
since $\<Psi>^{(n)}_{t,t+.}$ and $\<Psi>^{(n)}_{0,.}$ have the same distribution.
\end{proof}

\medskip

 We now pass to the limit $n\to \infty$ in the identity \eqref{lemmaMarkov}, which is justified by the almost sure convergence of ${X^{(n),u}_t}$ to ${X^{u}_t}$ (Proposition~\ref{prop:fixed-point}~$(ii)$), and we deduce the desired Markov identity: for every $\Phi \in \mathcal{C}_b(\B_{\infty}^{-\frac12-\eps})$,
\begin{equation*}
 \mathbb{E} \Big[\Phi\big({X^{u}_{t+h}}\big) | \mathcal{F}_t\Big] = \big(P_h\Phi\big) \big({X^{u}_t}\big),
 \end{equation*}

Observe also that the transition semigroup $(P_t)_{t\geq 0}$ satisfies the Feller property, as a result of the continuous dependence of the solutions on the initial data (see Proposition~\ref{prop:fixed-point}~$(ii)$). In other words, denoting by $\mathcal{C}_b(\B_{\infty}^{-\frac12-\eps})$ the space of continuous bounded functions on $\B_{\infty}^{-\frac12-\eps}$, one has $P_t \Phi \in \mathcal{C}_b(\B_{\infty}^{-\frac12-\eps})$ for all $\Phi \in \mathcal{C}_b(\B_{\infty}^{-\frac12-\eps})$ and $t\ge 0$.

\

\begin{proposition}
The family $(P_t)_{t\geq 0}$ introduced in \eqref{transi} defines a Feller semigroup.
\end{proposition}

\
 
\subsection{Existence of an invariant measure}

Thanks to the bounds established in the previous sections, we can now state our existence result for the invariant measure of the dynamics.

\begin{theorem} \label{ThmMeasure}
The semigroup $(P_t)_{t \geq 0} $ admits at least one invariant measure. 
\end{theorem}

\begin{proof} 
The proof relies on the Krylov-Bogoliubov method (see \cite[Section~11.2]{DPZ}; see also \cite[Section 4.2]{TW18} for a similar application in a related setting).

\smallskip

For any $A>0$ and $t\ge 1$, Corollary \ref{bound_final} yields, for every $u \in \B_{\infty}^{-\frac12 -\eps}$,
\begin{align*}
 \int_0^t  ds \,\mathbb{P}\Big( \big\|{X^{u}_s}\big\|_{\B_{\infty}^{-\frac12 -\eps}} \ge A\Big)  & \le  \frac{1}{A} \int_0^t  ds\, 
\mathbb{E} \Big[ \big\|{X^{u}_s}\big\|_{\B_{\infty}^{-\frac12 -\eps}}\Big]  \le \frac{C}{A} \bigg(\int_0^1 \frac{ds}{s^{\frac34}}  +\int_1^t ds\bigg) \le \frac{C t}{A},
\end{align*}
with $C$ independent of $t$. 
Consequently, for any $\eta>0$, taking $A=A_{\eta}=\frac{C}{\eta}$, we obtain
$$ \frac{1}{t} \int_0^t  ds\, \mathbb{P}\Big( \big\|{X^{u}_s}\big\|_{\B_{\infty}^{-\frac12 -\eps}} \ge A_\eta\Big)   \le  \eta. $$
The embedding $ \B_{\infty}^{-\frac12-\eps}\subset \B_{\infty}^{-\frac12-2\eps}$ is compact (see \cite[Proposition 1.1]{FI}), 
and thus setting 
$$\dis R_t=\frac{1}{t} \int_0^t ds \, P^{\ast}_s ,$$
 the sequence   $\dis \big\{R_t \delta_{u}\big\}_{t\ge 0}$ is tight in $\B_{\infty}^{-\frac12-2\eps}$. 
 The Krylov-Bogoliubov argument then implies that there exist a sequence $t_k \to \infty$ and a measure $\rho$ toward which~$R_{t_k}$ converges weakly as  $k \to \infty$. Moreover, $\rho$ is invariant for the transition semigroup $(P_t)_{t \geq 0} $.
\end{proof}

\

\subsection{Non-Gaussianity of the invariant measure}\label{subsec:non-gauss}

\begin{theorem}\label{theo:non-gauss}
If $\rho$ is an invariant measure for the dynamics, then it cannot be Gaussian. 
\end{theorem}

 Our argument relies on the consideration of the fourth-order cumulant, a classical tool to detect non-Gaussianity (see {\it e.g.}  \cite{NP} for a general overview of this approach in the context of Malliavin calculus, and~\cite{nualart-book} for background material).

\subsubsection{Reduction of the problem}
We introduce the family of operators $(\cl^{(r)})_{0<r\leq 1}$ defined by
$$\cl^{(r)}:=\big(1-e^{-rH}\big) e^{-rH}.$$
Then, for all (distributional-valued) random variables $Z_1,\ldots,Z_4$, set
\begin{align*}
&\ka^{(r)}\big(Z_1,\ldots,Z_4 \big):= \mathbb{E}\Big[ \big(\cl^{(r)}Z_1\big)(0)\big(\cl^{(r)}Z_2\big)(0)\big(\cl^{(r)}Z_3\big)(0)\big(\cl^{(r)}Z_4\big)(0)\Big]\\
&\hspace{1cm}-\mathbb{E}\big[ \big(\cl^{(r)}Z_1\big)(0)\big(\cl^{(r)}Z_2\big)(0)\big]\mathbb{E}\big[ \big(\cl^{(r)}Z_3\big)(0)\big(\cl^{(r)}Z_4\big)(0)\big]\\
&\hspace{2cm}-\mathbb{E}\big[ \big(\cl^{(r)}Z_1\big)(0)\big(\cl^{(r)}Z_3\big)(0)\big]\mathbb{E}\big[ \big(\cl^{(r)}Z_2\big)(0)\big(\cl^{(r)}Z_4\big)(0)\big]\\
&\hspace{3cm}-\mathbb{E}\big[ \big(\cl^{(r)}Z_1\big)(0)\big(\cl^{(r)}Z_4\big)(0)\big]\mathbb{E}\big[ \big(\cl^{(r)}Z_2\big)(0)\big(\cl^{(r)}Z_3\big)(0)\big].
\end{align*}

Consider an invariant measure $\rho$, so that for all $t\geq 0$, $X_t\sim \rho$. We proceed by contradiction and suppose that $X_1$ is Gaussian. In this case, the linear transformation $\cl^{(r)} X_1 (0) $ is also Gaussian and consequently we should have
\begin{equation}\label{ident-ka-r}
\ka^{(r)}\big(X_1, X_1, X_1, X_1 \big) =0.
\end{equation}
In fact, we will prove that for all $0<r\leq 1$ and $\eps>0$,
\begin{equation}\label{decomp-G}
\ka^{(r)}\big(X_1,X_1,X_1,X_1 \big) = { -4 \cm^{(r)}}+ \mathcal{O}\Big(\frac{1}{r^{\frac38+\eps}}\Big)
\end{equation}
with 
\begin{equation*}
 { \cm^{(r)}\geq \frac{c}{r^{\frac12}}},
\end{equation*}
for some strictly positive constant $c>0$, independent of $r$. As $r\to 0$, this contradicts~\eqref{ident-ka-r}, and hence $\rho$ cannot be Gaussian.  \medskip

To prove \eqref{decomp-G}, we will rely on the decomposition (see \eqref{decompos-x}, and set $\tau_t:=\tau_{0,t}$ for doubly-indexed objects)
$$X_1=\<Psi>_{1}-\<IPsi3>_{1}+U_1.$$
Substituting this decomposition into $\ka^{(r)}(X_1,X_1,X_1,X_1)$, we obtain, for some remainder $R^{(r)}$, 
\begin{align*}
&\ka^{(r)}\big(X_1,X_1,X_1,X_1 \big)=\ka^{(r)}\big(\<Psi>_1,\<Psi>_1,\<Psi>_1,\<Psi>_1 \big)-4 \ka^{(r)}\Big(\<IPsi3>_1,\<Psi>_1,\<Psi>_1,\<Psi>_1 \Big)+R^{(r)}.
\end{align*}
Since $\<Psi>_1$ is Gaussian, one has first $\ka^{(r)}(\<Psi>_1,\<Psi>_1,\<Psi>_1,\<Psi>_1) =0$. Then, as far as $R^{(r)}$ is concerned, recall that by Proposition \ref{prop:conv-arbre} and Theorem \ref{thm:global_Lp}, one has
\begin{equation*} 
\mathbb{E}\Big[ \big\| \<Psi>_1\big\|_{\cb_{\infty}^{-\frac12-\eps}}^p\Big]+ \mathbb{E}\Big[ \big\| \<IPsi3>_1\big\|_{\cb_{\infty}^{\frac12-\eps}}^p\Big]+\mathbb{E}\Big[ \big\| U_1\big\|_{\cb_{\infty}^{\frac34-\eps}}^p\Big] < \infty,
\end{equation*}
which yields for instance
\begin{align*}
\big|\ka^{(r)}\big(U_1,\<Psi>_1,\<Psi>_1,\<Psi>_1\big)\big|&\lesssim \mathbb{E}\Big[ \big\|\cl^{(r)} U_1\big\|_{L^\infty}^4\Big]^{\frac14} \bigg( \mathbb{E}\Big[ \big\|\cl^{(r)} \<Psi>_1\big\|_{L^\infty}^4\Big]^{\frac14} \bigg)^3\\
&\lesssim \bigg( r^{\frac38-\frac{\eps}{2}}\mathbb{E}\Big[ \big\| U_1\big\|_{\cb_{\infty}^{\frac34-\eps}}^4\Big]^{\frac14}\bigg) \bigg( \frac{1}{r^{\frac14+\frac{\eps}{2}}} \mathbb{E}\Big[ \big\| \<Psi>_1\big\|_{\cb_{\infty}^{-\frac12-\eps}}^4\Big]^{\frac14} \bigg)^3\lesssim \frac{1}{r^{\frac38+2\eps}}.
\end{align*}
The other terms contained in $R^{(r)}$ can be estimated in the same way, and we thus conclude that
\begin{align*}
&\ka^{(r)}\big(X_1,X_1,X_1,X_1 \big)=-4 \ka^{(r)}\Big(\<IPsi3>_1,\<Psi>_1,\<Psi>_1,\<Psi>_1 \Big)+\mathcal{O}\Big(\frac{1}{r^{\frac38+\eps}}\Big).
\end{align*} 
Once endowed with this decomposition, the desired expansion \eqref{decomp-G} immediately follows from the combination of Lemmas \ref{lemm-cm} and \ref{lemm-cm-2} below.

\medskip
  
\

 \subsubsection{Technical estimates} Recall that we denote by $A\gtrsim B$ any inequality of the form $A\geq cB$, where $c>0$ is a universal constant. In addition, we set 
$$K_{t,s}:=K_t-K_s,$$
so that in particular, for any function $f:\R^3\to \R$,
$$\big(\cl^{(r)}f\big)(x)=\int dy \, K_{r,2r}(x,y) f(y).$$
Finally, we denote by $(I_k)_{k\geq 0}$ the multiple integrals driven by $W$, as defined in \cite[Section~1.1.2]{nualart-book}. Recall in particular (see \cite[Section 4.1]{DFT}) that
\begin{equation}\label{represen-psi-0}
\<Psi>_t(x)=I_1\big(F_{t,x}\big)  \quad \text{with} \ \ F_{t,x}(s,w):=\1_{[0,t]}(s) K_{t-s}(x,w), \quad \<Psi3>_t(x)=I_3\big((F_{t,x})^{\otimes 3}\big),
\end{equation}
and for $t_1\neq t_2$
\begin{equation}\label{def-E}
\mathbb{E}\Big[ \<Psi>_{t_1}(y_1) \<Psi>_{t_2}(y_2)\Big]=\cac_{t_1,t_2}(y_1,y_2)=\frac12\int^{t_1+t_2}_{|t_2-t_1|}d\si\,  K_{\si}(y_1,y_2).
\end{equation}

\

\begin{lemma}\label{lemm-cm}
For all $0<r\leq 1$, one has
\begin{equation}\label{calcul-K1}
\ka^{(r)}\Big(\<IPsi3>_1,\<Psi>_1,\<Psi>_1,\<Psi>_1 \Big) =c\, \cm^{(r)}+\mathcal{O}\Big(\frac{1}{r^{\frac14}} \Big),
\end{equation} 
for some constant $c>0$, and where we have set
\begin{equation*}
\cm^{(r)}:=\int_{r}^{+\infty}ds\int_{s}^{s+r} d\si_1 \int_{s}^{s+r} d\si_2 \int_{s}^{s+r} d\si_3 \int dz\, K_{s,r+s}(0,z)  K_{\si_1}(0,z) K_{\si_2}(0,z)K_{\si_3}(0,z).
\end{equation*}
\end{lemma}

\begin{proof}

Observe first that for every $r\geq 0$,
\begin{equation}\label{prem-obs}
\mathbb{E}\Big[ \Big(\cl^{(r)}\<IPsi3>_1\Big)(0)\big(\cl^{(r)}\<Psi>_1\big)(0)\Big]=0.
\end{equation}
Indeed, using the representations in \eqref{represen-psi-0} and the multiplication rules in Wiener chaoses, we have
\begin{align*}
&\Big(\cl^{(r)}\<IPsi3>_1\Big)(0)\big(\cl^{(r)}\<Psi>_1\big)(0)=\\
&=\int dy_1  dy_2 \, K_{r,2r}(0,y_1)K_{r,2r}(0,y_2) \int_{0}^1 ds \int dz \, K_{1-s}(y_1,z)  I_3\big(F_{s,z}^{\otimes 3}\big)I_1\big(F_{1,y_2}\big)\\
&=\int dy_1  dy_2 \, K_{r,2r}(0,y_1)K_{r,2r}(0,y_2) \\
&\hspace{3cm}\int_{0}^1 ds \int dz \, K_{1-s}(y_1,z) \big[I_4\big(F_{s,z}^{\otimes 3}\otimes F_{1,y_2}\big)+3\,\cac_{s,1}(y_2,z)I_2\big(F_{s,z}\otimes F_{1,y_2}\big)\big],
\end{align*}
from which \eqref{prem-obs} immediately follows.

\smallskip

As a result of \eqref{prem-obs}, the quantity $\ka^{(r)}\big(\<IPsi3>_1,\<Psi>_1,\<Psi>_1,\<Psi>_1 \big)$ under consideration reduces to
\begin{multline}\label{calculM}
\ka^{(r)}\Big(\<IPsi3>_1,\<Psi>_1,\<Psi>_1,\<Psi>_1 \Big) =\mathbb{E}\Big[ \big(\cl^{(r)}\<IPsi3>_1\big)(0)\big(\cl^{(r)}\<Psi>_1\big)(0)\big(\cl^{(r)}\<Psi>_1\big)(0)\big(\cl^{(r)}\<Psi>_1\big)(0)\Big]= \\
=\int dy_1 \cdots dy_4 \, K_{r,2r}(0,y_1)K_{r,2r}(0,y_2)K_{r,2r}(0,y_3)K_{r,2r}(0,y_4)\mathbb{E}\Big[ \<IPsi3>_1(y_1)\<Psi>_1(y_2)\<Psi>_1(y_3)\<Psi>_1(y_4)\Big].
\end{multline}
Then 
\begin{equation*}
 \<IPsi3>_1(y_1)\<Psi>_1(y_2)\<Psi>_1(y_3)\<Psi>_1(y_4)
=\int_{0}^1 ds \int dz \, K_{1-s}(y_1,z)  I_3\big(F_{s,z}^{\otimes 3}\big)I_1\big(F_{1,y_2}\big)I_1\big(F_{1,y_3}\big)I_1\big(F_{1,y_4}\big),
\end{equation*}
and
\begin{multline*}
 \mathbb{E} \Big[ I_3\big(F_{s,z}^{\otimes 3}\big)I_1\big(F_{1,y_2}\big)I_1\big(F_{1,y_3}\big)I_1\big(F_{1,y_4}\big)\Big]=\\
\begin{aligned}
&=\mathbb{E} \Big[\Big( I_4\big(F_{s,z}^{\otimes 3} \otimes F_{1,y_2}\big)+c_1 \cac_{s,1}(z,y_2)I_2\big(F_{s,z}^{\otimes 2}\big)\Big) \Big(I_2\big(F_{1,y_3} \otimes F_{1,y_4}\big)+c_2 \cac_{1,1}(y_3,y_4)\Big)\Big]\\
&=c_1\, \cac_{s,1}(z,y_2) \mathbb{E} \Big[I_2\big(F_{s,z}^{\otimes 2}\big)I_2\big(F_{1,y_3} \otimes F_{1,y_4}\big)\Big]\\
&=c\, \cac_{s,1}(z,y_2)\cac_{s,1}(z,y_3)\cac_{s,1}(z,y_4),
\end{aligned}
\end{multline*}
for some combinatorial coefficient $c>0$. Therefore, going back to \eqref{calculM}, we obtain
\begin{align*}
&\mathbb{E}\Big[ \big(\cl^{(r)}\<IPsi3>_1\big)(0)\big(\cl^{(r)}\<Psi>_1\big)(0)\big(\cl^{(r)}\<Psi>_1\big)(0)\big(\cl^{(r)}\<Psi>_1\big)(0)\Big] =\\
&=c\int dy_1 \cdots dy_4 \, K_{r,2r}(0,y_1)K_{r,2r}(0,y_2)K_{r,2r}(0,y_3)K_{r,2r}(0,y_4) \\
&\hspace{2cm}\times \int_{0}^1 ds \int dz \, K_{1-s}(y_1,z)\cac_{s,1}(z,y_2)\cac_{s,1}(z,y_3)\cac_{s,1}(z,y_4)\\
&=c\int_{0}^1 ds\int dz\, \bigg(\int dy_1\, K_{r,2r}(0,y_1)K_{1-s}(y_1,z)\bigg) \bigg(\int_{1-s}^{1+s}d\si\bigg(\int dy_2\, K_{r,2r}(0,y_2)K_{\si}(z,y_2)\bigg)\bigg)^3 \\
&=c\int_{0}^1 ds \int dz\, K_{r+(1-s),2r+(1-s)}(0,z)  \bigg(\int_{1-s}^{1+s}d\si\,  K_{r+\si,2r+\si}(0,z)\bigg)^3\\
&=c\int_{0}^{1}ds \int dz\, K_{r+s,2r+s}(0,z) \bigg(\int_{s}^{2-s}d\si \, K_{r+\si,2r+\si}(0,z) \bigg)^3.
\end{align*}
Using the change of variables $(s',\sigma')=(s+r,\sigma+r)$, we obtain
\begin{multline*}
\mathbb{E}\Big[ \big(\cl^{(r)}\<IPsi3>_1\big)(0)\big(\cl^{(r)}\<Psi>_1\big)(0)\big(\cl^{(r)}\<Psi>_1\big)(0)\big(\cl^{(r)}\<Psi>_1\big)(0)\Big] =\\
=c\int_{r}^{1+r}ds  \int dz\, K_{s,r+s}(0,z) \bigg(\int_{s}^{2-s+2r}d\si\,  K_{\si,r+\si}(0,z)\bigg)^3 .
\end{multline*}
Let us set
$$\mathcal{R}^{(r)}:=\cm^{(r)}-\int_{r}^{1+r}ds  \int dz\, K_{s,r+s}(0,z) \bigg(\int_{s}^{2-s+2r}d\si\,  K_{\si,r+\si}(0,z)\bigg)^3.$$
By writing $\cm^{(r)}$ as
\begin{align*}
\cm^{(r)}&=\int_{r}^{+\infty}ds \int dz\, K_{s,r+s}(0,z) \bigg(\int_{s}^{s+r} d\si\,  K_{\si}(0,z)\bigg)^3\\
&=\int_{r}^{+\infty}ds \int dz\, K_{s,r+s}(0,z) \bigg(\int_{s}^{+\infty}d\si\, K_{\si,r+\si}(0,z)\bigg)^3,
\end{align*}
we deduce that
\begin{align*}
&\mathcal{R}^{(r)}=\int_{r}^{+\infty}ds \int dz\, K_{s,r+s}(0,z) \bigg(\int_{s}^{+\infty}d\si\, K_{\si,r+\si}(0,z)\bigg)^3\\
&\hspace{3cm}-\int_{r}^{1+r}ds  \int dz\, K_{s,r+s}(0,z) \bigg(\int_{s}^{2-s+2r}d\si\,  K_{\si,r+\si}(0,z)\bigg)^3\\
&=\int_{1+r}^{+\infty}ds \int dz\, K_{s,r+s}(0,z) \bigg(\int_{s}^{+\infty}d\si\, K_{\si,r+\si}(0,z)\bigg)^3\\
&\hspace{1cm}+\int_{r}^{1+r}ds\int dz\, K_{s,r+s}(0,z)\bigg[ \bigg(\int_{s}^{+\infty}d\si\, K_{\si,r+\si}(0,z)\bigg)^3-\bigg(\int_{s}^{2-s+2r}d\si\,  K_{\si,r+\si}(0,z)\bigg)^3 \bigg],
\end{align*}
and accordingly
\begin{align*}
&|\mathcal{R}^{(r)}|\lesssim \mathcal{R}_1^{(r)}+\mathcal{R}_2^{(r)},
\end{align*}
with
\begin{align*}
&\mathcal{R}_1^{(r)}:=\int_{1}^{+\infty}ds \int dz\, \big| K_{s,r+s}(0,z)\big| \bigg(\int_{s}^{+\infty}d\si\, \big| K_{\si,r+\si}(0,z)\big|\bigg)^3\\
\end{align*}
and
\begin{align*}
&\mathcal{R}_2^{(r)}:=\int_{r}^{1+r}ds\int dz\, \big| K_{s,r+s}(0,z)\big| \bigg(\int_{2-s+2r}^{+\infty}d\si\, \big| K_{\si,r+\si}(0,z)\big|\bigg) \bigg(\int_{s}^{+\infty}d\si\, \big| K_{\si,r+\si}(0,z)\big|\bigg)^2.
\end{align*}
For $\mathcal{R}^{(r)}_1$, it is readily checked that
\begin{align*}
&\mathcal{R}_1^{(r)}\lesssim \int_{1}^{+\infty}ds \int dz\, \big| K_{s,r+s}(0,z)\big| \bigg(\int_{s}^{+\infty}\frac{d\si}{\si^{\frac32}}\bigg)^3\lesssim \int_{1}^{+\infty}\frac{ds}{s^{\frac32}} \int dz\, \big| K_{s,r+s}(0,z)\big|\lesssim 1.
\end{align*}
As for $\mathcal{R}^{(r)}_2$, one has
\begin{align*}
\mathcal{R}_2^{(r)}&\lesssim \int_{r}^{1+r}ds\int dz\, \big| K_{s,r+s}(0,z)\big| \bigg(\int_{1}^{+\infty}d\si\, \big| K_{\si,r+\si}(0,z)\big|\bigg) \bigg(\int_{s}^{+\infty}d\si\, \big| K_{\si,r+\si}(0,z)\big|\bigg)^2\\
&\lesssim \int_{r}^{\infty}ds\int dz\, \big| K_{s,r+s}(0,z)\big| \bigg(\int_{1}^{+\infty}\frac{d\si}{\si^{\frac32}}\bigg) \bigg(\int_{s}^{+\infty}\frac{d\si}{\si^{\frac32}}\bigg)^2\\
&\lesssim \int_{r}^{\infty}\frac{ds}{s} \int dz\, \big| K_{s,r+s}(0,z)\big|\lesssim \int_{r}^{\infty}\frac{ds}{s} \lesssim \frac{1}{r^{\frac14}}.
\end{align*}

We have thus shown that
$$\mathcal{R}^{(r)}=\mathcal{O}\Big(\frac{1}{r^{\frac14}} \Big),$$
which completes the proof of \eqref{calcul-K1}.
\end{proof}

We now establish a lower bound for the quantity $\cm^{(r)}$ in \eqref{calcul-K1}. 

\begin{lemma}\label{lemm-cm-2}
With the notation of Lemma \ref{lemm-cm}, there exists $c>0$ such that for all $0<r\leq 1$
\begin{equation*} 
 \cm^{(r)}  \geq c\, r^{-\frac12}.
\end{equation*}
\end{lemma}

\begin{proof}  Let us set  
$$a_\si:=c \,\sinh(2\si)^{\frac32} \quad \text{and} \quad b_\si:=\frac14\Big(\tanh(\si)+\frac{1}{\tanh(\si)}\Big)$$
so that
$$K_\si(0,z) =\frac{1}{a_\si} \exp\big(-b_\si |z|^2\big).$$
Set also for simplicity
$$d_\si:=\tanh(\si), \quad e_\si:=1+\tanh^2(\si),$$
  and 
$$F_{\si_1,\si_2,\si_3}(s_1,s_2):= \Big(e_{s_1} d_{\si_1}d_{\si_2}d_{\si_3} +d_{s_2} \big(e_{ \si_1} d_{\si_2}d_{\si_3}+e_{\si_2} d_{\si_1}d_{\si_3}+e_{\si_3} d_{\si_1}d_{\si_2}\big)  \Big)^{-\frac32}.
$$  
With this notation, we have on the one hand
\begin{align*}
\int dz\, K_{s}(0,z)  K_{\si_1}(0,z) K_{\si_2}(0,z)K_{\si_3}(0,z) &=\frac{c_0}{a_s a_{\si_1}a_{\si_2}a_{\si_3}} \frac{1}{\big(b_s+b_{\si_1}+b_{\si_2}+b_{\si_3} \big)^{\frac32}}\\
&=c_0'\frac{d_s^{\frac32} }{a_s }\bigg( \frac{ (d_{\si_1}d_{\si_2}d_{\si_3})^{\frac32} }{ a_{\si_1}a_{\si_2}a_{\si_3}}\bigg)  F_{\si_1,\si_2,\si_3}(s,s),
\end{align*}
and in the same way
\begin{equation*}
\int dz\, K_{s+r}(0,z)  K_{\si_1}(0,z) K_{\si_2}(0,z)K_{\si_3}(0,z)=c_0'\frac{d_{s+r}^{\frac32} }{a_{s+r} }\bigg( \frac{ (d_{\si_1}d_{\si_2}d_{\si_3})^{\frac32} }{ a_{\si_1}a_{\si_2}a_{\si_3}}\bigg)  F_{\si_1,\si_2,\si_3}(s+r,s+r).
\end{equation*}
As a result, we have
\begin{equation*}
 {\cm^{(r)}} =c c_0'\int_{r}^{+\infty}ds\int_{s}^{s+r} d\si_1 \int_{s}^{s+r} d\si_2 \int_{s}^{s+r} d\si_3\, \Big[I^{(r)}_{s,\si_1,\si_2,\si_3}+II^{(r)}_{s,\si_1,\si_2,\si_3}+III^{(r)}_{s,\si_1,\si_2,\si_3}\Big] ,
\end{equation*}
with
\begin{equation*}
I^{(r)}_{s,\si_1,\si_2,\si_3}:=\bigg(\frac{d_{s}^{\frac32} }{a_{s} }-\frac{d_{s+r}^{\frac32} }{a_{s+r} }\bigg)\bigg( \frac{ (d_{\si_1}d_{\si_2}d_{\si_3})^{\frac32} }{ a_{\si_1}a_{\si_2}a_{\si_3}}\bigg)F_{\si_1,\si_2,\si_3}(s,s),
\end{equation*}
\begin{equation*}
II^{(r)}_{s,\si_1,\si_2,\si_3}:=\frac{d_{s+r}^{\frac32} }{a_{s+r} }\bigg( \frac{ (d_{\si_1}d_{\si_2}d_{\si_3})^{\frac32} }{ a_{\si_1}a_{\si_2}a_{\si_3}}\bigg)\Big( F_{\si_1,\si_2,\si_3}(s,s)-F_{\si_1,\si_2,\si_3}(s+r,s) \Big),
\end{equation*}
and
\begin{equation*}
III^{(r)}_{s,\si_1,\si_2,\si_3} :=\frac{d_{s+r}^{\frac32} }{a_{s+r} }\bigg( \frac{ (d_{\si_1}d_{\si_2}d_{\si_3})^{\frac32} }{ a_{\si_1}a_{\si_2}a_{\si_3}}\bigg)\Big( F_{\si_1,\si_2,\si_3}(s+r,s)-F_{\si_1,\si_2,\si_3}(s+r,s+r) \Big).
\end{equation*}
Observe first that
$$\frac{d_{s}^{\frac32} }{a_{s} }=\frac{c}{\cosh^3(s)} \geq \frac{c}{\cosh^3(s+r)}=\frac{d_{s+r}^{\frac32} }{a_{s+r} },$$
and so
$$I^{(r)}_{s,\si_1,\si_2,\si_3}\geq 0.$$
Then, since $e_s\leq e_{s+r}$, one has for all $A,B>0$,
$$\frac{1}{(A+e_s B)^{\frac32}}\geq \frac{1}{(A+e_{s+r} B)^{\frac32}},$$
which shows that $II^{(r)}_{s,\si_1,\si_2,\si_3} \geq 0$. It follows that
\begin{equation}\label{reduce-to-iii}
 {\cm^{(r)}} 
 \gtrsim\int_{r}^{+\infty}ds\int_{s}^{s+r} d\si_1 \int_{s}^{s+r} d\si_2 \int_{s}^{s+r} d\si_3\, III^{(r)}_{s,\si_1,\si_2,\si_3}.
\end{equation}
To exhibit a lower bound for the latter quantity, we will use the basic inequality: for all $A,B\geq 0$ and $0 \leq S\leq T$,
\begin{align*}
&\frac{1}{(A+SB)^{\frac32}}-\frac{1}{(A+TB)^{\frac32}}=\frac32 B\int_{S}^{T}\frac{du}{(A+uB)^{\frac52}} \gtrsim \frac{B(T-S)}{(A+TB)^{\frac52}}.
\end{align*}
This yields here
\begin{multline*}
F_{\si_1,\si_2,\si_3}(s+r,s)-F_{\si_1,\si_2,\si_3}(s+r,s+r) \gtrsim\\
\gtrsim  \Big(e_{ \si_1} d_{\si_2}d_{\si_3}+e_{\si_2} d_{\si_1}d_{\si_3}+e_{\si_3} d_{\si_1}d_{\si_2}\Big) \Big(F_{\si_1,\si_2,\si_3}(s+r,s+r)\Big)^{\frac53} (d_{s+r}-d_s) \geq 0.
\end{multline*}
Now,  for all $s \in [r,2r]$ and $\si_1,\si_2,\si_3\in [s,s+r] \subset [r,3r]$, if   $r>0$ is small enough
\begin{equation*}
III^{(r)}_{s,\si_1,\si_2,\si_3} \gtrsim  r^{-\frac92},
\end{equation*}
where we used that, up to a multiplicative factor,  $a_{\sigma}\sim \sigma^{\frac32}$, $d_{\sigma}\sim \sigma$, and  {$e_{\sigma}\sim 1$}. Going back to~\eqref{reduce-to-iii}, we obtain that for every $0<r\leq 1$,
\begin{equation*}
 {\cm^{(r)}} \gtrsim \int_{r}^{2r}ds\int_{s}^{s+r} d\si_1 \int_{s}^{s+r} d\si_2 \int_{s}^{s+r} d\si_3\, III^{(r)}_{s,\si_1,\si_2,\si_3}
\gtrsim  r^{-\frac12},
\end{equation*}
which yields the desired lower bound.
\end{proof}


\section{Uniqueness of the invariant measure at high temperature}\label{sec:uniq-inv}

In the preceding section, we have established the existence of an invariant measure for equation~\eqref{eq:intro}, for every $\la>0$. We now address the question of uniqueness in the high-temperature regime, that is, for $\la>0$ small enough.

\smallskip

To simplify the presentation, we will not write down the dependence on $n \geq 0$ in the (approximated) equation \eqref{appro-equ-introd}, and thus rephrase the model as 

\begin{equation} \label{eq-la}
\left\{
\begin{aligned}
&  \partial_t X^{(\la),u} + HX^{(\la),u}= -\la \big(X^{(\la),u}\big)^3+\big(3\la \frakc^{\mathbf{1}}-9\la^2 \frakc^{\mathbf{2}}\big) X^{(\la),u} + \xi, \quad t>0, \quad x \in \R^3, \\
& X^{(\la),u}_0 = u \in \cb_{\infty}^{-\frac58},
\end{aligned}
\right.
\end{equation}
for every $0<\la\leq 1$. Note that the above choice of regularity for the initial condition, namely  $u \in \cb_{\infty}^{-\frac58}$, is somewhat arbitrary: we could replace it with $u \in \cb_{\infty}^{-\frac12-\eps}$, for any $\eps>0$.

\

Our main technical result in this section (to be compared with the statement of \cite[Proposition~3.25]{DHYZ25}) reads as follows.

\begin{theorem}\label{thm-uniqueness}
There exist constants $\la_\star>0$ and $c_0,c_1>0$ such that for all $\la\in (0,\la_\star)$, all $T\geq 1$ and uniformly in $n \geq 0$, 
\begin{equation}\label{boun-to-uni}
\sup_{u_1,u_2\in \cb_{\infty}^{-\frac58}}\mathbb{E}\Big[\big\|X_T^{(\la),u_1}-X_T^{(\la),u_2}\big\|_{L^\infty}\Big]\leq c_0 e^{-c_1 T} .
\end{equation}
\end{theorem}
 
In the sequel, if $\rho$ is a measure on $\cb_{\infty}^{-\frac58}$ and $\Psi:\cb_{\infty}^{-\frac58}\to\R$ is a Lipschitz function, we  use the standard notation 
$$  \rho(\Psi):=   \int_{\cb_{\infty}^{-\frac58}} \Psi(v) \rho(dv). $$
Recall also the definition \eqref{transi}  of the transition semigroup $(P_t)_{t\ge 0}$, namely
 \begin{equation*}
\big(P^{(\lambda)}_t \Psi\big) (u): = \mathbb{E} \big[ \Psi (X_t^{(\lambda),u})\big], \quad u \in \B_{\infty}^{-\frac58}, \quad t\ge 0.
\end{equation*}
We also set, following the standard convention,
$$\big(P_t^{(\la)}\rho\big)(\Psi) :=\int_{\cb_{\infty}^{-\frac58}}  \big(P_t^{(\la)}\Psi\big)(v) \, \rho(dv). $$

Then we have the following result. 

\begin{corollary}\label{coro:uniquen-inva-m}
For every $\la \in (0,\la_\star)$, the transition semigroup $(P^{(\la)}_t)_{t\geq 0}$ generated by equation~\eqref{eq-la} admits a unique invariant measure $\rho^{(\lambda)}$, supported by $\cb_{\infty}^{-\frac58}$. Moreover, for all $u\in \cb_{\infty}^{-\frac58}$, all bounded Lipschitz function $\Psi:\cb_{\infty}^{-\frac58}\to\R$ and all $T\geq 1$, it holds that
\begin{equation}\label{expo-conve}
\big|  \big(P^{(\lambda)}_T \Psi\big) (u) -\rho^{(\lambda)}(\Psi)  \big|\leq c_0 \big\|\Psi\big\|_{\text{Lip}} e^{-c_1 T} .
\end{equation}
\end{corollary}

\begin{proof}[Proof of Corollary \ref{coro:uniquen-inva-m}]
Once endowed with \eqref{boun-to-uni}, the argument is standard. Namely, considering two measures $\rho,\tilde{\rho}$ on $\cb_{\infty}^{-\frac58}$ and a bounded Lipschitz function $\Psi:\cb_{\infty}^{-\frac58}\to\R$, we can write for all $T\geq 1$
\begin{multline*}
\big|\big(P_T^{(\la)}\rho\big)(\Psi)-\big(P_T^{(\la)}\tilde{\rho}\big)(\Psi)\big|=\Big|\int_{\cb_{\infty}^{-\frac58}} \big(P_T^{(\la)}\Psi\big)(u) \, \rho(du)-\int_{\cb_{\infty}^{-\frac58}} \big(P_T^{(\la)}\Psi\big)(\tilde{u}) \, \tilde{\rho}(d\tilde{u})\Big|\\
\begin{aligned}
&=\Big|\int_{\cb_{\infty}^{-\frac58}} \mathbb{E}\big[ \Psi\big(X^{(\la),u}_T\big)\big] \, \rho(du)-\int_{\cb_{\infty}^{-\frac58}} \mathbb{E}\big[ \Psi\big(X^{(\la),\tilde{u}}_T\big)\big] \, \tilde{\rho}(d\tilde{u})\Big|\\
&=\Big|\int_{\cb_{\infty}^{-\frac58}\times \cb_{\infty}^{-\frac58}} \Big(\mathbb{E}\big[ \Psi\big(X^{(\la),u}_T\big)\big] -\mathbb{E}\big[ \Psi\big(X^{(\la),\tilde{u}}_T\big)\big] \Big)\,\rho(du) \tilde{\rho}(d\tilde{u})\Big|\\
&\leq \int_{\cb_{\infty}^{-\frac58}\times \cb_{\infty}^{-\frac58}} \mathbb{E}\Big[ \big|\Psi\big(X^{(\la),u}_T\big) -\Psi\big(X^{(\la),\tilde{u}}_T\big)\big|\Big] \,\rho(du) \tilde{\rho}(d\tilde{u})\\
&\leq \big\|\Psi\big\|_{\text{Lip}} \int_{\cb_{\infty}^{-\frac58}\times \cb_{\infty}^{-\frac58}} \mathbb{E}\Big[ \big\|X^{(\la),u}_T -X^{(\la),\tilde{u}}_T\big\|_{L^\infty}\Big]\,  \rho(du) \tilde{\rho}(d\tilde{u})\leq c_0 \big\|\Psi\big\|_{\text{Lip}} e^{-c_1 T} .
\end{aligned}
\end{multline*}
Both conclusions of the statement are then immediate. First, if $\rho$ and $\tilde{\rho}$ are invariant measures, we obtain that for all $T\geq 1$
\begin{equation*}
\big|\rho(\Psi)-\tilde{\rho}(\Psi)\big|=\big|\big(P_T^{(\la)}\rho\big)(\Psi)-\big(P_T^{(\la)}\tilde{\rho}\big)(\Psi)\big|\leq c_0 \big\|\Psi\big\|_{\text{Lip}} e^{-c_1 T} 
\end{equation*}
and hence $\rho(\Psi)=\tilde{\rho}(\Psi)$, so that $\rho=\tilde{\rho}$. On the other hand, denoting by $\rho^{(\lambda)}$ the unique invariant measure of $P^{(\la)}$, we have for $\nu=\delta_{u}$,
\begin{equation*}
\big|\big(P^{(\lambda)}_T \Psi\big) (u)-\rho^{(\lambda)}(\Psi)\big|=\big|\big(P_T^{(\la)}\nu\big)(\Psi)-\big(P_T^{(\la)}\rho^{(\lambda)}\big)(\Psi)\big|\leq c_0 \big\|\Psi\big\|_{\text{Lip}} e^{-c_1 T} ,
\end{equation*}
which yields \eqref{expo-conve}.
\end{proof}

\subsection{Notation and scaling}

Let us start the proof of Theorem \ref{thm-uniqueness} by introducing a few technical tools and notation.

\subsubsection{Notation 1}\label{subsec:tree}

For every $\la\in (0,1]$, we set, following the pattern described in \eqref{def-luxost}-\eqref{ord5},
\begin{align*}
\big\{\<Psi>^{(\la)}_{s,.},\<Psi2>^{(\la)}_{s,.},\<IPsi3>^{(\la)}_{s,.},\<PsiIPsi3>^{(\la)}_{s,.},\<Psi2IPsi2>^{(\la)}_{s,.},\<Psi2IPsi3>^{(\la)}_{s,.}\big\}:=Z_{s,.}(\la^{\frac12}\xi).
\end{align*}

In other words, for all $t\geq s$,
\begin{equation*}
\<Psi>^{(\la)}_{s,t}=\la^{\frac12}\<Psi>_{s,t}, \quad\quad \<Psi2>^{(\la)}_{s,t}=\la\<Psi2>_{s,t}, \quad\quad  \<IPsi2>^{(\la)}_{s,t}=\la\<IPsi2>_{s,t}, \quad\quad \<IPsi3>^{(\la)}_{s,t}=\la^{\frac32}\<IPsi3>_{s,t},
\end{equation*}
\begin{equation*}
\<PsiIPsi3>^{(\la)}_{s,t}=\la^2\<PsiIPsi3>_{s,t}, \quad \quad \<Psi2IPsi2>^{(\la)}_{s,t}=\la^2\<Psi2IPsi2>_{s,t}, \quad \quad \<Psi2IPsi3>^{(\la)}_{s,t}=\la^{\frac52}\<Psi2IPsi3>_{s,t} .
\end{equation*}
Denoting by $\big\|Z_{s,.}(\la^{\frac12}\xi)\big\|_{\mathcal Z_{\frac{\eps}{2},[s,T]}}$ the quantity associated with $(\la^{\frac12}\xi)$ through \eqref{K}, one has for all $0 < \lambda \leq 1$
\begin{equation*} 
\big\|Z_{s,.}(\la^{\frac12}\xi)\big\|_{\mathcal Z_{\frac{\eps}{2},[s,T]}} \leq \lambda^{\frac12} \big\|Z_{s,.}(\xi)\big\|_{\mathcal Z_{\frac{\eps}{2},[s,T]}}.
\end{equation*}
In the same way, observe that with the notation in \eqref{constante-re}, one has trivially

\begin{equation*}
\frakc^{\mathbf{1},(\la)}_{s,t}(x):=\mathbb{E}\Big[ \big| \<Psi>^{(\la)}_{s,t}(x)\big|^2\Big]=\la \frakc^{\mathbf{1}}_{s,t}(x)
\end{equation*}
while
\begin{equation*}
\frakc^{\mathbf{2},(\la)}_{s,t}(x):=\mathbb{E}\Big[ \<Psi2>^{(\la)}_{s,t}(x) \<IPsi2>^{(\la)}_{s,t}(x)\Big]=\la^2 \frakc^{\mathbf{2}}_{s,t}(x),
\end{equation*}
which leads us to define the rescaled constants (see \eqref{cstts-rigor-1}-\eqref{cstts-rigor-2}) as
\begin{equation}\label{def:frakcla}
\frakc^{\mathbf{1},(\la)}(x):=\la \frakc^{\mathbf{1}}(x), \quad \frakc^{\mathbf{2},(\la)}(x):=\la^2 \frakc^{\mathbf{2}}(x), \quad  \frakc^{(\la)} := 3\la \frakc^{\mathbf{1}}-9\la^2 \frakc^{\mathbf{2}}  ,
\end{equation}
together with
\begin{equation}\label{def:frakclaoverline}
\overline{\frakc}^{(\la)}_{s,t} :=\big(3\frakc^{\mathbf{1},(\la)} - 9\frakc^{\mathbf{2},(\la)}\big)-\big(3\frakc^{\mathbf{1},(\la)}_{s,t} - 9\frakc^{\mathbf{2},(\la)}_{s,t}\big).
\end{equation}
With the notation in \eqref{llbrack}, it is readily checked that
\begin{equation*} 
\big\llbracket  \overline{\frakc}^{(\la)}_{s,.}\big\rrbracket_{3\eps,[s,T]}\lesssim \la\Big( \big\llbracket  \overline{\frakc}^{\mathbf{1}}_{s,.}\big\rrbracket_{3\eps,[s,T]}+ \big\llbracket  \overline{\frakc}^{\mathbf{2}}_{s,.}\big\rrbracket_{3\eps,[s,T]}\Big).
\end{equation*}
In the end, denoting by $M^{(\la)}_{s,T}$ the quantity associated with $(\la^{\frac12}\xi)$ and $\overline{\frakc}^{(\la)}_{s,t}$ through \eqref{notation-mst}, one has for all $0 < \lambda \leq 1$
\begin{equation*} 
  M^{(\la)}_{s,T} \leq   M_{s,T}.
\end{equation*}

\subsubsection{Notation 2}

First, given (say smooth) functions $\frakc^{\mathbf{1}},\frakc^{\mathbf{2}}$ on $\R^3$, a (say smooth) noise $\xi$ and an initial condition $u\in \cb_{\infty}^{-\frac58}$, let us denote by
$$\Phi\big(u;\frakc^{\mathbf{1}},\frakc^{\mathbf{2}};\xi\big)$$
the solution of the equation
\begin{equation} \label{eq-depart}
\left\{
\begin{aligned}
&  \partial_t X + HX= -X^3+\big(3 \frakc^{\mathbf{1}}-9 \frakc^{\mathbf{2}}\big) X + \xi, \quad t>0, \quad x \in \R^3, \\
& X_0 = u .
\end{aligned}
\right.
\end{equation}
On the other hand, for every $0<\la\leq 1$, denote by
$$\Phi^{(\la)}\big(u;\frakc^{\mathbf{1}},\frakc^{\mathbf{2}};\xi\big)$$
the solution of the generalized problem
\begin{equation*} 
\left\{
\begin{aligned}
&  \partial_t X+ HX= -\la X^3+\big(3\la \frakc^{\mathbf{1}}-9\la^2 \frakc^{\mathbf{2}}\big) X+ \xi, \quad t>0, \quad x \in \R^3, \\
& X_0 = u  .
\end{aligned}
\right.
\end{equation*}

\

With this notation, it can be checked that for every $0<\la\leq 1$,
\begin{equation}\label{identif}
\Phi^{(\la)}\big(u;\frakc^{\mathbf{1}},\frakc^{\mathbf{2}};\xi\big)=\la^{-\frac12} \Phi\big(\la^{\frac12}u;\la\frakc^{\mathbf{1}},\la^2 \frakc^{\mathbf{2}};\la^{\frac12}\xi\big).
\end{equation}
Indeed, by setting 
\begin{equation}\label{notation-x-la}
\widehat{X}^{(\la),u}:=\Phi\big(\la^{\frac12}u;\la\frakc^{\mathbf{1}},\la^2 \frakc^{\mathbf{2}};\la^{\frac12}\xi\big) \quad \text{and} \quad X^{(\la),u}:=\la^{-\frac12}\widehat{X}^{(\la),u},
\end{equation}
one has $X^{(\la),u}_0=\la^{-\frac12}\widehat{X}^{(\la),u}_0=u$ and
\begin{align*}
\partial_t X^{(\la),u} + HX^{(\la),u}&=\la^{-\frac12}\big(\partial_t \widehat{X}^{(\la),u} + H\widehat{X}^{(\la),u}\big)\\
&=\la^{-\frac12} \big(-(\widehat{X}^{(\la),u})^3+(3\la \frakc^{\mathbf{1}}-9\la^2 \frakc^{\mathbf{2}}) \widehat{X}^{(\la),u} + \la^{\frac12}\xi\big)\\
&=-\la \big(\la^{-\frac12}\widehat{X}^{(\la),u}\big)^3+(3\la \frakc^{\mathbf{1}}-9\la^2 \frakc^{\mathbf{2}})\big(\la^{-\frac12} \widehat{X}^{(\la),u}\big) +\xi\\
&=-\la \big(X^{(\la),u}\big)^3+(3\la \frakc^{\mathbf{1}}-9\la^2 \frakc^{\mathbf{2}})X^{(\la),u} +\xi,
\end{align*}
and so, by uniqueness of the solution,
$$X^{(\la),u}=\Phi^{(\la)}\big(u;\frakc^{\mathbf{1}},\frakc^{\mathbf{2}};\xi\big).$$

\

\subsubsection{Known controls}

Taking the notation in \eqref{notation-x-la} and the observations in Section \ref{subsec:tree} into account, we deduce from Theorem~\ref{thm:global_Lp} that for all $0\leq s\leq t$, $\widehat{X}^{(\la),u}_t$ can be decomposed as
\begin{equation}\label{decompo-x-hat-s-t}
\widehat{X}^{(\la),u}_t=\la^{\frac12}\<Psi>_{s,t}-\la^{\frac32} \<IPsi3>_{s,t}+\vtila_{s,t}+\wtila_{s,t},
\end{equation}
where the processes $\vtila_{s,t}$ and $\wtila_{s,t}$ are controlled through the following coming-down-from-infinity bounds:

\begin{equation}\label{bbb1}
\max\Big(\sup_{t\in [s,T]} \big\|\vtila_{s,t}\big\|_{\cb^{-\eps}_\infty}, \sup_{t\in [s,T]} (t-s)^{\frac38+\nu}\big\|\vtila_{s,t}\big\|_{\cb^{\frac34-\eps}_\infty} \Big)\lesssim   M^\ka_{s,T},
\end{equation}
as well as
\begin{equation}\label{bbb2}
\max\Big(\sup_{t\in [s,T]} (t-s)^{\frac12+\nu}\big\|\wtila_{s,t}\big\|_{\cb^{-\eps}_\infty}, \sup_{t\in [s,T]} (t-s)^{\frac32+\nu}\big\|\wtila_{s,t}\big\|_{\cb^{\frac32-\eps}_\infty} \Big)\lesssim   M^\ka_{s,T}.
\end{equation}

\

By setting
$$\vla_{s,t} :=\la^{-\frac12}\vtila_{s,t} \quad \text{and} \quad \wla_{s,t}:=\la^{-\frac12} \wtila_{s,t},$$
we immediately deduce from \eqref{identif} and \eqref{decompo-x-hat-s-t} that for all $0\leq s\leq t$,
\begin{equation}\label{decompo-x-s-t}
X^{(\la),u}_t=\<Psi>_{s,t}-\la \<IPsi3>_{s,t}+\vla_{s,t}+\wla_{s,t},
\end{equation}

Furthermore, the above bounds for $(\vtila,\wtila)$ can be readily transferred to $(\vla,\wla)$. More precisely, setting $U^{(\la),u}_{s,t}:=\vla_{s,t}+\wla_{s,t}$, it holds that

\begin{proposition}\label{cor:comes-down}
Let $\eps>0$ be small enough. There exists $\ka\geq 1$ such that for all $0\leq s<T\leq s+1$ and $0<\la\leq 1$, one has
\begin{equation*}
 \max\Big( \sup_{t\in [s,T]} (t-s)^{\frac12+\eps}\big\|U^{(\la),u}_{s,t}\big\|_{L^\infty},  \sup_{t\in [s,T]} (t-s)^{\frac56+2\eps}\big\|U^{(\la),u}_{s,t}\big\|_{\cb^{\frac12+2\eps}_\infty} \Big)\lesssim {\la^{-\frac12}}M^\ka_{s,T}.
\end{equation*}
\end{proposition}

\begin{proof}
By rescaling, the bounds \eqref{bbb1} and \eqref{bbb2} directly imply that for all $\nu>0$ 
 \begin{equation*}
\max\Big( \sup_{t\in [s,T]} \big\|\vla_{s,t}\big\|_{\cb^{-\eps}_\infty},  \sup_{t\in [s,T]} (t-s)^{\frac38+\nu}\big\|\vla_{s,t}\big\|_{\cb^{\frac34-\eps}_\infty} \Big)\lesssim {\la^{-\frac12}} M^\ka_{s,T},
\end{equation*}
as well as
\begin{equation*}
\max\Big( \sup_{t\in [s,T]} (t-s)^{\frac12+\nu}\big\|\wla_{s,t}\big\|_{\cb^{-\eps}_\infty},  \sup_{t\in [s,T]} (t-s)^{\frac32+\nu}\big\|\wla_{s,t}\big\|_{\cb^{\frac32-\eps}_\infty} \Big)\lesssim {\la^{-\frac12}} M^\ka_{s,T}.
\end{equation*}
Then, for $\eps>0$ we have  $\cb_\infty^{\eps}\subset L^\infty$, and thus
\begin{equation}\label{ddt}
 (t-s)^{\frac12+3\eps}\big\|U^{(\la),u}_{s,t}\big\|_{L^\infty}\lesssim (t-s)^{\frac12+2\eps}\Big(\big\|\vla_{s,t}\big\|_{\cb^{\eps}_\infty}+\big\|\wla_{s,t}\big\|_{\cb^{\eps}_\infty}\Big).
\end{equation}
We now use the interpolation bounds
$$\big\|\vla_{s,t}\big\|_{\cb^{\eps}_\infty} \leq \big\|\vla_{s,t}\big\|^{1-\frac83\eps}_{\cb^{-\eps}_\infty} \big\|\vla_{s,t}\big\|^{\frac83\eps}_{\cb^{\frac34-\eps}_\infty} \lesssim (t-s)^{-\eps-\frac83\eps \nu}{\la^{-\frac12}} M^\ka_{s,T}$$
and
$$\big\|\wla_{s,t}\big\|_{\cb^{\eps}_\infty} \leq \big\|\wla_{s,t}\big\|^{1-\frac43\eps}_{\cb^{-\eps}_\infty} \big\|\wla_{s,t}\big\|^{\frac43\eps}_{\cb^{\frac32-\eps}_\infty} \lesssim (t-s)^{-\frac12-2\eps-\nu}{\la^{-\frac12}} M^\ka_{s,T},$$
which, together with \eqref{ddt}, allows us to prove the first bound.

\smallskip

The second bound is obtained in a similar way.
\end{proof}

\subsubsection{Auxiliary system}

We introduce an auxiliary system that will prove to be closely related to the problem under consideration (see Section \ref{subsec:stra} below). To this end, we fix $s\geq 0$, $f\in L^\infty$, as well as two (regular) space-time functions $R,S$, with $S\geq 0$. Then for every $t\geq s$, we consider the solution 
$$ (\widetilde{\ell},\widetilde{m})$$
of  the following linear system: for every $t\geq s$, 
 \begin{equation}\label{systlm0}
\left\{
\begin{aligned}
&  \widetilde{\ell}_t  =-3\la  \int_s^t dr\, e^{-(t-r)H}\Big(  (\widetilde{\ell}+\widetilde{m}) \pl  \<Psi2>_{s,r}
\Big), \\
& \widetilde{m}_t  =e^{-(t-s)H} f+\int_s^t dr\, e^{-(t-r)H}\Big(\widetilde{\mathfrak{G}}^{(\la)}\big(\widetilde{\ell},\widetilde{m};R,S,Z_{s,.}\big)_{s,r}\Big),
\end{aligned}
\right.
\end{equation}
with 
\begin{multline}\label{sgma1}
\widetilde{\mathfrak{G}}^{(\la)}(\widetilde{\ell},\widetilde{m};R,S,Z_{s,.})_{s,r}:=-3\la (\widetilde{\ell}_r+\widetilde{m}_r) \pe  \<Psi2>_{s,r}-3\la (\widetilde{\ell}_r+\widetilde{m}_r) \pg  \<Psi2>_{s,r}\\
-9\la^2 \frakc^{\mathbf{2}}_{s,r}\cdot (\widetilde{\ell}_r+\widetilde{m}_r)-\la\, S_r \cdot  (\widetilde{\ell}_r+\widetilde{m}_r)+\la\, R_r \cdot  (\widetilde{\ell}_r+\widetilde{m}_r).
\end{multline}

In the sequel, we will crucially use the condition $S\geq 0$ (see Proposition \ref{PropS+}), which explains why we distinguish $S$ from $R$ in the above formulation. 

\medskip

Let us then set, as in Section \ref{Sect32},
\begin{align}
\mathrm{com}(\widetilde{\ell},\widetilde{m};Z_{s,.})_{s,t}&: = \mathrm{com}_1(\widetilde{\ell},\widetilde{m};Z_{s,.})_{s,t} \pe  \<Psi2>_{s,t} +\mathrm{com}_2 (\widetilde{\ell}+\widetilde{m};Z_{s,.})_{s,t} \nonumber \\
\mathrm{com}_1(\widetilde{\ell},\widetilde{m};Z_{s,.})_{s,t} &:=  \int_s^t e^{-(t-r)H} \big[(\widetilde{\ell}+\widetilde{m})_r \pl  \<Psi2>_{s,r}\big] \, dr -(\widetilde{\ell}+\widetilde{m})_t \pl  \<IPsi2>_{s,t},\nonumber\\ 
\mathrm{com}_2 (\widetilde{\ell}+\widetilde{m};Z_{s,.})_{s,t}&:= \big[ \pl, \pe\big] \Big(\widetilde{\ell}_t+\widetilde{m}_t, \<IPsi2>_{s,t}, \<Psi2>_{s,t}\Big), \nonumber
\end{align}
where we recall that
\begin{equation*} 
[\pl, \pe] (f,g,h) := (f \pl g) \pe h- f(g \pe h).
\end{equation*}

Using this notation, we can rephrase the expression \eqref{sgma1} in the following way:

\begin{lemma}
The expression \eqref{sgma1} can equivalently be written as
\begin{multline}\label{decc}
\widetilde{\mathfrak{G}}^{(\la)}(\widetilde{\ell},\widetilde{m};R,S,Z_{s,.})_{s,t}=\la\, R_t \cdot (\widetilde{\ell}_t+\widetilde{m}_t)-\la\, S_r \cdot  (\widetilde{\ell}_r+\widetilde{m}_r)-3\la (\widetilde{\ell}_t+\widetilde{m}_t) \pg  \<Psi2>_{s,t}\\
-3\la\,  \widetilde{m}_t\pe  \<Psi2>_{s,t}+9\la^2\,  \mathrm{com}\big(\widetilde{\ell},\widetilde{m};Z_{s,.}\big)_{s,t} +9\la^2\, (\widetilde{\ell}_t+\widetilde{m}_t )\, \<Psi2IPsi2>_{s,t} .
\end{multline}
\end{lemma}

\begin{proof}
Starting from \eqref{sgma1}, we can write
\begin{multline*}
\widetilde{\mathfrak{G}}^{(\la)}(\widetilde{\ell},\widetilde{m};R,S,Z_{s,.})_{s,t}=\la\, R_t \cdot (\widetilde{\ell}_t+\widetilde{m}_t)-\la\, S_r \cdot  (\widetilde{\ell}_r+\widetilde{m}_r)-3\la (\widetilde{\ell}_t+\widetilde{m}_t) \pg  \<Psi2>_{s,t}\\
\hspace{3cm}-3\la\,  \widetilde{m}_t\pe  \<Psi2>_{s,t}-3\Big[\la\, \widetilde{\ell}_t \pe  \<Psi2>_{s,t}+3 \la^2\, \frakc^{\mathbf{2}}_{s,t}\, (\widetilde{\ell}_t+\widetilde{m}_t)\Big]\\
\begin{aligned}
&=\la\, R_t \cdot (\widetilde{\ell}_t+\widetilde{m}_t)-\la\, S_r \cdot  (\widetilde{\ell}_r+\widetilde{m}_r)-3\la (\widetilde{\ell}_t+\widetilde{m}_t) \pg  \<Psi2>_{s,t}-3\la\,  \widetilde{m}_t\pe  \<Psi2>_{s,t}\\
&\hspace{2cm}+9\la^2\bigg[ \<Psi2>_{s,t} \pe \Big(\int_s^t e^{-(t-r)H} \big((\widetilde{\ell}+\widetilde{m} )_r\pl  \<Psi2>_{s,r}\big)   \Big)- \frakc^{\mathbf{2}}_{s,t}\, (\widetilde{\ell}_t+\widetilde{m}_t)\bigg]\\
&=\la\, R_t \cdot (\widetilde{\ell}_t+\widetilde{m}_t)-\la\, S_r \cdot  (\widetilde{\ell}_r+\widetilde{m}_r)-3\la (\widetilde{\ell}_t+\widetilde{m}_t) \pg  \<Psi2>_{s,t}-3\la\,  \widetilde{m}_t\pe  \<Psi2>_{s,t}\\
&\hspace{1cm}+9\la^2\, \<Psi2>_{s,t} \pe \mathrm{com}_1\big(\widetilde{\ell},\widetilde{m};Z_{s,.}\big)_{s,t} +9\la^2\bigg[\<Psi2>_{s,t} \pe \Big((\widetilde{\ell}_t+\widetilde{m}_t )\pl  \<IPsi2>_{s,t}   \Big)- \frakc^{\mathbf{2}}_{s,t}\, (\widetilde{\ell}_t+\widetilde{m}_t)\bigg]\\
&=\la\, R_t \cdot (\widetilde{\ell}_t+\widetilde{m}_t)-\la\, S_r \cdot  (\widetilde{\ell}_r+\widetilde{m}_r)-3\la (\widetilde{\ell}_t+\widetilde{m}_t) \pg  \<Psi2>_{s,t}-3\la\,  \widetilde{m}_t\pe  \<Psi2>_{s,t}\\
&+9\la^2\Big[ \<Psi2>_{s,t} \pe \mathrm{com}_1\big(\widetilde{\ell},\widetilde{m};Z_{s,.}\big)_{s,t} + \mathrm{com}_2\big(\widetilde{\ell}+\widetilde{m};Z_{s,.}\big)_{s,t }\Big] +9\la^2\bigg[(\widetilde{\ell}_t+\widetilde{m}_t ) \Big(\<Psi2>_{s,t} \pe  \<IPsi2>_{s,t}   \Big)- \frakc^{\mathbf{2}}_{s,t}\, (\widetilde{\ell}_t+\widetilde{m}_t)\bigg],
\end{aligned}
\end{multline*}
which, by \eqref{ordres4}, implies \eqref{decc}.
\end{proof}
 
\

Once endowed with the solution $(\widetilde{\ell},\widetilde{m})$ of \eqref{systlm0}, we set
\begin{equation}\label{not-l-ti}
\widetilde{L}^{(\la)}(f,s;R,S,Z_{s,.})_{t}:=\widetilde{\ell}_t+\widetilde{m}_t.
\end{equation}

\

\subsubsection{Dynamics of the difference}\label{subsec:stra}
Using the notation $\widetilde{L}^{(\la)}$ in \eqref{not-l-ti}, the auxiliary system \eqref{systlm0} and the original equation \eqref{eq-la} are linked through the following identity: 
\begin{proposition}\label{prop:egalite}
 For all $0<\la\leq 1$, $u_1,u_2\in \cb_{\infty}^{-\frac58}$ and $0 \leq s \leq t$, it holds that
\begin{equation} \label{egalite}
X^{(\la),u_1}_t-X^{(\la),u_2}_t=\widetilde{L}^{(\la)} \big(X^{(\la),u_1}_s-X^{(\la),u_2}_s,s; R^{(\la),u_1,u_2}_{s,.},S^{(\la),u_1,u_2}_{s,.},Z_{s,.}\big)_t,
\end{equation}
where we have set
\begin{align}\label{defi-s-la} 
&S^{(\la),u_1,u_2}_{s,r}:=\frac12 \big(\widetilde{U}^{(\la),u_1}_{s,r}\big)^2 +\frac12 \big(\widetilde{U}^{(\la),u_2}_{s,r}\big)^2 +\frac12 \big(\widetilde{U}^{(\la),u_1}_{s,r}+\widetilde{U}^{(\la),u_2}_{s,r}\big)^2   
\end{align}
and (recalling the definition \eqref{def:frakclaoverline} of $\overline{\frakc}^{(\la)}$)
\begin{align}
&R^{(\la),u_1,u_2}_{s,r}:=6\, \la \<PsiIPsi3nr>_{s,r}-3\,  \<Psi>_{s,r}\big(U^{(\la),u_1}_{s,r}+U^{(\la),u_2}_{s,r}\big)+ {\frac{1}{\la}\overline{\frakc}^{(\la)}_{s,r}}.\label{defi-r-la}
\end{align}
\end{proposition}

\begin{proof}
For convenience, let us define for all $0<\la\leq 1$, $u_1,u_2\in \cb_{\infty}^{-\frac58}$ and $t\geq 0$,
\begin{equation*} 
f^{(\la),u_1,u_2}_t:=X^{(\la),u_1}_t-X^{(\la),u_2}_t.
\end{equation*}
Going back to \eqref{eq-la} (and recalling the definition \eqref{def:frakcla} of $\frakc^{(\la)}$), one has for all $t\geq s$
\begin{align}
f^{(\la),u_1,u_2}_t&=e^{-(t-s)H} f^{(\la),u_1,u_2}_s \nonumber\\
&\hspace{1cm}+\int_s^t dr\,  e^{-(t-r)H}\Big(-\la \Big( \big(X^{(\la),u_1}_r\big)^3-\big(X^{(\la),u_2}_r\big)^3\Big) +\frakc^{(\la)}\big( X^{(\la),u_1}_r -X^{(\la),u_2}_r\big) \Big)\nonumber\\
&=e^{-(t-s)H} f^{(\la),u_1,u_2}_s+\int_s^t dr\,  e^{-(t-r)H}\big(A^{(\la),u_1,u_2}_r \cdot f^{(\la),u_1,u_2}_r \big), \label{meme-eq}
\end{align}
where we have set
\begin{align*}
 A^{(\la),u_1,u_2}_r &:=-\la \big(X^{(\la),u_1}_r\big)^2-\la X^{(\la),u_1}_r X^{(\la),u_2}_r-\la\big(X^{(\la),u_2}_r\big)^2+\frakc^{(\la)}.
\end{align*}
Using the decomposition \eqref{decompo-x-s-t}, and setting
\begin{equation}\label{defi:ula}
U^{(\la),u}_{s,t}:=\vla_{s,t}+\wla_{s,t}, \quad \quad \widetilde{U}^{(\la),u}_{s,t}=-\la \<IPsi3>_{s,t}+U^{(\la),u}_{s,t},
\end{equation}
we can expand the quantity $A^{(\la),u_1,u_2}_r$ as
\begin{align*}
A^{(\la),u_1,u_2}_r&=-\la \big(\<Psi>_{s,r}+\widetilde{U}^{(\la),u_1}_{s,r}\big)^2-\la \big(\<Psi>_{s,r}+\widetilde{U}^{(\la),u_1}_{s,r}\big) \big(\<Psi>_{s,r}+\widetilde{U}^{(\la),u_2}_{s,r}\big)-\la\big(\<Psi>_{s,r}+\widetilde{U}^{(\la),u_2}_{s,r}\big)^2+\frakc^{(\la)}\\
&=-3\la \big(\<Psi>_{s,r}\big)^2-3\la \<Psi>_{s,r}\big(\widetilde{U}^{(\la),u_1}_{s,r}+\widetilde{U}^{(\la),u_2}_{s,r}\big) +\frakc^{(\la)}-\la \big[\big(\widetilde{U}^{(\la),u_1}_{s,r}\big)^2 +\big(\widetilde{U}^{(\la),u_2}_{s,r}\big)^2 +\widetilde{U}^{(\la),u_1}_{s,r} \widetilde{U}^{(\la),u_2}_{s,r}  \big]\\
&=-3\la \<Psi2>_{s,r}-9 \la^2\frakc^{\mathbf{2}}_{s,r}-3\la \<Psi>_{s,r}\big(-2\la \<IPsi3>_{s,r}+U^{(\la),u_1}_{s,r}+U^{(\la),u_2}_{s,r}\big) \\
&\hspace{1cm}-\la \bigg[\frac12 \big(\widetilde{U}^{(\la),u_1}_{s,r}\big)^2 +\frac12 \big(\widetilde{U}^{(\la),u_2}_{s,r}\big)^2 +\frac12 \big(\widetilde{U}^{(\la),u_1}_{s,r}+\widetilde{U}^{(\la),u_2}_{s,r}\big)^2  \bigg]+ {\overline{\frakc}^{(\la)}_{s,r}},
\end{align*}
and finally
\begin{align*}
&A^{(\la),u_1,u_2}_r=-3\la \<Psi2>_{s,r}-9 \la^2\frakc^{\mathbf{2}}_{s,r}+\la\, R^{(\la),u_1,u_2}_{s,r}-\la\, S^{(\la),u_1,u_2}_{s,r}.
\end{align*}

\smallskip

Consequently, $f^{(\la),u_1,u_2}$ and $\widetilde{L}^{(\la)} \big(f^{(\la),u_1,u_2}_s,s; R^{(\la),u_1,u_2}_{s,.},S^{(\la),u_1,u_2}_{s,.},Z_{s,.}\big)$ satisfy the same equation \eqref{meme-eq}. By uniqueness of the solution, \eqref{egalite} follows. 
\end{proof}

\


\subsection{Fixed-point}

 Our aim here is to establish an a priori bound on the solution $ (\widetilde{\ell},\widetilde{m})$ of \eqref{systlm0}, which, combined with identity \eqref{egalite}, will play a fundamental role in the proof of Theorem \ref{thm-uniqueness}.  

\smallskip

We fix $s\geq 0$ and $f\in L^\infty$ for the whole section.

\subsubsection{A first simplification}

Consider the pair 
$$ (\ell,m)=(\ell^{(\la)}_{s,.},m^{(\la)}_{s,.})$$ 
which satisfies the following system: for every $t\geq s$,
 \begin{equation} \label{systlm}
\left\{
\begin{aligned}
&  \ell_t  =-3 \la \int_s^t dr\, e^{-(t-r)H}\Big(  (\ell_r+m_r) \pl  \<Psi2>_{s,r}\Big), \\
& m_t  =e^{-(t-s)H} f+\int_s^t dr\, e^{-(t-r)H}\Big(\mathfrak{G}^{(\la)}\big(\ell,m;R^{(\la),u_1,u_2}_{s,.},Z_{s,.}\big)_{s,r}\Big),
\end{aligned}
\right.
\end{equation}
with (recall the expression in \eqref{decc})
\begin{multline*}
\mathfrak{G}^{(\la)}(\ell,m;R^{(\la),u_1,u_2}_{s,.},Z_{s,.})_{s,t}=\widetilde{\mathfrak{G}}^{(\la)}(\ell,m;R^{(\la),u_1,u_2}_{s,.},0,Z_{s,.})_{s,t}\\
=\la\, R^{(\la),u_1,u_2}_{s,t} \cdot (\ell_t+m_t)-3\la (\ell_t+m_t) \pg  \<Psi2>_{s,t}-3\la\,  m_t\pe  \<Psi2>_{s,t}\\
\hspace{5cm}+9\la^2\,  \mathrm{com}\big(\ell,m;Z_{s,.}\big)_{s,t} +9\la^2\, \big(\ell_t+m_t\big)\, \<Psi2IPsi2>_{s,t} .
\end{multline*}

\smallskip

We also set
$$L^{(\la)}(f,s;R^{(\la),u_1,u_2}_{s,.},Z_{s,.})_{t}:=\ell_t+m_t.$$

The next result (inspired by \cite[Lemma 3.27]{DHYZ25}) will allow us to reduce the study of the system~\eqref{systlm0} to that of \eqref{systlm}.

\begin{proposition}\label{PropS+}
For all $s\leq t<s+1$, it holds that 
\begin{equation*}
\big\|\widetilde{L}^{(\la)} \big(f^{(\la),u_1,u_2}_s,s; R^{(\la),u_1,u_2}_{s,.},S^{(\la),u_1,u_2}_{s,.},Z_{s,.}\big)_t \big\|_{L^\infty} \leq \big\|L^{(\la)} \big(|f^{(\la),u_1,u_2}_s|,s;R^{(\la),u_1,u_2}_{s,.},Z_{s,.}\big)_t \big\|_{L^\infty}.
\end{equation*}
\end{proposition}

\begin{proof}
This is a straightforward consequence of the Feynman-Kac representation theorem (see {\it e.g.} \cite[Theorem 8.2.1, page 145]{Oksendal}), combined with the fact that $S^{(\la),u_1,u_2}_{s,r}\geq 0$ (see~\eqref{defi-s-la}). Observe first that for any $f\in L^\infty$, 
$$\tilde{d}:=\widetilde{L}^{(\la)} \big(f^{(\la),u_1,u_2}_s,s; R^{(\la),u_1,u_2}_{s,.},S^{(\la),u_1,u_2}_{s,.},Z_{s,.}\big)$$
satisfies the parabolic-type problem: $\tilde{d}_s=f$ and for all $t\geq s$,
$$(\partial_t-\Delta)\tilde{d}_t=\widetilde{V}_{s,t} \cdot \tilde{d}_t,$$
where $\widetilde{V}_{s,t}(x):=-|x|^2-3\la \<Psi2>_{s,t}(x)-9\la^2 \frakc^{\mathbf{2}}_{s,t}(x)+\la R^{(\la),u_1,u_2}_{s,t}(x)-\la S^{(\la),u_1,u_2}_{s,t}(x)$. Therefore, $\tilde{d}$ can be explicitly represented as
$$\tilde{d}_t(x):=\mathbb{E}_B \Big[ f\big(B_t\big) \exp\Big(\int_s^t \widetilde{V}_{s,t-r}(x+B_r) \, dr\Big)\Big],$$
for some Brownian motion $B$ independent from the underlying noise $W$. Since $S^{(\la),u_1,u_2}_{s,r}\geq 0$, we immediately deduce that
\begin{equation}\label{FKf}
\big|\tilde{d}_t(x)\big|\leq \mathbb{E}_B \Big[ \big|f\big(B_t\big)\big| \exp\Big(\int_s^t V_{s,t-r}(x+B_r) \, dr\Big)\Big],
\end{equation}
where $V_{s,t}(x):=\widetilde{V}_{s,t}(x)+\la S^{(\la),u_1,u_2}_{s,.}(x)$. It remains only to observe that, according to the Feynman-Kac formula again, the right-hand side of \eqref{FKf} coincides with the representation of
$$L^{(\la)} \big(|f^{(\la),u_1,u_2}_s|,s;R^{(\la),u_1,u_2}_{s,.},Z_{s,.}\big)_t(x).$$
\end{proof}

\begin{remark}
The elimination of the term $S^{(\la),u_1,u_2}$ between \eqref{systlm0} and \eqref{systlm} is particularly important from a technical standpoint, as it allows us to avoid potential non-integrable singularity issues (in time) related to the presence of squares in $S$: for instance, using Proposition \ref{cor:comes-down}, we only obtain the (non-integrable) bound $\big\|\big(\widetilde{U}^{(\la),u_1}_{s,t}\big)^2 \big\|_{L^\infty} \lesssim |t-s|^{-1-\eps}$, which would pose a serious difficulty for the subsequent fixed-point argument.
As mentioned earlier, this simplification had already been observed in \cite[Section 3.4]{DHYZ25}.
\end{remark}

\color{black}

\

\subsubsection{Topology and main local result}

\

\smallskip

For any given $\eps>0$, we fix four (finite) constants $C_0^{(\eps)},C_1,C_2^{(\eps)},C_3^{(\eps)}>0$ such that the following inequalities are satisfied for all (say smooth) functions $g$ and all $t\geq 0$:
\begin{equation}\label{c-0-c-1}
\big\|e^{-tH}g\big\|_{\cb^{1+4\eps}_{\infty}}\leq \frac{C_0^{(\eps)}}{t^{\frac12+2\eps}} \big\|g\big\|_{\cb^0_{\infty}}\quad , \quad \big\|g\big\|_{\cb^0_{\infty}} \leq C_1 \big\|g\big\|_{L^\infty},
\end{equation}
\begin{equation}
\big\|\big[e^{-tH}-\id\big]g\big\|_{\cb^\eps_{\infty}}\leq C_2^{(\eps)} t^{3\eps}\big\|g\big\|_{\cb^{7\eps}_{\infty}}\quad , \quad \big\|e^{-tH}g\big\|_{\cb^{7\eps}_{\infty}}\leq \frac{C_3^{(\eps)}}{t^{\frac{7\eps}{2}}} \big\|g\big\|_{\cb^0_{\infty}}.\label{c-2-c-3}
\end{equation}

\

Once endowed with these quantities, and for every $T\in (s,s+1]$, we define the space $X^{(\eps)}(s,T)$
 through the norm
\begin{align*}
&\big\|(\ell,m)\big\|_{X^{(\eps)}(s,T)}\\
&:=\max\bigg(\sup_{s\leq t\leq T} e^{3(t-s)}\big(\big\|\ell_t\big\|_{L^\infty}+\big\|m_t\big\|_{L^\infty}\big),\sup_{s\leq t\leq T} \big\|\ell_t\big\|_{\cb^{\frac12+2\eps}_{\infty}}, \sup_{s\leq t_1<t_2\leq T} \frac{\big\|\ell_{t_2}-\ell_{t_1}\big\|_{\cb^\eps_{\infty}}}{|t_2-t_1|^{3\eps}}, \\
&\hspace{1.5cm}\frac{1}{2C_0^{(\eps)}C_1} \sup_{s< t\leq T} |t-s|^{\frac23}\big\|m_t\big\|_{\cb^{1+4\eps}_{\infty}},\frac{1}{2C_1C_2^{(\eps)}C_3^{(\eps)}} \sup_{s< t_1<t_2\leq T} |t_1-s|^{\frac{7 \eps}2}\frac{\big\|m_{t_2}-m_{t_1}\big\|_{\cb^\eps_{\infty}}}{|t_2-t_1|^{3\eps}}\bigg).
\end{align*}

\

\

Recall the definition \eqref{notation-mst} of $M_{s,t}$. For all $\eta\geq 1$, we define
\begin{equation}\label{defTs}
T_s^{(\eta)}:=\inf\big\{t\geq s: \, M_{s,t}\geq \eta \big\} \wedge (s+1).
\end{equation}

\

\begin{proposition}\label{prop:l-m}
Fix $\eta\geq 1$. For every $\eps>0$ small enough, there exist constants $c_\eps, \kappa_\eps>0$ such that for all $\la\in(0,1)$,
\begin{equation*}
\big\|(\ell,m)\big\|_{X^{(\eps)}(s,T_s^{(\eta)})} \leq 9\big\|f\big\|_{L^\infty}+c_\eps \eta^{\ka_\eps} {\la^{\frac12}} \big\|(\ell,m)\big\|_{X^{(\eps)}(s,T_s^{(\eta)})}.
\end{equation*}
\end{proposition}

\

Assuming Proposition \ref{prop:l-m}, we deduce the following result:

\begin{corollary}\label{coro:la-star}
For all fixed $\eta\geq 1$ and $\eps>0$ small enough, set
\begin{equation}\label{la-star}
 \la_\star:= {\bigg(\frac{1}{10 c_\eps \eta^{\ka_\eps}}\bigg)^2}.
\end{equation}
Then for all $\la\in (0,\la_\star)$ and $t\in [s,T_s^{(\eta)}]$, one has
\begin{equation}\label{borne-expl}
\big\|\ell_t+m_t\big\|_{L^\infty} \leq 10\, e^{-3(t-s)}\big\|f\big\|_{L^\infty}.
\end{equation}
\end{corollary}

\

 We will see at the end of the proof of Theorem \ref{thm-uniqueness} that the explicit constants in \eqref{borne-expl} will play a key role (see Section \ref{section-concl}). Combined with the results of Propositions \ref{prop:egalite} and \ref{PropS+}, Corollary \ref{coro:la-star} then implies:

\begin{corollary}\label{coro:la-star2}
Let  $\la_\star>0$ be given by \eqref{la-star}. Then for  all  $u_1,u_2 \in \cb_{\infty}^{-\frac58}$, $\la\in (0,\la_\star)$,  $s\geq 0$,  $t\in [s,T_s^{(\eta)}]$
\begin{equation*} 
\big\|\big(X^{(\la),u_1}-X^{(\la),u_2}\big)_{t}\big\|_{L^\infty}\leq 10\, e^{-3(t-s)} \big\|\big(X^{(\la),u_1}-X^{(\la),u_2}\big)_{s}\big\|_{L^\infty}.
\end{equation*}

\end{corollary}

\

\subsection{Proof of Proposition \ref{prop:l-m}}\label{subsec:proofprolm}

\

\smallskip

This section is devoted to the proof of Proposition \ref{prop:l-m}. For clarity, we set throughout the proof $T:=T_s^{(\eta)}$ and
$$d_t:=\ell_t+m_t.$$

\subsubsection{Various bounds}

Note first that 
\begin{align*}
&\big\|d_t\big\|_{L^\infty} \lesssim \big\|(\ell,m)\big\|_{X^{(\eps)}(s,T)}, 
\end{align*}
while 
\begin{align}
\big\|d_t\big\|_{\cb^{\frac12+2\eps}_{\infty}}
&\lesssim\big\|\ell_t\big\|_{\cb^{\frac12+2\eps}_{\infty}}+\big\|m_t\big\|_{\cb^{0}_{\infty}}^{\frac12}\big\|m_t\big\|_{\cb^{1+4\eps}_{\infty}}^{\frac12}\lesssim \frac{1}{|t-s|^{\frac13 }} \big\|(\ell,m)\big\|_{X^{(\eps)}(s,T)}.\label{control-d-2}
\end{align}

\begin{lemma}\label{lem:com1}
For every $\varepsilon>0$ small enough, $Z\in  \mathcal Z_{\eps,T}$ and $(v,w) \in X^{(\eps)}(s,T)$, one has
\begin{align}\label{desi1}
{\big\|\mathrm{com}_1(\ell,m;Z_{s,.})_{s,t}\big\|_{\mathcal{B}_{\infty}^{1+2\varepsilon}}} & \lesssim \frac{1}{|t-s|^{\frac{1}{12}+\eps}}\eta \big\|(\ell,m)\big\|_{X^{(\eps)}(s,T)}
\end{align}
and 
\begin{align}\label{desi2}
\Big\|\big[ \pl, \pe\big] \Big(\ell_r+m_r, \<IPsi2>_{s,r}, \<Psi2>_{s,r}\Big) \Big\|_{\cb^\eps_{\infty}}\lesssim \frac{1}{|r-s|^{\frac{1}{3}}}\eta^2 \big\|(\ell,m)\big\|_{X^{(\eps)}(s,T)}.
 \end{align}
\end{lemma}

\begin{proof}  Recall that we have set $d_t=\ell_t+m_t$. We can rephrase the definition of $\mathrm{com}_1(\ell,m;Z_{s,.})$ as
\begin{align}
\mathrm{com}_1(\ell,m;Z_{s,.})_{s,t}&=\int_s^t dr\, e^{-(t-r)H} \big[d_r \pl  \<Psi2>_{s,r}\big]  -d_t \pl  \<IPsi2>_{s,t}\nonumber\\
&=\int_s^t dr \,\big[ e^{-(t-r)H},\pl \big]\big(d_r,\<Psi2>_{s,r}\big)-\int_s^t dr \, (d_t-d_r) \pl \big(e^{-(t-r)H}\<Psi2>_{s,r}\big).\label{decomp-ci}
\end{align}

Then, using the result of \cite[Lemma 13.26]{DFT} and \eqref{control-d-2}, we obtain on the one hand that for $\eps>0$ small enough,
\begin{align*}
\Big\|\int_s^t dr \,\big[ e^{-(t-r)H},\pl \big]\big(d_r,\<Psi2>_{s,r}\big) \Big\|_{ \B_\infty^{1+2\eps}}&\lesssim \int_s^t \frac{dr}{|t-r|^{\frac34+\eps}}\big\| d_r \big\|_{\B_\infty^{\frac12+2\eps}}\big\| \<Psi2>_{s,r} \big\|_{\B_\infty^{-1-\eps}}\\
&\lesssim \eta \big\|(\ell,m)\big\|_{X^{(\eps)}(s,T)}\int_s^t \frac{dr}{|t-r|^{\frac34+\eps}|r-s|^{\frac13}}\lesssim \frac{\eta \big\|(\ell,m)\big\|_{X^{(\eps)}(s,T)}}{|t-s|^{\frac{1}{12}+\eps}}.
\end{align*}
On the other hand,
\begin{align*}
\Big\|\int_s^t dr \, (d_t-d_r) \pl \big(e^{-(t-r)H}\<Psi2>_{s,r}\big)\Big\|_{ \B_\infty^{1+2\eps}}&\lesssim \int_s^t dr \, \big\|d_t-d_r\big\|_{L^\infty} \big\|e^{-(t-r)H}\<Psi2>_{s,r}\big\|_{ \B_\infty^{1+2\eps}}\\
& \lesssim  \big\|(\ell,m)\big\|_{X^{(\eps)}(s,T)}\int_s^t \frac{dr}{|t-r|^{1+\frac{3\eps}2} } \, \frac{|t-r|^{3\eps}}{|r-s|^{\frac{7\eps}2}} \big\|\<Psi2>_{s,r}\big\|_{ \B_\infty^{-1-\eps}}\\
& \lesssim \eta \big\|(\ell,m)\big\|_{X^{(\eps)}(s,T)}\int_s^t \frac{dr}{|t-r|^{1-\frac{3\eps}2}   |r-s|^{\frac{7\eps}2}   }\\
& \lesssim  \frac{\eta \big\|(\ell,m)\big\|_{X^{(\eps)}(s,T)}}{|t-s|^{2\eps}}.
\end{align*}
Inserting the above bounds into \eqref{decomp-ci} yields the desired estimate \eqref{desi1}. 

\

As for \eqref{desi2}, one has by \cite[Proposition 13.25]{DFT} and \eqref{control-d-2}
\begin{align*} 
  \Big\|\big[ \pl, \pe\big] \Big(\ell_r+m_r, \<IPsi2>_{s,r}, \<Psi2>_{s,r}\Big) \Big\|_{\cb^\eps_{\infty}} &\lesssim  \big\|\ell_r+m_r\big\|_{\cb^{4\eps}_\infty}  \big\|\<IPsi2>_{s,r}\big\|_{ \B_\infty^{1-\eps}} \big\|\<Psi2>_{s,r}\big\|_{ \B_\infty^{-1-\eps}} \\
 &\lesssim \frac{1}{|r-s|^{\frac{1}{3}}}\eta^2 \big\|(\ell,m)\big\|_{X^{(\eps)}(s,T)},
 \end{align*}
 which is exactly the claim.
\end{proof}

For further use, and with expression \eqref{defi-r-la} of $R^{(\la),u_1,u_2}_{s,.}$ in mind, we decompose the quantity $\mathfrak{G}^{(\la)}(\ell,m;R^{(\la),u_1,u_2}_{s,.},Z_{s,.})_{s,r}$ as
\begin{align}
&\mathfrak{G}^{(\la)}(\ell,m;R^{(\la),u_1,u_2}_{s,.},Z_{s,.})_{s,r}=\nonumber\\
&=-3 \la\, d_r\cdot \Big[ \<Psi>_{s,r}\big(U^{(\la),u_1}_{s,r}+U^{(\la),u_2}_{s,r}\big)-2 \la \<PsiIPsi3nr>_{s,r}- {\overline{\frakc}^{\mathbf{1}}_{s,r}+3\la \overline{\frakc}^{\mathbf{2}}_{s,r}}  \Big]\nonumber\\
&\hspace{2cm}-3\la\,  d_r \pg  \<Psi2>_{s,r}-3\la\,  m_r\pe  \<Psi2>_{s,r}+9\la^2\,  \mathrm{com}\big(\ell,m;Z_{s,.}\big)_{s,r} +9\la^2\, d_r\, \<Psi2IPsi2>_{s,r} \nonumber\\
&=-3 \la\,\big[ d_r\big(U^{(\la),u_1}_{s,r}+U^{(\la),u_2}_{s,r}\big) \big]\pl \<Psi>_{s,r}-3 \la\,\big[ d_r\big(U^{(\la),u_1}_{s,r}+U^{(\la),u_2}_{s,r}\big) \big]\pge \<Psi>_{s,r}\nonumber\\
&\hspace{1cm}+3 \la\, d_r\cdot \big[{\overline{\frakc}^{\mathbf{1}}_{s,r}-3\la \overline{\frakc}^{\mathbf{2}}_{s,r}}\big]-3\la\,  m_r\pe  \<Psi2>_{s,r}\nonumber\\
&\hspace{0.5cm}+9\la^2\,  \mathrm{com}\big(\ell,m;Z_{s,.}\big)_{s,r} -3\la\,  d_r \pg  \<Psi2>_{s,r}+3 \la^2\, d_r\cdot \Big[2\<PsiIPsi3nr>_{s,r} +3 \<Psi2IPsi2>_{s,r}\Big] =:\sum_{i=1,\ldots,7} \ci^{i,(\la)}_{s,r}.\label{decompo-sigma}
\end{align}

\

\begin{lemma}\label{lem:Ij}
Let $\ci^{i,(\la)}_{s,r}$ be defined along \eqref{decompo-sigma}. Then, for all {$s \leq t_1 <t_2 \leq T$} and all $1 \leq i \leq 7$ we have
\begin{equation}\label{normifJ}
{\int_{t_1}^{t_2}} dr \, \Big\|e^{-({t_2}-r)H}\Big(\ci^{i,(\la)}_{s,r}\Big)\Big\|_{\cb^{7\eps}_\infty} \lesssim  {\la^{\frac12}} \eta^\ka  \big\|(\ell,m)\big\|_{X^{(\eps)}(s,T)},
\end{equation}
\begin{equation}\label{normifJ2}
{\int_{t_1}^{t_2}} dr \, \Big\|e^{-({t_2}-r)H}\Big(\ci^{i,(\la)}_{s,r}\Big)\Big\|_{\cb^{\eps}_\infty} \lesssim  {\la^{\frac12}}  \eta^\ka  \big\|(\ell,m)\big\|_{X^{(\eps)}(s,T)}|{t_2-t_1}|^{\frac16-2\eps},
\end{equation}
and 
\begin{equation}\label{borneI4eps}
{\int_{s}^{t_2}} dr \, \Big\|e^{-({t_2}-r)H}\Big(\ci^{i,(\la)}_{s,r}\Big)\Big\|_{\cb^{1+4\eps}_\infty} \lesssim  {\la^{\frac12}}  \eta^\ka  \big\|(\ell,m)\big\|_{X^{(\eps)}(s,T)}{\frac{1}{|{t_2}-s|^{\frac13+4\eps}}}.
\end{equation}
\end{lemma}

\smallskip

\begin{proof}
Let $0 \leq \sigma \leq 1+4 \eps$.

\medskip

\noindent \underline{Case of $\ci^{1,(\la)}_{s,r}$}.   By Proposition \ref{cor:comes-down}, one has
\begin{align*}
{\int_{t_1}^{t_2}} dr \, \Big\|e^{-({t_2}-r)H}\Big(\ci^{1,(\la)}_{s,r}\Big)\Big\|_{\cb^\sigma_{\infty}}&\lesssim \la  {\int_{t_1}^{t_2}}dr \, \Big\|e^{-({t_2}-r)H}\Big(\big[ d_r\big( U^{(\la),u_1}_{s,r}+ U^{(\la),u_2}_{s,r}\big) \big]\pl \<Psi>_{s,r}\Big)\Big\|_{\cb^\sigma_{\infty}}\\
&\lesssim \la {\int_{t_1}^{t_2}}\frac{dr}{|{t_2}-r|^{\frac14+\frac{\sigma}2+\frac{\eps}2}} \, \Big\|\big[ d_r\big( U^{(\la),u_1}_{s,r}+ U^{(\la),u_2}_{s,r}\big) \big]\pl \<Psi>_{s,r}\Big\|_{\cb^{-\frac12-\eps}_\infty}\\
&\lesssim \la  \eta{\int_{t_1}^{t_2}} \frac{dr}{|{t_2}-r|^{\frac14+\frac{\sigma}2+\frac{\eps}2}} \, \big\|d_r \big\|_{L^\infty} \big\| U^{(\la),u_1}_{s,r}+ U^{(\la),u_2}_{s,r} \big\|_{L^\infty}\\
&\lesssim  {\la^{\frac12}}\eta^\ka  \big\|(\ell,m)\big\|_{X^{(\eps)}(s,T)}{\int_{t_1}^{t_2}} \frac{dr}{ |{t_2}-r|^{\frac14+\frac{\sigma}2+\frac{\eps}2}|r-s|^{\frac12+\eps}}.
\end{align*}
By taking the different values $\sigma \in \{\eps, 7\eps, 1+4\eps\}$ we obtain the announced bounds.
\medskip

\noindent
\underline{Case of $\ci^{2,(\la)}_{s,r}$}.  By Proposition \ref{cor:comes-down} and \eqref{control-d-2}, one has
\begin{multline*}
{\int_{t_1}^{t_2}} dr \, \Big\|e^{-({t_2}-r)H}\Big(\ci^{2,(\la)}_{s,r}\Big)\Big\|_{\cb^\sigma_{\infty}}\lesssim\\
\begin{aligned}
& \lesssim \la {\int_{t_1}^{t_2}} \frac{dr}{|{t_2}-r|^{\frac{\sigma}2-\frac{\eps}2}} \, \Big\|e^{-(t_2-r)H}\Big(\big[ d_r\big( U^{(\la),u_1}_{s,r}+ U^{(\la),u_2}_{s,r}\big) \big]\pge \<Psi>_{s,r}\Big)\Big\|_{\cb^\eps_{\infty}}\\
& \lesssim \la {\int_{t_1}^{t_2}} \frac{dr}{|{t_2}-r|^{\frac{\sigma}2-\frac{\eps}2}} \, \Big\| \big[ d_r\big( U^{(\la),u_1}_{s,r}+ U^{(\la),u_2}_{s,r}\big) \big]\pge \<Psi>_{s,r}\Big\|_{\cb^\eps_{\infty}}\\
&\lesssim \la \eta{\int_{t_1}^{t_2}} \frac{dr}{|{t_2}-r|^{\frac{\sigma}2-\frac{\eps}2}}\, \Big(\big\| d_r \big\|_{L^\infty}\big\|U^{(\la),u_1}_{s,r}+ U^{(\la),u_2}_{s,r} \big\|_{\cb^{\frac12+2\eps}_\infty}+\big\| d_r \big\|_{\cb^{\frac12+2\eps}_\infty}\big\|U^{(\la),u_1}_{s,r}+ U^{(\la),u_2}_{s,r} \big\|_{L^\infty}\Big)\\
&\lesssim {\la^{\frac12}} \eta^\ka\big\|(\ell,m)\big\|_{X^{(\eps)}(s,T)}{\int_{t_1}^{t_2}}  \frac{dr}{|{t_2}-r|^{\frac\sigma2-\frac\eps2}|r-s|^{\frac56+2\eps}}.
\end{aligned}
\end{multline*}
We now take  the different values $\sigma \in \{\eps, 7\eps, 1+4\eps\}$ and obtain the bounds.

\medskip

\noindent
\underline{Case of $\ci^{3,(\la)}_{s,r}$}. By Lemma \ref{lem:c}, one has
\begin{align*}
{\int_{t_1}^{t_2}} dr \, \Big\|e^{-(t_2-r)H}\Big(\ci^{3,(\la)}_{s,r}\Big)\Big\|_{\cb^{\sigma}_\infty}&\lesssim \la {\int_{t_1}^{t_2}} dr \, \Big\|e^{-({t_2}-r)H}\Big(d_r\cdot \big[ {\overline{\frakc}^{\mathbf{1}}_{s,r}-3\la \overline{\frakc}^{\mathbf{2}}_{s,r}} \big]\Big)\Big\|_{\cb^{\sigma}_\infty}\\
&\lesssim \la {\int_{t_1}^{t_2}} \frac{dr}{|{t_2}-r|^{\frac{\sigma}2}} \, \Big\|d_r\cdot \big[ {\overline{\frakc}^{\mathbf{1}}_{s,r}-3\la \overline{\frakc}^{\mathbf{2}}_{s,r}} \big]\Big\|_{\cb^{0}_\infty}\\
&\lesssim \la {\int_{t_1}^{t_2}} \frac{dr}{|{t_2}-r|^{\frac\sigma2}}\, \big\|d_r\big\|_{L^\infty}  \big\|{\overline{\frakc}^{\mathbf{1}}_{s,r}-3\la \overline{\frakc}^{\mathbf{2}}_{s,r}}\big\|_{L^\infty}\\
&\lesssim \la {\Big( \big\llbracket  \overline{\frakc}^{\mathbf{1}}_{s,.}\big\rrbracket_{3\eps,[s,T]}+ \big\llbracket  \overline{\frakc}^{\mathbf{2}}_{s,.}\big\rrbracket_{3\eps,[s,T]}\Big)} {\int_{t_1}^{t_2}} \frac{dr}{|{t_2}-r|^{\frac\sigma2}|r-s|^{{\frac{9}{16}}}}\big\|d_r\big\|_{L^\infty}\\
&\lesssim \la {\eta}\big\|(\ell,m)\big\|_{X^{(\eps)}(s,T)} {\int_{t_1}^{t_2}} \frac{dr}{|{t_2}-r|^{\frac\sigma2}|r-s|^{{\frac{9}{16}}}}.
\end{align*}
We are now able to conclude by taking  $\sigma \in \{\eps, 7\eps, 1+4\eps\}$.

\medskip

\noindent
\underline{Case of $\ci^{4,(\la)}_{s,r}$}. By Proposition \ref{Prop-est-para} $(i)$, one has
\begin{align*}
{\int_{t_1}^{t_2}} dr \, \Big\|e^{-({t_2}-r)H}\Big(\ci^{4,(\la)}_{s,r}\Big)\Big\|_{\cb^{\sigma}_\infty}&\lesssim \la {\int_{t_1}^{t_2}} dr \, \Big\|e^{-({t_2}-r)H}\big(m_r\pe  \<Psi2>_{s,r}\big)\Big\|_{\cb^{\sigma}_\infty}\\
&\lesssim \la {\int_{t_1}^{t_2}} \frac{dr}{|{t_2}-r|^{\frac\sigma2-\frac{\eps}2}} \, \big\|m_r\pe  \<Psi2>_{s,r}\big\|_{\cb^\eps_{\infty}}\\
&\lesssim \la \eta {\int_{t_1}^{t_2}} \frac{dr}{|{t_2}-r|^{\frac\sigma2-\frac{\eps}2}} \, \big\|m_r\big\|_{\cb^{1+4\eps}_\infty}\\
&\lesssim \la \eta \big\|(\ell,m)\big\|_{X^{(\eps)}(s,T)}{\int_{t_1}^{t_2}} \frac{dr}{|{t_2}-r|^{\frac\sigma2-\frac{\eps}2}|r-s|^{\frac23}}.
\end{align*}
We now take  the different values $\sigma \in \{\eps, 7\eps, 1+4\eps\}$ and obtain the bounds.

\medskip

\noindent
\underline{Case of $\ci^{5,(\la)}_{s,r}$}. By Lemma \ref{lem:com1}, one has
\begin{align*}
&{\int_{t_1}^{t_2}} dr \, \Big\|e^{-({t_2}-r)H}\Big(\ci^{5,(\la)}_{s,r}\Big)\Big\|_{\cb^{\sigma}_\infty}\lesssim \\
&\lesssim \la {\int_{t_1}^{t_2}} dr \, \Big\|e^{-({t_2}-r)H}\Big(\mathrm{com}_1(\ell,m;Z_{s,.})_{s,r} \pe  \<Psi2>_{s,r} \Big)\Big\|_{\cb^{\sigma}_\infty}+\la {\int_{t_1}^{t_2}} dr \, \Big\|e^{-({t_2}-r)H}\Big(\big[ \pl, \pe\big] \Big(d_r, \<IPsi2>_{s,r}, \<Psi2>_{s,r}\Big) \Big)\Big\|_{\cb^{\sigma}_\infty}\\
&\lesssim \la {\int_{t_1}^{t_2}}\frac{dr}{|{t_2}-r|^{\frac{\sigma}2-\frac{\eps}2}} \, \Big\|\mathrm{com}_1(\ell,m;Z_{s,.})_{s,r} \pe  \<Psi2>_{s,r} \Big\|_{\cb^\eps_{\infty}}+\la {\int_{t_1}^{t_2}}\frac{dr}{|{t_2}-r|^{\frac\sigma2-\frac\eps2}} \, \Big\|\big[ \pl, \pe\big] \Big(d_r, \<IPsi2>_{s,r}, \<Psi2>_{s,r}\Big) \Big\|_{\cb^\eps_{\infty}}\\
&\lesssim \la \eta{\int_{t_1}^{t_2}} \frac{dr}{|{t_2}-r|^{\frac{\sigma}2-\frac{\eps}2}} \, \big\|\mathrm{com}_1(\ell,m;Z_{s,.})_{s,r}  \Big\|_{\cb^{1+2\eps}_\infty}+\la \eta^2 \big\|(\ell,m)\big\|_{X^{(\eps)}(s,T)}{\int_{t_1}^{t_2}} \frac{dr}{|{t_2}-r|^{\frac{\sigma}2-\frac{\eps}2}|r-s|^{\frac13}}\\
&\lesssim \la \eta^2 \big\|(\ell,m)\big\|_{X^{(\eps)}(s,T)}{\int_{t_1}^{t_2}} \frac{dr}{|{t_2}-r|^{\frac{\sigma}2-\frac{\eps}2}|r-s|^{\frac{1}{12}+\eps}}+\la \eta^2 \big\|(\ell,m)\big\|_{X^{(\eps)}(s,T)}{\int_{t_1}^{t_2}} \frac{dr}{|{t_2}-r|^{\frac{\sigma}2-\frac{\eps}2}|r-s|^{\frac13}}.
\end{align*}
We now take  the different values $\sigma \in \{\eps, 7\eps, 1+4\eps\}$ and obtain the bounds.

\medskip

\noindent
\underline{Case of $\ci^{6,(\la)}_{s,r}$}. Thanks to \eqref{control-d-2}, one has
\begin{align*}
{\int_{t_1}^{t_2}} dr \, \Big\|e^{-({t_2}-r)H}\Big(\ci^{6,(\la)}_{s,r}\Big)\Big\|_{\cb^{\sigma}_\infty}&\lesssim \la {\int_{t_1}^{t_2}} dr \, \Big\|e^{-({t_2}-r)H}\big(d_r \pg  \<Psi2>_{s,r}\big)\Big\|_{\cb^{\sigma}_\infty}\\
&\lesssim \la {\int_{t_1}^{t_2}} \frac{dr}{|{t_2}-r|^{\frac\sigma2+\frac14-\frac\eps2}}\, \big\| d_r \pg  \<Psi2>_{s,r}\big\|_{\cb^{-\frac12+\eps}_\infty}\\
&\lesssim \la \eta {\int_{t_1}^{t_2}} \frac{dr}{|{t_2}-r|^{\frac\sigma2+\frac14-\frac\eps2}}\, \big\| d_r \big\|_{\cb^{\frac12+2\eps}_\infty}\\
&\lesssim \la \eta \big\|(\ell,m)\big\|_{X^{(\eps)}(s,T)} {\int_{t_1}^{t_2}} \frac{dr}{|{t_2}-r|^{\frac\sigma2+\frac14-\frac\eps2}|r-s|^{\frac13}}.
\end{align*}
We are now able to conclude by taking  $\sigma \in \{\eps, 7\eps, 1+4\eps\}$.

\medskip

\noindent
\underline{Case of $\ci^{7,(\la)}_{s,r}$}. Thanks to \eqref{control-d-2}, one has
\begin{align*}
{\int_{t_1}^{t_2}} dr \, \Big\|e^{-({t_2}-r)H}\Big(\ci^{7,(\la)}_{s,r}\Big)\Big\|_{\cb^{\sigma}_\infty}&\lesssim \la {\int_{t_1}^{t_2}}dr \, \Big\|e^{-({t_2}-r)H}\Big(d_r\cdot \Big[2\<PsiIPsi3nr>_{s,r} +3 \<Psi2IPsi2>_{s,r}\Big]\Big)\Big\|_{\cb^{\sigma}_\infty}\\
&\lesssim \la {\int_{t_1}^{t_2}}\frac{dr}{|{t_2}-r|^{\frac\sigma2+\frac14+\frac\eps2}} \, \Big\|d_r\cdot \Big[2\<PsiIPsi3nr>_{s,r} +3 \<Psi2IPsi2>_{s,r}\Big]\Big\|_{\cb^{-\frac12-\eps}_\infty}\\
&\lesssim \la \eta{\int_{t_1}^{t_2}} \frac{dr}{|{t_2}-r|^{\frac\sigma2+\frac14+\frac\eps2}} \, \big\|d_r\big\|_{\cb^{\frac12+2\eps}_\infty}\\
&\lesssim \la \eta \big\|(\ell,m)\big\|_{X^{(\eps)}(s,T)} {\int_{t_1}^{t_2}} \frac{dr}{|{t_2}-r|^{\frac\sigma2+\frac14+\frac\eps2}|r-s|^{\frac13}}.
\end{align*}
By taking the different values $\sigma \in \{\eps, 7\eps, 1+4\eps\}$ we obtain the announced bounds.
\end{proof}

\

\subsubsection{Control of $\big\|\ell_t \big\|_{L^\infty}$ and  $ \big\|\ell_t\big\|_{\cb^{\frac12+2\eps}_\infty}$}

One has
\begin{align}
\big\|\ell_t\big\|_{L^\infty} \lesssim \big\|\ell_t\big\|_{\cb^{\frac12+2\eps}_\infty} &\lesssim \la \int_s^t dr \, \Big\|e^{-(t-r)H}\big(d_r\pl  \<Psi2>_{s,r}\big)\Big\|_{\cb^{\frac12+2\eps}_\infty}\nonumber \\
&\lesssim \la \int_s^t \frac{dr}{|t-r|^{\frac34+2\eps}} \, \big\|d_r\pl  \<Psi2>_{s,r}\big\|_{\cb^{-1-\eps}_\infty}\nonumber \\
&\lesssim \la \eta \int_s^t \frac{dr}{|t-r|^{\frac34+2\eps}} \, \big\|d_r\big\|_{L^\infty}\lesssim \la \eta \big\|(\ell,m)\big\|_{X^{(\eps)}(s,T)}. \label{est01}
\end{align}

\

\subsubsection{Control of $\big\|\ell_{t_2}-\ell_{t_1}\big\|_{\cb_\infty^\eps}$}

One has
\begin{multline*}
\big\|\ell_{t_2}-\ell_{t_1}\big\|_{\cb_\infty^\eps} \lesssim \la \int_{t_1}^{t_2} dr \, \Big\|e^{-(t_2-r)H}\big(d_r\pl  \<Psi2>_{s,r}\big)\Big\|_{\cb^\eps_{\infty}}+\\
+\la \int_s^{t_1} dr \, \Big\|\big[e^{-(t_2-r)H}-e^{-(t_1-r)H}\big]\big(d_r\pl  \<Psi2>_{s,r}\big)\Big\|_{\cb^\eps_{\infty}}.
\end{multline*}

On the one hand,
\begin{align*}
  \int_{t_1}^{t_2} dr \, \Big\|e^{-(t_2-r)H}\big(d_r\pl  \<Psi2>_{s,r}\big)\Big\|_{\cb^\eps_{\infty}}&\lesssim  \int_{t_1}^{t_2} \frac{dr}{|t_2-r|^{\frac12+\eps}} \, \big\|d_r\pl  \<Psi2>_{s,r}\big\|_{\cb^{-1-\eps}_\infty}\\
&\lesssim  \eta \int_{t_1}^{t_2} \frac{dr}{|t_2-r|^{\frac12+\eps}} \, \big\|d_r\big\|_{L^\infty}\\
&\lesssim  \eta \big\|(\ell,m)\big\|_{X^{(\eps)}(s,T)}|t_2-t_1|^{\frac12-\eps}.
\end{align*}
On the other hand,
\begin{multline*}
  \int_s^{t_1} dr \,\Big\|\big[e^{-(t_2-t_1)H}-\id\big]\Big(e^{-(t_1-r)H}\big(d_r\pl  \<Psi2>_{s,r}\big)\Big)\Big\|_{\cb^\eps_{\infty}} \lesssim  |t_2-t_1|^{3\eps}\int_s^{t_1} dr \, \big\|e^{-(t_1-r)H}\big(d_r\pl  \<Psi2>_{s,r}\big)\big\|_{\cb^{7\eps}_\infty}\\
\begin{aligned}
&\lesssim  |t_2-t_1|^{3\eps} \int_s^{t_1}  \frac{dr}{|t_1-r|^{\frac12+4\eps}} \, \big\|d_r\pl  \<Psi2>_{s,r}\big\|_{\cb^{-1-\eps}_\infty}\\
&\lesssim  \eta |t_2-t_1|^{3\eps} \int_s^{t_1}  \frac{dr}{|t_1-r|^{\frac12+4\eps}} \, \big\|d_r\big\|_{L^\infty}\lesssim  \eta  \big\|(\ell,m)\big\|_{X^{(\eps)}(s,T)}|t_2-t_1|^{3\eps}.
\end{aligned}
\end{multline*}
The previous inequalities imply that for all $s\leq t_1<t_2\leq T$
\begin{equation} \label{est02}
\frac{\big\|\ell_{t_2}-\ell_{t_1}\big\|_{\cb^\eps_{\infty}}}{|t_2-t_1|^{3\eps}} \lesssim  \eta  \big\|(\ell,m)\big\|_{X^{(\eps)}(s,T)}.
\end{equation}

\

\subsubsection{Control of $\big\|m_t \big\|_{L^\infty}$}

One has by Lemma \ref{lem:dep},
\begin{align}
\big\|m_t \big\|_{L^\infty}&\leq \big\|e^{-(t-s)H}f\big\|_{L^\infty}+\int_s^t dr \, \Big\|e^{-(t-r)H}\Big(\mathfrak{G}^{(\la)}\big(\ell,m;R^{(\la),u_1,u_2}_{s,.},Z_{s,.}\big)_{s,r}\Big)\Big\|_{L^\infty}\nonumber\\
&\leq 8 e^{-3(t-s)} \big\|f\big\|_{L^\infty}+\int_s^t dr \, \Big\|e^{-(t-r)H}\Big(\mathfrak{G}^{(\la)}\big(\ell,m;R^{(\la),u_1,u_2}_{s,.},Z_{s,.}\big)_{s,r}\Big)\Big\|_{L^\infty}.\label{l+m}
\end{align}
With decomposition \eqref{decompo-sigma} in mind, we can then write
\begin{align*}
\int_s^t dr \, \Big\|e^{-(t-r)H}\Big(\mathfrak{G}^{(\la)}\big(\ell,m;R^{(\la),u_1,u_2}_{s,.},Z_{s,.}\big)_{s,r}\Big)\Big\|_{L^\infty}&\leq \sum_{i=1,\ldots,7} \int_s^t dr \, \Big\|e^{-(t-r)H}\Big(\ci^{i,(\la)}_{s,r}\Big)\Big\|_{L^\infty} \\
&\lesssim  \sum_{i=1,\ldots,7} \int_s^t dr \, \Big\|e^{-(t-r)H}\Big(\ci^{i,(\la)}_{s,r}\Big)\Big\|_{\cb^{\eps}_\infty}.
\end{align*}
Now we can use \eqref{normifJ} together with \eqref{l+m} to deduce the existence of a constant $c_\eps>0$ such that for all $t\in [s,T]$,
\begin{align}\label{est03}
e^{3(t-s)}\big\|m_t \big\|_{L^\infty}\leq 8  \big\|f\big\|_{L^\infty}+c_\eps  {\la^{\frac12}} \eta^\ka \big\|(\ell,m)\big\|_{X^{(\eps)}(s,T)}.
\end{align}

\

\subsubsection{Control of $\big\|m_t\big\|_{\cb^{1+4\eps}_\infty}$}

One has
\begin{align}
&\big\|m_t\big\|_{\cb^{1+4\eps}_\infty}\leq \big\|e^{-(t-s)H}f\big\|_{\cb^{1+4\eps}_\infty}+\int_s^t dr \, \Big\|e^{-(t-r)H}\Big(\mathfrak{G}^{(\la)}\big(\ell,m; R^{(\la),u_1,u_2}_{s,.},Z_{s,.}\big)_{s,r}\Big)\Big\|_{\cb^{1+4\eps}_\infty}.\label{m-1+2eps}
\end{align}
On the one hand, recalling the definition of $C^{(\eps)}_0$ and $C_1$ in \eqref{c-0-c-1}, one has
\begin{equation*}
\big\|e^{-(t-s)H}f\big\|_{\cb^{1+4\eps}_\infty}\leq \frac{C_0^{(\eps)}}{(t-s)^{\frac12+2\eps}} \big\|f\big\|_{\cb^0_\infty} \leq \frac{C_0^{(\eps)}C_1}{(t-s)^{\frac12+2\eps}} \big\|f\big\|_{L^\infty}.
\end{equation*}
On the other hand, with the decomposition \eqref{decompo-sigma} in mind, we can write
\begin{align*}
\int_s^t dr \, \Big\|e^{-(t-r)H}\Big(\mathfrak{G}^{(\la)}\big(\ell,m;R^{(\la),u_1,u_2}_{s,.},Z_{s,.}\big)_{s,r}\Big)\Big\|_{\cb^{1+4\eps}_\infty}\leq \sum_{i=1,\ldots,7} \int_s^t dr \, \Big\|e^{-(t-r)H}\Big(\ci^{i,(\la)}_{s,r}\Big)\Big\|_{\cb^{1+4\eps}_\infty}
\end{align*}
and then use \eqref{borneI4eps}. Going back to \eqref{m-1+2eps}, we have shown the existence of a constant $c_\eps>0$ such that for all $t\in [s,T]$,
\begin{align}\label{est04}
\frac{1}{2C_0^{(\eps)}C_1}|t-s|^{\frac23}\big\|m_t\big\|_{\cb^{1+4\eps}_\infty}\leq \frac12  \big\|f\big\|_{L^\infty}+c_\eps  \la^{\frac12} \eta^\ka \big\|(\ell,m)\big\|_{X^{(\eps)}(s,T)}.
\end{align}

\

\subsubsection{Control of $\big\|m_{t_2}-m_{t_1}\big\|_{\cb^\eps_{\infty}}$}
One has
\begin{align*}
&\big\|m_{t_2}-m_{t_1}\big\|_{\cb^\eps_{\infty}}\leq \Big\|\big[e^{-(t_2-s)H}-e^{-(t_1-s)H}\big]f\Big\|_{\cb^\eps_{\infty}}+\int_{t_1}^{t_2} dr \, \Big\|e^{-(t_2-r)H}\Big(\mathfrak{G}^{(\la)}\big(\ell,m;R,Z_{s,.}\big)_{s,r}\Big)\Big\|_{\cb^\eps_{\infty}}\nonumber\\
&\hspace{3cm}+\int_s^{t_1} dr \, \Big\|\big[e^{-(t_2-r)H}-e^{-(t_1-r)H}\big]\Big(\mathfrak{G}^{(\la)}\big(\ell,m;R,Z_{s,.}\big)_{s,r}\Big)\Big\|_{\cb^\eps_{\infty}}.
\end{align*}
On the one hand, recalling the definition of $C_1$, $C_2^{(\eps)}$ and $C_3^{(\eps)}$ in \eqref{c-0-c-1}-\eqref{c-2-c-3}, one has
\begin{align*}
\Big\|\big[e^{-(t_2-s)H}-e^{-(t_1-s)H}\big]f\Big\|_{\cb^{\eps}_{\infty}}& =\Big\|\big[e^{-(t_2-t_1)H}-\id\big]e^{-(t_1-s)H}f\Big\|_{\cb^{\eps}_{\infty}} \\
&\leq C_2^{(\eps)} |t_2-t_1|^{3\eps} \big\|e^{-(t_1-s)H}f\big\|_{\cb^{7\eps}_{\infty}}\\
&\leq C_2^{(\eps)} C_3^{(\eps)}\frac{|t_2-t_1|^{3\eps} }{|t_1-s|^{\frac{7\eps}{2}}} \big\|f\big\|_{\cb^0_{\infty}}\leq C_1 C_2^{(\eps)} C_3^{(\eps)}\frac{|t_2-t_1|^{3\eps} }{|t_1-s|^{\frac{7\eps}{2}}} \big\|f\big\|_{L^\infty}.
\end{align*}
On the other hand, using the decomposition in \eqref{decompo-sigma}, we can write
\begin{align*}
&\int_{t_1}^{t_2} dr \, \Big\|e^{-(t_2-r)H}\Big(\mathfrak{G}^{(\la)}\big(\ell,m;R,Z_{s,.}\big)_{s,r}\Big)\Big\|_{\cb^\eps_{\infty}}\\
&\hspace{3cm}+\int_s^{t_1} dr \, \Big\|\big[e^{-(t_2-r)H}-e^{-(t_1-r)H}\big]\Big(\mathfrak{G}^{(\la)}\big(\ell,m;R,Z_{s,.}\big)_{s,r}\Big)\Big\|_{\cb^\eps_{\infty}}\lesssim \\
&\leq \sum_{i=1,\ldots,7}\bigg[ \int_{t_1}^{t_2} dr \, \Big\|e^{-(t_2-r)H}\Big(\ci^{i,(\la)}_{s,r}\Big)\Big\|_{\cb^\eps_{\infty}}+\int_s^{t_1} dr \, \Big\|\big[e^{-(t_2-r)H}-e^{-(t_1-r)H}\big]\Big(\ci^{i,(\la)}_{s,r}\Big)\Big\|_{\cb^\eps_{\infty}}\bigg]\\
&\lesssim  \sum_{i=1,\ldots,7}\bigg[ \int_{t_1}^{t_2} dr \, \Big\|e^{-(t_2-r)H}\Big(\ci^{i,(\la)}_{s,r}\Big)\Big\|_{\cb^\eps_{\infty}}+|t_2-t_1|^{3\eps}\int_s^{t_1} dr \, \Big\|e^{-(t_1-r)H}\Big(\ci^{i,(\la)}_{s,r}\Big)\Big\|_{\cb^{7\eps}_\infty}\bigg].
\end{align*}
We now apply the bounds \eqref{normifJ} and \eqref{normifJ2}, which yields for all $s\leq t_1<t_2\leq T$
\begin{equation}\label{est05}
 \frac{1}{2C_1C_2^{(\eps)}C_3^{(\eps)}}  |t_1-s|^{\frac{7 \eps}2}\frac{\big\|m_{t_2}-m_{t_1}\big\|_{\cb^\eps_{\infty}}}{|t_2-t_1|^{3\eps}}    \leq \frac12 \big\|f\big\|_{L^\infty}+c_\eps  {\la^{\frac12}} \eta^\ka \big\|(\ell,m)\big\|_{X^{(\eps)}(s,T)}.
\end{equation}

\

Putting all the estimates \eqref{est01}, \eqref{est02}, \eqref{est03}, \eqref{est04} and \eqref{est05} together completes the proof of Proposition \ref{prop:l-m}.

\

\subsection{Conclusion: proof of Theorem \ref{thm-uniqueness}}\label{section-concl}

Let us now complete the proof of uniqueness of the invariant measure, as stated in Theorem \ref{thm-uniqueness}. The argument combines the local contraction estimates derived in Corollary \ref{coro:la-star2} with a probabilistic iteration scheme, following the approach of~\cite[Section~3]{DHYZ25}.

\medskip

Recall first that  thanks to Proposition \ref{prop:conv-arbre} and Lemma \ref{lem:c}, one has for every $p\geq 1$
\begin{equation}\label{general-moment-control}
\sup_{s\geq 0} \mathbb{E}\big[M_{s,s+1}^p\big] < \infty.
\end{equation}
Recall also the definition \eqref{defTs} of 
\begin{equation*} 
T_s^{(\eta)}=\inf\big\{t\geq s: \, M_{s,t}\geq \eta \big\} \wedge (s+1).
\end{equation*}
Then the bound \eqref{general-moment-control} yields (this is analogous to \cite[Lemma 3.24]{DHYZ25}):

\begin{lemma}\label{lem:eta}
There exists $\eta\geq 1$ such that for all $s\geq 0$,
\begin{equation*}
\mathbb{P}\big( T^{(\eta)}_s <s+1\big) \leq \frac{1}{20}.
\end{equation*}
\end{lemma}

\begin{proof}[Proof of Lemma \ref{lem:eta}]
We observe that $\big\{T^{(\eta)}_s <s+1\big\} \subset \big\{M_{s,s+1}\geq \eta\big\}$. Hence 
$$\mathbb{P}\big( T^{(\eta)}_s <s+1\big) \leq  \mathbb{P}\big( M_{s,s+1}\geq \eta\big)\leq \frac{\mathbb{E}\big[M_{s,s+1}\Big]}{\eta}.$$
The result then follows from \eqref{general-moment-control} upon taking $\eta\geq 1$ large enough.
\end{proof}

The remainder of the proof of Theorem \ref{thm-uniqueness} proceeds as follows. We fix $\eta\geq 1$ \textit{once and for all} as in Lemma~\ref{lem:eta}, and set from now on $T_s:=T_s^{(\eta)}$ for every $s\geq 0$. Also, for some (fixed) $\eps>0$ small enough, we define $\la_\star>0$ as in Corollary \ref{coro:la-star2}. In particular, we know that for all $\la\in (0,\la_\star)$ and $s\geq 0$, 
\begin{equation}\label{resul-local}
\big\|\big(X^{(\la),u_1}-X^{(\la),u_2}\big)_{T_s}\big\|_{L^\infty}\leq 10\, e^{-3(T_s-s)} \big\|\big(X^{(\la),u_1}-X^{(\la),u_2}\big)_{s}\big\|_{L^\infty}.
\end{equation}

\

Now, for all $k\geq 0$, we denote by $k=\tau_{k}^{(0)}< \tau_{k}^{(1)} < \ldots < \tau_{k}^{(J_{k})}=k+1$ the partition of $[k,k+1]$ defined through the iterative procedure
$$\tau_{k}^{(i+1)}:=\inf\big\{ t>\tau_{k}^{(i)}: \, M_{\tau_{k}^{(i)},t}\geq \eta\} \wedge (k+1).$$
Due to the independence and stationarity of Brownian increments, we can check that
\begin{equation}\label{stoc-prop}
\tau_{k}^{(i+1)}-\tau_{k}^{(i)}\, \bot \, \cf_{\tau_{k}^{(i)}} \quad , \quad J_k \stackrel{\text{(law)}}{=}J_1 \quad \text{and} \quad J_k\,  \bot \,  \cf_k.
\end{equation}
Moreover, by \eqref{resul-local}, we know that
$$\big\|\big(X^{(\la),u_1}-X^{(\la),u_2}\big)_{\tau^{(i+1)}_{k}}\big\|_{L^\infty}\leq 10\, e^{-3(\tau^{(i+1)}_{k}-\tau^{(i)}_{k})} \big\|\big(X^{(\la),u_1}-X^{(\la),u_2}\big)_{\tau^{(i)}_{k}}\big\|_{L^\infty}.$$
Therefore, for all $k\geq 0$,
\begin{equation*}
\big\|\big(X^{(\la),u_1}-X^{(\la),u_2}\big)_{k+1}\big\|_{L^\infty}\leq e^{-3}10^{J_k} \big\|\big(X^{(\la),u_1}-X^{(\la),u_2}\big)_{k}\big\|_{L^\infty},
\end{equation*}
which entails that for all $T\geq 1$,
\begin{equation*}
\big\|\big(X^{(\la),u_1}-X^{(\la),u_2}\big)_{T}\big\|_{L^\infty}\leq e^{-3(T-1)}10^{J_{T-1}+\ldots+J_1} \big\|\big(X^{(\la),u_1}-X^{(\la),u_2}\big)_{1}\big\|_{L^\infty}.
\end{equation*}

\smallskip

Using the stochastic properties in \eqref{stoc-prop}, we deduce that
\begin{equation}
\mathbb{E}\Big[\big\|\big(X^{(\la),u_1}-X^{(\la),u_2}\big)_{T}\big\|_{L^\infty}\Big]\leq \Big(e^{-3}\mathbb{E}\big[10^{J_1}\big]\Big)^{T-1} \mathbb{E}\Big[\big\|\big(X^{(\la),u_1}-X^{(\la),u_2}\big)_{1}\big\|_{L^\infty}\Big].\label{taking-expectation}
\end{equation}

\

We observe first that, with the notation in \eqref{defi:ula}, 
$$\mathbb{E}\Big[\big\|\big(X^{(\la),u_1}-X^{(\la),u_2}\big)_{1}\big\|_{L^\infty}\Big]\leq \mathbb{E}\Big[\big\|U^{(\la),u_1}_1\big\|_{L^\infty}\Big]+\mathbb{E}\Big[\big\|U^{(\la),u_2}_1\big\|_{L^\infty}\Big]$$
and so, by combining the estimate of Proposition \ref{cor:comes-down} with \eqref{general-moment-control} (for $s=0$), we obtain that
$$C_0:=\sup_{\la\in (0,1)}\sup_{u_1,u_2\in \cb_{\infty}^{-\frac58}} \mathbb{E}\Big[\big\|\big(X^{(\la),u_1}-X^{(\la),u_2}\big)_{1}\big\|_{L^\infty}\Big] <\infty.$$

On the other hand, using the first assertion in \eqref{stoc-prop}, we have for every $i\geq 2$,
\begin{align*}
\mathbb{P}\big(J_1\geq i\big) =\mathbb{P}\big(\tau_1^{(i-1)} <2\big)&\leq \mathbb{E}\big[\1_{\{\tau_1^{(i-2)} <2\}}\1_{\{\tau_1^{(i-1)}-\tau_1^{(i-2)} <1\}}\big] \\
&\leq \mathbb{E}\Big[\mathbb{E}\big[\1_{\{\tau_1^{(i-2)} <2\}}\1_{\{\tau_1^{(i-1)}-\tau_1^{(i-2)} <1\}}\big| \cf_{\tau_1^{(i-2)}}\big]\Big]  \\
&\leq \mathbb{P}\big(\tau_1^{(i-2)} <2\big) \mathbb{P}\big(\tau_1^{(i-1)}-\tau_1^{(i-2)} <1 \big)\\
&\leq \mathbb{P}\big(\tau_1^{(i-2)} <2\big) \mathbb{P}\big(T_0 <1 \big)\leq \frac{1}{20}\mathbb{P}\big(\tau_1^{(i-2)} <2\big) 
\end{align*}
and accordingly, for every $i\geq 1$
\begin{align*}
  \mathbb{P}\big(J_1\geq i\big) \leq \bigg(\frac{1}{20}\bigg)^{i-1} \mathbb{P}\big(\tau_1^{(0)} <2\big)\leq \bigg(\frac{1}{20}\bigg)^{i-1}. 
\end{align*}

\

As a result,
\begin{align*}
\mathbb{E}\big[10^{J_1}\big]=\sum_{i\geq 1} 10^i \mathbb{P}\big(J_1= i\big)  \leq 10\sum_{i\geq 1} \frac{1}{2^{i-1}} = 20.
\end{align*}
Going back to \eqref{taking-expectation}, we can conclude that for every $\la\in (0,\la_\star)$,
\begin{align*}
&\sup_{u_1,u_2\in \cb_{\infty}^{-\frac58}} \mathbb{E}\Big[\big\|\big(X^{(\la),u_1}-X^{(\la),u_2}\big)_{T}\big\|_{L^\infty}\Big]\leq C_0 e^{-(3-\log(20))(T-1)} ,
\end{align*}
and since $3>\log(20)$, this yields the desired bound.

\subsection{Transfer to the model with modulated confining potential}\label{subsec:transfer-alpha}

Throughout this subsection, we recast the main results of the paper in terms of the 
companion model in which the {\it confining harmonic potential} is modulated, instead 
of the nonlinearity. Beyond its conceptual appeal, this reformulation provides a direct physical 
interpretation of the high-temperature regime as a strong-confinement regime, and 
shows that all our previous results transfer with explicit dependence on the trap frequency.

\medskip

In what follows, $\mathbb{D}:=\{2^{2\ell}, \, \ell\geq 0\}$. For $\alpha \in \mathbb{D}$, set 
$$H_\alpha := -\Delta + \alpha^2 |x|^2$$
 on $\R^3$, so that $H = H_1$. 
We are interested in the equation
\begin{equation} \label{eq-alpha-renormalized}
\left\{
\begin{aligned}
&(\partial_t + H_\alpha)\, Y^{(\alpha),(n)} = -\bigl(Y^{(\alpha),(n)}\bigr)^{\!3} 
+ \widetilde{\frakc}^{(\alpha),(n)} Y^{(\alpha),(n)} + \xi^{(\alpha),(n)}, \quad t>0, \quad x \in \R^3, \\
& Y^{(\alpha),(n)}(0) = v,
\end{aligned}
\right.
\end{equation} 
where the noise regularization $\xi^{(\alpha),(n)}$ is the natural extension, for the operator 
$H_\alpha$, of the previous regularization $\xi^{(n)}$ (see \eqref{regu-noise-alpha}), 
and the renormalization sequence 
$\widetilde{\frakc}^{(\alpha),(n)}$ will be specified in \eqref{def-c-alpha}.

\subsubsection{The scaling correspondence}\label{subsubsec:scaling-correspondence}

For $\alpha \in \mathbb{D}$, define the spatial dilation $\cs_\alpha$ acting on functions 
$f : \R^3 \to \R$, and the space-time dilation $\ct_\alpha$ acting on trajectories 
$Z : \R_+ \times \R^3 \to \R$, by
\begin{equation}\label{def-S-T-alpha}
(\cs_\alpha f)(x) := \alpha^{1/4} f(\alpha^{1/2} x),
\qquad 
(\ct_\alpha Z)(t,x) := \alpha^{1/4} Z(\alpha t, \alpha^{1/2} x).
\end{equation}
A direct computation gives the intertwining relation
\begin{equation}\label{intertwining-H}
\cs_\alpha \circ H = \alpha^{-1} H_\alpha \circ \cs_\alpha.
\end{equation}
In particular, denoting by $(\vp_k)_{k\geq 0}$ the Hermite basis of $H$ from Section \ref{Sect2}, with eigenvalues $(\la_k)_{k\geq 0}$, the family
\begin{equation}\label{def-vp-alpha}
\vp^{(\alpha)}_k(x) := \alpha^{3/4} \vp_k(\alpha^{1/2} x)
\end{equation}
is a Hilbertian basis of $L^2(\R^3)$ made of eigenvectors of $H_\alpha$ with eigenvalues 
$\la^{(\alpha)}_k := \alpha\, \la_k$. The kernel of $e^{-\sigma H_\alpha}$ is therefore
\begin{equation}\label{K-alpha-scaling}
K^{(\alpha)}_\sigma(x,y) = \alpha^{3/2}\, K_{\alpha \sigma}(\alpha^{1/2} x, \alpha^{1/2} y).
\end{equation}

\subsubsection{Equivalence of harmonic Besov spaces}\label{subsubsec:besov-scaling}

For convenience, let us denote by $\cb^\sigma_{\infty,\infty}(H_\alpha)$ the harmonic Besov space 
associated to $H_\alpha$, defined as in \S\ref{sect-Besov} with $H$ replaced by 
$H_\alpha$. Notice that for $\alpha = 1$, this is the space $\cb^\sigma_{\infty,\infty}=\cb^\sigma_{\infty,\infty}(H)$ used throughout the paper. The multipliers $(\delta^{(\alpha)}_j)$ adapted to $H_\alpha$ are naturally defined by 
$$\delta^{(\alpha)}_j := \theta\Big(\frac{H_\alpha}{2^{2j}}\Big)  $$
and the norm reads 
$$ \|u\|_{\cb^\sigma_{\infty,\infty}(H_\alpha)} = \sup_{j \geq 0} 2^{j\sigma}\, 
\|\delta^{(\alpha)}_j u\|_{L^\infty(\R^3)}.$$

\begin{lemma}[Equivalence of harmonic Besov spaces under scaling]\label{lem:besov-scaling}
For all $\alpha\in \mathbb{D}$ and $\sigma \in \R$, one has
\begin{equation}\label{besov-equiv}
\|\cs_\alpha f\|_{\cb^\sigma_{\infty,\infty}(H_\alpha)} = \alpha^{1/4+ \sigma/2}\, \|f\|_{\cb^\sigma_{\infty,\infty}(H)}.
\end{equation}

\end{lemma}

\begin{proof}
Assume that $\al=2^{2\ell}$ ($\ell\geq 0$), and write, thanks to identity \eqref{def-vp-alpha},
\begin{align*}
\delta_j^{(\al)}\big(\cs_\al f\big)(x)&=\sum_{k\geq 0} \theta\Big(\frac{\la^{(\al)}_k}{2^{2j}}\Big) \langle \cs_\al f,\vp^{(\al)}_k\rangle\, \vp^{(\al)}_k(x)\\
&=\al^{\frac14}\sum_{k\geq 0} \theta\Big(\frac{\al \la_k}{2^{2j}}\Big) \langle  f,\vp_k\rangle\, \vp_k(\al^{\frac12}x)=\al^{\frac14}\sum_{k\geq 0} \theta\Big(\frac{\la_k}{2^{2(j-\ell)}}\Big) \langle  f,\vp_k\rangle\, \vp_k(\al^{\frac12}x).
\end{align*}
Since $\operatorname{Supp}\theta  \subset \Big\{ \xi \in \mathbb{R} : (\tfrac{3}{4})^2 \le |\xi| \le (\tfrac{8}{3})^2 \Big\}$ and $\la_k\geq 3$, the above quantity reduces to
\begin{align*}
\delta_j^{(\al)}\big(\cs_\al f\big)(x)
&=\al^{\frac14}\1_{\{j-\ell \geq 0\}} \delta_{j-\ell}(f)(\al^{\frac12}x)
\end{align*}
and accordingly
\begin{align*}
\|\cs_\alpha f\|_{\cb^\sigma_{\infty,\infty}(H_\alpha)}=\sup_{j\geq 0} 2^{j\si}\big\|\delta_j^{(\al)}\big(\cs_\al f\big)\big\|_{L^\infty}
&=\al^{\frac14}\sup_{j\geq \ell} 2^{j\si} \big\|\delta_{j-\ell}(f)\big\|_{L^\infty}\\
&=\al^{\frac14}2^{\ell \si}\sup_{j\geq 0} 2^{j\si} \big\|\delta_{j}(f)\big\|_{L^\infty}=\al^{\frac14+\frac{\si}{2}} \| f\|_{\cb^\sigma_{\infty,\infty}(H)}.
\end{align*}
\end{proof}

\subsubsection{Noise regularization for the $\alpha$-model}\label{subsubsec:noise-alpha}

Following \eqref{repres-noise} and \eqref{regu-noise}, we define the regularization 
adapted to $H_\alpha$ by
\begin{equation}\label{regu-noise-alpha}
\xi^{(\alpha),(n)}_t := \frac{dW^{(\alpha),(n)}_t}{dt}, 
\qquad 
W^{(\alpha),(n)}_t(x) := \sum_{k \geq 0} e^{-\varepsilon^{(\alpha)}_n \la^{(\alpha)}_k}
\beta^{(k)}_t \vp^{(\alpha)}_k(x),
\end{equation}
where $(\beta^{(k)})_{k\geq 0}$ is again the family of independent Brownian motions introduced above, and the regularization parameter is chosen as
\begin{equation}\label{def-eps-alpha-n}
\varepsilon^{(\alpha)}_n := \alpha^{-1} \varepsilon_n= \alpha^{-1} 2^{-n} .
\end{equation}

\

\begin{lemma}[Scaling identity for the regularized noise]\label{lem:noise-rescaling}
With the choice \eqref{def-eps-alpha-n}, the field
\begin{equation}\label{def-xi-tilde}
\widetilde{\xi}^{(n)}_t(y) := \alpha^{-5/4}\, \xi^{(\alpha),(n)}_{\alpha^{-1}t}(\alpha^{-1/2} y)
\end{equation}
satisfies $\widetilde{\xi}^{(n)} \stackrel{\mathrm{law}}{=} \xi^{(n)}$.
\end{lemma}

\begin{proof}
By \eqref{regu-noise-alpha}, \eqref{def-vp-alpha} and \eqref{def-eps-alpha-n}, one has
\[
W^{(\alpha),(n)}_t(x) 
= \alpha^{3/4}\, \sum_{k \geq 0} e^{-\varepsilon_n  \la_k}\, 
\beta^{(k)}_t\, \vp_k(\alpha^{1/2} x),
\]
and hence
\[
W^{(\alpha),(n)}_{\alpha^{-1} t}\bigl(\alpha^{-1/2} y\bigr) 
= \alpha^{3/4}\, \sum_{k \geq 0} e^{-\varepsilon_n  \la_k}\, 
\beta^{(k)}_{\alpha^{-1} t}\, \vp_k(y)=\alpha^{1/4}\, \sum_{k \geq 0} e^{-\varepsilon_n  \la_k}\, 
\big(\alpha^{1/2}\,\beta^{(k)}_{\alpha^{-1} t}\big)\, \vp_k(y).
\]
By the standard scaling property of Brownian motion, the family 
$\bigl(\alpha^{1/2}\, \beta^{(k)}_{\alpha^{-1}t}\bigr)_{k\geq 0}$ corresponds, in law, to a family 
of independent standard Brownian motions, and accordingly one has
\[
W^{(\alpha),(n)}_{\alpha^{-1} t}\bigl(\alpha^{-1/2} y\bigr) 
\,\stackrel{\mathrm{law}}{=}\, 
\alpha^{1/4}\, W^{(n)}_t(y).
\]
Taking the time derivative in $t$ and using 
$$\partial_t \bigl[ W^{(\alpha),(n)}_{\alpha^{-1}t}(\alpha^{-1/2}y) \bigr] 
= \alpha^{-1} \xi^{(\alpha),(n)}_{\alpha^{-1}t}( \alpha^{-1/2}y),$$
 together with 
$\partial_t W^{(n)}_t(y) = \xi^{(n)}_t(y)$, one immediately obtains the announced identity.
\end{proof}

\
 
\subsubsection{Scaling identity for solutions}
 
From now on, and with a slight abuse of notation, we still denote by $X^{(\la),(n)}$ the solution of \eqref{eq-la} associated with the white noise approximation $\widetilde{\xi}^{(n)}$ defined in \eqref{def-xi-tilde}. Thanks to the identification result in Lemma \ref{lem:noise-rescaling}, we know that $(X^{(\la),(n)})$ satisfies all the properties listed in Theorem \ref{thm:main}.
 
\smallskip
 
Using the notation in \eqref{constante-re}, we now define the renormalization sequence involved in the $\alpha$-model~\eqref{eq-alpha-renormalized} as
\begin{equation}\label{def-c-alpha}
\widetilde{\frakc}^{(\alpha),(n)} := 3\,\widetilde{\frakc}^{\mathbf{1},(\alpha),(n)} 
- 9\,\widetilde{\frakc}^{\mathbf{2},(\alpha),(n)},
\end{equation}
with
\begin{equation}\label{c-alpha-scaling}
\widetilde{\frakc}^{\mathbf{1},(\alpha),(n)}(x) 
= \alpha^{1/2}\, \frakc^{\mathbf{1},(n)}(\alpha^{1/2} x),
\qquad
\widetilde{\frakc}^{\mathbf{2},(\alpha),(n)}(x) 
= \frakc^{\mathbf{2},(n)}(\alpha^{1/2} x).
\end{equation}

\begin{lemma}[Scaling identity for solutions]\label{lem:scaling}
Fix $\alpha > 0$, $n \geq 1$, and set 
\[
\lambda := \alpha^{-1/2}.
\]
For every initial datum $Y_0 \in \cb^{-\frac12 - \eps}_{\infty,\infty}(H_{\alpha})$, set
\[
X_0 := \alpha^{-1/2}\, \cs_\alpha^{-1} Y_0 \,\in\, \cb^{-\frac12 - \eps}_{\infty,\infty}(H).
\]
Let $Y^{(\alpha),(n)}$ denote the unique solution of \eqref{eq-alpha-renormalized} 
with initial datum $Y_0$, and $X^{(\lambda),(n)}$ the unique solution of 
\eqref{eq-la} with parameter $\lambda$, noise $\widetilde{\xi}^{(n)}$ and initial datum $X_0$. Then
\begin{equation}\label{scaling-identity}
Y^{(\alpha),(n)} = \ct_\alpha\bigl( X^{(\lambda),(n)} \bigr).
\end{equation}
 
\end{lemma}
 
\smallskip
 
\begin{proof}
Set $\Xi_t(y) := (\ct_\alpha^{-1} Y^{(\alpha),(n)})_t(y) 
= \alpha^{-1/4}\, Y^{(\alpha),(n)}_{\alpha^{-1}t}( \alpha^{-1/2} y)$. Differentiating, we obtain
\[
(\partial_t \Xi)_t(y) 
= \alpha^{-5/4}\, (\partial_t Y^{(\alpha),(n)})_{\alpha^{-1}t}( \alpha^{-1/2}y),
\qquad
(H \Xi)_t(y) 
= \alpha^{-5/4}\, (H_\alpha Y^{(\alpha),(n)})_{\alpha^{-1}t}( \alpha^{-1/2}y),
\]
where the second identity can be derived from \eqref{intertwining-H}. As a result,
\begin{align*}
\big((\partial_t + H) \Xi\big)_t(y)&= \alpha^{-5/4}\big( (\partial_t Y^{(\alpha),(n)})+ (H_\alpha Y^{(\alpha),(n)})\big)_{\al^{-1}t}(\al^{-\frac12}y)\\
&= \alpha^{-5/4}\big( -\bigl(Y^{(\alpha),(n)}\bigr)^{3} + \widetilde{\frakc}^{(\alpha),(n)} Y^{(\alpha),(n)} + \xi^{(\alpha),(n)}\big)_{\al^{-1}t}(\al^{-\frac12}y)\\
&= -\alpha^{-1/2}\, \big(\Xi_t(y)\big)^3 
+ \alpha^{-1}\, \widetilde{\frakc}^{(\alpha),(n)}(\alpha^{-1/2} y)\cdot  \Xi_t(y) 
+ \widetilde{\xi}^{(n)}_t(y).
\end{align*}
Combining \eqref{def-c-alpha} and \eqref{c-alpha-scaling}, we obtain that
\begin{align*}
\alpha^{-1}\, \widetilde{\frakc}^{(\alpha),(n)}(\alpha^{-1/2} y) 
&= 3 \alpha^{-1}\, \alpha^{1/2}\, \frakc^{\mathbf{1},(n)}(y) - 9 \alpha^{-1}\, \frakc^{\mathbf{2},(n)}(y) = 3 \lambda\, \frakc^{\mathbf{1},(n)}(y) - 9 \lambda^2\, \frakc^{\mathbf{2},(n)}(y) ,
\end{align*}
which yields
\begin{align*}
\big((\partial_t + H) \Xi\big)_t(y)&= -\alpha^{-1/2}\, \big(\Xi_t(y)\big)^3 
+\big(3 \lambda\, \frakc^{\mathbf{1},(n)}(y) - 9 \lambda^2\, \frakc^{\mathbf{2},(n)}(y)\big)\cdot  \Xi_t(y) 
+ \widetilde{\xi}^{(n)}_t(y).
\end{align*}
In other words, $\Xi$ solves the equation 
\eqref{eq-la} with parameter $\lambda = \alpha^{-1/2}$, initial datum $X_0$, 
and driving noise $\widetilde{\xi}^{(n)}$. By uniqueness of the solution,
\[
\Xi =X^{(\lambda),(n)},
\]
and applying $\ct_\alpha$ to both sides gives \eqref{scaling-identity}.
\end{proof}

\smallskip

\subsubsection{Proof of Theorem \ref{thm:uniq-alpha}}

\begin{proof}
{\rm (i)} {\it Convergence.} Set $u:= \alpha^{-1/2}\, \cs_\alpha^{-1} v$. By Theorem~\ref{thm:main}~$(ii)$, applied with 
$\la = \alpha^{-1/2}$, the sequence $(X^{(\la),(n)})_{n \geq 1}$ converges almost surely 
in $\bigcap_{\eta>0} \cac\big([0,T); \cb^{-\frac12-\eta}_{\infty,\infty}(H)\big)$ to a 
limit $X^{(\la),u}$. Applying the map $\ct_\alpha$, whose continuity between the relevant 
Besov spaces is provided by Lemma~\ref{lem:besov-scaling}, and invoking the identity \eqref{scaling-identity}, we deduce that $(Y^{(\alpha),(n)})_{n \geq 1}$ converges 
 in $\bigcap_{\eta>0} \cac\big([0,T); \cb^{-\frac12-\eta}_{\infty,\infty}(H_\al)\big)$ to the limit 
\begin{equation}\label{keyidentity}
Y^{(\alpha),v} =\ct_\alpha\, X^{(\la),u}.
\end{equation}

\smallskip
 
\noindent
{\rm (ii)} {\it Pushforward formula and invariance.} Identity \eqref{keyidentity} can be rephrased as
\begin{equation}\label{keyidentity-bis}
(Y^{(\alpha), \cs_\al u}_t)_{t \geq 0}=(\cs_\alpha X^{(\la), u}_{\alpha t})_{t \geq 0},
\end{equation}
for every $u \in \cb^{-\frac12-\eps}_{\infty,\infty}(H)$. Recall  that the map
$\cs_\alpha : \cb^{-\frac12-\eps}_{\infty,\infty}(H) \to \cb^{-\frac12-\eps}_{\infty,\infty}(H_\alpha)$
is a linear isomorphism, as follows from Lemma~\ref{lem:besov-scaling}. Therefore, for any bounded measurable functional $\Phi : \cb^{-\frac12-\eps}_{\infty,\infty}(H_\alpha) \to \R$,
\begin{align*}
Q^{(\al)}_T\big( (\cs_\al)_\sharp \rho^{(\la)}\big)(\Phi)&=\int_{\cb^{-\frac12-\eps}_{\infty,\infty}(H_\al)} \big(Q^{(\al)}_T \Phi\big)(v) \, \big((\cs_\al)_\sharp \rho^{(\la)}\big)(dv)\\
&=\int_{\cb^{-\frac12-\eps}_{\infty,\infty}(H_\al)} \mathbb{E}\big[ \Phi\big( Y^{(\al),v}_T \big)\big] \, \big((\cs_\al)_\sharp \rho^{(\la)}\big)(dv)\\
&=\int_{\cb^{-\frac12-\eps}_{\infty,\infty}(H)} \mathbb{E}\big[ \Phi\big( Y^{(\al),\cs_\al u}_T \big)\big] \, \rho^{(\la)}(du)\\
&=\int_{\cb^{-\frac12-\eps}_{\infty,\infty}(H)} \mathbb{E}\big[ \big(\Phi\circ \cs_\al\big)\big( X^{(\la), u}_{\al T} \big)\big] \, \rho^{(\la)}(du) \qquad \text{(by \eqref{keyidentity-bis})}\\
&=P^{(\la)}_{\al T}\big(  \rho^{(\la)}\big)(\Phi\circ \cs_\al)= \rho^{(\la)}(\Phi\circ \cs_\al)=\big( (\cs_\al)_\sharp \rho^{(\la)}\big)(\Phi),
\end{align*} 
and therefore $(\cs_\al)_\sharp \rho^{(\la)}$ is invariant for $Q^{(\al)}$.

\smallskip
 
{\it Uniqueness.} Assume that $\al >\la_\star^{-2}$, where $\la_\star>0$ is the threshold introduced in Theorem \ref{thm-uniqueness}. Then, if $\nu, \widetilde{\nu}$ are two invariant measures for 
$Q^{(\alpha)}$, the above arguments show that $((\cs_\al)_\sharp)^{-1}\nu$ and $((\cs_\al)_\sharp)^{-1}\widetilde\nu$ are two invariant measures of $P^{(\la)}$ at $\la = \alpha^{-1/2} < \la_\star$. By the 
uniqueness statement of item~$(iv)$ of Theorem~\ref{thm:main}, $((\cs_\al)_\sharp)^{-1}\nu=((\cs_\al)_\sharp)^{-1}\widetilde\nu$,
hence $\nu = \widetilde\nu = (\cs_\al)_\sharp \rho^{(\la)}$.
 
{\it Non-Gaussianity.} Since $\cs_\alpha$ is a linear bijection between $\cb^{-\frac12-\eps}_{\infty,\infty}$ spaces and $\rho^{(\la)}$ is non-Gaussian (see item~$(iii)$ of Theorem~\ref{thm:main}), its pushforward $\nu^{(\alpha)}$ is also non-Gaussian.
 
\smallskip

\noindent
{\rm (iii)} Setting $u = (\cs_\alpha)^{-1} v$ and using \eqref{keyidentity}, we can write
\[
\mathbb{E}[\Phi(Y^{(\alpha),v}_T)] - \nu^{(\alpha)}(\Phi) 
\,=\, \mathbb{E}[(\Phi \circ \cs_\alpha)(X^{(\la),u}_{\alpha T})] - \rho^{(\la)}(\Phi \circ \cs_\alpha).
\]
Observe that $\|\Phi \circ \cs_\alpha\|_{\mathrm{Lip}} \leq \alpha^{1/4}\, \|\Phi\|_{\mathrm{Lip}}$, and thus, applying item~$(iv)$ of Theorem~\ref{thm:main} for the original model at time 
$\alpha T \geq 1$ yields \eqref{expo-conv-alpha}.
\end{proof}

\
 

\appendix
\section{Some functional inequalities}\label{appendix:functional}
In this section we gather some useful inequalities.   Recall that, by definition \eqref{def-theta},
\begin{equation}\label{def-bis}
 \big\|v\big\|_{\B^\alpha_{p,\infty}(\R^d)}= \sup_{j \geq 0} \big(2^{j \alpha} \big\| \delta_j(v) \big\|_{L^p(\R^d)}\big).
\end{equation}

  \begin{lemma} 
Let $d\geq 1$, $\alpha_0, \alpha_1 \in \R$, $1\leq p_0,p_1 \leq \infty$ and $0 \leq \nu \leq 1$. Define $\alpha \in \R$ and $1 \leq p \leq \infty$ by
$$  \alpha = (1- \nu) \alpha_0+ \nu \alpha_1, \qquad  \frac1p=\frac{1-\nu}{p_0}+ \frac{\nu}{p_1}.$$
Then one has 
\begin{equation}\label{interp}
\big\| v\big\|_{\B^\alpha_{p,\infty}(\R^d)} \lesssim  \big\| v\big\|^{1-\nu}_{\B^{\alpha_0}_{p_0,\infty}(\R^d)}  \big\| v\big\|^\nu_{\B^{\alpha_1}_{p_1,\infty}(\R^d)} .
\end{equation}
 \end{lemma}
 
  \begin{proof}
The result directly follows from \eqref{def-bis} and the H\"older inequality.
    \end{proof}

 \begin{lemma}\label{lem-bis}
Let $d\geq 1$, $\eta>0$ and $1\leq p \leq \infty$ be such that $\frac{d}p <{\eta}$. Then for all $\alpha \in \R$
\begin{equation}\label{lplq0}
\big\| v\big\|_{\B^\alpha_{p,\infty}(\R^d)} \lesssim  \big\| v\big\|_{\B^{\alpha+\eta}_{\infty,\infty}(\R^d)}.
\end{equation}
More generally, if $\eta>0$ and $1\leq p \leq q \leq \infty$ are such that $\dis d(\frac1p-\frac1q)< \eta$, then
\begin{equation}\label{lplq}
\big\| v\big\|_{\B^\alpha_{p,\infty}(\R^d)} \lesssim  \big\| v\big\|_{\B^{\alpha+\eta}_{q,\infty}(\R^d)}.
\end{equation}
 \end{lemma}

 \begin{proof}
We only prove \eqref{lplq0}; the argument for \eqref{lplq} is similar. 
 For $\eta >\frac{d}p $, we write 
 \begin{equation}\label{eqq1}
 \big\| \delta_j (v) \big\|_{L^p(\R^d)} \leq \big\| \langle x\rangle^{\eta} \delta_j (v) \big\|_{L^\infty(\R^d)}   \big\| \langle x\rangle^{-\eta} \big\|_{L^p(\R^d)} \lesssim \big\| \langle x\rangle^{\eta} \delta_j (v) \big\|_{L^\infty(\R^d)}. 
  \end{equation}
 Next, by the proof of \cite[Lemma 13.12]{DFT} we have
  \begin{equation}\label{eqq2}
  \big\| \langle x\rangle^{\eta} \delta_j (v) \big\|_{L^\infty(\R^d)}    \lesssim 2^{j \eta} \big\| \widehat{\delta}_j (v) \big\|_{L^\infty(\R^d)}, 
    \end{equation}
where $\widehat{\delta}_j$ is a spectral truncation which has a slightly larger support than $ \delta_j$ and satisfies $\widehat{\delta}_j  {\delta}_j= {\delta}_j$.  Combining the estimates \eqref{eqq1} and \eqref{eqq2} yields the result.
  \end{proof}
 
 Let us recall the Sobolev embeddings (see {\it e.g.} \cite[Proposition 1.1]{FI}).
 \begin{lemma} \label{lem:inclusion-besov}
Let $d\geq 1$, $\alpha  \in \mathbb R$ and  $1 \leq p \leq q \leq \infty$. Then 
$$\big\|u\big\|_{\mathcal{B}^{\alpha}_{q,\infty}(\R^d)}  \le C  \big\|u\big\|_{\mathcal{B}^{\alpha+d(\frac{1}{p}-\frac{1}{q})}_{p,\infty}(\R^d)} . $$
\end{lemma}

For the reader's convenience, we recall some bilinear estimates. For the proof, we refer to \cite[Proposition 13.19]{DFT}. 
\begin{proposition}\label{Prop-est-para}
Let $\alpha, \beta \in \R$ and $1 \leq p,p_1,p_2,q,q_1,q_2 \leq \infty$ be such that 
 $$\frac1{p_1}+\frac1{p_2}=\frac1{p} \quad \text{ and } \quad\frac1{q_1}+\frac1{q_2}=\frac1{q}.$$ 
 \begin{enumerate}[$(i)$]
 \item If $\alpha+\beta>0$, then the mapping $(f,g) \mapsto f \pe g$ extends to a continuous bilinear map from $\mathcal{B}^{\alpha}_{p_1,q_1}(\R^d) \times \mathcal{B}^{\beta}_{p_2,q_2}(\R^d) $ to $\mathcal{B}^{\alpha+\beta}_{p,q}(\R^d) $.
 \item The mapping $(f,g) \mapsto f \pl g$ extends to a continuous bilinear map from $L^{p_1}(\R^d) \times \mathcal{B}^{\beta}_{p_2,q}(\R^d)$ to $\mathcal{B}^{\beta}_{p,q}(\R^d)$.
 \item If $\alpha <0$, then the mapping $(f,g) \mapsto f \pl g$ extends to a continuous bilinear map from $\mathcal{B}^{\alpha}_{p_1,q_1}(\R^d)\times \mathcal{B}^{\beta}_{p_2,q_2}(\R^d)$ to $\mathcal{B}^{\alpha+\beta}_{p,q}(\R^d)$.
  \item If $\alpha <0< \beta$ and $\alpha+\beta>0$, then the mapping $(f,g) \mapsto f  g$ extends to a continuous bilinear map from $\mathcal{B}^{\alpha}_{p_1,q_1}(\R^d)\times \mathcal{B}^{\beta}_{p_2,q_2}(\R^d)$ to $\mathcal{B}^{\alpha}_{p,q}(\R^d)$.
    \item If $\alpha >0$, then the mapping $(f,g) \mapsto f  g$ extends to a continuous bilinear map from $\mathcal{B}^{\alpha}_{p_1,q}(\R^d)\times \mathcal{B}^{\alpha}_{p_2,q}(\R^d)$ to $\mathcal{B}^{\alpha}_{p,q}(\R^d)$. Moreover, for $1 \leq p_3, p_4 \leq \infty$ such that $\frac1{p_1}+\frac1{p_2}=\frac1{p_3}+\frac1{p_4}=\frac1{p}$, there exists $C>0$ satisfying
    \begin{equation*}
    \big\| fg \big\|_{ \mathcal{B}^{\alpha}_{p,q}(\R^d)} \leq C \Big( \big\| f \big\|_{L^{p_1}(\R^d)}   \big\| g \big\|_{ \mathcal{B}^{\alpha}_{p_2,q}(\R^d)}+ \big\| f \|_{ \mathcal{B}^{\alpha}_{p_3,q}(\R^d)}  \big\| g \big\| _{L^{p_4}(\R^d)}       \big).
    \end{equation*} 
 \end{enumerate}
  \end{proposition}

We now recall the following bounds on the heat flow (see \cite[Lemma 13.14]{DFT}).

 \begin{lemma}\label{lem:actisemi}
Let $d\geq 1$. Let $\alpha, \beta \in \R$ and $1 \leq p, q\leq \infty$.
 \begin{enumerate}[$(i)$]
 \item For all $t\geq 0$
  \begin{equation}\label{heat-Lp}
 \big\|  e^{-tH } u \big\|_{L^p(\R^d)} \leq C e^{- d t}  \big\|  u \big\|_{L^p(\R^d)}.
        \end{equation}
 \item If $\alpha \geq \beta$, then there exists $C>0$ such that for all $t>0$
 \begin{equation}\label{heat1}
 \big\| e^{-t H} u \big\|_{\B^\alpha_{p,q} (\R^d)}\leq C t^{-\frac{\alpha-\beta}2} e^{-\frac{dt}2} \big\|  u \big\|_{\B^\beta_{p,q}(\R^d) }.
       \end{equation}
      \item If $0 \leq \beta -\alpha \leq 2$, then there exists $C>0$ such that for all $t>0$
 \begin{equation}\label{heat2}
  \big\| (1-e^{-t H} )u \big\|_{\B^\alpha_{p,q} (\R^d)}\leq C t^{-\frac{\alpha-\beta}2}  \big\|  u \big\|_{\B^\beta_{p,q}(\R^d) }.
       \end{equation}  
  \end{enumerate}
 \end{lemma}

We also need the following bound, involving an explicit constant:

\begin{lemma}\label{lem:dep}
For every $t\geq 0$, it holds that
$$\big\|e^{-tH}f\big\|_{L^\infty}\leq 8e^{-3t}\big\|f\big\|_{L^\infty}.$$
\end{lemma}
\begin{proof}
One has 
$$\big\|e^{-tH}f\big\|_{L^\infty}\leq \big\|f\big\|_{L^\infty} \cdot \bigg(\sup_{x\in \R^3} \int dy \, K_t(x,y)\bigg).$$
Then 
\begin{align*}
\sup_{x\in \R^3} \int dy \, K_t(x,y)&\leq \frac{1}{\big(2\pi \sinh(2t)\big)^{\frac32}}\int dy \, \exp\bigg(-\frac{|y|^2}{4\tanh(t)}\bigg)\\
&\leq \bigg( \frac{2 \tanh(t)}{\sinh(2t)}\bigg)^{\frac32} \bigg(\frac{1}{\big(2\pi \big)^{\frac32}}\int dy \, \exp\bigg(-\frac{|y|^2}{2}\bigg)\bigg)\\
&\leq \bigg( \frac{1}{\cosh(t)}\bigg)^{3}\leq 8 \, e^{-3t}.
\end{align*}

\end{proof}

\

 
\section{About the fourth order diagram 1} \label{sec:diag-4th-order-1}

Recall that the functions $\{F^{(n)}_{t,x}, \, t\geq 0, \, x\in \R^3\}$, resp. $\{\cac^{(n)}_{s,t}(x,y), \, s,t\geq 0, \, x,y\in \R^3\}$, have been introduced in \eqref{defi-f-n-ell}, resp. \eqref{cova-c-n-l}.
Let us now introduce some additional notation that will be used in the next three sections.

\

\begin{notation}\label{notations}
For any function $g : \R^d \times \R^d  \longrightarrow \R$, we denote by $ g^\star$ the auxiliary function on $\R^d$ defined by 
\begin{equation}\label{star-not} 
g^\star(x) = \sup_{y \in \R^d}|g(y,y+x)|.
\end{equation}
Moreover, we will use the following notation for the increments of the kernel $K$ and the covariance function $\cac^{(n)}$: for $s,s',t,t'\geq 0$,
\begin{equation}\label{inc-k}
K_{s,t}:=K_t-K_s
\end{equation}
and
\begin{equation}\label{inc-cn}
\cac^{(n)}_{(s,s'),t}:=\cac^{(n)}_{s',t}-\cac^{(n)}_{s,t}, \quad \cac^{(n)}_{s,(t,t')}:=\cac^{(n)}_{s,t'}-\cac^{(n)}_{s,t}, \quad \cac^{(n)}_{(s,t),(s',t')}:=\cac^{(n)}_{t,t'}-\cac^{(n)}_{s,t'}-\cac^{(n)}_{t,s'}+\cac^{(n)}_{s,s'}.
\end{equation}
In particular, we have
$$\cac^{(n)}_{(s,t),(s',t')}=\cac^{(n)}_{(s,t),t'}-\cac^{(n)}_{(s,t),s'}=\cac^{(n)}_{t,(s',t')}-\cac^{(n)}_{s,(s',t')}.$$
\end{notation}

\

Let us also state the following technical property for further use.

\begin{lemma}\label{lem:techn}
Let $X:\Omega\to \cac^\infty(\R^d)$ be a random function, and set, for all $y_1,y_2\in \R^d$,
$$\cm(y_1,y_2):=\mathbb{E}\big[X(y_1)X(y_2)\big].$$
Then, for every $\al >\frac34$, there exists $p_\al \geq 1$ such that for every $p\geq p_\al$,
$$\int_{\R^d}dx \, \mathbb{E} \Big[ \big| H^{-\al}(X)(x)\big|^2\Big]^p \lesssim \bigg(\int_{\R^d} dy \, \cm^\star(y)\bigg)^p.$$
where the proportional constant only depends on $\al$.

\end{lemma}

\begin{proof}

Denote by $h_{-\al}(x,y)$ the kernel associated with $H^{-\al}$. Then one has
\begin{align*}
\mathbb{E}\bigg[ \Big| H^{-\al}\big(X\big)(x)\Big|^{2}\bigg]
&=\bigg|\int dy_1 dy_2\, h_{-\al}(x,y_1)h_{-\al}(x,y_2) \cm(y_1,y_2)\bigg|\\
&\lesssim \big\|h_{-\al}(x,\cdot)\big\|_{L^2_{y_2}} \bigg(\int dy_2\, \bigg| \int dy_1 \, h_{-\al}(x,y_1)\cm(y_1,y_2)\bigg|^2\bigg)^{\frac12}\\
&\lesssim \big\|h_{-\al}(x,\cdot)\big\|_{L^2_{y_2}} \bigg(\int dy_2\, \bigg| \int dy_1 \, h_{-\al}(x,y_1+y_2)\cm(y_1+y_2,y_2)\bigg|^2\bigg)^{\frac12}\\
&\lesssim \big\|h_{-\al}(x,\cdot)\big\|_{L^2_{y_2}} \bigg(\int dy_2\, \bigg| \int dy_1 \, \big|h_{-\al}(x,y_2+y_1)\big| \Big(\sup_{y_2}\big|\cm(y_1+y_2,y_2)\big| \Big)\bigg|^2\bigg)^{\frac12}\\
&\lesssim \big\|h_{-\al}(x,\cdot)\big\|_{L^2_{y}}^2 \bigg(\int dy_1 \,  \Big(\sup_{y_2}\big|\cm(y_1+y_2,y_2)\big| \Big)\bigg).
\end{align*}
Therefore, using that $\dis \big\|h_{-\alpha}(x , \cdot) \big\|^2_{L^2_y}=h_{-2\alpha}(x,x)$, we obtain 
\begin{equation*}
\int dx \, \mathbb{E}\bigg[ \Big| H^{-\al}\big(X\big)(x)\Big|^{2}\bigg]^p \lesssim \big(\int dx \, \big|h_{-2\al}(x,x)\big|^{p}\bigg) \bigg(\int dy_1 \,  \Big(\sup_{y_2}\big|\cm(y_1+y_2,y_2)\big| \Big)\bigg)^p,
\end{equation*}
which, since $\al>\frac34$, allows us to conclude (by \cite[Lemma 2.1]{DFT}) that
\begin{align*}
\int dx \, \mathbb{E}\bigg[ \Big| H^{-\al}\big(X\big)(x)\Big|^{2}\bigg]^p
&\lesssim \bigg(\int dy_1 \,  \Big(\sup_{y_2}\big|\cm(y_1+y_2,y_2)\big| \Big)\bigg)^p,
\end{align*}
for every $p$ large enough.
\end{proof}

\

In this section, we focus on the fourth-order diagram $\<PsiIPsi3>$, that is, the limit of
\begin{equation}\label{ord4}
\<PsiIPsi3>_{s,t}^{(n)}:=\<Psi>_{s,t}^{(n)}\pe \<IPsi3>_{s,t}^{(n)}.
\end{equation}
A first step toward the construction of this diagram was carried out in \cite[Proposition~8.1]{DFT}. Combined with the stationarity property of Lemma \ref{lem:stati-z}, this result immediately extends as follows:

\begin{proposition}\label{Prop-p41}
Fix $0\leq T_1<T_2$ and for all $t\in [T_1,T_2]$, set
\begin{equation*}
\widetilde{ \<PsiIPsi3>}^{(n)}_{T_1,t}:=\int_{T_1}^t \<PsiIPsi3>^{(n)}_{T_1,s} \, ds.
\end{equation*}
Then for all $0<\eps,\eta<\frac12$, there exists $\ka>0$ such that for every $p\geq 1$,
\begin{equation}\label{L3}
\sup_{\ell\geq 0}    \mathbb{E} \Big[ \Big\|\widetilde{ \<PsiIPsi3>}^{(n+1)}_{\ell,.}   - \widetilde{ \<PsiIPsi3>}^{(n)}_{\ell,.}\Big\|_{{\ov \cac}^{1-\eps}([\ell,\ell+2]; \cb^{-\eta}_{\infty})}^{2p} \Big]\lesssim 2^{-\ka n p }. 
\end{equation}
As a particular consequence, for all $0\leq T_1<T_2$, the sequence $(\widetilde{ \<PsiIPsi3>}^{(n)}_{T_1,.})$ converges almost surely to an element $\widetilde{\<PsiIPsi3>}_{T_1,.}$  in the space ${\cac}^{1-\eps}([T_1,T_2]; \cb^{-\eta}_{\infty})$, for all $0<\eps,\eta<\frac12$.
\end{proposition}

\medskip

In brief, our objective now is to derive a new estimate for $\<PsiIPsi3>$, which, combined with~\eqref{L3}, could yield both a proper definition and a sharp control of $\<PsiIPsi3>$ \textit{as a (continuous) function in time}. To be more specific, our main result about the diagram reads as follows.

 \begin{theorem}\label{coro-8.1}
For all $\eta>0$, there exist $\eps>0$ and $\ka>0$ such that for every $p\geq 1$
\begin{equation*}
\sup_{\ell\geq 0}   \mathbb{E} \Big[ \Big\|  \<PsiIPsi3>^{(n+1)}_{\ell,.} -  \<PsiIPsi3>^{(n)}_{\ell,.}\Big\|_{\cac^{\eps}([\ell,\ell+1]; \cb^{-\eta}_{\infty})}^{2p} \Big]\lesssim 2^{-\ka n p}. 
\end{equation*}
Consequently, for all $0\leq T_1<T_2$, the sequence $(  \<PsiIPsi3>^{(n)}_{T_1,.})$ converges almost surely to  an element $   \<PsiIPsi3>_{T_1,.}$ in  $\cac^{\varepsilon}\big([T_1,T_2];\cb^{-\eta}_{\infty}(\R^3)\big)$.

\end{theorem}

\begin{remark}  
By uniqueness of the limit, the element $ \<PsiIPsi3>$ defined above  is related to $\widetilde{ \<PsiIPsi3>}$ by 
\begin{equation*} 
\widetilde{ \<PsiIPsi3>}_{T_1,t}=\int_{T_1}^t  \<PsiIPsi3>_{T_1,s} \, ds.
\end{equation*}
\end{remark}

Let us now examine the transition from Proposition \ref{Prop-p41} to Theorem \ref{coro-8.1}. To this end, we consider the process
$$\<PsiIPsi3nr>^{(n)}_{s,t}:=\<Psi>^{(n)}_{s,t}\<IPsi3>^{(n)}_{s,t}.  $$
In other words, $\<PsiIPsi3nr>^{(n)}$ stands for the full product of the two diagrams $\<Psi>^{(n)}$ and $\<IPsi3>^{(n)}$, to be distinguished from the resonant product $\<PsiIPsi3>^{(n)}$ in \eqref{ord4}. The result of Theorem \ref{coro-8.1} will then be an easy consequence of the following property:

\begin{proposition}\label{pro-anex8}
There exists $\ka>0$ such that for every $p\geq 1$,
\begin{equation}\label{L30}
\sup_{\ell\geq 0}    \mathbb{E} \Big[ \Big\|\<PsiIPsi3nr>^{(n+1)}_{\ell,.}   -  \<PsiIPsi3nr>^{(n)}_{\ell,.}\Big\|_{{\cac}^{\frac1{10}}([\ell,\ell+2]; \cb^{-2}_\infty)}^{2p} \Big]\lesssim 2^{-\ka n p }. 
\end{equation}
\end{proposition}

\

Before we proceed with the proof of Proposition \ref{pro-anex8}, we show how it yields our main statement by interpolation.

\begin{proof}[Proof of Theorem \ref{coro-8.1}]  
Recall that
$$\<PsiIPsi3>^{(n)}_{s,t}=\<PsiIPsi3nr>^{(n)}_{s,t}-\<Psi>^{(n)}_{s,t}\pl \<IPsi3>^{(n)}_{s,t}-\<IPsi3>^{(n)}_{s,t}\pl\<Psi>^{(n)}_{s,t}.$$
Therefore, using the properties of the paraproduct (see Proposition \ref{Prop-est-para}) together with the results of \cite[Proposition 5.1]{DFT} and \cite[Proposition 7.1]{DFT}, it is readily checked that the estimate in \eqref{L30} also applies to the resonant product $\<PsiIPsi3>^{(n)}$, that is, one has
\begin{equation}\label{L3-bis}
\sup_{\ell\geq 0}    \mathbb{E} \Big[ \Big\|\<PsiIPsi3>^{(n+1)}_{\ell,.}   -  \<PsiIPsi3>^{(n)}_{\ell,.}\Big\|_{{\cac}^{\frac1{10}}([\ell,\ell+2]; \cb^{-2}_{\infty})}^{2p} \Big]\lesssim 2^{-\ka n p }. 
\end{equation}
Let $0<\eta<1$ and set $\theta=\frac{\eta}{2-\eta}<1$. Then for $\eps>0$ small enough,  one has $\lambda(\theta) :=\frac{\theta}{10}-\eps(1-\theta)>0$, and Lemma \ref{lemma-interpt} yields, for all $\ell\geq 0$,
\begin{multline*} 
\big\|  \<PsiIPsi3>^{(n+1)}_{\ell,.}-  \<PsiIPsi3>^{(n)}_{\ell,.}  \big\|_{\mathcal{C}^{\lambda(\theta)}([\ell,\ell+1]; \mathcal{B}_{\infty}^{-2\eta})} \lesssim \\
\lesssim \big\| \widetilde{ \<PsiIPsi3>}_{\ell,.}^{(n+1)}-\widetilde{ \<PsiIPsi3>}_{\ell,.}^{(n)} \big\|^{1-\theta}_{{\ov \cac}^{1-\eps}([\ell,\ell+2]; \mathcal{B}_{\infty}^{-\eta})} \big\|  \<PsiIPsi3>^{(n+1)}_{\ell,.}-  \<PsiIPsi3>^{(n)}_{\ell,.}  \big\|^{\theta}_{\mathcal{C}^{\frac1{10}}([\ell,\ell+2]; \mathcal{B}_{\infty}^{-2})}.
\end{multline*}
The conclusion now comes from the combination of \eqref{L3} and \eqref{L3-bis}.
\end{proof}

\medskip

We now turn to the proof of Proposition \ref{pro-anex8}. For the sake of conciseness, we will only focus on the uniform bound
\begin{equation*}
 \sup_{n\geq 1} \sup_{\ell\geq 0} \, \mathbb{E} \Big[ \big\| \<PsiIPsi3nr>^{(n)}_{\ell,.}  \big\|_{{\cac}^{\frac1{10}}([\ell,\ell+2];\cb_{\infty}^{-2})}^{2p} \Big] <\infty.
\end{equation*}
Note that, due to the stationarity property stated in Lemma \ref{lem:stati-z}, this reduces to proving that
\begin{equation*}
 \sup_{n\geq 1}  \, \mathbb{E} \Big[ \big\| \<PsiIPsi3nr>^{(n)}_{0,.}  \big\|_{{\cac}^{\frac1{10}}([0,2];\cb_{\infty}^{-2})}^{2p} \Big] <\infty.
\end{equation*}

\

 Fix  $\frac34<\al<1$. For every $p\geq 1$, one has, by \cite[Proposition B.6]{DFT},
 \begin{equation*}
\big\|\<PsiIPsi3nr>^{(n)}_{0,.}\big\|_{\cac^{\frac{1}{10}}([0,2];\cb_{\infty}^{-2})}^{2p}\lesssim \int_{0}^{2}\int_{0}^{2} dt_1 dt_2 \, \frac{\big\| \<PsiIPsi3nr>^{(n)}_{0,t_2}-\<PsiIPsi3nr>^{(n)}_{0,t_1}\big\|_{\cb_{\infty}^{-2}}^{2p}}{|t_2-t_1|^{\frac{p}{5} +2}},
\end{equation*}
and then we can use the Sobolev embedding 
 \begin{equation}\label{sobo-beso}
 \| u\|_{\mathcal{B}^{-2} _{\infty,\infty}(\R^3) } \leq C \| u\|_{\mathcal{W}^{-2\al,p}(\R^3) } 
  \end{equation}
to assert that for every $p\geq 1$ large enough,
 \begin{equation}
\big\|\<PsiIPsi3nr>^{(n)}_{0,.} \big\|_{\cac^{\frac{1}{10}}([0,2];\cb_{\infty}^{-2})}^{2p} \lesssim \int_{0}^{2}\int_{0}^{2} dt_1 dt_2 \, \frac{\big\|H^{-\al}\big(\<PsiIPsi3nr>^{(n)}_{0,t_2}-\<PsiIPsi3nr>^{(n)}_{0,t_1}\big)\big\|_{L^{2p}_x}^{2p}}{|t_2-t_1|^{\frac{ p}{5}+2}}.\label{normlip-40}
\end{equation}

Now, using the hypercontractivity property of Wiener chaoses, we obtain that
\begin{align}
\mathbb{E}\Big[\big\|\<PsiIPsi3nr>^{(n)}_{0,.} \big\|_{\cac^{\frac{1}{10}}([0,2];\cb_{\infty}^{-2})}^{2p}\Big]
& \lesssim  \int_{[0,2]^2} \frac{dt_1 dt_2 }{|t_2-t_1|^{\frac{ p}{5}+2}}\int dx \, \mathbb{E}\bigg[ \Big| H^{-\al}\big(\<PsiIPsi3nr>^{(n)}_{0,t_2}-\<PsiIPsi3nr>^{(n)}_{0,t_1}\big)(x)\Big|^{2p}\bigg]\nonumber\\
& \lesssim  \int_{[0,2]^2} \frac{dt_1 dt_2 }{|t_2-t_1|^{\frac{ p}{5}+2}}\int dx \, \mathbb{E}\bigg[ \Big| H^{-\al}\big(\<PsiIPsi3nr>^{(n)}_{0,t_2}-\<PsiIPsi3nr>^{(n)}_{0,t_1}\big)(x)\Big|^{2}\bigg]^p.\label{hypercon}
\end{align}
Thus, it suffices to show that
\begin{equation}\label{hyperc-40}
\sup_{n\geq 1}\int dx \, \mathbb{E}\bigg[ \Big| H^{-\al}\big(\<PsiIPsi3nr>^{(n)}_{0,t_2}-\<PsiIPsi3nr>^{(n)}_{0,t_1}\big)(x)\Big|^{2}\bigg]^p \lesssim |t_2-t_1|^{\ka p}  
\end{equation}
uniformly over $0\leq t_1,t_2\leq 2$, for some $\ka>\frac{1}{5}$. In fact, we will prove that~\eqref{hyperc-40} holds with $\kappa=\frac14$, uniformly over all $t_1,t_2\geq 0$.

\

\subsection{Proof of \eqref{hyperc-40}}

Using the notation in \eqref{inc-k}, we can write, for all $0 \leq t_1\leq t_2\leq 2$, 
\begin{multline*}
\big(\<PsiIPsi3nr>^{(n)}_{0,t_2}-\<PsiIPsi3nr>^{(n)}_{0,t_1}\big)(y)=\<Psi>^{(n)}_{0,t_2}(y)\<IPsi3>^{(n)}_{0,t_2}(y)-\<Psi>^{(n)}_{0,t_1}(y)\<IPsi3>^{(n)}_{0,t_1}(y) \\
\begin{aligned}
&=\big(\<Psi>^{(n)}_{0,t_2}-\<Psi>^{(n)}_{0,t_1}\big)(y)\<IPsi3>^{(n)}_{0,t_2}(y)+\<Psi>^{(n)}_{0,t_1}(y)\big(\<IPsi3>^{(n)}_{0,t_2}(y)-\<IPsi3>^{(n)}_{0,t_1}(y)\big)\\
&=\big(\<Psi>^{(n)}_{0,t_2}-\<Psi>^{(n)}_{0,t_1}\big)(y)\<IPsi3>^{(n)}_{0,t_2}(y)\\
&\hspace{1cm}+\<Psi>^{(n)}_{0,t_1}(y)\Big(\int_{0}^{t_2} ds \int dw \, K_{t_2-s}(y,w)  \<Psi3>^{(n)}_{0,s}(w)-\int_{0}^{t_1} ds \int dw \, K_{t_1-s}(y,w)  \<Psi3>^{(n)}_{0,s}(w)\Big)\\
&=\big(\<Psi>^{(n)}_{0,t_2}-\<Psi>^{(n)}_{0,t_1}\big)(y)\<IPsi3>^{(n)}_{0,t_2}(y)+\<Psi>^{(n)}_{0,t_1}(y)\Big(\int_{t_1}^{t_2} ds \int dw \, K_{t_2-s}(y,w)  \<Psi3>^{(n)}_{0,s}(w)\Big)\\
&\hspace{5cm}+\<Psi>^{(n)}_{0,t_1}(y)\Big(\int_{0}^{t_1} ds \int dw \,  K_{t_1-s,t_2-s}(y,w)  \<Psi3>^{(n)}_{0,s}(w)\Big)\\
&=:\mathfrak{I}^{(n)}_{t_1,t_2}(y)+\mathfrak{II}^{(n)}_{t_1,t_2}(y)+\mathfrak{III}^{(n)}_{t_1,t_2}(y).
\end{aligned}
\end{multline*}

\

We will successively show that
\begin{equation}\label{boun-frak-i}
\sup_{n\geq 1}\int dx \, \mathbb{E}\bigg[ \Big| H^{-\al}\big(\mathfrak{I}^{(n)}_{t_1,t_2}\big)(x)\Big|^{2}\bigg]^p \lesssim |t_2-t_1|^{\frac{p}{4}},
\end{equation}
\begin{equation}\label{boun-frak-ii}
\sup_{n\geq 1}\int dx \, \mathbb{E}\bigg[ \Big| H^{-\al}\big(\mathfrak{II}^{(n)}_{t_1,t_2}\big)(x)\Big|^{2}\bigg]^p \lesssim |t_2-t_1|^{\frac{p}{4}}
\end{equation}
and
\begin{equation}\label{boun-frak-iii}
\sup_{n\geq 1}\int dx \, \mathbb{E}\bigg[ \Big| H^{-\al}\big(\mathfrak{III}^{(n)}_{t_1,t_2}\big)(x)\Big|^{2}\bigg]^p \lesssim |t_2-t_1|^{\frac{p}{4}}
\end{equation}
uniformly over $0\leq t_1\leq t_2$. Observe that once endowed with \eqref{boun-frak-i}-\eqref{boun-frak-ii}-\eqref{boun-frak-iii}, the derivation of \eqref{hyperc-40} (with $\ka=\frac14$) is straightforward.

\

\noindent
\textbf{$\mathfrak{I}$. Proof of \eqref{boun-frak-i}.} One has
\begin{eqnarray*}
\mathfrak{I}^{(n)}_{t_1,t_2}(y)
& =&\int_{0}^{t_2}ds \int dw \, K_{t_2-s}(y,w)\big(\<Psi>^{(n)}_{0,t_2}-\<Psi>^{(n)}_{0,t_1}\big)(y)\<Psi3>^{(n)}_{0,s}(w)\\
&=&\int_{0}^{t_2}ds \int dw \, K_{t_2-s}(y,w)I^W_1\big(F^{(n)}_{t_2,y}-F^{(n)}_{t_1,y}\big)I^W_3\big(F^{(n)}_{s,w} \otimes F^{(n)}_{s,w} \otimes F^{(n)}_{s,w} \big),
\end{eqnarray*}
which, by applying the standard product rule for multiple integrals, yields the decomposition
\begin{align*}
&\mathfrak{I}^{(n)}_{t_1,t_2}(y)=\mathfrak{I}^{\mathbf{1},(n)}_{t_1,t_2} (y) +\mathfrak{I}^{\mathbf{2},(n)}_{t_1,t_2} (y) 
\end{align*}
with
\begin{align*}
\mathfrak{I}^{\mathbf{1},(n)}_{t_1,t_2} (y) :=\int_{0}^{t_2}ds \int dw \, K_{t_2-s}(y,w)I^W_4\Big(\big(F^{(n)}_{t_2,y}-F^{(n)}_{t_1,y}\big) \otimes  \big(F^{(n)}_{s,w} \otimes F^{(n)}_{s,w} \otimes F^{(n)}_{s,w} \big)\Big)
\end{align*}
and
\begin{align*}
\mathfrak{I}^{\mathbf{2},(n)}_{t_1,t_2} (y) :=3 \int_{0}^{t_2}ds \int dw \, K_{t_2-s}(y,w)\,\cac^{(n)}_{(t_1,t_2),s}(y,w) \,  I^W_2\big( F^{(n)}_{s,w} \otimes F^{(n)}_{s,w} \big).
\end{align*}
Recall that the above notation $\cac^{(n)}_{(t_1,t_2),s}$ has been introduced in \eqref{inc-cn}.

\smallskip

For $\mathbf{a}=1,2$, let us set
$$\cm^{\mathbf{a},(n)}_{t_1,t_2}(y_1,y_2):=\mathbb{E}\Big[\mathfrak{I}^{\mathbf{a},(n)}_{t_1,t_2}(y_1)\, \mathfrak{I}^{\mathbf{a},(n)}_{t_1,t_2}(y_2) \Big].$$
Then, thanks to Lemma \ref{lem:techn} and since $\al>\frac34$, we deduce that for $\mathbf{a}=1,2$,
\begin{equation}\label{contr-1-i}
\int dx \, \mathbb{E}\bigg[ \Big| H^{-\al}\big(\mathfrak{I}^{\mathbf{a},(n)}_{t_1,t_2}\big)(x)\Big|^{2}\bigg]^p\lesssim \bigg(\int dy_1 \,  \big|\cm^{\mathbf{a},(n),\star}_{t_1,t_2}(y_1)\big| \bigg)^p,
\end{equation}
for every $p$ large enough, where we have used the notation introduced in \eqref{star-not}.

\

\noindent
\textbf{$\mathfrak{I}$.1. Study of $\mathfrak{I}^{\mathbf{1},(n)}$.} One has here
\begin{multline*}
\cm^{\mathbf{1},(n)}_{t_1,t_2}(y_1,y_2)=\mathbb{E}\Big[\mathfrak{I}^{\mathbf{1},(n)}_{t_1,t_2}(y_1)\, \mathfrak{I}^{\mathbf{1},(n)}_{t_1,t_2}(y_2) \Big]=\\
\begin{aligned}
&=\int_{0}^{t_2}ds_1\int_{0}^{t_2}ds_2 \int dw_1 dw_2 \, K_{t_2-s_1}(y_1,w_1)K_{t_2-s_2}(y_2,w_2)\\
&\hspace{3cm} \times \mathbb{E}\Big[I^W_4\Big(\big(F^{(n)}_{t_2,y_1}-F^{(n)}_{t_1,y_1}\big) \otimes  \big(F^{(n)}_{s_1,w_1} \otimes F^{(n)}_{s_1,w_1} \otimes F^{(n)}_{s_1,w_1} \big)\Big)\\
&\hspace{5cm} \times I^W_4\Big(\big(F^{(n)}_{t_2,y_2}-F^{(n)}_{t_1,y_2}\big) \otimes  \big(F^{(n)}_{s_2,w_2} \otimes F^{(n)}_{s_2,w_2} \otimes F^{(n)}_{s_2,w_2} \big)\Big) \Big].
\end{aligned}
\end{multline*}
The latter expectation can be expanded as
\begin{align*}
&\mathbb{E}\Big[I^W_4\Big(\big(F^{(n)}_{t_2,y_1}-F^{(n)}_{t_1,y_1}\big) \otimes  \big(F^{(n)}_{s_1,w_1} \otimes F^{(n)}_{s_1,w_1} \otimes F^{(n)}_{s_1,w_1} \big)\Big)I^W_4\Big(\big(F^{(n)}_{t_2,y_2}-F^{(n)}_{t_1,y_2}\big) \otimes  \big(F^{(n)}_{s_2,w_2} \otimes F^{(n)}_{s_2,w_2} \otimes F^{(n)}_{s_2,w_2} \big)\Big) \Big]\\
& =c\, \Big\langle \text{Sym}\Big(\big(F^{(n)}_{t_2,y_1}-F^{(n)}_{t_1,y_1}\big) \otimes  \big(F^{(n)}_{s_1,w_1} \otimes F^{(n)}_{s_1,w_1} \otimes F^{(n)}_{s_1,w_1} \big)\Big), \\
 &\hspace{4cm}   \text{Sym}\Big(\big(F^{(n)}_{t_2,y_2}-F^{(n)}_{t_1,y_2}\big) \otimes  \big(F^{(n)}_{s_2,w_2} \otimes F^{(n)}_{s_2,w_2} \otimes F^{(n)}_{s_2,w_2} \big)\Big) \Big\rangle_{L^2((\R_+\times \R^3)^4)}\\
&= c_{\mathbf{1}} \cq^{\mathbf{1},\mathbf{1},(n)}_{t_1,t_2,s_1,s_2}({y,w})+c_{\mathbf{2}} \cq^{\mathbf{1},\mathbf{2},(n)}_{t_1,t_2,s_1,s_2}({y,w})
\end{align*}
for some combinatorial coefficients $c,c_{\mathbf{1}},c_{\mathbf{2}}\geq 0$, and with (recall the notation in \eqref{star-not})
\begin{align*}
 \cq^{\mathbf{1},\mathbf{1},(n)}_{t_1,t_2,s_1,s_2}(y,w)&:=\cac^{(n)}_{(t_1,t_2),(t_1,t_2)}(y_1,y_2)   \cac^{(n)}_{s_1,s_2}(w_1,w_2)^3,\\[2pt]
 \cq^{\mathbf{1},\mathbf{2},(n)}_{t_1,t_2,s_1,s_2}(y,w)&:=\cac^{(n)}_{(t_1,t_2),s_2}(y_1,w_2) \cac^{(n)}_{(t_1,t_2),s_1}(y_2,w_1)\cac^{(n)}_{s_1,s_2}(w_1,w_2)^2.
\end{align*}
Thus,
\begin{align*}
\cm^{\mathbf{1},(n)}_{t_1,t_2}(y_1,y_2)= c_{\mathbf{1}}\cm^{\mathbf{1,1},(n)}_{t_1,t_2}(y_1,y_2)+ c_{\mathbf{2}}\cm^{\mathbf{1,2},(n)}_{t_1,t_2}(y_1,y_2),
\end{align*}
with
\begin{multline*}
\cm^{\mathbf{1,b},(n)}_{t_1,t_2}(y_1,y_2):=\int_{0}^{t_2}ds_1\int_{0}^{t_2}ds_2\int dw_1 dw_2 \, K_{t_2-s_1}(y_1,w_1) K_{t_2-s_2}(y_2,w_2)\cq^{\mathbf{1},\mathbf{b},(n)}_{t_1,t_2,s_1,s_2}(y,w),
\end{multline*}
and accordingly
\begin{equation}\label{contr-2-i}
\int dy_1 \, \big|\cm^{\mathbf{1},(n),\star}_{t_1,t_2}(y_1)\big| \lesssim \int dy_1 \, \big|\cm^{\mathbf{1,1},(n),\star}_{t_1,t_2}(y_1)\big| +\int dy_1 \,  \big|\cm^{\mathbf{1,2},(n),\star}_{t_1,t_2}(y_1)\big| .
\end{equation}

\

\noindent
\textbf{$\mathfrak{I}$.1.1. \underline{Study of $\cm^{\mathbf{1,1},(n)}$}.} One has for all $y_1,y_2$,
\begin{align*} 
&\big|\cm^{\mathbf{1,1},(n)}_{t_1,t_2}(y_1+y_2,y_2)\big| \lesssim \\
&\lesssim {\cac}^{(n)}_{(t_1,t_2),(t_1,t_2)}(y_1+y_2,y_2)\int_{0}^{t_2}ds_1\int_{0}^{t_2}ds_2 \int dw_1\, \big|K_{t_2-s_1}(y_1+y_2,w_1) \big| \int dw_2 \, \big| K_{t_2-s_2}(y_2,w_2) \big| \cac^{(n)}_{s_1,s_2}(w_1,w_2)^3  \\
&\lesssim {\cac}^{(n),\star}_{(t_1,t_2),(t_1,t_2)}(y_1)\int_{0}^{t_2}ds_1\int_{0}^{t_2}ds_2 \int dw_1\, \big|K_{t_2-s_1}(y_1+y_2,w_1+y_2)\big|   \\
&\hspace{4cm}\times \int dw_2 \, \big|K_{t_2-s_2}(y_2,w_2+y_2)\big|  \cac^{(n)}_{s_1,s_2}(w_1+y_2,w_2+y_2)^3\\
&\lesssim {\cac}^{(n),\star}_{(t_1,t_2),(t_1,t_2)}(y_1)\int_{0}^{t_2}ds_1\int_{0}^{t_2}ds_2 \int dw_1\, K_{t_2-s_1}^\star(y_1-w_1) \int dw_2 \, K_{t_2-s_2}^\star(w_2)  \cac^{(n),\star}_{s_1,s_2}(w_1-w_2)^3\\
&\lesssim {\cac}^{(n),\star}_{(t_1,t_2),(t_1,t_2)}(y_1)\int_{0}^{t_2}ds_1\int_{0}^{t_2}ds_2\,  \Big[ K_{t_2-s_1}^\star  \ast \big[ K_{t_2-s_2}^\star \ast  (\cac^{(n),\star}_{s_1,s_2})^3\big]\Big](y_1).
\end{align*}
As a result,
\begin{align} 
\int dy_1 \, \big|\cm^{\mathbf{1,1},(n),\star}_{t_1,t_2}(y_1)\big|
&\lesssim \big\|{\cac}^{(n),\star}_{(t_1,t_2),(t_1,t_2)}\big\|_{L^2} \int_{0}^{t_2}ds_1\int_{0}^{t_2}ds_2\,  \Big\| K_{t_2-s_1}^\star  \ast \big[ K_{t_2-s_2}^\star \ast  (\cac^{(n),\star}_{s_1,s_2})^3\big]\Big\|_{L^2}\nonumber\\
&\lesssim \big\|{\cac}^{(n),\star}_{(t_1,t_2),(t_1,t_2)}\big\|_{L^2} \int_{0}^{t_2}ds_1\int_{0}^{t_2}ds_2\,  \big\| K_{t_2-s_1}^\star  \big\|_{L^2}\big\|  K_{t_2-s_2}^\star \ast  (\cac^{(n),\star}_{s_1,s_2})^3\big\|_{L^1}\nonumber\\
&\lesssim \big\|{\cac}^{(n),\star}_{(t_1,t_2),(t_1,t_2)}\big\|_{L^2}\int_{0}^{t_2}ds_1\int_{0}^{t_2}ds_2\,  \big\| K_{t_2-s_1}^\star  \big\|_{L^2}\big\|  K_{t_2-s_2}^\star\big\|_{L^1} \big\|\cac^{(n),\star}_{s_1,s_2}\big\|_{L^3}^3.\label{ref-1-1-papi}
\end{align}
We can then appeal to the estimates of Lemma \ref{lem:k-star-0} and Lemma \ref{lem:cac3}, which yield for all $\eps>0$,
\begin{align*} 
\int dy_1 \, \big|\cm^{\mathbf{1,1},(n),\star}_{t_1,t_2}(y_1)\big|
&\lesssim |t_2-t_1|^{\frac14}  \int_{0}^{t_2}\int_{0}^{t_2} \frac{ds_1ds_2}{|s_1-s_2|^\eps}\frac{e^{-(t_2-s_1)}}{|t_2-s_1|^{\frac34}}  e^{-(t_2-s_2)} \\
&\lesssim |t_2-t_1|^{\frac14}  \int_{0}^{+\infty}\int_{0}^{+\infty} \frac{dr_1dr_2}{|r_1-r_2|^\eps}\frac{e^{-r_1}}{|r_1|^{\frac34}}  e^{-r_2} ,
\end{align*}
and we can conclude that
\begin{equation}\label{cm-1-1-i} 
\int dy_1 \, \big|\cm^{\mathbf{1,1},(n),\star}_{t_1,t_2}(y_1)\big|\lesssim |t_2-t_1|^{\frac14}  ,
\end{equation}
uniformly over $t_1,t_2\geq 0$.

\

\noindent
\textbf{$\mathfrak{I}$.1.2. \underline{Study of $\cm^{\mathbf{1,2},(n)}$}.} Arguing as above, we obtain that
\begin{multline*}
\big|\cm^{\mathbf{1,2},(n)}_{t_1,t_2}(y_1+y_2,y_2)\big| \lesssim \\
\begin{aligned}
& \lesssim  \int_{0}^{t_2}ds_1\int_{0}^{t_2}ds_2\int dw_1 \, \big|K_{t_2-s_1}(y_1+y_2,w_1+y_2) \big|  \cac^{(n)}_{(t_1,t_2),s_1}(y_2,w_1+y_2)\\
&\hspace{1cm}\times \int dw_2\,  \big| K_{t_2-s_2}(y_2,w_2+y_2) \big| \cac^{(n)}_{(t_1,t_2),s_2}(y_1+y_2,w_2+y_2) \cac^{(n)}_{s_1,s_2}(w_1+y_2,w_2+y_2)^2\\
& \lesssim  \int_{0}^{t_2}ds_1\int_{0}^{t_2}ds_2\int dw_1 \, K_{t_2-s_1}^\star(y_1-w_1)  \cac^{(n),\star}_{(t_1,t_2),s_1}(w_1)\\
&\hspace{1cm}\times \int dw_2\,   K_{t_2-s_2}^\star(w_2)  \cac^{(n),\star}_{(t_1,t_2),s_2}(y_1-w_2) \cac^{(n),\star}_{s_1,s_2}(w_1-w_2)^2,
\end{aligned}
\end{multline*}
and so
\begin{multline}\label{ref-1-2-papi}
\int dy_1 \, \big|\cm^{\mathbf{1,2},(n),\star}_{t_1,t_2}(y_1)\big|\lesssim \\
\begin{aligned}
&\lesssim \int_{0}^{t_2}ds_1\int_{0}^{t_2}ds_2\, \big\|\cac^{(n),\star}_{(t_1,t_2),s_2}\big\|_{L^2} \big\| K_{t_2-s_1}^\star\big\|_{L^2}\int dw_1 \,  \cac^{(n),\star}_{(t_1,t_2),s_1}(w_1) \int dw_2\,   K_{t_2-s_2}^\star(w_2)  \cac^{(n),\star}_{s_1,s_2}(w_1-w_2)^2 \\
&\lesssim \int_{0}^{t_2}ds_1\int_{0}^{t_2}ds_2\, \big\|\cac^{(n),\star}_{(t_1,t_2),s_2}\big\|_{L^2} \big\| K_{t_2-s_1}^\star\big\|_{L^2}\int dw_1 \,  \cac^{(n),\star}_{(t_1,t_2),s_1}(w_1) \big[  K_{t_2-s_2}^\star \ast  (\cac^{(n),\star}_{s_1,s_2})^2\big](w_1)\\
&\lesssim \int_{0}^{t_2}ds_1\int_{0}^{t_2}ds_2\, \big\|\cac^{(n),\star}_{(t_1,t_2),s_2}\big\|_{L^2} \big\| K_{t_2-s_1}^\star\big\|_{L^2}\big\|\cac^{(n),\star}_{(t_1,t_2),s_1}\big\|_{L^2} \big\| K_{t_2-s_2}^\star\big\|_{L^2} \big\| \cac^{(n),\star}_{s_1,s_2}\big\|_{L^2}^2.
\end{aligned}
\end{multline}
Using the estimates contained in Lemma \ref{lem:k-star-0} and Lemma \ref{lem:cac3}, we deduce 
\begin{align*}
\int dy_1 \, \big|\cm^{\mathbf{1,2},(n),\star}_{t_1,t_2}(y_1)\big|&\lesssim |t_2-t_1|^{\frac14} \bigg(\int_{0}^{t_2}ds\, \frac{e^{-(t_2-s)}}{|t_2-s|^{\frac34}} \bigg)^2\lesssim |t_2-t_1|^{\frac14}\bigg(\int_0^{+\infty}dr\, \frac{e^{-r}}{ |r|^{\frac34}}\bigg)^2 ,
\end{align*}
and so we can conclude again that
\begin{equation}\label{cm-1-2-i} 
\int dy_1 \, \big|\cm^{\mathbf{1,2},(n),\star}_{t_1,t_2}(y_1)\big| \lesssim |t_2-t_1|^{\frac14} ,
\end{equation}
uniformly over $t_1,t_2\geq 0$.

\

\noindent
\textbf{$\mathfrak{I}$.2. Study of $\mathfrak{I}^{\mathbf{2},(n)}$.} The quantity $\cm^{\mathbf{2},(n)}_{t_1,t_2}$ can be computed as follows:
\begin{align*} 
&\cm^{\mathbf{2},(n)}_{t_1,t_2}(y_1,y_2)=\mathbb{E}\Big[\mathfrak{I}^{\mathbf{2},(n)}_{t_1,t_2}(y_1)\, \mathfrak{I}^{\mathbf{2},(n)}_{t_1,t_2}(y_2) \Big]=\\
&=9\int_{0}^{t_2}ds_1\int_{0}^{t_2}ds_2 \int dw_1 dw_2 \, K_{t_2-s_1}(y_1,w_1)K_{t_2-s_2}(y_2,w_2)\\
&\hspace{2cm}\times  \cac^{(n)}_{(t_1,t_2),s_1}(y_1,w_1) \cac^{(n)}_{(t_1,t_2),s_2}(y_2,w_2)\mathbb{E}\Big[ I^W_2\big( F^{(n)}_{s_1,w_1} \otimes F^{(n)}_{s_1,w_1} \big)I^W_2\big( F^{(n)}_{s_2,w_2} \otimes F^{(n)}_{s_2,w_2} \big)\Big] \\
&=c\int_{0}^{t_2}ds_1\int_{0}^{t_2}ds_2 \int dw_1 dw_2 \, K_{t_2-s_1}(y_1,w_1)K_{t_2-s_2}(y_2,w_2) \\
&\hspace{5cm}\times  \cac^{(n)}_{(t_1,t_2),s_1}(y_1,w_1) \cac^{(n)}_{(t_1,t_2),s_2}(y_2,w_2) \cac^{(n)}_{s_1,s_2}(w_1,w_2)^2.
\end{align*}
We can then write
\begin{align*} 
&\cm^{\mathbf{2},(n)}_{t_1,t_2}(y_1+y_2,y_2)=\\
&=\int_{0}^{t_2}ds_1\int_{0}^{t_2}ds_2 \int dw_1 dw_2 \, K_{t_2-s_1}(y_1+y_2,w_1+y_2)K_{t_2-s_2}(y_2,w_2+y_2) \\
&\hspace{3cm}\times  \cac^{(n)}_{(t_1,t_2),s_1}(y_1+y_2,w_1+y_2) \cac^{(n)}_{(t_1,t_2),s_2}(y_2,w_2+y_2) \cac^{(n)}_{s_1,s_2}(w_1+y_2,w_2+y_2)^2\\
&\lesssim \int_{0}^{t_2}ds_1\int_{0}^{t_2}ds_2 \int dw_1 \, K^\star_{t_2-s_1}(y_1-w_1)\cac^{(n),\star}_{(t_1,t_2),s_1}(y_1-w_1)  \\
&\hspace{3cm}\times\int dw_2\, \cac^{(n),\star}_{s_1,s_2}(w_1-w_2)^2 K^\star_{t_2-s_2}(w_2) \cac^{(n),\star}_{(t_1,t_2),s_2}(w_2) \\
&\lesssim \int_{0}^{t_2}ds_1\int_{0}^{t_2}ds_2\, \bigg[ \big( K^\star_{t_2-s_1}\cdot \cac^{(n),\star}_{(t_1,t_2),s_1}\big)\ast\Big[(\cac^{(n),\star}_{s_1,s_2})^2 \ast  \big(K^\star_{t_2-s_2}\cdot \cac^{(n),\star}_{(t_1,t_2),s_2}\big)\Big]\bigg](y_1) .
\end{align*}
As a result,
\begin{align} 
&\int dy_1\, \cm^{\mathbf{2},(n),\star}_{t_1,t_2}(y_1)\lesssim \int_{0}^{t_2}ds_1\int_{0}^{t_2}ds_2\,  \big\| K^\star_{t_2-s_1}\cdot \cac^{(n),\star}_{(t_1,t_2),s_1}\big\|_{L^1}\big\|(\cac^{(n),\star}_{s_1,s_2})^2 \ast  \big(K^\star_{t_2-s_2}\cdot \cac^{(n),\star}_{(t_1,t_2),s_2}\big)\big\|_{L^1}\nonumber\\
&\lesssim \int_{0}^{t_2}ds_1\int_{0}^{t_2}ds_2\,  \big\| K^\star_{t_2-s_1}\big\|_{L^2}\big\| \cac^{(n),\star}_{(t_1,t_2),s_1}\big\|_{L^2}\big\|(\cac^{(n),\star}_{s_1,s_2})^2 \big\|_{L^1}\big\|K^\star_{t_2-s_2}\big\|_{L^2}\big\| \cac^{(n),\star}_{(t_1,t_2),s_2}\big\|_{L^2}.\label{ref-2-papi}
\end{align}
Thanks to Lemma \ref{lem:k-star-0} and Lemma \ref{lem:cac3}, we obtain this time
\begin{equation}\label{cm-2-i} 
\int dy_1\, \cm^{\mathbf{2},(n),\star}_{t_1,t_2}(y_1)\lesssim |t_2-t_1|^{\frac14} \bigg(\int_{0}^{t_2}ds\, \frac{e^{-(t_2-s)}}{|t_2-s|^{\frac34}} \bigg)^2\lesssim |t_2-t_1|^{\frac14} ,
\end{equation}
uniformly over $t_1,t_2\geq 0$.

\

Finally, by inserting \eqref{cm-1-1-i}-\eqref{cm-1-2-i}-\eqref{cm-2-i} into \eqref{contr-1-i} and \eqref{contr-2-i}, we obtain \eqref{boun-frak-i}.

\

\

\noindent
\textbf{$\mathfrak{II}$. Proof of \eqref{boun-frak-ii}.} One has
\begin{eqnarray*}
\mathfrak{II}^{(n)}_{t_1,t_2}(y)
& =&\int_{t_1}^{t_2}ds \int dw \, K_{t_2-s}(y,w)\<Psi>^{(n)}_{0,t_1}(y)\<Psi3>^{(n)}_{0,s}(w)\\
&=&\int_{t_1}^{t_2}ds \int dw \, K_{t_2-s}(y,w)I^W_1\big(F^{(n)}_{t_1,y}\big)I^W_3\big(F^{(n)}_{s,w} \otimes F^{(n)}_{s,w} \otimes F^{(n)}_{s,w} \big),
\end{eqnarray*}
which, by applying the standard product rule for multiple integrals, yields the decomposition
\begin{align*}
&\mathfrak{II}^{(n)}_{t_1,t_2}(y)=\mathfrak{II}^{\mathbf{1},(n)}_{t_1,t_2} (y) +\mathfrak{II}^{\mathbf{2},(n)}_{t_1,t_2} (y) 
\end{align*}
with
\begin{align*}
\mathfrak{II}^{\mathbf{1},(n)}_{t_1,t_2} (y) :=\int_{t_1}^{t_2}ds \int dw \, K_{t_2-s}(y,w)I^W_4\Big(F^{(n)}_{t_1,y} \otimes  \big(F^{(n)}_{s,w} \otimes F^{(n)}_{s,w} \otimes F^{(n)}_{s,w} \big)\Big)
\end{align*}
and
\begin{align*}
\mathfrak{II}^{\mathbf{2},(n)}_{t_1,t_2} (y) :=3 \int_{t_1}^{t_2}ds \int dw \, K_{t_2-s}(y,w)\,\cac^{(n)}_{t_1,s}(y,w) \,  I^W_2\big( F^{(n)}_{s,w} \otimes F^{(n)}_{s,w} \big).
\end{align*}

\smallskip

For $\mathbf{a}=1,2$, let us set
$$\cm^{\mathbf{a},(n)}_{t_1,t_2}(y_1,y_2):=\mathbb{E}\Big[\mathfrak{II}^{\mathbf{a},(n)}_{t_1,t_2}(y_1)\, \mathfrak{II}^{\mathbf{a},(n)}_{t_1,t_2}(y_2) \Big].$$
Then, thanks to Lemma \ref{lem:techn} and since $\al>\frac34$, we deduce that for $\mathbf{a}=1,2$,
\begin{equation}\label{contr-1-ii}
\int dx \, \mathbb{E}\bigg[ \Big| H^{-\al}\big(\mathfrak{II}^{\mathbf{a},(n)}_{t_1,t_2}\big)(x)\Big|^{2}\bigg]^p
\lesssim \bigg(\int dy_1 \,  \big|\cm^{\mathbf{a},(n),\star}_{t_1,t_2}(y_1)\big| \bigg)^p,
\end{equation}
for every $p$ large enough, where we have used the notation introduced in \eqref{star-not}.

\

\noindent
\textbf{$\mathfrak{II}$.1. Study of $\mathfrak{II}^{\mathbf{1},(n)}$.} One has here
\begin{align*}
&\cm^{\mathbf{1},(n)}_{t_1,t_2}(y_1,y_2)=\mathbb{E}\Big[\mathfrak{II}^{\mathbf{1},(n)}_{t_1,t_2}(y_1)\, \mathfrak{II}^{\mathbf{1},(n)}_{t_1,t_2}(y_2) \Big]\\
&=\int_{t_1}^{t_2}ds_1\int_{t_1}^{t_2}ds_2 \int dw_1 dw_2 \, K_{t_2-s_1}(y_1,w_1)K_{t_2-s_2}(y_2,w_2)\\
&\times \mathbb{E}\Big[I^W_4\Big(F^{(n)}_{t_1,y_1} \otimes  \big(F^{(n)}_{s_1,w_1} \otimes F^{(n)}_{s_1,w_1} \otimes F^{(n)}_{s_1,w_1} \big)\Big)I^W_4\Big(F^{(n)}_{t_1,y_2} \otimes  \big(F^{(n)}_{s_2,w_2} \otimes F^{(n)}_{s_2,w_2} \otimes F^{(n)}_{s_2,w_2} \big)\Big) \Big].
\end{align*}
The latter expectation can be expanded as
\begin{align*}
&\mathbb{E}\Big[I^W_4\Big(F^{(n)}_{t_1,y_1} \otimes  \big(F^{(n)}_{s_1,w_1} \otimes F^{(n)}_{s_1,w_1} \otimes F^{(n)}_{s_1,w_1} \big)\Big)I^W_4\Big(F^{(n)}_{t_1,y_2} \otimes  \big(F^{(n)}_{s_2,w_2} \otimes F^{(n)}_{s_2,w_2} \otimes F^{(n)}_{s_2,w_2} \big)\Big) \Big]\\
& =c\, \Big\langle \text{Sym}\Big(F^{(n)}_{t_1,y_1} \otimes  \big(F^{(n)}_{s_1,w_1} \otimes F^{(n)}_{s_1,w_1} \otimes F^{(n)}_{s_1,w_1} \big)\Big), \\
 &\hspace{4cm}   \text{Sym}\Big(F^{(n)}_{t_1,y_2} \otimes  \big(F^{(n)}_{s_2,w_2} \otimes F^{(n)}_{s_2,w_2} \otimes F^{(n)}_{s_2,w_2} \big)\Big) \Big\rangle_{L^2((\R_+\times \R^3)^4)}\\
&= c_{\mathbf{1}} \cq^{\mathbf{1},\mathbf{1},(n)}_{t_1,s_1,s_2}({y,w})+c_{\mathbf{2}} \cq^{\mathbf{1},\mathbf{2},(n)}_{t_1,s_1,s_2}({y,w})
\end{align*}
for some combinatorial coefficients $c,c_{\mathbf{1}},c_{\mathbf{2}}\geq 0$, and with (recall the notation in \eqref{star-not})
\begin{align*}
 \cq^{\mathbf{1},\mathbf{1},(n)}_{t_1,s_1,s_2}(y,w)&:=\cac^{(n)}_{t_1,t_1}(y_1,y_2)   \cac^{(n)}_{s_1,s_2}(w_1,w_2)^3,\\[2pt]
 \cq^{\mathbf{1},\mathbf{2},(n)}_{t_1,s_1,s_2}(y,w)&:=\cac^{(n)}_{t_1,s_2}(y_1,w_2) \cac^{(n)}_{t_1,s_1}(y_2,w_1)\cac^{(n)}_{s_1,s_2}(w_1,w_2)^2.
\end{align*}
Thus,
\begin{align*}
\cm^{\mathbf{1},(n)}_{t_1,t_2}(y_1,y_2)= c_{\mathbf{1}}\cm^{\mathbf{1,1},(n)}_{t_1,t_2}(y_1,y_2)+ c_{\mathbf{2}}\cm^{\mathbf{1,2},(n)}_{t_1,t_2}(y_1,y_2),
\end{align*}
with
\begin{align*}
\cm^{\mathbf{1,b},(n)}_{t_1,t_2}(y_1,y_2):=\int_{t_1}^{t_2}ds_1\int_{t_1}^{t_2}ds_2\int dw_1 dw_2 \, K_{t_2-s_1}(y_1,w_1) K_{t_2-s_2}(y_2,w_2)\cq^{\mathbf{1},\mathbf{b},(n)}_{t_1,s_1,s_2}(y,w),
\end{align*}
and accordingly
\begin{equation}\label{contr-2-ii}
\int dy_1 \, \big|\cm^{\mathbf{1},(n),\star}_{t_1,t_2}(y_1)\big| \lesssim \int dy_1 \, \big|\cm^{\mathbf{1,1},(n),\star}_{t_1,t_2}(y_1)\big| +\int dy_1 \,  \big|\cm^{\mathbf{1,2},(n),\star}_{t_1,t_2}(y_1)\big| .
\end{equation}

\

\noindent
\textbf{$\mathfrak{II}$.1.1. \underline{Study of $\cm^{\mathbf{1,1},(n)}$}.} Arguing exactly as in the derivation of \eqref{ref-1-1-papi}, we obtain
\begin{align*} 
\int dy_1 \, \big|\cm^{\mathbf{1,1},(n),\star}_{t_1,t_2}(y_1)\big|
&\lesssim \big\|{\cac}^{(n),\star}_{t_1,t_1}\big\|_{L^2} \int_{t_1}^{t_2}ds_1\int_{t_1}^{t_2}ds_2\,  \big\| K_{t_2-s_1}^\star  \big\|_{L^2}\big\|  K_{t_2-s_2}^\star\big\|_{L^1} \big\|\cac^{(n),\star}_{s_1,s_2}\big\|_{L^3}^3.
\end{align*}
We can now lean on the estimates of Lemma \ref{lem:k-star-0} and Lemma \ref{lem:cac3}, which yield for every small $\eps>0$,
\begin{align*} 
\int dy_1 \, \big|\cm^{\mathbf{1,1},(n),\star}_{t_1,t_2}(y_1)\big|
&\lesssim \int_{t_1}^{t_2}\int_{t_1}^{t_2} \frac{ds_1ds_2}{|s_1-s_2|^\eps}\frac{e^{-(t_2-s_1)}}{|t_2-s_1|^{\frac34}}  e^{-(t_2-s_2)} \\
&\lesssim |t_2-t_1|^{\frac54-\eps}\int_0^1\int_0^1 \frac{dr_1dr_2}{|r_1-r_2|^\eps}\frac{e^{-(t_2-t_1)(1-r_1)}}{|1-r_1|^{\frac34}}  e^{-(t_2-t_1)(1-r_2)}\\
&\lesssim |t_2-t_1|^{\frac14}\int_0^1\int_0^1 \frac{dr_1dr_2}{|r_1-r_2|^\eps}\frac{1}{|1-r_1|^{\frac78}|1-r_2|^{\frac78-\eps}}  
\end{align*}
and we can conclude that
\begin{equation}\label{cm-1-1-ii} 
\int dy_1 \, \big|\cm^{\mathbf{1,1},(n),\star}_{t_1,t_2}(y_1)\big|\lesssim |t_2-t_1|^{\frac14} ,
\end{equation}
uniformly over $t_1,t_2\geq 0$.

\

\noindent
\textbf{$\mathfrak{II}$.1.2. \underline{Study of $\cm^{\mathbf{1,2},(n)}$}.} Repeating the reasoning that led to \eqref{ref-1-2-papi}, we find
\begin{align*}
&\int dy_1 \, \big|\cm^{\mathbf{1,2},(n),\star}_{t_1,t_2}(y_1)\big|\lesssim \int_{t_1}^{t_2}ds_1\int_{t_1}^{t_2}ds_2\, \big\|\cac^{(n),\star}_{t_1,s_2}\big\|_{L^2} \big\| K_{t_2-s_1}^\star\big\|_{L^2}\big\|\cac^{(n),\star}_{t_1,s_1}\big\|_{L^2} \big\| K_{t_2-s_2}^\star\big\|_{L^2} \big\| \cac^{(n),\star}_{s_1,s_2}\big\|_{L^2}^2.
\end{align*}
Using the estimates contained in Lemma \ref{lem:k-star-0} and Lemma \ref{lem:cac3}, we deduce 
\begin{align*}
\int dy_1 \, \big|\cm^{\mathbf{1,2},(n),\star}_{t_1,t_2}(y_1)\big|&\lesssim \bigg(\int_{t_1}^{t_2}ds\,  \frac{e^{-(t_2-s)}}{|t_2-s|^{\frac34}} \bigg)^2\\
&\lesssim |t_2-t_1|^{\frac12}\bigg( \int_0^1 dr\, \frac{e^{-(t_2-t_1)(1-r)}}{|1-r|^{\frac34}}\bigg)^2\lesssim |t_2-t_1|^{\frac14}\bigg( \int_0^1 \frac{dr}{|1-r|^{\frac78}}\bigg)^2
\end{align*}
and thus we can conclude again that
\begin{equation}\label{cm-1-2-ii}
\int dy_1 \, \big|\cm^{\mathbf{1,2},(n),\star}_{t_1,t_2}(y_1)\big| \lesssim |t_2-t_1|^{\frac14} ,
\end{equation}
uniformly over $t_1,t_2\geq 0$.

\

\noindent
\textbf{$\mathfrak{II}$.2. Study of $\mathfrak{II}^{\mathbf{2},(n)}$.} The quantity $\cm^{\mathbf{2},(n)}_{t_1,t_2}$ can be computed as follows:
\begin{align*} 
&\cm^{\mathbf{2},(n)}_{t_1,t_2}(y_1,y_2)=\mathbb{E}\Big[\mathfrak{II}^{\mathbf{2},(n)}_{t_1,t_2}(y_1)\, \mathfrak{II}^{\mathbf{2},(n)}_{t_1,t_2}(y_2) \Big]=\\
&=9\int_{t_1}^{t_2}ds_1\int_{t_1}^{t_2}ds_2 \int dw_1 dw_2 \, K_{t_2-s_1}(y_1,w_1)K_{t_2-s_2}(y_2,w_2)\\
&\hspace{3cm}\times  \cac^{(n)}_{t_1,s_1}(y_1,w_1) \cac^{(n)}_{t_1,s_2}(y_2,w_2)\mathbb{E}\Big[ I^W_2\big( F^{(n)}_{s_1,w_1} \otimes F^{(n)}_{s_1,w_1} \big)I^W_2\big( F^{(n)}_{s_2,w_2} \otimes F^{(n)}_{s_2,w_2} \big)\Big] \\
&=c\int_{t_1}^{t_2}ds_1\int_{t_1}^{t_2}ds_2 \int dw_1 dw_2 \, K_{t_2-s_1}(y_1,w_1)K_{t_2-s_2}(y_2,w_2) \\
&\hspace{5cm}\times  \cac^{(n)}_{t_1,s_1}(y_1,w_1) \cac^{(n)}_{t_1,s_2}(y_2,w_2) \cac^{(n)}_{s_1,s_2}(w_1,w_2)^2.
\end{align*}
Proceeding along the same lines as for \eqref{ref-2-papi}, we arrive at
\begin{align*} 
&\int dy_1\, \cm^{\mathbf{2},(n),\star}_{t_1,t_2}(y_1)\lesssim \int_{t_1}^{t_2}ds_1\int_{t_1}^{t_2}ds_2\,  \big\| K^\star_{t_2-s_1}\big\|_{L^2}\big\| \cac^{(n),\star}_{t_1,s_1}\big\|_{L^2}\big\|(\cac^{(n),\star}_{s_1,s_2})^2 \big\|_{L^1}\big\|K^\star_{t_2-s_2}\big\|_{L^2}\big\| \cac^{(n),\star}_{t_1,s_2}\big\|_{L^2}.
\end{align*}
Thanks to Lemma \ref{lem:k-star-0} and Lemma \ref{lem:cac3}, we obtain this time
\begin{equation}\label{cm-2-ii} 
\int dy_1\, \cm^{\mathbf{2},(n),\star}_{t_1,t_2}(y_1)\lesssim  \bigg(\int_{t_1}^{t_2}ds\, \frac{e^{-(t_2-s)}}{|t_2-s|^{\frac34}} \bigg)^2\lesssim  \bigg(\int_{t_1}^{t_2}\frac{ds}{|t_2-s|^{\frac78}} \bigg)^2\lesssim |t_2-t_1|^{\frac14} ,
\end{equation}
uniformly over $t_1,t_2\geq 0$.

\

Finally, by inserting \eqref{cm-1-1-ii}-\eqref{cm-1-2-ii}-\eqref{cm-2-ii} into \eqref{contr-1-ii} and \eqref{contr-2-ii}, we obtain \eqref{boun-frak-ii}.

\

\

\noindent
\textbf{$\mathfrak{III}$. Proof of \eqref{boun-frak-iii}.} One has
\begin{eqnarray*}
\mathfrak{III}^{(n)}_{t_1,t_2}(y)
& =&\int_{0}^{t_1}ds \int dw \, K_{t_1-s,t_2-s}(y,w)\<Psi>^{(n)}_{0,t_1}(y)\<Psi3>^{(n)}_{0,s}(w)\\
&=&\int_{0}^{t_1}ds \int dw \, K_{t_1-s,t_2-s}(y,w)I^W_1\big(F^{(n)}_{t_1,y}\big)I^W_3\big(F^{(n)}_{s,w} \otimes F^{(n)}_{s,w} \otimes F^{(n)}_{s,w} \big),
\end{eqnarray*}
which, by applying the standard product rule for multiple integrals, yields the decomposition
\begin{align*}
&\mathfrak{III}^{(n)}_{t_1,t_2}(y)=\mathfrak{III}^{\mathbf{1},(n)}_{t_1,t_2} (y) +\mathfrak{III}^{\mathbf{2},(n)}_{t_1,t_2} (y) 
\end{align*}
with
\begin{align*}
\mathfrak{III}^{\mathbf{1},(n)}_{t_1,t_2} (y) :=\int_{0}^{t_1}ds \int dw \, K_{t_1-s,t_2-s}(y,w)I^W_4\Big(F^{(n)}_{t_1,y} \otimes  \big(F^{(n)}_{s,w} \otimes F^{(n)}_{s,w} \otimes F^{(n)}_{s,w} \big)\Big)
\end{align*}
and
\begin{align*}
\mathfrak{III}^{\mathbf{2},(n)}_{t_1,t_2} (y) :=3 \int_{0}^{t_1}ds \int dw \, K_{t_1-s,t_2-s}(y,w)\,\cac^{(n)}_{t_1,s}(y,w) \,  I^W_2\big( F^{(n)}_{s,w} \otimes F^{(n)}_{s,w} \big).
\end{align*}
For $\mathbf{a}=1,2$, let us set
$$\cm^{\mathbf{a},(n)}_{t_1,t_2}(y_1,y_2):=\mathbb{E}\Big[\mathfrak{III}^{\mathbf{a},(n)}_{t_1,t_2}(y_1)\, \mathfrak{III}^{\mathbf{a},(n)}_{t_1,t_2}(y_2) \Big].$$
Then, thanks to Lemma \ref{lem:techn} and since $\al>\frac34$, we deduce that for $\mathbf{a}=1,2$,
\begin{equation}\label{contr-1-iii}
\int dx \, \mathbb{E}\bigg[ \Big| H^{-\al}\big(\mathfrak{III}^{\mathbf{a},(n)}_{t_1,t_2}\big)(x)\Big|^{2}\bigg]^p\lesssim \bigg(\int dy_1 \,  \big|\cm^{\mathbf{a},(n),\star}_{t_1,t_2}(y_1)\big| \bigg)^p,
\end{equation}
for every $p$ large enough, where we have used the notation introduced in \eqref{star-not}.

\

\noindent
\textbf{$\mathfrak{III}$.1. Study of $\mathfrak{III}^{\mathbf{1},(n)}$.} One has here
\begin{align*}
&\cm^{\mathbf{1},(n)}_{t_1,t_2}(y_1,y_2)=\mathbb{E}\Big[\mathfrak{III}^{\mathbf{1},(n)}_{t_1,t_2}(y_1)\, \mathfrak{III}^{\mathbf{1},(n)}_{t_1,t_2}(y_2) \Big]\\
&=\int_{0}^{t_1}ds_1\int_{0}^{t_1}ds_2 \int dw_1 dw_2 \, K_{t_1-s_1,t_2-s_1}(y_1,w_1)K_{t_1-s_2,t_2-s_2}(y_2,w_2)\\
&\hspace{1cm}\times \mathbb{E}\Big[I^W_4\Big(F^{(n)}_{t_1,y_1} \otimes  \big(F^{(n)}_{s_1,w_1} \otimes F^{(n)}_{s_1,w_1} \otimes F^{(n)}_{s_1,w_1} \big)\Big)I^W_4\Big(F^{(n)}_{t_1,y_2} \otimes  \big(F^{(n)}_{s_2,w_2} \otimes F^{(n)}_{s_2,w_2} \otimes F^{(n)}_{s_2,w_2} \big)\Big) \Big].
\end{align*}
The latter expectation can be expanded as
\begin{multline*}
\mathbb{E}\Big[I^W_4\Big(F^{(n)}_{t_1,y_1} \otimes  \big(F^{(n)}_{s_1,w_1} \otimes F^{(n)}_{s_1,w_1} \otimes F^{(n)}_{s_1,w_1} \big)\Big)I^W_4\Big(F^{(n)}_{t_1,y_2} \otimes  \big(F^{(n)}_{s_2,w_2} \otimes F^{(n)}_{s_2,w_2} \otimes F^{(n)}_{s_2,w_2} \big)\Big) \Big]\\
\begin{aligned}
& =c\, \Big\langle \text{Sym}\Big(F^{(n)}_{t_1,y_1} \otimes  \big(F^{(n)}_{s_1,w_1} \otimes F^{(n)}_{s_1,w_1} \otimes F^{(n)}_{s_1,w_1} \big)\Big), \\
 &\hspace{4cm}   \text{Sym}\Big(F^{(n)}_{t_1,y_2} \otimes  \big(F^{(n)}_{s_2,w_2} \otimes F^{(n)}_{s_2,w_2} \otimes F^{(n)}_{s_2,w_2} \big)\Big) \Big\rangle_{L^2((\R_+\times \R^3)^4)}\\
&= c_{\mathbf{1}} \cq^{\mathbf{1},\mathbf{1},(n)}_{t_1,s_1,s_2}({y,w})+c_{\mathbf{2}} \cq^{\mathbf{1},\mathbf{2},(n)}_{t_1,s_1,s_2}({y,w})
\end{aligned}
\end{multline*}
for some combinatorial coefficients $c,c_{\mathbf{1}},c_{\mathbf{2}}\geq 0$, and with
\begin{align*}
 \cq^{\mathbf{1},\mathbf{1},(n)}_{t_1,s_1,s_2}(y,w)&:=\cac^{(n)}_{t_1,t_1}(y_1,y_2)   \cac^{(n)}_{s_1,s_2}(w_1,w_2)^3,\\
 \cq^{\mathbf{1},\mathbf{2},(n)}_{t_1,s_1,s_2}(y,w)&:=\cac^{(n)}_{t_1,s_2}(y_1,w_2) \cac^{(n)}_{t_1,s_1}(y_2,w_1)\cac^{(n)}_{s_1,s_2}(w_1,w_2)^2.
\end{align*}
Thus,
\begin{align*}
\cm^{\mathbf{1},(n)}_{t_1,t_2}(y_1,y_2)= c_{\mathbf{1}}\cm^{\mathbf{1,1},(n)}_{t_1,t_2}(y_1,y_2)+ c_{\mathbf{2}}\cm^{\mathbf{1,2},(n)}_{t_1,t_2}(y_1,y_2),
\end{align*}
with
\begin{multline*}
\cm^{\mathbf{1,b},(n)}_{t_1,t_2}(y_1,y_2):=\\
=\int_{0}^{t_1}ds_1\int_{0}^{t_1}ds_2\int dw_1 dw_2 \, K_{t_1-s_1,t_2-s_1}(y_1,w_1) K_{t_1-s_2,t_2-s_2}(y_2,w_2)\cq^{\mathbf{1},\mathbf{b},(n)}_{t_1,s_1,s_2}(y,w),
\end{multline*}
and accordingly
\begin{equation}\label{contr-2-iii}
\int dy_1 \, \big|\cm^{\mathbf{1},(n),\star}_{t_1,t_2}(y_1)\big| \lesssim \int dy_1 \, \big|\cm^{\mathbf{1,1},(n),\star}_{t_1,t_2}(y_1)\big| +\int dy_1 \,  \big|\cm^{\mathbf{1,2},(n),\star}_{t_1,t_2}(y_1)\big| .
\end{equation}

\

\noindent
\textbf{$\mathfrak{III}$.1.1. \underline{Study of $\cm^{\mathbf{1,1},(n)}$}.} Following the same strategy as in the proof of \eqref{ref-1-1-papi}, we deduce
\begin{align*} 
\int dy_1 \, \big|\cm^{\mathbf{1,1},(n),\star}_{t_1,t_2}(y_1)\big|
&\lesssim \big\|{\cac}^{(n),\star}_{t_1,t_1}\big\|_{L^2} \int_{0}^{t_1}ds_1\int_{0}^{t_1}ds_2\,  \big\| K_{t_1-s_1,t_2-s_1}^\star  \big\|_{L^2}\big\|  K_{t_1-s_2,t_2-s_2}^\star\big\|_{L^1} \big\|\cac^{(n),\star}_{s_1,s_2}\big\|_{L^3}^3.
\end{align*}
We can then appeal to the estimates of Lemma \ref{lem:k-star-1} and Lemma \ref{lem:cac3}, which yield for all $\eps>0$ and $0\leq \theta_1,\theta_2 \leq 1$,
\begin{align*} 
\int dy_1 \, \big|\cm^{\mathbf{1,1},(n),\star}_{t_1,t_2}(y_1)\big|
&\lesssim |t_2-t_1|^{\theta_1+\theta_2} \int_{0}^{t_1}\int_{0}^{t_1} \frac{ds_1ds_2}{|s_1-s_2|^\eps}\frac{e^{-(t_1-s_1)}}{|t_1-s_1|^{\frac34+\frac74 \theta_1}}  \frac{e^{-(t_1-s_2)}}{|t_1-s_2|^{\frac52\theta_2}} \\
&\lesssim |t_2-t_1|^{\theta_1+\theta_2} \int_{0}^{+\infty}\int_{0}^{+\infty} \frac{ds_1ds_2}{|s_1-s_2|^\eps}\frac{e^{-s_1}}{|s_1|^{\frac34+\frac74 \theta_1}}  \frac{e^{-s_2}}{|s_2|^{\frac52\theta_2}} .
\end{align*}
By choosing for instance $\theta_1=\theta_2=\frac18$, we finally deduce
\begin{equation}\label{cm-1-1-iii} 
\int dy_1 \, \big|\cm^{\mathbf{1,1},(n),\star}_{t_1,t_2}(y_1)\big| \lesssim |t_2-t_1|^{\frac14},
\end{equation}
uniformly over $t_1,t_2\geq 0$.

\

\noindent
\textbf{$\mathfrak{III}$.1.2. \underline{Study of $\cm^{\mathbf{1,2},(n)}$}.} By an analogous computation to the one yielding \eqref{ref-1-2-papi}, we establish
\begin{multline*}
\int dy_1 \, \big|\cm^{\mathbf{1,2},(n),\star}_{t_1,t_2}(y_1)\big|\lesssim \\
\lesssim \int_{0}^{t_1}ds_1\int_{0}^{t_1}ds_2\, \big\|\cac^{(n),\star}_{t_1,s_2}\big\|_{L^2} \big\| K_{t_1-s_1,t_2-s_1}^\star\big\|_{L^2}\big\|\cac^{(n),\star}_{t_1,s_1}\big\|_{L^2} \big\| K_{t_1-s_2,t_2-s_2}^\star\big\|_{L^2} \big\| \cac^{(n),\star}_{s_1,s_2}\big\|_{L^2}^2.
\end{multline*}
Using the estimates contained in Lemma \ref{lem:k-star-1} and Lemma \ref{lem:cac3}, we deduce that for every $\theta\in [0,1]$,
\begin{align*}
\int dy_1 \, \big|\cm^{\mathbf{1,2},(n),\star}_{t_1,t_2}(y_1)\big|&\lesssim |t_2-t_1|^{2\theta}\bigg(\int_{0}^{t_1}ds\, \frac{e^{-(t_1-s)}}{|t_1-s|^{\frac34+\frac74 \theta}} \bigg)^2,
\end{align*}
and we can choose $\theta=\frac18$ to conclude that
\begin{equation}\label{cm-1-2-iii}
\int dy_1 \, \big|\cm^{\mathbf{1,2},(n),\star}_{t_1,t_2}(y_1)\big| \lesssim |t_2-t_1|^{\frac14},
\end{equation}
uniformly over $t_1,t_2\geq 0$.

\

\noindent
\textbf{$\mathfrak{III}$.2. Study of $\mathfrak{III}^{\mathbf{2},(n)}$.} The quantity $\cm^{\mathbf{2},(n)}_{t_1,t_2}$ can be computed as follows:
\begin{align*} 
&\cm^{\mathbf{2},(n)}_{t_1,t_2}(y_1,y_2)=\mathbb{E}\Big[\mathfrak{III}^{\mathbf{2},(n)}_{t_1,t_2}(y_1)\, \mathfrak{III}^{\mathbf{2},(n)}_{t_1,t_2}(y_2) \Big]=\\
&=9\int_{0}^{t_1}ds_1\int_{0}^{t_1}ds_2 \int dw_1 dw_2 \, K_{t_1-s_1,t_2-s_1}(y_1,w_1)K_{t_1-s_2,t_2-s_2}(y_2,w_2)\\
&\hspace{4cm}\times  \cac^{(n)}_{t_1,s_1}(y_1,w_1) \cac^{(n)}_{t_1,s_2}(y_2,w_2)\mathbb{E}\Big[ I^W_2\big( F^{(n)}_{s_1,w_1} \otimes F^{(n)}_{s_1,w_1} \big)I^W_2\big( F^{(n)}_{s_2,w_2} \otimes F^{(n)}_{s_2,w_2} \big)\Big] \\
&=c\int_{0}^{t_1}ds_1\int_{0}^{t_1}ds_2 \int dw_1 dw_2 \, K_{t_1-s_1,t_2-s_1}(y_1,w_1)K_{t_1-s_2,t_2-s_2}(y_2,w_2) \\
&\hspace{6cm}\times  \cac^{(n)}_{t_1,s_1}(y_1,w_1) \cac^{(n)}_{t_1,s_2}(y_2,w_2) \cac^{(n)}_{s_1,s_2}(w_1,w_2)^2.
\end{align*}
Applying the estimates in the same manner as for \eqref{ref-2-papi}, we see that
\begin{multline*} 
\int dy_1\, \cm^{\mathbf{2},(n),\star}_{t_1,t_2}(y_1)\lesssim \\
\lesssim \int_{0}^{t_1}ds_1\int_{0}^{t_1}ds_2\,  \big\| K^\star_{t_1-s_1,t_2-s_1}\big\|_{L^2}\big\| \cac^{(n),\star}_{t_1,s_1}\big\|_{L^2}\big\|(\cac^{(n),\star}_{s_1,s_2})^2 \big\|_{L^1}\big\|K^\star_{t_1-s_2,t_2-s_2}\big\|_{L^2}\big\| \cac^{(n),\star}_{t_1,s_2}\big\|_{L^2}.
\end{multline*}
Thanks to Lemma \ref{lem:k-star-1} and Lemma \ref{lem:cac3}, we obtain
\begin{equation}\label{cm-2-iii}
\int dy_1\, \cm^{\mathbf{2},(n),\star}_{t_1,t_2}(y_1)\lesssim |t_2-t_1|^{\frac14}\bigg(\int_{0}^{t_1}ds\frac{e^{-(t_1-s)}}{|t_1-s|^{\frac34+\frac{7}{32}}} \bigg)^2\lesssim |t_2-t_1|^{\frac14},
\end{equation}
uniformly over $t_1,t_2\geq 0$.

\

Finally, by inserting \eqref{cm-1-1-iii}-\eqref{cm-1-2-iii}-\eqref{cm-2-iii} into \eqref{contr-1-iii} and \eqref{contr-2-iii}, we obtain \eqref{boun-frak-iii}.

\

\


\section{About the fourth order diagram 2} \label{sec:diag-4th-order-2}

In this section, we focus on the second fourth-order diagram involved in the procedure, namely
$$\<Psi2IPsi2>^{(n)}_{s,t}:=\<Psi2>^{(n)}_{s,t} \pe \<IPsi2>^{(n)}_{s,t}- \frakc^{\mathbf{2},(n)}_{s,t}$$
where
$$\frakc^{\mathbf{2},(n)}_{s,t}(x) := \mathbb{E}\Big[ \<Psi2>^{(n)}_{s,t}(x) \<IPsi2>^{(n)}_{s,t}(x)\Big].$$

\smallskip

A first evaluation of the convergence and regularity of $\<Psi2IPsi2>^{(n)}$ can be found in \cite[Proposition~9.1]{DFT}. Combined with the stationarity property in Lemma \ref{lem:stati-z}, this regularity result can be summarized as follows.

\begin{proposition}\label{Prop-p51}
Fix $0\leq T_1<T_2$ and for all $t\in [T_1,T_2]$, set
\begin{equation*}
\widetilde{ \<Psi2IPsi2>}^{(n)}_{T_1,t}:=\int_{T_1}^t \<Psi2IPsi2>^{(n)}_{T_1,s} \, ds.
\end{equation*}
Then for all $0<\eps,\eta<\frac12$, there exists $\ka>0$ such that for every $p\geq 1$,
\begin{equation}\label{L3-2}
\sup_{\ell\geq 0}    \mathbb{E} \Big[ \Big\|\widetilde{ \<Psi2IPsi2>}^{(n+1)}_{\ell,.}   - \widetilde{ \<Psi2IPsi2>}^{(n)}_{\ell,.}\Big\|_{{\ov \cac}^{1-\eps}([\ell,\ell+2]; \cb^{-\eta}_{\infty})}^{2p} \Big]\lesssim 2^{-\ka n p }. 
\end{equation}
As a particular consequence, for all $0\leq T_1<T_2$, the sequence $(\widetilde{ \<Psi2IPsi2>}^{(n)}_{T_1,.})$ converges almost surely to an element $\widetilde{\<Psi2IPsi2>}_{T_1,.}$  in the space ${\cac}^{1-\eps}([T_1,T_2]; \cb^{-\eta}_{\infty})$, for all $0<\eps,\eta<\frac12$.
\end{proposition}

\medskip

Just as in the previous section,  we now plan to combine \eqref{L3-2} with an estimate of a different type, leading to a well-defined function $\<Psi2IPsi2>$ in a suitable space. The result of this (forthcoming) procedure can be stated as follows.

 \begin{theorem}\label{coro-9.1}
For all $\eta>0$, there exist $\eps>0$ and $\ka>0$ such that for every $p\geq 1$
\begin{equation}\label{boun-9.1}
\sup_{\ell\geq 0}   \mathbb{E} \Big[ \Big\|  \<Psi2IPsi2>^{(n+1)}_{\ell,.} -  \<Psi2IPsi2>^{(n)}_{\ell,.}\Big\|_{\cac^{\eps}([\ell,\ell+1]; \cb^{-\eta}_{\infty})}^{2p} \Big]\lesssim 2^{-\ka n p}. 
\end{equation}
Consequently, for all $0\leq T_1<T_2$, the sequence $(\<Psi2IPsi2>^{(n)}_{T_1,.})$ converges almost surely to  an element  $\<Psi2IPsi2>_{T_1,.}$ in  $\cac^{\varepsilon}\big([T_1,T_2];\cb^{-\eta}_{\infty}(\R^3)\big)$.

\end{theorem}

\

We introduce the intermediate \enquote{full-product} process
$$\<Psi2IPsi2nr>^{(n)}_{s,t}:=\<Psi2>^{(n)}_{s,t} \<IPsi2>^{(n)}_{s,t}- \frakc^{\mathbf{2},(n)}_{s,t},$$
noting that
\begin{equation}\label{difference}
\<Psi2IPsi2nr>^{(n)}_{s,t}-\<Psi2IPsi2>^{(n)}_{s,t}= \<Psi2>^{(n)}_{s,t} \pl \<IPsi2>^{(n)}_{s,t}+\<Psi2>^{(n)}_{s,t} \pg \<IPsi2>^{(n)}_{s,t}.
\end{equation}

Our main technical result in this section, leading to Theorem \ref{coro-9.1}, reads as follows:

\begin{proposition}\label{pro-anex9}
There exists $\ka>0$ such that for every $p\geq 1$,
\begin{equation}\label{L30-2}
\sup_{\ell\geq 0}    \mathbb{E} \Big[ \Big\|\<Psi2IPsi2nr>^{(n+1)}_{\ell,.}   -  \<Psi2IPsi2nr>^{(n)}_{\ell,.}\Big\|_{{\cac}^{\frac1{10}}([\ell,\ell+2]; \cb^{-2}_\infty)}^{2p} \Big]\lesssim 2^{-\ka n p }. 
\end{equation}
\end{proposition}

\

\begin{proof}[Proof of Theorem \ref{coro-9.1}] 
Starting from \eqref{difference} and using the regularity properties of $\<Psi2>^{(n)}$ contained in \cite[Proposition 6.1]{DFT}, we can transfer the estimate \eqref{L30-2} for $\<Psi2IPsi2nr>^{(n)}$ to the resonant part $\<Psi2IPsi2>^{(n)}$, that is, one has
\begin{equation}\label{L3-bis-2}
\sup_{\ell\geq 0}    \mathbb{E} \Big[ \Big\|\<Psi2IPsi2>^{(n+1)}_{\ell,.}   -  \<Psi2IPsi2>^{(n)}_{\ell,.}\Big\|_{{\cac}^{\frac1{10}}([\ell,\ell+2]; \cb^{-2}_{\infty})}^{2p} \Big]\lesssim 2^{-\ka n p }. 
\end{equation}
The bound \eqref{boun-9.1} is then obtained as in the proof of Theorem \ref{coro-8.1}, by replacing~\eqref{L3} with \eqref{L3-2} and \eqref{L3-bis} with \eqref{L3-bis-2}.
\end{proof}

\medskip

We are thus left with the proof of Proposition \ref{pro-anex9}. For the sake of conciseness, we will only focus on the uniform bound
\begin{equation*}
 \sup_{n\geq 1} \sup_{\ell\geq 0} \, \mathbb{E} \Big[ \big\| \<Psi2IPsi2nr>^{(n)}_{\ell,.}  \big\|_{{\cac}^{\frac1{10}}([\ell,\ell+2];\cb_{\infty}^{-2})}^{2p} \Big] <\infty.
\end{equation*}
Using the stationarity property contained in Lemma \ref{lem:stati-z}, it suffices to establish that
\begin{equation*}
 \sup_{n\geq 1}  \, \mathbb{E} \Big[ \big\| \<Psi2IPsi2nr>^{(n)}_{0,.}  \big\|_{{\cac}^{\frac1{10}}([0,2];\cb_{\infty}^{-2})}^{2p} \Big] <\infty.
\end{equation*}

\

Fix  $\frac34<\al<1$. For the same reasons as in Section \ref{sec:diag-4th-order-1} (see \eqref{sobo-beso}-\eqref{normlip-40}-\eqref{hypercon}), one has, for every $p\geq 1$ large enough,
\begin{equation*}
\mathbb{E}\Big[\big\|\<Psi2IPsi2nr>^{(n)}_{0,.} \big\|_{\cac^{\frac{1}{10}}([0,2];\cb_{\infty}^{-2})}^{2p}\Big]
 \lesssim \int_{[0,2]^2} \frac{dt_1 dt_2 }{|t_2-t_1|^{\frac{ p}{5}+2}}\int dx \, \mathbb{E}\bigg[ \Big| H^{-\al}\big(\<Psi2IPsi2nr>^{(n)}_{0,t_2}-\<Psi2IPsi2nr>^{(n)}_{0,t_1}\big)(x)\Big|^{2}\bigg]^p.
\end{equation*}
Thus, it suffices to show that
\begin{equation}\label{hyperc-40-1}
\sup_{n\geq 1}\int dx \, \mathbb{E}\bigg[ \Big| H^{-\al}\big(\<Psi2IPsi2nr>^{(n)}_{0,t_2}-\<Psi2IPsi2nr>^{(n)}_{0,t_1}\big)(x)\Big|^{2}\bigg]^p \lesssim |t_2-t_1|^{\ka p}
\end{equation}
uniformly over $0\leq t_1,t_2\leq 2$, for some $\ka>\frac{1}{5}$.

\

\

\subsection{Proof of \eqref{hyperc-40-1}}
Let us write
\begin{align*}
\<Psi2IPsi2nr>^{(n)}_{0,t}(y)&=\<Psi2>^{(n)}_{0,t}(y) \<IPsi2>^{(n)}_{0,t}(y)- \frakc^{\mathbf{2},(n)}_{0,t}(y)\\
& =\int_{0}^t ds \int dw \, K_{t-s}(y,w)\<Psi2>^{(n)}_{0,t}(y)\<Psi2>^{(n)}_{0,s}(w)- \frakc^{\mathbf{2},(n)}_{0,t}(y)\\
&=\int_{0}^t ds \int dw \, K_{t-s}(y,w)I^W_2\big(F^{(n)}_{t,y}\otimes F^{(n)}_{t,y}\big)I^W_2\big(F^{(n)}_{s,w} \otimes F^{(n)}_{s,w}  \big)- \frakc^{\mathbf{2},(n)}_{0,t}(y),
\end{align*}
which, by applying the product rule stated in \cite[Lemma 4.1]{DFT}, yields the decomposition
\begin{equation}\label{decompo-scretundeux}
\<Psi2IPsi2nr>^{(n)}_{0,t} (x) =\scret^{\mathbf{1},(n)}_t (x) +\scret^{\mathbf{2},(n)}_t (x) 
\end{equation}
with
\begin{align*}
\scret^{\mathbf{1},(n)}_t (y) :=\int_{0}^t ds \int dw \, K_{t-s}(y,w)I^W_4\Big(F^{(n)}_{t,y}\otimes F^{(n)}_{t,y}\otimes F^{(n)}_{s,w} \otimes F^{(n)}_{s,w} \Big) 
\end{align*}
and
\begin{align*}
\scret^{\mathbf{2},(n)}_t (y) :=3\int_{0}^t ds \int dw \, K_{t-s}(y,w) \,\cac^{(n)}_{t,s}(y,w) \,  I^W_2\big( F^{(n)}_{t,y} \otimes F^{(n)}_{s,w} \big).
\end{align*}
With this notation, we can write 
\begin{align}
&\<Psi2IPsi2nr>^{(n)}_{0,t_2}-\<Psi2IPsi2nr>^{(n)}_{0,t_1}  =\big(\scret^{\mathbf{1},(n)}_{t_2}-\scret^{\mathbf{1},(n)}_{t_1}\big) +\big(\scret^{\mathbf{2},(n)}_{t_2}- \scret^{\mathbf{2},(n)}_{t_1}\big) =:\scret^{\mathbf{1},(n)}_{t_1,t_2}+\scret^{\mathbf{2},(n)}_{t_1,t_2}.\label{decompo-scretun-scretdeux}
\end{align}

\

\

\noindent
\textbf{1. Study of $\scret^{\mathbf{1},(n)}_{t_1,t_2}$.}

\smallskip

One can decompose this quantity as
\begin{multline*}
\scret^{\mathbf{1},(n)}_{t_1,t_2}(y)=\int_{t_1}^{t_2} ds \int dw \, K_{t_2-s}(y,w)I^W_4\Big(F^{(n)}_{t_2,y}\otimes F^{(n)}_{t_2,y}\otimes F^{(n)}_{s,w} \otimes F^{(n)}_{s,w} \Big)\\
\hspace{-1cm}+\int_{0}^{t_1} ds \int dw \, K_{t_1-s,t_2-s}(y,w)I^W_4\Big(F^{(n)}_{t_2,y}\otimes F^{(n)}_{t_2,y}\otimes F^{(n)}_{s,w} \otimes F^{(n)}_{s,w} \Big)\\
+\int_{0}^{t_1} ds \int dw \, K_{t_1-s}(y,w)I^W_4\Big(\big(F^{(n)}_{t_2,y}-F^{(n)}_{t_1,y}\big)\otimes F^{(n)}_{t_2,y}\otimes F^{(n)}_{s,w} \otimes F^{(n)}_{s,w} \Big)\\
+\int_{0}^{t_1} ds \int dw \, K_{t_1-s}(y,w)I^W_4\Big(F^{(n)}_{t_1,y}\otimes \big(F^{(n)}_{t_2,y}-F^{(n)}_{t_1,y}\big)\otimes F^{(n)}_{s,w} \otimes F^{(n)}_{s,w} \Big)\\
=:\scret^{\mathbf{1,1},(n)}_{t_1,t_2}(y)+\scret^{\mathbf{1,2},(n)}_{t_1,t_2}(y)+\scret^{\mathbf{1,3},(n)}_{t_1,t_2}(y)+\scret^{\mathbf{1,4},(n)}_{t_1,t_2}(y).
\end{multline*}

\

For $\mathbf{a}=1,\ldots,4$, let us set
$$\cm^{\mathbf{1,a},(n)}_{t_1,t_2}(y_1,y_2):=\mathbb{E}\Big[\scret^{\mathbf{1,a},(n)}_{t_1,t_2}(y_1)\, \scret^{\mathbf{1,a},(n)}_{t_1,t_2}(y_2) \Big].$$
Then, thanks to Lemma \ref{lem:techn} and since $\al>\frac34$, we deduce that 
\begin{equation}\label{contr-1-bis}
\int dx \, \mathbb{E}\bigg[ \Big| H^{-\al}\big(\scret^{\mathbf{1},\mathbf{a},(n)}_{t_1,t_2}\big)(x)\Big|^{2}\bigg]^p\lesssim \sum_{\mathbf{a}=1,\ldots,4}\bigg(\int dy_1 \,  \big|\cm^{\mathbf{1,a},(n),\star}_{t_1,t_2}(y_1)\big| \bigg)^p,
\end{equation}
for every $p$ large enough, where we have used the notation $\star$ introduced in \eqref{star-not}.

\

\

\subsubsection{Study of $\cm^{\mathbf{1,1},(n)}$.} One has here
\begin{multline*}
\cm^{\mathbf{1,1},(n)}_{t_1,t_2}(y_1,y_2)=\mathbb{E}\Big[\scret^{\mathbf{1,1},(n)}_{t_1,t_2}(y_1)\, \scret^{\mathbf{1,1},(n)}_{t_1,t_2}(y_2) \Big]=\\
\begin{aligned}
&=\int_{t_1}^{t_2}ds_1\int_{t_1}^{t_2}ds_2 \int dw_1 dw_2 \, K_{t_2-s_1}(y_1,w_1)K_{t_2-s_2}(y_2,w_2)\\
&\hspace{1cm}\times \mathbb{E}\Big[I^W_4\Big(F^{(n)}_{t_2,y_1}\otimes F^{(n)}_{t_2,y_1}\otimes F^{(n)}_{s_1,w_1} \otimes F^{(n)}_{s_1,w_1} \Big)I^W_4\Big(F^{(n)}_{t_2,y_2}\otimes F^{(n)}_{t_2,y_2}\otimes F^{(n)}_{s_2,w_2} \otimes F^{(n)}_{s_2,w_2} \Big) \Big].
\end{aligned}
\end{multline*}
Recall that the latter expectation can be expanded as in \eqref{expand-expec}, which gives
\begin{align*}
\cm^{\mathbf{1,1},(n)}_{t_1,t_2}(y_1,y_2)= c_{\mathbf{1}}\cm^{\mathbf{1,1,1},(n)}_{t_1,t_2}(y_1,y_2)+ c_{\mathbf{2}}\cm^{\mathbf{1,1,2},(n)}_{t_1,t_2}(y_1,y_2)+ c_{\mathbf{3}}\cm^{\mathbf{1,1,3},(n)}_{t_1,t_2}(y_1,y_2),
\end{align*}
where (with the notation in \eqref{expand-expec})
\begin{align*}
\cm^{\mathbf{1,1,b},(n)}_{t_1,t_2}(y_1,y_2):=\int_{t_1}^{t_2}ds_1\int_{t_1}^{t_2}ds_2\int dw_1 dw_2 \, K_{t_2-s_1}(y_1,w_1) K_{t_2-s_2}(y_2,w_2)\cq^{\mathbf{1,b},(n)}_{t_2,s_1,s_2}(y,w).
\end{align*}
As a result,
\begin{multline}\label{contr-2-ii-bis-1-1}
\int dy_1 \, \big|\cm^{\mathbf{1,1},(n),\star}_{t_1,t_2}(y_1)\big|\lesssim\\
\lesssim \int dy_1 \, \big|\cm^{\mathbf{1,1,1},(n),\star}_{t_1,t_2}(y_1)\big| +\int dy_1 \,  \big|\cm^{\mathbf{1,1,2},(n),\star}_{t_1,t_2}(y_1)\big| +\int dy_1 \,  \big|\cm^{\mathbf{1,1,3},(n),\star}_{t_1,t_2}(y_1)\big|.
\end{multline}

\

\

\noindent
\textit{(i) \underline{Study of $\cm^{\mathbf{1,1,1},(n)}$.}} One has for all $y_1,y_2$,
\begin{multline*} 
\big|\cm^{\mathbf{1,1,1},(n)}_{t_1,t_2}(y_1+y_2,y_2)\big| \lesssim\\
\begin{aligned} 
&\lesssim \cac^{(n)}_{t_2,t_2}(y_1+y_2,y_2)^2\int_{t_1}^{t_2}ds_1\int_{t_1}^{t_2}ds_2 \int dw_1\, \big|K_{t_2-s_1}(y_1+y_2,w_1) \big| \\
&\hspace{4cm}\times \int dw_2 \, \big| K_{t_2-s_2}(y_2,w_2) \big| \cac^{(n)}_{s_1,s_2}(w_1,w_2)^2  \\
&\lesssim {\cac}^{(n),\star}_{t_2,t_2}(y_1)^2\int_{t_1}^{t_2}ds_1\int_{t_1}^{t_2}ds_2 \int dw_1\, \big| K_{t_2-s_1}(y_1+y_2,w_1+y_2)\big|    \\
&\hspace{4cm}\times \int dw_2 \, \big| K_{t_2-s_2}(y_2,w_2+y_2)\big|   \cac^{(n)}_{s_1,s_2}(w_1+y_2,w_2+y_2)^2\\
&\lesssim {\cac}^{(n),\star}_{t_2,t_2}(y_1)^2\int_{t_1}^{t_2}ds_1\int_{t_1}^{t_2}ds_2 \int dw_1\, K_{t_2-s_1}^\star(y_1-w_1) \int dw_2 \, K_{t_2-s_2}^\star(w_2)  \cac^{(n),\star}_{s_1,s_2}(w_1-w_2)^2\\
&\lesssim {\cac}^{(n),\star}_{t_2,t_2}(y_1)^2\int_{t_1}^{t_2}ds_1\int_{t_1}^{t_2}ds_2\,  \Big[ K_{t_2-s_1}^\star  \ast \big[ K_{t_2-s_2}^\star \ast  (\cac^{(n),\star}_{s_1,s_2})^2\big]\Big](y_1).
\end{aligned}
\end{multline*}
Thus,
\begin{multline}\label{ref-1-1-pap}
\int dy_1 \sup_{y_2} \big|\cm^{\mathbf{1,1,1},(n)}_{t_1,t_2}(y_1+y_2,y_2)\big|\lesssim\\
\begin{aligned}
&\lesssim \big\|(\cac^{(n),\star}_{t_2,t_2})^2\big\|_{L^1} \int_{t_1}^{t_2} ds_1 \int_{t_1}^{t_2} ds_2 \, \Big\|K^\star_{t_2-s_1}\ast\big[K^\star_{t_2-s_2}\ast (\cac^{(n),\star}_{s_1,s_2})^2\big]\Big\|_{L^\infty}\\
&\lesssim \big\|\cac^{(n),\star}_{t_2,t_2}\big\|_{L^2}^2 \int_{t_1}^{t_2}ds_1 \int_{t_1}^{t_2} ds_2 \, \big\|K^\star_{t_2-s_1}\big\|_{L^1} \big\| K^\star_{t_2-s_2}\big\|_{L^2} \big\|\cac^{(n),\star}_{s_1,s_2}\big\|_{L^4}^2.
\end{aligned}
\end{multline}
Then, thanks to Lemma \ref{lem:k-star-0} and Lemma \ref{lem:cac3}, we obtain
\begin{align*}
 \int dy_1 \sup_{y_2} \big|\cm^{\mathbf{1,1,1},(n)}_{t_1,t_2}(y_1+y_2,y_2)\big|
&\lesssim \int_{t_1}^{t_2} ds_1 \int_{t_1}^{t_2} ds_2 \, e^{-(t_2-s_1)}\frac{e^{- (t_2-s_2)}}{|t_2-s_2|^{\frac34}}   \frac{1}{|s_2-s_1|^{\frac14}}\\
&\lesssim \int_{t_1}^{t_2} ds_1 \int_{t_1}^{t_2} ds_2 \, \frac{1}{|t_2-s_2|^{\frac34}}   \frac{1}{|s_2-s_1|^{\frac14}} 
\end{align*}
and as a result,
\begin{align}\label{contr-2-ii-bis-1-1-1}
&\int dy_1 \sup_{y_2} \big|\cm^{\mathbf{1,1,1},(n)}_{t_1,t_2}(y_1+y_2,y_2)\big|\lesssim   |t_2-t_1|,
\end{align}
uniformly over $0\leq t_1,t_2\leq 2$.

\

\

\noindent
\textit{(ii) \underline{Study of $\cm^{\mathbf{1,1,2},(n)}$.}} One has for all $y_1,y_2$,
\begin{multline*} 
\big|\cm^{\mathbf{1,1,2},(n)}_{t_1,t_2}(y_1+y_2,y_2)\big| \lesssim\\
\begin{aligned} 
&\lesssim \cac^{(n)}_{t_2,t_2}(y_1+y_2,y_2)\int_{t_1}^{t_2} ds_1\int_{t_1}^{t_2} ds_2 \int dw_1\, \big|K_{t_2-s_1}(y_1+y_2,w_1+y_2) \big|\cac^{(n)}_{t_2,s_1}(y_2,w_1+y_2)\\
&\hspace{3cm}\times \int dw_2 \, \big| K_{t_2-s_2}(y_2,w_2+y_2) \big| \cac^{(n)}_{t_2,s_2}(y_1+y_2,w_2+y_2) \cac^{(n)}_{s_1,s_2}(w_1+y_2,w_2+y_2)  \\
&\lesssim {\cac}^{(n),\star}_{t_2,t_2}(y_1)\int_{t_1}^{t_2} ds_1\int_{t_1}^{t_2} ds_2 \int dw_1\, K_{t_2-s_1}^\star(y_1-w_1)\cac^{(n),\star}_{t_2,s_1}(w_1) \\
&\hspace{3cm}\times \int dw_2 \, K_{t_2-s_2}^\star(w_2) \cac^{(n),\star}_{t_2,s_2}(y_1-w_2) \cac^{(n),\star}_{s_1,s_2}(w_1-w_2)\\
&\lesssim {\cac}^{(n),\star}_{t_2,t_2}(y_1)\int_{t_1}^{t_2} ds_1\, \big\|\cac^{(n),\star}_{t_2,s_1}\big\|_{L^\infty}\int_{t_1}^{t_2} ds_2\, \big\|\cac^{(n),\star}_{t_2,s_2}\big\|_{L^\infty}\Big[ K_{t_2-s_1}^\star \ast \big[K_{t_2-s_2}^\star \ast  \cac^{(n),\star}_{s_1,s_2}\big]\Big](y_1).
\end{aligned}    
 \end{multline*}                
Thus, 
\begin{multline} \label{ref-1-2-pap}
 \int dy_1 \sup_{y_2} \big|\cm^{\mathbf{1,1,2},(n)}_{t_1,t_2}(y_1+y_2,y_2)\big| \lesssim\\
\begin{aligned}
&\lesssim \big\|\cac^{(n),\star}_{t_2,t_2}\big\|_{L^2} \int_{t_1}^{t_2} ds_1\, \big\|\cac^{(n),\star}_{t_2,s_1}\big\|_{L^\infty}\int_{t_1}^{t_2} ds_2\, \big\|\cac^{(n),\star}_{t_2,s_2}\big\|_{L^\infty} \, \Big\|K^\star_{t_2-s_1}\ast\big[K^\star_{t_2-s_2}\ast \cac^{(n),\star}_{s_1,s_2}\big]\Big\|_{L^2} \\
&\lesssim \big\|\cac^{(n),\star}_{t_2,t_2}\big\|_{L^2} \int_{t_1}^{t_2} ds_1\, \big\|\cac^{(n),\star}_{t_2,s_1}\big\|_{L^\infty}\int_{t_1}^{t_2} ds_2\, \big\|\cac^{(n),\star}_{t_2,s_2}\big\|_{L^\infty} \, \big\|K^\star_{t_2-s_1}\big\|_{L^1}\big\|K^\star_{t_2-s_2}\big\|_{L^1} \big\|\cac^{(n),\star}_{s_1,s_2}\big\|_{L^2},
\end{aligned}
 \end{multline}
Now using estimates in Lemma \ref{lem:k-star-0} and Lemma \ref{lem:cac3}, we can conclude that
\begin{align}\label{contr-2-ii-bis-1-1-2}
&\int dy_1 \sup_{y_2} \big|\cm^{\mathbf{1,1,2},(n)}_{t_1,t_2}(y_1+y_2,y_2)\big|\lesssim \int_{t_1}^{t_2} \frac{ds_1}{|t_2-s_1|^{\frac12}} \int_{t_1}^{t_2} \frac{ds_2}{|t_2-s_2|^{\frac12}}\lesssim   |t_2-t_1|,
\end{align}
uniformly over $0\leq t_1,t_2\leq 2$.

\

\

\noindent
\textit{(iii) \underline{Study of $\cm^{\mathbf{1,1,3},(n)}$.}} One has for all $y_1,y_2$,
\begin{multline*} 
\big|\cm^{\mathbf{1,1,3},(n)}_{t_1,t_2}(y_1+y_2,y_2)\big|\lesssim \\
\begin{aligned} 
&\lesssim \int_{t_1}^{t_2} ds_1\int_{t_1}^{t_2} ds_2 \int dw_1\, \big|K_{t_2-s_1}(y_1+y_2,w_1+y_2) \big| \cac^{(n)}_{t_2,s_1}(y_2,w_1+y_2)^2 \\
&\hspace{3cm}\times \int dw_2 \, \big| K_{t_2-s_2}(y_2,w_2+y_2) \big| \cac^{(n)}_{t_2,s_2}(y_1+y_2,w_2+y_2)^2  \\
&\lesssim \int_{t_1}^{t_2} ds_1\int_{t_1}^{t_2} ds_2 \int dw_1\, K_{t_2-s_1}^\star(y_1-w_1)\cac^{(n),\star}_{t_2,s_1}(w_1)^2 \int dw_2 \, K_{t_2-s_2}^\star(w_2)  \cac^{(n),\star}_{t_2,s_2}(y_1-w_2)^2\\
&\lesssim \int_{t_1}^{t_2} ds_1\int_{t_1}^{t_2} ds_2\,  \Big[ K_{t_2-s_1}^\star  \ast   (\cac^{(n),\star}_{t_2,s_1})^2\Big](y_1)\Big[ K_{t_2-s_2}^\star  \ast   (\cac^{(n),\star}_{t_2,s_2})^2\Big](y_1).
\end{aligned}
\end{multline*} 
Thus
\begin{align}
 \int dy_1 \sup_{y_2} \big|\cm^{\mathbf{1,1,3},(n)}_{t_1,t_2}(y_1+y_2,y_2)\big|
&\lesssim  \bigg(\int_{t_1}^{t_2} ds \, \big\| K_{t_2-s}^\star \big\|_{L^1}   \big\|\cac^{(n),\star}_{t_2,s}\big\|_{L^4}^2\bigg)^2.\label{ref-1-3-pap}
\end{align}
Using again Lemma \ref{lem:k-star-0} and Lemma \ref{lem:cac3}, we obtain that
\begin{align}
 \int dy_1 \sup_{y_2} \big|\cm^{\mathbf{1,1,3},(n)}_{t_1,t_2}(y_1+y_2,y_2)\big|
&\lesssim  \bigg(\int_{t_1}^{t_2} \frac{ds}{|t_2-s|^{\frac14}}  \bigg)^2 \lesssim |t_2-t_1|^{\frac32},\label{contr-2-ii-bis-1-1-3}
\end{align}
uniformly over $0\leq t_1,t_2\leq 2$.

\

Finally, by inserting estimates \eqref{contr-2-ii-bis-1-1-1} to \eqref{contr-2-ii-bis-1-1-3} into \eqref{contr-2-ii-bis-1-1}, we can conclude that
\begin{align}\label{cm11}
&\int dy_1 \, \big|\cm^{\mathbf{1,1},(n),\star}_{t_1,t_2}(y_1)\big|\lesssim |t_2-t_1|,
\end{align}
uniformly over $0\leq t_1,t_2\leq 2$.

\

\

\subsubsection{Study of $\cm^{\mathbf{1,2},(n)}$.} One has here
\begin{multline*}
\cm^{\mathbf{1,2},(n)}_{t_1,t_2}(y_1,y_2)=\mathbb{E}\Big[\scret^{\mathbf{1,2},(n)}_{t_1,t_2}(y_1)\, \scret^{\mathbf{1,2},(n)}_{t_1,t_2}(y_2) \Big]\\
\hspace{-2cm}=\int_{0}^{t_1}ds_1\int_{0}^{t_1}ds_2 \int dw_1 dw_2 \, K_{t_1-s_1,t_2-s_1}(y_1,w_1)K_{t_1-s_2,t_2-s_2}(y_2,w_2)\\
\hspace{2cm}\times \mathbb{E}\Big[I^W_4\Big(F^{(n)}_{t_2,y_1}\otimes F^{(n)}_{t_2,y_1}\otimes F^{(n)}_{s_1,w_1} \otimes F^{(n)}_{s_1,w_1} \Big)I^W_4\Big(F^{(n)}_{t_2,y_2}\otimes F^{(n)}_{t_2,y_2}\otimes F^{(n)}_{s_2,w_2} \otimes F^{(n)}_{s_2,w_2} \Big) \Big].
\end{multline*}

The latter expectation can be expanded as
\begin{multline}\label{expand-expec}
\mathbb{E}\Big[I^W_4\Big(F^{(n)}_{t_2,y_1}\otimes F^{(n)}_{t_2,y_1}\otimes F^{(n)}_{s_1,w_1} \otimes F^{(n)}_{s_1,w_1} \Big)I^W_4\Big(F^{(n)}_{t_2,y_2}\otimes F^{(n)}_{t_2,y_2}\otimes F^{(n)}_{s_2,w_2} \otimes F^{(n)}_{s_2,w_2} \Big) \Big]\\
\begin{aligned}
& =c\, \Big\langle \text{Sym}\Big(F^{(n)}_{t_2,y_1}\otimes F^{(n)}_{t_2,y_1}\otimes F^{(n)}_{s_1,w_1} \otimes F^{(n)}_{s_1,w_1} \Big),  \\
 &\hspace{3.5cm}   \text{Sym}\Big(F^{(n)}_{t_2,y_2}\otimes F^{(n)}_{t_2,y_2}\otimes F^{(n)}_{s_2,w_2} \otimes F^{(n)}_{s_2,w_2} \Big) \Big\rangle_{L^2((\R_+\times \R^3)^4)} \\
&= c_{\mathbf{1}} \cq^{\mathbf{1,1},(n)}_{t_2,s_1,s_2}({y,w})+c_{\mathbf{2}} \cq^{\mathbf{1,2},(n)}_{t_2,s_1,s_2}({y,w})+c_{\mathbf{3}} \cq^{\mathbf{1,3},(n)}_{t_2,s_1,s_2}({y,w}),
\end{aligned}
\end{multline}
for some combinatorial coefficients $c,c_{\mathbf{1}},c_{\mathbf{2}},c_{\mathbf{3}}\geq 0$, and with 
\begin{align*}
 \cq^{\mathbf{1,1},(n)}_{t_2,s_1,s_2}(y,w)&:=\cac^{(n)}_{t_2,t_2}(y_1,y_2)^2   \cac^{(n)}_{s_1,s_2}(w_1,w_2)^2,\\[2pt]
 \cq^{\mathbf{1,2},(n)}_{t_2,s_1,s_2}(y,w)&:=\cac^{(n)}_{t_2,t_2}(y_1,y_2)\cac^{(n)}_{t_2,s_2}(y_1,w_2) \cac^{(n)}_{t_2,s_1}(y_2,w_1)\cac^{(n)}_{s_1,s_2}(w_1,w_2)\\[2pt]
\cq^{\mathbf{1,3},(n)}_{t_2,s_1,s_2}(y,w)&:=\cac^{(n)}_{t_2,s_2}(y_1,w_2)^2 \cac^{(n)}_{t_2,s_1}(y_2,w_1)^2.
\end{align*}

\

As a result, we obtain the decomposition
\begin{align*}
\cm^{\mathbf{1,2},(n)}_{t_1,t_2}(y_1,y_2)= c_{\mathbf{1}}\cm^{\mathbf{1,2,1},(n)}_{t_1,t_2}(y_1,y_2)+ c_{\mathbf{2}}\cm^{\mathbf{1,2,2},(n)}_{t_1,t_2}(y_1,y_2)+ c_{\mathbf{3}}\cm^{\mathbf{1,2,3},(n)}_{t_1,t_2}(y_1,y_2),
\end{align*}
with
\begin{multline*}
\cm^{\mathbf{1,2,b},(n)}_{t_1,t_2}(y_1,y_2):=\int_{0}^{t_1}ds_1\int_{0}^{t_1}ds_2\int dw_1 dw_2 \,\\
 \times K_{t_1-s_1,t_2-s_1}(y_1,w_1)K_{t_1-s_2,t_2-s_2}(y_2,w_2)\cq^{\mathbf{1,b},(n)}_{t_2,s_1,s_2}(y,w),
\end{multline*}
and accordingly
\begin{multline}\label{contr-2-ii-bis-1-2}
\int dy_1 \, \big|\cm^{\mathbf{1,2},(n),\star}_{t_1,t_2}(y_1)\big|\lesssim\\
\lesssim \int dy_1 \, \big|\cm^{\mathbf{1,2,1},(n),\star}_{t_1,t_2}(y_1)\big| +\int dy_1 \,  \big|\cm^{\mathbf{1,2,2},(n),\star}_{t_1,t_2}(y_1)\big| +\int dy_1 \,  \big|\cm^{\mathbf{1,2,3},(n),\star}_{t_1,t_2}(y_1)\big|.
\end{multline}

\

\

\noindent
\textit{(i) \underline{Study of $\cm^{\mathbf{1,2,1},(n)}$.}} Arguing exactly as in the derivation of \eqref{ref-1-1-pap}, we readily obtain
\begin{multline*}
\int dy_1 \sup_{y_2} \big|\cm^{\mathbf{1,2,1},(n)}_{t_1,t_2}(y_1+y_2,y_2)\big|\lesssim\\
\lesssim \big\|\cac^{(n),\star}_{t_2,t_2}\big\|_{L^2}^2 \int_{0}^{t_1} ds_1 \int_{0}^{t_1} ds_2 \, \big\|K^\star_{t_1-s_1,t_2-s_1}\big\|_{L^1} \big\| K^\star_{t_1-s_2,t_2-s_2}\big\|_{L^2} \big\|\cac^{(n),\star}_{s_1,s_2}\big\|_{L^4}^2.
\end{multline*}
Then, using the results in Lemma \ref{lem:k-star-1} and Lemma \ref{lem:cac3}, we obtain for all $\theta_1,\theta_2\in [0,1]$,
\begin{multline*}
\int dy_1 \sup_{y_2} \big|\cm^{\mathbf{1,2,1},(n)}_{t_1,t_2}(y_1+y_2,y_2)\big|\lesssim\\
\begin{aligned}
&\lesssim |t_2-t_1|^{\theta_1+\theta_2} \int_{0}^{t_1} ds_1 \int_{0}^{t_1} ds_2 \, \frac{e^{-(t_1-s_1)}}{|t_1-s_1|^{\frac52\theta_1}}\frac{e^{- (t_1-s_2)}}{|t_1-s_2|^{\frac34+\frac74\theta_2}}   \frac{1}{|s_2-s_1|^{\frac14}}\\
&\lesssim |t_2-t_1|^{\theta_1+\theta_2} \int_{0}^{+\infty} ds_1 \int_{0}^{+\infty} ds_2 \, \frac{e^{-s_1}}{|s_1|^{\frac52\theta_1}}\frac{e^{- s_2}}{|s_2|^{\frac34+\frac74\theta_2}}   \frac{1}{|s_2-s_1|^{\frac14}}
\end{aligned}
\end{multline*}
and finally, by choosing for instance $\theta_1=\theta_2=\frac18$, we obtain 
\begin{multline}
\int dy_1 \sup_{y_2} \big|\cm^{\mathbf{1,2,1},(n)}_{t_1,t_2}(y_1+y_2,y_2)\big|\lesssim\\
\lesssim |t_2-t_1|^{\frac14} \int_{0}^{+\infty} ds_1 \int_{0}^{+\infty} ds_2 \, \frac{e^{-s_1}}{|s_1|^{\frac{5}{16}}}\frac{e^{- s_2}}{|s_2|^{\frac{31}{32}}}   \frac{1}{|s_2-s_1|^{\frac14}}\lesssim |t_2-t_1|^{\frac14},\label{contr-2-ii-bis-1-2-1}
\end{multline}
uniformly over $0\leq t_1,t_2\leq 2$.

\

\

\noindent
\textit{(ii) \underline{Study of $\cm^{\mathbf{1,2,2},(n)}$.}} Repeating the reasoning that led to \eqref{ref-1-2-pap}, we find
\begin{multline*}
 \int dy_1 \sup_{y_2} \big|\cm^{\mathbf{1,2,2},(n)}_{t_1,t_2}(y_1+y_2,y_2)\big| \lesssim   \big\|\cac^{(n),\star}_{t_2,t_2}\big\|_{L^2} \int_{0}^{t_1}ds_1\, \big\|\cac^{(n),\star}_{t_2,s_1}\big\|_{L^\infty}\int_{0}^{t_1}ds_2\, \big\|\cac^{(n),\star}_{t_2,s_2}\big\|_{L^\infty} \,\\
\times  \big\|K^\star_{t_1-s_1,t_2-s_1}\big\|_{L^1}\big\|K^\star_{t_1-s_2,t_2-s_2}\big\|_{L^1} \big\|\cac^{(n),\star}_{s_1,s_2}\big\|_{L^2},
\end{multline*}
Now we can use the estimates in Lemma \ref{lem:k-star-1} and Lemma \ref{lem:cac3} to assert that for all $\theta\in [0,1]$,
\begin{align*}
\int dy_1 \sup_{y_2} \big|\cm^{\mathbf{1,2,2},(n)}_{t_1,t_2}(y_1+y_2,y_2)\big|&\lesssim |t_2-t_1|^{2\theta}\bigg( \int_{0}^{t_1}ds\, \frac{e^{-(t_1-s)}}{|t_1-s|^{\frac12+\frac52\theta}} \bigg)^2\\
&\lesssim |t_2-t_1|^{2\theta}\bigg( \int_{0}^{+\infty}ds\, \frac{e^{-s}}{|s|^{\frac12+\frac52\theta}} \bigg)^2.
\end{align*}
Picking $\theta=\frac16$, we immediately obtain that
\begin{align}
\int dy_1 \sup_{y_2} \big|\cm^{\mathbf{1,2,2},(n)}_{t_1,t_2}(y_1+y_2,y_2)\big|&\lesssim |t_2-t_1|^{\frac13}\label{contr-2-ii-bis-1-2-2}
\end{align}
uniformly over $0\leq t_1,t_2\leq 2$.

\

\

\noindent
\textit{(iii) \underline{Study of $\cm^{\mathbf{1,2,3},(n)}$.}} Proceeding along the same lines as for \eqref{ref-1-3-pap}, we deduce
\begin{align*}
& \int dy_1 \sup_{y_2} \big|\cm^{\mathbf{1,2,3},(n)}_{t_1,t_2}(y_1+y_2,y_2)\big|
\lesssim  \bigg(\int_{0}^{t_1} ds \, \big\| K_{t_1-s,t_2-s}^\star \big\|_{L^1}   \big\|\cac^{(n),\star}_{t_2,s}\big\|_{L^4}^2\bigg)^2.
\end{align*}
Using again Lemma \ref{lem:k-star-1} and Lemma \ref{lem:cac3}, we obtain that for all $\theta\in [0,1]$,
\begin{eqnarray*}
 \int dy_1 \sup_{y_2} \big|\cm^{\mathbf{1,2,3},(n)}_{t_1,t_2}(y_1+y_2,y_2)\big|
&\lesssim&  |t_2-t_1|^{2\theta}\bigg(\int_{0}^{t_1}ds\, \frac{e^{-(t_1-s)}}{|t_1-s|^{\frac14+\frac52\theta}}  \bigg)^2 \\
&\lesssim & |t_2-t_1|^{2\theta}\bigg(\int_{0}^{+\infty}ds\, \frac{e^{-s}}{|s|^{\frac14+\frac52\theta}}  \bigg)^2,
\end{eqnarray*}
and by choosing $\theta=\frac14$, this yields
\begin{align}
 \int dy_1 \sup_{y_2} \big|\cm^{\mathbf{1,2,3},(n)}_{t_1,t_2}(y_1+y_2,y_2)\big| &\lesssim  |t_2-t_1|^{\frac12},\label{contr-2-ii-bis-1-2-3}
\end{align}
uniformly over $0\leq t_1,t_2\leq 2$.

\

Finally, by substituting estimates \eqref{contr-2-ii-bis-1-2-1}-\eqref{contr-2-ii-bis-1-2-2}-\eqref{contr-2-ii-bis-1-2-3} into \eqref{contr-2-ii-bis-1-2}, we can conclude that
\begin{align}\label{cm12}
&\int dy_1 \, \big|\cm^{\mathbf{1,2},(n),\star}_{t_1,t_2}(y_1)\big|\lesssim |t_2-t_1|^{\frac14},
\end{align}
uniformly over $0\leq t_1,t_2\leq 2$.

\

\subsubsection{Study of $\cm^{\mathbf{1,3},(n)}$.} One has here
\begin{multline*}
\cm^{\mathbf{1,3},(n)}_{t_1,t_2}(y_1,y_2)=\mathbb{E}\Big[\scret^{\mathbf{1,3},(n)}_{t_1,t_2}(y_1)\, \scret^{\mathbf{1,3},(n)}_{t_1,t_2}(y_2) \Big]\\
\begin{aligned}
&=\int_{0}^{t_1}ds_1\int_{0}^{t_1}ds_2 \int dw_1 dw_2 \, K_{t_1-s_1}(y_1,w_1)K_{t_1-s_2}(y_2,w_2)\\
&\hspace{2cm} \times \mathbb{E}\Big[I^W_4\Big(\big(F^{(n)}_{t_2,y_1}-F^{(n)}_{t_1,y_1}\big)\otimes F^{(n)}_{t_2,y_1}\otimes F^{(n)}_{s_1,w_1} \otimes F^{(n)}_{s_1,w_1} \Big) \\
&\hspace{4cm}\times I^W_4\Big(\big(F^{(n)}_{t_2,y_2}-F^{(n)}_{t_1,y_2}\big)\otimes F^{(n)}_{t_2,y_2}\otimes F^{(n)}_{s_2,w_2} \otimes F^{(n)}_{s_2,w_2} \Big) \Big].
\end{aligned}
\end{multline*}
The latter expectation can be expanded as
\begin{align*}
&\mathbb{E}\Big[I^W_4\Big(\big(F^{(n)}_{t_2,y_1}-F^{(n)}_{t_1,y_1}\big)\otimes F^{(n)}_{t_2,y_1}\otimes F^{(n)}_{s_1,w_1} \otimes F^{(n)}_{s_1,w_1} \Big)I^W_4\Big(\big(F^{(n)}_{t_2,y_2}-F^{(n)}_{t_1,y_2}\big)\otimes F^{(n)}_{t_2,y_2}\otimes F^{(n)}_{s_2,w_2} \otimes F^{(n)}_{s_2,w_2} \Big) \Big]\nonumber\\
& =c\, \Big\langle \text{Sym}\Big(\big(F^{(n)}_{t_2,y_1}-F^{(n)}_{t_1,y_1}\big)\otimes F^{(n)}_{t_2,y_1}\otimes F^{(n)}_{s_1,w_1} \otimes F^{(n)}_{s_1,w_1} \Big), \nonumber\\
 &\hspace{4cm}   \text{Sym}\Big(\big(F^{(n)}_{t_2,y_2}-F^{(n)}_{t_1,y_2}\big)\otimes F^{(n)}_{t_2,y_2}\otimes F^{(n)}_{s_2,w_2} \otimes F^{(n)}_{s_2,w_2} \Big) \Big\rangle_{L^2((\R_+\times \R^3)^4)}\nonumber\\
&=\sum_{\mathbf{b}=1,...,7} c_{\mathbf{b}} \cq^{\mathbf{1,3,b},(n)}_{t_1,t_2,s_1,s_2}({y,w}),
\end{align*}
for some combinatorial coefficients $c,c_{\mathbf{b}}\geq 0$, and with (recall the notation in \eqref{inc-cn})
\begin{align*}
 \cq^{\mathbf{1,3,1},(n)}_{t_1,t_2,s_1,s_2}(y,w)&:=\cac^{(n)}_{(t_1,t_2),(t_1,t_2)}(y_1,y_2) \cac^{(n)}_{t_2,t_2}(y_1,y_2)   \cac^{(n)}_{s_1,s_2}(w_1,w_2)^2,\\[2pt]
 \cq^{\mathbf{1,3,2},(n)}_{t_1,t_2,s_1,s_2}(y,w)&:=\cac^{(n)}_{(t_1,t_2),t_2}(y_1,y_2) \cac^{(n)}_{(t_1,t_2),t_2}(y_2,y_1)   \cac^{(n)}_{s_1,s_2}(w_1,w_2)^2\\[2pt]
 \cq^{\mathbf{1,3,3},(n)}_{t_1,t_2,s_1,s_2}(y,w)&:=\cac^{(n)}_{(t_1,t_2),(t_1,t_2)}(y_1,y_2) \cac^{(n)}_{t_2,s_2}(y_1,w_2)  \cac^{(n)}_{t_2,s_1}(y_2,w_1)   \cac^{(n)}_{s_1,s_2}(w_1,w_2)\\[2pt]
 \cq^{\mathbf{1,3,4},(n)}_{t_1,t_2,s_1,s_2}(y,w)&:=\cac^{(n)}_{(t_1,t_2),t_2}(y_1,y_2) \cac^{(n)}_{t_2,s_2}(y_1,w_2)  \cac^{(n)}_{(t_1,t_2),s_1}(y_2,w_1)  \cac^{(n)}_{s_1,s_2}(w_1,w_2)\\[2pt]
 \cq^{\mathbf{1,3,5},(n)}_{t_1,t_2,s_1,s_2}(y,w)&:=\cac^{(n)}_{(t_1,t_2),s_2}(y_1,w_2) \cac^{(n)}_{(t_1,t_2),t_2}(y_2,y_1)  \cac^{(n)}_{t_2,s_1}(y_2,w_1)  \cac^{(n)}_{s_1,s_2}(w_1,w_2)\\[2pt]
 \cq^{\mathbf{1,3,6},(n)}_{t_1,t_2,s_1,s_2}(y,w)&:=\cac^{(n)}_{(t_1,t_2),s_2}(y_1,w_2) \cac^{(n)}_{t_1,t_2}(y_1,y_2)  \cac^{(n)}_{(t_1,t_2),s_1}(y_2,w_1)    \cac^{(n)}_{s_1,s_2}(w_1,w_2)\\[2pt]
 \cq^{\mathbf{1,3,7},(n)}_{t_1,t_2,s_1,s_2}(y,w)&:=\cac^{(n)}_{(t_1,t_2),s_2}(y_1,w_2) \cac^{(n)}_{t_2,s_2}(y_1,w_2)   \cac^{(n)}_{(t_1,t_2),s_1}(y_2,w_1) \cac^{(n)}_{t_2,s_1}(y_2,w_1).
\end{align*}
As a consequence,
\begin{align*}
\cm^{\mathbf{1,3},(n)}_{t_1,t_2}(y_1,y_2)=\sum_{\mathbf{b}=1,...,7} c_{\mathbf{b}}\cm^{\mathbf{1,3,b},(n)}_{t_1,t_2}(y_1,y_2),
\end{align*}
with
\begin{align*}
\cm^{\mathbf{1,3,b},(n)}_{t_1,t_2}(y_1,y_2):=\int_{0}^{t_1}ds_1\int_{0}^{t_1}ds_2\int dw_1 dw_2 \,K_{t_1-s_1}(y_1,w_1)K_{t_1-s_2}(y_2,w_2)\cq^{\mathbf{1,3,b},(n)}_{t_1,t_2,s_1,s_2}(y,w),
\end{align*}
and therefore
\begin{align}
&\int dy_1 \, \big|\cm^{\mathbf{1,3},(n),\star}_{t_1,t_2}(y_1)\big|\lesssim \sum_{\mathbf{b}=1,...,7}\int dy_1 \, \big|\cm^{\mathbf{1,3,b},(n),\star}_{t_1,t_2}(y_1)\big|.\label{contr-2-ii-bis-1-3}
\end{align}

\

\

\noindent
\textit{(i) \underline{Study of $\cm^{\mathbf{1,3,1},(n)}$.}} Following the same strategy as in the proof of \eqref{ref-1-1-pap}, we obtain
\begin{multline*}
\int dy_1 \sup_{y_2} \big|\cm^{\mathbf{1,3,1},(n)}_{t_1,t_2}(y_1+y_2,y_2)\big|\lesssim \\
\lesssim \big\|\cac^{(n),\star}_{(t_1,t_2),(t_1,t_2)}\big\|_{L^2} \big\|{\cac}^{(n),\star}_{t_2,t_2}\big\|_{L^2}\int_{0}^{t_1} ds_1 \int_{0}^{t_1} ds_2 \, \big\|K^\star_{t_1-s_1}\big\|_{L^1} \big\| K^\star_{t_1-s_2}\big\|_{L^2} \big\|\cac^{(n),\star}_{s_1,s_2}\big\|_{L^4}^2.
\end{multline*}
Then, using the results in Lemma \ref{lem:k-star-1} and Lemma \ref{lem:cac3}, we obtain that
\begin{eqnarray*}
\int dy_1 \sup_{y_2} \big|\cm^{\mathbf{1,3,1},(n)}_{t_1,t_2}(y_1+y_2,y_2)\big|&\lesssim &
  |t_2-t_1|^{\frac14} \int_{0}^{t_1} ds_1 \int_{0}^{t_1} ds_2 \, e^{-(t_1-s_1)}\frac{e^{- (t_1-s_2)}}{|t_1-s_2|^{\frac34}}   \frac{1}{|s_2-s_1|^{\frac14}}\\
&\lesssim &|t_2-t_1|^{\frac14} \int_{0}^{+\infty} ds_1 \int_{0}^{+\infty} ds_2 \, e^{-s_1}\frac{e^{- s_2}}{|s_2|^{\frac34}}   \frac{1}{|s_2-s_1|^{\frac14}}
\end{eqnarray*}
and thus we have
\begin{align}
\int dy_1 \sup_{y_2} \big|\cm^{\mathbf{1,3,1},(n)}_{t_1,t_2}(y_1+y_2,y_2)\big|&\lesssim |t_2-t_1|^{\frac14},\label{contr-2-ii-bis-1-3-1}
\end{align}
uniformly over $0\leq t_1,t_2\leq 2$.

\

\

\noindent
\textit{(ii) \underline{Study of $\cm^{\mathbf{1,3,2},(n)}$.}} By an analogous computation to the one yielding \eqref{ref-1-1-pap}, we deduce
\begin{multline*}
\int dy_1 \sup_{y_2} \big|\cm^{\mathbf{1,3,2},(n)}_{t_1,t_2}(y_1+y_2,y_2)\big|\lesssim \\
\lesssim \big\|\cac^{(n),\star}_{(t_1,t_2),t_2}\big\|_{L^2}^2 \int_{0}^{t_1} ds_1 \int_{0}^{t_1} ds_2 \, \big\|K^\star_{t_1-s_1}\big\|_{L^1} \big\| K^\star_{t_1-s_2}\big\|_{L^2} \big\|\cac^{(n),\star}_{s_1,s_2}\big\|_{L^4}^2.
\end{multline*}
Using the results in Lemma \ref{lem:k-star-1} and Lemma \ref{lem:cac3}, we obtain that
\begin{align*}
\int dy_1 \sup_{y_2} \big|\cm^{\mathbf{1,3,2},(n)}_{t_1,t_2}(y_1+y_2,y_2)\big|&\lesssim |t_2-t_1|^{\frac12} \int_{0}^{t_1} ds_1 \int_{0}^{t_1} ds_2 \, e^{-(t_1-s_1)}\frac{e^{- (t_1-s_2)}}{|t_1-s_2|^{\frac34}}   \frac{1}{|s_2-s_1|^{\frac14}}\\
&\lesssim |t_2-t_1|^{\frac12} \int_{0}^{+\infty} ds_1 \int_{0}^{+\infty} ds_2 \, e^{-s_1}\frac{e^{- s_2}}{|s_2|^{\frac34}}   \frac{1}{|s_2-s_1|^{\frac14}}
\end{align*}
and thus we have
\begin{align}
\int dy_1 \sup_{y_2} \big|\cm^{\mathbf{1,3,2},(n)}_{t_1,t_2}(y_1+y_2,y_2)\big|&\lesssim |t_2-t_1|^{\frac12},\label{contr-2-ii-bis-1-3-2}
\end{align}
uniformly over $0\leq t_1,t_2\leq 2$.

\

\

\noindent
\textit{(iii) \underline{Study of $\cm^{\mathbf{1,3,3},(n)}$.}} Applying the estimates in the same manner as for \eqref{ref-1-2-pap}, we deduce
\begin{multline*}
 \int dy_1 \sup_{y_2} \big|\cm^{\mathbf{1,3,3},(n)}_{t_1,t_2}(y_1+y_2,y_2)\big|\lesssim  \\
\lesssim \big\|\cac^{(n),\star}_{(t_1,t_2),(t_1,t_2)}\big\|_{L^2} \int_{0}^{t_1}ds_1\, \big\|\cac^{(n),\star}_{t_2,s_1}\big\|_{L^\infty}\int_{0}^{t_1}ds_2\, \big\|\cac^{(n),\star}_{t_2,s_2}\big\|_{L^\infty} \, \big\|K^\star_{t_1-s_1}\big\|_{L^1}\big\|K^\star_{t_1-s_2}\big\|_{L^1} \big\|\cac^{(n),\star}_{s_1,s_2}\big\|_{L^2},
\end{multline*}
and using estimates in Lemma \ref{lem:k-star-1} and Lemma \ref{lem:cac3}, we obtain that 
\begin{align*}
\int dy_1 \sup_{y_2} \big|\cm^{\mathbf{1,3,3},(n)}_{t_1,t_2}(y_1+y_2,y_2)\big|&\lesssim |t_2-t_1|^{\frac14}\bigg( \int_{0}^{t_1}ds\, \frac{e^{-(t_1-s)}}{|t_1-s|^{\frac12}} \bigg)^2\lesssim |t_2-t_1|^{\frac14}\bigg( \int_{0}^{+\infty}ds\, \frac{e^{-s}}{|s|^{\frac12}} \bigg)^2.
\end{align*}
We have thus shown that
\begin{align}
\int dy_1 \sup_{y_2} \big|\cm^{\mathbf{1,3,3},(n)}_{t_1,t_2}(y_1+y_2,y_2)\big|&\lesssim |t_2-t_1|^{\frac14}\label{contr-2-ii-bis-1-3-3}
\end{align}
uniformly over $0\leq t_1,t_2\leq 2$.

\

\

\noindent
\textit{(iv) \underline{Study of $\cm^{\mathbf{1,3,4},(n)}$.}} A similar analysis to the one carried out for \eqref{ref-1-2-pap} yields
\begin{multline*}
 \int dy_1 \sup_{y_2} \big|\cm^{\mathbf{1,3,4},(n)}_{t_1,t_2}(y_1+y_2,y_2)\big| \lesssim \\
\lesssim \big\|\cac^{(n),\star}_{(t_1,t_2),t_2}\big\|_{L^2} \int_{0}^{t_1}ds_1\, \big\|\cac^{(n),\star}_{(t_1,t_2),s_1}\big\|_{L^\infty}\int_{0}^{t_1}ds_2\, \big\|\cac^{(n),\star}_{t_2,s_2}\big\|_{L^\infty} \, \big\|K^\star_{t_1-s_1}\big\|_{L^1}\big\|K^\star_{t_1-s_2}\big\|_{L^1} \big\|\cac^{(n),\star}_{s_1,s_2}\big\|_{L^2},
\end{multline*}
and using estimates in Lemma \ref{lem:k-star-1} and Lemma \ref{lem:cac3}, we derive that
\begin{align*}
\int dy_1 \sup_{y_2} \big|\cm^{\mathbf{1,3,4},(n)}_{t_1,t_2}(y_1+y_2,y_2)\big|&\lesssim |t_2-t_1|^{\frac14}\bigg( \int_{0}^{t_1}ds\, \frac{e^{-(t_1-s)}}{|t_1-s|^{\frac12}} \bigg)^2\lesssim |t_2-t_1|^{\frac14}\bigg( \int_{0}^{+\infty}ds\, \frac{e^{-s}}{|s|^{\frac12}} \bigg)^2.
\end{align*}
We have thus shown that
\begin{align}
\int dy_1 \sup_{y_2} \big|\cm^{\mathbf{1,3,4},(n)}_{t_1,t_2}(y_1+y_2,y_2)\big|&\lesssim |t_2-t_1|^{\frac14}\label{contr-2-ii-bis-1-3-4}
\end{align}
uniformly over $0\leq t_1,t_2\leq 2$.

\

\

\noindent
\textit{(v) \underline{Study of $\cm^{\mathbf{1,3,5},(n)}$.}} Running the same estimates as in the proof of \eqref{ref-1-2-pap}, we obtain
\begin{multline*}
 \int dy_1 \sup_{y_2} \big|\cm^{\mathbf{1,3,5},(n)}_{t_1,t_2}(y_1+y_2,y_2)\big| \lesssim \\
\lesssim \big\|\cac^{(n),\star}_{(t_1,t_2),t_2}\big\|_{L^2} \int_{0}^{t_1}ds_1\, \big\|\cac^{(n),\star}_{t_2,s_1}\big\|_{L^\infty}\int_{0}^{t_1}ds_2\, \big\|\cac^{(n),\star}_{(t_1,t_2),s_2}\big\|_{L^\infty} \, \big\|K^\star_{t_1-s_1}\big\|_{L^1}\big\|K^\star_{t_1-s_2}\big\|_{L^1} \big\|\cac^{(n),\star}_{s_1,s_2}\big\|_{L^2},
\end{multline*}
and using estimates in Lemma \ref{lem:k-star-1} and Lemma \ref{lem:cac3}, we derive that
\begin{align*}
\int dy_1 \sup_{y_2} \big|\cm^{\mathbf{1,3,5},(n)}_{t_1,t_2}(y_1+y_2,y_2)\big|&\lesssim |t_2-t_1|^{\frac14}\bigg( \int_{0}^{t_1}ds\, \frac{e^{-(t_1-s)}}{|t_1-s|^{\frac12}} \bigg)^2\lesssim |t_2-t_1|^{\frac14}\bigg( \int_{0}^{+\infty}ds\, \frac{e^{-s}}{|s|^{\frac12}} \bigg)^2.
\end{align*}
We have thus shown that
\begin{align}
\int dy_1 \sup_{y_2} \big|\cm^{\mathbf{1,3,5},(n)}_{t_1,t_2}(y_1+y_2,y_2)\big|&\lesssim |t_2-t_1|^{\frac14}\label{contr-2-ii-bis-1-3-5}
\end{align}
uniformly over $0\leq t_1,t_2\leq 2$.

\

\

\noindent
\textit{(vi) \underline{Study of $\cm^{\mathbf{1,3,6},(n)}$.}} Mimicking the computation that produced \eqref{ref-1-2-pap}, we obtain
\begin{multline*}
 \int dy_1 \sup_{y_2} \big|\cm^{\mathbf{1,3,6},(n)}_{t_1,t_2}(y_1+y_2,y_2)\big| \lesssim \\
\lesssim \big\|\cac^{(n),\star}_{t_2,t_2}\big\|_{L^2} \int_{0}^{t_1}ds_1\, \big\|\cac^{(n),\star}_{(t_1,t_2),s_1}\big\|_{L^\infty}\int_{0}^{t_1}ds_2\, \big\|\cac^{(n),\star}_{(t_1,t_2),s_2}\big\|_{L^\infty} \, \big\|K^\star_{t_1-s_1}\big\|_{L^1}\big\|K^\star_{t_1-s_2}\big\|_{L^1} \big\|\cac^{(n),\star}_{s_1,s_2}\big\|_{L^2},
\end{multline*}
and using estimates in Lemma \ref{lem:k-star-1} and Lemma \ref{lem:cac3}, we derive that
\begin{align*}
\int dy_1 \sup_{y_2} \big|\cm^{\mathbf{1,3,6},(n)}_{t_1,t_2}(y_1+y_2,y_2)\big|&\lesssim |t_2-t_1|^{\frac12}\bigg( \int_{0}^{t_1}ds\, \frac{e^{-(t_1-s)}}{|t_1-s|^{\frac34}} \bigg)^2\lesssim |t_2-t_1|^{\frac12}\bigg( \int_{0}^{+\infty}ds\, \frac{e^{-s}}{|s|^{\frac34}} \bigg)^2.
\end{align*}
We have thus shown that
\begin{align}
\int dy_1 \sup_{y_2} \big|\cm^{\mathbf{1,3,6},(n)}_{t_1,t_2}(y_1+y_2,y_2)\big|&\lesssim |t_2-t_1|^{\frac12}\label{contr-2-ii-bis-1-3-6}
\end{align}
uniformly over $0\leq t_1,t_2\leq 2$.

\

\

\noindent
\textit{(vii) \underline{Study of $\cm^{\mathbf{1,3,7},(n)}$.}} The same line of reasoning as for \eqref{ref-1-3-pap} gives
\begin{align*}
\int dy_1 \sup_{y_2} \big|\cm^{\mathbf{1,3,7},(n)}_{t_1,t_2}(y_1+y_2,y_2)\big|
&\lesssim  \bigg(\int_{0}^{t_1} ds \, \big\| K_{t_1-s}^\star \big\|_{L^1}   \big\|\cac^{(n),\star}_{(t_1,t_2),s}\big\|_{L^4}\big\|\cac^{(n),\star}_{t_2,s}\big\|_{L^4}\bigg)^2.
\end{align*}
Using again Lemma \ref{lem:k-star-1} and Lemma \ref{lem:cac3}, we obtain that
\begin{align*}
\int dy_1 \sup_{y_2} \big|\cm^{\mathbf{1,3,7},(n)}_{t_1,t_2}(y_1+y_2,y_2)\big|
&\lesssim  |t_2-t_1|\bigg(\int_{0}^{t_1}ds\, \frac{e^{-(t_1-s)}}{|t_1-s|^{\frac34}}  \bigg)^2\lesssim  |t_2-t_1|\bigg(\int_{0}^{+\infty}ds\, \frac{e^{-s}}{|s|^{\frac34}}  \bigg)^2,
\end{align*}
which yields
\begin{align}
 \int dy_1 \sup_{y_2} \big|\cm^{\mathbf{1,3,7},(n)}_{t_1,t_2}(y_1+y_2,y_2)\big| &\lesssim  |t_2-t_1|,\label{contr-2-ii-bis-1-3-7}
\end{align}
uniformly over $0\leq t_1,t_2\leq 2$.

\

Finally, by substituting estimates \eqref{contr-2-ii-bis-1-3-1} to \eqref{contr-2-ii-bis-1-3-7} into \eqref{contr-2-ii-bis-1-3}, we can conclude that
\begin{align}\label{cm13}
&\int dy_1 \, \big|\cm^{\mathbf{1,3},(n),\star}_{t_1,t_2}(y_1)\big|\lesssim |t_2-t_1|^{\frac14},
\end{align}
uniformly over $0\leq t_1,t_2\leq 2$.

\

\

\subsubsection{Study of $\cm^{\mathbf{1,4},(n)}$.}

The same arguments as those used for $\scret^{\mathbf{1,3},(n)}$ yield
\begin{align}\label{cm14}
&\int dy_1 \, \big|\cm^{\mathbf{1,4},(n),\star}_{t_1,t_2}(y_1)\big|\lesssim |t_2-t_1|^{\frac14},
\end{align}
uniformly over $0\leq t_1,t_2\leq 2$.

\

We are finally in a position to substitute \eqref{cm11}, \eqref{cm12}, \eqref{cm13} and \eqref{cm14} into \eqref{contr-1-bis}: as a result, we obtain that
\begin{equation}\label{scretun}
\int dx \, \mathbb{E}\bigg[ \Big| H^{-\al}\big(\scret^{\mathbf{1},(n)}_{t_1,t_2}\big)(x)\Big|^{2}\bigg]^p\lesssim |t_2-t_1|^{\frac{p}{4}},
\end{equation}
uniformly over $0\leq t_1,t_2\leq 2$.

\

\

\noindent
\textbf{2. Study of $\scret^{\mathbf{2},(n)}_{t_1,t_2}$.}

\smallskip

Using the notation introduced in \eqref{inc-k}, we can write 
\begin{multline*}
\scret^{\mathbf{2},(n)}_{t_1,t_2} (y) =\int_{t_1}^{t_2} ds \int dw \, K_{t_2-s}(y,w) \,\cac^{(n)}_{t_2,s}(y,w) \,  I^W_2\big( F^{(n)}_{t_2,y} \otimes F^{(n)}_{s,w} \big)\\
\begin{aligned}
&+\int_{0}^{t_1} ds \int dw \, K_{t_1-s,t_2-s}(y,w) \,\cac^{(n)}_{t_2,s}(y,w) \,  I^W_2\big( F^{(n)}_{t_2,y} \otimes F^{(n)}_{s,w} \big)\\
&+\int_{0}^{t_1} ds \int dw \, K_{t_1-s}(y,w) \,\cac^{(n)}_{(t_1,t_2),s}(y,w) \,  I^W_2\big( F^{(n)}_{t_2,y} \otimes F^{(n)}_{s,w} \big)\\
&+\int_{0}^{t_1} ds \int dw \, K_{t_1-s}(y,w) \,\cac^{(n)}_{t_1,s}(y,w) \,  I^W_2\big( (F^{(n)}_{t_2,y}-F^{(n)}_{t_1,y}) \otimes F^{(n)}_{s,w} \big)\\
&\hspace{2cm}=:\scret^{\mathbf{2,1},(n)}_{t_1,t_2} (y)+\scret^{\mathbf{2,2},(n)}_{t_1,t_2} (y)+\scret^{\mathbf{2,3},(n)}_{t_1,t_2} (y)+\scret^{\mathbf{2,4},(n)}_{t_1,t_2} (y).
\end{aligned}
\end{multline*}

\

For $\mathbf{a}=1,\ldots,4$, let us set
$$\cm^{\mathbf{2,a},(n)}_{t_1,t_2}(y_1,y_2):=\mathbb{E}\Big[\scret^{\mathbf{2,a},(n)}_{t_1,t_2}(y_1)\, \scret^{\mathbf{2,a},(n)}_{t_1,t_2}(y_2) \Big].$$
Then, thanks to Lemma \ref{lem:techn} and since $\al>\frac34$, we obtain that 
\begin{equation}\label{contr-1-bis-2}
\int dx \, \mathbb{E}\bigg[ \Big| H^{-\al}\big(\scret^{\mathbf{2},(n)}_{t_1,t_2}\big)(x)\Big|^{2}\bigg]^p\lesssim \sum_{\mathbf{a}=1,\ldots,4}\bigg(\int dy_1 \,  \big|\cm^{\mathbf{2,a},(n),\star}_{t_1,t_2}(y_1)\big| \bigg)^p,
\end{equation}
for every $p$ large enough, where we have used the notation $\star$ introduced in \eqref{star-not}.

\

\subsubsection{Study of $\cm^{\mathbf{2,1},(n)}$.} One has here
\begin{multline*}
\cm^{\mathbf{2,1},(n)}_{t_1,t_2}(y_1,y_2)=\mathbb{E}\Big[\scret^{\mathbf{2,1},(n)}_{t_1,t_2}(y_1)\, \scret^{\mathbf{2,1},(n)}_{t_1,t_2}(y_2) \Big]=\\
\begin{aligned}
&=\int_{t_1}^{t_2}ds_1\int_{t_1}^{t_2}ds_2 \int dw_1 dw_2 \, K_{t_2-s_1}(y_1,w_1)K_{t_2-s_2}(y_2,w_2) \cac^{(n)}_{t_2,s_1}(y_1,w_1) \cac^{(n)}_{t_2,s_2}(y_2,w_2)\\
&\hspace{2cm}\times \mathbb{E}\Big[I^W_2\big( F^{(n)}_{t_2,y_1} \otimes F^{(n)}_{s_1,w_1} \big)I^W_2\big( F^{(n)}_{t_2,y_2} \otimes F^{(n)}_{s_2,w_2} \big)\Big]\\
&=\int_{t_1}^{t_2}ds_1\int_{t_1}^{t_2}ds_2 \int dw_1 dw_2 \, K_{t_2-s_1}(y_1,w_1)K_{t_2-s_2}(y_2,w_2) \cac^{(n)}_{t_2,s_1}(y_1,w_1) \cac^{(n)}_{t_2,s_2}(y_2,w_2)\\
&\hspace{2cm}\times \Big[c_{\mathbf{1}}\, \cac^{(n)}_{t_2,t_2}(y_1,y_2) \cac^{(n)}_{s_1,s_2}(w_1,w_2) +c_{\mathbf{2}}\, \cac^{(n)}_{t_2,s_2}(y_1,w_2) \cac^{(n)}_{t_2,s_1}(y_2,w_1) \Big],
\end{aligned}
\end{multline*}
for some combinatorial coefficients $c_{\mathbf{1}},c_{\mathbf{2}}\geq 0$. Thus,
\begin{align*}
\cm^{\mathbf{2,1},(n)}_{t_1,t_2}(y_1,y_2)= c_{\mathbf{1}}\cm^{\mathbf{2,1,1},(n)}_{t_1,t_2}(y_1,y_2)+ c_{\mathbf{2}}\cm^{\mathbf{2,1,2},(n)}_{t_1,t_2}(y_1,y_2),
\end{align*}
with
\begin{multline*}
\cm^{\mathbf{2,1,1},(n)}_{t_1,t_2}(y_1,y_2):=\int_{t_1}^{t_2}ds_1\int_{t_1}^{t_2}ds_2\int dw_1 dw_2 \, K_{t_2-s_1}(y_1,w_1) K_{t_2-s_2}(y_2,w_2)\\
\hspace{3cm}\times \cac^{(n)}_{t_2,s_1}(y_1,w_1) \cac^{(n)}_{t_2,s_2}(y_2,w_2)\cac^{(n)}_{t_2,t_2}(y_1,y_2) \cac^{(n)}_{s_1,s_2}(w_1,w_2),
\end{multline*}
and 
\begin{multline*}
\cm^{\mathbf{2,1,2},(n)}_{t_1,t_2}(y_1,y_2):=\int_{t_1}^{t_2}ds_1\int_{t_1}^{t_2}ds_2\int dw_1 dw_2 \, K_{t_2-s_1}(y_1,w_1) K_{t_2-s_2}(y_2,w_2)\\
\hspace{3cm}\times \cac^{(n)}_{t_2,s_1}(y_1,w_1) \cac^{(n)}_{t_2,s_2}(y_2,w_2)\cac^{(n)}_{t_2,s_2}(y_1,w_2) \cac^{(n)}_{t_2,s_1}(y_2,w_1).
\end{multline*}
In this way,
\begin{align}
\int dy_1 \, \big|\cm^{\mathbf{2,1},(n),\star}_{t_1,t_2}(y_1)\big|&\lesssim \int dy_1 \, \big|\cm^{\mathbf{2,1,1},(n),\star}_{t_1,t_2}(y_1)\big| +\int dy_1 \,  \big|\cm^{\mathbf{2,1,2},(n),\star}_{t_1,t_2}(y_1)\big| .\label{contr-2-ii-bis-2-1}
\end{align}

\

\noindent
\textit{(i) \underline{Study of $\cm^{\mathbf{2,1,1},(n)}$.}} One has for all $y_1,y_2$,
\begin{align*} 
&\big|\cm^{\mathbf{2,1,1},(n)}_{t_1,t_2}(y_1+y_2,y_2)\big| \lesssim\\
&\lesssim \cac^{(n)}_{t_2,t_2}(y_1+y_2,y_2)\int_{t_1}^{t_2} ds_1 \int_{t_1}^{t_2} ds_2 \int dw_1\, \big|K_{t_2-s_1}(y_1+y_2,w_1) \big| \cac^{(n)}_{t_2,s_1}(y_1+y_2,w_1)\\
&\hspace{4cm}\times \int dw_2 \, \big| K_{t_2-s_2}(y_2,w_2) \big|\cac^{(n)}_{t_2,s_2}(y_2,w_2)  \cac^{(n)}_{s_1,s_2}(w_1,w_2)  \\
&\lesssim \cac^{(n),\star}_{t_2,t_2}(y_1)\int_{t_1}^{t_2} ds_1 \int_{t_1}^{t_2} ds_2 \int dw_1\, \big|K_{t_2-s_1}(y_1+y_2,w_1+y_2) \big| \cac^{(n)}_{t_2,s_1}(y_1+y_2,w_1+y_2)\\
&\hspace{4cm}\times \int dw_2 \, \big| K_{t_2-s_2}(y_2,w_2+y_2) \big|\cac^{(n)}_{t_2,s_2}(y_2,w_2+y_2)  \cac^{(n)}_{s_1,s_2}(w_1+y_2,w_2+y_2)  \\
&\lesssim \cac^{(n),\star}_{t_2,t_2}(y_1)\int_{t_1}^{t_2} ds_1 \int_{t_1}^{t_2} ds_2 \int dw_1\, K^\star_{t_2-s_1}(y_1-w_1) \cac^{(n),\star}_{t_2,s_1}(y_1-w_1) \\
&\hspace{4cm}\times\int dw_2 \,  K^\star_{t_2-s_2}(w_2) \cac^{(n),\star}_{t_2,s_2}(w_2)  \cac^{(n),\star}_{s_1,s_2}(w_1-w_2)  \\
&\lesssim \cac^{(n),\star}_{t_2,t_2}(y_1)\int_{t_1}^{t_2} ds_1 \int_{t_1}^{t_2} ds_2 \,\Big[ \big( K^\star_{t_2-s_1}\cdot \cac^{(n),\star}_{t_2,s_1} \big)\ast\Big(\cac^{(n),\star}_{s_1,s_2} \ast \big(K^\star_{t_2-s_2}\cdot \cac^{(n),\star}_{t_2,s_2}\big)\Big) \Big](y_1).
\end{align*}
Thus,
\begin{multline}\label{ref-2-1-pap}
 \int dy_1 \sup_{y_2} \big|\cm^{\mathbf{2,1,1},(n)}_{t_1,t_2}(y_1+y_2,y_2)\big|\lesssim\\
\begin{aligned}
&\lesssim \big\|\cac^{(n),\star}_{t_2,t_2}\big\|_{L^2} \int_{t_1}^{t_2} ds_1 \int_{t_1}^{t_2} ds_2 \, \Big\|\big( K^\star_{t_2-s_1}\cdot \cac^{(n),\star}_{t_2,s_1} \big)\ast\Big(\cac^{(n),\star}_{s_1,s_2} \ast \big(K^\star_{t_2-s_2}\cdot \cac^{(n),\star}_{t_2,s_2}\big)\Big)\Big\|_{L^2} \\
&\lesssim \big\|\cac^{(n),\star}_{t_2,t_2}\big\|_{L^2} \int_{t_1}^{t_2} ds_1 \int_{t_1}^{t_2} ds_2 \, \big\| K^\star_{t_2-s_1}\cdot \cac^{(n),\star}_{t_2,s_1} \big\|_{L^1}\big\|\cac^{(n),\star}_{s_1,s_2}\big\|_{L^2} \big\|K^\star_{t_2-s_2}\cdot \cac^{(n),\star}_{t_2,s_2}\big\|_{L^1} \\
&\lesssim \big\|\cac^{(n),\star}_{t_2,t_2}\big\|_{L^2} \int_{t_1}^{t_2} ds_1 \int_{t_1}^{t_2} ds_2\, \big\| K^\star_{t_2-s_1}\big\|_{L^1} \big\|\cac^{(n),\star}_{t_2,s_1} \big\|_{L^\infty}\big\|\cac^{(n),\star}_{s_1,s_2}\big\|_{L^2} \big\|K^\star_{t_2-s_2}\big\|_{L^1} \big\|\cac^{(n),\star}_{t_2,s_2}\big\|_{L^\infty}.
\end{aligned}
\end{multline}
Then, thanks to Lemma \ref{lem:k-star-0} and Lemma \ref{lem:cac3}, we deduce
\begin{align}\label{contr-2-ii-bis-2-1-1}
&\int dy_1 \sup_{y_2} \big|\cm^{\mathbf{2,1,1},(n)}_{t_1,t_2}(y_1+y_2,y_2)\big|\lesssim \bigg(\int_{t_1}^{t_2} \frac{ds}{|t_2-s|^{\frac12}}\bigg)^2 \lesssim   |t_2-t_1|,
\end{align}
uniformly over $0\leq t_1,t_2\leq 2$.

\

\noindent
\textit{(ii) \underline{Study of $\cm^{\mathbf{2,1,2},(n)}$.}} 
In this case, for all $y_1,y_2$,
\begin{multline*} 
\big|\cm^{\mathbf{2,1,2},(n)}_{t_1,t_2}(y_1+y_2,y_2)\big| \lesssim\\
\begin{aligned} 
&\lesssim \int_{t_1}^{t_2}ds_1\int_{t_1}^{t_2}ds_2 \int dw_1\, \big|K_{t_2-s_1}(y_1+y_2,w_1+y_2) \big| \cac^{(n)}_{t_2,s_1}(y_1+y_2,w_1+y_2) \cac^{(n)}_{t_2,s_1}(y_2,w_1+y_2)\\
&\hspace{4cm}\times \int dw_2 \, \big| K_{t_2-s_2}(y_2,w_2+y_2) \big|  \cac^{(n)}_{t_2,s_2}(y_2,w_2+y_2)\cac^{(n)}_{t_2,s_2}(y_1+y_2,w_2+y_2) \\
&\lesssim \int_{t_1}^{t_2}ds_1 \int_{t_1}^{t_2}ds_2 \int dw_1\, K^\star_{t_2-s_1}(y_1-w_1) \cac^{(n),\star}_{t_2,s_1}(y_1-w_1) \cac^{(n),\star}_{t_2,s_1}(w_1)\\
&\hspace{4cm}\times \int dw_2 \,  K^\star_{t_2-s_2}(w_2)   \cac^{(n),\star}_{t_2,s_2}(w_2)\cac^{(n),\star}_{t_2,s_2}(y_1-w_2) \\
&\lesssim \int_{t_1}^{t_2}ds_1 \int_{t_1}^{t_2}ds_2 \, \Big[ \big(K^\star_{t_2-s_1}\cdot \cac^{(n),\star}_{t_2,s_1}\big)\ast \cac^{(n),\star}_{t_2,s_1}\Big](y_1)\Big[ \big(  K^\star_{t_2-s_2}\cdot   \cac^{(n),\star}_{t_2,s_2}\big)\ast \cac^{(n),\star}_{t_2,s_2}\Big](y_1) .
\end{aligned}
\end{multline*} 
Therefore,
\begin{align}
\int dy_1 \sup_{y_2} \big|\cm^{\mathbf{2,1,2},(n)}_{t_1,t_2}(y_1+y_2,y_2)\big|
&\lesssim \bigg(\int_{t_1}^{t_2}ds\,  \Big\| \big(K^\star_{t_2-s}\cdot \cac^{(n),\star}_{t_2,s}\big)\ast \cac^{(n),\star}_{t_2,s}\Big\|_{L^2} \bigg)^2\nonumber\\
&\lesssim \bigg(\int_{t_1}^{t_2}ds\,  \big\| K^\star_{t_2-s}\cdot \cac^{(n),\star}_{t_2,s}\big\|_{L^1} \big\|\cac^{(n),\star}_{t_2,s}\big\|_{L^2} \bigg)^2\nonumber\\
&\lesssim \bigg(\int_{t_1}^{t_2}ds\,  \big\| K^\star_{t_2-s}\big\|_{L^1} \big\|\cac^{(n),\star}_{t_2,s}\big\|_{L^\infty} \big\|\cac^{(n),\star}_{t_2,s}\big\|_{L^2} \bigg)^2. \label{ref-2-2-pap}
\end{align}
Using the results in Lemma \ref{lem:k-star-0} and Lemma \ref{lem:cac3}, this yields
\begin{align}\label{contr-2-ii-bis-2-1-2}
&\int dy_1 \sup_{y_2} \big|\cm^{\mathbf{2,1,2},(n)}_{t_1,t_2}(y_1+y_2,y_2)\big|\lesssim \bigg(\int_{t_1}^{t_2} \frac{ds}{|t_2-s|^{\frac12}}\bigg)^2 \lesssim   |t_2-t_1|,
\end{align}
uniformly over $0\leq t_1,t_2\leq 2$.

\

By inserting the two estimates \eqref{contr-2-ii-bis-2-1-1} and \eqref{contr-2-ii-bis-2-1-2} into \eqref{contr-2-ii-bis-2-1}, we can conclude that
\begin{align}\label{cm21}
&\int dy_1 \, \big|\cm^{\mathbf{2,1},(n),\star}_{t_1,t_2}(y_1)\big|\lesssim |t_2-t_1|^{\frac14},
\end{align}
uniformly over $0\leq t_1,t_2\leq 2$.

\

\

\subsubsection{Study of $\cm^{\mathbf{2,2},(n)}$.} One has here
\begin{align*}
&\cm^{\mathbf{2,2},(n)}_{t_1,t_2}(y_1,y_2)=\mathbb{E}\Big[\scret^{\mathbf{2,2},(n)}_{t_1,t_2}(y_1)\, \scret^{\mathbf{2,2},(n)}_{t_1,t_2}(y_2) \Big]\\
&=\int_{0}^{t_1}ds_1\int_{0}^{t_1}ds_2 \int dw_1 dw_2 \, K_{t_1-s_1,t_2-s_1}(y_1,w_1)K_{t_1-s_2,t_2-s_2}(y_2,w_2) \cac^{(n)}_{t_2,s_1}(y_1,w_1) \cac^{(n)}_{t_2,s_2}(y_2,w_2)\\
&\hspace{2cm}\times \mathbb{E}\Big[I^W_2\big( F^{(n)}_{t_2,y_1} \otimes F^{(n)}_{s_1,w_1} \big)I^W_2\big( F^{(n)}_{t_2,y_2} \otimes F^{(n)}_{s_2,w_2} \big)\Big]\\
&=\int_{0}^{t_1}ds_1\int_{0}^{t_1}ds_2 \int dw_1 dw_2 \, K_{t_1-s_1,t_2-s_1}(y_1,w_1)K_{t_1-s_2,t_2-s_2}(y_2,w_2) \cac^{(n)}_{t_2,s_1}(y_1,w_1) \cac^{(n)}_{t_2,s_2}(y_2,w_2)\\
&\hspace{2cm}\times \Big[c_{\mathbf{1}}\, \cac^{(n)}_{t_2,t_2}(y_1,y_2) \cac^{(n)}_{s_1,s_2}(w_1,w_2) +c_{\mathbf{2}}\, \cac^{(n)}_{t_2,s_2}(y_1,w_2) \cac^{(n)}_{t_2,s_1}(y_2,w_1) \Big],
\end{align*}
for some combinatorial coefficients $c_{\mathbf{1}},c_{\mathbf{2}}\geq 0$. Thus,
\begin{align*}
\cm^{\mathbf{2,2},(n)}_{t_1,t_2}(y_1,y_2)= c_{\mathbf{1}}\cm^{\mathbf{2,2,1},(n)}_{t_1,t_2}(y_1,y_2)+ c_{\mathbf{2}}\cm^{\mathbf{2,2,2},(n)}_{t_1,t_2}(y_1,y_2),
\end{align*}
with
\begin{align*}
&\cm^{\mathbf{2,2,1},(n)}_{t_1,t_2}(y_1,y_2):=\int_{0}^{t_1}ds_1\int_{0}^{t_1}ds_2\int dw_1 dw_2 \,  K_{t_1-s_1,t_2-s_1}(y_1,w_1)K_{t_1-s_2,t_2-s_2}(y_2,w_2)\\
&\hspace{5cm}\times\cac^{(n)}_{t_2,s_1}(y_1,w_1) \cac^{(n)}_{t_2,s_2}(y_2,w_2)\cac^{(n)}_{t_2,t_2}(y_1,y_2) \cac^{(n)}_{s_1,s_2}(w_1,w_2),
\end{align*}
and 
\begin{align*}
&\cm^{\mathbf{2,2,2},(n)}_{t_1,t_2}(y_1,y_2):=\int_{0}^{t_1}ds_1\int_{0}^{t_1}ds_2\int dw_1 dw_2 \,  K_{t_1-s_1,t_2-s_1}(y_1,w_1)K_{t_1-s_2,t_2-s_2}(y_2,w_2)\\
&\hspace{5cm}\times\cac^{(n)}_{t_2,s_1}(y_1,w_1) \cac^{(n)}_{t_2,s_2}(y_2,w_2)\cac^{(n)}_{t_2,s_2}(y_1,w_2) \cac^{(n)}_{t_2,s_1}(y_2,w_1).
\end{align*}
In this way,
\begin{align}
\int dy_1 \, \big|\cm^{\mathbf{2,2},(n),\star}_{t_1,t_2}(y_1)\big|&\lesssim \int dy_1 \, \big|\cm^{\mathbf{2,2,1},(n),\star}_{t_1,t_2}(y_1)\big| +\int dy_1 \,  \big|\cm^{\mathbf{2,2,2},(n),\star}_{t_1,t_2}(y_1)\big| .\label{contr-2-ii-bis-2-2}
\end{align}

\

\

\noindent
\textit{(i) \underline{Study of $\cm^{\mathbf{2,2,1},(n)}$.}} Arguing as in the derivation of \eqref{ref-2-1-pap}, we obtain
\begin{multline*}
 \int dy_1 \sup_{y_2} \big|\cm^{\mathbf{2,2,1},(n)}_{t_1,t_2}(y_1+y_2,y_2)\big|\lesssim \\
\lesssim \big\|\cac^{(n),\star}_{t_2,t_2}\big\|_{L^2} \int_{0}^{t_1}ds_1\int_{0}^{t_1}ds_2 \, \big\| K^\star_{t_1-s_1,t_2-s_1}\big\|_{L^1} \big\|\cac^{(n),\star}_{t_2,s_1} \big\|_{L^\infty}\big\|\cac^{(n),\star}_{s_1,s_2}\big\|_{L^2} \big\|K^\star_{t_1-s_2,t_2-s_2}\big\|_{L^1} \big\|\cac^{(n),\star}_{t_2,s_2}\big\|_{L^\infty},
\end{multline*}
and thanks to Lemma \ref{lem:k-star-1} and Lemma \ref{lem:cac3}, we deduce
\begin{align}
\int dy_1 \sup_{y_2} \big|\cm^{\mathbf{2,2,1},(n)}_{t_1,t_2}(y_1+y_2,y_2)\big|&\lesssim |t_2-t_1|^{\frac13}\bigg(\int_{0}^{t_1} \frac{ds}{|t_1-s|^{\frac12+\frac{5}{12}}}e^{-(t_1-s)}\bigg)^2\nonumber\\
&\lesssim |t_2-t_1|^{\frac13}\bigg(\int_{0}^{\infty} \frac{ds}{|s|^{\frac{11}{12}}}e^{-s}\bigg)^2\nonumber\\
& \lesssim   |t_2-t_1|^{\frac13},\label{contr-2-ii-bis-2-2-1}
\end{align}
uniformly over $0\leq t_1,t_2\leq 2$.

\

\

\noindent
\textit{(ii) \underline{Study of $\cm^{\mathbf{2,2,2},(n)}$.}} Repeating the argument that led to \eqref{ref-2-2-pap}, we find
\begin{align*}
\int dy_1 \sup_{y_2} \big|\cm^{\mathbf{2,2,2},(n)}_{t_1,t_2}(y_1+y_2,y_2)\big|
&\lesssim \bigg(\int_{0}^{t_1}ds\,  \big\| K^\star_{t_1-s,t_2-s}\big\|_{L^1} \big\|\cac^{(n),\star}_{t_2,s}\big\|_{L^\infty} \big\|\cac^{(n),\star}_{t_2,s}\big\|_{L^2} \bigg)^2
\end{align*}
and using the results in Lemma \ref{lem:k-star-1} and Lemma \ref{lem:cac3}, we deduce
\begin{align}
\int dy_1 \sup_{y_2} \big|\cm^{\mathbf{2,2,2},(n)}_{t_1,t_2}(y_1+y_2,y_2)\big|&\lesssim |t_2-t_1|^{\frac13}\bigg(\int_{0}^{t_1} \frac{ds}{|t_1-s|^{\frac12+\frac{5}{12}}}e^{-(t_1-s)}\bigg)^2\nonumber\\
&\lesssim |t_2-t_1|^{\frac13}\bigg(\int_{0}^{\infty} \frac{ds}{|s|^{\frac{11}{12}}}e^{-s}\bigg)^2\nonumber\\
& \lesssim   |t_2-t_1|^{\frac13},\label{contr-2-ii-bis-2-2-2}
\end{align}
uniformly over $0\leq t_1,t_2\leq 2$.

\

By inserting the two estimates \eqref{contr-2-ii-bis-2-2-1} and \eqref{contr-2-ii-bis-2-2-2} into \eqref{contr-2-ii-bis-2-2}, we can conclude that
\begin{align}\label{cm22}
&\int dy_1 \, \big|\cm^{\mathbf{2,2},(n),\star}_{t_1,t_2}(y_1)\big|\lesssim |t_2-t_1|^{\frac14},
\end{align}
uniformly over $0\leq t_1,t_2\leq 2$.

\

\

\subsubsection{Study of $\cm^{\mathbf{2,3},(n)}$.} One has here
\begin{align*}
&\cm^{\mathbf{2,3},(n)}_{t_1,t_2}(y_1,y_2)=\mathbb{E}\Big[\scret^{\mathbf{2,3},(n)}_{t_1,t_2}(y_1)\, \scret^{\mathbf{2,3},(n)}_{t_1,t_2}(y_2) \Big]\\
&=\int_{0}^{t_1}ds_1\int_{0}^{t_1}ds_2 \int dw_1 dw_2 \, K_{t_1-s_1}(y_1,w_1)K_{t_1-s_2}(y_2,w_2) \cac^{(n)}_{(t_1,t_2),s_1}(y_1,w_1) \cac^{(n)}_{(t_1,t_2),s_2}(y_2,w_2)\\
&\hspace{2cm}\mathbb{E}\Big[I^W_2\big( F^{(n)}_{t_2,y_1} \otimes F^{(n)}_{s_1,w_1} \big)I^W_2\big( F^{(n)}_{t_2,y_2} \otimes F^{(n)}_{s_2,w_2} \big)\Big]\\
&=\int_{0}^{t_1}ds_1\int_{0}^{t_1}ds_2 \int dw_1 dw_2 \, K_{t_1-s_1}(y_1,w_1)K_{t_1-s_2}(y_2,w_2) \cac^{(n)}_{(t_1,t_2),s_1}(y_1,w_1) \cac^{(n)}_{(t_1,t_2),s_2}(y_2,w_2)\\
&\hspace{2cm}\Big[c_{\mathbf{1}}\, \cac^{(n)}_{t_2,t_2}(y_1,y_2) \cac^{(n)}_{s_1,s_2}(w_1,w_2) +c_{\mathbf{2}}\, \cac^{(n)}_{t_2,s_2}(y_1,w_2) \cac^{(n)}_{t_2,s_1}(y_2,w_1) \Big],
\end{align*}
for some combinatorial coefficients $c_{\mathbf{1}},c_{\mathbf{2}}\geq 0$. Thus,
\begin{align*}
\cm^{\mathbf{2,3},(n)}_{t_1,t_2}(y_1,y_2)= c_{\mathbf{1}}\cm^{\mathbf{2,3,1},(n)}_{t_1,t_2}(y_1,y_2)+ c_{\mathbf{2}}\cm^{\mathbf{2,3,2},(n)}_{t_1,t_2}(y_1,y_2),
\end{align*}
with
\begin{align*}
&\cm^{\mathbf{2,3,1},(n)}_{t_1,t_2}(y_1,y_2):=\int_{0}^{t_1}ds_1\int_{0}^{t_1}ds_2\int dw_1 dw_2 \,  K_{t_1-s_1}(y_1,w_1)K_{t_1-s_2}(y_2,w_2)\\
&\hspace{5cm}\times \cac^{(n)}_{(t_1,t_2),s_1}(y_1,w_1) \cac^{(n)}_{(t_1,t_2),s_2}(y_2,w_2)\cac^{(n)}_{t_2,t_2}(y_1,y_2) \cac^{(n)}_{s_1,s_2}(w_1,w_2),
\end{align*}
and 
\begin{align*}
&\cm^{\mathbf{2,3,2},(n)}_{t_1,t_2}(y_1,y_2):=\int_{0}^{t_1}ds_1\int_{0}^{t_1}ds_2\int dw_1 dw_2 \,  K_{t_1-s_1}(y_1,w_1)K_{t_1-s_2}(y_2,w_2)\\
&\hspace{5cm}\times\cac^{(n)}_{(t_1,t_2),s_1}(y_1,w_1) \cac^{(n)}_{(t_1,t_2),s_2}(y_2,w_2)\cac^{(n)}_{t_2,s_2}(y_1,w_2) \cac^{(n)}_{t_2,s_1}(y_2,w_1).
\end{align*}
In this way,
\begin{align}
\int dy_1 \, \big|\cm^{\mathbf{2,3},(n),\star}_{t_1,t_2}(y_1)\big|&\lesssim \int dy_1 \, \big|\cm^{\mathbf{2,3,1},(n),\star}_{t_1,t_2}(y_1)\big| +\int dy_1 \,  \big|\cm^{\mathbf{2,3,2},(n),\star}_{t_1,t_2}(y_1)\big| .\label{contr-2-ii-bis-2-3}
\end{align}

\

\

\noindent
\textit{(i) \underline{Study of $\cm^{\mathbf{2,3,1},(n)}$.}} By the same arguments as for \eqref{ref-2-1-pap}, we obtain
\begin{multline*}
 \int dy_1 \sup_{y_2} \big|\cm^{\mathbf{2,3,1},(n)}_{t_1,t_2}(y_1+y_2,y_2)\big| \lesssim \\
\lesssim \big\|\cac^{(n),\star}_{t_2,t_2}\big\|_{L^2} \int_{0}^{t_1}ds_1\int_{0}^{t_1}ds_2 \, \big\| K^\star_{t_1-s_1}\big\|_{L^1} \big\|\cac^{(n),\star}_{(t_1,t_2),s_1} \big\|_{L^\infty}\big\|\cac^{(n),\star}_{s_1,s_2}\big\|_{L^2} \big\|K^\star_{t_1-s_2}\big\|_{L^1} \big\|\cac^{(n),\star}_{(t_1,t_2),s_2}\big\|_{L^\infty},
\end{multline*}
and thanks to Lemma \ref{lem:k-star-0} and Lemma \ref{lem:cac3}, we deduce
\begin{align}
\int dy_1 \sup_{y_2} \big|\cm^{\mathbf{2,3,1},(n)}_{t_1,t_2}(y_1+y_2,y_2)\big|&\lesssim |t_2-t_1|^{\frac12}\bigg(\int_{0}^{t_1} \frac{ds}{|t_1-s|^{\frac34}}e^{-(t_1-s)}\bigg)^2\nonumber\\
&\lesssim |t_2-t_1|^{\frac12}\bigg(\int_{0}^{\infty} \frac{ds}{|s|^{\frac34}}e^{-s}\bigg)^2\nonumber\\
& \lesssim   |t_2-t_1|^{\frac12},\label{contr-2-ii-bis-2-3-1}
\end{align}
uniformly over $0\leq t_1,t_2\leq 2$.

\

\

\noindent
\textit{(ii) \underline{Study of $\cm^{\mathbf{2,3,2},(n)}$.}} With the same arguments as those leading to \eqref{ref-2-2-pap}, we obtain
\begin{align*}
\int dy_1 \sup_{y_2} \big|\cm^{\mathbf{2,3,2},(n)}_{t_1,t_2}(y_1+y_2,y_2)\big|
&\lesssim \bigg(\int_{0}^{t_1}ds\,  \big\| K^\star_{t_1-s}\big\|_{L^1} \big\|\cac^{(n),\star}_{(t_1,t_2),s}\big\|_{L^\infty} \big\|\cac^{(n),\star}_{t_2,s}\big\|_{L^2} \bigg)^2
\end{align*}
and using the results in Lemma \ref{lem:k-star-0} and Lemma \ref{lem:cac3}, we deduce
\begin{align}
\int dy_1 \sup_{y_2} \big|\cm^{\mathbf{2,3,2},(n)}_{t_1,t_2}(y_1+y_2,y_2)\big|&\lesssim |t_2-t_1|^{\frac12}\bigg(\int_{0}^{t_1} \frac{ds}{|t_1-s|^{\frac34}}e^{-(t_1-s)}\bigg)^2\nonumber\\
&\lesssim |t_2-t_1|^{\frac12}\bigg(\int_{0}^{\infty} \frac{ds}{|s|^{\frac34}}e^{-s}\bigg)^2\nonumber\\
& \lesssim   |t_2-t_1|^{\frac12},\label{contr-2-ii-bis-2-3-2}
\end{align}
uniformly over $0\leq t_1,t_2\leq 2$.

\

By inserting the two estimates \eqref{contr-2-ii-bis-2-3-1} and \eqref{contr-2-ii-bis-2-3-2} into \eqref{contr-2-ii-bis-2-3}, we can conclude that
\begin{align}\label{cm23}
&\int dy_1 \, \big|\cm^{\mathbf{2,3},(n),\star}_{t_1,t_2}(y_1)\big|\lesssim |t_2-t_1|^{\frac14},
\end{align}
uniformly over $0\leq t_1,t_2\leq 2$.

\

\

\subsubsection{Study of $\cm^{\mathbf{2,4},(n)}$.} One has here
\begin{multline*}
\cm^{\mathbf{2,4},(n)}_{t_1,t_2}(y_1,y_2)=\mathbb{E}\Big[\scret^{\mathbf{2,4},(n)}_{t_1,t_2}(y_1)\, \scret^{\mathbf{2,4},(n)}_{t_1,t_2}(y_2) \Big]=\\
\begin{aligned}
&=\int_{0}^{t_1}ds_1\int_{0}^{t_1}ds_2 \int dw_1 dw_2 \, K_{t_1-s_1}(y_1,w_1)K_{t_1-s_2}(y_2,w_2) \cac^{(n)}_{t_1,s_1}(y_1,w_1) \cac^{(n)}_{t_1,s_2}(y_2,w_2)\\
&\hspace{2cm}\times \mathbb{E}\Big[I^W_2\big( (F^{(n)}_{t_2,y_1}-F^{(n)}_{t_1,y_1}) \otimes F^{(n)}_{s_1,w_1} \big)I^W_2\big( (F^{(n)}_{t_2,y_2}-F^{(n)}_{t_1,y_2}) \otimes F^{(n)}_{s_2,w_2} \big)\Big]\\
&=\int_{0}^{t_1}ds_1\int_{0}^{t_1}ds_2 \int dw_1 dw_2 \, K_{t_1-s_1}(y_1,w_1)K_{t_1-s_2}(y_2,w_2) \cac^{(n)}_{t_1,s_1}(y_1,w_1) \cac^{(n)}_{t_1,s_2}(y_2,w_2)\\
&\hspace{1cm}\times \Big[c_{\mathbf{1}}\, \cac^{(n)}_{(t_1,t_2),(t_1,t_2)}(y_1,y_2) \cac^{(n)}_{s_1,s_2}(w_1,w_2) +c_{\mathbf{2}}\, \cac^{(n)}_{(t_1,t_2),s_2}(y_1,w_2) \cac^{(n)}_{(t_1,t_2),s_1}(y_2,w_1) \Big],
\end{aligned}
\end{multline*}
for some combinatorial coefficients $c_{\mathbf{1}},c_{\mathbf{2}}\geq 0$. Thus,
\begin{align*}
\cm^{\mathbf{2,4},(n)}_{t_1,t_2}(y_1,y_2)= c_{\mathbf{1}}\cm^{\mathbf{2,4,1},(n)}_{t_1,t_2}(y_1,y_2)+ c_{\mathbf{2}}\cm^{\mathbf{2,4,2},(n)}_{t_1,t_2}(y_1,y_2),
\end{align*}
with
\begin{multline*}
\cm^{\mathbf{2,4,1},(n)}_{t_1,t_2}(y_1,y_2):=\int_{0}^{t_1}ds_1\int_{0}^{t_1}ds_2\int dw_1 dw_2 \,  K_{t_1-s_1}(y_1,w_1)K_{t_1-s_2}(y_2,w_2)\\
\hspace{3cm}\times \cac^{(n)}_{t_1,s_1}(y_1,w_1) \cac^{(n)}_{t_1,s_2}(y_2,w_2)\cac^{(n)}_{(t_1,t_2),(t_1,t_2)}(y_1,y_2) \cac^{(n)}_{s_1,s_2}(w_1,w_2),
\end{multline*}
and 
\begin{multline*}
\cm^{\mathbf{2,4,2},(n)}_{t_1,t_2}(y_1,y_2):=\int_{0}^{t_1}ds_1\int_{0}^{t_1}ds_2\int dw_1 dw_2 \,  K_{t_1-s_1}(y_1,w_1)K_{t_1-s_2}(y_2,w_2)\\
\hspace{3cm}\times\cac^{(n)}_{t_1,s_1}(y_1,w_1) \cac^{(n)}_{t_1,s_2}(y_2,w_2)\cac^{(n)}_{(t_1,t_2),s_2}(y_1,w_2) \cac^{(n)}_{(t_1,t_2),s_1}(y_2,w_1).
\end{multline*}
In this way,
\begin{align}
\int dy_1 \, \big|\cm^{\mathbf{2,4},(n),\star}_{t_1,t_2}(y_1)\big|&\lesssim \int dy_1 \, \big|\cm^{\mathbf{2,4,1},(n),\star}_{t_1,t_2}(y_1)\big| +\int dy_1 \,  \big|\cm^{\mathbf{2,4,2},(n),\star}_{t_1,t_2}(y_1)\big| .\label{contr-2-ii-bis-2-4}
\end{align}

\

\

\noindent
\textit{(i) \underline{Study of $\cm^{\mathbf{2,4,1},(n)}$.}} The same calculation as for \eqref{ref-2-1-pap} shows that
\begin{multline*}
 \int dy_1 \sup_{y_2} \big|\cm^{\mathbf{2,4,1},(n)}_{t_1,t_2}(y_1+y_2,y_2)\big|\lesssim \\
\lesssim \big\|\cac^{(n),\star}_{(t_1,t_2),(t_1,t_2)}\big\|_{L^2} \int_{0}^{t_1}ds_1\int_{0}^{t_1}ds_2 \, \big\| K^\star_{t_1-s_1}\big\|_{L^1} \big\|\cac^{(n),\star}_{t_1,s_1} \big\|_{L^\infty}\big\|\cac^{(n),\star}_{s_1,s_2}\big\|_{L^2} \big\|K^\star_{t_1-s_2}\big\|_{L^1} \big\|\cac^{(n),\star}_{t_1,s_2}\big\|_{L^\infty},
\end{multline*}
and thanks to Lemma \ref{lem:k-star-0} and Lemma \ref{lem:cac3}, we deduce
\begin{align}
\int dy_1 \sup_{y_2} \big|\cm^{\mathbf{2,4,1},(n)}_{t_1,t_2}(y_1+y_2,y_2)\big|&\lesssim |t_2-t_1|^{\frac14}\bigg(\int_{0}^{t_1} \frac{ds}{|t_1-s|^{\frac12}}e^{-(t_1-s)}\bigg)^2\nonumber\\
&\lesssim |t_2-t_1|^{\frac14}\bigg(\int_{0}^{\infty} \frac{ds}{|s|^{\frac12}}e^{-s}\bigg)^2\nonumber\\
& \lesssim   |t_2-t_1|^{\frac14},\label{contr-2-ii-bis-2-4-1}
\end{align}
uniformly over $0\leq t_1,t_2\leq 2$.

\

\

\noindent
\textit{(ii) \underline{Study of $\cm^{\mathbf{2,4,2},(n)}$.}} Paralleling the derivation of \eqref{ref-2-2-pap}, we arrive at
\begin{align*}
\int dy_1 \sup_{y_2} \big|\cm^{\mathbf{2,4,2},(n)}_{t_1,t_2}(y_1+y_2,y_2)\big|
&\lesssim \bigg(\int_{0}^{t_1}ds\,  \big\| K^\star_{t_1-s}\big\|_{L^1} \big\|\cac^{(n),\star}_{t_1,s}\big\|_{L^\infty} \big\|\cac^{(n),\star}_{(t_1,t_2),s}\big\|_{L^2} \bigg)^2
\end{align*}
and using the results in Lemma \ref{lem:k-star-0} and Lemma \ref{lem:cac3}, we obtain that
\begin{align}
\int dy_1 \sup_{y_2} \big|\cm^{\mathbf{2,4,2},(n)}_{t_1,t_2}(y_1+y_2,y_2)\big|&\lesssim |t_2-t_1|^{\frac12}\bigg(\int_{0}^{t_1} \frac{ds}{|t_1-s|^{\frac12}}e^{-(t_1-s)}\bigg)^2\nonumber\\
&\lesssim |t_2-t_1|^{\frac12}\bigg(\int_{0}^{\infty} \frac{ds}{|s|^{\frac12}}e^{-s}\bigg)^2\nonumber\\
& \lesssim   |t_2-t_1|^{\frac12},\label{contr-2-ii-bis-2-4-2}
\end{align}
uniformly over $0\leq t_1,t_2\leq 2$.

\

By inserting the two estimates \eqref{contr-2-ii-bis-2-4-1} and \eqref{contr-2-ii-bis-2-4-2} into \eqref{contr-2-ii-bis-2-4}, we can conclude that
\begin{align}\label{cm24}
&\int dy_1 \, \big|\cm^{\mathbf{2,4},(n),\star}_{t_1,t_2}(y_1)\big|\lesssim |t_2-t_1|^{\frac14},
\end{align}
uniformly over $0\leq t_1,t_2\leq 2$.

\

We can finally substitute \eqref{cm21}, \eqref{cm22}, \eqref{cm23} and \eqref{cm24} into \eqref{contr-1-bis-2}, which yields
\begin{equation}\label{scretdeux}
\int dx \, \mathbb{E}\bigg[ \Big| H^{-\al}\big(\scret^{\mathbf{2},(n)}_{t_1,t_2}\big)(x)\Big|^{2}\bigg]^p\lesssim |t_2-t_1|^{\frac{p}{4}},
\end{equation}
uniformly over $0\leq t_1,t_2\leq 2$.

\

By combining \eqref{scretun} and \eqref{scretdeux}, we deduce the desired bound \eqref{hyperc-40-1} (recall the decomposition in \eqref{decompo-scretun-scretdeux}).
 

\

\

\section{About the fifth-order diagram} \label{sec:fifth-o-i}

The fifth-order diagram was introduced in \eqref{ord5} as
$$\<Psi2IPsi3>^{(n)}_{s,t}(x):=\Big(\<Psi2>^{(n)}_{s,t} \pe \<IPsi3>^{(n)}_{s,t}\Big)(x) - 3\, \frakc^{\mathbf{2},(n)}_{s,t}(x) \, \<Psi>_{s,t}^{(n)}{(x)} ,$$
where the deterministic sequence $\frakc^{\mathbf{2},(n)}$ is given by 
$$\frakc^{\mathbf{2},(n)}_{s,t}(x):=\mathbb{E}\Big[ \<Psi2>^{(n)}_{s,t}(x) \<IPsi2>^{(n)}_{s,t}(x)\Big].$$

\smallskip

A first construction and regularity result for this object was obtained in \cite[Proposition~10.1]{DFT}. The following extension is an immediate consequence of Lemma \ref{lem:stati-z}:

\begin{proposition}
Fix $0\leq T_1<T_2$ and for all $t\in [T_1,T_2]$, set
\begin{equation*}
\widetilde{\<Psi2IPsi3>}^{(n)}_{T_1,t}:=\int_{T_1}^t \<Psi2IPsi3>^{(n)}_{T_1,s} \, ds.
\end{equation*}
Then for all $0<\eps,\eta<\frac12$, there exists $\ka>0$ such that
\begin{equation}\label{contro-5-1} 
\sup_{\ell\geq 0}    \mathbb{E} \Big[ \Big\|\widetilde{ \<Psi2IPsi3>}^{(n+1)}_{\ell,.}   - \widetilde{ \<Psi2IPsi3>}^{(n)}_{\ell,.}\Big\|_{{\ov \cac}^{\frac34-\eps}([\ell,\ell+2];  \cb_\infty^{-\eta})}^{2p} \Big]\lesssim 2^{-\ka n p }. 
\end{equation}
As a particular consequence, for all $0\leq T_1<T_2$, the sequence $(\widetilde{\<Psi2IPsi3>}^{(n)}_{T_1,.})_{n\geq 1}$ converges almost surely to an element $\widetilde{\<Psi2IPsi3>}_{T_1,.}$ in $\cac^{\frac34-\eps}([T_1,T_2]; \cb_{\infty}^{-\eta})$, for all $\varepsilon,\eta>0$. 

\end{proposition}

We now aim to construct $\<Psi2IPsi3>_{\ell,.}$ as a continuous function of time. Our final result in this direction reads as follows:

\begin{theorem}\label{theo-10.1}
Let $T>0$. For all $0<\eta<\frac12$, there exist $\eps>0$ and $\ka>0$ such that
\begin{eqnarray*}
&\sup_{\ell\geq 0}\mathbb{E} \Big[ \Big\|{ \<Psi2IPsi3>}^{(n+1)}_{\ell,.} -{ \<Psi2IPsi3>}^{(n)}_{\ell,.}\Big\|_{\cac^{\eps}([\ell,\ell+2]; \cb_{\infty}^{-\frac56-\eta})}^{2p} \Big]\lesssim 2^{-\ka n p} .
\end{eqnarray*}
Consequently, for all $0\leq T_1<T_2$, the sequence $( {\<Psi2IPsi3>}^{(n)}_{T_1,.})$ converges almost surely to an element $ {\<Psi2IPsi3>_{T_1,.}}$ in  $\cac^{\eps}([T_1,T_2]; \cb_{\infty}^{-\frac56-\eta})$, for all $\eta>0$ and $\varepsilon>0$ small enough.

\end{theorem}

Our main technical result towards Theorem \ref{theo-10.1} relies on the study of the complete fifth-order diagram  $\<Psi2IPsi3nr> $, that is,
 \begin{equation*}
 \<Psi2IPsi3nr>^{(n)}_{s,t}(x):=\Big(\<Psi2>^{(n)}_{s,t}  \<IPsi3>^{(n)}_{s,t}\Big)(x) - 3\, \frakc^{\mathbf{2},(n)}_{s,t}(x) \, \<Psi>_{s,t}^{(n)}{(x)}.
 \end{equation*}
 Then the following statement holds: 
 
 \begin{proposition}\label{prop-10.1}
Let $T>0$. For all $0<\eta<\frac12$, there exists $\ka>0$ such that
\begin{equation}\label{boun-10.1}
\sup_{\ell\geq 0}\mathbb{E} \bigg[ \Big\|{ \<Psi2IPsi3nr>}^{(n+1)}_{\ell,.} -{ \<Psi2IPsi3nr>}^{(n)}_{\ell,.}\Big\|_{\cac^{\frac15}([\ell,\ell+2]; \cb_{\infty}^{-\frac32-\eta})}^{2p} \bigg]\lesssim 2^{-\ka n p} .
\end{equation}

\end{proposition}
 
We now show that the result of Proposition \ref{prop-10.1} implies that of Theorem \ref{theo-10.1}.

 \begin{proof}[Proof of Theorem \ref{theo-10.1}]

Observe first that 
\begin{equation}\label{decomp-1}
\<Psi2IPsi3nr>^{(n)}- \<Psi2IPsi3>^{(n)}= \<Psi2>^{(n)} \pl\<IPsi3>^{(n)}+\<IPsi3>^{(n)} \pl \<Psi2>^{(n)}.
\end{equation}
Therefore, using the paraproduct rules together with the regularity results contained in \cite[Proposition 6.1]{DFT} and  \cite[Proposition 7.1]{DFT}, it is readily checked that the estimate \eqref{boun-10.1} also applies to the resonant part $\<Psi2IPsi3>^{(n)}$. In other words, one has, for all $0<\eta<\frac12$
\begin{equation}\label{contro-5-2} 
\sup_{\ell\geq 0}\mathbb{E} \bigg[ \Big\|  \<Psi2IPsi3>^{(n+1)}_{\ell,.}-\<Psi2IPsi3>^{(n)}_{\ell,.}\Big\|_{\cac^{\frac15}([\ell,\ell+2]; \cb_{\infty}^{-\frac32-\eta})}^{2p} \bigg]\lesssim 2^{-\ka n p} .
\end{equation}
We can then implement an interpolation procedure based on Lemma \ref{lemma-interpt}. Namely, we fix $0<\eta<\frac12$ and pick $\eps>0$ small enough so that
$$\eta':=\eta+\frac{\frac38+\frac56\eps}{\frac14+2\eps}-\frac32>0 \quad \text{and} \quad \theta:=\frac{\frac14+2\eps}{\frac{9}{20}+\eps}\in (0,1).$$
By applying Lemma \ref{lemma-interpt}, we obtain that for all $\ell\geq 0$,
\begin{multline*} 
\big\|  \<Psi2IPsi3>^{(n+1)}_{\ell,.}-  \<Psi2IPsi3>^{(n)}_{\ell,.}  \big\|_{\cac^{\eps}([\ell,\ell+2]; \cb_{\infty}^{-\frac56-\eta})} \lesssim \\
\lesssim \big\| \widetilde{ \<Psi2IPsi3>}_{\ell,.}^{(n+1)}-\widetilde{ \<Psi2IPsi3>}_{\ell,.}^{(n)} \big\|^{1-\theta}_{{\ov \cac}^{\frac34-\eps}([\ell,\ell+2];  \cb_{\infty}^{-\eta})} \big\|  \<Psi2IPsi3>^{(n+1)}_{\ell,.}-\<Psi2IPsi3>^{(n)}_{\ell,.}\big\|_{\cac^{\frac15}([0,T]; \cb_{\infty}^{-\frac32-\eta'})}^\theta.
\end{multline*}
The conclusion now follows from the two controls \eqref{contro-5-1} and \eqref{contro-5-2}.
 \end{proof}

\

 The remaining part of this section is thus devoted to the proof of Proposition \ref{prop-10.1}. For the sake of conciseness, we will restrict our attention to the uniform bound
\begin{equation*} 
\sup_{n\geq 1}\sup_{\ell\geq 0}\mathbb{E} \bigg[ \Big\|{ \<Psi2IPsi3nr>}^{(n)}_{\ell,.}\Big\|_{\cac^{\frac15}([\ell,\ell+2]; \cb_{\infty}^{-\frac32-\eps})}^{2p} \bigg]<\infty .
\end{equation*}
Thanks to the stationarity property stated in Lemma \ref{lem:stati-z}, the above objective reduces to proving that
\begin{equation*} 
\sup_{n\geq 1}\mathbb{E} \bigg[ \Big\|{ \<Psi2IPsi3nr>}^{(n)}_{0,.}\Big\|_{\cac^{\frac15}([0,2]; \cb_{\infty}^{-\frac32-\eps})}^{2p} \bigg]<\infty .
\end{equation*}

\

To start with, recall that, as in the previous sections, one has, for every $p\geq 1$ large enough,
\begin{equation*}
\mathbb{E}\Big[\big\|\<Psi2IPsi3nr>^{(n)}_{0,.} \big\|_{\cac^{\frac15}([0,2];\cb_{\infty}^{-\frac32-\eps})}^{2p}\Big] \lesssim \int_{[0,2]^2} \frac{dt_1 dt_2 }{|t_2-t_1|^{\frac25 p+2}}\int dx \, \mathbb{E}\bigg[ \Big| H^{-\al}\big(\<Psi2IPsi3nr>^{(n)}_{0,t_2}-\<Psi2IPsi3nr>^{(n)}_{0,t_1}\big)(x)\Big|^{2}\bigg]^p.
\end{equation*}
As a result, it remains for us to establish that
\begin{equation}\label{hold-coucou}
\sup_{n\geq 1}\int dx \, \mathbb{E}\bigg[ \Big| H^{-\al}\big(\<Psi2IPsi3nr>^{(n)}_{0,t_2}-\<Psi2IPsi3nr>^{(n)}_{0,t_1}\big)(x)\Big|^{2}\bigg]^p \lesssim |t_2-t_1|^{\theta p},
\end{equation}
uniformly over $0\leq t_1,t_2\leq 2$, for some $\theta>\frac25 $.

\

\

\subsection{Proof of \eqref{hold-coucou}}

Observe first that $\frakc^{\mathbf{2},(n)}(y)$ can be recast as
\begin{align}
\frakc^{\mathbf{2},(n)}_{0,t}(y)&= \int_{0}^t ds \int dw \, K_{t-s}(y,w)\mathbb{E}\Big[ \<Psi2>^{(n)}_{0,t}(y) \<Psi2>^{(n)}_{0,s}(w)\Big] =2\int_{0}^t ds \int dw \, K_{t-s}(y,w)\cac^{(n)}_{t,s}(y,w)^2. \label{identif-c2}
\end{align}

\smallskip

Moreover, one has
\begin{eqnarray*}
\Big(\<Psi2>^{(n)}_{0,t}  \<IPsi3>^{(n)}_{0,t}\Big)(y)
&=& \<Psi2>^{(n)}_{0,t}(y) \int_{0}^t ds \int dw \, K_{t-s}(y,w) \<Psi3>^{(n)}_{0,s}(w) \\
&=& \int_{0}^t ds \int dw \, K_{t-s}(y,w) I^W_2\big(F^{(n)}_{t,y}\otimes F^{(n)}_{t,y}\big)I^W_3\big(F^{(n)}_{s,w}\otimes F^{(n)}_{s,w}\otimes F^{(n)}_{s,w}\big).
\end{eqnarray*}
Using the multiplication rule in \cite[Lemma 4.1]{DFT} and the identity \eqref{identif-c2}, we obtain the decomposition
\begin{equation}\label{decomp-arbre-ordre-5}
\<Psi2IPsi3nr>^{(n)}_{0,t}(y)=\calt^{\mathbf{1},(n)}_t(y)+6\, \calt^{\mathbf{2},(n)}_t(y)+6\, \calt^{\mathbf{3},(n)}_t(y)+6\, \calt^{\mathbf{4},(n)}_t(y),
\end{equation}
with
\begin{align*}
\calt^{\mathbf{1},(n)}_t(y):= \int_{0}^t ds \int dw \, K_{t-s}(y,w) I^W_5\big(F^{(n)}_{t,y}\otimes F^{(n)}_{t,y}\otimes F^{(n)}_{s,w}\otimes F^{(n)}_{s,w}\otimes F^{(n)}_{s,w}\big),
\end{align*}
\begin{align*}
\calt^{\mathbf{2},(n)}_t(y):= \int_{0}^t ds \int dw \, K_{t-s}(y,w) \cac^{(n)}_{t,s}(y,w)I^W_3\big(F^{(n)}_{t,y}\otimes F^{(n)}_{s,w}\otimes F^{(n)}_{s,w}\big),
\end{align*}
\begin{align*}
\calt^{\mathbf{3},(n)}_t(y)&:=\int_{0}^t ds \int dw \, K_{t-s}(y,w)\cac^{(n)}_{t,s}(y,w)^2 \big(\<Psi>^{(n)}_{0,s}(w)-\<Psi>^{(n)}_{0,t}(w)\big),
\end{align*}
\begin{align*}
\calt^{\mathbf{4},(n)}_t(y)&:= \int_{0}^t ds \int dw \, K_{t-s}(y,w)\cac^{(n)}_{t,s}(y,w)^2 \big(\<Psi>^{(n)}_{0,t}(w)-\<Psi>^{(n)}_{0,t}(y)\big).
\end{align*}

\smallskip

With expression \eqref{decomp-arbre-ordre-5} in mind, we will decompose the increments of each diagram $\calt^{\mathbf{a},(n)}$ as a sum (see {\it e.g.} \eqref{exemple-decompos})
$$\calt^{\mathbf{a},(n)}_{t_2}-\calt^{\mathbf{a},(n)}_{t_1} =\sum_{\mathbf{b}=1,\ldots, 4}\calt^{\mathbf{a,b},(n)}_{t_1,t_2},$$
with corresponding covariances
$$\cm^{\mathbf{a,b},(n)}_{t_1,t_2}(y_1,y_2):=\mathbb{E}\Big[ \calt^{\mathbf{a,b},(n)}_{t_1,t_2} \calt^{\mathbf{a,b},(n)}_{t_1,t_2} \Big].$$
In turn, following the Wiener chaoses computation rules, each quantity $\cm^{\mathbf{a,b},(n)}$ will be further decomposed as a finite sum $\cm^{\mathbf{a,b},(n)}:= \sum_{\mathbf{d}}c_{\mathbf{d}} \cm^{\mathbf{a,b,d},(n)}$ (see {\it e.g.}~\eqref{exemple-decompos-1}). With this notation in hand, Lemma \ref{lem:techn} yields, for all $\al>\frac34$ and $p\geq 1$,
\begin{align}
&\int dx \, \mathbb{E}\bigg[ \Big| H^{-\al}\big(\<Psi2IPsi3nr>^{(n)}_{0,t_2}-\<Psi2IPsi3nr>^{(n)}_{0,t_1}\big)(x)\Big|^{2}\bigg]^p
\lesssim \sum_{\mathbf{a}, \mathbf{b},\mathbf{d}} \bigg(\int dy_1 \,  \Big(\sup_{y_2}\big|\cm^{\mathbf{a,b,d},(n)}_{t_1,t_2}(y_1+y_2,y_2)\big| \Big)\bigg)^p.\label{transition-coucou}
\end{align}
In brief, to reach the H{\"o}lder-type control \eqref{hold-coucou}, we will need to show that, for all $\mathbf{a},\mathbf{b},\mathbf{d}$, one has
\begin{equation}\label{hold-m-a-b-c}
\sup_{n\geq 1}\sup_{t\geq 0}\int dy_1 \,  \Big(\sup_{y_2}\big|\cm^{\mathbf{a,b,d},(n)}_{t_1,t_2}(y_1+y_2,y_2)\big| \Big) \lesssim |t_2-t_1|^{\theta}
\end{equation}
uniformly over $0\leq t_1,t_2\leq 2$, for some $\theta>\frac25 $.

\

\subsubsection{Increments: first diagram}

We decompose the increments of $\calt^{\mathbf{1},(n)}$ as
\begin{align}
&\calt^{\mathbf{1},(n)}_{t_2}(y)-\calt^{\mathbf{1},(n)}_{t_1}(y)=\calt^{\mathbf{1,1},(n)}_{t_1,t_2}(y)+\calt^{\mathbf{1,2},(n)}_{t_1,t_2}(y)+\calt^{\mathbf{1,3},(n)}_{t_1,t_2}(y)+\calt^{\mathbf{1,4},(n)}_{t_1,t_2}(y),\label{exemple-decompos}
\end{align}
with
\begin{align*}
&\calt^{\mathbf{1,1},(n)}_{t_1,t_2}(y):= \int_{t_1}^{t_2} ds \int dw \, K_{t_2-s}(y,w) I^W_5\Big(F^{(n)}_{t_2,y}\otimes F^{(n)}_{t_2,y}\otimes F^{(n)}_{s,w}\otimes F^{(n)}_{s,w}\otimes F^{(n)}_{s,w}\Big),
\end{align*}

\begin{align*}
&\calt^{\mathbf{1,2},(n)}_{t_1,t_2}(y):= \int_{0}^{t_1} ds \int dw \, K_{t_1-s,t_2-s}(y,w) I^W_5\Big(F^{(n)}_{t_2,y}\otimes F^{(n)}_{t_2,y}\otimes F^{(n)}_{s,w}\otimes F^{(n)}_{s,w}\otimes F^{(n)}_{s,w}\Big),
\end{align*}

\begin{align*}
&\calt^{\mathbf{1,3},(n)}_{t_1,t_2}(y):= \int_{0}^{t_1} ds \int dw \, K_{t_1-s}(y,w) I^W_5\Big(\big(F^{(n)}_{t_2,y}-F^{(n)}_{t_1,y}\big)\otimes F^{(n)}_{t_2,y}\otimes F^{(n)}_{s,w}\otimes F^{(n)}_{s,w}\otimes F^{(n)}_{s,w}\Big),
\end{align*}

\begin{align*}
&\calt^{\mathbf{1,4},(n)}_{t_1,t_2}(y):= \int_{0}^{t_1} ds \int dw \, K_{t_1-s}(y,w) I^W_5\Big(F^{(n)}_{t_1,y} \otimes \big(F^{(n)}_{t_2,y}-F^{(n)}_{t_1,y}\big) \otimes F^{(n)}_{s,w}\otimes F^{(n)}_{s,w}\otimes F^{(n)}_{s,w}\Big).
\end{align*}

\

\noindent
\underline{\textit{Case of $\calt^{\mathbf{1,1},(n)}$}}.
One has
\begin{align}
\cm^{\mathbf{1,1},(n)}_{t_1,t_2}(y_1,y_2):=\mathbb{E}\Big[\calt^{\mathbf{1,1},(n)}_{t_1,t_2}(y_1)\calt^{\mathbf{1,1},(n)}_{t_1,t_2}(y_2)\Big]=\sum_{\mathbf{a}=1}^3 c_{\mathbf{a}}\cm^{\mathbf{1,1,a},(n)}_{t_1,t_2}(y_1,y_2),\label{exemple-decompos-1}
\end{align}
with
\begin{multline*}
\cm^{\mathbf{1,1,1},(n)}_{t_1,t_2}(y_1,y_2):=\int_{t_1}^{t_2} ds_1 \int dw_1 \, K_{t_2-s_1}(y_1,w_1)\\
\hspace{2cm}\times \int_{t_1}^{t_2} ds_2 \int dw_2 \, K_{t_2-s_2}(y_2,w_2)\cac^{(n)}_{t_2,t_2}(y_1,y_2)^2 \cac^{(n)}_{s_1,s_2}(w_1,w_2)^3,
\end{multline*}

\begin{multline*}
\cm^{\mathbf{1,1,2},(n)}_{t_1,t_2}(y_1,y_2):=\int_{t_1}^{t_2} ds_1 \int dw_1 \, K_{t_2-s_1}(y_1,w_1)\\
\times \int_{t_1}^{t_2} ds_2 \int dw_2 \, K_{t_2-s_2}(y_2,w_2)\cac^{(n)}_{t_2,t_2}(y_1,y_2)\cac^{(n)}_{t_2,s_2}(y_1,w_2)\cac^{(n)}_{t_2,s_1}(y_2,w_1) \cac^{(n)}_{s_1,s_2}(w_1,w_2)^2,
\end{multline*}

\begin{multline*}
\cm^{\mathbf{1,1,3},(n)}_{t_1,t_2}(y_1,y_2):=\int_{t_1}^{t_2} ds_1 \int dw_1 \, K_{t_2-s_1}(y_1,w_1)\int_{t_1}^{t_2} ds_2\\
\times  \int dw_2 \, K_{t_2-s_2}(y_2,w_2)\cac^{(n)}_{t_2,s_2}(y_1,w_2)^2\cac^{(n)}_{t_2,s_1}(y_2,w_1)^2 \cac^{(n)}_{s_1,s_2}(w_1,w_2).
\end{multline*}

\

\paragraph{Case of $\cm^{\mathbf{1,1,1},(n)}$}

Observe that
\begin{align*}
&\big|\cm^{\mathbf{1,1,1},(n)}_{t_1,t_2}(y_1+y_2,y_2)\big| \lesssim \\
&\lesssim \cac^{(n)}_{t_2,t_2}(y_1+y_2,y_2)^2 \int_{t_1}^{t_2} ds_1 \int_{t_1}^{t_2} ds_2\int dw_1 dw_2 \, K_{t_2-s_1}(y_1+y_2,w_1)  K_{t_2-s_2}(y_2,w_2)\cac^{(n)}_{s_1,s_2}(w_1,w_2)^3\\
&\lesssim \cac^{(n),\star}_{t_2,t_2}(y_1)^2 \int_{t_1}^{t_2} ds_1 \int_{t_1}^{t_2} ds_2\int dw_1 dw_2  \, K_{t_2-s_1}(y_1+y_2,w_1+y_2) K_{t_2-s_2}(y_2,w_2+y_2)\cac^{(n)}_{s_1,s_2}(w_1+y_2,w_2+y_2)^3\\
&\lesssim \cac^{(n),\star}_{t_2,t_2}(y_1)^2 \int_{t_1}^{t_2} ds_1 \int_{t_1}^{t_2} ds_2\int dw_1 dw_2  \, K^\star_{t_2-s_1}(y_1-w_1)K^\star_{t_2-s_2}(w_2)\cac^{(n),\star}_{s_1,s_2}(w_1-w_2)^3\\
&\lesssim \cac^{(n),\star}_{t_2,t_2}(y_1)^2 \int_{t_1}^{t_2} ds_1 \int_{t_1}^{t_2} ds_2\, \Big[K^\star_{t_2-s_1}\ast\big[K^\star_{t_2-s_2}\ast (\cac^{(n),\star}_{s_1,s_2})^3\big]\Big](y_1).
\end{align*}

Observe that we have
\begin{align*}
 \Big\|K^\star_{t_2-s_1}\ast\big[K^\star_{t_2-s_2}\ast (\cac^{(n),\star}_{s_1,s_2})^3\big]\Big\|_{L^\infty} 
&\lesssim  \big\|K^\star_{t_2-s_1}\big\|_{L^1} \big\| K^\star_{t_2-s_2}\ast (\cac^{(n),\star}_{s_1,s_2})^3\big\|_{L^\infty}\\
&\lesssim   \big\|K^\star_{t_2-s_1}\big\|_{L^1} \big\| K^\star_{t_2-s_2}\big\|_{L^2} \big\|(\cac^{(n),\star}_{s_1,s_2})^3\big\|_{L^2} \\
&\lesssim  \big\|K^\star_{t_2-s_1}\big\|_{L^1} \big\| K^\star_{t_2-s_2}\big\|_{L^2} \big\|\cac^{(n),\star}_{s_1,s_2}\big\|_{L^6}^3,
\end{align*}
then
\begin{align}
 \int dy_1 \sup_{y_2} \big|\cm^{\mathbf{1,1,1},(n)}_{t_1,t_2}(y_1+y_2,y_2)\big|
&\lesssim \big\|(\cac^{(n),\star}_{t_2,t_2})^2\big\|_{L^1} \int_{t_1}^{t_2} ds_1 \int_{t_1}^{t_2} ds_2\, \Big\|K^\star_{t_2-s_1}\ast\big[K^\star_{t_2-s_2}\ast (\cac^{(n),\star}_{s_1,s_2})^3\big]\Big\|_{L^\infty}\nonumber\\
&\lesssim \big\|\cac^{(n),\star}_{t_2,t_2}\big\|_{L^2}^2 \int_{t_1}^{t_2} ds_1 \int_{t_1}^{t_2} ds_2 \, \big\|K^\star_{t_2-s_1}\big\|_{L^1} \big\| K^\star_{t_2-s_2}\big\|_{L^2} \big\|\cac^{(n),\star}_{s_1,s_2}\big\|_{L^6}^3.\label{ref-1-1-papier}
\end{align}

Now, thanks to Lemma \ref{lem:k-star-0} and Lemma \ref{lem:cac3}, one has
\begin{align*}
\int dy_1 \sup_{y_2} \big|\cm^{\mathbf{1,1,1},(n)}_{t_1,t_2}(y_1+y_2,y_2)\big|
&\lesssim \int_{t_1}^{t_2} ds_1 \int_{t_1}^{t_2} ds_2 \, e^{-(t_2-s_1)}\frac{e^{- (t_2-s_2)}}{|t_2-s_2|^{\frac34}}   \frac{1}{|s_2-s_1|^{\frac34}}\\
&\lesssim \int_{t_1}^{t_2} ds_1 \int_{t_1}^{t_2} ds_2 \, \frac{1}{|t_2-s_2|^{\frac34}}   \frac{1}{|s_2-s_1|^{\frac34}} ,
\end{align*}
and as result
\begin{align*}
&\int dy_1 \sup_{y_2} \big|\cm^{\mathbf{1,1,1},(n)}_{t_1,t_2}(y_1+y_2,y_2)\big|\lesssim   |t_2-t_1|^{\frac12},
\end{align*}
uniformly over $0\leq t_1,t_2\leq 2$.

\

\paragraph{Case of $\cm^{\mathbf{1,1,2},(n)}$}

One has here
\begin{align*}
&\cm^{\mathbf{1,1,2},(n)}_{t_1,t_2}(y_1+y_2,y_2)=\cac^{(n)}_{t_2,t_2}(y_1+y_2,y_2)\int_{t_1}^{t_2} ds_1 \int dw_1 \, K_{t_2-s_1}(y_1+y_2,w_1)\cac^{(n)}_{t_2,s_1}(y_2,w_1) \\
&\hspace{1.5cm}\times \int_{t_1}^{t_2}ds_2 \int dw_2 \, K_{t_2-s_2}(y_2,w_2)\cac^{(n)}_{t_2,s_2}(y_1+y_2,w_2)\cac^{(n)}_{s_1,s_2}(w_1,w_2)^2\\
&=\cac^{(n)}_{t_2,t_2}(y_1+y_2,y_2)\int_{t_1}^{t_2} ds_1 \int dw_1 \, K_{t_2-s_1}(y_1+y_2,w_1+y_2)\cac^{(n)}_{t_2,s_1}(y_2,w_1+y_2)\\
&\hspace{0.5cm}\times \int_{t_1}^{t_2} ds_2 \int dw_2 \, K_{t_2-s_2}(y_2,w_2+y_2)\cac^{(n)}_{t_2,s_2}(y_1+y_2,w_2+y_2) \cac^{(n)}_{s_1,s_2}(w_1+y_2,w_2+y_2)^2.
\end{align*}
This implies that 
\begin{align*}
&\cm^{\mathbf{1,1,2},(n)}_{t_1,t_2}(y_1+y_2,y_2)\lesssim \cac^{(n),\star}_{t_2,t_2}(y_1)\int_{t_1}^{t_2} ds_1 \int_{t_1}^{t_2} ds_2  \int dw_1 \, K^\star_{t_2-s_1}(y_1-w_1)\cac^{(n),\star}_{t_2,s_1}(w_1)\\
& \hspace{3cm}       \times   \int dw_2 \, K^\star_{t_2-s_2}(w_2)\cac^{(n),\star}_{t_2,s_2}(y_1-w_2) \cac^{(n),\star}_{s_1,s_2}(w_1-w_2)^2\\
&\lesssim \cac^{(n),\star}_{t_2,t_2}(y_1)\int_{t_1}^{t_2}ds_1 \int_{t_1}^{t_2} ds_2\, \big\|\cac^{(n),\star}_{t_2,s_2}\big\|_{L^\infty}\int dw_1 \, K^\star_{t_2-s_1}(y_1-w_1)\cac^{(n),\star}_{t_2,s_1}(w_1)\\
& \hspace{3cm}     \times \int dw_2 \, K^\star_{t_2-s_2}(w_2)\cac^{(n),\star}_{s_1,s_2}(w_1-w_2)^2\\
&\lesssim \cac^{(n),\star}_{t_2,t_2}(y_1)\int_{t_1}^{t_2} ds_1 \int_{t_1}^{t_2} ds_2\, \big\|\cac^{(n),\star}_{t_2,s_2}\big\|_{L^\infty} \bigg[ K^\star_{t_2-s_1}\ast\Big[\cac^{(n),\star}_{t_2,s_1}\cdot  \big[ K^\star_{t_2-s_2}\ast (\cac^{(n),\star}_{s_1,s_2})^2\big]\Big]\bigg](y_1).
\end{align*}
Next, we can bound
 \begin{align*}
  \bigg\| K^\star_{t_2-s_1}\ast\Big[\cac^{(n),\star}_{t_2,s_1}\cdot  \big[ K^\star_{t_2-s_2}\ast (\cac^{(n),\star}_{s_1,s_2})^2\big]\Big]\bigg\|_{L^2} 
&\lesssim    \big\| K^\star_{t_2-s_1}\big\|_{L^2} \Big\|\cac^{(n),\star}_{t_2,s_1}\cdot  \big[ K^\star_{t_2-s_2}\ast (\cac^{(n),\star}_{s_1,s_2})^2\big]\Big\|_{L^1} \\
&\lesssim  \big\| K^\star_{t_2-s_1}\big\|_{L^2} \big\|\cac^{(n),\star}_{t_2,s_1}\big\|_{L^2}  \big\| K^\star_{t_2-s_2}\ast (\cac^{(n),\star}_{s_1,s_2})^2\big\|_{L^2} \\
&\lesssim   \big\| K^\star_{t_2-s_1}\big\|_{L^2} \big\|\cac^{(n),\star}_{t_2,s_1}\big\|_{L^2}  \big\| K^\star_{t_2-s_2}\big\|_{L^1} \big\|\cac^{(n),\star}_{s_1,s_2}\big\|_{L^4}^2,
\end{align*}
 and therefore
\begin{multline}
\int dy_1 \sup_{y_2} \big|\cm^{\mathbf{1,1,2},(n)}_{t_1,t_2}(y_1+y_2,y_2)\big| \lesssim \\
\begin{aligned}
&\lesssim \big\|\cac^{(n),\star}_{t_2,t_2}\big\|_{L^2}\int_{t_1}^{t_2} ds_1 \int_{t_1}^{t_2} ds_2\, \big\|\cac^{(n),\star}_{t_2,s_2}\big\|_{L^\infty}  \bigg\| K^\star_{t_2-s_1}\ast\Big[\cac^{(n),\star}_{t_2,s_1}\cdot  \big[ K^\star_{t_2-s_2}\ast (\cac^{(n),\star}_{s_1,s_2})^2\big]\Big]\bigg\|_{L^2} \\
&\lesssim \big\|\cac^{(n),\star}_{t_2,t_2}\big\|_{L^2}\int_{t_1}^{t_2} ds_1 \int_{t_1}^{t_2} ds_2\, \big\|\cac^{(n),\star}_{t_2,s_2}\big\|_{L^\infty}  \big\| K^\star_{t_2-s_1}\big\|_{L^2} \big\|\cac^{(n),\star}_{t_2,s_1}\big\|_{L^2}  \big\| K^\star_{t_2-s_2}\big\|_{L^1} \big\|\cac^{(n),\star}_{s_1,s_2}\big\|_{L^4}^2.  \label{ref-1-2-papier}
\end{aligned}
\end{multline}

As a result, using Lemma \ref{lem:k-star-0} and Lemma \ref{lem:cac3},
\begin{align*}
\int dy_1 \sup_{y_2} \big|\cm^{\mathbf{1,1,2},(n)}_{t_1,t_2}(y_1+y_2,y_2)\big|&\lesssim \int_{t_1}^{t_2} \frac{ds_1}{|t_2-s_1|^{\frac34}}e^{- (t_2-s_1)} \int_{t_1}^{t_2} \frac{ds_2}{|t_2-s_2|^{\frac12}}e^{- (t_2-s_2)}\frac{1}{|s_2-s_1|^{\frac14}}\\
&\lesssim \int_{t_1}^{t_2} \frac{ds_1}{|t_2-s_1|^{\frac34}} \int_{t_1}^{t_2} \frac{ds_2}{|t_2-s_2|^{\frac12}}\frac{1}{|s_2-s_1|^{\frac14}},
\end{align*}
and from here we deduce that
\begin{align*}
&\int dy_1 \sup_{y_2} \big|\cm^{\mathbf{1,1,2},(n)}_{t_1,t_2}(y_1+y_2,y_2)\big|\lesssim |t_2-t_1|^{\frac12},
\end{align*}
uniformly over $0\leq t_1,t_2\leq 2$.

\

\paragraph{Case of $\cm^{\mathbf{1,1,3},(n)}$}

With the same arguments as above,
\begin{align*}
&\cm^{\mathbf{1,1,3},(n)}_{t_1,t_2}(y_1+y_2,y_2)\lesssim \\
&\lesssim \int_{t_1}^{t_2} ds_1 \int_{t_1}^{t_2}  ds_2\int dw_1 \, K^\star_{t_2-s_1}(y_1-w_1)\cac^{(n),\star}_{t_2,s_1}(w_1)^2\int dw_2 \, K^\star_{t-s_2}(w_2)\cac^{(n),\star}_{t_2,s_2}(y_1-w_2)^2 \cac^{(n),\star}_{s_1,s_2}(w_1-w_2)\\
&\lesssim \int_{t_1}^{t_2} ds_1 \int_{t_1}^{t_2}  ds_2\, \big\|\cac^{(n),\star}_{s_1,s_2}\big\|_{L^\infty} \big[K^\star_{t_2-s_1}\ast (\cac^{(n),\star}_{t_2,s_1})^2\big](y_1) \big[K^\star_{t_2-s_2}\ast (\cac^{(n),\star}_{t_2,s_2})^2\big](y_1) 
\end{align*}
and thus
\begin{align}
&\int dy_1 \sup_{y_2} \big|\cm^{\mathbf{1,1,3},(n)}_{t_1,t_2}(y_1+y_2,y_2)\big|\nonumber \\
&\lesssim \int_{t_1}^{t_2} ds_1 \int_{t_1}^{t_2}  ds_2\, \big\|\cac^{(n),\star}_{s_1,s_2}\big\|_{L^\infty} \int dy_1 \,  \big[K^\star_{t_2-s_1}\ast (\cac^{(n),\star}_{t_2,s_1})^2\big](y_1) \big[K^\star_{t_2-s_2}\ast (\cac^{(n),\star}_{t_2,s_2})^2\big](y_1)\nonumber\\
&\lesssim \int_{t_1}^{t_2} ds_1 \int_{t_1}^{t_2}  ds_2\, \big\|\cac^{(n),\star}_{s_1,s_2}\big\|_{L^\infty} \big\|K^\star_{t_2-s_1}\ast (\cac^{(n),\star}_{t_2,s_1})^2\big\|_{L^2} \big\|K^\star_{t_2-s_2}\ast (\cac^{(n),\star}_{t_2,s_2})^2\big\|_{L^2}\nonumber\\
&\lesssim \int_{t_1}^{t_2} ds_1 \int_{t_1}^{t_2}  ds_2\, \big\|\cac^{(n),\star}_{s_1,s_2}\big\|_{L^\infty}  \big\|K^\star_{t_2-s_1}\big\|_{L^{\frac65}}\big\|(\cac^{(n),\star}_{t_2,s_1})^2\big\|_{L^{\frac32}} \big\|K^\star_{t_2-s_2}\big\|_{L^{\frac65}}\big\|(\cac^{(n),\star}_{t_2,s_2})^2\big\|_{L^{\frac32}}.  \label{ref-1-3-papier}
\end{align}
Then, using Lemma \ref{lem:k-star-0} and Lemma \ref{lem:cac3},
\begin{align*}
\int dy_1 \sup_{y_2} \big|\cm^{\mathbf{1,1,3},(n)}_{t_1,t_2}(y_1+y_2,y_2)\big|
&\lesssim \int_{t_1}^{t_2} \frac{ds_1}{|t_2-s_1|^{\frac14+\eps}}e^{-(t_2-s_1)} \int_{t_1}^{t_2} \frac{ds_2}{|t_2-s_2|^{\frac14+\eps}} e^{-(t_2-s_2)} \frac{1}{|s_2-s_1|^{\frac12}}\\
&\lesssim \int_{t_1}^{t_2} \frac{ds_1}{|t_2-s_1|^{\frac14+\eps}} \int_{t_1}^{t_2} \frac{ds_2}{|t_2-s_2|^{\frac14+\eps}}  \frac{1}{|s_2-s_1|^{\frac12}},
\end{align*}
and in this way
\begin{align*}
&\int dy_1 \sup_{y_2} \big|\cm^{\mathbf{1,1,3},(n)}_{t_1,t_2}(y_1+y_2,y_2)\big|\lesssim |t_2-t_1|^{1-2\eps},
\end{align*}
uniformly over $0\leq t_1,t_2\leq 2$.

\

\noindent
\underline{\textit{Case of $\calt^{\mathbf{1,2},(n)}$}}.
One has
\begin{align*}
\cm^{\mathbf{1,2},(n)}_{t_1,t_2}(y_1,y_2):=\mathbb{E}\Big[\calt^{\mathbf{1,2},(n)}_{t_1,t_2}(y_1)\calt^{\mathbf{1,2},(n)}_{t_1,t_2}(y_2)\Big]=\sum_{\mathbf{a}=1}^3 c_{\mathbf{a}}\cm^{\mathbf{1,2,a},(n)}_{t_1,t_2}(y_1,y_2)
\end{align*}
with
\begin{multline*}
\cm^{\mathbf{1,2,1},(n)}_{t_1,t_2}(y_1,y_2):=\int_{0}^{t_1} ds_1 \int dw_1 \,K_{t_1-s_1,t_2-s_1}(y_1,w_1)\\
\hspace{2cm}\times \int_{0}^{t_1} ds_2 \int dw_2 \, K_{t_1-s_2,t_2-s_2}(y_2,w_2)\cac^{(n)}_{t_2,t_2}(y_1,y_2)^2 \cac^{(n)}_{s_1,s_2}(w_1,w_2)^3,
\end{multline*}

\begin{multline*}
\cm^{\mathbf{1,2,2},(n)}_{t_1,t_2}(y_1,y_2):=\int_{0}^{t_1} ds_1 \int dw_1 \, K_{t_1-s_1,t_2-s_1}(y_1,w_1)\\
\times\int_{0}^{t_1} ds_2 \int dw_2 \, K_{t_1-s_2,t_2-s_2}(y_2,w_2)\cac^{(n)}_{t_2,t_2}(y_1,y_2)\cac^{(n)}_{t_2,s_2}(y_1,w_2)\cac^{(n)}_{t_2,s_1}(y_2,w_1) \cac^{(n)}_{s_1,s_2}(w_1,w_2)^2,
\end{multline*}

\begin{multline*}
\cm^{\mathbf{1,2,3},(n)}_{t_1,t_2}(y_1,y_2):=\int_{0}^{t_1} ds_1 \int dw_1 \, \\
\times\int_{0}^{t_1} ds_2 \int dw_2 \, K_{t_1-s_2,t_2-s_2}(y_2,w_2)\cac^{(n)}_{t_2,s_2}(y_1,w_2)^2\cac^{(n)}_{t_2,s_1}(y_2,w_1)^2 \cac^{(n)}_{s_1,s_2}(w_1,w_2).
\end{multline*}

\

\paragraph{Case of $\cm^{\mathbf{1,2,1},(n)}$}
One has in this situation
\begin{multline*}
\big|\cm^{\mathbf{1,2,1},(n)}_{t_1,t_2}(y_1+y_2,y_2)\big| \lesssim\\
\lesssim \cac^{(n),\star}_{t_2,t_2}(y_1)^2\int_{0}^{t_1} ds_1 \int dw_1 \,  K^\star_{t_1-s_1,t_2-s_1}(y_1-w_1)\int_{0}^{t_1} ds_2 \int dw_2 \, K^\star_{t_1-s_2,t_2-s_2}(w_2) \cac^{(n),\star}_{s_1,s_2}(w_1-w_2)^3,
\end{multline*}
and so
\begin{multline*}
\int dy_1\sup_{y_2} \big|\cm^{\mathbf{1,2,1},(n)}_{t_1,t_2}(y_1+y_2,y_2)\big|\lesssim \\
\begin{aligned}
&\lesssim \int_{0}^{t_1} ds_1 \int_{0}^{t_1} ds_2\int dw_1 \, \bigg(\int dy_1\, K^\star_{t_1-s_1,t_2-s_1}(y_1-w_1)\cac^{(n),\star}_{t_2,t_2}(y_1)^2\bigg) \\
&\hspace{6cm}\times \bigg(\int dw_2 \, K^\star_{t_1-s_2,t_2-s_2}(w_2) \cac^{(n),\star}_{s_1,s_2}(w_1-w_2)^3\bigg)\\
&\lesssim \int_{0}^{t_1} ds_1\int_{0}^{t_1} ds_2  \,  \int dw_1 \, \big[ K^\star_{t_1-s_1,t_2-s_1}\ast(\cac^{(n),\star}_{t_2,t_2})^2\big](w_1)\big[ K^\star_{t_1-s_2,t_2-s_2}\ast (\cac^{(n),\star}_{s_1,s_2})^3 \big](w_1)\\
&\lesssim \int_{0}^{t_1} ds_1\int_{0}^{t_1} ds_2  \,  \big\| K^\star_{t_1-s_1,t_2-s_1}\ast(\cac^{(n),\star}_{t_2,t_2})^2\big\|_{L^2}\big\| K^\star_{t_1-s_2,t_2-s_2}\ast (\cac^{(n),\star}_{s_1,s_2})^3 \big\|_{L^2}\\
&\lesssim \big\| (\cac^{(n),\star}_{t_2,t_2})^2\big\|_{L^{\frac{6}{4+\eps}}}\int_{0}^{t_1} ds_1\int_{0}^{t_1} ds_2  \,   \big\| K^\star_{t_1-s_1,t_2-s_1}\big\|_{L^{\frac{6}{5-\eps}}}\big\| K^\star_{t_1-s_2,t_2-s_2} \big\|_{L^2}\big\| (\cac^{(n),\star}_{s_1,s_2})^3 \big\|_{L^1}.
\end{aligned}
\end{multline*}
Using both Lemma \ref{lem:k-star-1} and Lemma \ref{lem:cac3}, we deduce that
\begin{multline*}
\int dy_1\sup_{y_2} \big|\cm^{\mathbf{1,2,1},(n)}_{t_1,t_2}(y_1+y_2,y_2)\big|\lesssim \\
\begin{aligned}
&\lesssim |t_2-t_1|^{\frac13+\frac17-\eps}\int_{0}^{t_1} ds_1\int_{0}^{t_1} ds_2  \,   \frac{e^{-(t_1-s_1)}}{|t_1-s_1|^{1-\frac{\eps}{6}}}\frac{e^{-(t_1-s_2)}}{|t_1-s_2|^{1-\frac{7\eps}{4}}}\frac{1}{|s_1-s_2|^\eps}\\
&\lesssim |t_2-t_1|^{\frac13+\frac17-\eps}\int_{0}^{+\infty} ds_1\int_{0}^{+\infty} ds_2  \,   \frac{e^{-s_1}}{|s_1|^{1-\frac{\eps}{6}}}\frac{e^{-s_2}}{|s_2|^{1-\frac{7\eps}{4}}}\frac{1}{|s_1-s_2|^\eps}\lesssim |t_2-t_1|^{\frac{10}{21}-\eps},
\end{aligned}
\end{multline*}
uniformly over $0\leq t_1,t_2\leq 2$.

\

\paragraph{Case of $\cm^{\mathbf{1,2,2},(n)}$}

One has in this case
\begin{multline*}
\big|\cm^{\mathbf{1,2,2},(n)}_{t_1,t_2}(y_1+y_2,y_2)\big|\lesssim \\
\begin{aligned}
&\lesssim \cac^{(n),\star}_{t_2,t_2}(y_1)\int_{0}^{t_1} ds_1\int_{0}^{t_1} ds_2 \,  \int dw_1 \, K_{t_1-s_1,t_2-s_1}^\star(y_1-w_1)\cac^{(n),\star}_{t_2,s_1}(w_1)\\
&\hspace{5cm}\int dw_2 \,  K_{t_1-s_2,t_2-s_2}^\star(w_2)\cac^{(n),\star}_{t_2,s_2}(y_1-w_2) \cac^{(n),\star}_{s_1,s_2}(w_1-w_2)^2\\
&\lesssim \cac^{(n),\star}_{t_2,t_2}(y_1)\int_{0}^{t_1} ds_1\, \big\|\cac^{(n),\star}_{t_2,s_1}\big\|_{L^\infty}\int_{0}^{t_1}ds_2\, \big\|\cac^{(n),\star}_{t_2,s_2}\big\|_{L^\infty}  \int dw_1 \, K_{t_1-s_1,t_2-s_1}^\star(y_1-w_1)\\
&\hspace{5cm}\int dw_2 \,  K_{t_1-s_2,t_2-s_2}^\star(w_2) \cac^{(n),\star}_{s_1,s_2}(w_1-w_2)^2
\end{aligned}
\end{multline*}
and therefore
\begin{align*}
&\int dy_1 \sup_{y_2}\big|\cm^{\mathbf{1,2,2},(n)}_{t_1,t_2}(y_1+y_2,y_2)\big|\lesssim \int_{0}^{t_1} ds_1\, \big\|\cac^{(n),\star}_{t_2,s_1}\big\|_{L^\infty}\int_{0}^{t_1}ds_2\, \big\|\cac^{(n),\star}_{t_2,s_2}\big\|_{L^\infty}  \\
&\hspace{1cm} \times \int dw_1 \bigg[ \int dy_1\, \cac^{(n),\star}_{t_2,t_2}(y_1) K_{t_1-s_1,t_2-s_1}^\star(y_1-w_1)\bigg]\bigg[\int dw_2 \,  K_{t_1-s_2,t_2-s_2}^\star(w_2) \cac^{(n),\star}_{s_1,s_2}(w_1-w_2)^2\bigg]\\
&\lesssim \int_{0}^{t_1} ds_1\, \big\|\cac^{(n),\star}_{t_2,s_1}\big\|_{L^\infty}\int_{0}^{t_1}ds_2\, \big\|\cac^{(n),\star}_{t_2,s_2}\big\|_{L^\infty}   \big\| \cac^{(n),\star}_{t_2,t_2}\ast K_{t_1-s_1,t_2-s_1}^\star\big\|_{L^{\frac{3}{1+\eps}}} \big\|  K_{t_1-s_2,t_2-s_2}^\star\ast (\cac^{(n),\star}_{s_1,s_2})^2\big\|_{L^{\frac{3}{2-\eps}}}\\
&\lesssim \big\| \cac^{(n),\star}_{t_2,t_2}\big\|_{L^{\frac{3}{1+\eps}}}\int_{0}^{t_1} ds_1\, \big\|\cac^{(n),\star}_{t_2,s_1}\big\|_{L^\infty}\int_{0}^{t_1}ds_2\, \big\|\cac^{(n),\star}_{t_2,s_2}\big\|_{L^\infty}   \big\| K_{t_1-s_1,t_2-s_1}^\star\big\|_{L^1} \big\|  K_{t_1-s_2,t_2-s_2}^\star\big\|_{L^1}\big\|(\cac^{(n),\star}_{s_1,s_2})^2\big\|_{L^{\frac{3}{2-\eps}}}.
\end{align*}
As a result, using Lemma \ref{lem:k-star-1} and Lemma \ref{lem:cac3}, one has for every $\theta\in [0,1]$,
\begin{multline*}
\int dy_1 \sup_{y_2}\big|\cm^{\mathbf{1,2,2},(n)}_{t_1,t_2}(y_1+y_2,y_2)\big|\lesssim \\
\begin{aligned}
&\lesssim |t_2-t_1|^{2\theta}\int_{0}^{t_1} \frac{ds_1}{|t_2-s_1|^{\frac12}}\int_{0}^{t_1} \frac{ds_2}{|t_2-s_2|^{\frac12}}  \frac{|t_2-s_1|^{\frac32\theta}}{|t_1-s_1|^{\frac52\theta}}e^{-(t_1-s_1)}\frac{|t_2-s_2|^{\frac32\theta}}{|t_1-s_2|^{\frac52\theta}}\frac{e^{-(t_1-s_2)}}{|s_2-s_1|^{\frac{\eps}{2}}}\\
&\lesssim |t_2-t_1|^{2\theta}\int_{0}^{t_1} \frac{ds_1\, e^{-(t_1-s_1)}}{|t_2-s_1|^{\frac12-\frac32\theta} |t_1-s_1|^{\frac52 \theta}}\int_{0}^{t_1} \frac{ds_2\, e^{-(t_1-s_2)}}{|t_2-s_2|^{\frac12-\frac32\theta}|t_1-s_2|^{\frac52\theta}} \frac{1}{|s_2-s_1|^{\frac{\eps}{2}}}.
\end{aligned}
\end{multline*}
By choosing $\theta=\frac13$, we obtain, for $\eps>0$ small enough,
\begin{align*}
\int dy_1 \sup_{y_2}\big|\cm^{\mathbf{1,2,2},(n)}_{t_1,t_2}(y_1+y_2,y_2)\big|
&\lesssim |t_2-t_1|^{\frac23}\int_{0}^{t_1} \frac{ds_1 \, e^{-(t_1-s_1)} }{|t_1-s_1|^{\frac56}}\int_{0}^{t_1} \frac{ds_2 \, e^{-(t_1-s_2)}}{|t_1-s_2|^{\frac56}|s_2-s_1|^{\frac{\eps}{2}}}   \\
&\lesssim |t_2-t_1|^{\frac23}\int_{0}^{+\infty} \frac{ds_1\, e^{-s_1}}{|s_1|^{\frac56}}\int_{0}^{+\infty} \frac{ds_2\, e^{-s_2} }{|s_2|^{\frac56}|s_2-s_1|^{\frac{\eps}{2}}}  \lesssim |t_2-t_1|^{\frac23},
\end{align*}
uniformly over $0\leq t_1,t_2\leq 2$.

\

\paragraph{Case of $\cm^{\mathbf{1,2,3},(n)}$}

In this situation, with the same arguments as those leading to \eqref{ref-1-3-papier}, we obtain that
\begin{multline*}
\int dy_1 \sup_{y_2} \big|\cm^{\mathbf{1,2,3},(n)}_{t_1,t_2}(y_1+y_2,y_2)\big|\lesssim\\
\begin{aligned}
&\lesssim \int_{0}^{t_1}ds_1\int_{0}^{t_1} ds_2 \,\big\|\cac^{(n),\star}_{s_1,s_2}\big\|_{L^\infty} \big\| K_{t_1-s_1,t_2-s_1}^\star\big\|_{L^{\frac65}}\big\| (\cac^{(n),\star}_{t_2,s_1})^2\big\|_{L^{\frac32}}\big\| K_{t_1-s_2,t_2-s_2}^\star\big\|_{L^{\frac65}} \big\|  (\cac^{(n),\star}_{t_2,s_2})^2\big\|_{L^{\frac32}} .
\end{aligned}
\end{multline*}
Thanks to Lemma \ref{lem:k-star-1} and Lemma \ref{lem:cac3}, we obtain that for every $\theta\in [0,1]$,
\begin{multline*}
\int dy_1 \sup_{y_2} \big|\cm^{\mathbf{1,2,3},(n)}_{t_1,t_2}(y_1+y_2,y_2)\big|\lesssim\\
\begin{aligned}
&\lesssim |t_2-t_1|^{2\theta}\int_{0}^{t_1}ds_1\int_{0}^{t_1} ds_2 \, \frac{1}{|s_2-s_1|^{\frac12}}  \frac{e^{-(t_1-s_1)}}{|t_1-s_1|^{\frac14+\frac{9}{4}\theta}}\frac{1}{|t_2-s_1|^\eps}\frac{e^{-(t_1-s_2)}}{|t_1-s_2|^{\frac14+\frac{9}{4}\theta}} \frac{1}{|t_2-s_2|^\eps}\\
&\lesssim |t_2-t_1|^{2\theta}\int_{0}^{t_1}ds_1\int_{0}^{t_1} ds_2 \, \frac{1}{|s_2-s_1|^{\frac12}}  \frac{e^{-(t_1-s_1)}}{|t_1-s_1|^{\frac14+\frac{9}{4}\theta+\eps}}\frac{e^{-(t_1-s_2)}}{|t_1-s_2|^{\frac14+\frac{9}{4}\theta+\eps}}.
\end{aligned}
\end{multline*}
By choosing $\theta=\frac{2}{9}-\eps$, we obtain
\begin{multline*}
\int dy_1 \sup_{y_2} \big|\cm^{\mathbf{1,2,3},(n)}_{t_1,t_2}(y_1+y_2,y_2)\big|\lesssim\\
\begin{aligned}
&\lesssim |t_2-t_1|^{\frac{4}{9}-2\eps}\int_{0}^{t_1}ds_1\int_{0}^{t_1} ds_2 \, \frac{1}{|s_2-s_1|^{\frac12}}  \frac{e^{-(t_1-s_1)}}{|t_1-s_1|^{\frac34-\frac{5}{4}\eps}}\frac{e^{-(t_1-s_2)}}{|t_1-s_2|^{\frac34-\frac{5}{4}\eps}}\\
&\lesssim |t_2-t_1|^{\frac{4}{9}-2\eps}\int_{0}^{+\infty}ds_1\int_{0}^{+\infty} ds_2 \, \frac{1}{|s_2-s_1|^{\frac12}}  \frac{e^{-s_1}}{|s_1|^{\frac34-\frac{5}{4}\eps}}\frac{e^{-s_2}}{|s_2|^{\frac34-\frac{5}{4}\eps}}\lesssim |t_2-t_1|^{\frac{4}{9}-2\eps},
\end{aligned}
\end{multline*}
uniformly over $0\leq t_1,t_2\leq 2$.

\

\

\noindent
\underline{\textit{Case of $\calt^{\mathbf{1,3},(n)}$}}. One has
\begin{align*}
\cm^{\mathbf{1,3},(n)}_{t_1,t_2}(y_1,y_2):=\mathbb{E}\Big[\calt^{\mathbf{1,3},(n)}_{t_1,t_2}(y_1)\calt^{\mathbf{1,3},(n)}_{t_1,t_2}(y_2)\Big]=\sum_{\mathbf{a}=1}^{7} c_{\mathbf{a}}\cm^{\mathbf{1,3,a},(n)}_{t_1,t_2}(y_1,y_2),
\end{align*}
with
\begin{align*}
&\cm^{\mathbf{1,3,1},(n)}_{t_1,t_2}(y_1,y_2):=\int_{0}^{t_1} ds_1 \int dw_1 \, K_{t_2-s_1}(y_1,w_1)\int_{0}^{t_1} ds_2 \int dw_2 \, K_{t_2-s_2}(y_2,w_2)\\
&\hspace{4cm}\times \cac^{(n)}_{(t_1,t_2),(t_1,t_2)}(y_1,y_2) \cac^{(n)}_{t_2,t_2}(y_1,y_2)\cac^{(n)}_{s_1,s_2}(w_1,w_2)^3,
\end{align*}

\begin{align*}
&\cm^{\mathbf{1,3,2},(n)}_{t_1,t_2}(y_1,y_2):=\int_{0}^{t_1} ds_1 \int dw_1 \, K_{t_2-s_1}(y_1,w_1)\int_{0}^{t_1} ds_2 \int dw_2 \, K_{t_2-s_2}(y_2,w_2)\\
&\hspace{4cm}\times\cac^{(n)}_{(t_1,t_2),(t_1,t_2)}(y_1,y_2) \cac^{(n)}_{t_2,s_2}(y_1,w_2)\cac^{(n)}_{t_2,s_1}(y_2,w_1)\cac^{(n)}_{s_1,s_2}(w_1,w_2)^2,
\end{align*}

\begin{align*}
&\cm^{\mathbf{1,3,3},(n)}_{t_1,t_2}(y_1,y_2):=\int_{0}^{t_1} ds_1 \int dw_1 \, K_{t_2-s_1}(y_1,w_1)\int_{0}^{t_1} ds_2 \int dw_2 \, K_{t_2-s_2}(y_2,w_2)\\
&\hspace{4cm}\times\cac^{(n)}_{(t_1,t_2),t_2}(y_1,y_2) \cac^{(n)}_{(t_1,t_2),t_2}(y_2,y_1) \cac^{(n)}_{s_1,s_2}(w_1,w_2)^3,
\end{align*}

\begin{align*}
&\cm^{\mathbf{1,3,4},(n)}_{t_1,t_2}(y_1,y_2):=\int_{0}^{t_1} ds_1 \int dw_1 \, K_{t_2-s_1}(y_1,w_1)\int_{0}^{t_1} ds_2 \int dw_2 \, K_{t_2-s_2}(y_2,w_2)\\
&\hspace{4cm}\times\cac^{(n)}_{(t_1,t_2),t_2}(y_1,y_2)  \cac^{(n)}_{(t_1,t_2),s_1}(y_2,w_1)\cac^{(n)}_{t_2,s_2}(y_1,w_2)\cac^{(n)}_{s_1,s_2}(w_1,w_2)^2,
\end{align*}

\begin{align*}
&\cm^{\mathbf{1,3,5},(n)}_{t_1,t_2}(y_1,y_2):=\int_{0}^{t_1} ds_1 \int dw_1 \, K_{t_2-s_1}(y_1,w_1)\int_{0}^{t_1} ds_2 \int dw_2 \, K_{t_2-s_2}(y_2,w_2)\\
&\hspace{4cm}\times\cac^{(n)}_{(t_1,t_2),s_2}(y_1,w_2)\cac^{(n)}_{(t_1,t_2),t_2}(y_2,y_1)\cac^{(n)}_{t_2,s_1}(y_2,w_1)\cac^{(n)}_{s_1,s_2}(w_1,w_2)^2,
\end{align*}

\begin{align*}
&\cm^{\mathbf{1,3,6},(n)}_{t_1,t_2}(y_1,y_2):=\int_{0}^{t_1} ds_1 \int dw_1 \, K_{t_2-s_1}(y_1,w_1)\int_{0}^{t_1} ds_2 \int dw_2 \, K_{t_2-s_2}(y_2,w_2)\\
&\hspace{4cm}\times\cac^{(n)}_{(t_1,t_2),s_2}(y_1,w_2)  \cac^{(n)}_{(t_1,t_2),s_1}(y_2,w_1) \cac^{(n)}_{t_2,t_2}(y_1,y_2)\cac^{(n)}_{s_1,s_2}(w_1,w_2)^2,
\end{align*}

\begin{align*}
&\cm^{\mathbf{1,3,7},(n)}_{t_1,t_2}(y_1,y_2):=\int_{0}^{t_1} ds_1 \int dw_1 \, K_{t_2-s_1}(y_1,w_1)\int_{0}^{t_1} ds_2 \int dw_2 \, K_{t_2-s_2}(y_2,w_2)\\
&\hspace{3cm}\times\cac^{(n)}_{(t_1,t_2),s_2}(y_1,w_2) \cac^{(n)}_{(t_1,t_2),s_1}(y_2,w_1)  \cac^{(n)}_{t_2,s_1}(y_2,w_1)\cac^{(n)}_{t_2,s_2}(y_1,w_2)\cac^{(n)}_{s_1,s_2}(w_1,w_2).
\end{align*}

\

\paragraph{Case of $\cm^{\mathbf{1,3,1},(n)}$}

One has
\begin{multline*}
\big|\cm^{\mathbf{1,3,1},(n)}_{t_1,t_2}(y_1+y_2,y_2)\big|\lesssim \cac^{(n),\star}_{(t_1,t_2),(t_1,t_2)}(y_1)  \cac^{(n),\star}_{t_2,t_2}(y_1)\\
\begin{aligned}
&\hspace{1cm}\times\int_{0}^{t_1} ds_1 \int_{0}^{t_1} ds_2 \int dw_1 dw_2 \, K^\star_{t_2-s_1}(y_1-w_1) K^\star_{t_2-s_2}(w_2)\cac^{(n),\star}_{s_1,s_2}(w_1-w_2)^3\\
&\lesssim \cac^{(n),\star}_{(t_1,t_2),(t_1,t_2)}(y_1)  \cac^{(n),\star}_{t_2,t_2}(y_1)\int_{0}^{t_1} ds_1 \int_{0}^{t_1} ds_2\, \Big[ K^\star_{t_2-s_1}\ast \big[K^\star_{t_2-s_2}\ast (\cac^{(n),\star}_{s_1,s_2})^3\big]\Big](y_1),
\end{aligned}
\end{multline*}
and so
\begin{multline*}
\int dy_1 \sup_{y_2} \big|\cm^{\mathbf{1,3,1},(n)}_{t_1,t_2}(y_1+y_2,y_2)\big| \lesssim  \\
\begin{aligned}
&\lesssim \big\|\cac^{(n),\star}_{(t_1,t_2),(t_1,t_2)}  \cac^{(n),\star}_{t_2,t_2}\big\|_{L^1}\int_{0}^{t_1} ds_1 \int_{0}^{t_1} ds_2\,  \Big\| K^\star_{t_2-s_1}\ast \big[K^\star_{t_2-s_2}\ast (\cac^{(n),\star}_{s_1,s_2})^3\big]\Big\|_{L^\infty}\\
&\lesssim \big\|\cac^{(n),\star}_{(t_1,t_2),(t_1,t_2)}\big\|_{L^{\frac{3}{2-\eps}}} \big\|\cac^{(n),\star}_{t_2,t_2}\big\|_{L^{\frac{3}{1+\eps}}}\int_{0}^{t_1} ds_1 \int_{0}^{t_1} ds_2\,  \big\| K^\star_{t_2-s_1}\big\|_{L^1} \big\| K^\star_{t_2-s_2}\big\|_{L^2} \big\|(\cac^{(n),\star}_{s_1,s_2})^3\big\|_{L^2}.
\end{aligned}
\end{multline*}
Using Lemma \ref{lem:k-star-0} and Lemma \ref{lem:cac3},
\begin{multline*}
\int dy_1 \sup_{y_2} \big|\cm^{\mathbf{1,3,1},(n)}_{t_1,t_2}(y_1+y_2,y_2)\big| \lesssim  \\
\begin{aligned}
&\lesssim \big\|\cac^{(n),\star}_{(t_1,t_2),(t_1,t_2)}\big\|_{L^{\frac{3}{2-\eps}}} \int_{0}^{t_1} ds_1\, e^{- (t_2-s_1)} \int_{0}^{t_1} \frac{ds_2}{|t_2-s_2|^{\frac34}}e^{- (t_2-s_2)}   \big\|\cac^{(n),\star}_{s_1,s_2}\big\|_{L^6}^3\\
&\lesssim \big\|\cac^{(n),\star}_{(t_1,t_2),(t_1,t_2)}\big\|_{L^1}^{\frac13-\frac23\eps}\big\|\cac^{(n),\star}_{(t_1,t_2),(t_1,t_2)}\big\|_{L^2}^{\frac23+\frac23\eps} \int_{0}^{t_1} ds_1 \, e^{- (t_2-s_1)} \int_{0}^{t_1} \frac{ds_2\, e^{- (t_2-s_2)}}{|t_2-s_2|^{\frac34}|s_2-s_1|^{\frac34} } \\
&\lesssim |t_2-t_1|^{\frac13-\frac23\eps}|t_2-t_1|^{\frac14(\frac23+\frac23\eps)}\int_{0}^{+\infty} ds_1 \, e^{- s_1} \int_{0}^{+\infty} \frac{ds_2}{|s_2|^{\frac34}}e^{- s_2}\frac{1}{|s_2-s_1|^{\frac34}}\\
& \lesssim |t_2-t_1|^{\frac12-\frac12\eps},
\end{aligned}
\end{multline*}
uniformly over $0\leq t_1,t_2\leq 2$.

\

\paragraph{Case of $\cm^{\mathbf{1,3,2},(n)}$}

One has
\begin{multline*}
\big| \cm^{\mathbf{1,3,2},(n)}_{t_1,t_2}(y_1+y_2,y_2)\big|\lesssim \cac^{(n),\star}_{(t_1,t_2),(t_1,t_2)}(y_1)\int_{0}^{t_1} ds_1 \int dw_1 \, K^\star_{t_2-s_1}(y_1-w_1) \\
\begin{aligned}
&\hspace{3cm}\times \cac^{(n),\star}_{t_2,s_1}(w_1)\int_{0}^{t_1} ds_2 \int dw_2 \, K^\star_{t_2-s_2}(w_2)\cac^{(n),\star}_{t_2,s_2}(y_1-w_2)\cac^{(n),\star}_{s_1,s_2}(w_1-w_2)^2\\
&\lesssim \cac^{(n),\star}_{(t_1,t_2),(t_1,t_2)}(y_1)\int_{0}^{t_1} ds_1\, \big\|\cac^{(n),\star}_{t_2,s_1}\big\|_{L^\infty}\int_{0}^{t_1} ds_2\, \big\|\cac^{(n),\star}_{t_2,s_2}\big\|_{L^\infty}\\
& \hspace{4cm} \times \int dw_1 \, K^\star_{t_2-s_1}(y_1-w_1) \int dw_2 \, K^\star_{t_2-s_2}(w_2)\cac^{(n),\star}_{s_1,s_2}(w_1-w_2)^2\\
&\lesssim \cac^{(n),\star}_{(t_1,t_2),(t_1,t_2)}(y_1)\int_{0}^{t_1} ds_1\, \big\|\cac^{(n),\star}_{t_2,s_1}\big\|_{L^\infty}\int_{0}^{t_1} ds_2\, \big\|\cac^{(n),\star}_{t_2,s_2}\big\|_{L^\infty}\Big[  K^\star_{t_2-s_1}\ast\big[ K^\star_{t_2-s_2}\ast (\cac^{(n),\star}_{s_1,s_2})^2\big]\Big](y_1),
\end{aligned}
\end{multline*}
and so
\begin{multline*}
\int dy_1 \sup_{y_2} \big| \cm^{\mathbf{1,3,2},(n)}_{t_1,t_2}(y_1+y_2,y_2)\big| \lesssim \\
\begin{aligned}
&\lesssim \big\|\cac^{(n),\star}_{(t_1,t_2),(t_1,t_2)}\big\|_{L^{\frac{1}{1-\eps}}}\int_{0}^{t_1} ds_1\, \big\|\cac^{(n),\star}_{t_2,s_1}\big\|_{L^\infty}\int_{0}^{t_1} ds_2\, \big\|\cac^{(n),\star}_{t_2,s_2}\big\|_{L^\infty}\Big\|  K^\star_{t_2-s_1}\ast\big[ K^\star_{t_2-s_2}\ast (\cac^{(n),\star}_{s_1,s_2})^2\big]\Big\|_{L^{\frac{1}{\eps}}}\\
&\lesssim \big\|\cac^{(n),\star}_{(t_1,t_2),(t_1,t_2)}\big\|_{L^1}^{1-2\eps}\big\|\cac^{(n),\star}_{(t_1,t_2),(t_1,t_2)}\big\|_{L^2}^{2\eps}\\
&\hspace{3cm}\times \int_{0}^{t_1} ds_1\, \big\|\cac^{(n),\star}_{t_2,s_1}\big\|_{L^\infty}\int_{0}^{t_1} ds_2\, \big\|\cac^{(n),\star}_{t_2,s_2}\big\|_{L^\infty}\big\|K^\star_{t_2-s_1}\big\|_{L^1}\big\| K^\star_{t_2-s_2}\big\|_{L^1} \big\|\cac^{(n),\star}_{s_1,s_2}\big\|_{L^{\frac{2}{\eps}}}^2.
\end{aligned}
\end{multline*}
By using Lemma \ref{lem:k-star-0} and Lemma \ref{lem:cac3}, we obtain
\begin{multline*}
\int dy_1 \sup_{y_2} \big| \cm^{\mathbf{1,3,2},(n)}_{t_1,t_2}(y_1+y_2,y_2)\big| \lesssim \\
\begin{aligned}
&\lesssim |t_2-t_1|^{1-2\eps}|t_2-t_1|^{\frac{\eps}{2}}\int_{0}^{t_1} \frac{ds_1}{|t_2-s_1|^{\frac12}}e^{- (t_2-s_1)} \int_{0}^{t_1} \frac{ds_2}{|t_2-s_2|^{\frac12}}e^{- (t_2-s_2)}\frac{1}{|s_2-s_1|^{1-\frac32 \eps}}\\
&\lesssim |t_2-t_1|^{1-\frac32\eps}\int_{0}^{+\infty} \frac{ds_1}{|s_1|^{\frac12}}e^{- s_1} \int_{0}^{+\infty} \frac{ds_2}{|s_2|^{\frac12}}e^{- s_2}\frac{1}{|s_2-s_1|^{1-\frac32 \eps}}\\
&\lesssim |t_2-t_1|^{1-\frac32\eps},
\end{aligned}
\end{multline*}
uniformly over $0\leq t_1,t_2\leq 2$.

\

\paragraph{Case of $\cm^{\mathbf{1,3,3},(n)}$}

An identical calculation to the one behind \eqref{ref-1-1-papier} leads to
\begin{multline*}
\int dy_1 \sup_{y_2} \big|\cm^{\mathbf{1,3,3},(n)}_{t_1,t_2}(y_1+y_2,y_2)\big|  \lesssim \\
\begin{aligned}
&\lesssim \big\|\cac^{(n),\star}_{(t_1,t_2),t_2}\big\|_{L^2}^2\int_{0}^{t_1} ds_1 \int_{0}^{t_1} ds_2\,  \big\| K^\star_{t_2-s_1}\big\|_{L^1} \big\| K^\star_{t_2-s_2}\big\|_{L^2} \big\|\cac^{(n),\star}_{s_1,s_2}\big\|_{L^6}^3.
\end{aligned}
\end{multline*}
Then, by combining Lemma \ref{lem:k-star-0} and Lemma \ref{lem:cac3},
\begin{align*}
\int dy_1 \sup_{y_2} \big|\cm^{\mathbf{1,3,3},(n)}_{t_1,t_2}(y_1+y_2,y_2)\big|  
&\lesssim |t_2-t_1|^{\frac12}\int_{0}^{t_1} ds_1\, e^{- (t_2-s_1)} \int_{0}^{t_1} \frac{ds_2}{|t_2-s_2|^{\frac34}}e^{- (t_2-s_2)}\frac{1}{|s_2-s_1|^{\frac34}}\\
&\lesssim |t_2-t_1|^{\frac12}\int_{0}^{+\infty} ds_1\, e^{- s_1} \int_{0}^{+\infty} \frac{ds_2}{|s_2|^{\frac34}}e^{- s_2}\frac{1}{|s_2-s_1|^{\frac34}}\\
&\lesssim |t_2-t_1|^{\frac12},
\end{align*}
uniformly over $0\leq t_1,t_2\leq 2$.

\

\paragraph{Case of $\cm^{\mathbf{1,3,4},(n)}$}

One has
\begin{align}
&\big|\cm^{\mathbf{1,3,4},(n)}_{t_1,t_2}(y_1+y_2,y_2)\big|\lesssim \cac^{(n),\star}_{(t_1,t_2),t_2}(y_1)\int_{0}^{t_1} ds_1  \int_{0}^{t_1} ds_2 \int dw_1 \, K^{\star}_{t_2-s_1}(y_1-w_1)\nonumber \\
&\hspace{2cm}\times \cac^{(n),\star}_{(t_1,t_2),s_1}(w_1)\int dw_2 \, K^{\star}_{t_2-s_2}(w_2) \cac^{(n),\star}_{t_2,s_2}(y_1-w_2)\cac^{(n),\star}_{s_1,s_2}(w_1-w_2)^2 \nonumber\\
&\lesssim \cac^{(n),\star}_{(t_1,t_2),t_2}(y_1)\int_{0}^{t_1} ds_1  \int_{0}^{t_1} ds_2\, \big\|\cac^{(n),\star}_{t_2,s_2}\big\|_{L^\infty}  \bigg[K^{\star}_{t_2-s_1}\ast \Big[\cac^{(n),\star}_{(t_1,t_2),s_1}\cdot\big[ K^{\star}_{t_2-s_2}\ast (\cac^{(n),\star}_{s_1,s_2})^2\big]\Big]\bigg](y_1),\label{same-quant}
\end{align}
and so, using the results of Lemma \ref{lem:k-star-0} and Lemma \ref{lem:cac3},
\begin{multline*}
\int dy_1 \sup_{y_2} \big|\cm^{\mathbf{1,3,4},(n)}_{t_1,t_2}(y_1+y_2,y_2)\big|  \lesssim \\
\begin{aligned}
&\lesssim \big\|\cac^{(n),\star}_{(t_1,t_2),t_2}\big\|_{L^2}\int_{0}^{t_1} ds_1  \int_{0}^{t_1} \frac{ds_2}{|t_2-s_2|^{\frac12}}  \bigg\|K^{\star}_{t_2-s_1}\ast \Big[\cac^{(n),\star}_{(t_1,t_2),s_1}\cdot\big[ K^{\star}_{t_2-s_2}\ast (\cac^{(n),\star}_{s_1,s_2})^2\big]\Big]\bigg\|_{L^2}\\
&\lesssim |t_1-t_2|^{\frac14}\int_{0}^{t_1} ds_1  \int_{0}^{t_1} \frac{ds_2}{|t_2-s_2|^{\frac12}}  \big\|K^{\star}_{t_2-s_1}\big\|_{L^2} \Big\|\cac^{(n),\star}_{(t_1,t_2),s_1}\cdot\big[ K^{\star}_{t_2-s_2}\ast (\cac^{(n),\star}_{s_1,s_2})^2\big]\Big\|_{L^1}\\
&\lesssim |t_1-t_2|^{\frac14}\int_{0}^{t_1} \frac{ds_1}{|t_2-s_1|^{\frac34}}e^{- (t_2-s_1)}  \int_{0}^{t_1} \frac{ds_2}{|t_2-s_2|^{\frac12}}   \big\|\cac^{(n),\star}_{(t_1,t_2),s_1}\big\|_{L^2} \big\| K^{\star}_{t_2-s_2}\ast (\cac^{(n),\star}_{s_1,s_2})^2\big\|_{L^2}\\
&\lesssim |t_1-t_2|^{\frac12}\int_{0}^{t_1} \frac{ds_1}{|t_2-s_1|^{\frac34}} e^{- (t_2-s_1)}  \int_{0}^{t_1} \frac{ds_2}{|t_2-s_2|^{\frac12}} \  \big\| K^{\star}_{t_2-s_2}\big\|_{L^1} \big\|\cac^{(n),\star}_{s_1,s_2}\big\|_{L^4}^2\\
&\lesssim |t_1-t_2|^{\frac12}\int_{0}^{t_1} \frac{ds_1}{|t_2-s_1|^{\frac34}} e^{- (t_2-s_1)} \int_{0}^{t_1} \frac{ds_2}{|t_2-s_2|^{\frac12}}e^{- (t_2-s_2)}  \frac{1}{|s_2-s_1|^{\frac14}}\\
&\lesssim |t_1-t_2|^{\frac12}\int_{0}^{+\infty} \frac{ds_1}{|s_1|^{\frac34}} e^{- s_1} \int_{0}^{+\infty} \frac{ds_2}{|s_2|^{\frac12}}e^{- s_2}  \frac{1}{|s_2-s_1|^{\frac14}}\lesssim |t_1-t_2|^{\frac12},
\end{aligned}
\end{multline*}
uniformly over $0\leq t_1,t_2\leq 2$.

\

\paragraph{Case of $\cm^{\mathbf{1,3,5},(n)}$}

One has
\begin{align*}
&\big|\cm^{\mathbf{1,3,5},(n)}_{t_1,t_2}(y_1+y_2,y_2)\big| \lesssim \cac^{(n),\star}_{(t_1,t_2),t_2}(y_1) \int_{0}^{t_1} ds_1  \int_{0}^{t_1} ds_2\\
&\hspace{1cm}\times \int dw_1 \, K^\star_{t_2-s_1}(y_1-w_1)\cac^{(n),\star}_{t_2,s_1}(w_1) \int dw_2 \,  K^\star_{t_2-s_2}(w_2)\cac^{(n),\star}_{(t_1,t_2),s_2}(y_1-w_2)\cac^{(n),\star}_{s_1,s_2}(w_1-w_2)^2\\
&\lesssim \cac^{(n),\star}_{(t_1,t_2),t_2}(y_1) \int_{0}^{t_1} \frac{ds_1}{|t_2-s_1|^{\frac12}}  \int_{0}^{t_1} ds_2\\
&\hspace{1cm}\times \int dw_2 \,  K^\star_{t_2-s_2}(w_2)\cac^{(n),\star}_{(t_1,t_2),s_2}(y_1-w_2)\int dw_1 \, K^\star_{t_2-s_1}(y_1-w_1)\cac^{(n),\star}_{s_1,s_2}(w_1-w_2)^2 \\
&\lesssim \cac^{(n),\star}_{(t_1,t_2),t_2}(y_1) \int_{0}^{t_1} \frac{ds_2}{|t_2-s_2|^{\frac12}}  \int_{0}^{t_1} ds_1\\
&\hspace{1cm}\times \int dw_1 \,  K^\star_{t_2-s_1}(w_1)\cac^{(n),\star}_{(t_1,t_2),s_1}(y_1-w_1)\int dw_2 \, K^\star_{t_2-s_2}(y_1-w_2)\cac^{(n),\star}_{s_1,s_2}(w_1-w_2)^2 \\
&\lesssim \cac^{(n),\star}_{(t_1,t_2),t_2}(y_1) \int_{0}^{t_1} ds_1\int_{0}^{t_1} \frac{ds_2}{|t_2-s_2|^{\frac12}}  \\
&\hspace{1cm}\times \int dw_1 \,  K^\star_{t_2-s_1}(y_1-w_1)\cac^{(n),\star}_{(t_1,t_2),s_1}(w_1)\int dw_2 \, K^\star_{t_2-s_2}(w_2)\cac^{(n),\star}_{s_1,s_2}(w_1-w_2)^2 ,
\end{align*}
and we are dealing here with the same quantity as in \eqref{same-quant}.

\

\paragraph{Case of $\cm^{\mathbf{1,3,6},(n)}$}

One has
\begin{align*}
&\big|\cm^{\mathbf{1,3,6},(n)}_{t_1,t_2}(y_1+y_2,y_2)\big|\lesssim \cac^{(n),\star}_{t_2,t_2}(y_1)\int_{0}^{t_1} ds_1\int_{0}^{t_1} ds_2 \\
&\hspace{1cm}\times\int dw_1 \, K^\star_{t_2-s_1}(y_1-w_1)\cac^{(n),\star}_{(t_1,t_2),s_1}(w_1) \int dw_2\,  K^{\star}_{t_2-s_2}(w_2)\cac^{(n),\star}_{(t_1,t_2),s_2}(y_1-w_2)   \cac^{(n),\star}_{s_1,s_2}(w_1-w_2)^2\\
&\lesssim \cac^{(n),\star}_{t_2,t_2}(y_1)\int_{0}^{t_1} ds_1\int_{0}^{t_1} \frac{ds_2}{|t_2-s_2|^{\frac12}} \\
&\hspace{1cm}\times\int dw_1 \, K^\star_{t_2-s_1}(y_1-w_1)\cac^{(n),\star}_{(t_1,t_2),s_1}(w_1) \int dw_2\,  K^{\star}_{t_2-s_2}(w_2)   \cac^{(n),\star}_{s_1,s_2}(w_1-w_2)^2\\
&\lesssim \cac^{(n),\star}_{t_2,t_2}(y_1)\int_{0}^{t_1} ds_1  \int_{0}^{t_1} \frac{ds_2}{|t_2-s_2|^{\frac12}}  \bigg[K^{\star}_{t_2-s_1}\ast \Big[\cac^{(n),\star}_{(t_1,t_2),s_1}\cdot\big[ K^{\star}_{t_2-s_2}\ast (\cac^{(n),\star}_{s_1,s_2})^2\big]\Big]\bigg](y_1)
\end{align*}
and so, thanks to Lemma \ref{lem:k-star-0} and Lemma \ref{lem:cac3},
\begin{align*}
&\int dy_1 \sup_{y_2} \big|\cm^{\mathbf{1,3,6},(n)}_{t_1,t_2}(y_1+y_2,y_2)\big|  \lesssim\\
&\lesssim \big\|\cac^{(n),\star}_{t_2,t_2}\big\|_{L^{\frac{3}{1+\eps}}}\int_{0}^{t_1} ds_1  \int_{0}^{t_1} \frac{ds_2}{|t_2-s_2|^{\frac12}}  \bigg\|K^{\star}_{t_2-s_1}\ast \Big[\cac^{(n),\star}_{(t_1,t_2),s_1}\cdot\big[ K^{\star}_{t_2-s_2}\ast (\cac^{(n),\star}_{s_1,s_2})^2\big]\Big]\bigg\|_{L^{\frac{3}{2-\eps}}}\\
&\lesssim \int_{0}^{t_1} ds_1  \int_{0}^{t_1} \frac{ds_2}{|t_2-s_2|^{\frac12}}  \big\|K^{\star}_{t_2-s_1}\big\|_{L^{\frac{3}{2-\eps}}} \Big\|\cac^{(n),\star}_{(t_1,t_2),s_1}\cdot\big[ K^{\star}_{t_2-s_2}\ast (\cac^{(n),\star}_{s_1,s_2})^2\big]\Big\|_{L^{1}}\\
&\lesssim \int_{0}^{t_1} \frac{ds_1}{|t_2-s_1|^{\frac12+\frac{\eps}{2}}} e^{- (t_2-s_1)}\big\|\cac^{(n),\star}_{(t_1,t_2),s_1}\big\|_{L^{\frac{3}{3-2\eps}}} \int_{0}^{t_1} \frac{ds_2}{|t_2-s_2|^{\frac12}}    \big\| K^{\star}_{t_2-s_2}\ast (\cac^{(n),\star}_{s_1,s_2})^2\big\|_{L^{\frac{3}{2\eps}}}\\
&\lesssim \int_{0}^{t_1} \frac{ds_1}{|t_2-s_1|^{\frac12+\frac{\eps}{2}}}e^{-(t_2-s_1)}\big\|\cac^{(n),\star}_{(t_1,t_2),s_1}\big\|_{L^{1}}^{1-\frac43\eps}\big\|\cac^{(n),\star}_{(t_1,t_2),s_1}\big\|_{L^{2}}^{\frac43\eps}  \int_{0}^{t_1} \frac{ds_2}{|t_2-s_2|^{\frac12}}   \big\| K^{\star}_{t_2-s_2}\big\|_{L^{\frac{3}{2(1+\eps)}}} \big\|\cac^{(n),\star}_{s_1,s_2}\big\|_{L^{6}}^2\\
&\lesssim |t_1-t_2|^{1-\eps}\int_{0}^{t_1} \frac{ds_1}{|t_2-s_1|^{\frac12+\frac{\eps}{2}}} e^{-(t_2-s_1)} \int_{0}^{t_1} \frac{ds_2}{|t_2-s_2|^{1-\eps}} e^{-(t_2-s_2)} \frac{1}{|s_2-s_1|^{\frac12}}\\
&\lesssim |t_1-t_2|^{1-\eps}\int_{0}^{+\infty} \frac{ds_1}{|s_1|^{\frac12+\frac{\eps}{2}}} e^{-s_1} \int_{0}^{+\infty} \frac{ds_2}{|s_2|^{1-\eps}} e^{- s_2} \frac{1}{|s_2-s_1|^{\frac12}}\lesssim |t_1-t_2|^{1-\eps},
\end{align*}
uniformly over $0\leq t_1,t_2\leq 2$.

\

\paragraph{Case of $\cm^{\mathbf{1,3,7},(n)}$}

One has
\begin{align*}
&\big|\cm^{\mathbf{1,3,7},(n)}_{t_1,t_2}(y_1+y_2,y_2)\big|\lesssim \int_{0}^{t_1} ds_1 \int_{0}^{t_1} ds_2 \int dw_1  \, K^\star_{t_2-s_1}(y_1-w_1)\cac^{(n),\star}_{(t_1,t_2),s_1}(w_1)  \cac^{(n),\star}_{t_2,s_1}(w_1)\\
&\hspace{3cm}\times \int dw_2\,  K^\star_{t_2-s_2}(w_2)\cac^{(n),\star}_{(t_1,t_2),s_2}(y_1-w_2) \cac^{(n),\star}_{t_2,s_2}(y_1-w_2)\cac^{(n),\star}_{s_1,s_2}(w_1-w_2)\\
&\lesssim \int_{0}^{t_1} \frac{ds_1}{|t_2-s_1|^{\frac12}} \int_{0}^{t_1}\frac{ds_2}{|t_2-s_2|^{\frac12}} \frac{1}{|s_2-s_1|^{\frac12}} \\
&\hspace{2cm}\times \int dw_1  \, K^\star_{t_2-s_1}(y_1-w_1)\cac^{(n),\star}_{(t_1,t_2),s_1}(w_1) \int dw_2\,  K^\star_{t_2-s_2}(w_2)\cac^{(n),\star}_{(t_1,t_2),s_2}(y_1-w_2) \\
&\lesssim \int_{0}^{t_1} \frac{ds_1}{|t_2-s_1|^{\frac12}} \int_{0}^{t_1}\frac{ds_2}{|t_2-s_2|^{\frac12}} \frac{1}{|s_2-s_1|^{\frac12}} \big[ K^\star_{t_2-s_1}\ast\cac^{(n),\star}_{(t_1,t_2),s_1}\big](y_1) \big[  K^\star_{t_2-s_2}\ast \cac^{(n),\star}_{(t_1,t_2),s_2}\big](y_1), 
\end{align*}
and so, by Lemma \ref{lem:k-star-0} and Lemma \ref{lem:cac3},
\begin{multline*}
\int dy_1 \sup_{y_2}\big|\cm^{\mathbf{1,3,7},(n)}_{t_1,t_2}(y_1+y_2,y_2)\big| \lesssim \\
\begin{aligned}
&\lesssim \int_{0}^{t_1} \frac{ds_1}{|t_2-s_1|^{\frac12}} \int_{0}^{t_1}\frac{ds_2}{|t_2-s_2|^{\frac12}} \frac{1}{|s_2-s_1|^{\frac12}} \big\| K^\star_{t_2-s_1}\ast\cac^{(n),\star}_{(t_1,t_2),s_1}\big\|_{L^2} \big\|  K^\star_{t_2-s_2}\ast \cac^{(n),\star}_{(t_1,t_2),s_2}\big\|_{L^2}\\ 
&\lesssim \int_{0}^{t_1} \frac{ds_1}{|t_2-s_1|^{\frac12}} \int_{0}^{t_1}\frac{ds_2}{|t_2-s_2|^{\frac12}} \frac{1}{|s_2-s_1|^{\frac12}} \big\| K^\star_{t_2-s_1}\big\|_{L^1}\big\| \cac^{(n),\star}_{(t_1,t_2),s_1}\big\|_{L^2} \big\|  K^\star_{t_2-s_2}\big\|_{L^1}\big\|  \cac^{(n),\star}_{(t_1,t_2),s_2}\big\|_{L^2} \\
&\lesssim |t_2-t_1|^{\frac12} \int_{0}^{t_1} \frac{ds_1}{|t_2-s_1|^{\frac12}}e^{- (t_2-s_1)} \int_{0}^{t_1}\frac{ds_2}{|t_2-s_2|^{\frac12}}e^{-(t_2-s_2)} \frac{1}{|s_2-s_1|^{\frac12}} \\
&\lesssim |t_2-t_1|^{\frac12} \int_{0}^{+\infty} \frac{ds_1}{|s_1|^{\frac12}}e^{-s_1}  \int_{0}^{+\infty}\frac{ds_2}{|s_2|^{\frac12}}e^{- s_2} \frac{1}{|s_2-s_1|^{\frac12}} \lesssim |t_2-t_1|^{\frac12},
\end{aligned}
\end{multline*}
uniformly over $0\leq t_1,t_2\leq 2$.

\

\

\subsubsection{Increments: second diagram}

We deal here with
\begin{align*}
\calt^{\mathbf{2},(n)}_t(y):= \int_{0}^t ds \int dw \, K_{t-s}(y,w) \cac^{(n)}_{t,s}(y,w)I^W_3\big(F^{(n)}_{t,y}\otimes F^{(n)}_{s,w}\otimes F^{(n)}_{s,w}\big).
\end{align*}
For $0\leq t_1<t_2$, we decompose the increment of $\calt^{\mathbf{2},(n)}$ as
\begin{align*}
\calt^{\mathbf{2},(n)}_{t_1,t_2}(y)&:=\calt^{\mathbf{2},(n)}_{t_2}(y)-\calt^{\mathbf{2},(n)}_{t_1}(y)\\
&=\calt^{\mathbf{2,1},(n)}_{t_1,t_2}(y)+\calt^{\mathbf{2,2},(n)}_{t_1,t_2}(y)+\calt^{\mathbf{2,3},(n)}_{t_1,t_2}(y)+\calt^{\mathbf{2,4},(n)}_{t_1,t_2}(y),
\end{align*}
with
\begin{align*}
\calt^{\mathbf{2,1},(n)}_{t_1,t_2}(y)&:=\int_{t_1}^{t_2} ds \int dw \, K_{t_2-s}(y,w) \cac^{(n)}_{t_2,s}(y,w)I^W_3\big(F^{(n)}_{t_2,y}\otimes F^{(n)}_{s,w}\otimes F^{(n)}_{s,w}\big),
\end{align*}

\begin{align*}
\calt^{\mathbf{2,2},(n)}_{t_1,t_2}(y)&:=\int_{0}^{t_1} ds \int dw \, K_{t_1-s,t_2-s}(y,w) \cac^{(n)}_{t_2,s}(y,w)I^W_3\big(F^{(n)}_{t_2,y}\otimes F^{(n)}_{s,w}\otimes F^{(n)}_{s,w}\big),
\end{align*}

\begin{align*}
\calt^{\mathbf{2,3},(n)}_{t_1,t_2}(y)&:=\int_{0}^{t_1} ds \int dw \, K_{t_1-s}(y,w) \cac^{(n)}_{(t_1,t_2),s}(y,w)I^W_3\big(F^{(n)}_{t_2,y}\otimes F^{(n)}_{s,w}\otimes F^{(n)}_{s,w}\big),
\end{align*}

\begin{align*}
\calt^{\mathbf{2,4},(n)}_{t_1,t_2}(y)&:=\int_{0}^{t_1} ds \int dw \, K_{t_1-s}(y,w) \cac^{(n)}_{t_1,s}(y,w)I^W_3\big((F^{(n)}_{t_2,y}-F^{(n)}_{t_1,y})\otimes F^{(n)}_{s,w}\otimes F^{(n)}_{s,w}\big).
\end{align*}

\

\noindent
\underline{\textit{Case of $\calt^{\mathbf{2,1},(n)}$}}. One has
\begin{align*}
\cm^{\mathbf{2,1},(n)}_{t_1,t_2}(y_1,y_2):=\mathbb{E}\Big[\calt^{\mathbf{2,1},(n)}_{t_1,t_2}(y_1)\calt^{\mathbf{2,1},(n)}_{t_1,t_2}(y_2)\Big]=\sum_{\mathbf{a}=1}^2 c_{\mathbf{a}}\cm^{\mathbf{2,1,a},(n)}_{t_1,t_2}(y_1,y_2),
\end{align*}
with
\begin{multline*}
\cm^{\mathbf{2,1,1},(n)}_{t_1,t_2}(y_1,y_2):=\int_{t_1}^{t_2}ds_1\int_{t_1}^{t_2}ds_2\int dw_1 dw_2 K_{t_2-s_1}(y_1,w_1) \cac^{(n)}_{t_2,s_1}(y_1,w_1) \\
\times K_{t_2-s_2}(y_2,w_2) \cac^{(n)}_{t_2,s_2}(y_2,w_2)\cac^{(n)}_{t_2,t_2}(y_1,y_2) \cac^{(n)}_{s_1,s_2}(w_1,w_2)^2
\end{multline*}
and
\begin{multline*}
\cm^{\mathbf{2,1,2},(n)}_{t_1,t_2}(y_1,y_2):=\int_{t_1}^{t_2}ds_1\int_{t_1}^{t_2}ds_2\int dw_1 dw_2 K_{t_2-s_1}(y_1,w_1) \cac^{(n)}_{t_2,s_1}(y_1,w_1)\\
\times K_{t_2-s_2}(y_2,w_2) \cac^{(n)}_{t_2,s_2}(y_2,w_2)\cac^{(n)}_{t_2,s_2}(y_1,w_2) \cac^{(n)}_{t_2,s_1}(y_2,w_1)\cac^{(n)}_{s_1,s_2}(w_1,w_2).
\end{multline*}

\

\paragraph{Case of $\cm^{\mathbf{2,1,1},(n)}$}

One has
\begin{multline*}
\cm^{\mathbf{2,1,1},(n)}_{t_1,t_2}(y_1+y_2,y_2)=\cac^{(n)}_{t_2,t_2}(y_1+y_2,y_2)\int_{t_1}^{t_2}ds_1\int_{t_1}^{t_2}ds_2\int dw_1 dw_2K_{t_2-s_1}(y_1+y_2,w_1+y_2) \\
\begin{aligned}
&\hspace{1.5cm}\times  \cac^{(n)}_{t_2,s_1}(y_1+y_2,w_1+y_2)K_{t_2-s_2}(y_2,w_2+y_2) \cac^{(n)}_{t_2,s_2}(y_2,w_2+y_2) \cac^{(n)}_{s_1,s_2}(w_1+y_2,w_2+y_2)^2\\
&\lesssim \cac^{(n),\star}_{t_2,t_2}(y_1)\int_{t_1}^{t_2}\frac{ds_1}{|t_2-s_1|^{\frac12}}\int_{t_1}^{t_2}\frac{ds_2}{|t_2-s_2|^{\frac12}}\int dw_1 \, K_{t_2-s_1}^\star(y_1-w_1)\int dw_2\,  K^\star_{t_2-s_2}(w_2) \cac^{(n),\star}_{s_1,s_2}(w_1-w_2)^2\\
&\lesssim \cac^{(n),\star}_{t_2,t_2}(y_1)\int_{t_1}^{t_2}\frac{ds_1}{|t_2-s_1|^{\frac12}}\int_{t_1}^{t_2}\frac{ds_2}{|t_2-s_2|^{\frac12}}\Big[ K_{t_2-s_1}^\star \ast\big[ K^\star_{t_2-s_2}\ast (\cac^{(n),\star}_{s_1,s_2})^2\big]\Big](y_1),
\end{aligned}
\end{multline*}
and so
\begin{multline*}
\int dy_1 \sup_{y_2}\, \big|\cm^{\mathbf{2,1,1},(n)}_{t_1,t_2}(y_1+y_2,y_2)\big| \lesssim \\
\begin{aligned}
&\lesssim \int_{t_1}^{t_2}\frac{ds_1}{|t_2-s_1|^{\frac12}}\int_{t_1}^{t_2}\frac{ds_2}{|t_2-s_2|^{\frac12}}\int dy_1\, \cac^{(n),\star}_{t_2,t_2}(y_1)\Big[ K_{t_2-s_1}^\star \ast\big[ K^\star_{t_2-s_2}\ast (\cac^{(n),\star}_{s_1,s_2})^2\big]\Big](y_1)\\
&\lesssim \big\|\cac^{(n),\star}_{t_2,t_2}\big\|_{L^2}\int_{t_1}^{t_2}\frac{ds_1}{|t_2-s_1|^{\frac12}}\int_{t_1}^{t_2}\frac{ds_2}{|t_2-s_2|^{\frac12}}\Big\| K_{t_2-s_1}^\star \ast\big[ K^\star_{t_2-s_2}\ast (\cac^{(n),\star}_{s_1,s_2})^2\big]\Big\|_{L^2}\\
&\lesssim \big\|\cac^{(n),\star}_{t_2,t_2}\big\|_{L^2}\int_{t_1}^{t_2}\frac{ds_1}{|t_2-s_1|^{\frac12}}\int_{t_1}^{t_2}\frac{ds_2}{|t_2-s_2|^{\frac12}}\big\| K_{t_2-s_1}^\star \big\|_{L^1}\big\|  K^\star_{t_2-s_2}\big\|_{L^1} \big\| (\cac^{(n),\star}_{s_1,s_2})^2\big\|_{L^2}\\
&\lesssim \big\|\cac^{(n),\star}_{t_2,t_2}\big\|_{L^2}\int_{t_1}^{t_2}\frac{ds_1}{|t_2-s_1|^{\frac12}}e^{-(t_2-s_1)}\int_{t_1}^{t_2}\frac{ds_2}{|t_2-s_2|^{\frac12}}e^{-(t_2-s_2)}\big\| \cac^{(n),\star}_{s_1,s_2}\big\|_{L^4}^2\\
&\lesssim \int_{t_1}^{t_2}\frac{ds_1}{|t_2-s_1|^{\frac12}}e^{-(t_2-s_1)}\int_{t_1}^{t_2}\frac{ds_2}{|t_2-s_2|^{\frac12}}e^{-(t_2-s_2)}\frac{1}{|s_1-s_2|^{\frac14}}.
\end{aligned}
\end{multline*}
From here we deduce
\begin{align*}
\int dy_1 \sup_{y_2}\, \big|\cm^{\mathbf{2,1,1},(n)}_{t_1,t_2}(y_1+y_2,y_2)\big|&\lesssim \int_{t_1}^{t_2}\frac{ds_1}{|t_2-s_1|^{\frac12}}\int_{t_1}^{t_2}\frac{ds_2}{|t_2-s_2|^{\frac12}}\frac{1}{|s_1-s_2|^{\frac14}}\\
&\lesssim |t_2-t_1|^{\frac34}\int_{0}^{1}\frac{dr_1}{|1-r_1|^{\frac12}}\int_{0}^{1}\frac{dr_2}{|1-r_2|^{\frac12}}\frac{1}{|r_1-r_2|^{\frac14}} \lesssim |t_2-t_1|^{\frac34},
\end{align*}
uniformly over $0\leq t_1,t_2\leq 2$.

\

\paragraph{Case of $\cm^{\mathbf{2,1,2},(n)}$}

One has
\begin{align*}
&\cm^{\mathbf{2,1,2},(n)}_{t_1,t_2}(y_1+y_2,y_2)=\int_{t_1}^{t_2}ds_1\int_{t_1}^{t_2}ds_2\int dw_1 dw_2\, K_{t_2-s_1}(y_1+y_2,w_1+y_2) \cac^{(n)}_{t_2,s_1}(y_1+y_2,w_1+y_2)\\
&\times  K_{t_2-s_2}(y_2,w_2+y_2) \cac^{(n)}_{t_2,s_2}(y_2,w_2+y_2)\cac^{(n)}_{t_2,s_2}(y_1+y_2,w_2+y_2) \cac^{(n)}_{t_2,s_1}(y_2,w_1+y_2)\cac^{(n)}_{s_1,s_2}(w_1+y_2,w_2+y_2)\\
&\lesssim \int_{t_1}^{t_2}\frac{ds_1}{|t_2-s_1|^{\frac12}}\int_{t_1}^{t_2}\frac{ds_2}{|t_2-s_2|^{\frac12}}\frac{1}{|s_1-s_2|^{\frac12}}\int dw_1 \, K^\star_{t_2-s_1}(y_1-w_1) \cac^{(n),\star}_{t_2,s_1}(w_1)\int dw_2\,K^\star_{t_2-s_2}(w_2) \cac^{(n),\star}_{t_2,s_2}(y_1-w_2) \\
&\lesssim \int_{t_1}^{t_2}\frac{ds_1}{|t_2-s_1|^{\frac12}}\int_{t_1}^{t_2}\frac{ds_2}{|t_2-s_2|^{\frac12}}\frac{1}{|s_1-s_2|^{\frac12}}\big[K^\star_{t_2-s_1} \ast \cac^{(n),\star}_{t_2,s_1}\big](y_1) \big[K^\star_{t_2-s_2}\ast \cac^{(n),\star}_{t_2,s_2}\big](y_1) ,
\end{align*}
and so
\begin{multline*}
\int dy_1 \sup_{y_2}\, \big|\cm^{\mathbf{2,1,2},(n)}_{t_1,t_2}(y_1+y_2,y_2)\big|\lesssim \\
\begin{aligned}
&\lesssim \int_{t_1}^{t_2}\frac{ds_1}{|t_2-s_1|^{\frac12}}\int_{t_1}^{t_2}\frac{ds_2}{|t_2-s_2|^{\frac12}}\frac{1}{|s_1-s_2|^{\frac12}}\big\|K^\star_{t_2-s_1} \ast \cac^{(n),\star}_{t_2,s_1}\big\|_{L^2} \big\|K^\star_{t_2-s_2}\ast \cac^{(n),\star}_{t_2,s_2}\big\|_{L^2} \\
&\lesssim \int_{t_1}^{t_2}\frac{ds_1}{|t_2-s_1|^{\frac12}}\int_{t_1}^{t_2}\frac{ds_2}{|t_2-s_2|^{\frac12}}\frac{1}{|s_1-s_2|^{\frac12}}\big\|K^\star_{t_2-s_1} \big\|_{L^1} \big\| \cac^{(n),\star}_{t_2,s_1}\big\|_{L^2} \big\|K^\star_{t_2-s_2}\big\|_{L^1} \big\| \cac^{(n),\star}_{t_2,s_2}\big\|_{L^2} \\
&\lesssim \int_{t_1}^{t_2}\frac{ds_1}{|t_2-s_1|^{\frac12}}e^{-(t_2-s_1)}\int_{t_1}^{t_2}\frac{ds_2}{|t_2-s_2|^{\frac12}}\frac{1}{|s_1-s_2|^{\frac12}}e^{-(t_2-s_2)}.
\end{aligned}
\end{multline*}
This allows us to deduce that
\begin{align*}
\int dy_1 \sup_{y_2}\, \big|\cm^{\mathbf{2,1,2},(n)}_{t_1,t_2}(y_1+y_2,y_2)\big|&\lesssim \int_{t_1}^{t_2}\frac{ds_1}{|t_2-s_1|^{\frac12}}\int_{t_1}^{t_2}\frac{ds_2}{|t_2-s_2|^{\frac12}}\frac{1}{|s_1-s_2|^{\frac12}}\\
&\lesssim |t_2-t_1|^{\frac12}\int_{0}^{1}\frac{dr_1}{|1-r_1|^{\frac12}}\int_{0}^{1}\frac{dr_2}{|1-r_2|^{\frac12}}\frac{1}{|r_1-r_2|^{\frac12}}\lesssim |t_2-t_1|^{\frac12},
\end{align*}
uniformly over $0\leq t_1,t_2\leq 2$.

\

\noindent
\underline{\textit{Case of $\calt^{\mathbf{2,2},(n)}$}}. One has
\begin{align*}
\cm^{\mathbf{2,2},(n)}_{t_1,t_2}(y_1,y_2):=\mathbb{E}\Big[\calt^{\mathbf{2,2},(n)}_{t_1,t_2}(y_1)\calt^{\mathbf{2,2},(n)}_{t_1,t_2}(y_2)\Big]=\sum_{\mathbf{a}=1}^2 c_{\mathbf{a}}\cm^{\mathbf{2,2,a},(n)}_{t_1,t_2}(y_1,y_2),
\end{align*}
with
\begin{multline*}
\cm^{\mathbf{2,2,1},(n)}_{t_1,t_2}(y_1,y_2):=\int_{0}^{t_1}ds_1\int_{0}^{t_1}ds_2\int dw_1 dw_2 K_{t_1-s_1,t_2-s_1}(y_1,w_1) \cac^{(n)}_{t_2,s_1}(y_1,w_1)\\
\times K_{t_1-s_2,t_2-s_2}(y_2,w_2) \cac^{(n)}_{t_2,s_2}(y_2,w_2)\cac^{(n)}_{t_2,t_2}(y_1,y_2) \cac^{(n)}_{s_1,s_2}(w_1,w_2)^2
\end{multline*}
and
\begin{multline*}
\cm^{\mathbf{2,2,2},(n)}_{t_1,t_2}(y_1,y_2):=\int_{0}^{t_1}ds_1\int_{0}^{t_1}ds_2\int dw_1 dw_2 K_{t_1-s_1,t_2-s_1}(y_1,w_1) \cac^{(n)}_{t_2,s_1}(y_1,w_1)\\
 \times K_{t_1-s_2,t_2-s_2}(y_2,w_2) \cac^{(n)}_{t_2,s_2}(y_2,w_2)\cac^{(n)}_{t_2,s_2}(y_1,w_2) \cac^{(n)}_{t_2,s_1}(y_2,w_1)\cac^{(n)}_{s_1,s_2}(w_1,w_2).
\end{multline*}

\

\paragraph{Case of $\cm^{\mathbf{2,2,1},(n)}$}

One has
\begin{align*}
&\cm^{\mathbf{2,2,1},(n)}_{t_1,t_2}(y_1+y_2,y_2)\lesssim \cac^{(n),\star}_{t_2,t_2}(y_1)\int_{0}^{t_1}ds_1\int_{0}^{t_1}ds_2 \int dw_1\, K_{t_1-s_1,t_2-s_1}^\star(y_1-w_1) \cac^{(n),\star}_{t_2,s_1}(y_1-w_1) \\
& \hspace{2cm}\times \int dw_2\, K^\star_{t_1-s_2,t_2-s_2}(w_2) \cac^{(n),\star}_{t_2,s_2}(w_2) \cac^{(n),\star}_{s_1,s_2}(w_1-w_2)^2\\
&\lesssim \cac^{(n),\star}_{t_2,t_2}(y_1)\int_{0}^{t_1}\frac{ds_1}{|t_2-s_1|^{\frac12}}\int_{0}^{t_1}\frac{ds_2}{|t_2-s_2|^{\frac12}}  \int dw_1\, K_{t_1-s_1,t_2-s_1}^\star(y_1-w_1)\\
& \hspace{5cm} \times \int dw_2\, K^\star_{t_1-s_2,t_2-s_2}(w_2)  \cac^{(n),\star}_{s_1,s_2}(w_1-w_2)^2\\
&\lesssim \cac^{(n),\star}_{t_2,t_2}(y_1)\int_{0}^{t_1}\frac{ds_1}{|t_2-s_1|^{\frac12}}\int_{0}^{t_1}\frac{ds_2}{|t_2-s_2|^{\frac12}}  \Big[ K_{t_1-s_1,t_2-s_1}^\star \ast \big[K^\star_{t_1-s_2,t_2-s_2}\ast (\cac^{(n),\star}_{s_1,s_2})^2\big]\Big](y_1),
\end{align*}
and so
\begin{align*}
&\int dy_1 \sup_{y_2}\, \cm^{\mathbf{2,2,1},(n)}_{t_1,t_2}(y_1+y_2,y_2)\lesssim \\
&\lesssim \int_{0}^{t_1}\frac{ds_1}{|t_2-s_1|^{\frac12}}\int_{0}^{t_1}\frac{ds_2}{|t_2-s_2|^{\frac12}} \big\|\cac^{(n),\star}_{t_2,t_2}\big\|_{L^{\frac{3}{1+\eps}}} \Big\| K_{t_1-s_1,t_2-s_1}^\star \ast \big[K^\star_{t_1-s_2,t_2-s_2}\ast (\cac^{(n),\star}_{s_1,s_2})^2\big]\Big\|_{L^{\frac{3}{2-\eps}}}\\
&\lesssim \int_{0}^{t_1}\frac{ds_1}{|t_2-s_1|^{\frac12}}\int_{0}^{t_1}\frac{ds_2}{|t_2-s_2|^{\frac12}} \big\| K_{t_1-s_1,t_2-s_1}^\star\big\|_{L^1} \big\|K^\star_{t_1-s_2,t_2-s_2}\big\|_{L^1} \big\|(\cac^{(n),\star}_{s_1,s_2})^2\big\|_{L^{\frac{3}{2-\eps}}}\\
&\lesssim |t_2-t_1|^{2\theta}\int_{0}^{t_1}\frac{ds_1}{|t_2-s_1|^{\frac12}}\int_{0}^{t_1}\frac{ds_2}{|t_2-s_2|^{\frac12}} \frac{1}{|s_2-s_1|^{\frac{\eps}{2}}}\big\| K_{t_1-s_1,t_2-s_1}^\star\big\|_{L^1} \big\|K^\star_{t_1-s_2,t_2-s_2}\big\|_{L^1} \\
&\lesssim |t_2-t_1|^{2\theta} \int_{0}^{t_1}\frac{ds_1\, e^{-(t_1-s_1)}}{|t_2-s_1|^{\frac12-\frac32 \theta}|t_1-s_1|^{\frac52\theta}}\int_{0}^{t_1}\frac{ds_2\, e^{-(t_1-s_2)}}{|t_2-s_2|^{\frac12-\frac32 \theta}|t_1-s_2|^{\frac52\theta}|s_2-s_1|^{\frac{\eps}{2}}}.
\end{align*}
By picking $\theta:=\frac{1}{3}-\eps$, we obtain that
\begin{multline*}
\int dy_1 \sup_{y_2}\, \cm^{\mathbf{2,2,1},(n)}_{t_1,t_2}(y_1+y_2,y_2)\lesssim \\
\begin{aligned}
&\lesssim |t_2-t_1|^{\frac{2}{3}-2\eps}  \int_{0}^{t_1}\frac{ds_1\, e^{-(t_1-s_1)} }{|t_2-s_1|^{\frac32 \eps}|t_1-s_1|^{\frac{5}{6}-\frac52\eps}}\int_{0}^{t_1}\frac{ds_2\, e^{-(t_1-s_2)} }{|t_2-s_2|^{\frac32 \eps}|t_1-s_2|^{\frac{5}{6}-\frac52\eps}|s_2-s_1|^{\frac{\eps}{2}}} \\
&\lesssim |t_2-t_1|^{\frac{2}{3}-2\eps}  \int_{0}^{t_1}\frac{ds_1 \, e^{-(t_1-s_1)}  }{|t_1-s_1|^{\frac{5}{6}-\eps}}\int_{0}^{t_1}\frac{ds_2 \, e^{-(t_1-s_2)}  }{|t_1-s_2|^{\frac{5}{6}-\eps} |s_2-s_1|^{\frac{\eps}{2}}} \\
&\lesssim |t_2-t_1|^{\frac{2}{3}-2\eps} |t_2-t_1|^{\frac{2}{3}-2\eps}  \int_{0}^{+\infty}\frac{ds_1}{|s_1|^{\frac{5}{6}-\eps}}e^{-s_1}\int_{0}^{+\infty}\frac{ds_2\, e^{-s_2} }{|s_2|^{\frac{5}{6}-\eps}|s_2-s_1|^{\frac{\eps}{2}}}   \\
& \lesssim |t_2-t_1|^{\frac{2}{3}-2\eps},
\end{aligned}
\end{multline*}
uniformly over $0\leq t_1,t_2\leq 2$.

\

\paragraph{Case of $\cm^{\mathbf{2,2,2},(n)}$}

One has
\begin{align*}
&\cm^{\mathbf{2,2,2},(n)}_{t_1,t_2}(y_1+y_2,y_2)\lesssim \int_{0}^{t_1}ds_1\int_{0}^{t_1}ds_2\int dw_1\, K^\star_{t_1-s_1,t_2-s_1}(y_1-w_1) \cac^{(n),\star}_{t_2,s_1}(y_1-w_1)\cac^{(n),\star}_{t_2,s_1}(w_1)\\
&\hspace{2cm} \times \int dw_2\,  K^\star_{t_1-s_2,t_2-s_2}(w_2) \cac^{(n),\star}_{t_2,s_2}(w_2)\cac^{(n),\star}_{t_2,s_2}(y_1-w_2) \cac^{(n),\star}_{s_1,s_2}(w_1-w_2)\\
&\lesssim \int_{0}^{t_1}\frac{ds_1}{|t_2-s_1|^{\frac12}}\int_{0}^{t_1}\frac{ds_2}{|t_2-s_2|^{\frac12}}\int dw_1\, K^\star_{t_1-s_1,t_2-s_1}(y_1-w_1)\cac^{(n),\star}_{t_2,s_1}(w_1) \\
&\hspace{5cm}\times \int dw_2\,  K^\star_{t_1-s_2,t_2-s_2}(w_2) \cac^{(n),\star}_{t_2,s_2}(y_1-w_2)\cac^{(n),\star}_{s_1,s_2}(w_1-w_2) \\
&\lesssim \int_{0}^{t_1}\frac{ds_1}{|t_2-s_1|^{\frac12}}\int_{0}^{t_1}\frac{ds_2}{|t_2-s_2|^{\frac12}}\frac{1}{|s_2-s_1|^{\frac12}}\big[K^\star_{t_1-s_1,t_2-s_1}\ast \cac^{(n),\star}_{t_2,s_1}\big](y_1) \big[K^\star_{t_1-s_2,t_2-s_2}\ast \cac^{(n),\star}_{t_2,s_2}\big](y_1),
\end{align*}
and so, for every $\theta\in [0,1]$,
\begin{multline*}
\int dy_1 \sup_{y_2}\, \big|\cm^{\mathbf{2,2,2},(n)}_{t_1,t_2}(y_1+y_2,y_2)\big|\lesssim \\
\begin{aligned}
&\lesssim \int_{0}^{t_1}\frac{ds_1}{|t_2-s_1|^{\frac12}}\int_{0}^{t_1}\frac{ds_2}{|t_2-s_2|^{\frac12}}\frac{1}{|s_1-s_2|^{\frac12}}\big\|K^\star_{t_1-s_1,t_2-s_1} \ast \cac^{(n),\star}_{t_2,s_1}\big\|_{L^2} \big\|K^\star_{t_1-s_2,t_2-s_2}\ast \cac^{(n),\star}_{t_2,s_2}\big\|_{L^2} \\
&\lesssim \int_{0}^{t_1}\frac{ds_1}{|t_2-s_1|^{\frac12}}\int_{0}^{t_1}\frac{ds_2}{|t_2-s_2|^{\frac12}}\frac{1}{|s_1-s_2|^{\frac12}}\big\|K^\star_{t_1-s_1,t_2-s_1} \big\|_{L^1} \big\| \cac^{(n),\star}_{t_2,s_1}\big\|_{L^2} \big\|K^\star_{t_1-s_2,t_2-s_2}\big\|_{L^1} \big\| \cac^{(n),\star}_{t_2,s_2}\big\|_{L^2} \\
&\lesssim |t_2-t_1|^{2\theta}\int_{0}^{t_1}\frac{ds_1 \, e^{-(t_1-s_1)}}{|t_2-s_1|^{\frac12-\frac32\theta}|t_1-s_1|^{\frac52\theta}}\int_{0}^{t_1}\frac{ds_2 \, e^{-(t_1-s_2)} }{|t_2-s_2|^{\frac12-\frac32\theta}|t_1-s_2|^{\frac52\theta}}\frac{1}{|s_1-s_2|^{\frac12}}.
\end{aligned}
\end{multline*}
By picking $\theta:=\frac14-\eps$, we obtain that
\begin{multline*}
\int dy_1 \sup_{y_2}\, \cm^{\mathbf{2,2,2},(n)}_{t_1,t_2}(y_1+y_2,y_2)\lesssim \\
\begin{aligned}
&\lesssim |t_2-t_1|^{\frac12-2\eps} \int_{0}^{t_1}\frac{ds_1 \, e^{-(t_1-s_1)} }{|t_2-s_1|^{\frac18+\frac32\eps}|t_1-s_1|^{\frac58-\frac52\eps}}\int_{0}^{t_1}\frac{ds_2 \, e^{-(t_1-s_2)} }{|t_2-s_2|^{\frac18+\frac32\eps}|t_1-s_2|^{\frac58-\frac52\eps}|s_1-s_2|^{\frac12}}\\
&\lesssim |t_2-t_1|^{\frac12-2\eps} \int_{0}^{t_1}\frac{ds_1 \, e^{-(t_1-s_1)} }{|t_1-s_1|^{\frac34-\eps}}\int_{0}^{t_1}\frac{ds_2 \,e^{-(t_1-s_2)} }{|t_1-s_2|^{\frac34-\eps}|s_1-s_2|^{\frac12} }\\
&\lesssim |t_2-t_1|^{\frac12-2\eps} \int_{0}^{+\infty}\frac{ds_1\, e^{-s_1}}{|s_1|^{\frac34-\eps}}\int_{0}^{+\infty}\frac{ds_2\, }{|s_2|^{\frac34-\eps}|s_1-s_2|^{\frac12}}  \lesssim |t_2-t_1|^{\frac12-2\eps},
\end{aligned}
\end{multline*}
uniformly over $0\leq t_1,t_2\leq 2$.

\

\

\noindent
\underline{\textit{Case of $\calt^{\mathbf{2,3},(n)}$}}. One has

\begin{align*}
\cm^{\mathbf{2,3},(n)}_{t_1,t_2}(y_1,y_2):=\mathbb{E}\Big[\calt^{\mathbf{2,3},(n)}_{t_1,t_2}(y_1)\calt^{\mathbf{2,3},(n)}_{t_1,t_2}(y_2)\Big]=\sum_{\mathbf{a}=1}^2 c_{\mathbf{a}}\cm^{\mathbf{2,3,a},(n)}_{t_1,t_2}(y_1,y_2),
\end{align*}
with
\begin{multline*}
\cm^{\mathbf{2,3,1},(n)}_{t_1,t_2}(y_1,y_2):=\int_{0}^{t_1}ds_1\int_{0}^{t_1}ds_2\int dw_1 dw_2 K_{t_1-s_1}(y_1,w_1) \cac^{(n)}_{(t_1,t_2),s_1}(y_1,w_1)  \\
\times K_{t_1-s_2}(y_2,w_2) \cac^{(n)}_{(t_1,t_2),s_2}(y_2,w_2)\cac^{(n)}_{t_2,t_2}(y_1,y_2) \cac^{(n)}_{s_1,s_2}(w_1,w_2)^2
\end{multline*}
and
\begin{multline*}
\cm^{\mathbf{2,3,2},(n)}_{t_1,t_2}(y_1,y_2):=\int_{0}^{t_1}ds_1\int_{0}^{t_1}ds_2\int dw_1 dw_2 K_{t_1-s_1}(y_1,w_1) \cac^{(n)}_{(t_1,t_2),s_1}(y_1,w_1)  \\
 \times K_{t_1-s_2}(y_2,w_2) \cac^{(n)}_{(t_1,t_2),s_2}(y_2,w_2)\cac^{(n)}_{t_2,s_2}(y_1,w_2) \cac^{(n)}_{t_2,s_1}(y_2,w_1)\cac^{(n)}_{s_1,s_2}(w_1,w_2).
\end{multline*}

\

\paragraph{Case of $\cm^{\mathbf{2,3,1},(n)}$}

One has
\begin{multline*}
\cm^{\mathbf{2,3,1},(n)}_{t_1,t_2}(y_1+y_2,y_2)\lesssim \cac^{(n),\star}_{t_2,t_2}(y_1)\int_{0}^{t_1}ds_1\, \big\|\cac^{(n)}_{(t_1,t_2),s_1}\big\|_{L^\infty}\int_{0}^{t_1}ds_2\, \big\|\cac^{(n)}_{(t_1,t_2),s_2}\big\|_{L^\infty} \\
\begin{aligned}
&\hspace{2cm} \times \int dw_1\, K_{t_1-s_1}^\star(y_1-w_1) \int dw_2\, K^\star_{t_1-s_2}(w_2) \cac^{(n),\star}_{s_1,s_2}(w_1-w_2)^2\\
&\lesssim |t_2-t_1|^{2\theta}\cac^{(n),\star}_{t_2,t_2}(y_1)\int_{0}^{t_1}\frac{ds_1}{|t_1-s_1|^{\frac12+\theta}}\int_{0}^{t_1}\frac{ds_2}{|t_1-s_2|^{\frac12+\theta}}\,  \Big[ K_{t_1-s_1}^\star \ast \big[K^\star_{t_1-s_2}\ast (\cac^{(n),\star}_{s_1,s_2})^2\big]\Big](y_1),
\end{aligned}
\end{multline*}
for every $\theta\in [0,1]$. As a result,
\begin{multline*}
\int dy_1 \sup_{y_2}\, \cm^{\mathbf{2,3,1},(n)}_{t_1,t_2}(y_1+y_2,y_2)\lesssim \\
\begin{aligned}
&\lesssim |t_2-t_1|^{2\theta}\int_{0}^{t_1}\frac{ds_1}{|t_1-s_1|^{\frac12+\theta}}\int_{0}^{t_1}\frac{ds_2}{|t_1-s_2|^{\frac12+\theta}} \big\|\cac^{(n),\star}_{t_2,t_2}\big\|_{L^2} \Big\| K_{t_1-s_1}^\star \ast \big[K^\star_{t_1-s_2}\ast (\cac^{(n),\star}_{s_1,s_2})^2\big]\Big\|_{L^2}\\
&\lesssim |t_2-t_1|^{2\theta}\int_{0}^{t_1}\frac{ds_1}{|t_1-s_1|^{\frac12+\theta}}\int_{0}^{t_1}\frac{ds_2}{|t_1-s_2|^{\frac12+\theta}} \big\| K_{t_1-s_1}^\star\big\|_{L^1} \big\|K^\star_{t_1-s_2}\big\|_{L^1} \big\|(\cac^{(n),\star}_{s_1,s_2})^2\big\|_{L^2}\\
&\lesssim |t_2-t_1|^{2\theta}\int_{0}^{t_1}\frac{ds_1}{|t_1-s_1|^{\frac12+\theta}}e^{-(t_1-s_1)}\int_{0}^{t_1}\frac{ds_2}{|t_1-s_2|^{\frac12+\theta}} e^{-(t_1-s_2)}\frac{1}{|s_2-s_1|^{\frac14}} .
\end{aligned}
\end{multline*}
By picking $\theta:=\frac{3}{8}-\eps$, we obtain that
\begin{multline*}
\int dy_1 \sup_{y_2}\, \cm^{\mathbf{2,3,1},(n)}_{t_1,t_2}(y_1+y_2,y_2)\lesssim \\
\begin{aligned}
&\lesssim |t_2-t_1|^{\frac{3}{4}-2\eps} \int_{0}^{t_1}\frac{ds_1}{|t_1-s_1|^{\frac78-\eps}}e^{-(t_1-s_1)}\int_{0}^{t_1}\frac{ds_2}{|t_1-s_2|^{\frac78-\eps}} e^{-(t_1-s_2)}\frac{1}{|s_2-s_1|^{\frac14}}\\
&\lesssim |t_2-t_1|^{\frac{3}{4}-2\eps} \int_{0}^{+\infty}\frac{ds_1}{|s_1|^{\frac78-\eps}}e^{-s_1}\int_{0}^{+\infty}\frac{ds_2}{|s_2|^{\frac78-\eps}} e^{-s_2}\frac{1}{|s_2-s_1|^{\frac14}} \lesssim |t_2-t_1|^{\frac{3}{4}-2\eps},
\end{aligned}
\end{multline*}
uniformly over $0\leq t_1,t_2\leq 2$.

\

\paragraph{Case of $\cm^{\mathbf{2,3,2},(n)}$}

One has
\begin{align*}
&\cm^{\mathbf{2,3,2},(n)}_{t_1,t_2}(y_1+y_2,y_2)\lesssim \int_{0}^{t_1}ds_1\int_{0}^{t_1}ds_2\int dw_1 dw_2\,  K^\star_{t_1-s_1}(y_1-w_1) \cac^{(n),\star}_{(t_1,t_2),s_1}(y_1-w_1)\\
&\hspace{1cm} \times K^\star_{t_1-s_2}(w_2) \cac^{(n),\star}_{(t_1,t_2),s_2}(w_2)\cac^{(n),\star}_{t_2,s_2}(y_1-w_2) \cac^{(n),\star}_{t_2,s_1}(w_1)\cac^{(n),\star}_{s_1,s_2}(w_1-w_2)\\
&\lesssim |t_2-t_1|^{2\theta}\int_{0}^{t_1}\frac{ds_1}{|t_1-s_1|^{\frac12+\theta}}\int_{0}^{t_1}\frac{ds_2}{|t_1-s_2|^{\frac12+\theta}}\frac{1}{|s_2-s_1|^{\frac12}}\int dw_1\, K^\star_{t_1-s_1}(y_1-w_1) \cac^{(n),\star}_{t_2,s_1}(w_1)\\
&\hspace{4cm}\times \int dw_2\,   K^\star_{t_1-s_2}(w_2) \cac^{(n),\star}_{t_2,s_2}(y_1-w_2)\\
&\lesssim |t_2-t_1|^{2\theta}\int_{0}^{t_1}\frac{ds_1}{|t_1-s_1|^{\frac12+\theta}}\int_{0}^{t_1}\frac{ds_2}{|t_1-s_2|^{\frac12+\theta}}\frac{1}{|s_2-s_1|^{\frac12}}\big[ K^\star_{t_1-s_1}\ast \cac^{(n),\star}_{t_2,s_1}\big](y_1) \big[K^\star_{t_1-s_2}\ast\cac^{(n),\star}_{t_2,s_2}\big](y_1),
\end{align*}
for every $\theta\in [0,1]$. As a result,
\begin{align*}
&\int dy_1 \sup_{y_2}\, \big|\cm^{\mathbf{2,3,2},(n)}_{t_1,t_2}(y_1+y_2,y_2)\big|\lesssim \\
&\lesssim |t_2-t_1|^{2\theta}\int_{0}^{t_1}\frac{ds_1}{|t_1-s_1|^{\frac12+\theta}}\int_{0}^{t_1}\frac{ds_2}{|t_1-s_2|^{\frac12+\theta}}\frac{1}{|s_2-s_1|^{\frac12}}\big\| K^\star_{t_1-s_1}\ast \cac^{(n),\star}_{t_2,s_1}\big\|_{L^2} \big\| K^\star_{t_1-s_2}\ast\cac^{(n),\star}_{t_2,s_2}\big\|_{L^2}\\
&\lesssim |t_2-t_1|^{2\theta}\int_{0}^{t_1}\frac{ds_1}{|t_1-s_1|^{\frac12+\theta}}\int_{0}^{t_1}\frac{ds_2}{|t_1-s_2|^{\frac12+\theta}}\frac{1}{|s_2-s_1|^{\frac12}}\big\| K^\star_{t_1-s_1}\big\|_{L^1}\big\|  \cac^{(n),\star}_{t_2,s_1}\big\|_{L^2}  \big\| K^\star_{t_1-s_2}\big\|_{L^1}\big\| \cac^{(n),\star}_{t_2,s_2}\big\|_{L^2}\\
&\lesssim |t_2-t_1|^{2\theta}\int_{0}^{t_1}\frac{ds_1}{|t_1-s_1|^{\frac12+\theta}}e^{-(t_1-s_1)}\int_{0}^{t_1}\frac{ds_2}{|t_1-s_2|^{\frac12+\theta}}e^{-(t_1-s_2)}\frac{1}{|s_2-s_1|^{\frac12}}.
\end{align*}
By choosing $\theta:=\frac14-\eps$, we deduce
\begin{multline*}
\int dy_1 \sup_{y_2}\, \big|\cm^{\mathbf{2,3,2},(n)}_{t_1,t_2}(y_1+y_2,y_2)\big|\lesssim \\
\begin{aligned}
&\lesssim  |t_2-t_1|^{\frac12-2\eps}\int_{0}^{t_1}\frac{ds_1}{|t_1-s_1|^{\frac34-\eps}}e^{-(t_1-s_1)}\int_{0}^{t_1}\frac{ds_2}{|t_1-s_2|^{\frac34-\eps}}e^{-(t_1-s_2)}\frac{1}{|s_2-s_1|^{\frac12}}\\
&\lesssim |t_2-t_1|^{\frac12-2\eps}\int_{0}^{+\infty}\frac{ds_1}{|s_1|^{\frac34-\eps}}e^{-s_1}\int_{0}^{+\infty}\frac{ds_2}{|s_2|^{\frac34-\eps}}e^{-s_2}\frac{1}{|s_2-s_1|^{\frac12}}\\
&\lesssim |t_2-t_1|^{\frac12-2\eps},
\end{aligned}
\end{multline*}
uniformly over $0\leq t_1,t_2\leq 2$.

\

\

\noindent
\underline{\textit{Case of $\calt^{\mathbf{2,4},(n)}$}}. One has

\begin{align*}
\cm^{\mathbf{2,4},(n)}_{t_1,t_2}(y_1,y_2):=\mathbb{E}\Big[\calt^{\mathbf{2,4},(n)}_{t_1,t_2}(y_1)\calt^{\mathbf{2,4},(n)}_{t_1,t_2}(y_2)\Big]=\sum_{\mathbf{a}=1}^2 c_{\mathbf{a}}\cm^{\mathbf{2,4,a},(n)}_{t_1,t_2}(y_1,y_2),
\end{align*}
with (recall the notation introduced in \eqref{inc-cn})
\begin{multline*}
\cm^{\mathbf{2,4,1},(n)}_{t_1,t_2}(y_1,y_2):=\int_{0}^{t_1}ds_1\int_{0}^{t_1}ds_2\int dw_1 dw_2  K_{t_1-s_1}(y_1,w_1) \cac^{(n)}_{t_1,s_1}(y_1,w_1)\\
\times K_{t_1-s_2}(y_2,w_2) \cac^{(n)}_{t_1,s_2}(y_2,w_2)\cac^{(n)}_{(t_1,t_2),(t_1,t_2)}(y_1,y_2) \cac^{(n)}_{s_1,s_2}(w_1,w_2)^2
\end{multline*}
and
\begin{multline*}
 \cm^{\mathbf{2,4,2},(n)}_{t_1,t_2}(y_1,y_2):=\int_{0}^{t_1}ds_1\int_{0}^{t_1}ds_2\int dw_1 dw_2 K_{t_1-s_1}(y_1,w_1) \cac^{(n)}_{t_1,s_1}(y_1,w_1)K_{t_1-s_2}(y_2,w_2) \\
  \times \cac^{(n)}_{t_1,s_2}(y_2,w_2)\cac^{(n)}_{(t_1,t_2),s_2}(y_1,w_2) \cac^{(n)}_{(t_1,t_2),s_1}(y_2,w_1)\cac^{(n)}_{s_1,s_2}(w_1,w_2).
\end{multline*}

\

\paragraph{Case of $\cm^{\mathbf{2,4,1},(n)}$}

One has
\begin{align*}
&\big| \cm^{\mathbf{2,4,1},(n)}_{t_1,t_2}(y_1+y_2,y_2)\big|\lesssim \cac^{(n),\star}_{(t_1,t_2),(t_1,t_2)}(y_1)\int_{0}^{t_1}ds_1\int_{0}^{t_1}ds_2\int dw_1 dw_2\,  K^\star_{t_1-s_1}(y_1-w_1) \\
&\hspace{1cm}\times \cac^{(n),\star}_{t_1,s_1}(y_1-w_1)K^\star_{t_1-s_2}(w_2) \cac^{(n),\star}_{t_1,s_2}(w_2) \cac^{(n),\star}_{s_1,s_2}(w_1-w_2)^2\\
&\lesssim \cac^{(n),\star}_{(t_1,t_2),(t_1,t_2)}(y_1)\int_{0}^{t_1}\frac{ds_1}{|t_1-s_1|^{\frac12}}\int_{0}^{t_1}\frac{ds_2}{|t_1-s_2|^{\frac12}}\Big[K^\star_{t_1-s_1}\ast \big[K^\star_{t_1-s_2}\ast  (\cac^{(n),\star}_{s_1,s_2})^2\big]\Big](y_1).
\end{align*}
Thus,
\begin{align*}
&\int dy_1 \sup_{y_2} \big| \cm^{\mathbf{2,4,1},(n)}_{t_1,t_2}(y_1+y_2,y_2)\big|\lesssim \\
&\lesssim \big\|\cac^{(n),\star}_{(t_1,t_2),(t_1,t_2)}\big\|_{L^{\frac{1}{1-\eps}}} \int_{0}^{t_1}\frac{ds_1}{|t_1-s_1|^{\frac12}}\int_{0}^{t_1}\frac{ds_2}{|t_1-s_2|^{\frac12}}\Big\|K^\star_{t_1-s_1}\ast \big[K^\star_{t_1-s_2}\ast  (\cac^{(n),\star}_{s_1,s_2})^2\big]\Big\|_{L^{\frac{1}{\eps}}}\\
&\lesssim \big\|\cac^{(n),\star}_{(t_1,t_2),(t_1,t_2)}\big\|_{L^1}^{1-2\eps} \big\|\cac^{(n),\star}_{(t_1,t_2),(t_1,t_2)}\big\|_{L^2}^{2\eps}\int_{0}^{t_1}\frac{ds_1}{|t_1-s_1|^{\frac12}}\int_{0}^{t_1}\frac{ds_2}{|t_1-s_2|^{\frac12}}\big\|K^\star_{t_1-s_1}\big\|_{L^1}\big\|K^\star_{t_1-s_2}\big\|_{L^1}\big\| \cac^{(n),\star}_{s_1,s_2}\big\|_{L^{\frac{2}{\eps}}}^2\\
&\lesssim |t_2-t_1|^{1-2\eps} \int_{0}^{t_1}\frac{ds_1}{|t_1-s_1|^{\frac12}}e^{-(t_1-s_1)}\int_{0}^{t_1}\frac{ds_2}{|t_1-s_2|^{\frac12}}e^{-(t_1-s_2)}\frac{1}{|s_1-s_2|^{1-\frac32\eps}}\\
&\lesssim |t_2-t_1|^{1-2\eps} \int_{0}^{+\infty}\frac{ds_1}{|s_1|^{\frac12}}e^{-s_1}\int_{0}^{+\infty}\frac{ds_2}{|s_2|^{\frac12}}e^{-s_2}\frac{1}{|s_1-s_2|^{1-\frac32\eps}}\\
&\lesssim |t_2-t_1|^{1-2\eps},
\end{align*}
uniformly over $0\leq t_1,t_2\leq 2$.

\

\paragraph{Case of $\cm^{\mathbf{2,4,2},(n)}$}

One has
\begin{multline*}
\big|\cm^{\mathbf{2,4,2},(n)}_{t_1,t_2}(y_1+y_2,y_2)\big|\lesssim \int_{0}^{t_1}ds_1\int_{0}^{t_1}ds_2\int dw_1 dw_2\\
  \times K^\star_{t_1-s_1}(y_1-w_1) \cac^{(n),\star}_{t_1,s_1}(y_1-w_1)K^\star_{t_1-s_2}(w_2) \\
 \hspace{4cm} \times \cac^{(n),\star}_{t_1,s_2}(w_2)\cac^{(n),\star}_{(t_1,t_2),s_2}(y_1-w_2) \cac^{(n),\star}_{(t_1,t_2),s_1}(w_1)\cac^{(n),\star}_{s_1,s_2}(w_1-w_2)\\
\lesssim \int_{0}^{t_1}\frac{ds_1}{|t_1-s_1|^{\frac12}}\int_{0}^{t_1}\frac{ds_2}{|t_1-s_2|^{\frac12}}\frac{1}{|s_2-s_1|^{\frac12}}\big[ K^\star_{t_1-s_1}\ast\cac^{(n),\star}_{(t_1,t_2),s_1}\big](y_1)\big[ K^\star_{t_1-s_2}\ast \cac^{(n),\star}_{(t_1,t_2),s_2}\big](y_1) ,
\end{multline*}
and so
\begin{multline*}
\int dy_1 \sup_{y_2} \big|\cm^{\mathbf{2,4,2},(n)}_{t_1,t_2}(y_1+y_2,y_2)\big|\lesssim \\
\begin{aligned}
&\lesssim \int_{0}^{t_1}\frac{ds_1}{|t_1-s_1|^{\frac12}}\int_{0}^{t_1}\frac{ds_2}{|t_1-s_2|^{\frac12}}\frac{1}{|s_2-s_1|^{\frac12}}\big\| K^\star_{t_1-s_1}\ast\cac^{(n),\star}_{(t_1,t_2),s_1}\big\|_{L^2}\big\| K^\star_{t_1-s_2}\ast \cac^{(n),\star}_{(t_1,t_2),s_2}\big\|_{L^2}\\ 
&\lesssim \int_{0}^{t_1}\frac{ds_1}{|t_1-s_1|^{\frac12}}e^{-(t_1-s_1)}\int_{0}^{t_1}\frac{ds_2}{|t_1-s_2|^{\frac12}}e^{-(t_1-s_2)}\frac{1}{|s_2-s_1|^{\frac12}}\big\|\cac^{(n),\star}_{(t_1,t_2),s_1}\big\|_{L^2}\big\| \cac^{(n),\star}_{(t_1,t_2),s_2}\big\|_{L^2}\\
&\lesssim |t_2-t_1|^{\frac12} \int_{0}^{t_1}\frac{ds_1}{|t_1-s_1|^{\frac12}}e^{-(t_1-s_1)}\int_{0}^{t_1}\frac{ds_2}{|t_1-s_2|^{\frac12}}e^{-(t_1-s_2)}\frac{1}{|s_2-s_1|^{\frac12}}\\
&\lesssim |t_2-t_1|^{\frac12} \int_{0}^{+\infty}\frac{ds_1}{|s_1|^{\frac12}}e^{-s_1}\int_{0}^{+\infty}\frac{ds_2}{|s_2|^{\frac12}}e^{-s_2}\frac{1}{|s_2-s_1|^{\frac12}}\\
&\lesssim |t_2-t_1|^{\frac12},
\end{aligned}
\end{multline*}
uniformly over $0\leq t_1,t_2\leq 2$.


\

\subsubsection{Increments: third diagram}

Recall the definition
\begin{align*}
\calt^{\mathbf{3},(n)}_t(y)&:=\int_{0}^t ds \int dw \, K_{t-s}(y,w)\cac^{(n)}_{t,s}(y,w)^2 \big(\<Psi>^{(n)}_{0,s}(w)-\<Psi>^{(n)}_{0,t}(w)\big).
\end{align*}

\

\

For $0\leq t_1<t_2$, one has
\begin{align*}
\calt^{\mathbf{3},(n)}_{t_1,t_2}(y)&:=\calt^{\mathbf{3},(n)}_{t_2}(y)-\calt^{\mathbf{3},(n)}_{t_1}(y)\\
&=\calt^{\mathbf{3,1},(n)}_{t_1,t_2}(y)+\calt^{\mathbf{3,2},(n)}_{t_1,t_2}(y)+\calt^{\mathbf{3,3},(n)}_{t_1,t_2}(y)+\calt^{\mathbf{3,4},(n)}_{t_1,t_2}(y)+\calt^{\mathbf{3,5},(n)}_{t_1,t_2}(y),
\end{align*}
with
\begin{align*}
\calt^{\mathbf{3,1},(n)}_{t_1,t_2}(y)&:=\int_{t_1}^{t_2} ds \int dw \, K_{t_2-s}(y,w)\cac^{(n)}_{s,t_2}(y,w)^2 \big(\<Psi>^{(n)}_{0,s}(w)-\<Psi>^{(n)}_{0,t_2}(w)\big),
\end{align*}

\begin{align*}
\calt^{\mathbf{3,2},(n)}_{t_1,t_2}(y)&:=\int_{0}^{t_1} ds \int dw \, K_{t_2-s}(y,w)\cac^{(n)}_{s,t_2}(y,w)^2 \big((\<Psi>^{(n)}_{0,s}(w)-\<Psi>^{(n)}_{0,t_2}(w))-(\<Psi>^{(n)}_{0,s}(w)-\<Psi>^{(n)}_{0,t_1}(w))\big),
\end{align*}

\begin{align*}
\calt^{\mathbf{3,3},(n)}_{t_1,t_2}(y)&:=\int_{0}^{t_1} ds \int dw \, K_{t_1-s,t_2-s}(y,w)\cac^{(n)}_{s,t_2}(y,w)^2 \big(\<Psi>^{(n)}_{0,s}(w)-\<Psi>^{(n)}_{0,t_1}(w)\big),
\end{align*}

\begin{align*}
\calt^{\mathbf{3,4},(n)}_{t_1,t_2}(y)&:=\int_{0}^{t_1} ds \int dw \, K_{t_1-s}(y,w)\cac^{(n)}_{s,(t_1,t_2)}(y,w)\cac^{(n)}_{s,t_2}(y,w) \big(\<Psi>^{(n)}_{0,s}(w)-\<Psi>^{(n)}_{0,t_1}(w)\big),
\end{align*}

\begin{align*}
\calt^{\mathbf{3,5},(n)}_{t_1,t_2}(y)&:=\int_{0}^{t_1} ds \int dw \, K_{t_1-s}(y,w)\cac^{(n)}_{s,t_1}(y,w)\cac^{(n)}_{s,(t_1,t_2)}(y,w) \big(\<Psi>^{(n)}_{0,s}(w)-\<Psi>^{(n)}_{0,t_1}(w)\big).
\end{align*}

\

\noindent
\underline{\textit{Case of $\calt^{\mathbf{3,1},(n)}$}}. One has
\begin{multline*}
\cm^{\mathbf{3,1},(n)}_{t_1,t_2}(y_1,y_2)=\int_{t_1}^{t_2} ds_1 \int_{t_1}^{t_2} ds_2 \int dw_1 \int dw_2 \, K_{t_2-s_1}(y_1,w_1)\cac^{(n)}_{t_2,s_1}(y_1,w_1)^2 \\
\times K_{t_2-s_2}(y_2,w_2)\cac^{(n)}_{t_2,s_2}(y_2,w_2)^2\cac^{(n)}_{(t_2,s_1),(t_2,s_2)}(w_1,w_2).
\end{multline*}
As a result,
\begin{multline*}
\int dy_1 \sup_{y_2}\big|\cm^{\mathbf{3,1},(n)}_{t_1,t_2}(y_1+y_2,y_2)\big|\lesssim \int_{t_1}^{t_2} ds_1 \int_{t_1}^{t_2} ds_2 \int dy_1\int dw_1 \int dw_2 \, K^\star_{t_2-s_1}(y_1-w_1) \\
\begin{aligned}
&\hspace{3cm} \times \big(\cac^{(n),\star}_{t_2,s_1}(y_1-w_1)\big)^2 K^\star_{t_2-s_2}(w_2)\big(\cac^{(n),\star}_{t_2,s_2}(w_2)\big)^2\cac^{(n),\star}_{(t_2,s_1),(t_2,s_2)}(w_1-w_2)\\
&\lesssim \int_{t_1}^{t_2} ds_1 \int_{t_1}^{t_2} ds_2 \, \big\|K^\star_{t_2-s_1}(\cac^{(n),\star}_{t_2,s_1})^2\big\|_{L^1} \big\|\cac^{(n),\star}_{(t_2,s_1),(t_2,s_2)} \ast (K^\star_{t_2-s_2}(\cac^{(n),\star}_{t_2,s_2})^2)\big\|_{L^1}\\
&\lesssim \int_{t_1}^{t_2} ds_1 \int_{t_1}^{t_2} ds_2 \, \big\|K^\star_{t_2-s_1}\big\|_{L^1}\big\|(\cac^{(n),\star}_{t_2,s_1})^2\big\|_{L^\infty} \big\|\cac^{(n),\star}_{(t_2,s_1),(t_2,s_2)} \big\|_{L^1}\big\|K^\star_{t_2-s_2}\big\|_{L^1} \big\|(\cac^{(n),\star}_{t_2,s_2})^2\big\|_{L^\infty}\\
&\lesssim \int_{t_1}^{t_2} \frac{ds_1}{|t_2-s_1|} e^{-(t_2-s_1)}\int_{t_1}^{t_2} \frac{ds_2}{|t_2-s_2|} e^{-(t_2-s_2)}  \big\|\cac^{(n),\star}_{(t_2,s_1),(t_2,s_2)} \big\|_{L^1}\\
&\lesssim \int_{t_1}^{t_2} \frac{ds_1}{|t_2-s_1|^{\frac12}} e^{-(t_2-s_1)}\int_{t_1}^{t_2} \frac{ds_2}{|t_2-s_2|^{\frac12}}e^{-(t_2-s_2)},
\end{aligned}
\end{multline*}
from which we obtain 
\begin{align*}
\int dy_1 \sup_{y_2}\big|\cm^{\mathbf{3,1},(n)}_{t_1,t_2}(y_1+y_2,y_2)\big|&\lesssim \int_{t_1}^{t_2} \frac{ds_1}{|t_2-s_1|^{\frac12}} \int_{t_1}^{t_2} \frac{ds_2}{|t_2-s_2|^{\frac12}}\lesssim |t_2-t_1|,
\end{align*}
uniformly over $0\leq t_1,t_2\leq 2$.

\

\noindent
\underline{\textit{Case of $\calt^{\mathbf{3,2},(n)}$}}. One has
\begin{multline*}
\cm^{\mathbf{3,2},(n)}_{t_1,t_2}(y_1,y_2):=\int_{0}^{t_1} ds_1 \int_{0}^{t_1} ds_2 \int dw_1 \int dw_2 \, K_{t_2-s_1}(y_1,w_1)\\
\times \cac^{(n)}_{t_2,s_1}(y_1,w_1)^2 K_{t_2-s_2}(y_2,w_2)\cac^{(n)}_{t_2,s_2}(y_2,w_2)^2\cac^{(n)}_{(t_1,t_2),(t_1,t_2)}(w_1,w_2).
\end{multline*}
Therefore
\begin{multline*}
\int dy_1 \sup_{y_2}\big|\cm^{\mathbf{3,2},(n)}_{t_1,t_2}(y_1+y_2,y_2)\big|\lesssim \int_{0}^{t_1} ds_1 \int_{0}^{t_1} ds_2  \int dy_1\int dw_1 \int dw_2 \, K^\star_{t_2-s_1}(y_1-w_1)\\
\begin{aligned}
&\hspace{3.5cm}\times \big(\cac^{(n),\star}_{t_2,s_1}(y_1-w_1)\big)^2 K^\star_{t_2-s_2}(w_2)\big(\cac^{(n),\star}_{t_2,s_2}(w_2)\big)^2\cac^{(n),\star}_{(t_1,t_2),(t_1,t_2)}(w_1-w_2)\\
&\lesssim \int_{0}^{t_1} ds_1 \int_{0}^{t_1} ds_2\, \big\|K^\star_{t_2-s_1}(\cac^{(n),\star}_{t_2,s_1})^2\big\|_{L^1} \big\|\cac^{(n),\star}_{(t_1,t_2),(t_1,t_2)} \ast (K^\star_{t_2-s_2}(\cac^{(n),\star}_{t_2,s_2})^2)\big\|_{L^1}\\
&\lesssim \int_{0}^{t_1} ds_1 \int_{0}^{t_1} ds_2 \, \big\|K^\star_{t_2-s_1}\big\|_{L^1}\big\|(\cac^{(n),\star}_{t_2,s_1})^2\big\|_{L^\infty} \big\|\cac^{(n),\star}_{(t_1,t_2),(t_1,t_2)} \big\|_{L^1}\big\|K^\star_{t_2-s_2}\big\|_{L^1} \big\|(\cac^{(n),\star}_{t_2,s_2})^2\big\|_{L^\infty}\\
&\lesssim \big\|\cac^{(n),\star}_{(t_1,t_2),(t_1,t_2)} \big\|_{L^1}\int_{0}^{t_1} \frac{ds_1}{|t_2-s_1|}e^{-(t_2-s_1)} \int_{0}^{t_1} \frac{ds_2}{|t_2-s_2|}e^{-(t_2-s_2)}  \\
&\lesssim |t_2-t_1| \bigg(\int_{0}^{t_1} \frac{ds_1}{|t_2-t_1|^{\eps} |t_1-s_1|^{1-\eps}}e^{-(t_1-s_1)}\bigg)^2 \\
&\lesssim |t_2-t_1|^{1-2\eps}\bigg(\int_{0}^{t_1} \frac{ds_1}{ |t_1-s_1|^{1-\eps}}e^{-(t_1-s_1)}\bigg)^2\\
&\lesssim |t_2-t_1|^{1-2\eps}\bigg(\int_{0}^{+\infty} \frac{ds_1}{ |s_1|^{1-\eps}}e^{-s_1}\bigg)^2\lesssim |t_2-t_1|^{1-2\eps}
\end{aligned}
\end{multline*}
for any small $\eps>0$, uniformly over $0\leq t_1,t_2\leq 2$.

\

\noindent
\underline{\textit{Case of $\calt^{\mathbf{3,3},(n)}$}}. One has
\begin{multline*}
\cm^{\mathbf{3,3},(n)}_{t_1,t_2}(y_1,y_2):=\int_{0}^{t_1} ds_1 \int_{0}^{t_1} ds_2 \int dw_1 \int dw_2\,  K_{t_1-s_1,t_2-s_1}(y_1,w_1) \cac^{(n)}_{t_2,s_1}(y_1,w_1)^2\\
 \times K_{t_1-s_2,t_2-s_2}(y_2,w_2)\cac^{(n)}_{t_2,s_2}(y_2,w_2)^2\cac^{(n)}_{(t_1,s_1),(t_1,s_2)}(w_1,w_2).
\end{multline*}
Therefore, for any $\theta\in [0,1]$,
\begin{multline*}
\int dy_1 \sup_{y_2}\big|\cm^{\mathbf{3,3},(n)}_{t_1,t_2}(y_1+y_2,y_2)\big|\lesssim \int_{0}^{t_1} ds_1 \int_{0}^{t_1} ds_2 \int dy_1\int dw_1 \int dw_2  K_{t_1-s_1,t_2-s_1}^\star(y_1-w_1)\\
\begin{aligned}
&\hspace{3cm}\times \big(\cac^{(n),\star}_{t_2,s_1}(y_1-w_1)\big)^2 K_{t_1-s_2,t_2-s_2}^\star(w_2)\big(\cac^{(n),\star}_{t_2,s_2}(w_2)\big)^2\cac^{(n),\star}_{(t_1,s_1),(t_1,s_2)}(w_1-w_2)\\
&\lesssim \int_{0}^{t_1} ds_1 \int_{0}^{t_1} ds_2\, \big\|K_{t_1-s_1,t_2-s_1}^\star(\cac^{(n),\star}_{t_2,s_1})^2\big\|_{L^1} \big\|\cac^{(n),\star}_{(t_1,s_1),(t_1,s_2)} \ast (K_{t_1-s_2,t_2-s_2}^\star(\cac^{(n),\star}_{t_2,s_2})^2)\big\|_{L^1}\\
&\lesssim \int_{0}^{t_1} ds_1 \int_{0}^{t_1} ds_2 \, \big\|K_{t_1-s_1,t_2-s_1}^\star\big\|_{L^1}\big\|(\cac^{(n),\star}_{t_2,s_1})^2\big\|_{L^\infty} \big\|\cac^{(n),\star}_{(t_1,s_1),(t_1,s_2)} \big\|_{L^1}\big\|K_{t_1-s_2,t_2-s_2}^\star\big\|_{L^1} \big\|(\cac^{(n),\star}_{t_2,s_2})^2\big\|_{L^\infty}\\
&\lesssim |t_2-t_1|^{2\theta}\int_{0}^{t_1} \frac{ds_1}{|t_2-s_1|^{1-\frac32\theta}}|t_1-s_1|^{\frac12-\frac52 \theta}e^{-(t_1-s_1)} \int_{0}^{t_1} \frac{ds_2}{|t_2-s_2|^{1-\frac32\theta}}  |t_1-s_2|^{\frac12-\frac52 \theta}e^{-(t_1-s_2)} \\
&\lesssim |t_2-t_1|^{2\theta}\bigg(\int_{0}^{t_1} \frac{ds}{|t_2-s|^{1-\frac32\theta}}|t_1-s|^{\frac12-\frac52 \theta}e^{-(t_1-s)}\bigg)^2
\end{aligned}
\end{multline*}
and from here we can pick for instance $\theta=\frac12-\eps$ to assert that
\begin{align*}
\int dy_1 \sup_{y_2}\big|\cm^{\mathbf{3,3},(n)}_{t_1,t_2}(y_1+y_2,y_2)\big|&\lesssim  |t_2-t_1|^{1-2\eps}\bigg(\int_{0}^{t_1} \frac{ds}{|t_2-s|^{\frac14+\frac32\eps}}|t_1-s|^{-\frac34+\frac52 \eps}e^{-(t_1-s)}\bigg)^2\\
&\lesssim  |t_2-t_1|^{1-2\eps}\bigg(\int_{0}^{t_1} \frac{ds}{|t_1-s|^{1-\eps}}e^{-(t_1-s)}\bigg)^2\\
&\lesssim  |t_2-t_1|^{1-2\eps}\bigg(\int_{0}^{+\infty} \frac{ds}{|s|^{1-\eps}}e^{-s}\bigg)^2 \lesssim |t_2-t_1|^{1-2\eps},
\end{align*}
for every small $\eps>0$, uniformly over $0\leq t_1,t_2\leq 2$.

\

\noindent
\underline{\textit{Case of $\calt^{\mathbf{3,4},(n)}$}}. One has
\begin{multline*}
\cm^{\mathbf{3,4},(n)}_{t_1,t_2}(y_1,y_2):=\int_{0}^{t_1} ds_1 \int_{0}^{t_1} ds_2\int dw_1 \int dw_2   K_{t_1-s_1}(y_1,w_1)\cac^{(n)}_{s_1,(t_1,t_2)}(y_1,w_1)\cac^{(n)}_{s_1,t_2}(y_1,w_1)\\
\times K_{t_1-s_2}(y_2,w_2)\cac^{(n)}_{s_2,(t_1,t_2)}(y_2,w_2)\cac^{(n)}_{s_2,t_2}(y_2,w_2)\cac^{(n)}_{(t_1,s_1),(t_1,s_2)}(w_1,w_2).
\end{multline*}
As a result,
\begin{multline*}
\int dy_1 \sup_{y_2}\big|\cm^{\mathbf{3,4},(n)}_{t_1,t_2}(y_1+y_2,y_2)\big|\lesssim \int_{0}^{t_1} ds_1 \int_{0}^{t_1} ds_2 \int dy_1\int dw_1 \int dw_2 K^\star_{t_1-s_1}(y_1-w_1)\cac^{(n),\star}_{s_1,(t_1,t_2)}(y_1-w_1)\\
\begin{aligned}
&\hspace{3cm}\times  \cac^{(n),\star}_{s_1,t_2}(y_1-w_1) K^\star_{t_1-s_2}(w_2)\cac^{(n),\star}_{s_2,(t_1,t_2)}(w_2)\cac^{(n),\star}_{s_2,t_2}(w_2)\cac^{(n),\star}_{(t_1,s_1),(t_1,s_2)}(w_1-w_2)\\
&\lesssim \int_{0}^{t_1} ds_1 \int_{0}^{t_1} ds_2\, \big\|K^\star_{t_1-s_1}\cac^{(n),\star}_{s_1,(t_1,t_2)}\cac^{(n),\star}_{s_1,t_2}\big\|_{L^1} \big\|\cac^{(n),\star}_{(t_1,s_1),(t_1,s_2)} \ast (K^\star_{t_1-s_2}\cac^{(n),\star}_{s_2,(t_1,t_2)}\cac^{(n),\star}_{s_2,t_2})\big\|_{L^1}\\
&\lesssim \int_{0}^{t_1} ds_1 \int_{0}^{t_1} ds_2 \, \big\|K^\star_{t_1-s_1}\big\|_{L^1}\big\|\cac^{(n),\star}_{s_1,(t_1,t_2)}\big\|_{L^\infty} \big\|\cac^{(n),\star}_{s_1,t_2}\big\|_{L^\infty} \big\|\cac^{(n),\star}_{(t_1,s_1),(t_1,s_2)} \big\|_{L^1}\big\|K^\star_{t_1-s_2}\big\|_{L^1} \big\|\cac^{(n),\star}_{s_2,(t_1,t_2)}\big\|_{L^\infty}\big\|\cac^{(n),\star}_{s_2,t_2}\big\|_{L^\infty}\\
&\lesssim |t_2-t_1|^{2\theta}\bigg(\int_{0}^{t_1} ds  \, \frac{e^{-(t_1-s)}}{|t_1-s|^{\frac12+\theta}}\bigg)^2
\end{aligned}
\end{multline*}
for every $\theta\in [0,1]$. By choosing $\theta:=\frac12-\eps$, we obtain
\begin{align*}
&\int dy_1 \sup_{y_2}\big|\cm^{\mathbf{3,4},(n)}_{t_1,t_2}(y_1+y_2,y_2)\big|\lesssim |t_2-t_1|^{1-2\eps}\bigg(\int_{0}^{+\infty}ds\, \frac{e^{-s}}{|s|^{1-\eps}}\bigg)^2\lesssim |t_2-t_1|^{1-2\eps},
\end{align*}
for every $\eps >0$, uniformly over $0\leq t_1,t_2\leq 2$.

\

\noindent
\underline{\textit{Case of $\calt^{\mathbf{3,5},(n)}$}}. This quantity can be controlled with the same arguments as for $\calt^{\mathbf{3,4},(n)}$, leading to the same bound
\begin{align*}
&\int dy_1 \sup_{y_2}\big|\cm^{\mathbf{3,5},(n)}_{t_1,t_2}(y_1+y_2,y_2)\big|\lesssim |t_2-t_1|^{1-2\eps},
\end{align*}
for every $\eps >0$, uniformly over $0\leq t_1,t_2\leq 2$.

\

\

\subsubsection{Increments: fourth diagram}

Recall the definition
\begin{align*}
\calt^{\mathbf{4},(n)}_t(y)&:= \int_{0}^t ds \int dw \, K_{t-s}(y,w)\cac^{(n)}_{s,t}(y,w)^2 \big(\<Psi>^{(n)}_{0,t}(w)-\<Psi>^{(n)}_{0,t}(y)\big),
\end{align*}
and so, for $0\leq t_1<t_2$,
\begin{align*}
\calt^{\mathbf{4},(n)}_{t_1,t_2}(y)&:=\calt^{\mathbf{4},(n)}_{t_2}(y)-\calt^{\mathbf{4},(n)}_{t_1}(y)\\
&=\calt^{\mathbf{4,1},(n)}_{t_1,t_2}(y)+\calt^{\mathbf{4,2},(n)}_{t_1,t_2}(y)+\calt^{\mathbf{4,3},(n)}_{t_1,t_2}(y)+\calt^{\mathbf{4,4},(n)}_{t_1,t_2}(y)+\calt^{\mathbf{4,5},(n)}_{t_1,t_2}(y),
\end{align*}
with
\begin{align*}
\calt^{\mathbf{4,1},(n)}_{t_1,t_2}(y)&:=\int_{t_1}^{t_2} ds \int dw \, K_{t_2-s}(y,w)\cac^{(n)}_{s,t_2}(y,w)^2 \big(\<Psi>^{(n)}_{0,t_2}(w)-\<Psi>^{(n)}_{0,t_2}(y)\big),
\end{align*}

\begin{align*}
\calt^{\mathbf{4,2},(n)}_{t_1,t_2}(y)&:=\int_{0}^{t_1} ds \int dw \, K_{t_2-s}(y,w)\cac^{(n)}_{s,t_2}(y,w)^2 \big[\big(\<Psi>^{(n)}_{0,t_2}(w)-\<Psi>^{(n)}_{0,t_2}(y)\big)-\big(\<Psi>^{(n)}_{0,t_1}(w)-\<Psi>^{(n)}_{0,t_1}(y)\big)\big],
\end{align*}

\begin{align*}
\calt^{\mathbf{4,3},(n)}_{t_1,t_2}(y)&:=\int_{0}^{t_1} ds \int dw \, K_{t_1-s,t_2-s}(y,w)\cac^{(n)}_{s,t_2}(y,w)^2 \big(\<Psi>^{(n)}_{0,t_1}(w)-\<Psi>^{(n)}_{0,t_1}(y)\big),
\end{align*}

\begin{align*}
\calt^{\mathbf{4,4},(n)}_{t_1,t_2}(y)&:=\int_{0}^{t_1} ds \int dw \, K_{t_1-s}(y,w)\cac^{(n)}_{s,(t_1,t_2)}(y,w)\cac^{(n)}_{s,t_2}(y,w) \big(\<Psi>^{(n)}_{0,t_1}(w)-\<Psi>^{(n)}_{0,t_1}(y)\big),
\end{align*}

\begin{align*}
\calt^{\mathbf{4,5},(n)}_{t_1,t_2}(y)&:=\int_{0}^{t_1} ds \int dw \, K_{t_1-s}(y,w)\cac^{(n)}_{s,t_1}(y,w)\cac^{(n)}_{s,(t_1,t_2)}(y,w)\big(\<Psi>^{(n)}_{0,t_1}(w)-\<Psi>^{(n)}_{0,t_1}(y)\big).
\end{align*}

\

For clarity, we introduce the following notation.

\begin{notation}
We set
\begin{equation*} 
\cd^{(n)}_{t,t}\big((y_1,w_1),(y_2,w_2)\big):=\mathbb{E} \Big[ \big( \<Psi>^{(n)}_{0,t}(w_1)-\<Psi>^{(n)}_{0,t}(y_1)\big)\big(\<Psi>^{(n)}_{0,t}(w_2)-\<Psi>^{(n)}_{0,t}(y_2)\big)\Big].
\end{equation*}
\end{notation}

\

\noindent
\underline{\textit{Case of $\calt^{\mathbf{4,1},(n)}$}}. One has, for $0\leq t_1<t_2$,
\begin{multline*}
\cm^{\mathbf{4,1},(n)}_{t_1,t_2}(y_1,y_2)=\int_{t_1}^{t_2} ds_1\int_{t_1}^{t_2} ds_2 \int dw_1\int dw_2\, K_{t_2-s_1}(y_1,w_1)\cac^{(n)}_{s_1,t_2}(y_1,w_1)^2 \\
\hspace{2cm}\times K_{t_2-s_2}(y_2,w_2)\cac^{(n)}_{s_2,t_2}(y_2,w_2)^2 \cd^{(n)}_{t_2,t_2}\big((y_1,w_1),(y_2,w_2)\big).
\end{multline*}
Therefore
\begin{multline}
 \int dy_1 \sup_{y_2} \big| \cm^{\mathbf{4,1},(n)}_{t_1,t_2}(y_1+y_2,y_2)\big|  \lesssim   \\
 \lesssim \int_{t_1}^{t_2} \frac{ds_1}{|t_2-s_1|}\int_{t_1}^{t_2} \frac{ds_2}{|t_2-s_2|}\int dy_1\int dw_1\int dw_2\, K^\star_{t_2-s_1}(y_1-w_1) K^\star_{t_2-s_2}(w_2) F^{(n)}_{t_2}(y_1,w_1,w_2),\label{goback}
\end{multline}
where we have set
\begin{align}\label{deffn}
F^{(n)}_{t}(y_1,w_1,w_2):=\sup_{y_2}\big| \cd^{(n)}_{t,t}\big((y_1+y_2,w_1+y_2),(y_2,w_2+y_2)\big)\big|.
\end{align}
Observe that for all $t_2\geq 0$ and $\la\in [0,\frac12]$,
\begin{align}
&\big| \cd^{(n)}_{t_2,t_2}\big((y_1+y_2,w_1+y_2),(y_2,w_2+y_2)\big)\big|\lesssim \nonumber\\
&\lesssim \int_{0}^{+\infty}d\si\, \big| K_\si(w_1+y_2,w_2+y_2)-K_\si(w_1+y_2,y_2)-K_\si(y_1+y_2,w_2+y_2)+K_\si(y_1+y_2,y_2)\big|\nonumber\\
&\lesssim \int_{0}^{+\infty}d\si\, \Big[|w_1-y_1| \big\| \nabla_1 K_\si\big\|_{L^\infty(\R^6)}  \Big]^\la \Big[ |w_2| \big\| \nabla_2 K_\si\big\|_{L^\infty(\R^6)}\Big]^\la \nonumber\\
& \hspace{1cm}\times \Big[| K_\si(w_1+y_2,w_2+y_2)|+|K_\si(w_1+y_2,y_2)|+|K_\si(y_1+y_2,w_2+y_2)|+|K_\si(y_1+y_2,y_2)|\Big]^{1-2\la}\nonumber\\
&\lesssim |w_1-y_1|^\la |w_2|^\la \int_{0}^{+\infty}\frac{d\si}{\si^{4\la}}\, \Big[| K_\si^\star(w_1-w_2)|^{1-2\la}+|K_\si^\star(w_1)|^{1-2\la}+|K_\si^\star(y_1-w_2)|^{1-2\la}+|K_\si^\star(y_1)|^{1-2\la}\Big],\label{justas}
\end{align}
which, going back to \eqref{goback}, yields for every $\la\in [0,\frac12]$
\begin{multline*}
\int dy_1 \sup_{y_2} \big| \cm^{\mathbf{4,1},(n)}_{t_1,t_2}(y_1+y_2,y_2)\big|\lesssim \\
\begin{aligned}
&\lesssim \int_{t_1}^{t_2} \frac{ds_1}{|t_2-s_1|}\int_{t_1}^{t_2} \frac{ds_2}{|t_2-s_2|}\int_{0}^{+\infty}\frac{d\si}{\si^{4\la}}\int dy_1\int dw_1\int dw_2\, \big[K^\star_{t_2-s_1}(y_1-w_1) |w_1-y_1|^\la \big]\\
&\hspace{.5cm} \times \big[ K^\star_{t_2-s_2}(w_2)  |w_2|^\la\big] \, \big[| K_\si^\star(w_1-w_2)|^{1-2\la}+|K_\si^\star(w_1)|^{1-2\la}+|K_\si^\star(y_1-w_2)|^{1-2\la}+|K_\si^\star(y_1)|^{1-2\la}\big]\\
&\lesssim \int_{t_1}^{t_2} \frac{ds_1}{|t_2-s_1|}\big\|K^\star_{t_2-s_1}\cdot |.|^\la \big\|_{L^1}\int_{t_1}^{t_2} \frac{ds_2}{|t_2-s_2|}\big\|K^\star_{t_2-s_2}\cdot |.|^\la \big\|_{L^1}\int_{0}^{+\infty}\frac{d\si}{\si^{4\la}}\,\big\|| K_\si^\star|^{1-2\la}\big\|_{L^1}.
\end{aligned}
\end{multline*}
It is readily checked that for all $\la\in [0,\frac12)$,
\begin{align}
\big\|| K_\si^\star|^{1-2\la}\big\|_{L^1}&\lesssim \frac{1}{\sinh(\si)^{\frac32(1-2\la)}}\int dy \, \exp\Big(-\frac{(1-2\la)}{4\tanh(\si)}|y|^2\Big)\nonumber\\
&\lesssim \frac{\tanh(\si)^{\frac32}}{\sinh(\si)^{\frac32}}\sinh(\si)^{3\la}\lesssim \frac{\sinh(\si)^{3\la}}{\cosh(\si)^{\frac32}}\lesssim \si^{3\la}e^{-3\si(\frac12-\la)}\label{using}
\end{align}
due to the elementary estimates $\sinh(\si)\lesssim \si e^{\si}$ and $\cosh(\si)\gtrsim e^\si$.
Combining this bound with the result of Lemma \ref{lem:k-star-0}, we deduce that for every $\la\in [0,\frac12)$,
\begin{align*}
\int dy_1 \sup_{y_2} \big| \cm^{\mathbf{4,1},(n)}_{t_1,t_2}(y_1+y_2,y_2)\big|&\lesssim \bigg(\int_{t_1}^{t_2} \frac{ds}{|t_2-s|^{1-\frac{\la}{2}}}e^{-(t_2-s)}\bigg)^2\int_{0}^{+\infty}\frac{d\si}{\si^{\la}}e^{-3\si(\frac12-\la)}.
\end{align*}
By choosing $\la:=\frac12-\eps$, we deduce that
\begin{align*}
\int dy_1 \sup_{y_2} \big| \cm^{\mathbf{4,1},(n)}_{t_1,t_2}(y_1+y_2,y_2)\big|&\lesssim \bigg(\int_{t_1}^{t_2} \frac{ds}{|t_2-s|^{\frac34+\frac{\eps}{2}}}\bigg)^2\int_{0}^{+\infty}\frac{d\si}{\si^{\frac12-\eps}}e^{-3\si\eps}\lesssim |t_2-t_1|^{\frac12-\eps},
\end{align*}
for every $\eps>0$, uniformly over $0\leq t_1,t_2\leq 2$.

\

\

\noindent
\underline{\textit{Case of $\calt^{\mathbf{4,2},(n)}$}}. One has, for $0\leq t_1<t_2$,
\begin{align*}
&\big|\cm^{\mathbf{4,2},(n)}_{t_1,t_2}(y_1,y_2)\big|=\bigg|\int_{0}^{t_1} ds_1\int_{0}^{t_1} ds_2 \int dw_1\int dw_2\, 
\\
&\hspace{2cm}\times  K_{t_2-s_1}(y_1,w_1)\cac^{(n)}_{s_1,t_2}(y_1,w_1)^2 K_{t_2-s_2}(y_2,w_2)\cac^{(n)}_{s_2,t_2}(y_2,w_2)^2 \\
&\hspace{4cm}\times \mathbb{E}\Big[\big[\big(\<Psi>^{(n)}_{0,t_2}(w_1)-\<Psi>^{(n)}_{0,t_2}(y_1)\big)-\big(\<Psi>^{(n)}_{0,t_1}(w_1)-\<Psi>^{(n)}_{0,t_1}(y_1)\big)\big] \\
 & \hspace{6cm}  \times   \big[\big(\<Psi>^{(n)}_{0,t_2}(w_2)-\<Psi>^{(n)}_{0,t_2}(y_2)\big)-\big(\<Psi>^{(n)}_{0,t_1}(w_2)-\<Psi>^{(n)}_{0,t_1}(y_2)\big)\big] \Big]\bigg|\\
&\lesssim \int_{0}^{t_1} \frac{ds_1}{|t_2-s_1|}\int_{0}^{t_1} \frac{ds_2}{|t_2-s_2|} \int dw_1\int dw_2\, K^\star_{t_2-s_1}(y_1-w_1) K^\star_{t_2-s_2}(y_2-w_2) \\
&\hspace{2cm}\times \Big|\mathbb{E}\Big[\big[\big(\<Psi>^{(n)}_{0,t_2}(w_1)-\<Psi>^{(n)}_{0,t_2}(y_1)\big)-\big(\<Psi>^{(n)}_{0,t_1}(w_1)-\<Psi>^{(n)}_{0,t_1}(y_1)\big)\big]\\
    & \hspace{6cm}  \times \big[\big(\<Psi>^{(n)}_{0,t_2}(w_2)-\<Psi>^{(n)}_{0,t_2}(y_2)\big)-\big(\<Psi>^{(n)}_{0,t_1}(w_2)-\<Psi>^{(n)}_{0,t_1}(y_2)\big)\big] \Big]\Big|.
\end{align*}
Observe that by symmetry
\begin{multline*}
\Big|\mathbb{E}\Big[\big[\big(\<Psi>^{(n)}_{0,t_2}(w_1)-\<Psi>^{(n)}_{0,t_2}(y_1)\big)-\big(\<Psi>^{(n)}_{0,t_1}(w_1)-\<Psi>^{(n)}_{0,t_1}(y_1)\big)\big] \\
  \hspace{3cm}  \times \big[\big(\<Psi>^{(n)}_{0,t_2}(w_2)-\<Psi>^{(n)}_{0,t_2}(y_2)\big)-\big(\<Psi>^{(n)}_{0,t_1}(w_2)-\<Psi>^{(n)}_{0,t_1}(y_2)\big)\big] \Big]\Big|\\
  \begin{aligned}
&=\Big|\mathbb{E}\Big[\big(\<Psi>^{(n)}_{0,t_2}(w_1)-\<Psi>^{(n)}_{0,t_2}(y_1)\big)\big[\big(\<Psi>^{(n)}_{0,t_2}(w_2)-\<Psi>^{(n)}_{0,t_2}(y_2)\big)-\big(\<Psi>^{(n)}_{0,t_1}(w_2)-\<Psi>^{(n)}_{0,t_1}(y_2)\big)\big] \Big]\\
&\hspace{3cm}+\mathbb{E}\Big[\big(\<Psi>^{(n)}_{0,t_1}(w_1)-\<Psi>^{(n)}_{0,t_1}(y_1)\big) \big[\big(\<Psi>^{(n)}_{0,t_1}(w_2)-\<Psi>^{(n)}_{0,t_1}(y_2)\big)-\big(\<Psi>^{(n)}_{0,t_2}(w_2)-\<Psi>^{(n)}_{0,t_2}(y_2)\big)\big] \Big]\Big|\\
&\lesssim \int_0^{|t_2-t_1|} d\si \, \Big| K_\si(w_1,w_2) - K_\si(y_1,w_2)- K_\si(w_1,y_2)+ K_\si (y_1,y_2) \Big|\\
&\hspace{3cm}+\int_{2(t_1+\eps_n)}^{2(t_2+\eps_n)} d\si \, \Big| K_\si(w_1,w_2) - K_\si(y_1,w_2)- K_\si(w_1,y_2)+ K_\si (y_1,y_2) \Big|.
  \end{aligned}
\end{multline*}
Then, just as in \eqref{justas}, one has for all $\si>0$ and $\la\in [0,\frac12]$,
\begin{align*}
&\Big| K_\si(w_1+y_2,w_2+y_2) - K_\si(y_1+y_2,w_2+y_2)- K_\si(w_1+y_2,y_2)+ K_\si (y_1+y_2,y_2) \Big| \lesssim\\
&\lesssim   \frac{|w_1-y_1|^\la |w_2|^\la }{\si^{4\la}} \Big(K_\si^\star(w_1-w_2)^{1-2\la}+ K_\si^\star(y_1-w_2)^{1-2\la}+ K_\si^\star(w_1)^{1-2\la}+ K_\si^\star(y_1)^{1-2\la} \Big).
\end{align*}
Therefore, for every $\la\in [0,\frac12)$, one has, using both \eqref{using} and Lemma \ref{lem:k-star-0},
\begin{align*}
 \int dy_1 \sup_{y_2} \big| \cm^{\mathbf{4,2},(n)}_{t_1,t_2}(y_1+y_2,y_2)\big|&\lesssim  \int_{0}^{t_1} \frac{ds_1}{|t_2-s_1|}\big\|K^\star_{t_2-s_1}\cdot |.|^\la \big\|_{L^1}\int_{0}^{t_1} \frac{ds_2}{|t_2-s_2|}\big\|K^\star_{t_2-s_2}\cdot |.|^\la \big\|_{L^1} \\
&\hspace{1cm}\cdot \bigg[\int_0^{|t_2-t_1|}\frac{d\si}{\si^{4\la}}\,\big\|| K_\si^\star|^{1-2\la}\big\|_{L^1}+\int_{2(t_1+\eps_n)}^{2(t_2+\eps_n)}\frac{d\si}{\si^{4\la}}\,\big\|| K_\si^\star|^{1-2\la}\big\|_{L^1}\bigg]\\
&\lesssim \bigg(\int_{0}^{t_1} \frac{ds}{|t_2-s|^{1-\frac{\la}{2}}}e^{-(t_2-s)}\bigg)^2\bigg[\int_0^{|t_2-t_1|}\frac{d\si}{\si^{\la}}+\int_{2(t_1+\eps_n)}^{2(t_2+\eps_n)}\frac{d\si}{\si^{\la}} \bigg]\\
&\lesssim \bigg(\int_{0}^{t_1} \frac{ds}{|t_2-s|^{1-\frac{\la}{2}}}e^{-(t_2-s)}\bigg)^2 |t_2-t_1|^{1-\la}.
\end{align*}
By choosing $\la:=\eps>0$, we obtain 
\begin{align*}
\int dy_1 \sup_{y_2} \big| \cm^{\mathbf{4,2},(n)}_{t_1,t_2}(y_1+y_2,y_2)\big|&\lesssim \bigg(\int_{0}^{+\infty} \frac{ds}{|s|^{1-\frac{\eps}{2}}}e^{-s}\bigg)^2  |t_2-t_1|^{1-\eps} \lesssim  |t_2-t_1|^{1-\eps}, 
\end{align*}
for every small $\eps>0$, uniformly over $0\leq t_1,t_2\leq 2$.

\

\noindent
\underline{\textit{Case of $\calt^{\mathbf{4,3},(n)}$}}. One has, for $0\leq t_1<t_2$,
\begin{multline*}
\cm^{\mathbf{4,3},(n)}_{t_1,t_2}(y_1,y_2)=\int_{0}^{t_1} ds_1\int_{0}^{t_1} ds_2  \int dw_1\int dw_2\, K_{t_1-s_1,t_2-s_1}(y_1,w_1)\cac^{(n)}_{s_1,t_2}(y_1,w_1)^2\\
\times  K_{t_1-s_2,t_2-s_2}(y_2,w_2)\cac^{(n)}_{s_2,t_2}(y_2,w_2)^2 \cd^{(n)}_{t_1,t_1}\big((y_1,w_1),(y_2,w_2)\big).
\end{multline*}
Therefore, using the notation $F^{(n)}$ introduced in \eqref{deffn},
\begin{multline*}
\int dy_1 \sup_{y_2} \big| \cm^{\mathbf{4,3},(n)}_{t_1,t_2}(y_1+y_2,y_2)\big| \lesssim
 \int_{0}^{t_1} \frac{ds_1}{|t_2-s_1|}\int_{0}^{t_1} \frac{ds_2}{|t_2-s_2|}\int dy_1\int dw_1\int dw_2\, \\
\hspace{4cm} \times K^\star_{t_1-s_1,t_2-s_1}(y_1-w_1) K^\star_{t_1-s_2,t_2-s_2}(w_2) F^{(n)}_{t_1}(y_1,w_1,w_2).
\end{multline*}
Recall that for every $\la\in [0,\frac12]$, one has by \eqref{justas}
\begin{align}
&\big|F^{(n)}_{t_1}(y_1,w_1,w_2)\big|\lesssim  |w_1-y_1|^\la |w_2|^\la \nonumber\\
& \times \int_{0}^{+\infty}\frac{d\si}{\si^{4\la}}\, \Big[| K_\si^\star(w_1-w_2)|^{1-2\la}+|K_\si^\star(w_1)|^{1-2\la}+|K_\si^\star(y_1-w_2)|^{1-2\la}+|K_\si^\star(y_1)|^{1-2\la}\Big],\label{again}
\end{align}
and so, for every $\la\in [0,\frac12]$, 
\begin{multline*}
\int dy_1 \sup_{y_2} \big| \cm^{\mathbf{4,3},(n)}_{t_1,t_2}(y_1+y_2,y_2)\big|\lesssim \\
\begin{aligned}
&\lesssim \int_{0}^{t_1} \frac{ds_1}{|t_2-s_1|}\big\|K^\star_{t_1-s_1,t_2-s_1}\cdot |.|^\la \big\|_{L^1}\int_{0}^{t_1} \frac{ds_2}{|t_2-s_2|}\big\|K^\star_{t_1-s_2,t_2-s_2}\cdot |.|^\la \big\|_{L^1}\int_{0}^{+\infty}\frac{d\si}{\si^{4\la}}\,\big\|| K_\si^\star|^{1-2\la}\big\|_{L^1}\\
&\lesssim \bigg(\int_{0}^{t_1} \frac{ds}{|t_2-s|}\big\|K^\star_{t_1-s,t_2-s}\cdot |.|^\la \big\|_{L^1}\bigg)^2\int_{0}^{+\infty}\frac{d\si}{\si^{4\la}}\,\big\|| K_\si^\star|^{1-2\la}\big\|_{L^1}.
\end{aligned}
\end{multline*}
Using both \eqref{using} and the result of Lemma \ref{lem:k-star-1}, we deduce that for all $\theta\in [0,1]$ and $\la\in [0,\frac12)$,
\begin{multline*}
\int dy_1 \sup_{y_2} \big| \cm^{\mathbf{4,3},(n)}_{t_1,t_2}(y_1+y_2,y_2)\big|\lesssim\\
\lesssim |t_2-t_1|^{2\theta}\bigg(\int_{0}^{t_1} \frac{ds\, e^{-(t_1-s)}}{|t_1-s|^{\frac52 \theta}|t_2-s|^{1-(\frac32\theta +\frac12 \la)}}\bigg)^2\int_{0}^{+\infty}\frac{d\si}{\si^{\la}}e^{-3\si(\frac12-\la)}.
\end{multline*}
By choosing $\la:=\frac12-\eps$ and $\theta:=\frac14-\eps$, we obtain 
\begin{align*}
\int dy_1 \sup_{y_2} \big| \cm^{\mathbf{4,3},(n)}_{t_1,t_2}(y_1+y_2,y_2)\big|&\lesssim |t_2-t_1|^{\frac12-2\eps}\bigg(\int_{0}^{t_1} \frac{ds\, e^{-(t_1-s)}}{|t_1-s|^{\frac58-\frac52 \eps}|t_2-s|^{\frac38+2\eps}}\bigg)^2\int_{0}^{+\infty}\frac{d\si}{\si^{\la}}e^{-3\si\eps}\\
&\lesssim |t_2-t_1|^{\frac12-2\eps}\bigg(\int_{0}^{t_1} \frac{ds\, e^{-(t_1-s)}}{|t_1-s|^{1-\frac{\eps}{2}}}\bigg)^2\int_{0}^{+\infty}\frac{d\si}{\si^{\la}}e^{-3\si\eps}\\
&\lesssim |t_2-t_1|^{\frac12-2\eps}\bigg(\int_{0}^{+\infty} \frac{ds}{|s|^{1-\frac{\eps}{2}}}e^{-s}\bigg)^2\int_{0}^{+\infty}\frac{d\si}{\si^{\la}}e^{-3\si\eps} \lesssim  |t_2-t_1|^{\frac12-2\eps},
\end{align*}
uniformly over $0\leq t_1,t_2\leq 2$.

\

\noindent
\underline{\textit{Case of $\calt^{\mathbf{4,4},(n)}$}}. One has, for $0\leq t_1<t_2$,
\begin{align*}
&\cm^{\mathbf{4,4},(n)}_{t_1,t_2}(y_1,y_2)=\int_{0}^{t_1} ds_1\int_{0}^{t_1} ds_2 \int dw_1\int dw_2\, K_{t_1-s_1}(y_1,w_1)\cac^{(n)}_{s_1,(t_1,t_2)}(y_1,w_1)\cac^{(n)}_{s_1,t_2}(y_1,w_1) \\
&\hspace{3.5cm}  \times K_{t_1-s_2}(y_2,w_2)\cac^{(n)}_{s_2,(t_1,t_2)}(y_2,w_2)\cac^{(n)}_{s_2,t_2}(y_2,w_2) \cd^{(n)}_{t_1,t_1}\big((y_1,w_1),(y_2,w_2)\big).
\end{align*}
Therefore, with the notation $F^{(n)}$ introduced in \eqref{deffn},
\begin{multline*}
 \int dy_1 \sup_{y_2} \big| \cm^{\mathbf{4,4},(n)}_{t_1,t_2}(y_1+y_2,y_2)\big|\lesssim \int_{0}^{t_1} \frac{ds_1}{|t_1-s_1|^{\frac12}}\int_{0}^{t_1} \frac{ds_2}{|t_1-s_2|^{\frac12}}\big\|\cac^{(n),\star}_{s_1,(t_1,t_2)}\big\|_{L^\infty} \big\|\cac^{(n),\star}_{s_2,(t_1,t_2)}\big\|_{L^\infty}\\
\hspace{1cm} \times \int dy_1\int dw_1\int dw_2\, K^\star_{t_1-s_1}(y_1-w_1) K^\star_{t_1-s_2}(w_2) F^{(n)}_{t_1}(y_1,w_1,w_2).
\end{multline*}
Using again \eqref{again}, we deduce for every $\la\in [0,\frac12]$,
\begin{multline*}
\int dy_1 \sup_{y_2} \big| \cm^{\mathbf{4,4},(n)}_{t_1,t_2}(y_1+y_2,y_2)\big|\lesssim \\
\begin{aligned}
&\lesssim \int_{0}^{t_1} \frac{ds_1}{|t_1-s_1|^{\frac12}}\big\|\cac^{(n),\star}_{s_1,(t_1,t_2)}\big\|_{L^\infty}\big\|K^\star_{t_1-s_1}\cdot |.|^\la \big\|_{L^1}\\
 & \hspace{3cm} \times  \int_{0}^{t_1} \frac{ds_2}{|t_1-s_2|^{\frac12}}\big\|\cac^{(n),\star}_{s_2,(t_1,t_2)}\big\|_{L^\infty}\big\|K^\star_{t_1-s_2}\cdot |.|^\la \big\|_{L^1}\int_0^{+\infty}\frac{d\si}{\si^{4\la}}\,\big\|| K_\si^\star|^{1-2\la}\big\|_{L^1}\\
&\lesssim \bigg(\int_{0}^{t_1} \frac{ds}{|t_1-s|^{\frac12}}\big\|\cac^{(n),\star}_{s,(t_1,t_2)}\big\|_{L^\infty}\big\|K^\star_{t_1-s}\cdot |.|^\la \big\|_{L^1}\bigg)^2 \int_0^{+\infty}\frac{d\si}{\si^{4\la}}\,\big\|| K_\si^\star|^{1-2\la}\big\|_{L^1}.
\end{aligned}
\end{multline*}
By combining the results of Lemma \ref{lem:k-star-0}, Lemma \ref{lem:cac3} and \eqref{using}, we derive that for all $\la\in [0,\frac12)$ and $\theta\in [0,1]$,
\begin{equation*}
\int dy_1 \sup_{y_2} \big| \cm^{\mathbf{4,4},(n)}_{t_1,t_2}(y_1+y_2,y_2)\big| 
\lesssim |t_2-t_1|^{2\theta}\bigg(\int_{0}^{t_1} \frac{ds}{|t_1-s|^{1-\frac{\la}{2}+\theta}}e^{-(t_1-s)}\bigg)^2\int_0^{+\infty}\frac{d\si}{\si^{\la}}e^{-3\si(\frac12-\la)}.
\end{equation*}
By choosing $\la:=\frac12-\eps$ and $\theta:=\frac14-\eps$, we obtain 
\begin{align*}
\int dy_1 \sup_{y_2} \big| \cm^{\mathbf{4,4},(n)}_{t_1,t_2}(y_1+y_2,y_2)\big|&\lesssim |t_2-t_1|^{\frac12-2\eps}\bigg(\int_{0}^{t_1} \frac{ds}{|t_1-s|^{1-\frac{\eps}{2}}}e^{-(t_1-s)}\bigg)^2\int_0^{+\infty}\frac{d\si}{\si^{\la}}e^{-3\si\eps}\\
&\lesssim |t_2-t_1|^{\frac12-2\eps}\bigg(\int_{0}^{+\infty} \frac{ds}{|s|^{1-\frac{\eps}{2}}}e^{-s}\bigg)^2 \int_0^{+\infty}\frac{d\si}{\si^{\la}}e^{-3\si\eps}\lesssim |t_2-t_1|^{\frac12-2\eps}, 
\end{align*}
uniformly over $0\leq t_1,t_2\leq 2$.

\

\

\noindent
\underline{\textit{Case of $\calt^{\mathbf{4,5},(n)}$}}. This quantity can be controlled with the same arguments as for $\calt^{\mathbf{4,4},(n)}$, leading to the same bound
\begin{align*}
\int dy_1 \sup_{y_2} \big| \cm^{\mathbf{4,5},(n)}_{t_1,t_2}(y_1+y_2,y_2)\big|&\lesssim |t_2-t_1|^{\frac12-2\eps}, 
\end{align*}
uniformly over $0\leq t_1,t_2\leq 2$.

\

We have finally shown that condition \eqref{hold-m-a-b-c} is indeed satisfied for all $\mathbf{a}$, $\mathbf{b}$, $\mathbf{d}$, and by \eqref{transition-coucou}, this completes the proof of \eqref{hold-coucou}.


\

\section{Technical lemmas}\label{sect:techlem}

\subsection{About $K^\star$}

Recall that, according to \eqref{star-not}, we have, for $\sigma>0$, 
\begin{equation*}
K_{\sigma}^\star(x) = \sup_{y \in \R^3}\big|K_{\sigma}(y,y+x)\big|, \quad x \in \R^3.
\end{equation*}
Similarly, for $\sigma_1, \sigma_2>0$, 
\begin{equation} 
K_{\sigma_1,\sigma_2}^\star(x) =    \sup_{y \in \R^3}\big|K_{\sigma_2}(y,y+x)-K_{\sigma_1}(y,y+x)\big|      , \quad x \in \R^3.
\end{equation}

\begin{lemma}\label{lem:k-star-0}
For all $p\geq 1$ and $\si,\la\geq 0$, one has
$$\big\|K^\star_\si \cdot |.|^\la\big\|_{L^p(\R^3)} \lesssim \frac{1\wedge \si^{\frac{\la}{2}}}{\si^{\frac32-\frac{3}{2p}}}e^{-\si}.$$
\end{lemma}

\begin{proof}
It is readily checked that
\begin{align*}
\big\|K^\star_\si\cdot |.|^\la\big\|_{L^p(\R^3)}&=\frac{c}{\sinh(2\si)^{\frac32}}\bigg(\int dy \, |y|^{\la p}\exp\Big(-\frac{p}{4 \tanh(\si)}|y|^2 \Big)\bigg)^{\frac1p}\\
&=c\,\frac{\tanh(\si)^{\frac{3}{2p}+\frac{\la}{2}}}{\sinh(2\si)^{\frac32}}=c\, \frac{\tanh(\si)^{\frac{\la}{2}}\sinh(\si)^{\frac{3}{2p}}}{\sinh(2\si)^{\frac32}\cosh(\si)^{\frac{3}{2p}}}\lesssim \frac{\tanh(\si)^{\frac{\la}{2}}}{\sinh(\si)^{\frac32-\frac{3}{2p}}\cosh(\si)^{\frac{3}{2p}}},
\end{align*}
which, combined with the fact that $\tanh(\si)\lesssim 1\wedge \si$, $\sinh(\si) \gtrsim \si e^{\frac{2\si}{3}}$ and $\cosh(\si)\gtrsim e^{ \si}$, entails
$$\big\|K^\star_\si \cdot |.|^\la\big\|_{L^p(\R^3)}\lesssim \frac{1\wedge \si^{\frac{\la}{2}}}{\si^{\frac32-\frac{3}{2p}}}\exp\Big(-\frac{2\si}{3}(\frac32-\frac{3}{2p})-\frac{3\si}{2p}\Big)\lesssim \frac{1\wedge \si^{\frac{\la}{2}}}{\si^{\frac32-\frac{3}{2p}}}e^{-\si}.$$
\end{proof}

\begin{lemma}\label{lem:k-star-1}
For all $p\geq 1$, $0< \si_1<\si_2$, $\la\geq 0$ and $0\leq \theta\leq 1$, one has
\begin{equation}\label{boun-k-star-incr}
\big\|K_{\si_1,\si_2}^\star \cdot |.|^\la\big\|_{L^p(\R^3)} \lesssim |\si_2-\si_1|^\theta \frac{(1\wedge \si_2^{\frac{3}{2p}\theta+\frac{\la}{2}})}{\si_1^{\frac32-\frac{3}{2p}+\theta(1+\frac{3}{2p})}}e^{-\si_1},
\end{equation}
where the proportional constant does not depend on $\si_1,\si_2$.
\end{lemma}

\begin{proof}
We first establish the two extreme cases $\theta=0$ and $\theta=1$ of \eqref{boun-k-star-incr}, and then derive the general case by interpolation.

\smallskip

\noindent
$\bullet$ Case $\theta=0$. Since $K^\star_{\si_1,\si_2}\leq K^\star_{\si_1}+K^\star_{\si_2}$, we know by Lemma \ref{lem:k-star-0} that
\begin{equation}\label{ksta1}
\big\| K^\star_{\si_1,\si_2}\cdot |.|^\la\big\|_{L^p(\R^3)} \lesssim \frac{1\wedge \si_1^{\frac{\la}{2}}}{\si_1^{\frac32-\frac{3}{2p}}}e^{-\si_1}+\frac{1\wedge \si_2^{\frac{\la}{2}}}{\si_2^{\frac32-\frac{3}{2p}}}e^{-\si_2}\lesssim  \frac{1\wedge \si_2^{\frac{\la}{2}}}{\si_1^{\frac32-\frac{3}{2p}}}e^{-\si_1},
\end{equation}
where the second inequality follows from $0<\si_1<\si_2$ together with the fact that $\frac32-\frac{3}{2p}\geq 0$. This yields the claimed bound when $\theta=0$.

\smallskip

\noindent
$\bullet$ Case $\theta=1$.  By the Mehler formula \eqref{mehler2}, one has, for all $\si>0$ and $y,w\in \R^3$,
$$K_\si(y,y+w)=\frac{1}{(2\pi \sinh(2\si))^{\frac32}}\exp\Big(-\frac{|w|^2}{4\tanh(\si)}-\frac{\tanh(\si)}{4}|2y+w|^2\Big),$$
and therefore, differentiating with respect to $\si$,
\begin{equation}\label{derivK}
\partial_\si \big[K_\si(y,y+w)\big]=K_\si(y,y+w)\bigg[-3\,\frac{\cosh(2\si)}{\sinh(2\si)}+\frac{|w|^2}{4\sinh(\si)^2}-\frac{|2y+w|^2}{4\cosh(\si)^2}\bigg].
\end{equation}
On the one hand, we have the bound
$$\frac{|2y+w|^2}{4\cosh(\si)^2}\,\exp\Big(-\frac{\tanh(\si)}{4}|2y+w|^2\Big) \lesssim \frac{1}{\sinh(2\si)} \lesssim \frac{\cosh(2\si)}{\sinh(2\si)},$$
and, on the other hand,
$$\frac{|w|^2}{4\sinh(\si)^2}\,\exp\Big(-\frac{|w|^2}{4\tanh(\si)}\Big) \lesssim\frac{1}{\sinh(2\si)}\exp\Big(-\frac{|w|^2}{8\tanh(\si)}\Big)\lesssim\frac{\cosh(2\si)}{\sinh(2\si)}\exp\Big(-\frac{|w|^2}{8\tanh(\si)}\Big).$$

Inserting these two bounds into  \eqref{derivK}, we deduce that, uniformly over $y\in \R^3$,
$$\big|\partial_\si \big[K_\si(y,y+w)\big]\big|\lesssim  \frac{\cosh(2\si)}{\sinh(2\si)^{\frac52}}\,\exp\Big(-\frac{|w|^2}{8\tanh(\si)}\Big).$$
 Since $\cosh(2\si)\leq e^{2\si}$ and $\sinh(2\si)\gtrsim \si\, e^{\frac65 \si}$, it follows that
$$\sup_{y\in \R^3}\big|\partial_\si \big[K_\si(y,y+w)\big]\big|\lesssim \si^{-\frac52}\, e^{-\si}\,\exp\Big(-\frac{|w|^2}{8\tanh(\si)}\Big).$$
As a consequence, for all $0<\si_1<\si_2$ and $y,w\in \R^3$,  
$$\big|K_{\si_2}(y,y+w)-K_{\si_1}(y,y+w)\big|\leq \int_{\si_1}^{\si_2}d\si \, \big|\partial_\si \big[K_\si(y,y+w)\big]\big|\lesssim |\si_2-\si_1|\, \si_1^{-\frac52}\,e^{-\si_1}\exp\Big(-\frac{|w|^2}{8\tanh(\si_2)}\Big).$$
Taking the supremum over $y$ in the left-hand side, we obtain
$$\big|K^\star_{\si_1,\si_2}(w)\big|\lesssim |\si_2-\si_1|\, \si_1^{-\frac52}\,e^{-\si_1} \exp\Big(-\frac{|w|^2}{8\tanh(\si_2)}\Big),$$
which yields
\begin{equation}\label{ksta2}
\big\| K^\star_{\si_1,\si_2}\cdot |.|^\la\big\|_{L^p(\R^3)}\lesssim \tanh(\si_2)^{\frac{3}{2p}+\frac{\la}{2}}\si_1^{-\frac52} |\si_2-\si_1| e^{-\si_1}\lesssim (1\wedge \si_2^{\frac{3}{2p}+\frac{\la}{2}})\si_1^{-\frac52} |\si_2-\si_1| e^{-\si_1}.
\end{equation}

\

Interpolating between \eqref{ksta1} and \eqref{ksta2} gives the desired bound \eqref{boun-k-star-incr}.
\end{proof}

\

\subsection{About $\cac^{(n),\star}$}

In the following result, we gather some  estimates of $\cac^{(n),\star}$ in $L^p$ spaces.

\begin{lemma}\label{lem:cac3}
For all $0<\eps<1$ and  $t_1,t_2\geq 0$, it holds that
\begin{equation*} 
\sup_{n\geq 1} \big\|  \cac^{(n),\star}_{t_1,t_2}\big\|_{L^p(\R^3)}\lesssim 1 \quad \text{if} \ p<3,
\end{equation*}  
\begin{equation*} 
\sup_{n\geq 1} \big\|  \cac^{(n),\star}_{t_1,t_2}\big\|_{L^3(\R^3)}\lesssim \frac{1}{|t_2-t_1|^{\varepsilon}},
\end{equation*}   
while for $3<p\leq +\infty$,
\begin{equation*} 
\sup_{n\geq 1} \big\|  \cac^{(n),\star}_{t_1,t_2}\big\|_{L^p(\R^3)} \lesssim \frac{1}{|t_2-t_1|^{\frac12-\frac{3}{2p}}},
\end{equation*}                      
where the proportional constants do not depend on $t_1,t_2\geq 0$. Moreover, for $0\leq s,t_1\leq t_2$, one has 
\begin{equation*}
\sup_{n\geq 1} \big\|  \cac^{(n),\star}_{s,(t_1,t_2)}\big\|_{L^1(\R^3)} \lesssim |t_2-t_1|,
\end{equation*} 
\begin{equation*}
\sup_{n\geq 1} \big\|  \cac^{(n),\star}_{s,(t_1,t_2)}\big\|_{L^2(\R^3)} \lesssim |t_2-t_1|^{\frac14},
\end{equation*} 
while for all $\theta\in [0,1]$ and $3<p\leq +\infty$,
\begin{equation*}
\sup_{n\geq 1} \big\|  \cac^{(n),\star}_{s,(t_1,t_2)}\big\|_{L^p(\R^3)} \lesssim \frac{|t_2-t_1|^\theta}{\big(|t_1-s|\wedge |t_2-s|\big)^{\frac12-\frac{3}{2p}+\theta}}.
\end{equation*}
\end{lemma}

\

\subsection{An interpolation result}

 By \cite[Theorem 2.36 and Lemma 3.31]{BCD}, if $0<\lambda<1$, we have the equivalence 
 $$\mathcal{C}_b^\lambda(\R) = B^{\lambda} _{\infty,\infty}(\R)$$
 where $ \mathcal{C}_b^\lambda(\R)$ is the space of  bounded H\"older functions with exponent $\lambda$ defined by the norm
 \begin{equation*} 
 \| u\|_{ \mathcal{C}_b^\lambda(\R)}= \| u\|_{L^\infty(\R)}+ \sup_{x, y \in \R, x \neq y} \frac{|u(x)-u(y)|}{|x-y|^\lambda}.
  \end{equation*}
  Then, for $\gamma \in \R$, we denote by $\mathcal{C}_b^\lambda(\R ; \mathcal{B}_{\infty}^{\gamma}) $ the space of $\mathcal{B}_{\infty}^{\gamma}-$valued, bounded H\"older functions, equipped with the norm
   \begin{equation}\label{HolB1}
 \| u\|_{ \mathcal{C}_b^\lambda(\R ; \mathcal{B}_{\infty}^{\gamma})}= \| u\|_{L^\infty(\R;  \mathcal{B}_{\infty}^{\gamma})}+ \sup_{x, y \in \R, x \neq y} \frac{\|u(x)-u(y)\|_{\mathcal{B}_{\infty}^{\gamma}}}{|x-y|^\lambda}.
  \end{equation}
 Using this identification, we prove the following interpolation result.

 	\begin{lemma}\label{lemma-interpt}
Fix $0\leq T_1<T_2$, $\beta \leq \alpha$ and $0<\ka, \eta <1$. Let $f\in \mathcal{C}^\eta([T_1,T_2+1]; \B^\beta_{\infty})$ and set
$$\tilde{f}_{T_1,t}:=\int_{T_1}^t f_s \, ds.$$
Assume that $\tilde{f}_{T_1,.} \in  \mathcal{C}^{1-\ka}([T_1,T_2+1]; \B^\alpha_{\infty})$. Then, for all $\dis \theta\in (\frac{\ka}{\ka+\eta},1)$, setting
$$\lambda(\theta):= -\ka+\theta(\ka+\eta)>0,\quad  \quad \gamma(\theta):=\alpha -\theta (\alpha-\beta),$$
it holds that
\begin{equation*} 
\big\| f\big\|_{\mathcal{C}^{\lambda(\theta)}([T_1,T_2]; \mathcal{B}_{\infty}^{\gamma(\theta)})} \lesssim  \big\| \tilde{f}_{T_1,.}\big\|^{1-\theta}_{{\ov \cac}^{1-\ka}([T_1,T_2+1]; \mathcal{B}_{\infty}^\alpha)} \big\| f\big\|^{\theta}_{\mathcal{C}^{\eta}([T_1,T_2+1]; \mathcal{B}_{\infty}^\beta)},
\end{equation*}   
where the proportional constant does not depend on $T_1,T_2$.
\end{lemma}
 
 \begin{proof}
 To begin with, we extend $f$ by $f(T_1)$ on $(-\infty,T_1)$ in the variable $t$, which still defines a H\"older function with the same exponent. Let us introduce $t \mapsto \phi(t)  \in \mathcal{C}_0^\infty(\R)$ such that $0 \leq \phi \leq 1$, $0\leq |\phi'|\leq 2$, $\phi \equiv 1 $ on $[T_1,T_2]$ and $\phi=0$ on $(-\infty,T_1-1] \cup [T_2+1,\infty)$. Then, for $\gamma \in \R$ and $0<\la <1$, we have
\begin{equation}\label{eqqq}
\big\| f\big\|_{\mathcal{C}^{\lambda}([T_1,T_2]; \mathcal{B}_{\infty}^{\gamma})} \leq            \big\| \phi f\big\|_{\mathcal{C}_b^{\lambda}(\R; \mathcal{B}_{\infty}^{\gamma})}  \leq 6 \big\| f\big\|_{\mathcal{C}^{\lambda}([T_1,T_2+1]; \mathcal{B}_{\infty}^{\gamma})},
\end{equation}  
where $\mathcal{C}_b^{\lambda}(\R ; \mathcal{B}_{\infty}^{\gamma})$ is defined in \eqref{HolB1}. 
				
\smallskip
				
Recall that we have the equivalence (see  \cite[Theorem~2.36 and Lemma 3.31]{BCD})
\begin{equation}\label{eqqq-cucu}
\mathcal{C}_b^\lambda(\R ; \mathcal{B}_{\infty}^{\gamma})   \cong  B^{\lambda}(\R; \mathcal{B}_{\infty}^{\gamma}) ,
\end{equation} 
where $B^{\lambda}(\R; \mathcal{B}_{\infty}^{\gamma}):= B^{\lambda} _{\infty,\infty}(\R; \mathcal{B}_{\infty}^{\gamma})$ is defined through a straightforward extension of~\eqref{def-besov} (see \cite[Section 5]{Amann}). Namely, consider $\chi$ as in \eqref{partition} and set 
$$\Delta_k u =  \chi_k(\sqrt{1-\partial^2_t}) u=  \chi\big( 2^{-k} \sqrt{1-\partial^2_t}\big)u,$$
so that the definition in \cite[(5.1)]{Amann} reads as
$$ \| u\| _{B^{\lambda}(\R; \mathcal{B}_{\infty}^{\gamma})} = \sup_{\substack {j \geq 0 \\ k \geq -1}  } 2^{j\gamma }2^{k\lambda}    \| \Delta_k  {\delta}_j u \|_{L^\infty(\R \times \R^3)} .
$$   
It is now readily checked that the latter norm satisfies the interpolation property: for all $0 \leq \theta \leq 1$, 
 \begin{equation}\label{interPO}
\| u\| _{B^{\lambda(\theta)}(\R; \mathcal{B}_{\infty}^{\gamma(\theta)})}   \leq \| u\| ^{1-\theta}_{B^{-\ka}(\R; \mathcal{B}_{\infty}^{\alpha})}  \| u\|^\theta _{B^{\eta}(\R; \mathcal{B}_{\infty}^{\beta})},
\end{equation}  
where $\lambda(\theta):= -\ka+\theta(\ka+\eta)>0$ and $\gamma(\theta):=\alpha -\theta (\alpha-\beta)$.

\smallskip

At this point, we deduce that for all $\dis \theta\in (\frac{\ka}{\ka+\eta},1)$ and by setting $\lambda= -\ka+\theta(\ka+\eta)\in (0,1)$, $ \gamma=\alpha -\theta (\alpha-\beta)$,
\begin{align*}
\big\| f\big\|_{\mathcal{C}^{\lambda}([T_1,T_2]; \mathcal{B}_{\infty}^{\gamma})} &\lesssim \big\| \phi f\big\|_{\mathcal{C}_b^{\lambda}(\R; \mathcal{B}_{\infty}^{\gamma})} \quad \quad \text{(by \eqref{eqqq})}\\
&\lesssim \big\| \phi f\big\|_{B^{\lambda}(\R; \mathcal{B}_{\infty}^{\gamma})} \quad \quad \text{(by \eqref{eqqq-cucu})}\\
&\lesssim \|\phi f \| ^{1-\theta}_{B^{-\ka}(\R; \mathcal{B}_{\infty}^{\alpha})}  \| \phi f\|^\theta _{B^{\eta}(\R; \mathcal{B}_{\infty}^{\beta})}  \quad \quad \text{(by \eqref{interPO})}\\
&\lesssim \|\phi f \| ^{1-\theta}_{B^{-\ka}(\R; \mathcal{B}_{\infty}^{\alpha})} \big\| f\big\|_{\mathcal{C}^{\eta}([T_1,T_2+1]; \mathcal{B}_{\infty}^{\beta})}^\theta  \quad \quad \text{(by \eqref{eqqq})},
\end{align*}
and so it remains to show that
\begin{equation}\label{lastpap}
\big\| \phi f\big\|_{{B}^{-\ka}(\R ; \mathcal{B}_{\infty}^{\alpha})}  \lesssim     \big\| \tilde{f}_{T_1,.}\big\|_{{\ov \cac}^{1-\ka}([T_1,T_2+1]; \mathcal{B}_{\infty}^\alpha)}.
\end{equation}
To this end, write $\phi f = \partial_t (\phi \tilde{f}_{T_1,.}) - (\partial_t \phi ) \tilde{f}_{T_1,.}$, which gives
$$\big\| \phi f\big\|_{{B}^{-\ka}(\R ; \mathcal{B}_{\infty}^{\alpha})} \lesssim \big\| \partial_t (\phi \tilde{f}_{T_1,.})\big\|_{{B}^{-\ka}(\R ; \mathcal{B}_{\infty}^{\alpha})}+\big\| (\partial_t \phi ) \tilde{f}_{T_1,.}\big\|_{{B}^{-\ka}(\R ; \mathcal{B}_{\infty}^{\alpha})}.$$
Then, by using the definition of the Besov spaces, we obtain that 
$$\big\| \partial_t ( \phi  \tilde{f}_{T_1,.} )\big\|_{{B}^{-\ka}(\R ; \mathcal{B}_{\infty}^{\alpha})} \lesssim \big\|  \phi  \tilde{f}_{T_1,.} \big\|_{{B}^{1-\ka}(\R ; \mathcal{B}_{\infty}^{\alpha})}$$
and hence
\begin{align*}
\big\| \phi f\big\|_{{B}^{-\ka}(\R ; \mathcal{B}_{\infty}^{\alpha})} &\lesssim  \big\| \phi \tilde{f}_{T_1,.}\big\|_{{B}^{1-\ka}(\R ; \mathcal{B}_{\infty}^{\alpha})}+ \big\| (\partial_t \phi)  \tilde{f}_{T_1,.} \big\|_{{B}^{1-\ka}(\R ; \mathcal{B}_{\infty}^{\alpha})}\\
& \lesssim  \big\| \phi \tilde{f}_{T_1,.}\big\|_{\mathcal{C}_b^{1-\ka}(\R ; \mathcal{B}_{\infty}^{\alpha})}+ \big\| (\partial_t \phi)  \tilde{f}_{T_1,.} \big\|_{\mathcal{C}_b^{1-\ka}(\R ; \mathcal{B}_{\infty}^{\alpha})}  \quad \quad \text{(by \eqref{eqqq-cucu})}\\
& \lesssim  \big\|  \tilde{f}_{T_1,.}\big\|_{{\mathcal{C}}^{1-\ka}([T_1,T_2+1] ; \mathcal{B}_{\infty}^{\alpha})}, 
\end{align*}
where the last estimate can be proved with the same arguments as \eqref{eqqq}. The derivation of \eqref{lastpap} is now immediate: since $\tilde{f}_{T_1,T_1}=0$, one has, indeed,
$$\big\| \tilde{f}_{T_1,.}\big\|_{\cac^{1-\ka}([T_1,T_2+1]; \mathcal{B}_{\infty}^\alpha)}= \big\| \tilde{f}_{T_1,.}\big\|_{{\ov \cac}^{1-\ka}([T_1,T_2+1]; \mathcal{B}_{\infty}^\alpha)}.$$
\end{proof}


\end{document}